\documentclass[oneside, 12pt]{book}              
\usepackage{makeidx}
\usepackage{amsmath, amsfonts, amsthm, amssymb, mathrsfs, pb-diagram, graphics}

\usepackage{verbatim, xcolor, tikz, fullpage}
\usepackage{enumerate}
\usepackage{hyperref}
\hypersetup{colorlinks,%
citecolor=red,%
linkcolor=blue,%
}
\usepackage{color}
\usepackage[shortlabels]{enumitem}
\usepackage{mathtools}
\usepackage{graphicx}
\usepackage{wrapfig}
\usepackage{tikz}
\usepackage{tabularx}

\usepackage{exercise,chngcntr}

\newtheorem{thm}{Theorem}[chapter]
\newtheorem{lem}{Lemma}[chapter]
\newtheorem{cor}{Corollary}[chapter]
\newtheorem{defi}{Definition}[chapter]
\newtheorem{example}{Example}[chapter]

\newtheorem{prob}{}[chapter]
\newtheorem{conjecture}{Conjecture}
\newtheorem{question}{Question}
\newtheorem{claim}{Claim}

\newcommand{\floor}[1]{\left\lfloor #1 \right\rfloor}
\newcommand{\eps}{\varepsilon}
\newcommand{\om}{\omega}

\newcommand{\gT}{\mathbf{T}}
\newcommand{\cE}{\mathcal{E}}
\newcommand{\cR}{\mathcal{R}}
\newcommand{\One}{\mathbf{1}}

\newcommand{\notimplies}{%
  \mathrel{{\ooalign{\hidewidth$\not\phantom{=}$\hidewidth\cr$\implies$}}}}
  
\newcommand{\red}[1]{\textcolor{red}{#1}}
\newcommand{\blue}[1]{\textcolor{blue}{#1}}
\newcommand{\green}[1]{\textcolor{teal}{#1}}

\newcommand{\icol}[1]{             
  \left(\begin{smallmatrix}#1\end{smallmatrix}\right)
}

\DeclareMathOperator{\vol}{vol}
\DeclareMathOperator{\proj}{Proj}

\DeclareMathOperator{\conv}{conv}
\DeclareMathOperator{\interior}{int}
\DeclareMathOperator{\closure}{clos}

\DeclareMathOperator{\supp}{support}

\def\sinc{{\rm{sinc}}}
\def\SL_d{{\rm SL_d}}
\def\torus{{\mathbb T^d}}

\def\F{{\mathcal F}}

\def\K{{\mathcal K}}
\def\L{{\mathcal L}}

\def\P{{\mathcal P}}

\def\S{{\mathcal S}}

\def\b{{\bf b}}

\def\n{{\bf n}}

\def\x{{\bf x}}

\def\z{{\bf z}}

\def\0{{\bf 0}}
\def\1{{\bf 1}}

\def\Q{\mathbb{Q}}
\def\R{\mathbb{R}}
\def\Z{\mathbb{Z}}
\def\C{\mathbb{C}}

\def\RightTriangle
{   
\begin{tikzpicture}
\draw[blue, thick] (0,0) -- (0, 0.3) -- (0.3, 0) -- cycle;
\end{tikzpicture}
}

\def\rt
{   
\begin{tikzpicture}
\draw[blue] (0,0) -- (0, 0.23) -- (0.23, 0) -- cycle;
\end{tikzpicture}
}

\title{
{\bf A friendly introduction}\\
{\bf  to Fourier analysis on polytopes}  \\  \bigskip \bigskip
{and the geometry of numbers}  \\ \bigskip
}
\bigskip
\bigskip

\author{Sinai Robins}          
\date{\today}

\makeindex

\begin{document} 

\maketitle               

\thispagestyle{empty}

\tableofcontents

\chapter*{Acknowledgements}
The famous saying ``no man is an island'' is doubly-true in Mathematics, 
and indeed I've had the good fortune to know and learn from many interesting people, 
concerning the contents of this book.   Special thanks goes to Ricardo Diaz, my first collaborator along these topics.
I would like to thank the following people, from the bottom of my heart, for their valuable input and interesting 
discussions about some of these topics over the years:    

Ian Alevy, Artur~Andr\'e, Christine~Bachoc, Tamar~Bar, 
Imre~B\'ar\'any,  Alexander~Barvinok, Matthias~Beck, Dori~Bejleri,
Luca~Brandolini, Michel~Brion, Sunil~Chetty, 
Henry~Cohn, Leonardo~Colzani,  Amalia~Culiuc,  Pierre~Deligne,
Jes\'us A. De Loera, Holley~Friedlander,
Michel~Faleiros,  Brett~Frankel, 
Lenny~Fukshansky, Nick~Gravin, Tom~Hagedorn,
Martin~Henk,  Didier~Henrion,   Roberto Hirata Junior,  
Jeffrey~Hoffstein, Judy~Holdener,
Alex~Iosevich, Michael Joswig, Gil Kalai,
Marvin~Knopp, Mihalis~Kolountzakis, Matthias~K\"oppe, Greg~Kuperberg, 
Jean~Bernard Lasserre,   
Nhat~Le~Quang, Rafael~Zuolo~Coppini~Lima,  Sameer~Iyer, 
Fabr\'icio~Caluza~Machado,   
Romanos~Malikiosis, M\'at\'e Matolci, 
Tyrrell~McAllister, Nathan~McNew, Paul~Melvin, 
Victor Moll, Mel~Nathanson,   
James~Pommersheim,  Jim~Propp, 
Thales~Paiva,   Jill~Pipher,    Geremias~Polanco,
 Jorge Luis Ram\'irez Alfons\'in,  Ethan~Reiner,  Bruce~Reznick, 
Tiago~Royer,     Nicolas~Salter, 
  Gerv\'asio~Santos, Richard~Schwartz,
 Dima~Shiryaev,  
   Joseph~Silverman, Richard~Stanley,    Irena~Swanson,  Stephanie~Treneer,
   Christophe~Vignat,   Sergei~Tabachnikov, 
   Karen~Taylor,
 Giancarlo~Travaglini,   Mckenzie~West,  Ian~Whitehead, Kevin~Woods,  Ren~Yi,  
 G\"unter~Ziegler,  Chuanming~Zong.

\mainmatter                          
              

\chapter{\blue{Once upon a time.....an introduction}  }            

\begin{wrapfigure}{R}{0.4\textwidth}
\centering
\includegraphics[width=0.34\textwidth]{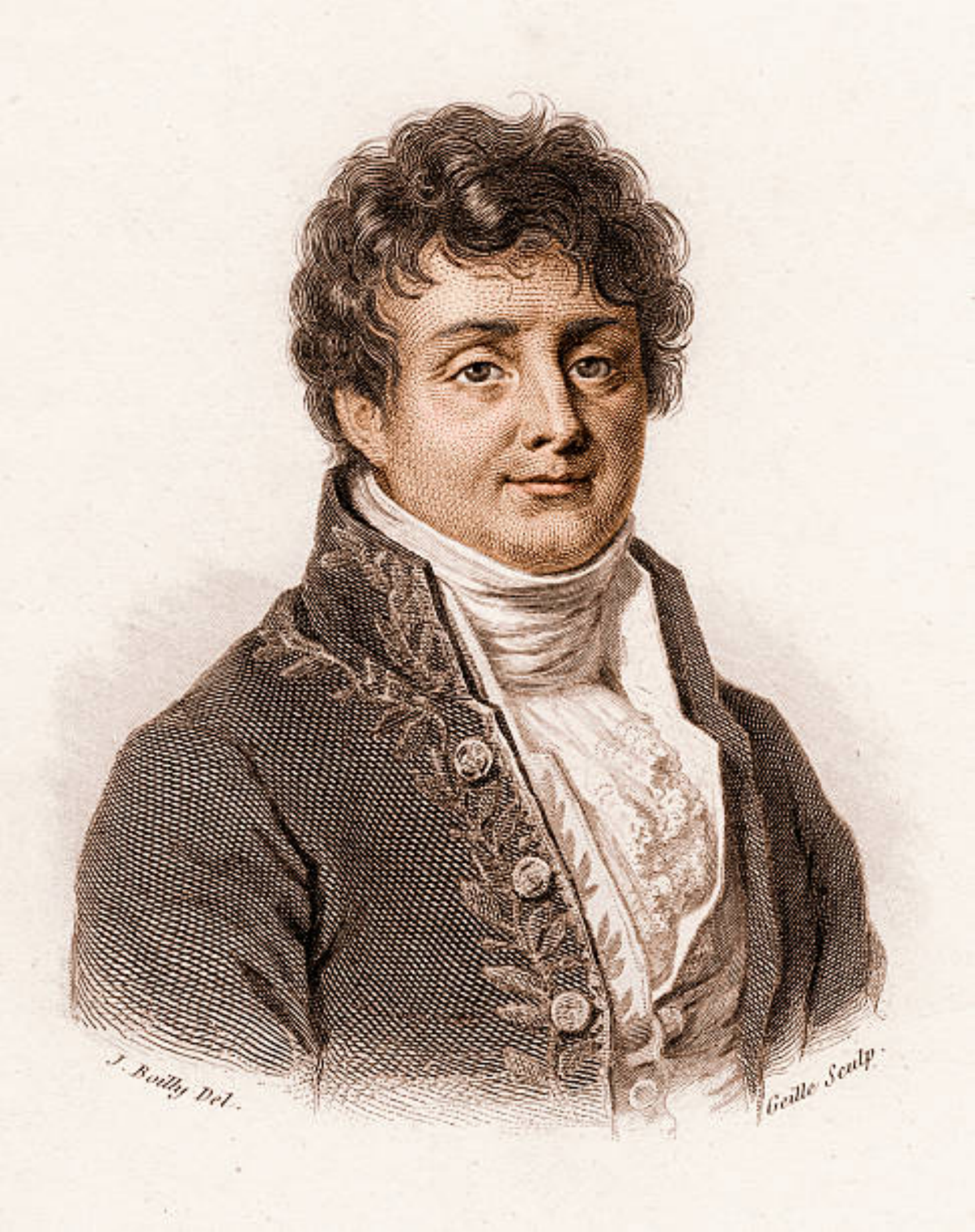}
\caption{Joseph Fourier}  \label{Fourier}
\end{wrapfigure}

What is a Fourier transform? Why is it so useful?  How can we apply Fourier transforms and Fourier series - which were originally used by Fourier to study heat diffusion - in order  to better understand topics in
 discrete and combinatorial geometry, number theory, and sampling theory?   

To begin, there are some useful analogies: imagine that you are drinking a milk-shake (lactose-free), and you want to know the ingredients of your tasty drink.   You would need to filter out the shake into some of its most basic components.   This decomposition into its basic ingredients  may be thought of as a sort of 
``Fourier transform of the milk-shake''.  Once we understand each of the ingredients, we will also be able to
 restructure these ingredients in new ways, to form many other types of tasty goodies.   To move the analogy back into mathematical language, the milkshake represents a function, and each of its basic ingredients represents for us the basis of sines and cosines; we may also think of a basic ingredient more compactly as a complex exponential $e^{2\pi i nx}$, for some $n\in \Z$.   Composing these basic ingredients together in a new way represents a Fourier series.  

\section{Introduction}
Mathematically, one of the most basic kinds of milk-shakes is the indicator function of the unit interval, and to break it down into its basic components, Mathematicians, Engineers, Computer scientists, and Physicists have used the 
sinc function (since  the $1800$'s):
\[
\sinc(z):= \frac{\sin(\pi z)}{\pi z}
\]
with great success, because it happens to be the Fourier transform of the unit interval 
$[-\frac{1}{2},  \frac{1}{2}]$:
\[
\int_{-\frac{1}{2}}^\frac{1}{2}    e^{-2\pi i z x} dx = \sinc(z),
\]
as we will compute shortly in identity \eqref{sinc function formula}. 
  Somewhat surprisingly, comparatively little energy has been given to some of its higher dimensional 
extensions, namely those extensions that arise naturally as Fourier transforms of polytopes. 

One motivation for this book is to better understand how this $1$-dimensional function -- which has proved to be extremely powerful in applications -- extends to higher dimensions.   Namely, we will build various mathematical structures that are motivated by the question:
\[
\text{ {\bf 
What is the Fourier transform of a polytope}? } 
\]
Of course, we will ask ``how can we apply it"?   An alternate title for this book might have been:

\centerline{ {\bf  We're taking Poisson summation and Fourier transforms of polytopes}} 
\centerline{ {\bf  for a very long ride....}}

Historically, sinc functions were used by Shannon (as well as Hardy, Kotelnikov, and Whittaker) when he published his seminal work on sampling theory and information theory.

In the first part of this book, we will learn how to use the technology of Fourier transforms of polytopes in order to  prove some of Minkowski's basic theorems in the geometry of numbers, to build the (Ehrhart) theory of integer point enumeration in polytopes, and to understand when a polytope tiles Euclidean space by translations.

In the second portion of this book, we give some applications to active research areas which are sometimes considered more applied, including 
the sphere-packing problem,  and the sampling of signals in higher dimensions.

There are also current research developments of the material developed here, to
 the learning of deep neural networks. 
In many applied scientific areas, in particular radio astronomy, computational tomography, and magnetic resonance imaging, a frequent theme is the reconstruction of a  function from knowledge of its Fourier transform.   Somewhat surprisingly, 
in various applications we only require very partial/sparse knowledge of its Fourier transform in order to reconstruct the required function, which may represent an image or a signal.

One of the goals here is to allow the general mathematical reader to approach the forefront of modern research in this expanding area, and even tackle some of its unsolved problems.
 There is a rapidly increasing amount of research focused in these directions in recent years, and it is therefore time to put some of these new findings in one place, making them much more  accessible to a general scientific reader.

The fact that the sinc function is indeed the Fourier transform of the $1$-dimensional line segment 
$[-\frac{1}{2}, \frac{1}{2}]$, which is a $1$-dimensional polytope, \index{polytope} 
 gives us a first hint that 
there is a deeper link between the geometry of a polytope and the analysis of its Fourier transform.

Indeed one reason that sampling and information theory, as initiated by Claude Shannon, \index{Shannon, Claude} 
works so well  is precisely because the Fourier transform of the unit interval has this nice form, and even more-so  because of the existence of the Poisson summation formula.    

The approach we take here is to gain insight into  how the Fourier transform of a polytope \index{polytope}
can be used to solve various specific problems in discrete geometry, combinatorics, optimization, approximation theory,  and the Shannon-Whittaker sampling theory in higher dimensions:
\begin{enumerate}[(a)]
\item Analyze tilings of Euclidean space by translations of a polytope 
\item   Give wonderful formulas for volumes of polytopes
\item   Compute discrete volumes of polytopes, which are combinatorial approximations to the continuous volume
\item   Introduce and develop the geometry of numbers, via Poisson summation
\item    Optimize sphere packings, and get bounds on their optimal densities
\item   Study the Shannon-Whittaker sampling theorem and its higher-dimensional siblings
\end{enumerate}

\medskip


Let's see at least one direction that quickly motivates the study of Fourier transforms.   In particular, we often begin with
simple-sounding problems that arise naturally in combinatorial enumeration, discrete and computational geometry, and number theory.

Throughout, an {\bf integer point} \index{integer point}  is any vector $v:=(v_1, \dots, v_d)\in \R^d$, all of whose coordinates $v_j$ are integers.   In other words, $v$ belongs to the integer lattice $\Z^d$. 
A {\bf rational point}   
\index{rational point}  is a point $m$ whose coordinates are rational numbers, in other words 
$m \in \Q^d$. 
We define the {\bf Fourier transform} of a function $f(x)$:
\begin{align} \label{Fourier transform}    \index{Fourier transform} 
\hat f(\xi) := \int_{\R^d} f(x) e^{-2\pi i \langle \xi, x \rangle} dx,
\end{align}
defined for all $\xi \in \R^d$ for which the latter integral converges, and where we use the standard inner product
$\langle a, b \rangle:= a_1 b_1 + \cdots + a_d b_d$.  We will also use the notation
$\F(f)$ for the Fourier transform of $f$, which is useful in some typographical contexts, for example when considering $\F^{-1}(f)$.

We introduce one of the main objects of study in this book,  the {\bf Fourier transform of a polytope} 
\index{Fourier transform of a polytope} 
$\P$, defined by:
\begin{align} \label{Fourier transform of P}
\hat 1_\P(\xi) 
 := \int_{\R^d} 1_\P(x) e^{-2\pi i \langle \xi, x \rangle} dx =  \int_{\P} e^{-2\pi i \langle \xi, x \rangle} dx,
\end{align}
where the function $1_\P(x)$ is the {\bf indicator function} of $\P$, defined by 
\[
1_\P(x):=
\begin{cases}  
1    &      \mbox{if } x\in \P \\ 
0  &        \mbox{if not}.
\end{cases}
\]
Thus, the words ``Fourier transform of a polytope $\P$'' will always mean the Fourier transform of the 
indicator function  \index{indicator function} of $\P$.

\section{The Poisson summation formula}

\begin{wrapfigure}{R}{0.4\textwidth}
\centering
\includegraphics[width=0.32\textwidth]{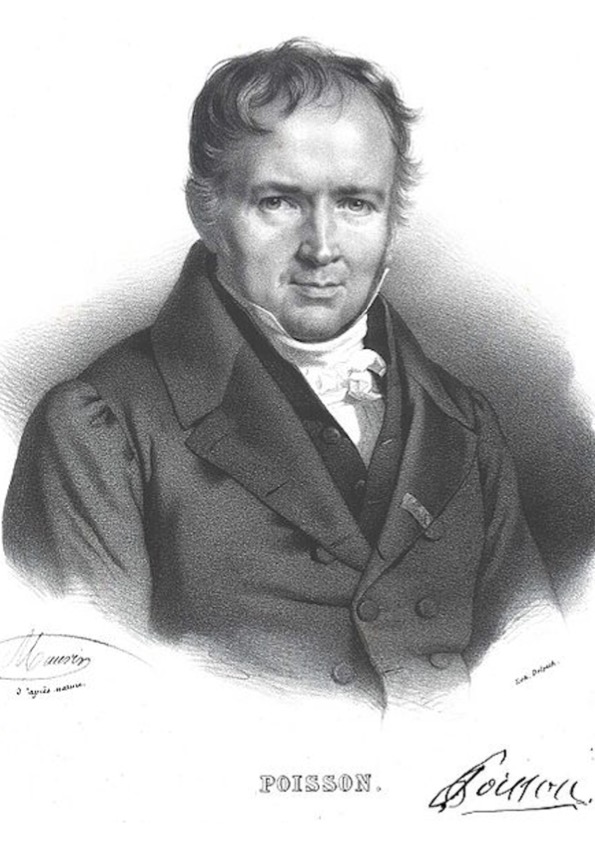}
\caption{Sim\'eon Denis Poisson}    \label{PoissonHimself}
\end{wrapfigure}

The {\bf Poisson summation formula}, named after  Sim\'eon Denis Poisson, 
  \index{Poisson summation formula}
 tells us that for any ``sufficiently nice"
 function $f : \R^d \rightarrow \C$ we have:
\begin{equation}\label{first appearance of Poisson summation}
\sum_{n \in \Z^d} f(n) = \sum_{\xi \in \Z^d} \hat f(\xi).
\end{equation}
In particular, if we were to naively set $f(n) :=  1_{\P}(n)$, the indicator function of a polytope $\P$, then we would get:
\begin{align} \label{PoissonSummation1}
\sum_{n \in \Z^d} 1_{\P}(n) = \sum_{\xi \in \Z^d} \hat 1_{\P}(\xi), 
\end{align}
which is technically false for functions, due to the fact that the indicator function $1_\P$ is discontinuous on $\R^d$.  

However,  this technically false statement is very useful!  We make this claim because it helps us build intuition for the more rigorous statements that are true, and which we study in later chapters.   
For applications to discrete geometry, we are interested in the number of integer points in a closed convex polytope $\P$, namely $|\P \cap \Z^d|$.  
The combinatorial-geometric quantity $|\P \cap \Z^d|$
may be regarded as a {\bf discrete volume}
 \index{discrete volume} 
for $\P$. 
From the definition of the indicator function of a polytope, the left-hand-side of \eqref{PoissonSummation1}
  counts the number of integer points in $\P$, namely 
we have by definition

\begin{equation}
\sum_{n \in \Z^d} 1_{\P}(n)  =  |\P \cap \Z^d|.
\end{equation}
On the other hand, the right-hand-side of  \eqref{PoissonSummation1} allows us to compute this discrete volume of $\P$ in a new way.   This is great, because it opens a wonderful window of computation for us in the following sense: 
\begin{align} \label{PoissonSummation2}
|\P \cap \Z^d| = \sum_{\xi \in \Z^d} \hat 1_{\P}(\xi).
\end{align}
We notice that for the $\xi = 0$ term, we have 
\begin{align} \label{Fourier transform at 0}
\hat 1_\P(0) := \int_{\R^d} 1_{\P}(x) e^{-2\pi i \langle 0, x \rangle} dx =   
 \int_{\P}    dx = \vol(\P),
\end{align}

and therefore the {\bf discrepancy} \index{discrepancy}
 between the continuous volume of $\P$ and the discrete volume of $\P$ is
\begin{align} \label{PoissonSummation3}
|\P \cap \Z^d| -  \vol(\P)  = \sum_{\xi \in   \Z^d-\{0\}} \hat 1_\P(\xi), 
\end{align}

showing us very quickly that indeed $|\P \cap \Z^d|$ is a discrete approximation to the classical Lebesgue volume $\vol(\P)$, and pointing us to the task of finding ways to evaluate the transform $\hat 1_P(\xi)$.  
From the trivial but often very useful identity 
\[
\hat 1_\P(0) = \vol(\P),
\]
we see another important motivation for this book:  the Fourier transform of a polytope is a very {\bf natural extension of volume}.   
\index{volume} 
Computing the volume of a polytope $\P$ captures a bit of information about $\P$, but we also lose a lot of information.  

On the other hand, computing the Fourier transform of a polytope 
$\hat 1_\P(\xi)$ uniquely determines $\P$, so we do not lose any information at all.   Another way of saying this is that the Fourier transform of a polytope is a {\bf complete invariant}.   
\index{complete invariant}
  In other words, 
it is a fact of life that 
\[
\hat 1_\P(\xi) = \hat 1_{\mathcal Q}(\xi) \text{ for all } \xi \in \R^d   \   \iff  \  \P = \mathcal Q.
\]

Combinatorially, there are brilliant identities (notably the Brion identities) that emerge between 
 the Fourier and Laplace transforms of a given polytope,  and its facets and vertex tangent cones.

In Statistics, the moment generating function of any probability distribution is given by a Fourier transform of the indicator function of the distribution, hence Fourier transforms arise very naturally in Statistical applications.
At this point, a natural glaring question naturally comes up: 
\begin{equation}
 \text{   How do we {\bf compute} the Fourier transform of a polytope }
\hat 1_P(\xi)?  
\end{equation}
And how do we use such computations to help us understand the important 
``error''  term 
\[
\sum_{\xi \in   \Z^d-\{0\}} \hat 1_\P(\xi)
\]
that came up naturally in \eqref{PoissonSummation3} above? 

There are many applications of the theory that we will build-up.   Often, we find it instructive to 
sometimes give an informal proof first, because it 
 brings the  intuitive ideas to the 
foreground, allowing the reader to gain an overview of the steps. Later on, we revisit the same intuitive proof again,  making all of the steps rigorous.

The Poisson summation formula \index{Poisson summation formula} is one of our main stars,
and some of its variations have relatively easy proofs.   But it constitutes a very first step for many of our explorations.

\section{Possible course outlines for teachers}

There is enough content here for  $2$ semesters, so it may be useful to outline some possible trajectories for a one-semester course:

{\bf [First option] }  Here the professor may follow an introduction to very basic Fourier analysis and then focus on the classical geometry of numbers, which may work well even for advanced undergraduates.
To save time, one can begin with Chapters \ref{Chapter.Tiling.A.Rectangle}, 
\ref{Chapter.Examples}, 
and in Chapter \ref{Fourier analysis basics} only cover 
Sections \ref{Fourier analysis basics}.1 - \ref{Fourier analysis basics}.15, which includes the important Poisson summation II, \'a la Poisson himself, as well as convolutions.  Then, one can cover:

Chapter \ref{Geometry of numbers}:  Sections 5.1 - 5.8 offer a novel introduction to the geometry of numbers.

Chapter \ref{chapter.lattices}:  introduces the bread-and-butter of lattices, but covering only Sections \ref{chapter.lattices}.1 - \ref{chapter.lattices}-7 may be sufficient.

Chapter \ref{Chapter.geometry of numbers II}:  this is a brief chapter that gives more theorems in the classical geometry of numbers, mainly the elegant theorems of Blichfeldt and Remak.

Chapter \ref{chapter.Brion}: basic Fourier transforms techniques are applied to recover a classical formula of Brion, for the Fourier transform of a simple polytope, via its vertex description.

\bigskip

{\bf [Second option] }  As a more advanced course, it's also possible to go more deeply into the Fourier transforms of polytopes.
I would recommend commencing with 
Chapters \ref{Chapter.Tiling.A.Rectangle}, \ref{Chapter.Examples}, and parts of 
\ref{Fourier analysis basics}, 
and a bit of 
Chapter \ref{chapter.lattices} on lattices, to give the basics of Fourier analysis, as well as sufficient intuition and motivation for what follows.   Then one can cover:

Chapter \ref{chapter.Brion}: gives some detailed and complete formulas for the Fourier transform of a polytope, given its vertex description (a formula of Brion).    

Chapter \ref{Angle polynomial}:  introduces and develops the extension of a $2$-dimensional angle (often called a solid angle) to higher dimensions.  

Chapter \ref{chapter:Discrete Brion}: gives a discretized version of Brion's theorem, for the integer point transform of a polytope, a discretized version of the Fourier transform of a polytope.  
Chapter \ref{chapter:Discrete Brion} sets the stage for the Ehrhart theory of Chapter \ref{Ehrhart theory}, which is an important theory that studies discrete volumes of polytopes.

The good news is that Chapters \ref{Ehrhart theory}, \ref{Stokes' formula and transforms}, 
\ref{Chapter.geometry of numbers III}, 
\ref{Sphere packings}, 
and \ref{Chapter:Shannon sampling} are essentially independent of each other. 
So there is a choice of ending the course with:

Chapter \ref{Stokes' formula and transforms}:  Stokes' theorem and its application to the Fourier transform of a polytope,

or  Chapter \ref{Chapter.geometry of numbers III}: more advanced topics in the geometry of numbers that use theta functions,
 
or Chapter \ref{Sphere packings}: Sphere packings, with upper bounds given by Poisson summation,

or Chapter \ref{Chapter:Shannon sampling}: Shannon's sampling theory in one and several variables, using Poisson summation.

\section{Prerequisites}

A word about {\bf prerequisites} for this book:  {\bf  Linear Algebra} is always very useful!  
A couple of calculus courses are required as well, with some real analysis.  In particular, familiarity with infinite series is assumed.   We give new proofs for
some of the main theorems in this theory,
including  Theorem \ref{Siegel for general lattices}, Theorem \ref{brion, continuous form}, Theorem \ref{brion2}, and Theorem \ref{brion, discrete form}.  Corollary \ref{THE NEW RESULT, vanishing of the FT} 
is one of the new results that appear in this book, which may also prove useful in extending the study of zero sets of the Fourier transform. 
These new Fourier-type proofs help streamline the theory, unifying sporadic results in the literature.  This
 unifying thread will hopefully help the reader put the various results - from antiquity to modernity - into context.  

We will assume some familiarity with the basic definitions of polytopes and their faces, although at places we will remind the reader of some of these definitions.
 There are many excellent texts that introduce the student to the classical language of polytopes, in particular the two classics: G\"unter Ziegler's ``Lectures on Polytopes" \cite{Ziegler},   and Branko Gr\"unbaum's ``Convex Polytopes" \cite{Grunbaum}.
For an easy introduction to the interactions between polytopes and lattice point enumeration, the reader is invited to consult ``Computing the continuous discretely: integer point enumeration in polytopes", by Beck and Robins \cite{BeckRobins}.  But the contents of the latter book are not necessary for the study of the current book. 
 
 The level of the current book is aimed at {\bf advanced undergraduates} and {\bf beginning graduate students} in various fields, and in particular Mathematics, Computer Science, Electrical Engineering, and Physics. 
 But I've included some goodies here and there for researchers as well.  Indeed, one of the goals of this book is to allow the reader to rapidly reach the forefront of research in this area.
 
 Because of the large number of exercises, with solutions to many of them in the back, this book can also be used effectively for self-study.  
If an exercise is marked with a $\clubsuit$ symbol, it means that we've mentioned this exercise in the body of the text, for that chapter. 
\red{Finally, this book is still in draft form, and in particular Chapters 10, 11, 12, and 15, 
are still under revision}. 
 
We proceed by developing an intuitive understanding first, using many examples and analogies, and this intuition  then points us to a rigorous path for the details of the ensuing proofs.   
\bigskip

Sinai Robins  \hfill February 2023

IME, University of S\~ao Paulo



 \chapter{
 \blue{
 A motivating problem: \\ tiling a rectangle with rectangles
 }
 }
 \label{Chapter.Tiling.A.Rectangle}
 \index{Tiling a rectangle}

 \begin{quote}                         
``Ripping up carpet is easy -- {\it tiling} is the issue''.

-- Douglas Wilson
 \end{quote}

\begin{wrapfigure}{R}{0.5\textwidth}
\centering
\includegraphics[width=0.45\textwidth]{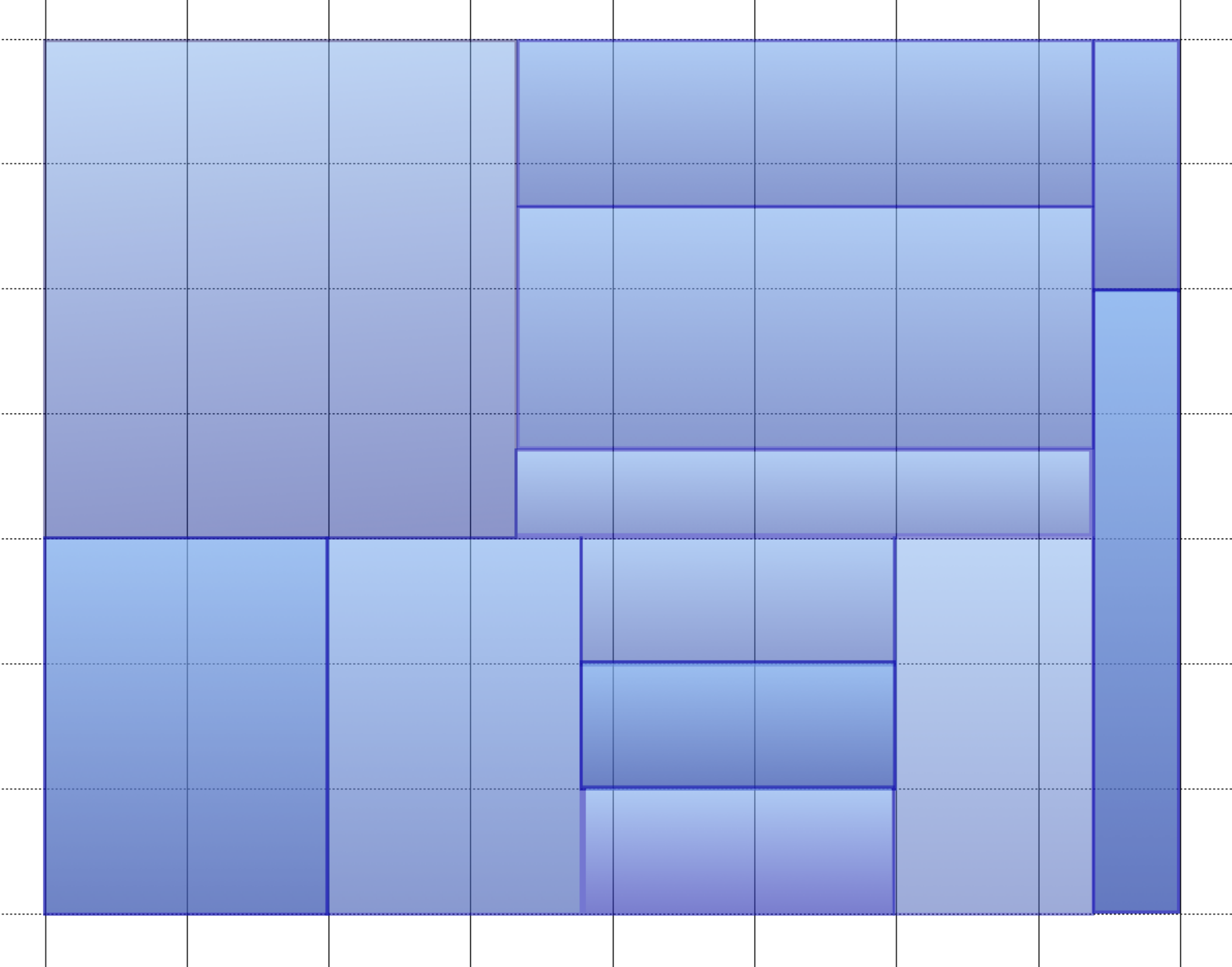}
\caption{A rectangle tiled by nice rectangles}    \label{nice rectangle}
\end{wrapfigure}

  \section{Intuition}
 
 To warm up, we begin with a simple  tiling problem in the plane.   A rectangle will be called 
  {\bf nice}  if at least one of its sides is an integer.  
We prove a now-classical fact about tiling a rectangle with nice rectangles, namely Theorem \ref{Integer.Side.Rectangle}, and we focus on the {\bf method} of the straightforward proof.

This proof brings to the foreground an important idea: by simply taking a Fourier transform of a body $B$, we immediately get interesting geometric consequences for $B$.
In particular, we will see throughout this book various ways in which the Fourier transform of a geometric body is a natural extension of its volume, sometimes in a continuous way, and sometimes in a discrete way.  
 So in order to study relationships between volumes of bodies, it is very natural and useful to play with their Fourier transforms.

\section{Nice rectangles}

The tilings that we focus on, in this small chapter, are tilings that are composed of smaller rectangles, all of which have their sides parallel to the axes, and all of which are nice.
There are at least $14$ different known proofs \cite{StanWagon} of Theorem \ref{Integer.Side.Rectangle}.
Here we give the proof that uses very basic Fourier tools, from first principles, motivating the chapters that follow.
The idea for this proof goes back to Nicolaas Govert De Bruijn \cite{DeBruijn.Book}.
\index{De Bruijn, Nicolaas Govert}

\begin{thm}[De Bruijn]  \label{Integer.Side.Rectangle}
 Suppose we tile a fixed rectangle $\cR$ with smaller, nice rectangles.   \\
 Then $\cR$ is a nice rectangle. 
 \end{thm}

\begin{proof}
Suppose that the rectangle $\cR$ is tiled with smaller rectangles $\cR_1, \dots, \cR_N$, as in
 Figure~\ref{nice rectangle}.
Due to our tiling
\index{tiling}
 hypothesis, we have
\begin{equation} \label{first identity of rectangles}
1_{\cR}(x) = \sum_{k=1}^N 1_{\cR_k}(x) + \sum (\pm \text{ indicator functions of lower-dimensional polytopes}),
\end{equation}
where the notation $1_S(x)$ always means we are using indicator functions.  
To ease the reader into the computations, we 
recall that the Fourier transform of the indicator function of any rectangle $R:=[a, b] \times [c, d]$ is defined by:
\begin{equation}
\hat 1_{\cR}(\xi) :=     \int_{\R^2}  1_{\cR}(x) e^{-2\pi i \langle \xi, x \rangle} dx 
=\int_a^b \int_c^d   e^{-2\pi i (\xi_1 x_1 + \xi_2 x_2)}dx_1 dx_2.
\end{equation}
Now we may formally take the Fourier transform of both sides of  \eqref{first identity of rectangles}.
In other words we
simply multiply both sides of \eqref{first identity of rectangles} 
 by the  exponential function  $e^{-2\pi i \langle \xi, x \rangle} $
 and then integrate both sides over $\R^2$, to get:
\begin{equation}   \label{sum.of.little.transforms}
\hat 1_{\cR}(\xi) = \sum_{k=1}^N \hat 1_{\cR_k}(\xi).
\end{equation}

In \eqref{sum.of.little.transforms}, we have used the fact that a $2$-dimensional integral over a $1$-dimensional line segment always vanishes, due to the fact that
a line segment has measure $0$ relative to the $2$-dimensional measure of the $2$-dimensional transform.
Let's compute one of these integrals, over a generic rectangle $\cR_k := [a_1, a_2] \times [b_1, b_2]$:
\begin{align} \label{transform.of.a.rectangle}
\hat 1_{\cR_k}(\xi) &:=  \int_{\R^2}   1_{\cR_k}(x)  e^{-2\pi i \langle x, \xi \rangle}  dx =  
\int_{\cR_k} e^{-2\pi i \langle x, \xi \rangle}  dx \\
&= \int_{b_1}^{b_2}   \int_{a_1}^{a_2}    e^{-2\pi i \langle x, \xi \rangle}  dx \\
&=   \int_{a_1}^{a_2} e^{-2\pi i \xi_1 x_1}  dx_1      \int_{b_1}^{b_2}    e^{-2\pi i\xi_2 x_2}  dx_2 \\
&=    \frac{    e^{-2\pi i \xi_1 a_2} - e^{-2\pi i \xi_1 a_1}   }{-2\pi i \xi_1} 
\cdot \frac{    e^{-2\pi i \xi_2 b_2} - e^{-2\pi i \xi_2 b_1}   }{-2\pi i \xi_2}\\  \label{last one}
&=   \frac{1}{(-2\pi i)^2}    \frac{     e^{-2\pi i (\xi_1 a_1 + \xi_2 b_1)}   }{\xi_1 \xi_2} 
(e^{-2\pi i \xi_1 (a_2-a_1)} - 1 ) (e^{-2\pi i \xi_2 (b_2-b_1)} -1),
\end{align}
valid for all $(\xi_1, \xi_2) \in \R^2$ except for the union of the two lines $\xi_1 = 0$ and $\xi_2 = 0$. 
Considering the latter formula for the Fourier transform of a rectangle, we make the following leap of faith:

{\bf Claim}. \ Suppose that $\cR$ is a rectangle whose sides are parallel to the axes.  Then 
\begin{equation}      \label{First.case.of.Fourier.tiling.criterion} 
\index{tiling}
\cR  \text{  is a nice rectangle  } \iff  \hat 1_{\cR}\Big(\icol{1\\1} \Big) = 0.
\end{equation}
Proof of the claim.   \ We consider the last equality \eqref{last one}.
We see that   
\begin{equation}   \label{twofactors}
\hat 1_{\cR_k}(\xi) =0   \iff  (e^{-2\pi i \xi_1 (a_2-a_1)} - 1 ) (e^{-2\pi i \xi_2 (b_2-b_1)} -1)=0,
\end{equation}
which is equivalent to having either $e^{-2\pi i \xi_1 (a_2-a_1)} =1$, or   $e^{-2\pi i \xi_2 (b_2-b_1)} =1$.
But we know that due to Euler, 
$e^{2\pi i \theta} = 1$ if and only if $\theta \in \Z$ (Exercise \ref{TrivialExponential}),  so we have
\begin{equation}       \label{the last bit}
\hat 1_{\cR}(\xi) = 0   \  \iff        \xi_1 (a_2-a_1) \in \Z     \  \text{   or   }  \   \xi_2 (b_2-b_1)   \in \Z.
\end{equation}
Now, if $\cR$ is a nice rectangle, then one of its sides is an integer, say $a_1 - a_2 \in \Z$ without loss of generality.   Therefore  
$\xi_1 (a_2-a_1) \in \Z$ for $\xi_1 = 1$,  and by \eqref{the last bit}, we see that 
$ \hat 1_{\cR}\Big(\icol{1\\1} \Big) = 0$.
 Conversely, if we assume that  $ \hat 1_{\cR}\Big(\icol{1\\1} \Big) = 0$, then by 
\eqref{the last bit}  either  \\
$1\cdot (a_2-a_1) \in \Z    \text{ or }   1\cdot (b_2-b_1)   \in \Z$,  proving the claim.

\medskip
To finish the proof of the theorem, by hypothesis each little rectangle $\cR_k$ is a nice rectangle, so by the claim above it satisfies
$ \hat 1_{\cR_k}\Big(\icol{1\\1} \Big) = 0$.
Returning to \eqref{sum.of.little.transforms}, we see that therefore
$\hat 1_{\cR}(\xi) = \sum_{k=1}^N \hat 1_{\cR_k}(\xi) = 0$, for  $\xi =  \icol{1\\1}$, and using the claim again (the converse part of it this time), we see that $\cR$ must be nice.
\end{proof}
The proof of Theorem \ref{Integer.Side.Rectangle}  was simple and elegant, motivating the use of Fourier transforms of polytopes in the ensuing chapters.   The claim, namely equation \eqref{First.case.of.Fourier.tiling.criterion},  offers an intriguing springboard for deeper investigations - it tells us that we can convert a geometric statement about tiling into a purely analytic statement about the vanishing of a certain integral transform.   Later, when we learn about
 Theorem \ref{zero set of the FT of a polytope}, we will see that this small initial success of  \eqref{First.case.of.Fourier.tiling.criterion} is part of a larger theory.  This is the beginning of a beautiful friendship.......

\bigskip
\section{Conventions, and some definitions}
We mention some conventions that we use throughout the book. 
First, we note that whenever we are given a complex-valued function $f:\R^d \rightarrow \C$,  we may write $f$ in terms of its
real and imaginary parts: $f(x):= u(x) + i v(x)$.   The  {\bf integral of such an $f$} is defined by 
\begin{equation}
\int_{\R^d}  f(x) dx := \int_{\R^d}  u(x) dx + i  \int_{\R^d}  v(x) dx,
\end{equation}
so that all of our Fourier transforms are really reduced to the usual integration of real-valued functions on Euclidean space
(see Exercise \ref{definition of complex integral}).   This is good news for the reader, because even though we see complex functions in the integrand,  elementary calculus suffices.
\medskip

Let $S\subset \R^d$ be a set.    
For our purposes, we may call $S$ a {\bf measurable} set if the integral
$
\int_S dx \text{ exists},
$
and in this case we define
\[
{\rm measure}(S) := \int_S dx.
\]
Equivalently, we may call $S$ measurable if the indicator function $1_S$ is an integrable function, by definition of the (Lebesgue) integral. 
We'll use the fact that every open set, every closed set (and hence every compact set) is measurable \cite{RudinGreenBook}.
A set $S$ is said to have {\bf measure zero} if 
\[
\int_S dx = 0.
\]
In $\R$, for example, we may alternatively define a set $S$ of measure $0$ as follows.
 Given any $\varepsilon >0$,
 there exists a countable collection of open intervals $I_n$ that cover all of $S$, and whose total length satisfies
  $\sum_{n=1}^\infty |I_n|< \varepsilon$.  But we will assume the reader knows the definition(s) of an integral 
  (either the Riemann integral or the Lebesgue integral),
circumventing discussions about $\sigma$-algebras of sets, so that the background required of the reader is kept to a minimum.

 The point we want to make here is that most things are in fact easier than the reader may have previously thought. 
 
We say that a statement $A(x)$ concerning points $x\in \R^d$ 
{\bf holds for almost every} $x\in \R^d$ (we also use the words {\bf almost everywhere})
 if the set of $x\in \R^d$ for which $A(x)$ is false is a set of measure $0$.  
For example, we have the following fact from real analysis:
\[
\int_{\R^d}  \left| f(x) - g(x) \right| dx =0  \ \iff \ f = g  \text{ almost everywhere},
\]
which means that $f(x) = g(x)$ for almost every $x\in \R^d$.  

We also mention our convention/notation for some definitions.   Whenever we want to define a new object called $N$, in terms of some combination of previously known mathematical objects called $K$, we will use the standard notation
\[
N:= K.
\]

For any set $A\subset \R^d$, we define the {\bf closure} of $A$ as the   
 the smallest (w.r.t containment) closed set that contains $A$, written as $\closure A$.  
 We define the {\bf interior} of $A$ as the set of all points $x \in A$  such that
   there exists a ball of some positive radius $\varepsilon$, centered at $x$, with $B_\varepsilon(x)  \subset A$. 
We define the {\bf boundary} of $A$, written as $\partial A$, by
\[
\partial A:=  \closure A \setminus \interior A. 
\]
An important concept is that of the support of a function $f:\R^d\rightarrow \C$, defined by
\begin{equation}  \label{def of support}
\supp(f):= \closure \{  x \in \R^d \bigm |  f(x) \not=0 \}.
\end{equation}
With this definition, we have for example:
\[
\supp(1_{[0, 1]}) = \supp(1_{(0, 1)}) = [0, 1]. 
\]
We will also say that a function $f$ is {\bf compactly supported} if the support of $f$ is a 
compact set $C$.  In particular this means that $f$ vanishes outside of $C$.


\bigskip

\section*{Notes}
\begin{enumerate}[(a)]
\item This little chapter was motivated by the lovely article written by Stan Wagon \cite{StanWagon}, which gives $14$ different proofs of Theorem \ref{Integer.Side.Rectangle}.
   The article \cite{StanWagon} is important because it shows -  in a concrete manner - how tools from one field can leak into another field, and may therefore lead to important discoveries in the future. 

\item In a related direction, we might wonder which polygons, and more generally which polytopes, tile Euclidean space by translations with a lattice.  It turns out (Theorem~\ref{zero set of the FT of a polytope}) that this question is equivalent to  the statement that the Fourier transform of $\P$ vanishes on a (dual) lattice.

\item  In the context of the Hilbert space of functions $L^2([0,1])$, Exercise   \ref{orthogonality for exponentials} 
 is one step towards showing that the set of exponentials $\{ e_n(x) \}_{n \in \Z}$ forms an orthonormal  basis for $L^2([0,1])$.  Namely, the identity above shows that these basis elements are orthogonal to each other - their inner product 
$\langle e_a, e_b \rangle := \int_0^1 e_a(x)  \overline{e_b(x)} dx $ vanishes 
for integers $a \not= b$.  Thus, the identity of Exercise  \ref{orthogonality for exponentials} 
 is often called the orthogonality relations
for exponentials, over   $L^2([0,1])$.    To show that they {\it span}  the space of functions in $L^2([0,1])$ is a bit harder, but see \cite{Travaglini} for details.

\item   The question in Exercise   \ref{Erdos lattice partition problem}  for $\Z$ was originally asked by Paul  Erd\H os
\index{Erd\"os, Paul} in $1951$, and has an affirmative answer.   This question also has
higher-dimensional analogues:
\begin{quote}
 Suppose we give a partition of the integer lattice  $\Z^d$
into a finite, disjoint union of translated sublattices.
Is it always true that at least two of these sublattices are translates of each other?
\end{quote}
The answer is known to be false for $d \geq 3$, but  is still unsolved for $d=2$ (see \cite{FeldmanProppRobins},\cite{BorodzikNguyenRobins}). 
\end{enumerate}


\bigskip

\section*{Exercises}
\addcontentsline{toc}{section}{Exercises}
\markright{Exercises}


\begin{quote}
``The game's afoot"

-- Arthur Conan Doyle (in his book Sherlock Holmes)
\end{quote}

\begin{prob}   $\clubsuit$  \label{TrivialExponential} 
Show that if $x \in \C$, then  $e^{2\pi i x} = 1$ if and only if $x \in \Z$.
\end{prob}

\medskip
\begin{prob}    \label{bound of the exponential function}
Show that $|e^z| \leq e^{|z|}$, for all complex numbers $z \in \C$.
\end{prob}

\medskip
\begin{prob}   $\clubsuit$  \label{orthogonality for exponentials}   
\index{orthogonality of exponentials in $L^2([0,1])$}
Here we prove the {\bf orthogonality relations for the exponential functions} defined by
$e_n(x) := e^{2\pi i n x}$, for each integer $n$.
Recall that the complex conjugate of any complex number $x + iy$ is defined by 
\[
\overline{x+ iy} := x - iy,
\]
so that $\overline{e^{i \theta}} := e^{-i \theta}$ for all real $\theta$. 
Prove that for all integers $a,b$:
\begin{equation}
\int_0^1 e_a(x)  \overline{e_b(x)} dx =  
\begin{cases}  
1    &      \mbox{if } a=b \\ 
0  &        \mbox{if not}.
\end{cases}
\end{equation}
\end{prob}

\medskip
\begin{prob}    \label{definition of complex integral}
Here the reader may gain some practice with the definitions of integrals that use complex-valued integrands $f(x) := u(x) + iv(x)$.
We recall for the reader the following definition:
\begin{equation} \label{real.and.imaginary.parts}
\int_{\R^d} f(x) dx := \int_{\R^d} \left( u(x) + i v(x) \right) dx :=  \int_{\R^d} u(x) dx + i  \int_{\R^d} v(x) dx,
\end{equation}
a linear combination of two real-valued integrals.  
Recalling that by definition, 
\[
\hat 1_{[0,1]}(\xi) := \int_{[0,1] }   e^{-2\pi i \xi x}  dx, 
\]
show directly from definition \ref{real.and.imaginary.parts} and from Euler's identity
$e^{i\theta} = \cos \theta + i \sin \theta$,
 that for any nonzero $\xi \in \R$, we have 
\begin{equation*}
\int_{[0,1]}   e^{-2\pi i \xi x}  dx =  \frac{e^{-2\pi i \xi } -1}{ -2\pi i \xi}.
\end{equation*}

{\rm Notes. Another way of thinking about this exercise is that it extends the `Fundamental theorem of calculus' to complex-valued functions in a rather easy way.    The anti-derivative of the integrand $f(x):= e^{-2\pi i \xi x}$ is $F(x):= \frac{e^{-2\pi i \xi x}}{-2\pi i \xi}$, and  we are saying that it is ok to use it in place of the usual anti-derivative in Calculus $1$ - it is consistent with definition \ref{real.and.imaginary.parts}. 
In the future, we generally do not have to break up complex integrals into their real and imaginary parts, 
because we can make use of the fact that antiderivatives of complex-valued functions are often simple, such as the one in this example.
 
We also note that this is {\bf not} calculus with a complex variable, because the {\bf domains of our integrands}, as well as the measures we are using throughout this book,  are defined over real Euclidean space $\R^d$.   This means we are still using basic Calculus.    
}
\end{prob}

 \begin{figure}[htb]
\begin{center}
\includegraphics[totalheight=2in]{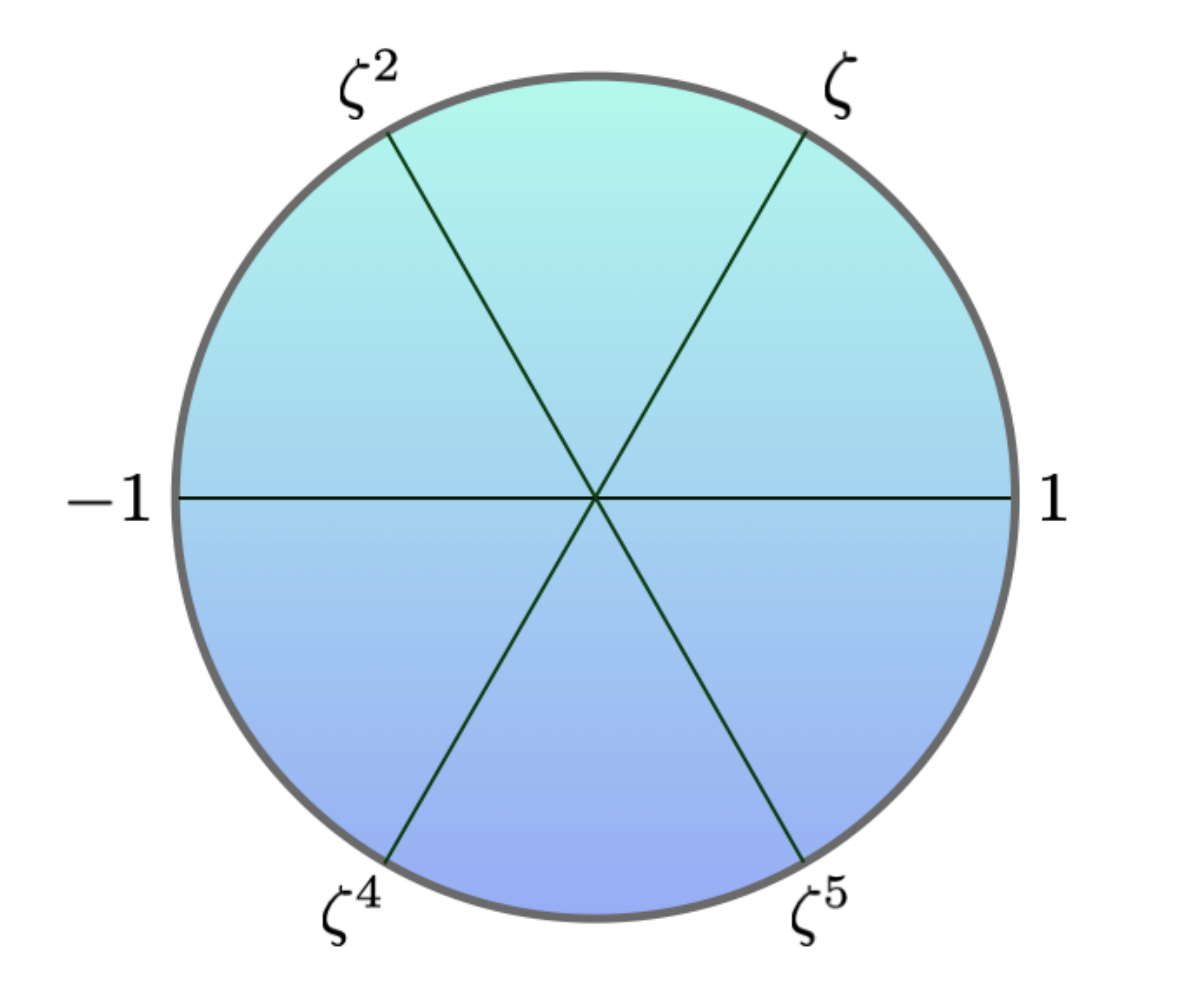}
\end{center}
\caption{The $6$'th roots of unity, with  $\zeta:= e^{\frac{2\pi i}{6}}$.  Geometrically, Exercise \ref{SumOfRootsOfUnity}  tells
us that their center of mass is the origin. }    \label{6th roots of unity}
\end{figure}

\medskip
\begin{prob}  $\clubsuit$  \label{SumOfRootsOfUnity} 
We recall that the $N$'th roots of unity are by definition the set of $N$ complex solutions to $z^N =1$, and are given by the set 
$\{e^{2\pi i k/N} \mid k = 0, 1, 2, \dots, N-1 \}$ of points on the unit circle.  Prove that the sum of all of the $N$'th roots of unity vanishes.  Precisely, fix any positive integer $N\geq 2$, 
and show that 
\[
\sum_{k = 0}^{N-1}    e^{\frac{2\pi i k}{N}} = 0.
\]
\end{prob}

\medskip
\begin{prob}    \label{DivisibilityUsingExponentials} 
Prove that, given positive integers $M, N$, we have 
\[
\frac{1}{N} \sum_{k = 0}^{N-1}    e^{\frac{2\pi i kM}{N}} = 
\begin{cases}  
1    &      \mbox{if } N \mid M \\ 
0  &        \mbox{if not}.
\end{cases}
\]
\end{prob} 

Notes.  This result is sometimes referred to as {\bf ``the harmonic detector"} for detecting when a rational
number $\frac{M}{N}$ is an integer;  that is, it assigns a value of $1$ to the sum if $\frac{M}{N} \in \Z$, and it assigns a value 
of $0$ to the sum if 
 $\frac{M}{N} \not\in \Z$.

\medskip
\begin{prob}  $\clubsuit$  \label{Orthogonality.for.roots.of.unity}   \index{orthogonality, roots of unity}
\index{root of unity}
Here we prove the {\bf orthogonality relations for roots of unity}.  Namely,  fix any two nonnegative integers $a,b$,
and prove that 
\begin{equation}  \label{12345}
\frac{1}{N} \sum_{k = 0}^{N-1}    e^{\frac{2\pi i ka}{N}} e^{-\frac{2\pi i kb}{N}} = 
\begin{cases}  
1    &      \mbox{if } a \equiv b \mod N \\ 
0  &        \mbox{if not}.
\end{cases}
\end{equation}
\end{prob} 

Notes.  In a later chapter on Euclidean lattices (Chapter \ref{chapter.lattices}), we will see that the identity
\ref{12345}  is a special case of the 
more general orthogonality relations for characters on lattices.  From this perspective, this exercise gives the orthogonality relations on the finite cyclic group $\Z/{N\Z}$.  There are more general orthogonality relations 
for characters of group representations, which play an important role in Number Theory.

\medskip
\begin{prob} \label{trick-write an integer as a product with roots of unity}
Show that for any positive integer $n$, we have
\[
n = \prod_{k=1}^{n-1} (1-\zeta^k),
\]
where $\zeta:= e^{2\pi i / n}$.
\end{prob} 


\medskip
\begin{prob}    \label{PrimitiveRootsOfUnity}  \index{root of unity, primitive}
An $N$'th root of unity is called a {\bf primitive root of unity}  if it is not a $k$'th root of unity for some smaller
positive integer $k < N$.  
Show that the primitive $N$'th roots of unity are precisely the numbers $e^{2\pi i k/N}$ for which 
$\gcd(k, N) = 1$. 
\end{prob}

\medskip
\begin{prob}    \label{SumOfPrimitiveRootsOfUnity} 
 The M\"obius $\mu$-function 
 \index{M\"obius $\mu$-function}
 is defined by:
\[
\mu(n) := \begin{cases}  
(-1)^{\text{ number of distinct prime factors of } n}     &      \mbox{if } n > 1 \mbox{ is a product of distinct primes }  \\ 
0   &        \mbox{if } n   \mbox{ is divisible by a square} \\
1   &        \mbox{if } n=1. 
\end{cases}
\]
Prove that the sum of all of the primitive $N$'th roots of unity is equal to the 
M\"obius $\mu$-function, evaluated at $N$:
\begin{equation}
\sum_{1\leq k < N \atop \gcd(k, N) = 1}   e^{\frac{2\pi i k}{N}} = \mu(N).
\end{equation}
\end{prob} 

Notes. \ See problem \ref{Mobius meets Poisson}, as a way of intertwining the M\"obius $\mu$-function with Poisson summation.

\medskip
\begin{prob}  $\clubsuit$  \label{extension of exponential}
We follow the Weierstrassian approach to defining the complex exponential $e^{z}$
 for all complex $z \in \C$:
\begin{equation}
e^{z} := \sum_{n=0}^\infty \frac{1}{n!} z^n,
\end{equation}
which converges absolutely for all $z\in \C$.    We also have the (Weierstrassian) definitions
of $\cos z$ and $\sin z$:
\[
\cos z:=  \sum_{n=0}^\infty \frac{1}{(2n)!} (-1)^n  z^{2n}, \quad  
\sin z:= \sum_{n=1}^\infty \frac{1}{(2n-1)!}  (-1)^{n-1}  z^{2n-1},
\]
 both converging absolutely again for all $z \in \C$.   Using these three Taylor series in $z$, prove that Euler's formula 
has the extension:
\[
e^{iz} = \cos z + i \sin z,
\]
valid for all $z \in \C$. 
\end{prob}

Notes.   \ Beginning with such a power series approach to many of the standard functions, Karl Weierstrass developed a rigorous and beautiful theory of real and complex functions.

\medskip
\begin{prob}
Here the reader needs to know a little bit about the quotient of two groups (this is one of the few exercises that assumes group theory).  We
 prove that the group of `real numbers mod $1$'  under addition, is
  isomorphic to the unit circle, under multiplication of complex numbers.  Precisely, we can define
$h: \R \rightarrow S^1$ by $h(x) := e^{2\pi i x}$.
\begin{enumerate}[(a)]
\item We recall the definition of the kernel of a map, namely $ker(h):= \{ x\in \R \mid h(x) = 1\}$.  Show that $ker(h) = \Z$.
\item Using the first isomorphism Theorem for groups, show that 
$\R/\Z$ is isomorphic to the unit circle $S^1$.
 \end{enumerate}
\end{prob}

\medskip
\begin{prob} 
Using gymnastics with roots of unity, we recall here a very classical solution to the problem of finding the roots of a cubic polynomial.
\begin{enumerate}[(a)]
\item    Let $\omega:= e^{2\pi i/3}$, and show that we have the polynomial identity:
 \[
 (x + a + b)( x + \omega a + \omega^2 b)   (x + \omega^2 a + \omega b)  
= x^3 - 3abx + a^3 + b^3.
\]
\item    Using the latter identity, solve the cubic polynomial:
$x^3 - px + q = 0$ by substituting $p = 3ab$ and $q= a^3 + b^3$. 
\end{enumerate}
\end{prob}

\medskip
\begin{prob}\label{zeros of the sin function}
Thinking of the function $\sin(\pi z)$ as a function of a complex variable $z\in \mathbb C$, show that its zeros are precisely the set of integers $\Z$.  
\end{prob}

\medskip
\begin{prob}
Here we give another equivalent condition for a rectangle in Theorem  \ref{Integer.Side.Rectangle}  to be a nice rectangle, using the same
definitions as before.   

Let's call $\xi \in \Z^2$ a {\bf generic} integer point
if $\xi$ is not orthogonal to any of the edges of $\cR$. In other words, a generic integer vector satisfies
$\langle \xi, p \rangle \not=0$, for all $p\in \cR$, and in particular $p=0$ is not generic, nor is any
point $p$ on the $x$-axis or the $y$-axis. 
Then 
\begin{equation}     
\cR  \text{  is a nice rectangle  } \iff    \hat 1_{\cR}(\xi) = 0,   \text{  for all generic points } \xi \in \Z^2. 
\end{equation}
\end{prob}

\medskip
\begin{prob}[ Erd\H os, 1951] \label{Erdos lattice partition problem}
\rm{
 Erd\H os asked: ``Can the set $\Z_{>0}$ of all positive integers be partitioned (that is, written as a disjoint union) into a finite number of
 arithmetic progressions, such that no two of the arithmetic progressions  will have the same common difference?'' 
 
Precisely, suppose that we have 
\begin{equation}\label{disjoint union of arithmetic progressions}
\Z= \{ a_1 n + b_1  \mid n \in \Z\}   \cup   \cdots  \cup   \{ a_N n + b_N  \mid n \in \Z\},
\end{equation}
for some positive integers $a_1 \leq a_2 \leq \cdots \leq a_N$, and $N\geq 2$, and where the arithmetic progressions are pairwise disjoint.

Prove that in any such partitioning of the integers,   $a_N = a_{N-1}$ (that is, the largest common difference must appear at least twice).
  
Notes.  For example, if we write $\Z = \{ 4 n + 1  \mid n \in \Z\} \cup  \{ 2 n   \mid n \in \Z\} \cup  \{ 4 n + 3  \mid n \in \Z\}$,  a disjoint union of
$3$ arithmetic progressions, then we see
that the largest common difference of $4$ appears twice.   Erd\H os  noticed that such a phenomenon must always occur. 
(See also Exercise \ref{Extension of Erdos to dimension d} for an extension to lattices in $\R^d$).
}
\end{prob}

\medskip
\begin{prob}
\rm{Continuing with the ideas of Exercise \ref{Erdos lattice partition problem}, suppose we are given a 
disjoint union of arithmetic progressions such as \eqref{disjoint union of arithmetic progressions} above.  
\begin{enumerate}[(a)]
\item  Prove that:
\begin{equation}\label{density identity}
1 = \frac{1}{a_1} + \cdots + \frac{1}{a_N}.
\end{equation}
\item  Show further that  $\gcd(a_i, a_j)>1$ for all indices $i, j$.
\item
Conversely, suppose that we are given positive integers 
$a_1, \cdots, a_N$, with $\gcd(a_i, a_j)>1$ for all indices $i, j$, and such that 
 $1= \tfrac{1}{a_1} + \cdots + \tfrac{1}{a_N}$.
 Can we always find integers $b_1, \dots, b_N$, such that we have the following disjoint union of arithmetic progressions:
$\Z=      \{ a_1 n + b_1  \mid n \in \Z\}   \cup   \cdots  \cup   \{ a_N n + b_N  \mid n \in \Z\}$ ?
(as in  \eqref{disjoint union of arithmetic progressions}) 
\end{enumerate}
}
\end{prob}


\chapter{
\blue{
Examples nourish the theory
}
}

\begin{quote}   
``To many, mathematics is a collection of theorems. For me, mathematics is a collection of examples; a theorem is a statement about a collection of examples and the purpose of proving theorems is to classify and explain the examples...''

-- John B. Conway
 \end{quote}

\label{Chapter.Examples} 
        
\begin{quote}                         
``A pint of example is worth a gallon of advice.'' 

-- Anonymous
 \end{quote}

\smallskip




\begin{figure}[htb]
\begin{center}
\includegraphics[totalheight=2.5in]{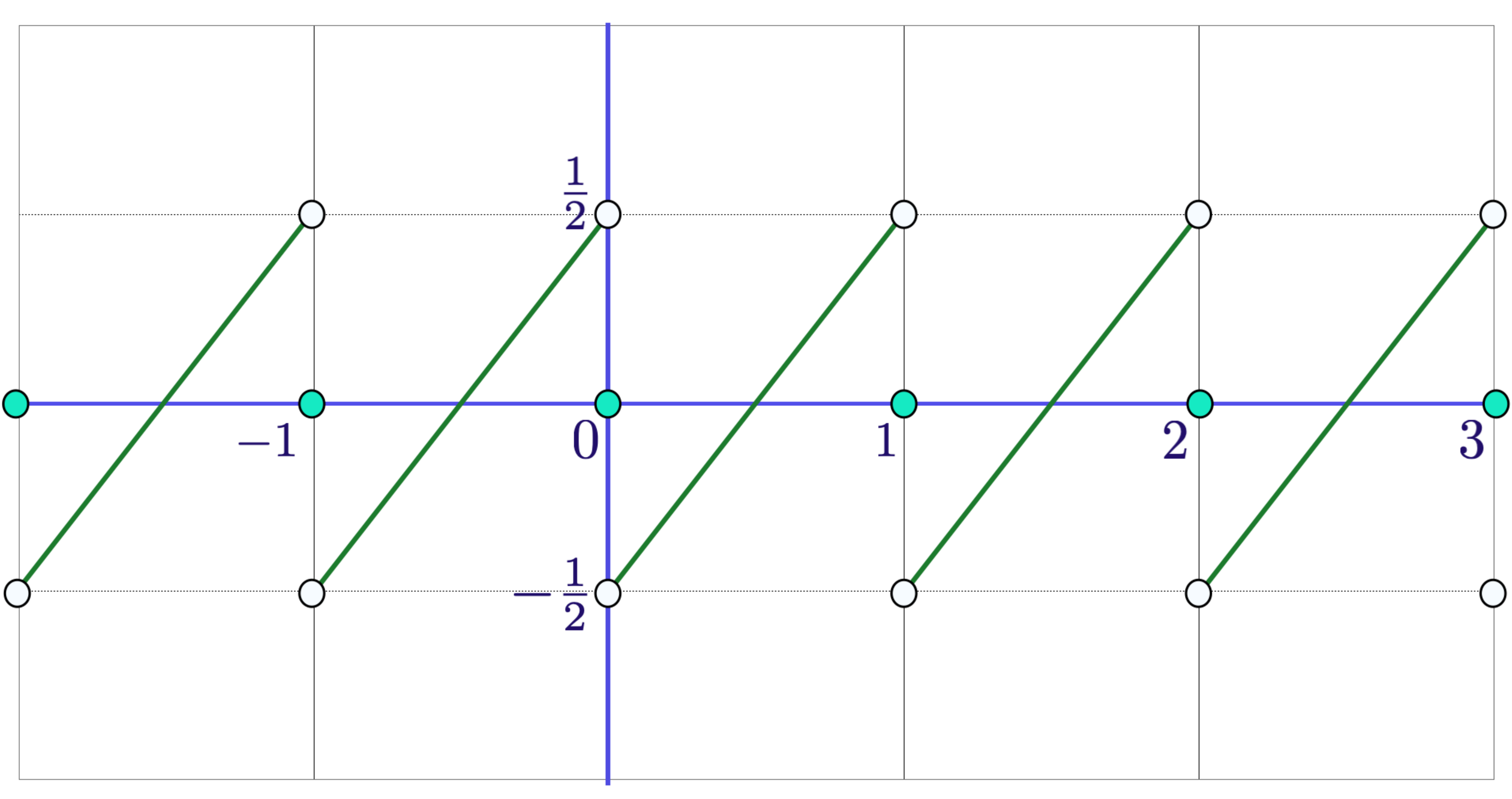}
\end{center}
\caption{The first periodic Bernoulli polynomial $P_1(x)$, sometimes called the sawtooth function, 
which turns out to be one of the building blocks of integer point enumeration in polytopes }
   \label{FirstBernoulli}
\end{figure}

\section{Intuition}
One way to think about the Fourier transform of a polytope $\P \subset \R^d$ is that it simultaneously captures all of the moments of $\P$, thereby uniquely defining $\P$.  Here we begin concretely by computing some Fourier transforms of various polytopes in dimensions $1$ and $2$, as well as the Fourier transforms of some simple families of polytopes in dimension $d$ as well.  

The $2$-dimensional computations will get the reader more comfortable with the basics.   In later chapters, once we learn a little more theory, we will return to these families of polytopes and compute some of their Fourier transforms in general.    

We also see, from small examples, that the Bernoulli polynomials immediately enter into the picture, forming natural building blocks.   
In this chapter we compute Fourier transforms without thinking too much about convergence issues, to let the reader run with the ideas.   
But commencing with the next chapter, we will be more rigorous when using Poisson summation, and with convergence issues.

\section{Dimension $1$ - the classical sinc function}    

\bigskip

We begin by computing the classical $1$-dimensional example of the Fourier transform \index{Fourier transform} 
 of the symmetrized unit interval $\P:= [-\frac{1}{2}, \frac{1}{2}]$:
 \begin{align*} \label{ClassicalExample}
\hat 1_{\P}(\xi)  := \int_{\R}  1_\P(x)  \ e^{-2\pi i x \xi } dx  
=  \int_{[-\frac{1}{2}, \frac{1}{2}]}  e^{-2\pi i x \xi } dx.
\end{align*}
 For all $\xi\not=0$, we have:
 \begin{align}
 \int_{[-\frac{1}{2}, \frac{1}{2}]}  e^{-2\pi i x \xi } dx
 &= \frac{e^{-2\pi i \left( \frac{1}{2} \right) \xi}  - e^{-2\pi i  \left( \frac{-1}{2} \xi\right) }    }{-2\pi i \xi} \\ \label{sinc}
& = \frac{ \cos (-\pi \xi) + i \sin(-\pi \xi)   -  (\cos(\pi \xi) + i \sin(\pi \xi))    }{-2\pi i \xi} \\  \label{sinc function formula}
&= \frac{\sin(\pi \xi)}{\pi \xi}.
\end{align}
Noticing that $\xi = 0$ is a removable singularity, we define the continuous {\bf sinc-function} \index{Sinc function} 
by
\begin{equation}\label{SincFunction}
\sinc(x):=  \begin{cases}  
\frac{\sin(\pi x)}{\pi x},     &\mbox{if } x \not= 0 \\ 
1  & \mbox{if } x= 0,
\end{cases}
\end{equation}
which is in fact infinitely smooth, via Lemma \ref{FT of a polytope is entire} below.

\begin{figure}[htb]
\begin{center}
\includegraphics[totalheight=1.6in]{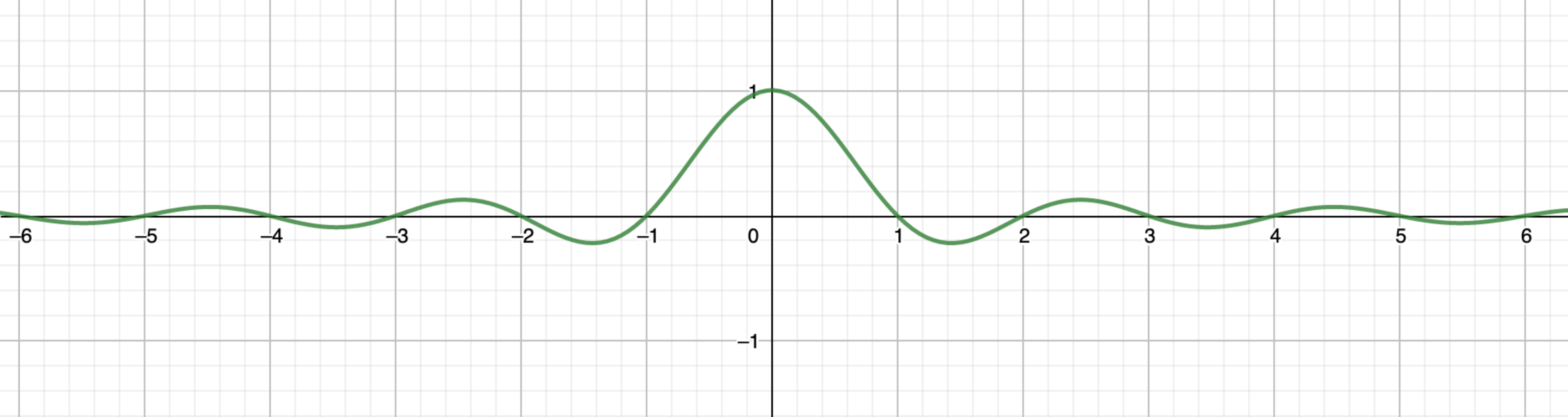}
\end{center}
\caption{The function $\sinc(x)$, which is Fourier transform of the $1$-dimensional 
polytope  $\P = [-\frac{1}{2}, \frac{1}{2}]$.  }
   \label{sinc.pic} \index{sinc function}
\end{figure}

 \bigskip
 \section{The Fourier transform of $\P$ as a complete invariant}  
\label{Fourier inversion}

The main goal of this section is to state Lemma \ref{complete invariance of the FT}, which tells us that all of the information
about a polytope is contained in its Fourier transform.  To that end, we 
 introduce 
the inverse Fourier transform, \index{inverse Fourier transform} often called the {\bf Fourier inversion formula}.
We'd like to see the fundamental fact that under certain conditions, 
the Fourier transform is invertible. First, we call a function $f:\R^d \rightarrow \C$ {\bf absolutely integrable} if 
$\int_{\R^d} | f(x) | dx < \infty$, and we write this as $f \in L^1(\R^d)$.
\begin{thm}  \label{thm:Inverse Fourier transform}
Given  a function $f$ such that both $f \in L^1(\R^d)$ and $\hat f \in L^1(\R^d)$, we have 
\begin{equation}\label{second version of Fourier inversion}
f(x) = \int_{\R^d}  \hat f(\xi)  e^{2\pi i \langle \xi, x \rangle} d\xi,
\end{equation}
for all $x \in \R^d$.
\hfill $\square$
\end{thm}
(see \cite{EinsiedlerWardBook} for a proof). 
We'll also use the notation $\F(f) := \hat f$.  Equation \eqref{second version of Fourier inversion} tells us that the inverse Fourier transform $\F^{-1}$ exists, and is almost equal to $\F$ itself.   A moment's thought reveals that we may rewrite \eqref{second version of Fourier inversion}  in the following useful form:
\begin{equation}   \label{InverseFourierTransform2}
(\F     \circ  \F)f (x)  = f(-x).
\end{equation}

\bigskip
\begin{example}\label{Integral.of.sinc}
\rm{
A famous and historically somewhat tricky integral formula for the sinc function is the following fact:
\begin{equation}  \label{area under sinc}
\int_{-\infty}^\infty   \sinc(x) dx := \int_{-\infty}^\infty    \frac{\sin(\pi x)}{\pi x} dx = 1,
\end{equation}
also known as the Dirichlet integral.  \index{Dirichlet integral}
The careful reader might notice that the latter integrand is not absolutely convergent, which means that 
$\int_{-\infty}^\infty    \Big| \frac{\sin(\pi x)}{\pi x} \Big| dx = \infty$ (Exercise \ref{divergence of |sinc|}).  So we have to specify what we really mean 
by the identity \eqref{area under sinc}.   The rigorous claim is:   
\begin{equation*}
\lim_{N\rightarrow \infty}
\int_{0}^N    \frac{\sin(\pi x)}{\pi x} dx =\frac{1}{2}.
\end{equation*}
Let's see an intuitive derivation of \eqref{area under sinc}, where we will be fast-and-loose for the moment.   Using  \eqref{ClassicalExample}, 
we've seen above that the Fourier transform of the indicator function of the interval $\P := [-\frac{1}{2}, \frac{1}{2}]$ is:
\begin{equation}
\F(1_{\P})(\xi) =  \frac{\sin(\pi \xi)}{\pi \xi},
\end{equation}
so that 
\begin{equation}
 \F \left(  \frac{\sin(\pi \xi)}{\pi \xi}        \right) 
= (\F \circ \F)(1_{\P})(\xi) 
= 1_{    \P        }(-\xi).
\end{equation}
Using the definition of the Fourier transform, the latter identity is:   
\begin{equation}
\int_{\R}        \frac{\sin(\pi x)}{\pi x}      e^{-2\pi i \xi x} dx 
= 1_{\P}(\xi),
\end{equation}
and now evaluating both sides at $\xi = 0$ gives us \eqref{area under sinc}.  
Although this derivation appears very convincing, it would not make it past the rigor police 
(see also note \ref{Ogsood's book}).
  So why not?  It is because we applied the Fourier inversion formula 
to a function that was {\bf not} in $L^1(\R)$, namely the sinc function.   So we owe it to ourselves to pursue a rigorous
approach by showing that 
\begin{equation}
\lim_{N\rightarrow \infty} \int_{-N}^N   \frac{\sin(\pi \xi)}{\pi \xi}  e^{-2\pi i \langle \xi, x \rangle} d\xi = 1_ {[-\frac{1}{2}, \frac{1}{2}]}(x),
\end{equation}
whose validity would give us a variation on Fourier inversion,  for a function that is not in $L^1(\R)$, namely  $\hat 1_ {[-\frac{1}{2}, \frac{1}{2}]}(\xi)= \sinc(\xi)$.
This is tricky business, but such an endeavor is taken up in Exercise \ref{rigorous inversion formula for sinc}.
}
\hfill   $\square$
\end{example}

We can extend Example \ref{Integral.of.sinc} in a natural way to all Fourier pairs of functions, $\{f(x), \hat f(\xi)\}$, provided 
that we may apply Fourier inversion, as follows. 
Simply let $x=0$ in \eqref{second version of Fourier inversion}, to get:
\begin{equation}  \label{IntegralTrick}
f(0)=\int_{\R^d} \hat f(x) dx.
\end{equation}
To summarize, Example \ref{Integral.of.sinc} is simply identity \eqref{IntegralTrick} with 
$f(x) := 1_{[-\frac{1}{2}, \frac{1}{2}]}(x)$.

Another nice and very useful  fact about the Fourier transform
 of a polytope is that it is an entire function,  meaning that it is differentiable everywhere.   This differentiability
 is already observable in the sinc function above, with its removable singularity at the origin.
 
\begin{lem} \label{FT of a polytope is entire}
Let $\P \subset \R^d$ be a $d$-dimensional polytope.  Then $\hat 1_\P(\xi)$ is an entire function of $\xi \in \C^d$.  
\end{lem}
\begin{proof}
Because $\P$ is compact, we can safely differentiate under the integral sign (this is a special case of Lebesgue's Dominated Convergence Theorem). Namely, for the coordinate variable $\xi_1$, we have:
$\frac{d}{d\xi_1} \int_{\P} e^{-2\pi i \langle \xi, x \rangle} dx=  \int_{\P} \frac{d}{d\xi_1}e^{-2\pi i \langle \xi, x \rangle} dx
= 2\pi i \int_{\P} x_1 e^{-2\pi i \langle \xi, x \rangle} dx$.  Since one complex derivative of  $\hat 1_\P(\xi)$ now exists  (in the complex variable $\xi_1$), 
the function $\hat 1_\P(\xi)$ is analytic in $\xi_1$, and using the same reasoning it is also analytic in each of the  variables $\xi_2, \xi_3, \dots, \xi_d$.
\end{proof}

We also have the very fortuitous fact that the Fourier transform of any polytope  $\P \subset \R^d$ is a complete invariant, in the following sense.
We recall that  by definition a polytope is in particular a closed set. 
\begin{lem} \label{complete invariance of the FT}
Let   $\P\subset \R^d$ be a polytope. Then
$\hat 1_\P(\xi)$ uniquely determines $\P$.  Precisely, given any two  $d$-dimensional polytopes $P, Q\subset \R^d$, 
we have
\[
\hat 1_\P(\xi) = \hat 1_{Q}(\xi) \text{ for all } \xi \in \R^d   \   \iff  \  \P = Q.
\]
In other words, for any polytope $\P$, its Fourier transform $\hat 1_\P$ uniquely determines the polytope. 
\end{lem} \label{FT.complete invariant}
\begin{proof} (outline)
If $\P=Q$, it is clear that $\hat 1_\P(\xi) = \hat 1_{Q}(\xi)$ for all $\xi \in \R^d$.  Conversely, suppose that 
$\hat 1_\P(\xi) = \hat 1_{Q}(\xi)$ for all $\xi \in \R^d$.  Using Fourier inversion, 
namely Theorem \ref{thm:Inverse Fourier transform},  
 we may take the Fourier transform of both sides of the latter equation to get $1_\P(-\xi) = 1_Q(-\xi)$, for all $\xi \in \R^d$. 
\end{proof}
The reason that the proof above is only an outline - at this point - 
 is due to the fact that we have applied the Fourier inversion formula \eqref{second version of Fourier inversion}
 to $\hat 1_\P$, which is not an absolutely integrable function 
(as we'll see in  Corollary \ref{the FT of an indicator function is not in L^1} in even greater generality).

We'll revisit Lemma \ref{complete invariance of the FT} in Chapter \ref{Fourier analysis basics}, as Theorem \ref{the FT of a convex set determines the set}  for a rigorous proof.  In fact much more is true - see note \ref{identity thm for two compact sets whose FT's agree on a convergent sequence}.
There is also a nice version of the Fourier inversion formula, due to Podkorytov and Minh, which is related \cite{Podkorytov}.  
The reason we've put Lemma \ref{complete invariance of the FT}
so early in the text is because it offers an extremely strong motivation for the study of Fourier transforms of polytopes, showing 
that they are complete invariants.

A fascinating consequence of Lemma \ref{complete invariance of the FT} is that when we take the Fourier transform of a polytope, then {\bf all of the combinatorial and geometric information} of $\P$ is contained in the formula of its Fourier transform......somehow.  
So we may begin to create a complete dictionary between the geometry and combinatorics of a polytope in the space domain, and its Fourier transform in the frequency domain.

\bigskip
\section{Bernoulli polynomials}   \index{Bernoulli polynomial}

We introduce the Bernoulli polynomials, which turn out to be a sort of ``glue'' between
discrete geometry, number theory, and Fourier analysis, as we will see throughout the book.
Historically, Jacob Bernoulli was considering the formulas
\[
1+ 2 + \cdots + n = \frac{n(n+1)}{2}, 
\]
\[
1^2+ 2^2 + \cdots + n^2 = \frac{n(n+1)(2n+1)}{6}, 
\]
\[
1^3+ 2^3 + \cdots + n^3 = \frac{n^2(n+1)^2}{4}, 
\]
and so on.   Jacob was wondering how to find a general formula for the sums:
\[
1^d+ 2^d + \cdots + n^d = \text{ some polynomial in the variable } n?
\]
With hindsight giving us slightly better vision, the modern approach to the latter polynomials
begins with  the following generating function:
\begin{equation}  \label{generating function for Bernoulli polynomials}
\frac{te^{xt}}{e^t-1}   =  \sum_{k =0}^\infty   B_k(x)  \frac{t^k}{k!}.
\end{equation}
It follows from this definition \eqref{generating function for Bernoulli polynomials} that each coefficient $B_k(x)$ is a polynomial in $x$, of degree $k$ 
(Exercise \ref{brute force Bernoulli polys}). 
These polynomials $B_k(x)$ are called 
{\bf Bernoulli polynomials},   \index{Bernoulli polynomial}
and Bernoulli was able to show that in general:
\[
\sum_{k=0}^{n-1} k^{d-1} = \frac{    B_d(n) - B_d(0)  }{    d   },
\]
for all integers $d \geq 1$ and $n \geq 2$ (Exercise \ref{historical origin of Bernoulli poly}).
The reader can develop her skills by proving some of the surprising and important properties of Bernoulli polynomials in Exercises \ref{brute force Bernoulli polys} through \ref{asymptotics for Beroulli numbers}.

\begin{example} \label{first few Bernoulli polys}
\rm{
 The first few Bernoulli polynomials are:
\begin{align}
B_0(x) &= 1 \\
B_1(x) &= x - \frac{1}{2} \\
B_2(x) &= x^2 - x + \frac{1}{6} \\
B_3(x) &= x^3 - \frac{3}{2} x^2 +\frac{1}{2} x  \\
B_4(x) &= x^4 -2x^3 + x^2 - \frac{1}{30}    \\
B_5(x) &= x^5 - \frac{5}{2} x^4 + \frac{5}{3} x^3 - \frac{1}{6} x   \\
B_6(x) &= x^6 - 3x^5 +  \frac{5}{2} x^4 - \frac{1}{2} x^2 + \frac{1}{42}    \label{B_6}
\end{align}
}
 \hfill $\square$
 \end{example}

It turns out that it's very useful to periodize the Bernoulli polynomials, in the following sense. 
We first define:
\[
\{x\} := x - \lfloor x  \rfloor,
\]
the fractional part of $x$.  \index{fractional part}
Now we define the $n$'th {\bf periodic Bernoulli polynomial}:  \index{periodic Bernoulli polynomial}
\begin{equation}  \label{definition of periodic Bernoulli polys}
P_n(x) := B_n(\{x\}),
\end{equation}
for $n\geq 2$.      Since $P_n(x)$  is periodic on $\R$ with period $1$, 
it has a Fourier series, and it turns out that
\begin{equation}
P_n(x) = -\frac{n!}{(2\pi i)^n} \sum_{k \in \Z - \{0\}}   \frac{e^{2\pi i k x}}{k^n},
\end{equation}
valid for $x \in \R$  (Exercise \ref{Bernoulli Polynomials}).
When $n=1$, we have the first Bernoulli polynomial 
\[
P_1(x):=
\begin{cases}
 x - \lfloor x  \rfloor - \frac{1}{2} & \text{ if } x \notin \Z, \\
 0 & \text{ if } x \in \Z.
\end{cases}
\]
which
 is very special (see Figure \ref{FirstBernoulli}).   For one thing, 
$P_1(x)$ is the only periodic Bernoulli polynomial that is not  continuous everywhere, 
and we note that its Fourier series does not converge absolutely, although it is quite appealing:
\begin{equation} \label{FirstBernoulliPolynomial}
P_1(x) = -\frac{1}{2\pi i} \sum_{k \in \Z - \{0\}}  \frac{e^{2\pi i k x}}{k},
\end{equation}
valid for all $x \notin \Z$.    But how are we supposed to sum up a conditionally convergent series such as \eqref{FirstBernoulliPolynomial}?  A common way to define it
rigorously is to prove that
\[
\lim_{N\rightarrow \infty}
 \sum_{k =1}^N  \frac{e^{2\pi i k x}}{k} \ \text{exists}.
\]
As we can see, special care must be taken  with $P_1(x)$, and
Exercise \ref{rigorous convergence of P_1(x)} provides a rigorous proof 
of the convergence of \eqref{FirstBernoulliPolynomial}.  
The {\bf Bernoulli numbers} are defined to be the constant terms of the Bernoulli polynomials:
\[
B_k := B_k(0). 
\]
Perusing Example \ref{first few Bernoulli polys}, we see that the 
first few Bernoulli numbers are:  \\
\[
B_0 = 1,  \ B_1 = -\frac{1}{2},  \ B_2 =  \frac{1}{6}, \ B_3 = 0, \ B_4 = - \frac{1}{30}, \ B_5 = 0, \ B_6 = \frac{1}{42}.
\]

It follows quickly from definition  \ref{generating function for Bernoulli polynomials}
 above that for odd $k \geq 3$, $B_k = 0$ (Exercise \ref{odd Bernoulli numbers}).   
 Using  the generating function \ref{generating function for Bernoulli polynomials}, 
 the Bernoulli numbers are defined via
 \begin{equation} \label{Def. of Bernoulli numbers}
 \frac{t}{e^t-1}   =  \sum_{k =0}^\infty   B_k \frac{t^k}{k!}.
 \end{equation}
 
An interesting identity that allows us to compute the Bernoulli numbers recursively  is:
\[
\sum_{k=0}^n  {n+1 \choose k}B_k = 0,
\]
valid for all $n \geq 1$ (Exercise \ref{vanishing identity for Beroulli numbers}). 
Some of the most natural, and beautiful, Fourier series arise naturally from the periodized Bernoulli polynomials.  

Recalling the statement of Poisson summation \eqref{first appearance of Poisson summation} from the Introduction, we now give a fast-and-loose application in dimension $d=1$. 
The following intuitive application of the Poisson summation formula already suggests an initial connection between periodized Bernoulli polynomials and Fourier transforms of polytopes.

\begin{example} [Intuitive Poisson summation] \label{intuitive Poisson example}
\index{Poisson summation formula}
 \rm{
 In this example we allow ourselves to be completely intuitive, and unrigorous at this moment, 
 but often such arguments are useful in pointing us to their rigorous counterparts. 
Consider the $1$-dimensional polytope $\P:= [a,b]$, and restrict attention to the case 
of $a, b \not\in \Z$.   If we could use the Poisson summation formula
\[
\sum_{n \in \Z^d} f(n) = \sum_{\xi \in \Z^d} \hat f(\xi),
\]
applied to the function $f(x):= 1_\P(x)$, then we would get:
\begin{align*}
\sum_{n \in \Z} 1_\P(n) &``=\text{''} \sum_{\xi \in \Z} \hat 1_\P(\xi)\\
&``=\text{''}  \ 
\hat 1_\P(0)+\sum_{\xi \in \Z - \{0\}}   \frac{e^{-2\pi i \xi b}    -  e^{-2\pi i \xi a}    }{-2\pi i \xi} \\
&``=\text{''} \ 
 (b-a) -\frac{1}{2\pi i}   \sum_{\xi \in \Z - \{0\}}   \frac{e^{-2\pi i \xi b}}{\xi}    
      +  \frac{1}{2\pi i}  \sum_{\xi \in \Z-\{0\}}  \frac{e^{-2\pi i \xi a}}{\xi} \\
&``=\text{''} \ 
 (b-a) +    \frac{1}{2\pi i}   \sum_{\xi \in \Z - \{0\}}   \frac{e^{2\pi i \xi b}}{\xi}    
      -  \frac{1}{2\pi i}  \sum_{\xi \in \Z-\{0\}}  \frac{e^{2\pi i \xi a}}{\xi} \\
&``=\text{''} \ 
 (b-a) -\left(     \{b\}- \frac{1}{2}     \right)     + \left(     \{a\} - \frac{1}{2}    \right)    \\
&``=\text{''} \ 
  b -   \{b\} - ( a - \{a\} ) = \lfloor b  \rfloor  - \lfloor a  \rfloor.
\end{align*}
We've used quotation marks around the latter string of equalities because the sums are formally divergent.  But we already know how to evaluate the left-hand side of Poisson summation above, namely
$\sum_{n \in \Z} 1_\P(n) = \#\left\{ \Z \cap \P   \right\} = \lfloor b  \rfloor  - \lfloor a  \rfloor$.  So we've confirmed that Poisson summation has given us  the correct formula here, in spite of the lack of rigor at this point.  We also see rather quickly why the first periodic Bernoulli polynomial $P_1(x)$ 
appears so naturally in integer point enumeration in polytopes, from this perspective.

Why is the intuitive argument above not rigorous yet?  In order to plug a function $f$ into  
Poisson summation,  and consider convergence at each point of the domain, 
$f$ and its Fourier transform $\hat f$ must both satisfy some growth conditions at infinity, at the very least 
ensuring proper convergence of both sides of the Poisson summation formula. 
We will see such conditions later,  in Chapter~\ref{Fourier analysis basics},  Theorem \ref{nice2}.  After
 we learn how to use Poisson summation, we will return to this example, which will become rigorous in Section \ref{sec:1-dim polytopes}.
}
\hfill $\square$
\end{example}

We recall that a series $\sum_{n\in \Z} a_n$ is said to {\bf converge absolutely} if $\sum_{n\in \Z} |a_n|$ converges.   It's easy to see that the series in \eqref{FirstBernoulliPolynomial} for $P_1(x)$ does not converge absolutely.  Such convergent series that do not converge absolutely are called {\bf conditionally convergent}. 

To prove rigorously that the conditionally convergent series 
\eqref{FirstBernoulliPolynomial} does in fact converge, see 
Exercises \ref{Abel summation by parts}, \ref{Dirichlet's convergence test},  \ref{exponential sum bound}, and
  \ref{rigorous convergence of P_1(x)}, which include the Abel summation formula, and the Dirichlet convergence test.


\bigskip
\section{The cube, and its Fourier transform}

Perhaps the easiest way to extend the Fourier transform of the unit interval is to consider the 
 $d$-dimensional unit cube
\[
\square := \left[-\frac{1}{2}, \frac{1}{2} \right]^d.
\]
What is its Fourier transform?  When we compute a Fourier transform of a function $f$, we will
say that $\{ f, \hat f\}$ is a {\bf  Fourier pair}.  We have seen that
$\left\{  1_{[-\frac{1}{2}, \frac{1}{2}]}(x),   \sinc(\xi) \right\}$ is a Fourier pair in dimension $1$. 

\bigskip
\begin{example} \label{Example, unit cube}
\rm{
Due to the fact that the  cube is the direct product of line segments, it follows that the ensuing integral can be separated into a product of integrals, and so it is the product of $1$-dimensional transforms:
\begin{align}
\hat 1_{\square}(\xi) &= \int_{\R^d} 1_\square(x)   e^{-2\pi i \langle x, \xi \rangle} dx \\
 &= \int_{\square}    e^{-2\pi i(x_1 \xi_1 + \cdots + x_d \xi_d)}   dx \\
&= \prod_{k=1}^d    \int_{-\frac{1}{2}}^{\frac{1}{2}}   e^{-2\pi i x_k  \xi_k} dx_k \\
&= \prod_{k=1}^d   \frac{\sin(\pi \xi_k)}{\pi \xi_k},
\end{align}
valid for all $\xi \in \R^d$ such that none of their coordinates vanishes.  So here we have the Fourier pair
\[
\left\{ 1_\square(x),   \,  \prod_{k=1}^d   \frac{\sin(\pi \xi_k)}{\pi \xi_k}     \right\}.
\]
 In general, though, 
polytopes are not a direct product of lower-dimensional polytopes, so we will need to develop more tools to compute their Fourier transforms.
 }
 \hfill $\square$
\end{example}

\bigskip
\section{The simplex, and its Fourier transform}   

Another basic building block for polytopes is the  {\bf standard simplex}, 
\index{standard simplex}
defined by 
\begin{equation}
\RightTriangle := \left\{ 
x \in \R^d  \bigm |    \,   x_1 + \cdots + x_d  \leq 1, \text{ and all }  x_k \geq 0 
\right\} .
\end{equation}

 \begin{figure}[!h] 
		\centering
		\begin{tikzpicture}[scale=1]
			\draw (0,0) node[below left] {$0$};
			\draw[loosely dotted] (-1,-1) grid (2,2);
			\draw[->] (-1.25,0) -- (2.25,0) node[right] {$x$};
			\draw[->] (0,-1.25) -- (0,2.25) node[above] {$y$};
			\draw[thick] (0,0) -- (1,0) -- (0,1) -- cycle;
			\filldraw[nearly transparent, blue] (0,0) -- (1,0) -- (0,1) -- cycle;					
		\end{tikzpicture}  
		\caption{The standard simplex in $\R^2$} 
		\label{standard simplex in two dimensions}  \index{standard simplex}
	\end{figure}
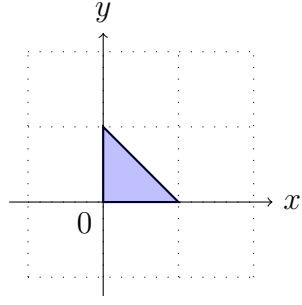

\begin{example}\label{standard simplex FT}
\index{standard simplex}	
\rm{
Just for fun, let's compute the Fourier transform of     $\triangle$
for $d=2$, via brute-force.
We may use the following parametrization (called a hyperplane description) 
for this standard triangle:
\[
\RightTriangle = \left\{   (x, y)   \bigm |        x+y  \leq 1, \text{ and }   x\geq 0, y\geq0   \right\}. 
\]
    Hence, we have: 
\begin{align*}
&\hat 1_{\, \rt}(\xi_1, \xi_2) := \int_{\, \rt}    e^{-2\pi i \big(x \xi_1 + y  \xi_2\big)} dx dy \\
&=   \int_0^1  \int_{y=0}^{y=1-x}   e^{-2\pi i \big(x \xi_1 + y \xi_2\big)} dy dx \\
&=   \int_0^1    e^{-2\pi i x \xi_1}  
\left[
   \frac{   e^{-2\pi i y \xi_2 }   }{-2\pi i  \xi_2 }   \Big|_{y=0}^{y=1-x}
   \right]
      dx  \\
&=  \frac{1}{-2\pi i  \xi_2 }       \int_0^1    e^{-2\pi i x \xi_1}   
                           \left(    e^{-2\pi i (1-x) \xi_2 }     - 1  \right)  dx  \\
&=  \frac{1}{-2\pi i  \xi_2 }       \int_0^1    
\left( 
e^{-2\pi i x (\xi_1 -\xi_2)}   e^{-2\pi i  \xi_2 }   - e^{-2\pi i x \xi_1}  
\right)  dx    \\
&=  \frac{1}{(-2\pi i)^2}
\frac{     e^{-2\pi i  \xi_2 }    }{      \xi_2(\xi_1-\xi_2)     }   
(e^{-2\pi i  (\xi_1 -\xi_2)}  -1) 
-
\frac{1}{(-2\pi i)^2}  \frac{     e^{-2\pi i  \xi_1 }   -1 }{        \xi_1 \xi_2     }      \\
&=  \frac{1}{(-2\pi i)^2}   \left[
\frac{     e^{-2\pi i  \xi_1} - e^{-2\pi i  \xi_2  }    }{      \xi_2(\xi_1-\xi_2)     }   
-
 \frac{     e^{-2\pi i  \xi_1 }   -1 }{        \xi_1 \xi_2     }   
\right].
\end{align*}
We may simplify further by noticing the rational function identity 
\[
\frac{     e^{-2\pi i  \xi_1 }   }{   \xi_2 (       \xi_1 - \xi_2 )    }  
-\frac{     e^{-2\pi i  \xi_1 }   }{        \xi_1 \xi_2     }  
=  \frac{     e^{-2\pi i  \xi_1 }   }{   \xi_1 (       \xi_1 - \xi_2 )    }, 
\]
giving us the symmetric function of $(\xi_1, \xi_2)$:
\begin{equation}\label{actual FT of the standard simplex}
\hat 1_{\rt}(\xi_1, \xi_2) = 
 \frac{1}{(-2\pi i)^2}   \left[
 \frac{     e^{-2\pi i  \xi_1 }   }{   \xi_1 (       \xi_1 - \xi_2 )    }
 + \frac{     e^{-2\pi i  \xi_2  }    }{      \xi_2(\xi_2-\xi_1)     }   
+
 \frac{     1 }{        \xi_1 \xi_2     }
 \right].
\end{equation}
}
\hfill $\square$
\end{example}

\medskip
\section{Convex sets and polytopes}

We need the concept of a {\bf convex set} $X\subset \R^d$, defined by the property that for any two points 
$x, y \in X$, the line segment joining them also lies in $X$.  In other words, the line segment
$\left\{   \lambda x + (1-\lambda)y \bigm |   0\leq \lambda \leq 1  \right\}   \subset X$, $\forall x, y \in X$.

Given any finite set of points $S:= \{ v_1, v_2, \dots, v_N\} \subset  \R^d$, we can also form the set of all
 {\bf convex linear combinations}  of $S$ by defining
\begin{equation}
\conv(S):= 
\left\{   
\lambda_1 v_1 + \lambda_2 v_2+ \dots +\lambda_{N} v_{N}
  \bigm |   
\sum_{k=1}^N \lambda_k = 1,  \text{ where all }  \lambda_k  \geq  0
\right\}.
\end{equation}
Given any set $U\subset \R^d$ (which is not restricted to be finite, or bounded), we define the  {\bf convex hull} of $U$
\index{convex hull}
as the set of convex linear combinations, taken over all finite subsets of $U$, and denoted by $\conv(U)$.

We define a {\bf polytope} as the convex hull of any finite set of points in $\R^d$.  This definition of a polytope is called its {\bf vertex description}.
\index{vertex description of a polytope}
We define a {\bf $k$-simplex}  \index{simplex} 
$\Delta$ as the convex hull of a finite set of vectors
$\{ v_1, v_2, \dots, v_{k+1} \}$:
\[
\Delta :=  \conv\{  v_1, v_2, \dots, v_{k+1} \},
\]  
where $0 \leq k \leq d$, and $v_2-v_1, v_3-v_1, \dots, v_{k+1} - v_1$ are linearly independent vectors in $\R^d$.
The points $v_1, v_2, \dots, v_{k+1}$ are called the vertices of $\Delta$, and this  
object is one of the basic building-blocks of polytopes, especially when triangulating a polytope. 


%

\begin{figure}[htb]
\begin{center}
\includegraphics[totalheight=3.2in]{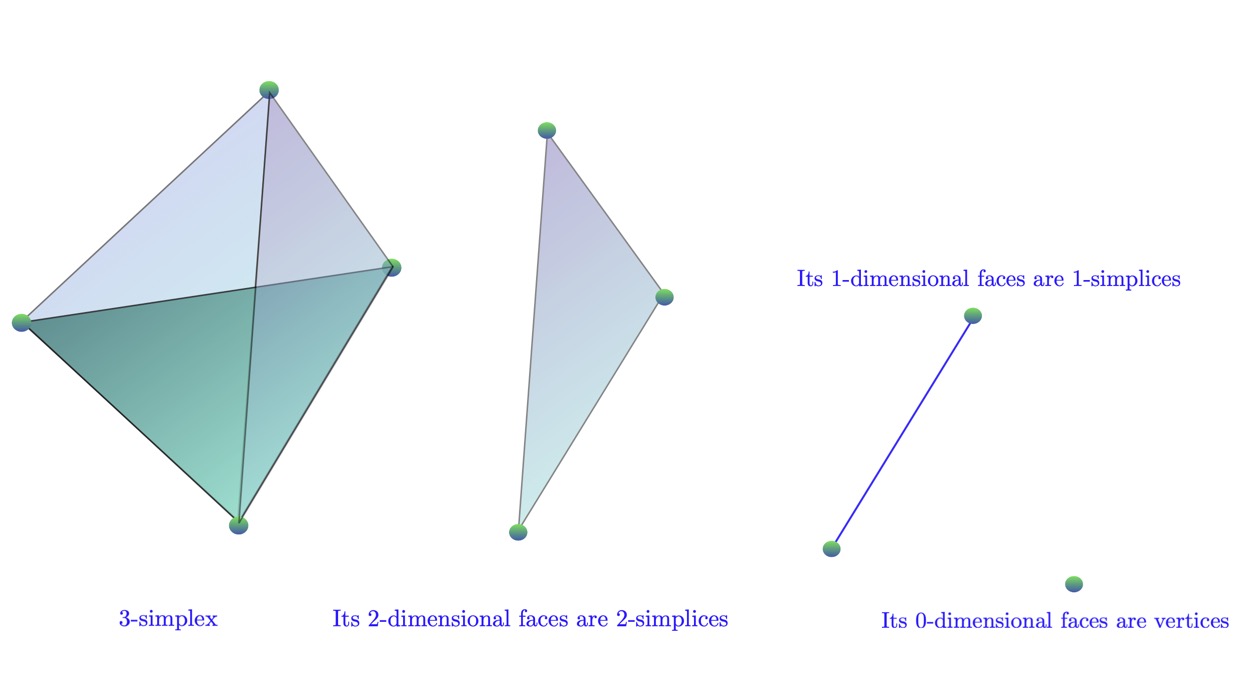}
\end{center}
\caption{A $3$-simplex and its faces, which are lower-dimensional simplices as well}    
\label{simplex}
\end{figure}

The simplex $\Delta$ is a $k$-dimensional polytope, sitting in $\R^d$.  When $k=d$, the dimension of $\Delta$ equals the dimension of the ambient space $\R^d$ - see Figure \ref{simplex}.
We have already  computed the Fourier transform of a particular $2$-simplex, 
in \eqref{actual FT of the standard simplex}.

How do we define a face of a polytope $\P$ more precisely? To begin, 
a  {\bf hyperplane} is defined by 
\[
H:=\{ x \in \R^d \mid \langle x, n \rangle = b \},
\]
for any fixed vector $n\in \R^d$, and any $b\in \R$. A hyperplane is called a
 {\bf supporting hyperplane for} $\P$ if $\P$ lies on one side of $H$, in the precise sense that:
\[
\P \subset \{ x \in \R^d \mid \langle x, n \rangle \leq b \}   \  \text{ or } 
\P \subset \{ x \in \R^d \mid \langle x, n \rangle \geq b \}.
\]
We now call $F\subseteq \P$ a {\bf face of } $\P$ if $F= H \cap \P$, for some supporting hyperplane $H$ of $\P$.  
As a consequence of the latter definition, the empty set is also a face of $\P$,  because we may pick a hyperplane very far from $\P$, which does not intersect $\P$.  
As a separate definition, we define  $\P$ to be a face of $\P$ itself.
 

With these preliminaries, we're now ready to compute the Fourier transform of any $2$-simplex in 
$\R^2$. 
In order to handle a general triangle, let $\Delta$ be any triangle in the plane, with vertices 
\[
v_1:= \icol{ a_1 \\ b_1},  v_2:=\icol{ a_2 \\ b_2} , v_3:= \icol{ a_3 \\ b_3}. 
\]
Can we reduce the computation of $\hat 1_{\Delta}$ to our already known formula 
for  $\hat 1_{\rt}$, given by \eqref{actual FT of the standard simplex}?
We first notice (after a brief cup of coffee) that we can map any triangle in the plane to the standard triangle,  by using a linear transformation followed by a translation: 
\begin{equation}\label{M followed by T}
\Delta = M ( \,  \RightTriangle  ) + v_3,
\end{equation}
where $M$ is the $2\times2$ matrix whose columns are $v_1-v_3$ and $v_2-v_3$.  
We are now ready to compute the Fourier transform of a general triangle $\Delta$:  
\[
\hat 1_{\Delta}(\xi) = \int_{\Delta} e^{-2\pi i \langle \xi, x \rangle} dx  
                          =    \int_{M(\rt)+ v_3} e^{-2\pi i \langle \xi, x \rangle} dx.
\]
Making the substitution $x := My+v_3$, with $y \in \rt$, we have $dx = |\det M| dy$, and so
\begin{align*}
& \int_{M(\rt)+v_3} e^{-2\pi i \langle \xi, x \rangle} dx  
=   |\det M| \int_{\, \rt} e^{-2\pi i \langle \xi, M y +v_3  \rangle} dy   \\
&=   |\det M|   e^{-2\pi i \langle \xi, v_3  \rangle}
               \int_{\,   \rt} e^{-2\pi i \langle   M^{T}  \xi,  y   \rangle} dy  \\
&=  |\det M|  e^{-2\pi i \langle \xi, v_3  \rangle}  \hat 1_{\, \rt}(M^{T} \xi) \\
&=  |\det M|  e^{-2\pi i \langle \xi, v_3  \rangle}  
    \hat 1_{\, \rt}   \big(  \langle v_1-v_3,   \xi  \rangle,    \langle v_2-v_3,   \xi  \rangle    \big) \\
 &=
 |\det M|  e^{-2\pi i \langle \xi, v_3  \rangle}  
 \frac{1}{(-2\pi i)^2}   \left[
 \frac{     e^{-2\pi i  z_1 }   }{   z_1 (       z_1 - z_2 )    }
 + \frac{     e^{-2\pi i  z_2  }    }{      z_2(z_2-z_1)     }   
+
 \frac{     1 }{        z_1 z_2     }
 \right],
 \end{align*}
 where we've used our formula \eqref{actual FT of the standard simplex} for the FT of the standard triangle (thereby bootstrapping out way to the general case) with 
 $z_1:=\langle v_1-v_3,   \xi  \rangle$, and $z_2:= \langle v_2-v_3,   \xi  \rangle$.  Substituting these values into the latter expression, we finally arrive at the FT of our general triangle $\Delta$:
 \begin{align}\label{FT of a general triangle}
\hat 1_{\Delta}(\xi) 
= \tfrac{  |\det M|  }{(-2\pi i)^2}
  \left[
 \frac{     e^{-2\pi i  \langle v_1,   \xi  \rangle }   }{   \langle v_1-v_3,   \xi  \rangle 
                                                                                     \langle v_1-v_2,   \xi  \rangle    }
 + \frac{  e^{-2\pi i  \langle v_2,   \xi  \rangle }   }{   \langle v_2-v_3,   \xi  \rangle 
                                                                                     \langle v_2-v_1,   \xi  \rangle    }  
+
 \frac{    e^{-2\pi i \langle \xi, v_3  \rangle}   }{   \langle v_3-v_1,   \xi  \rangle 
                                                                                     \langle v_3-v_2,   \xi  \rangle    }  
   \right].
\end{align}

We can notice in equation \eqref{FT of a general triangle} many of the same patterns that had already occurred in Example 
\ref{cross-polytope example in R^2}.  Namely,
the Fourier transform of a triangle has denominators that are products of linear forms in $\xi$, 
and it is a finite linear combination of rational functions multiplied by complex exponentials.

Also, in the particular case of 
equation  \eqref{FT of a general triangle}, $ \hat 1_\Delta(\xi) $ is a symmetric function of $v_1, v_2, v_3$, as we might have expected.

Using exactly the same ideas that were used in equation \eqref{FT of a general triangle}, it is possible to prove (by induction on the dimension) that  the Fourier transform of a general $d$-dimensional simplex $\Delta \subset \R^d$  is:
\begin{equation}\label{FT of a d-dimensional simplex}
\hat 1_{\Delta}(\xi) = (\vol \Delta) d!  \sum_{j=1}^N   
\frac{e^{-2\pi i \langle v_j, \xi \rangle}}{\prod_{k=1}^d \langle v_j-v_k, \xi \rangle   }[k \not= j],
\end{equation}
where the vertex set of $\P$ is $\{ v_1, \dots, v_N\}$  (Exercise \ref{FT of a general simplex, brute-force}), and in fact the same formula persists for all complex $\xi \in \C^d$ such that the products of linear forms in the denominators do not vanish.

However, looking back at the computation leading to \eqref{FT of a general triangle}, 
and the corresponding computation which would give \eqref{FT of a d-dimensional simplex},
the curious reader might be thinking: 

\medskip
\centerline{  ``There must be an easier way!'' }

But never fear - indeed there is.
So even though at this point the computation of $\hat 1_\Delta(\xi)$ may be a bit laborious (but still interesting),  
computing the Fourier transform of a general simplex will become quite easy once we will revisit it in a later chapter (see Theorem \ref{brion, continuous form}).

\bigskip
\section{Stretching and translating}
\bigskip
The perspicacious reader may have noticed that in order to arrive at the formula 
\eqref{FT of a general triangle} above for the FT of a general triangle, we exploited the fact that the Fourier transform interacted peacefully with the linear transformation $M$, and with the translation by the vector $v$.   Is this true in general? 

Indeed it is, and we record these thoughts in the following two lemmas, which will become our bread and butter for future computations. 
In general, given any invertible linear transformation $M :\R^d \rightarrow \R^d$, and any function $f:\R^d \rightarrow \C$ whose FT (Fourier transform) exists, we have the following useful interaction between Fourier transforms and linear transformations.

\begin{lem}[Stretch]    \index{stretch lemma}
\label{FT under linear maps}
\begin{equation}\label{The FT under streching}
(\widehat{f \circ M})(\xi)= \frac{1}{|\det M|}  \hat f\left(M^{-T}\xi \right) 
\end{equation}
\end{lem}
\begin{proof}
By definition, we have
$
(\widehat{f \circ M})(\xi) :=\int_{\R^d} f(Mx) e^{-2\pi i \langle \xi, x \rangle} dx.
$
We perform the change of variable $y:= Mx$, implying that $dy = |\det M| dx$, so that:
\begin{align*}
(\widehat{f \circ M})(\xi) &=   
\frac{1}{|\det M|}  \int_{\R^d} f(y) e^{-2\pi i \langle \xi, M^{-1}y \rangle} dy \\
&=\frac{1}{|\det M|}  \int_{\R^d} f(y) e^{-2\pi i \langle M^{-T}\xi, y \rangle} dy \\
&= \frac{1}{|\det M|}  \hat f\left(M^{-T} \xi \right).
\end{align*}
\end{proof}

What about translations?  They are even simpler.
\begin{lem}[Translate]     \index{translate lemma}
\label{FT under translations}
For any translation $T(x):= x + v$, where $v\in \R^d$ is a fixed vector, we have
\begin{equation}\label{The FT under translations}
(\widehat{f \circ T})(\xi)= e^{2\pi i \langle \xi, v \rangle}    \hat f(\xi).        
\end{equation}
\end{lem}
\begin{proof} 
Again, by definition we have
$
(\widehat{f \circ T})(\xi) := \int_{\R^d}    f( Tx) e^{-2\pi i \langle \xi, x \rangle} dx,
$
so that performing the simple change of variable $y = Tx := x + v$, we have $dy = dx$. 
The latter integral becomes
\begin{align*}
(\widehat{f \circ T})(\xi)
&=   \int_{\R^d} f(y) e^{-2\pi i \langle \xi, y-v \rangle} dy \\
&=  e^{2\pi i \langle \xi, v \rangle}    
                \int_{\R^d} f(y) e^{-2\pi i \langle \xi, y \rangle} dy :=  
           e^{2\pi i \langle \xi, v \rangle}    \hat f(\xi). 
\end{align*}
\end{proof}

 



In general, any function $\phi:\R^d \rightarrow \C$ of the form
\begin{equation}
\phi(x) = Mx+v,
\end{equation}
where $M$ is a fixed linear transformation and $v\in \R^d$ is a fixed vector, is
called an {\bf affine transformation}. \index{affine transformation}  For example, we've already seen
in \eqref{M followed by T} that the right triangle  $\RightTriangle$ was mapped to the more general triangle
$\Delta$ by an affine transformation.
So the latter two lemmas allow us to compose Fourier transforms very easily with affine transformations.

\begin{example}
\rm{
The simplest example of the Stretch Lemma \ref{FT under linear maps} is obtained
in $\R$, where the matrix $M = r$, a positive real number.  So we have $M^{-T} = \frac{1}{r}$.
Considering $f(rx)$ as a function of $x \in \R$, we have by \eqref{The FT under translations}:
\begin{equation} \label{simple example of stretching}
\widehat{f(rx)} := (\widehat{f \circ M})(\xi) = \tfrac{1}{r}  \hat f\left( \tfrac{1}{r}  \xi \right).
\end{equation}
As an interesting sub-example, let's take $f(x) := 1_{\left[-\tfrac{c}{2}, \tfrac{c}{2} \right]}(x) $, 
for a fixed constant $c>0$.  What's the easy way to use the Stretch lemma to compute $\hat f(\xi)$? First, we have to make a slight conversion:  $1_{\left[-\tfrac{c}{2}, \tfrac{c}{2} \right]}(x) 
= 1_{\left[-\tfrac{1}{2}, \tfrac{1}{2} \right]}(\tfrac{1}{c} x)$.  Using the FT of the unit interval, equation \eqref{sinc function formula}, 
together with \eqref{simple example of stretching}, we have:
\begin{equation}\label{Stretch lemma for the sinc function}
\hat  1_{\left[-\tfrac{c}{2}, \tfrac{c}{2} \right]}(\xi) = c\,  \hat 1_{\left[-\tfrac{1}{2}, \tfrac{1}{2} \right]}(c \xi)
=c\, \sinc(c\xi) = \frac{\sin(c\pi \xi)}{\pi \xi}.
\end{equation}
}
\hfill $\square$
\end{example}

\bigskip
\begin{example}
\rm{
Consider any set $B\subset \R^d$, for which $1_B$ is integrable, and let's translate $B$ by a fixed vector $v \in \R^d$, and compute
$\hat 1_{B+v}(\xi)$. 

We note that because $1_{B+v}(\xi) = 1_{B}(\xi-v)$, the translate lemma applies, but with a minus sign.  That is, we can use
$T(x):= x-v$ and  $f:= 1_B$ to get:
\begin{equation} \label{The FT of a translate of B}
\hat 1_{B+v}(\xi) = \widehat{   (1_B \circ T)        }(\xi) = e^{-2\pi i  \langle  \xi, v \rangle} \hat 1_B(\xi).
\end{equation}
}
\hfill $\square$
\end{example}

\bigskip
\section{The parallelepiped, and its Fourier transform}   \index{parallelepiped}

Now that we know how to compose the FT with affine transformations (translations and linear transformations), we can easily find the FT of
any parallelepiped in $\R^d$ by using our formula for the Fourier transform of the unit cube 
$\square := \left[-\frac{1}{2}, \frac{1}{2} \right]^d$, which we derived in Example \ref{Example, unit cube}:
\begin{align}\label{first FT of a cube}
\hat 1_{\square}(\xi) = \prod_{k=1}^d   \frac{\sin(\pi \xi_k)}{\pi \xi_k},
\end{align}
for all $\xi \in \R^d$ such that all the coordinates of $\xi$ do not vanish.  
First, we translate the cube $\square$
 by the vector $(\frac{1}{2}, \cdots, \frac{1}{2})$, to obtain 
 \[
 C:= \square + \left(\frac{1}{2}, \,  \cdots, \frac{1}{2}  \right) = [0, 1]^d.
 \]
It's straightforward to compute its FT as well (Exercise \ref{transform.of.unit.cube}), by using  Lemma \ref{FT under translations}, the `translate' lemma:
\begin{equation}\label{second appearance of FT of the cube}
\hat 1_{C}(\xi) =  \frac{1}{(2\pi i)^d}  \prod_{k=1}^d \frac{    1- e^{-2\pi i \xi_k}   }{   \xi_k         }.
\end{equation}

Next, we define a $d$-dimensional {\bf parallelepiped} $\P \subset \R^d$ as an affine image of the unit cube.  In other words, 
any parallelepiped has the description 
\[
\P = M(C) + v, 
\]
for some linear transformation $M$, and some translation vector $v$.    Geometrically, the cube is stretched and translated 
into a parallelepiped. 

\begin{figure}[htb]
\begin{center}
\includegraphics[totalheight=3.2in]{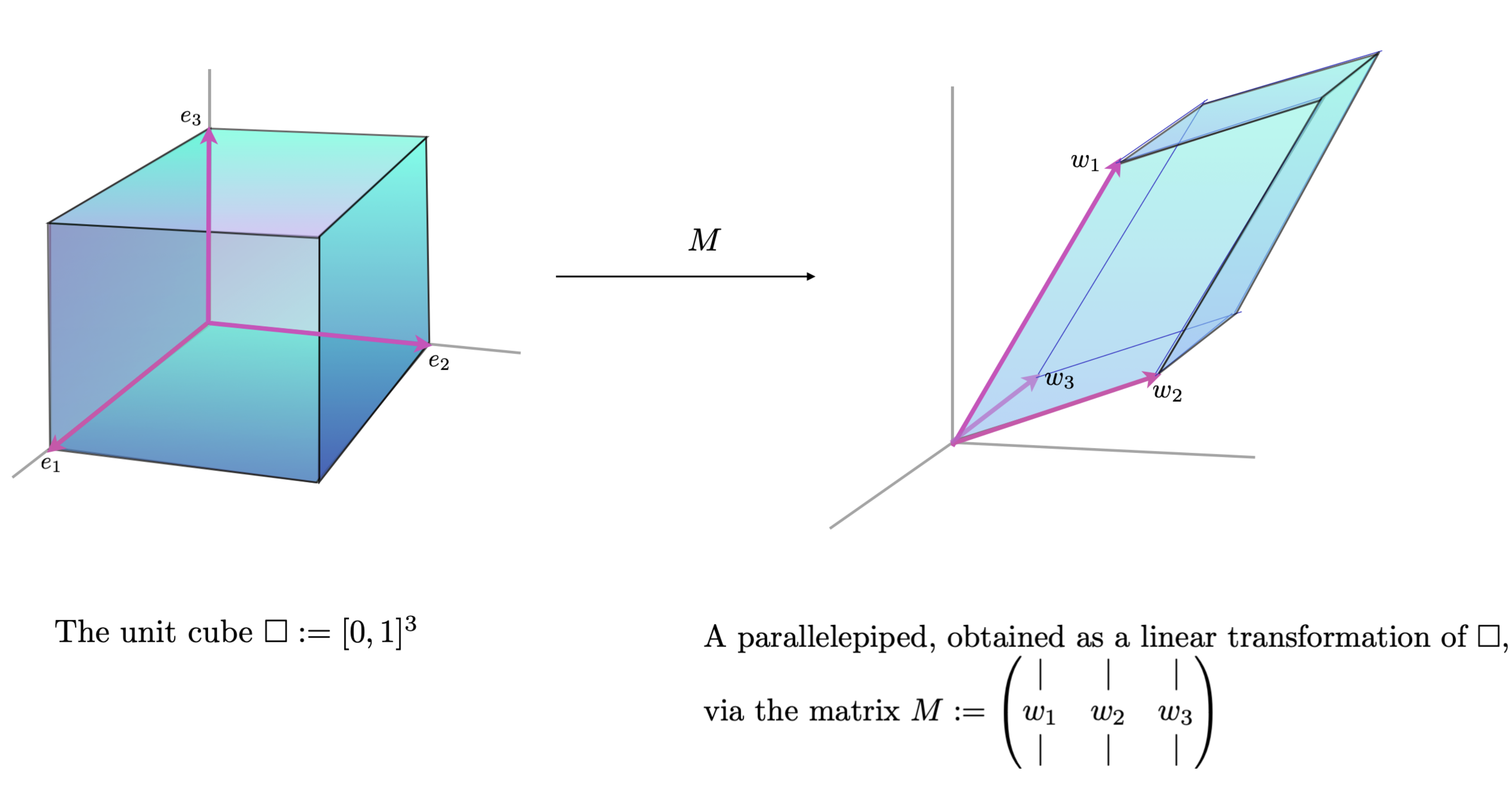}
\end{center}
\caption{Mapping the unit cube to a parallelepiped}    
\label{3Dparallelepiped}
\end{figure}

For the sake of concreteness,  will will first set $v:= 0$ and compute the Fourier transform of  $\P:= M(C)$, where we now give $M$ as a  $d \times d$ invertible matrix  whose columns are
$w_1, w_2, \dots, w_d$.   Because the cube $C$ may be written as a convex linear combination of the basis vectors $e_j$, 
we see that $\P$ may be written as a convex linear combination of  $M e_j = w_j$.
In other words, we see that the (closed) parallelepiped $\P$ has the equivalent vertex description:
\[
\P = \left\{  \sum_{k=1}^d   \lambda_k w_k \bigm |    \text{ all }  \lambda_k \in [0, 1] \right\}.
\]
To review the basics, let's compute the FT of our parallelepiped $\P$ from first principles:
\begin{align} \label{cube composed with M}
\hat 1_\P(\xi) &:= \int_{\P} e^{-2\pi i \langle \xi, x \rangle} dx = \int_{M(C)} e^{-2\pi i \langle \xi, x \rangle} dx\\
&= |\det M|   \int_{C} e^{-2\pi i \langle \xi, My \rangle} dy\\
&= |\det M|   \int_{C} e^{-2\pi i \langle M^T \xi, y \rangle} dy := |\det M|  \, \hat 1_C\left(M^T \xi  \right) \\  \label{for Q}
&=    \frac{   |\det M|    }{(2\pi i)^d}  \prod_{k=1}^d \frac{    1- e^{-2\pi i \langle w_k, \xi \rangle}   }{   \langle w_k, \xi \rangle     }.
\end{align}
where in the third equality we used the substitution $x:= My$, with $y\in C$, yielding $dx = |\det M| dy$.  In the last equality, we used our known formula \eqref{second appearance of FT of the cube} for the FT of the cube $C$, together with the elementary linear algebra fact that the $k$'th coordinate of $M^T \xi$ is given by
$\langle w_k, \xi \rangle$. 

Finally, for a general parallelepiped, we have $Q:= \P + v$, so that by definition
\[
Q = \left\{ v+  \sum_{k=1}^d   \lambda_k w_k \bigm |      \text{ all }  \lambda_k \in [0, 1]  \right\}.
\]

Noting that $1_{\P+ v}(\xi) = 1_{\P}(\xi -v)$,  
we compute the Fourier transform of $Q$ by using the `translate lemma' (Lemma \ref{FT under translations}), 
together with formula \eqref{for Q} for the Fourier transform of $\P$:
\begin{equation}\label{FT of a general parallelepiped} 
\hat  1_Q(\xi)    =   e^{-2\pi i \langle \xi, v \rangle} 
     \frac{   |\det M|    }{(2\pi i)^d}  
         \prod_{k=1}^d \frac{    1- e^{-2\pi i \langle w_k, \xi \rangle}   }{   \langle w_k, \xi \rangle     },
\end{equation}
for all $\xi \in \R^d$, except for those $\xi$ that are orthogonal to one of the $w_k$ (which are edge vectors for $Q$).
\begin{example}
\rm{
A straightforward computation shows that if we let $v:= -\frac{w_1 + \cdots + w_d}{2}$, then
$Q:= \{ v+  \sum_{k=1}^d   \lambda_k w_k \mid  \text{ all }  \lambda_k \in [0, 1] \}$ is symmetric about the origin, in the sense that
$x \in Q \iff -x \in Q$ (Exercise \ref{symmetrized parallelepiped}).  In other words, the center of mass of this new $Q$ is now the origin. 
Geometrically, we've translated the previous parallelepiped by using half its `body diagonal'. 
 For such a parallelepiped $Q$, centered at the origin, formula \eqref{FT of a general parallelepiped}  above gives the more pleasing expression: 
\begin{align}
\hat  1_Q(\xi) &= e^{2\pi i \langle \xi,    \frac{w_1 + \cdots + w_d}{2}     \rangle} 
     \frac{   |\det M|    }{(2\pi i)^d}  
         \prod_{k=1}^d \frac{    1- e^{-2\pi i \langle w_k, \xi \rangle}   }{   \langle w_k, \xi \rangle     } \\
 &=       
     \frac{   |\det M|    }{(2\pi i)^d}  
         \prod_{k=1}^d \frac{    e^{\pi i \langle w_k, \xi \rangle}- e^{-\pi i \langle w_k, \xi \rangle}   }{   \langle w_k, \xi \rangle     } \\
         &= \frac{   |\det M|    }{(2\pi i)^d}  
         \prod_{k=1}^d \frac{   (2i) \sin( \pi  \langle w_k, \xi \rangle)   }{   \langle w_k, \xi \rangle     } \\
       &=   |\det M|  
         \prod_{k=1}^d \frac{  \sin( \pi  \langle w_k, \xi \rangle)   }{ \pi  \langle w_k, \xi \rangle }.
\end{align}
To summarize, for a parallelepiped that is symmetric about the origin, we have the Fourier pair
\[
\left\{ 1_Q(x),   \ \   | \det M |  
         \prod_{k=1}^d \frac{  \sin( \pi  \langle w_k, \xi \rangle)   }{ \pi  \langle w_k, \xi \rangle }     \right\}.
\]
We could have also computed the latter FT by beginning with our known 
Fourier transform \eqref{first FT of a cube} of the cube $\square$, composing the FT 
with the same linear transformation $M$ of \eqref{cube composed with M}, and using the `stretch' lemma, so everything is consistent. 
}
\hfill $\square$
\end{example}



\bigskip
\section{The cross-polytope}   
\index{cross-polytope}

\begin{figure}[htb]
\begin{center}
\includegraphics[totalheight=2.8in]{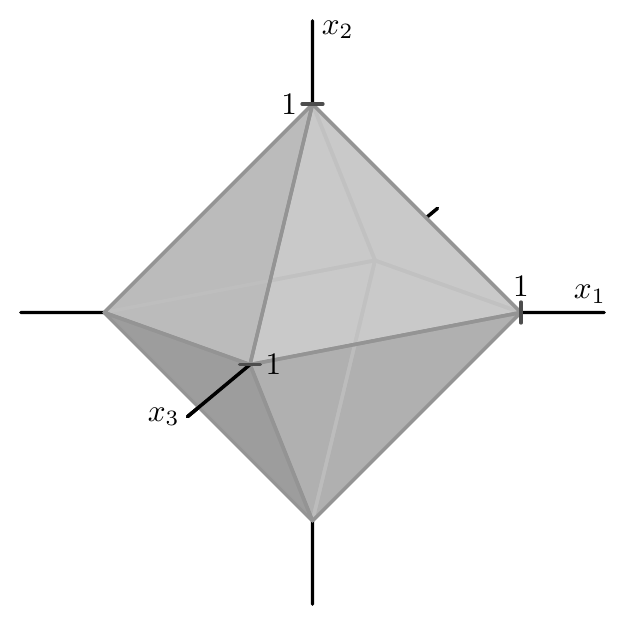}
\end{center}
\caption{The cross-polytope $\Diamond$ in $\R^3$ (courtesy of David Austin)}    \label{crosspic}
\end{figure}

Another natural convex body in $\R^2$ is the cross-polytope
\begin{equation}  \label{2dim.crosspolytope}
\Diamond_2 := \left\{ \left( x_1, x_2 \right) \in \R^2    \bigm |    \,    \left| x_1 \right| + \left| x_2 \right|   \leq 1 \right\} .
\end{equation}

In dimension $d$, the {\bf cross-polytope} $\Diamond_d$      \label{cross polytope}
 can be defined similarly by its {\bf  hyperplane description}
\begin{equation}  \label{crosspolytopehyperplanes}
\Diamond_d := \left\{ \left( x_1, x_2, \dots, x_d \right) \in \R^d  \bigm |     \, \left| x_1 \right| + \left| x_2 \right| + \dots + \left| x_d \right| \leq 1 \right\} .
\end{equation}
The cross-polytope is also, by definition, the unit ball in the $L_1$-norm on Euclidean space, and from this perspective a very natural object.  In $\mathbb R^3$, 
the cross-polytope $\Diamond_3$ is often called an {\bf {octahedron}}.  \index{octahedron}

In this section we only work out the $2$-dimensional case of the Fourier transfrom of the crosspolytope, 
In Chapter \ref{chapter.Brion}, we will work out the 
Fourier transform of any $d$-dimensional cross-polytope, $\hat 1_{\Diamond_d}$, 
because we will have more tools at our disposal.

Nevertheless, it's instructive to compute  $\hat 1_{\Diamond_2}$ via brute-force for $d=2$ here, 
 in order to gain some facility with the computation of Fourier transforms.

\begin{example}   \label{cross-polytope example in R^2}
\rm{
Using the definition of the Fourier transform, we first compute  the FT of the  $2$-dimensional 
cross polytope:
\begin{align}
\hat 1_{\Diamond_2}(\xi)    &:= \int_{\Diamond_2} e^{-2\pi i \langle   \xi, x     \rangle} dx.
\end{align}
In $\R^2$, we may write $\Diamond_2$ as a union of the following $4$ triangles:  
\begin{align*}
\Delta_1&:= \conv (  \icol{0\\0},   \icol{1\\0},  \icol{0\\1}                      )\\ 
\Delta_2&:= \conv (          \icol{0\\0},   \icol{-1\\0},  \icol{0\\1}        )\\
\Delta_3&:= \conv (        \icol{0\\0},   \icol{-1\\0},  \icol{0\\-1}           )\\
\Delta_4&:= \conv (        \icol{0\\0},   \icol{1\\0},  \icol{0\\-1}           ).
\end{align*}
Since these four triangles only intersect in lower-dimensional subsets of $\R^2$, the $2$-dimensional integral vanishes on such lower dimensional subsets, and we have:
\begin{equation}  \label{transform of 2d crosspolytope}
\hat 1_{\Diamond_2}(\xi)  = \hat 1_{\Delta_1}(\xi) + \hat 1_{\Delta_2}(\xi) 
+ \hat 1_{\Delta_3}(\xi)  + \hat 1_{\Delta_4}(\xi).
\end{equation}
Recalling from equation \eqref{actual FT of the standard simplex} 
of example \ref{standard simplex FT} 
that the Fourier transform of the standard simplex
\index{standard simplex}
 $\Delta_1$ is
\begin{equation}\label{simplex transform}
 \hat 1_{\Delta_1}(\xi)  =  \left( \frac{1}{2\pi i} \right)^2 
     \left(
       \frac{1}{\xi_1 \xi_2}
    +  \frac{\ e^{-2\pi i  \xi_1} }{(-\xi_1 +  \xi_2) \xi_1}
    + \frac{  \ e^{-2\pi i  \xi_2} }{( \xi_1 -  \xi_2) \xi_2}
     \right),
\end{equation}
we can compute  $\hat 1_{\Delta_2}(\xi)$,  by reflecting $\Delta_2$ about the $x_2-axis$ (the Jacobian of this transformation is $1$), and
using the already-computed transform \eqref{simplex transform} of $\Delta_1$:
\begin{align*}
\hat 1_{\Delta_2}(\xi_1, \xi_2) &:= \int_{\Delta_2}    e^{-2\pi i (x_1 \xi_1 + x_2 \xi_2)} dx \\
&=  \int_{\Delta_1}                                     e^{-2\pi i (-x_1 \xi_1 + x_2 \xi_2)} dx \\
&=  \int_{\Delta_1}                                     e^{-2\pi i (x_1 (-\xi_1) + x_2 \xi_2)} dx \\
&= \hat 1_{\Delta_1}(-\xi_1, \xi_2)).
\end{align*}
Similarly, we have $\hat 1_{\Delta_3}(\xi_1, \xi_2) = \hat 1_{\Delta_1}(-\xi_1, -\xi_2)$, and
$\hat 1_{\Delta_4}(\xi_1, \xi_2) = \hat 1_{\Delta_1}(\xi_1, -\xi_2)$.

Hence we may continue the computation from equation  \ref{transform of 2d crosspolytope} above, putting all the pieces back together:
\begin{align} \label{Fourier transform of 2d crosspolytope}
\hat 1_{\Diamond_2}(\xi)  &= \hat 1_{\Delta_1}(\xi_1, \xi_2) + \hat 1_{\Delta_1}(-\xi_1, \xi_2) 
+ \hat 1_{\Delta_1}(-\xi_1, -\xi_2) + \hat 1_{\Delta_1}(\xi_1, -\xi_2) \\
&=  \left( \frac{1}{2\pi i} \right)^2 
     \left(
       \frac{1}{\xi_1 \xi_2}
    +  \frac{-\ e^{2\pi i  \xi_1} }{(-\xi_1 +  \xi_2) \xi_1}
    + \frac{ - \ e^{2\pi i  \xi_2} }{( \xi_1 -  \xi_2) \xi_2}
     \right)  \\
&+ \left( \frac{1}{2\pi i} \right)^2 
     \left(
       \frac{-1}{\xi_1 \xi_2}
    +  \frac{\ e^{-2\pi i  \xi_1} }{(\xi_1 +  \xi_2) \xi_1}
    + \frac{  \ e^{2\pi i  \xi_2} }{( \xi_1 +  \xi_2) \xi_2}
     \right) \\
&+ \left( \frac{1}{2\pi i} \right)^2 
     \left(
       \frac{1}{\xi_1 \xi_2}
    +  \frac{e^{-2\pi i  \xi_1} }{(\xi_1 -  \xi_2) \xi_1}
    + \frac{e^{-2\pi i  \xi_2} }{( -\xi_1 + \xi_2) \xi_2}
     \right) \\
&+ \left( \frac{1}{2\pi i} \right)^2 
     \left(
       \frac{-1}{\xi_1 \xi_2}
    +  \frac{e^{2\pi i  \xi_1} }{(\xi_1 +  \xi_2) \xi_1}
    + \frac{e^{-2\pi i  \xi_2} }{( \xi_1 +  \xi_2) \xi_2}
     \right) \\
 &=  -\frac{1}{2\pi^2}  \left(
     \frac{\cos(2\pi \xi_1) }{(\xi_1 -  \xi_2) \xi_1}
 +  \frac{\cos(2\pi \xi_2) }{(-\xi_1 +  \xi_2) \xi_2}
 +  \frac{\cos(2\pi \xi_1) }{(\xi_1 +  \xi_2) \xi_1}
  +  \frac{\cos(2\pi \xi_2) }{(\xi_1 +  \xi_2) \xi_2}
  \right)  \\   \label{formula 1 for the FT of the 2-d crosspolytope}
  &=  -\frac{1}{\pi^2}  \left(
     \frac{     \cos(2\pi \xi_1)  -     \cos(2\pi \xi_2)     }{  (\xi_1 +  \xi_2)(   \xi_1 -  \xi_2)   }
  \right).
\end{align}
}
\hfill $\square$
\end{example}
It's time to mention another important relationship between the cross-polytope $\Diamond$ 
\index{cross-polytope}
and the cube $\P~:= ~[-1, 1]^d$.
To see this relationship, we define, for any polytope $\P \subset \R^d$, its {\bf polar polytope}:
\begin{equation}\label{polar polytope, definition}
\P^o := \left\{ x\in \R^d   \bigm |   \,  \langle x, y\rangle \leq 1, \text{ for all } y \in  \P    \right\}.
\end{equation}

\begin{figure}[htb]
\begin{center}
\includegraphics[totalheight=3.0in]{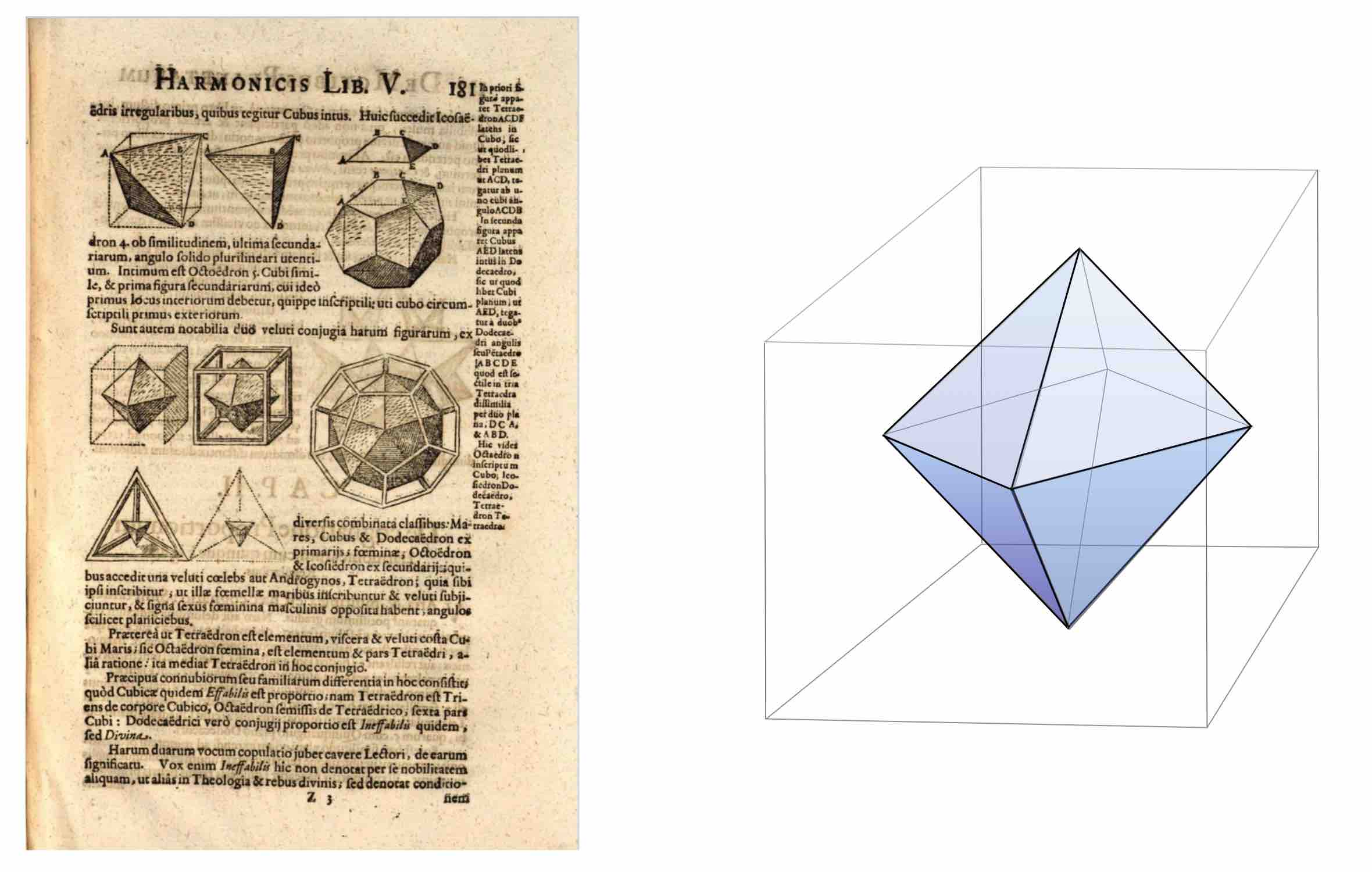}
\end{center}
\caption{Left:  a page from 
Kepler's book, \emph{Harmonices Mundi} ($1619$),  showing the author's interest in various polar polytopes, over $400$ years ago.  Right: The cube and the cross-polytope as polars of each other.  }    \label{duals}
\end{figure}

It is an easy fact (Exercise  \ref{polars of each other}) that in $\R^d$,  the cross-polytope 
$\Diamond_d$ and the cube $\P:= [-1, 1]^d$ are polar to each other, as in Figure \ref{duals}.

In many applications, it turns out the the volume of the cross-polytope plays an important role.  So we compute it here, for a generalized cross-polytope defined by the 
image of $\Diamond_d$ under any invertible linear transformation $M$, whose $k$'th olumn is defined by $v_k$:
\begin{equation}
Q:= M(\Diamond_d) = \conv\left( \pm v_1, \pm v_2, \dots, \pm v_d      \right).
\end{equation}

\begin{lem}
With the notation above, we have
\begin{enumerate}[(a)]
\item \label{first part of volume of crosspolytope}
\[
\vol \Diamond_d = \frac{2^d}{d!}
\]
\item   \label{second part of volume of crosspolytope}
\[
\vol Q  = |\det M| \frac{2^d}{d!}.
\]
\item  \label{third part of volume of crosspolytope}
In particular, if $v_k:= \alpha_k e_k$, then we have
\[
\vol Q = \alpha_1 \alpha_2 \cdots \alpha_d \frac{2^d}{d!}.
\]
\end{enumerate}
\end{lem}
\begin{proof}
To prove part \ref{first part of volume of crosspolytope}, we can simply triangulate the crosspolytope into $2^d$ isometric simplices by using the coordinate planes.  Each such simplex has volume $\frac{1}{d!}$, so we're done. 
For part \ref{second part of volume of crosspolytope}, we have:
\begin{equation}
\vol Q:= \int_{M(\Diamond_d)} dx =
|\det M|  \int_{\Diamond_d} dy = |\det M| \vol \Diamond_d = |\det M| \frac{2^d}{d!},
\end{equation}
where we used the change of variable $x:= M y$ and its ensuing Jacobian transformation
$dx = |\det M| dy$. The last equality above follows from 
part \ref{first part of volume of crosspolytope}.
 Part \ref{third part of volume of crosspolytope} follows trivially from part 
  \ref{second part of volume of crosspolytope}, using the determinant of a diagonal matrix.
\end{proof}


\bigskip
\section{Observations and questions}

Now we can make several observations about all of the formulas that we found so far,
 for the Fourier transforms of various polytopes.   For the $2$-dimensional cross-polytope, we found that
\begin{equation} \label{again the FT of a 2d crosspolytope}
\hat 1_{\Diamond_2}(\xi)  = -\frac{1}{\pi^2}  \left(
     \frac{     \cos(2\pi \xi_1)  -     \cos(2\pi \xi_2)     }{  (\xi_1 +  \xi_2)(   \xi_1 -  \xi_2)   }
  \right).
  \end{equation}
\begin{enumerate}[(a)]
\item \  It is real-valued for all $\xi \in \R^2$, and this is due to the fact that $\Diamond_2$ is symmetric about the origin (see section \ref{Centrally symmetric polytopes}).   
\begin{question}
Is it true that {\emph any} symmetric property of a polytope $\P$ is somehow mirrored by a corresponding symmetric property  of its Fourier transform?
\end{question}
Although this question is not well-defined at the moment (it depends on how we define `symmetric property'), 
it does sound exciting, and we can morph it into a few well-defined questions later.

\item  \ The only apparent singularities of the FT in \eqref{again the FT of a 2d crosspolytope}
(though they are in fact removable singularities) are the two lines 
$\xi_1 - \xi_2=0$ and $\xi_1 + \xi_2=0$, and these two
lines are {\it perpendicular} to the facets of $\Diamond_2$,
which is not a coincidence (see Chapter \ref{Stokes' formula and transforms}).

\item  \ It is always true that the Fourier transform of a polytope is an entire function, by  
Lemma \ref{FT of a polytope is entire}, so that the singularities
in the denominator $ (\xi_1 +  \xi_2)(   \xi_1 -  \xi_2)$  of 
\eqref{again the FT of a 2d crosspolytope} must be removable singularities!

\item  The denominators of all of the FT's so far are always products of {\bf linear forms} in $\xi$.   
\begin{question} \rm{[Rhetorical]}
Is  it true that the Fourier transform of any polytope is always a finite sum of rational functions times an exponential, where the denominators of the rational functions are always products of linear forms?
\end{question}

{\bf Answer}: (spoiler alert) Yes!  It's too early to prove this here, but we'll do it in the course of proving Theorem \ref{brion2}.

\item   \ We may retrieve the volume of $\Diamond_2$ by letting $\xi_1$ and $\xi_2$ tend to zero 
(Exercise \ref{retreiving volume of 2d.crosspolytope}), as always.   Doing so, we obtain 
\[
\lim_{\xi \rightarrow 0}  \hat 1_{\Diamond_2}(\xi) = 2 = \text{Area}(\Diamond_2).
\]
\end{enumerate}


\section*{Notes}
\begin{enumerate}[(a)]
\item  Another way to compute $1_{\Diamond}(\xi)$ for the $2$-dimensional cross-polytope $\Diamond$ is by starting with the square $[-\frac{1}{2}, \frac{1}{2}]^2$ and applying a rotation of the plane by $\pi/4$, followed by a simple dilation.  Because we know that  linear transformations interact in a very elegant way with the FT,  this method gives an alternate approach
 for the Example  \ref{cross-polytope example in R^2}   in $\R^2$.   

However, this method no longer works for the cross-polytope in dimensions $d \geq 3$, 
where it is not (yet) known if there is a simple way to go from the FT of the cube to the FT of the cross-polytope. 
\index{cross-polytope}

More generally, one may ask: 
\begin{question}
is there a nice relationship between the FT of a polytope $\P$ and the FT of its polar? 
\end{question}

\item \label{identity thm for two compact sets whose FT's agree on a convergent sequence}
 With regards to Lemma \ref{complete invariance of the FT}, much more is true.  If the Fourier transforms of any two
compact sets $A, B\subset \R^d$ agree on any convergent sequence (with a finite limit point), then $A= B$. The reason is that here $\hat 1_A(\xi)$ and $\hat 1_B(\xi)$ 
are both entire functions of $\xi \in \C^d$, so the proof follows from the identity theorem in complex variables.

\item We note that $P_1({x})$ is defined to be
equal to $0$ at the integers, because its Fourier series naturally converges to the mean of
the discontinuity of the function, at each integer.

\item  It has been known since the work of Riemann that the Bernoulli numbers occur as special values of the Riemann zeta function  (see Exercise \ref{Riemann zeta function, and Bernoulli numbers}).   Similarly, the 
Hurwitz zeta function, defined for each fixed $x >0$ by 
\[
\zeta(s, x) := \sum_{n =0}^\infty \frac{1}{(n+x)^s},
\]
has a meromorphic continuation to all of $\C$, and its special values at the negative integers are the Bernoulli polynomials $B_n(x)$ (up to a multiplicative constant). 

\item  There are sometimes very unusual (yet useful) formulations for the Fourier transform of certain functions.  
Ramanujan  (\cite{Ramanujan1}, eq. (2)) discovered the following remarkable formula
 for the Fourier transform of the Gamma function:
\begin{equation}
\int_{\R}   |\Gamma(a + iy)|^2 e^{-2\pi i \xi y} dy = 
\frac{ \sqrt{\pi} \  \Gamma(a) \Gamma( a + \frac{1}{2})}{\cosh(\pi \xi)^{2a}},
\end{equation}
valid for $a>0$.  

For example, with $a:= \frac{1}{2}$, in the language of this chapter we have the Fourier pair 
$\{   |\Gamma(\frac{1}{2} + iy)|^2,   \frac{\pi}{\cosh(\pi \xi)}   \}$.  But from the 
$\Gamma$-function identity \eqref{Gamma-function identity} below (extended to a  complex variable $s$), it quickly follows that 
$|\Gamma(\frac{1}{2} + iy)|^2= \frac{\pi}{\cosh(\pi y)}$. So this special case of $a:= \frac{1}{2}$ allows us to conclude the interesting fact that $f(y):= \frac{1}{\cosh(\pi x)}$ is a fixed point of the Fourier transform.

\bigskip

\item  \label{Ogsood's book}
I borrowed this joke from \cite{OsgoodBook}, a nice and informal introduction to Fourier analysis.

\end{enumerate}

\section*{Exercises}
\addcontentsline{toc}{section}{Exercises}
\markright{Exercises}

\begin{quote}    
Problems worthy of attack prove their worth by fighting back.

--  Paul  Erd\H os  
 \end{quote}
\medskip

\begin{prob}  \label{transform.of.interval.a.to.b}  $\clubsuit$
Show that the Fourier transform of the closed interval $[a, b]$ is: 
\[
\hat 1_{[a,b]}(\xi) =\frac{ e^{-2\pi i \xi a}   -  e^{-2\pi i \xi b}         }{2\pi i \xi},
\]
for $\xi \not=0$.
\end{prob} 

\medskip
\begin{prob}  \label{transform.of.unit.cube} 
Show that the Fourier transform of the unit cube $C:= [0,1]^d \subset \R^d$ is:
\begin{equation}
\hat 1_{C}(\xi) =  \frac{1}{(2\pi i)^d}  \prod_{k=1}^d \frac{    1- e^{-2\pi i \xi_k}   }{   \xi_k         },
\end{equation}
valid for all $\xi \in \R^d$, except for the union of hyperplanes defined by \\
$H := \left\{  x \in \R^d \bigm |    \xi_1 = 0 \text{ or } \xi_2 = 0  \dots  \text{ or }  \xi_d = 0  \right\}$.
\end{prob}


\medskip
\begin{prob}
Suppose we are given two polynomials $p(x)$ and $q(x)$, of degree $d$.  If there are $d+1$ distinct points
$\{z_1, \dots, z_{d+1}\}$ in the complex plane such that $p(z_k) = q(z_k)$ for $k = 1, \dots, d+1$, show that the two polynomials are identical. (Hint: consider $(p-q)(z_k)$)
\end{prob}

\medskip
\begin{prob} \label{brute force Bernoulli polys}
To gain some facility with generating functions, show by a brute-force computation with Taylor series that
the coefficients on the right-hand-side of equation \eqref{generating function for Bernoulli polynomials}, which are 
called $B_n(x)$ by definition, must in fact be polynomials in $x$.  

In fact, your direct computations will show that for all $n \geq 1$, we have
\[
B_n(x) = \sum_{k=0}^n  {n \choose k} B_{n-k} \ x^k,
\]
where $B_j$ is the $j$'th Bernoulli number.  \index{Bernoulli number}
\end{prob}

\medskip
\begin{prob}    $\clubsuit$ \label{Reflection property for B_n(x)}
Show that for all $n \geq 1$, we have
\[
B_n(1-x) = (-1)^n B_n(x).
\]
\end{prob}

\medskip
\begin{prob}    $\clubsuit$ \label{difference of Bernoulli polys}
Show that for all $n \geq 1$, we have
\[
B_n(x+1) - B_n(x) = n x^{n-1}.
\]
\end{prob}

\medskip
\begin{prob}    $\clubsuit$ \label{derivative of Bernoulli polys}
Show that for all $n \geq 1$, we have
\[
\frac{d}{dx}  B_n(x) = n B_{n-1}(x).
\]
\end{prob}

\medskip
\begin{prob}    $\clubsuit$ \label{historical origin of Bernoulli poly}
Prove that:
\[
\sum_{k=0}^{n-1} k^{d-1} = \frac{    B_d(n) - B_d   }{    d   },
\]
for all integers $d \geq 1$ and $n \geq 2$. 
\end{prob}

\medskip
\begin{prob} $\clubsuit$  
 \label{Bernoulli Polynomials}  \index{Bernoulli polynomial}
Show that the periodic Bernoulli polynomials $P_n(x) := B_n(\{ x \})$, for all $n\geq 2$, have the following Fourier series:
\begin{equation}
P_n(x) = -\frac{n!}{(2\pi i)^n} \sum_{k \not=0}   \frac{e^{2\pi i k x}}{k^n},
\end{equation}
valid for all $x \in \R$.   For $n \geq 2$, these series are absolutely convergent.  We note that from the definition above, $B_n(x) = P_n(x)$ when  $x\in (0,1)$.
\end{prob}

\medskip
\begin{prob}  \label{Raabe's identity for { } via Fourier series}
Show that the greatest integer function $\floor{x}$  (often called the `floor function') enjoys the
property:
\[
\sum_{k=0}^{N-1} \floor{ x +  \frac{k}{N}       }  = \floor{Nx},
\]
for all $x\in \R$, and all positive integers $N$.
\end{prob}

\medskip
\begin{prob}
Show that the Bernoulli polynomials enjoy the following identity, proved by Joseph Ludwig Raabe in 1851:
\[
B_n(Nx) = N^{n-1} \sum_{k=0}^{N-1}  B_n\left(  x + \frac{k}{N} \right),
\]
for all $x\in \R$, all positive integers $N$, and for each $n \geq 1$.

Notes.  Such formulas, in these last two exercises, are also called ``multiplication Theorems'', and they hold
for many other functions, including the Gamma function, the dilogarithm, the Hurwitz zeta function, the cotangent, and many more.
\end{prob}

\medskip
\begin{prob}    $\clubsuit$  \label{another definition for Bernoulli polynomials}
\index{Bernoulli polynomial}
Here we give a different method for defining the Bernoulli polynomials, based on the following three properties that they enjoy:
\begin{enumerate}
\item $B_0(x) = 1$.
\item For all $n \geq 1, \frac{d}{dx} B_n(x) = n B_{n-1}(x)$.
\item  For all $n \geq 1$, we have $\int_0^1 B_n(x) dx = 0$.
\end{enumerate}
Show that the latter three properties imply the original defining property of the Bernoulli polynomials
\eqref{generating function for Bernoulli polynomials}.
\end{prob}

\medskip
\begin{prob}   
Here is a more explicit, useful recursion for computing the Bernoulli polynomials.  Show that 
\[
\sum_{k=0}^{n-1}  {n \choose k} B_k(x)  = n x^{n-1},
\]
for all $n \geq 2$.
\end{prob}

\medskip
\begin{prob}    \label{B_7} 
Use the previous exercise, together with the known list the first $6$ Bernoulli polynomials that appear
in equation \ref{B_6}, to compute $B_7(x)$. 
\end{prob}

\medskip
\begin{prob}    \label{odd Bernoulli numbers}
Show that for odd $k \geq 3$, we have $B_k = 0$. 
\end{prob}

\medskip
\begin{prob}    \label{ Bernoulli numbers alternate in sign}
Show that  the even Bernoulli numbers alternate in sign.   More precisely,  show that 
\[
(-1)^{n+1} B_{2n}  \geq 0,
\]
for each positive integer $n$. 
\end{prob}

\medskip
\begin{prob} \label{vanishing identity for Beroulli numbers}
Show that the Bernoulli numbers enjoy the recursive property:
\[
\sum_{k=0}^n  {n+1 \choose k}B_k = 0,
\]
for all $n \geq 1$.
\end{prob}

\medskip
\begin{prob} \label{asymptotics for Beroulli numbers}
Show that the Bernoulli numbers enjoy the following asymptotics:
\[
B_{2n} \sim  2 \frac{(2n)!}{(2\pi)^{2n}}
\]
as $n\rightarrow \infty$.  Here we are using the usual notation for asymptotic functions, namely that
$f(n) \sim g(n)$ as $n\rightarrow \infty$ if $\lim_{n\rightarrow \infty}  \frac{f(n)}{g(n)} \rightarrow 1$.
\end{prob}

\medskip
\begin{prob}    \label{Fresnel} $\clubsuit$
Show that the following  integrals converge and have the closed forms: 
\begin{align}
\int_{-\infty}^\infty \cos(x^2) dx &= \sqrt{\frac{\pi}{2}},   \\ 
\int_{-\infty}^\infty \sin(x^2) dx &= \sqrt{\frac{\pi}{2}}.
\end{align}
Notes.   These integrals are called Fresnel integrals, and they are related to the Cornu spiral, which was created by Marie Alfred Cornu.    Marie used the spiral
as a tool for computing diffraction patterns that arise naturally in optics.
\end{prob}

\medskip
\begin{prob}     
 Prove the following Gamma function identity, using the sinc function:
\begin{equation}\label{Gamma-function identity}
 \frac{\sin(\pi x)}{\pi x} =  \frac{1}{\Gamma(1+x)\Gamma(1-x)},
\end{equation}
 for all $x \not\in \Z$. 
 
 Notes.  This identity is often called Euler's reflection formula.  $\Gamma(x):= \int_0^\infty e^{-t} t^{x-1} dt$ is by definition 
  the Gamma function, 
 where the integral converges for all $x > 0$ (see Section \ref{Volume of the ball, the Gamma function} for more on the $\Gamma$ function). 
\end{prob}

\medskip
\begin{prob}
\label{retreiving volume of 2d.crosspolytope}
$\clubsuit$ Using the formula for the Fourier transform of the   $2$-dimensional cross-polytope $\Diamond$, derived in the text, 
namely   
\[
\hat 1_{\Diamond}(\xi) = 
-\frac{1}{\pi^2}  \left(
     \frac{     \cos(2\pi \xi_1)  -     \cos(2\pi \xi_2)     }{  \xi_1^2 -  \xi_2^2}
  \right), 
  \]
find the area of $\Diamond$ by letting $\xi \rightarrow 0$ in the latter formula.
\end{prob}

\medskip
\begin{prob}
\label{Elementary bounds for sin(x), sinc(x)}
Some elementary but very useful bounds for trig functions are developed here.  
\begin{enumerate}[(a)]
\item   Prove that
\[
\frac{2}{\pi} < \frac{\sin x}{x} \leq 1,
\]
where the left inequality holds for $ 0 < x < \frac{\pi}{2}$, and the right inequality holds for $x\in \R$.
\item   \label{Elementary trig bounds, part b}
Prove that
\[
     \frac{2x}{\pi} \leq |1-e^{ix} | \leq |x|,
\]
where the left inequality holds for $|x| \leq \pi$, and the right inequality holds for $x\in \R$. 
\item  Prove that 
\[
     \frac{2x^2}{\pi^2} \leq    |1-\cos x | \leq \frac{x^2}{2},
\]
where the left inequality holds for $|x| \leq \pi$, and the right inequality holds for $x\in \R$.
\end{enumerate}
\end{prob}

\medskip
\begin{prob}
\label{divergence of |sinc|}
$\clubsuit$
Show that 
$\int_{-\infty}^\infty    \Big| \frac{\sin(\pi x)}{\pi x} \Big| dx = \infty$.
\end{prob}

Notes.  Once we have the inverse Fourier transform and its consequences at our disposal, this exercise will become trivial, and much more general - see Corollary \ref{the FT of an indicator function is not in L^1}.

\medskip
\begin{prob}
\rm{
There are (at least) two different ways of periodizing a given function $f:\R \rightarrow \C$ 
with respect to $\Z$. 
First, we can define $F_1(x) := f(\{x\})$, so that $F_1$ is periodic on $\R$ with period~$1$.
Second, we may also define $F_2(x) := \sum_{n\in \Z} f(x+n)$, which is also a periodic function on $\R$ with period $1$. 

Find an absolutely integrable (meaning that $\int_\R |f(x)| dx$ converges) function $f$ for which these two functions are not equal: $F_1 \not= F_2$.

Notes. \ In Chapter \ref{Fourier analysis basics}, we will see that the latter function $F_2(x) := \sum_{n\in \Z} f(x+n)$ captures a lot more information about $f$, and often captures all of $f$ as well.  
}
\end{prob}

\medskip
\begin{prob} \label{symmetrized parallelepiped}  \index{parallelepiped}
Given linearly independent vectors $w_1, \dots, w_d \in \R^d$, let $v:= -\frac{w_1 + \cdots + w_d}{2}$, and define
$Q:= \{ v+  \sum_{k=1}^d   \lambda_k w_k \mid  \text{ all }  \lambda_k \in [0, 1] \}$, a parallelepiped. 
Show that $Q$ is symmetric about the origin, in the precise sense that
$x \in Q \iff -x \in Q$.
\end{prob}


\medskip
\begin{prob}
\label{polars of each other}
$\clubsuit$ Show that the $d$-dimensional cross-polytope $\Diamond$ and the cube
$\square:= [-1, 1]^d$ are polar to each other.
\end{prob}

\medskip
\begin{prob} \label{polytope with 5 vertices}
\begin{enumerate}[(a)]
\item
Suppose $C\subset \R^3$ is a convex polytope with $5$ vertices.   Prove that at least one of the vertices of $C$ has degree $4$.  \\
\item Construct a convex polytope $\P\subset \R^3$ with $5$ vertices, such that all of its 
 vertices have degree $4$. 
 \end{enumerate}
\end{prob}

\medskip
\begin{prob}
\label{2-dim'l formula for triangle}
Prove the following  $2$-dimensional integral formula:
\begin{align}
 \int_{\lambda_1, \lambda_2 \geq 0    \atop   \lambda_1 + \lambda_2 \leq 1}
                     e^{a \lambda_1  }  e^{b \lambda_2  } d\lambda_1 d\lambda_2
                     = \frac{  b e^a -  a e^b }{ab(a-b)} +\frac{1}{ab},
\end{align}
valid for all $a, b \in \C$ such that $ab(a-b) \not=0$. 
\end{prob}

\medskip
\begin{prob}
\label{FT of a general simplex, brute-force}
Using the ideas of Example  \ref{FT of a general triangle}, prove (by induction on the dimension) that  the Fourier transform of a general $d$-dimensional simplex $\Delta \subset \R^d$  is given by:
\begin{equation}
\hat 1_{\Delta}(\xi) = (\vol \Delta) d!
\sum_{j=1}^N   
\frac{e^{-2\pi i \langle v_j, \xi \rangle}}{\prod_{1\leq  k \leq d} \langle v_j-v_k, \xi \rangle   }
[k\not=j],
\end{equation}
for all $\xi \in \R^d$, where the vertex set of $\P$ is $\{ v_1, \dots, v_N\}$.
\end{prob}

\medskip
\begin{prob}[Abel summation by parts]  $\clubsuit$    \label{Abel summation by parts}
Here we prove the straightforward but very useful technique of Niels Abel, called {\bf Abel summation by parts}. 
\index{Abel, Niels}   \index{Abel summation by parts}
Suppose we are given two sequences $\{a_n\}_{n=1}^\infty$, and $\{b_n\}_{n=1}^\infty$.   We define the finite partial sums
$B_n:= \sum_{k=1}^n b_k$.   Then we have
\begin{equation} \label{actual Abel summation}
 \sum_{k=1}^n a_k b_k     = a_n B_n + \sum_{k=1}^{n-1} B_k(a_k - a_{k+1}),
\end{equation}
for all $n\geq 2$.
\end{prob}

Notes. \ Using the forward difference operator, it's easy to recognize identity \eqref{actual Abel summation}
as a discrete version of integration by parts. 

\medskip
\begin{prob}[{\bf Dirichlet's convergence test}]   $\clubsuit$    \label{Dirichlet's convergence test}
\index{Dirichlet's convergence test}
Suppose we are given a real sequence $\{a_n\}_{n=1}^\infty$, and a complex sequence 
$\{b_n\}_{n=1}^\infty$,
such that 
 \begin{enumerate}[(a)]
 \item    $\{a_n\}$ is monotonically decreasing to $0$, and  
 \item  $| \sum_{k=1}^n  b_k |  \leq M$, for some positive  constant $M$, and all $n \geq 1$.  
\end{enumerate}
Then $\sum_{k=1}^\infty a_k b_k$ converges.
\end{prob}

\medskip
\begin{prob} \label{first Dirichlet kernel}
Prove that  for all $x \in \R-\Z$, we have the following important identity, called the ``Dirichlet kernel'', 
\index{Dirichlet kernel}
named after Peter Gustav Lejeune Dirichlet:
\begin{equation}\label{Dirichlet Kernel, identity}
  \sum_{k= -n}^n     e^{2\pi i k x} =  \frac{\sin \left( \pi x(2n + 1) \right) }{\sin(\pi x)}.
\end{equation}
\end{prob}

Notes. \ An equivalent way to write \eqref{Dirichlet Kernel, identity} is clearly:
\[
1+ 2 \sum_{k=1}^n \cos\left( 2\pi  k x \right) =  \frac{\sin \left( \pi x(2n + 1) \right) }{\sin(\pi x)}.
\]

\medskip
\begin{prob} \label{exponential sum bound}
Prove that we have the bound on the following exponential sum:
\begin{equation}
 \left |  \sum_{k= 1}^n     e^{2\pi i k x} \right |  \leq     \frac{1}{  |  \sin(\pi x)  |  },
\end{equation}
for any fixed  $x \in \R-\Z$, and for all $n \in \Z_{>0}$.
\end{prob}

\medskip
\begin{prob} $\clubsuit$    \label{rigorous convergence of P_1(x)}
Prove that $\sum_{m = 1}^\infty  \frac{e^{2\pi i m a}}{m}$ converges, given any fixed $a \in \R - \Z$.

Notes. \ We see that, although $\sum_{m = 1}^\infty  \frac{e^{2\pi i m a}}{m}$ does not converge absolutely, 
Abel's  summation formula \eqref{actual Abel summation}    gives us
\[
\sum_{k = 1}^n  \frac{e^{2\pi i k a}}{k} = \frac{1}{n}\sum_{r=1}^n e^{2\pi i r a} 
+  \sum_{k=1}^{n-1}      \Big(    \sum_{r=1}^k   e^{2\pi i r a}   \Big)     \frac{1}{k(k+1)},
\]
and the latter series {\bf does converge absolutely}, as $n \rightarrow + \infty$.  So we see that Abel summation transforms one series (that barely converges at all)  into 
another series that converges more rapidly. 
\end{prob}


\medskip
\begin{prob} $\clubsuit$
\label{integral of sinc is 1}
{\rm
 Here we'll prove that
\begin{equation} \label{exercise:the Dirichlet integral}
   \int_{-\infty}^\infty      \frac{ \sin(\pi t) }{\pi t} dt = 1,
\end{equation}
in the sense that $\int_0^ \infty      \frac{ \sin t }{ t} dt =\frac{\pi}{2}$. 
The integral \eqref{exercise:the Dirichlet integral} 
is sometimes called ``the Dirichlet integral''.  Comparing this Dirichlet integral with
Exercise \ref{divergence of |sinc|}, we see that there is something subtle going on here.   
We'll end up proving something slightly more general here:
\[
 \int_0^\infty e^{-st}   \frac{ \sin t }{ t} dt = \frac{\pi}{2} - \tan^{-1} s, 
\]
for all $s>0$.

\begin{enumerate}[(a)]
\item Define 
\begin{equation} \label{Laplace transform of sinc}
F(s):= \int_0^\infty e^{-st}   \frac{ \sin t }{ t} dt, 
\end{equation}
for each $s>0$.   Justify differentiation under the integral sign, and show that 
\[
\frac{dF}{ds}  = -\int_0^\infty e^{-st} \sin t dt, 
\]
\item Show that  $\int_0^\infty e^{-st} \sin t dt = \frac{1}{1+s^2}$.
\item Show that $F(s) = C - \tan^{-1} s$, and then show that the constant $C= \frac{\pi}{2}$.
\item Prove that $F$ is a continuous function  of $s\in \R_{>0}$, and finally prove that 
\[
\lim_{s\rightarrow 0} F(s) = \frac{\pi}{2}, 
\]
which is the desired result (Here you might want to integrate by parts first, and then use the Dominated convergence theorem).
\end{enumerate}
}
\end{prob}

Notes.  There are many proofs of this famous identity \eqref{exercise:the Dirichlet integral}, and although the method of contour integration is arguably the most straightforward, here we are only assuming 
knowledge of some real analysis.  The expression in \eqref{Laplace transform of sinc} is also known as the Laplace transform of the sinc function, and it is a variation of the Fourier transform that we will revisit when studying similar transforms of cones in 
Section \ref{Fourier Laplace transforms of cones}.

\medskip
\begin{prob} $\clubsuit$
\label{rigorous inversion formula for sinc}
\rm{
Here we give a rigorous proof of the tricky fact that for all $x\in \R$, we have
\begin{equation*}
\lim_{N\rightarrow \infty} \int_{-N}^N   \frac{\sin(\pi \xi)}{\pi \xi}  e^{-2\pi i \xi x} d\xi = 1_ {[-\frac{1}{2}, \frac{1}{2}]}(x),
\end{equation*}
following an approach taken by S. Bochner \cite{BochnerBook}.
We begin by noticing that this integral can be easily reduced to a real-valued integral:
\begin{equation*}
 \int_{-N}^N   \frac{\sin(\pi \xi)}{\pi \xi}  e^{-2\pi i \xi x} d\xi
  = \int_{-N}^N   \frac{\sin(\pi \xi)}{\pi \xi}  \cos( 2\pi \xi x) d\xi,
\end{equation*}
because for each  $x\in \R$,  $\int_{-N}^N   \frac{\sin(\pi \xi)}{\pi \xi}  \sin( 2\pi  \xi x) d\xi = 0$,  owing to the oddness of the integrand.
\begin{enumerate}[(a)]
\item  \label{tricky middle part with alpha}
Using the result from Exercise \ref{integral of sinc is 1}, prove that
\[
\lim_{N\rightarrow \infty} \int_{-N}^N    \frac{\sin(\pi   \alpha t)}{\pi t} dt= 
 \begin{cases}  
\ \  1    &      \mbox{if } \alpha >0,  \\ 
\  \ 0  &        \mbox{if }  \alpha =0, \\
-1  &        \mbox{if }  \alpha <0.
\end{cases}
\]
\item  Finish up by using \  $2\sin t \cos(\alpha t) = \sin(1- \alpha)t   + \sin(1+\alpha)t$, thereby showing that the desired
integral 
\[
\lim_{N\rightarrow \infty}  \int_{-N}^N   \frac{\sin(\pi t)}{\pi t}  \cos( 2\pi t x) dt
\]
 reduces  to part \ref{tricky middle part with alpha}.
\end{enumerate}
}
\end{prob}


\chapter{\blue{The basics of Fourier analysis} }
 \label{Fourier analysis basics}
 \index{Fourier analysis}

\begin{quote}
``If a function is periodic, then we should try to expand it into its Fourier series, and wonderful things will begin to happen....."

-- Erich Hecke
\end{quote}

 \begin{quote}                           
    ``. . .   Fourier's great mathematical poem.''
    
  -- William Thomson Kelvin    \index{Kelvin, William Thomson}
  
   [Referring to Fourier's mathematical theory of the conduction of heat]

   \end{quote}

\begin{figure}[htb]
\begin{center}
\includegraphics[totalheight=1.4in]{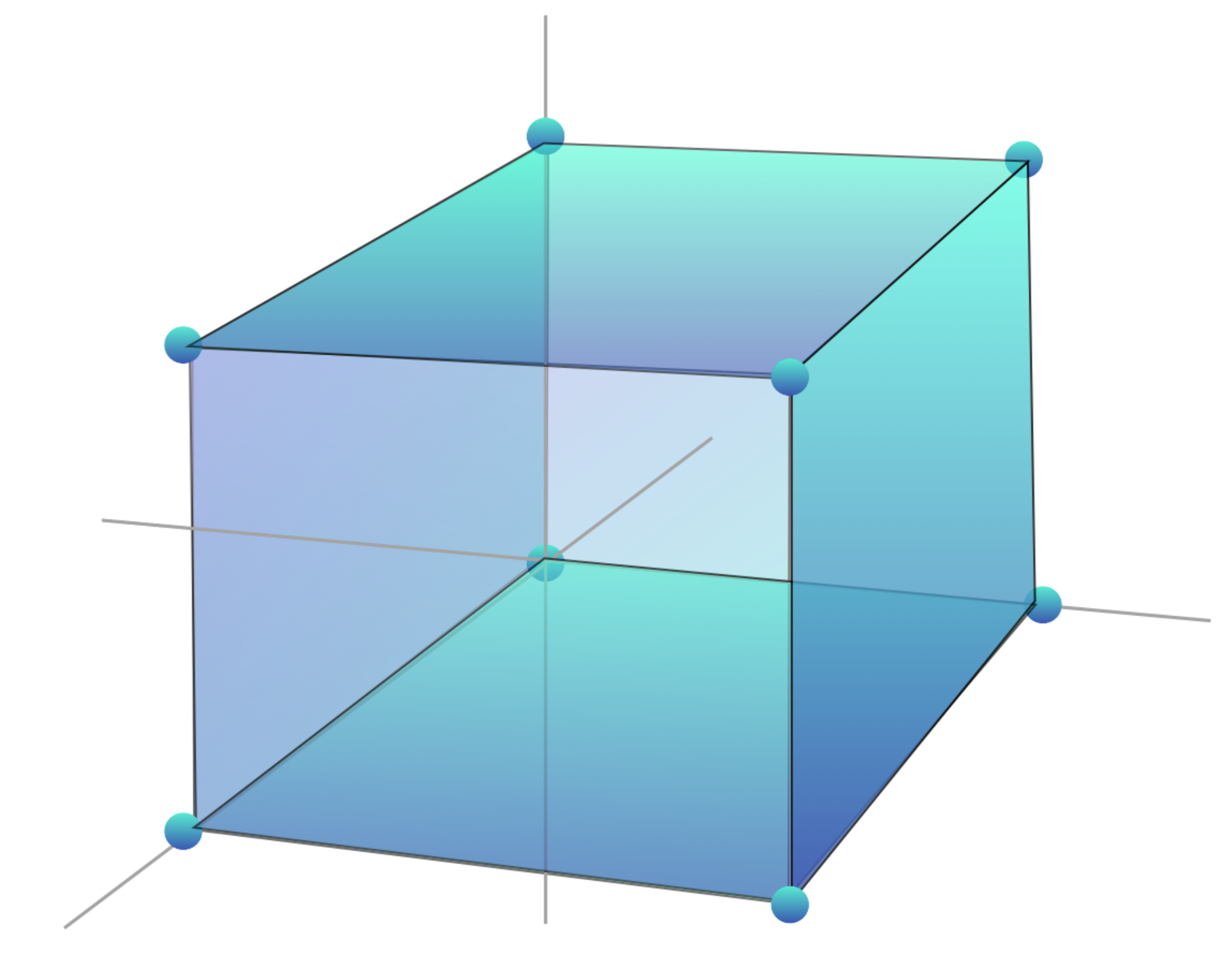}
\end{center}
\caption{The unit cube $\square:= [0,1]^3$, in $\R^3$, which tiles the space by translations.
Which other polytopes tile by translations?  How can we make mathematical use of such tilings?
In particular, can we give an explicit basis of exponentials for functions defined on $\square$?}    \label{crosspic}
\end{figure}

 \bigskip
 \section{Intuition}
 Because we will use tools from Fourier analysis throughout, we introduce them here as an outline of the field, 
 with the goal of {\it applying} them to the discrete geometry of polytopes, lattices, and their interactions.
 We will sometimes introduce a concept by using an intuitive argument, which we call ``fast and loose", but after such 
 an intuitive argument, we state the precise version of the corresponding theorem. 
In this chapter, we will sometimes point to the literature for some of the proofs.    
  
Our goal is to use the necessary tools of Fourier analysis in order to tackle problems in the enumerative combinatorics of polytopes, 
  in number theory, discrete geometry, and in some other fields.
 We emphasize that the Poisson summation formula allows us to {\bf discretize  integrals}, 
 in a sense that will be made precise in later chapters.
 
 One pattern that the reader may have already noticed, among all of the examples of Fourier transforms of polytopes computed thus far, is that each of them is a linear combination of a very special kind of rational function of $\xi$, multiplied by a complex exponential that involves a vertex of the polytope:
 \begin{equation}\label{structure 1}
 \hat 1_\P(\xi) = 
 \sum_{k=1}^M  \frac{1}{\prod_{j=1}^d \left\langle   \omega_{j,k}(v_k),  \xi  \right\rangle} \, 
 e^{2\pi i \langle v_k, \xi \rangle},
 \end{equation}
 where the vertices of $\P$ are $v_1, \dots, v_N$, and where $M \geq N$.   We observed that in all of our examples thus far, the  denominators  
 are in fact products of linear forms, as in \eqref{structure 1}.
 We will  be able to see some of the more precise geometric structure for  these products of linear forms, which 
come from the edges of the polytope, once we learn more about Fourier-Laplace transforms of cones.

 It is rather astounding that every single fact about a given polytope $\P$  is somehow 
 hiding inside these {\bf rational-exponential} functions given by \eqref{structure 1}, due to the fact that the Fourier transform $\hat 1_\P$ is a complete invariant (Lemma \ref{complete invariance of the FT}).
 
 Finally, it is worth mentioning that not every fact in this chapter is necessary for the comprehension of the rest of the book.  The reader is advised to learn just some of this chapter, and as she/he reads the rest of the book, periodically revisit this chaper.

 \bigskip
 \section{Introducing the Fourier transform on $L^1(\R^d)$}
 
In the spirit of bringing the reader very quickly up to speed, regarding the applications of Fourier analytic tools, we outline the basics 
of the field, and prove some of them.   Nowadays, there are many good texts on Fourier analysis, and the reader is encouraged to peruse some of these books (see  Note \ref{Fourier books}). 

Unless otherwise stated, all of our functions will have the form $f:\R^d \rightarrow \C$. 
One of the most useful tools for us is the Poisson summation formula. 
We provide several versions of Poisson summation, each of 
which uses a different set of sufficient conditions. 
 
 As we will see, the Fourier transform \index{Fourier transform}  is a very friendly creature, allowing us to travel back and forth between the ``space domain''  and the ``frequency domain'' to obtain many useful results. 
 The readers who are already familiar with basics of Fourier analysis may easily skip 
this chapter without impeding their understanding of the rest of the book.  
Although we enjoy thinking about the warm and cozy Hilbert spaces $L^2(\R^d)$ and $L^2([0, 1]^d)$,  there are exotic Fourier series that are pointwise divergent, and yet represent continuous functions, 
a whole field onto itself.   We won't go there. 
However, the very basic convergence issues are still important for us as well, and we will study them because we want to get the reader up and running.


The function space that immediately come up very naturally is the
the space of {\bf absolutely integrable functions} on $\R^d$:
 \[
 L^1(\R^d) :=\left\{  f: \R^d \rightarrow \C \bigm |   \  \int_{\R^d} |f(x)| dx < \infty \right\}.
 \]
Secondly, the space of
 {\bf square-integrable functions} on $\R^d$  is defined by:
 \[
 L^2(\R^d) := \left\{  f: \R^d \rightarrow \C \bigm |   \  \int_{\R^d} |f(x)|^2 dx < \infty \right\}.
 \]
The usual theory of Fourier transforms progresses by first defining the Fourier transform for functions belonging to $L^1(\R^d)$, which is quite a natural condition, and then later extending the Fourier transform to the $L^2(\R^d)$ space by taking appropriate limits.  We initially restrict attention to functions $f \in L^1(\R^d)$.

 \begin{figure}[htb]
\begin{center}
\includegraphics[totalheight=2.8in]{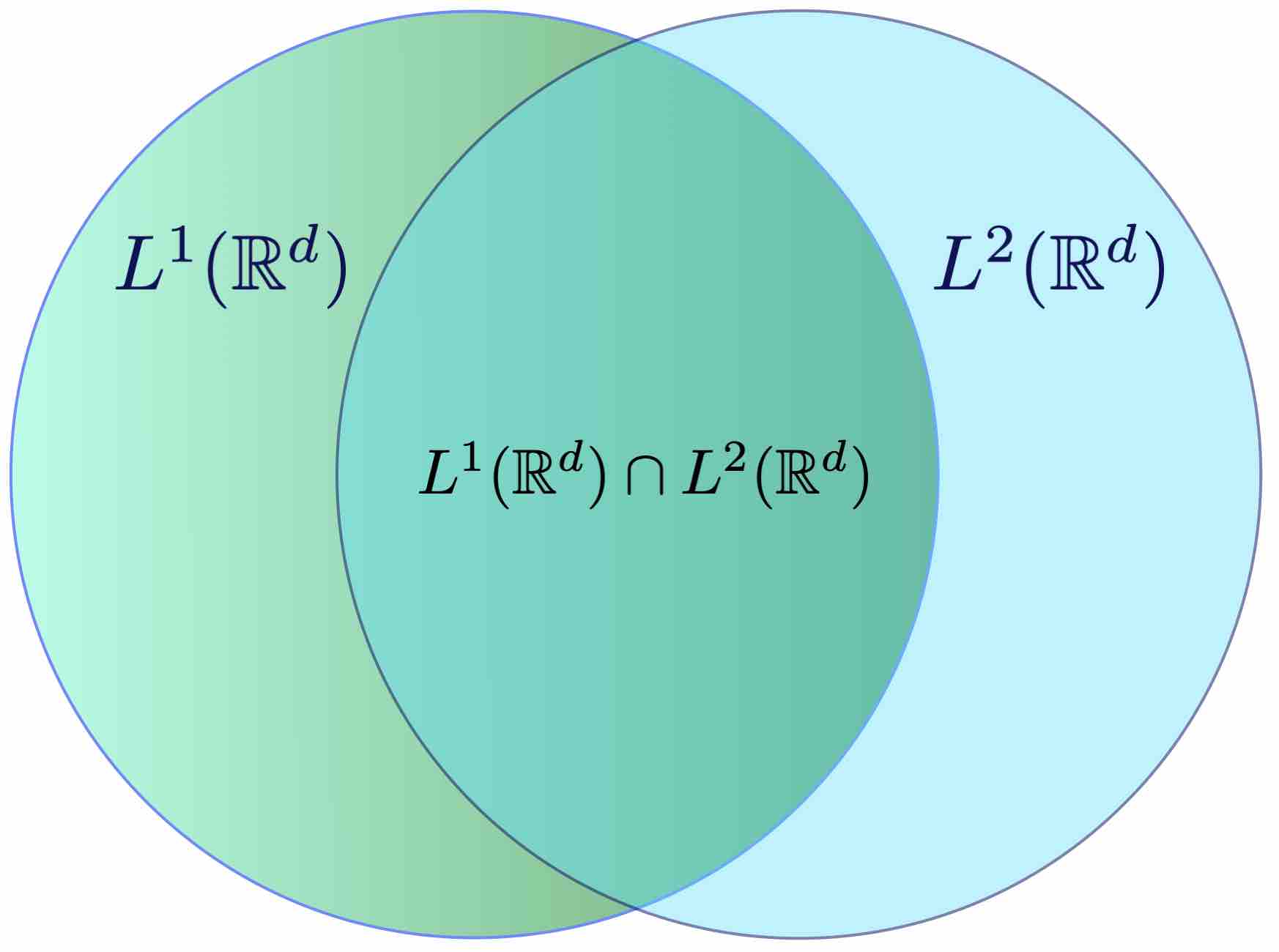}
\end{center}
\caption{Neither of the function spaces $L^1(\R^d)$ and $L^2(\R^d)$ is contained in the other, as in Example \ref{example:L1andL2}.}
    \label{L1 and L2}
\end{figure}
 
There are many fascinating facts  about all of these functions spaces.   For practice, let's ask:
\begin{question} \rm{[Rhetorical]} 
\label{Question: product of two L^2 functions}
Given two functions  $f, g \in L^2(\R^d)$, is their product always in $L^1(\R^d)$?
\end{question}
   Well, we have the {\bf Cauchy-Schwartz inequality} for the 
 Hilbert space $L^2(\R^d)$:
\begin{equation}\label{product of two L^2 functions is L^1}
\int_{\R^d} |f(x) g(x)|dx \leq  
      \left(\int_{\R^d} |f(x)|^2 dx  \right)^{\frac{1}{2}}   \left(\int_{\R^d} |g(x)|^2 dx\right)^{\frac{1}{2}} < \infty,
\end{equation}
the latter inequality holding by the assumption  $f, g \in L^2(\R^d)$.  So the product $f(x) g(x)$ is indeed in $L^1(\R^d)$, answering Question \ref{Question: product of two L^2 functions}
 in the affirmative.  This is the first sign that there  are fascinating links between $L^1$ functions and $L^2$ functions.  In fact, this metaphorical dance between $L^1(\R^d)$ and $L^2(\R^d)$ is simply too useful to ignore, so we will use study some of the interactions between these two spaces, from first principles (see Section \ref{Section: more relations between L^1 and L^2} below, for example). 

The utility of the Cauchy-Schwarz inequality should never be underestimated, and holds in greater generality.
\begin{lem}   \label{Cauchy-Schwartz}
Let $V$ be an inner product space, with the inner product $\langle x, y\rangle$. 
Then the following 
Cauchy-Schwarz inequality holds:
\[
\left|   \langle x, y   \rangle  \right|  \leq \|x\| \|y\|
\]
for all $x, y \in V$.   Moreover, equality holds $\iff$  $x$ and $y$ are linearly dependent. 
In addition, the function 
\[
\|  x  \|:=  \sqrt{ \langle x, x \rangle}
\]
  is a norm on $V$.
\hfill $\square$
\end{lem}
(For a proof see \cite{EinsiedlerWardBook}, Prop. 3.2.     
Appendix \ref{Appendix:bump functions and inner products}   has some related material) 
So we can conclude from Lemma \ref{Cauchy-Schwartz}  that every inner product space is also a normed vector space.  
The converse is false, though, in the sense that there are normed vector spaces, whose norm does not arise from any inner product.   
One fascinating example of such a space is $L^1(\R^d)$, which is not a Hilbert space, as we now easily show by exhibiting a counter-example
to the Cauchy-Schwarz inequality.

\bigskip
\begin{example}  \label{cool CS counterexample}
\rm{
We claim that the Cauchy-Schwarz inequality
is  false in $L^1(\R)$.   If the Cauchy-Schwarz inequality was true here, then \eqref{product of two L^2 functions is L^1} would be valid 
for all functions
$f, g \in L^1(\R)$.  But as a counterexample,  let 
\[
f(x):= 1_{(0,1)}(x) \frac{1}{\sqrt x}.
\]
 It's easy to see that $f \in L^1(\R)$:
 \[
\int_{\R}   1_{(0,1)}(x) \frac{1}{\sqrt x}  dx  =  \int_0^1  \frac{1}{\sqrt x}  dx = \frac{1}{2}.
\]
But $ \int_{\R} f(x)  \cdot   f(x)  dx  = \int_0^1 \frac{1}{x} dx$ diverges, so that we do not have a Cauchy-Schwarz inequality in $L^1(\R)$, because here both the left-hand-side and the right-hand-side of such an inequality do not even converge.  
}
\hfill $\square$
\end{example}

We say that $f$ is bounded on a measurable set $S\subset \R^d$ by a constant $M>0$,
 if $|f(x)| < M$, for all $x \in S$.
In the opposite direction of Example \ref{cool CS counterexample}, if two functions $f, g$ are bounded on $\R^d$, and absolutely integrable on $\R^d$, then we do have a Cauchy-Schwartz inequality
for the pair $f, g$, and we let the reader enjoy its verification.

We'll see that despite the fact that $L^1(\R^d)$ is not a Hilbert space, it does have a very beautiful structure, namely it is a Banach algebra (Lemma \ref{convolution theorem}).

\bigskip
\section{The triangle inequality for integrals}

An easy but extremely important inequality is the  {\bf triangle inequality for integrals}, as follows. 
\begin{thm}   \label{triangle inequality for integrals}  
 For any $f\in L^1(\R^d)$, and any measurable subset $S \subset \R^d$,
we have:
\begin{equation}  \index{triangle inequality for integrals}
\Big|  \int_S f(x) dx  \Big|  \leq       \int_S | f(x) | dx.
\end{equation}
\end{thm}
\begin{proof}
Letting $z:=  \int_S f(x) dx \in \C$, we may write  $|z| = \alpha z$, for a (unique) complex $\alpha$ on the unit circle.  
We let $u$ be the real part of $\alpha f:= u + iv $, so that $u \leq \sqrt{u^2 + v^2} = |\alpha f| = |f|$.
Altogether, we have:
\begin{equation} \label{string of inequalities in the triangle inequality}
\Big|  \int_S f(x) dx  \Big| = \alpha  \int_S f(x) dx  =    \int_S \alpha f(x) dx  =  \int_S u(x) dx 
\leq  \int_S   |f(x)| dx. 
\end{equation}
In the third equality, we used the fact that  $\int_S \alpha f(x) dx$ is real, which follows
from the first two equalities:
 $\int_S \alpha f(x) dx =  \Big|  \int_S f(x) dx  \Big|$. 
\end{proof}
Although Theorem \ref{triangle inequality for integrals} seems innocently trivial, it is sometimes quite powerful.
\begin{cor}   If $f$ is bounded on a measurable set $S\subset \R^d$ by a constant $M>0$, then:
\begin{equation} 
\Big|  \int_S f(x) dx  \Big|  \leq    \vol (S) \cdot M.
\end{equation}
\begin{proof}
\[
\Big|  \int_S f(x) dx  \Big|  \leq       \int_S | f(x) | dx  \leq   \int_S M dx  =      \vol (S) \cdot M,
\]
where the first inequality uses the triangle inequality for integrals, namely Theorem \eqref{triangle inequality for integrals}, and the second inequality uses
the boundedness assumption on $f$.
\end{proof}
\end{cor}
What about the equality case in Theorem \ref{triangle inequality for integrals}? 
Luckily there is a very satisfying answer, which turns out to be so useful that it merits its own Corollary.  We'll keep using the same notation as in the proof of 
Theorem \ref{triangle inequality for integrals}.

\begin{cor}[Equality conditions for the triangle inequality]
\label{cor:equality case of the triangle inequality}
Let $f \in L^1(\R^d)$, fix any measurable subset $S\subset \R^d$, and suppose that 
\[
\left | \int_S f(x) dx \right | = \int_S | f(x) | dx.
\]
Then we have $| f(x) | = \alpha f(x)$  for some complex number $\alpha$ on the unit circle,
and for almost all $x\in \R^d$.
\end{cor}
\begin{proof}
Returning to \eqref{string of inequalities in the triangle inequality} in the proof of Theorem \ref{triangle inequality for integrals}, our assumption of equality now gives us
\begin{equation}\label{vanishing of triangle inequality integral}
\int_S \big( |f(x)| - u(x) \big) dx = 0,
\end{equation}
We also have $|f(x)| - u(x) \geq 0$, so by \eqref{vanishing of triangle inequality integral} we now have $|f(x)| - u(x) =0$, almost everywhere.  Now we recall that $u$ is the real part of $\alpha f$, with $\alpha$ a complex number on the unit circle.  So we have
$|f(x)| = \Re(\alpha f(x))$ almost everywhere.  Since $|f(x)| = |\alpha f(x)| $, we see that
$ \Re \left (\alpha f(x) \right ) = |\alpha f(x)|$ almost everywhere.  In other words,
 $\alpha f(x) = | \alpha f(x) | = |f(x)|$ almost everywhere.
\end{proof}

Next, let's show that
\[
L^1(\R^d) \not\subset L^2(\R^d), \text{ and } L^2(\R^d) \not\subset L^1(\R^d),
\]
confirming the validity of set intersections in Figure \ref{L1 and L2}.    We'll do it for $d=1$, but the idea works for any dimension.  
\begin{example}\label{example:L1andL2}
\rm{
Let's define
\[
f(x):= \begin{cases}
x^{-\frac{2}{3}} & \text{ if } 0 < x < 1, \\
0  &  otherwise. 
\end{cases}
\]
Then $f \in L^1(\R)$, because $\int_\R |f(x)|dx := \int_0^1 |x^{-\frac{2}{3}}| dx = 3 x^{\frac{1}{3}} \Big |_0^1 = 3$. 
But  $f \not\in L^2(\R)$, because 
$\int_\R |f(x)|^2 dx   =  \int_0^1 |x^{-\frac{4}{3}}| dx =  -3 x^{-\frac{1}{3}} \Big |_0^1 = \infty$.  So $L^1(\R) \not\subset L^2(\R)$.

On the other hand, if we consider 
\[
g(x):= \begin{cases}
x^{-\frac{2}{3}} & \text{ if } x>1, \\
0  &  otherwise, 
\end{cases}
\]
then  $g \in L^2(\R)$, because 
$\int_\R |g(x)|^2 dx   =  \int_1^\infty |x^{-\frac{4}{3}}| dx =  -3x^{-\frac{1}{3}} \Big |_1^\infty = 3$.
But  $g \not\in L^1(\R)$ because
 $\int_\R |g(x)| dx = \int_1^\infty x^{-\frac{2}{3}} dx =3 x^{\frac{1}{3}} \Big |_1^\infty = \infty$. So $L^2(\R) \not\subset L^1(\R)$.
}
\hfill $\square$
\end{example}

\medskip
We've defined the Fourier transform before, and we remind the reader that for any function  $f\in L^1(\R^d)$,
the {\bf Fourier transform} of $f$ \index{Fourier transform} is
\begin{equation}
\hat f(\xi) := \int_{\R^d} f(x) e^{-2\pi i \langle x, \xi \rangle} dx.
\end{equation}
Where does this definition really come from?  One motivation comes from the inner product for functions
(in $L^2(\R^d)$),
 where we project a function $f$ onto each exponential function:
\[
\langle f, e^{2\pi i \langle x, \xi \rangle}  \rangle := \int_{\R^d} f(x) e^{-2\pi i \langle x, \xi \rangle} dx.
\]
Another motivation comes from the proof of the Poisson summation formula - eq. \eqref{finally, the FT} below, which
shows a crucial connection between the Fourier transform of $f$ and the Fourier coefficients of the periodized function
$\sum_{n\in \Z^d} f(x+n)$.  

One of the first things we might notice is:
\begin{claim}
The Fourier transform is a bounded linear operator.  
\end{claim}
The Fourier transform is a linear operator, by the linearity
of the integral: $\widehat{(f+g)} = \hat f + \hat g$, and it is a bounded operator 
due to the elementary estimate in Lemma \ref{lem:first bound for the FT} below.

A natural question is:  where does the Fourier transform take a function $f\in L^1(\R^d)$? 
An immediate partial answer is that for any $f\in L^1(\R^d)$, we have:
\[
\hat f \in B(\R^d),
\]
where $B(\R^d):= 
\{f:\R^d \rightarrow \C  \bigm |  \, \exists M>0  \text{ such that } |f(x)| < M, \text{ for all } x \in \R^d \}$ is the space of bounded functions on $\R^d$.  Here the constant $M$ depends only on $f$. 
To see this, consider:
\begin{align}
 | \hat f(\xi) |  &:= \left |   \int_{\R^d} f(x) e^{-2\pi i \langle x, \xi \rangle} dx    \right | 
 \leq 
   \int_{\R^d}   \left |  f(x) e^{-2\pi i \langle x, \xi \rangle}    \right | dx  \\  \label{boundedness of FT}
   &=  \int_{\R^d}   \left |  f(x)   \right | dx := \| f \|_{L^1(\R^d)},   
\end{align}
where we used Theorem \ref{triangle inequality for integrals}, the triangle inequality for integrals, together with the fact that
$ \left | e^{-2\pi i \langle x, \xi \rangle}    \right | =1$.   So we've just proved the following fact.

\begin{lem} \label{lem:first bound for the FT}
Given $f \in L^1(\R^d)$, its Fourier transform is uniformly bounded, with the following bound:
\begin{equation}\label{first bound for the FT}
 | \hat f(\xi) |  \leq \| f \|_{L^1(\R^d)},
\end{equation}
for all $\xi \in \R^d$.
\hfill $\square$
\end{lem}

\begin{example}
\rm{
Let's bound the Fourier transform of an indicator function $1_S$, for any bounded measurable set $S\subset \R^d$:
\[
| \hat 1_S(\xi) |:= \left |     \int_{S} e^{-2\pi i \langle x, \xi \rangle} dx \right | 
\leq   \int_{S}   \left |     e^{-2\pi i \langle x, \xi \rangle} \right |  dx = \int_{S}   dx =  {\rm measure} (S).
\]
In particular, for any polytope $\P\subset \R^d$, 
\[
|\hat 1_\P(\xi) | \leq \vol \P, \text{ for all } \xi \in \R^d.
\]
We already know that  $\hat 1_\P(0) = \vol \P$, so it's natural to ask whether the maximum allowed value of $\vol \P$ can also be achieved by a nonzero $\xi \in \R^d$;  or perhaps it may be the case that we always have the strict inequality
$|\hat 1_\P(\xi) | < \vol \P, \text{ for all nonzero } \xi \in \R^d$? 
(See Exercise \ref{strictly less than the FT at zero}).
}
\hfill $\square$
\end{example}

But a lot more is true for absolutely integrable functions.   
\medskip
\begin{lem}   \label{uniform continuity}
 If $f\in L^1(\R^d)$, then $\hat f$ is uniformly continuous on $\R^d$.
\begin{proof}
We fix any $\xi \in \R^d$, and $h \in \R^d$, and we compute:
\begin{align*}
\hat f(\xi + h) - \hat f(\xi) &:= 
\int_{\R^d} f(x) \Big(  e^{-2\pi i \langle x, \xi + h \rangle}  -  e^{-2\pi i \langle x, \xi  \rangle} \Big) dx  \\
&= \int_{\R^d} f(x)e^{-2\pi i \langle x, \xi  \rangle} \Big(    e^{-2\pi i \langle x, h \rangle} -1  \Big) dx,
\end{align*}
so by the triangle inequality for integrals, we have
\begin{equation} \label{pre-uniform conv}
|  \hat f(\xi + h) - \hat f(\xi)|       \leq       \int_{\R^d} |f(x)| |   e^{-2\pi i \langle x, h \rangle} -1 | dx.
\end{equation}
Letting $g_h(x):= f(x) \Big(    e^{-2\pi i \langle x, h \rangle} -1  \Big)$, we see that 
\[
  |g_h(x)| \leq 2 |f(x)|,       \text{ and }         \lim_{h \rightarrow 0} |g_h(x)| =0, 
\]
using $ |e^{-2\pi i \langle x, h \rangle} -1| \leq 2$.  We may now use the dominated convergence theorem, because the functions
$g_h$ are dominated by the absolutely integrable function $2 f$.  So we get:
\[
   \lim_{h\rightarrow 0} \int_{\R^d} |f(x)| |   e^{-2\pi i \langle x, h \rangle} -1 | dx = \int_{\R^d} \lim_{h\rightarrow 0}  |f(x)| |   e^{-2\pi i \langle x, h \rangle} -1 | dx
   = 0.
\]
Because the latter limit is independent of $\xi$, \eqref{pre-uniform conv} tells us that $|  \hat f(\xi + h) - f(\xi)|  \rightarrow 0$, as $h\rightarrow 0$, uniformly  
 in $\xi \in \R^d$.
\end{proof}
\end{lem}

It turns out that sometimes we need to measure distance between functions in a manner different than just pointwise convergence.  We therefore introduce convergence in the $L^2$~norm.   We say that a sequence of functions $f_n:\R^d \rightarrow \C$ converges to a function $f$ 
{\bf in the $L^2$ norm} if
\begin{equation} \label{convergence in the L^2 norm}
\index{convergence in the $L^2$ norm}
\int_{\R^d}  \left|  f_n(x) - f(x) \right |^2 dx \rightarrow 0,  \text{  as  }  n \rightarrow \infty,
\end{equation}
for which we also use the notation  $\lim_{n\rightarrow \infty} \| f_n - f\|_2 =0$.
It is also very useful to define the $L^p(\R^d)$ spaces, for each $1\leq p < \infty$:
\begin{equation}
L^p(\R^d):= \{  f:\R^d \rightarrow \C \bigm |     \int_{\R^d} |f(x)|^p dx < \infty \},
\end{equation}
 which naturally extend the $L^1$ and $L^2$ spaces.    In fact  for functions $f \in L^p(\R^d)$, the function $\| f \|_{L^p(\R^d)} :=   \left( \int_{\R^d} |f(x)|^p dx \right)^{\frac{1}{p}}$ is a norm;  it's also a fact that 
 for $p\not=2$, this norm does not arise from an inner product.  But of course, for $p=2$ this norm does arise from an inner product, via Lemma \ref{Cauchy-Schwartz}.
 It is well-known that among all of the $L^p(\R^d)$ spaces, the only one that is a Hilbert space is $L^2(\R^d)$.  For the curious reader, the other $L^p(\R^d)$ spaces, for $p \not=2$, also possess some additional structure, namely they are Banach spaces, after identifying two functions that are equal a.e.  (see \cite{EinsiedlerWardBook} for details).    The development of $L^p$ spaces is very important for Fourier analysis; for the sake of simplicity of exposition, here we will mostly work with $p=1$ and $p=2$.
 
 Similarly to \eqref{convergence in the L^2 norm}, we define {\bf convergence in the $L^p$ norm}, 
 for $1\leq p < \infty$ by
 \begin{equation} \label{convergence in the L^p norm}
\int_{\R^d} \left|  f_n(x) - f(x) \right |^p dx \rightarrow 0,  \text{  as  }  n \rightarrow \infty,
\end{equation}
for which we also use the notation 
\[
\lim_{n\rightarrow \infty} \| f_n - f \|_p = 0.
\]
For a review of some of these various forms of convergence, see the Appendix - Chapter \ref{Appendix A}.

\bigskip
\section{The Riemann--Lebesgue lemma}\index{Riemann-Lebesgue lemma}

The celebrated {\bf Riemann--Lebesgue lemma} 
gives us the basic decay property of the Fourier transform $\hat f(\xi)$ as $|\xi| \rightarrow \infty$.
To prove it, we will use the fact that we can approximate any function $f\in L^1(\R^d)$ with arbitrary precision by using
`step functions' in $\R^d$.  More precisely, let a {\bf box in $\R^d$} \index{box}  
be defined by $\P:= [a_1, b_1]\times \cdots \times [a_d, b_d]$, and consider 
the indicator function $1_\P$ of this box.   If we consider the set of all finite sums, taken over all such indicator functions (varying over all boxes), with arbitrary real coefficients, then this set turns out to be dense in $L^1(\R^d)$, in the $L^1$ norm. We record this fact as a lemma.
 \begin{lem}  \label{box functions dense in L^1}
If $f \in L^1(\R^d)$, then there is a finite sum of indicator functions of boxes that approaches $f$, in the $L^1$ norm.
\hfill $\square$
\end{lem}


 \index{Riemann-Lebesgue lemma}
 \begin{lem}[Riemann-Lebesgue lemma]      \label{Riemann--Lebesgue lemma}
If $f \in L^1(\R^d)$, then:
 \[
 \lim_{|\xi| \rightarrow \infty}  \hat f(\xi) = 0.
 \]
 \begin{proof}
 We first show the result in the case that $f$ is the indicator function of a box.  We already know, via
  Exercise \ref{transform.of.interval.a.to.b}, that if
  $\P:= [a_1, b_1]\times \cdots \times [a_d, b_d]$, then 
  \begin{equation} \label{FT of boxes...}
  \hat 1_\P(\xi) = \prod_{k=1}^d  
\frac{ e^{-2\pi i \xi_k a_k}   -  e^{-2\pi i \xi_k b_k}         }{2\pi i \xi_k}.
  \end{equation}
As $|\xi| \rightarrow \infty$ through a sequence of $\xi$'s with nonvanishing coordinates, we see that while the numerator of \eqref{FT of boxes...} stays bounded, the denominator satisfies
 $\prod_{k=1}^d |\xi_k| \rightarrow \infty$.  Hence we've proved
the lemma for indicator functions of boxes.  Since $f \in L^1(\R^d)$, we know by Lemma~\ref{box functions dense in L^1} that 
there exists a sequence of functions $g_n \in L^1(\R^d)$ such that  $\|f -g_n\|_1 \rightarrow 0$, as $n\rightarrow \infty$.
Also, by \eqref{FT of boxes...} we know that this sequence already satisfies $\lim_{|\xi| \rightarrow \infty} \hat g_n(\xi) =0$.
  Using the elementary inequality \eqref{first bound for the FT}, 
we get: 
\[
 \big| \hat f(\xi)    -   \hat g_n(\xi) \big|   
= \big| \widehat{(f - g_n)}(\xi) \big|  
\leq   \|f -g_n\|_1  \rightarrow 0,
\]
as $n\rightarrow \infty$.  Therefore $\lim_{|\xi| \rightarrow \infty} \hat f(\xi) = 0$.
 \end{proof}
 \end{lem}

With all of the above properties, it is now natural to consider the space of all uniformly continuous functions on $\R^d$ that go to $0$ at infinity:
\begin{equation}
C_0(\R^d):= \{ f: \R^d \rightarrow \C \bigm |  f \text{ is uniformly continuous on } \R^d, \text{ and } \lim_{|x| \rightarrow \infty} |f| = 0 \}.
\end{equation}
So although the Fourier transform does not map the space $L^1(\R^d)$ into itself, all of the above results may be summarized as follows. 
\begin{lem}
If $f \in L^1(\R^d)$, then $\hat f \in C_0(\R^d)$.
\end{lem}
\begin{proof}
The boundedness of $\hat f$ was given by the inequality 
$ | \hat f(\xi) | \leq \|f\|_1$ \eqref{first bound for the FT}, the uniform continuity by 
Lemma \ref{uniform continuity}, 
and the decay to zero at infinity by Lemma \ref{Riemann--Lebesgue lemma}. 
\index{Riemann-Lebesgue lemma} 
 \end{proof}

\medskip
 \section{The inverse Fourier transform}
 
 To invert the Fourier transform, we already mentioned briefly, 
  in Theorem \ref{thm:Inverse Fourier transform} and in
  example \ref{Integral.of.sinc}, an intuitive description of this process.  
 Now we state things more formally. 
 
 \begin{thm}[The inverse Fourier transform]
   \label{thm:Inverse Fourier transform, second showing}
If  $f \in L^1(\R^d)$ and $\hat f \in L^1(\R^d)$, then 
\begin{equation}\label{Fourier inversion formula, formal}
f(x) = \int_{\R^d}  \hat f(\xi)  e^{2\pi i \langle \xi, x \rangle} d\xi,
\end{equation}
for all $x \in \R^d$.
\hfill $\square$
\end{thm}
The reader is invited to see \cite{EinsiedlerWardBook} for a proof. 
Almost all proofs proceed by introducing a Gaussian approximate identity inside the integrand of the inversion formula, then recognizing
 the integrand as a convolution with an approximate identity, and finally removing the approximate identity by invoking a limit, such as \eqref{basic smoothing} below, at each point of continuity of $f$.

The inverse Fourier transform is sometimes called `Fourier inversion'.  Let's see an interesting application, showing in particular that the Fourier transform of a polytope is \emph{not} absolutely integrable.
 
\begin{cor} \label{the FT of an indicator function is not in L^1}
 Let $C\subset \R^d$ be a compact set.  Then $\hat 1_C \notin L^1(\R^d)$.
\end{cor}
\begin{proof}
Suppose to the contrary that $\hat 1_C \in L^1(\R^d)$.  We may apply Fourier inversion, namely formula \eqref{Fourier inversion formula, formal}, because we also (trivially) have  
$1_C \in L^1(\R^d)$:
\begin{equation}\label{Fourier inversion for the indicator function, 1}
\F\left( \hat 1_C \right)(x) = 1_C(-x).
\end{equation}
But by Lemma \ref{uniform continuity}, we also know that $\F\left( \hat 1_C \right)$ is a continuous function, giving us the contradiction that $1_C$ is a continuous function. 
\end{proof}
If $C$ is a convex set, for example, then we see that $\hat 1_C$ is not absolutely integrable, and we did not have to make any messy computations to see it. 
Of course, one of the most basic consequences of Fourier inversion is the uniqueness of transforms, as follows.

\begin{cor} \label{equality of FT's implies equality of functions}
Suppose that $f, g \in L^1(\R^d)$, and that $\hat f(\xi) = \hat g(\xi)$ for all $\xi \in \R^d$.
Then $f=g$ almost everywhere. 
\end{cor}
\begin{proof}
Letting $h:= f-g$, we clearly have $h \in L^1(\R^d)$.   Let's compute:  $\hat h(\xi) = (\hat f- \hat g)(\xi) = 0$ for all $x\in \R^d$. 
In particular $\hat h \in L^1(\R^d)$, so that we may apply Fourier inversion:
\[
h(x) = \int_{\R^d} \hat h(\xi) e^{2\pi i \langle x, \xi \rangle} dx =
 \int_{\R^d}0 \, dx  = 0,
 \]
almost everywhere.
\end{proof}
Now we can revisit our intuitive Lemma \ref{complete invariance of the FT} and give a rigorous proof of a more general statement 
(see also note
 \ref{identity thm for two compact sets whose FT's agree on a convergent sequence}
in Chapter \ref{Chapter.Examples}).

\begin{thm} \label{the FT of a convex set determines the set}
Let   $\P\subset \R^d$ be a compact set. Then
$\hat 1_\P(\xi)$ uniquely determines $\P$.  Precisely, given any two  $d$-dimensional compact
sets $\P, Q\subset \R^d$, 
we have
\begin{equation}\label{difficult identity}
\hat 1_\P(\xi) = \hat 1_{Q}(\xi) \text{ for all } \xi \in \R^d   \   \iff  \  \P = Q.
\end{equation}
In particular, for any polytope $\P$, its Fourier transform $\hat 1_\P$ uniquely determines the polytope. 
\end{thm} 
\begin{proof} 
Suppose that $\hat 1_\P(\xi) = \hat 1_{Q}(\xi) \text{ for all } \xi \in \R^d$.  
We apply Corollary \ref{equality of FT's implies equality of functions} to the $L^1$ functions
$1_\P$ and $1_Q$, to conclude that $1_\P(x) = 1_Q(x)$ for almost all $x \in \R^d$.  In other words,
$1_\P-1_Q=0$ almost everywhere.   

The latter statement implies that $1_\P - 1_Q$ vanishes at each of its points of continuity.  
But $1_\P$ is continuous on the whole interior of $\P$, and similarly for $1_Q$ (they are both identically $1$ there).  
Therefore  
$1_\P(x) = 1_Q(x)=1$ for all $x$ in the interior of $\P$ and for all $x$ in the interior of $Q$. 
Also, $1_\P(x)= 1_Q(x) = 0$  for each $x \notin \P$, and for each $x \notin Q$.  
Therefore $1_\P = 1_Q$, and hence $\P=Q$.
\end{proof}

\begin{example}\label{assuming only a lattice of equalities does not imply P=Q}
{\rm
What would happen if we assume less, and replace $\R^d$ by a lattice, say $\Z^d$? 
Is it possible for the following phenomenon to occur:
\begin{equation}\label{distinct polytopes, yet their FT's agree on a lattice}
\hat 1_\P(\xi) = \hat 1_{Q}(\xi) \text{ for all } \xi \in \Z^d \notimplies \P=Q?
\end{equation}
Indeed this scenario can happen, but we need to learn about extremal bodies first
 (see  Section \ref{extremal bodies}).
\hfill $\square$
}
\end{example}

 \medskip
 \section{The torus $\R^d/\Z^d$}

Suppose a function $f:\R \rightarrow \C$ is {\bf periodic on the real line}, with period $1$:  $f(x + 1) = f(x)$, for all $x\in \R$.  Then we may think
of $f$ as `living' on the unit circle, via the map $x\rightarrow e^{2\pi i x}$ which wraps the real line onto the unit  circle.  In this setting, we may also think of the circle as the quotient group $\R/\Z$.  As we promised, group theory will not be assumed of the reader, but it will be developed a little bit in the
concrete context of lattices - see Section \ref{alternate def of a lattice - subgroups} as well. 

We may also traverse these ideas in the other direction:  commencing with any function $g$ whose domain is just $[0, 1)$, we can always extend $g$ by periodicity to the whole real line by
defining  $G(x):= \{ x\}$, the fractional part of $x$, for all $x \in \R$.   Then $G(x) = g(x)$ for all $x \in \mathbb T$, $G$ is periodic on $\R$, 
and therefore we may think of
$g$ as living on the circle $\mathbb T$.

More generally, 
 we may think of  a {\bf periodic function $f:\R^d \rightarrow \C$} as living on the cube $\square := [0,1]^d$, 
if we insist that $f$ is periodic in the following sense:
\[
f(x) = f(x + e_k), \text{ for all }  x\in \square, \text{ and all } 1\leq k \leq d.
\]
In this case, the $1$-dimensional circle is replaced by the $d$-dimensional torus 
\[
\torus:= \R^d/\Z^d,
\]
which we may also think of as the unit cube $[0, 1]^d$, but with opposite facets `glued together'. 
Here we define another infinite-dimensional vector space, namely:
\begin{equation}
 L^2(\torus):= \{  f:\torus \rightarrow \C \bigm |  \int_{[0, 1]^d} |f(x)|^2 dx < \infty \}.
\end{equation}
We notice that the domains of the integrals in  $L^2(\torus)$ are cubes, and hence always compact.  So we may therefore expect nicer phenomena to occur in this space.

We also have the space of absolutely integrable functions on the torus:
\begin{equation}
 L^1(\torus):= \{  f:\torus \rightarrow \C \bigm |  \int_{[0, 1]^d} |f(x)| dx < \infty \},
\end{equation}
which plays a simpler role than the analogous $L^1(\R^d)$ space we had before.  And finally 
we also define the useful space of $k$-differentiable functions on the torus:
\begin{equation} \label{k derivatives on the torus}
 C^k(\torus):= \{  f:\torus \rightarrow \C \bigm |  f \text{ has $k$ continuous derivatives} \}.
\end{equation}
As a special case, we'll simply denote by $C(\torus)$ the space of all continuous functions on the torus. 
We emphasize that by definition, all of the latter function spaces,  $C^k(\torus), L^1(\torus), L^2(\torus)$, consist
of \emph{periodic functions} on the cube $[0, 1]^d$.

Similarly to the inner product on $L^2(\R^d)$, we also have in this new context a natural inner product for the space of square-integrable functions 
 $f\in L^2(\torus)$, defined by:
 \begin{equation}
 \langle f, g \rangle := \int_{[0,1]^d} f(x) \overline{g(x)} dx,
 \end{equation}
making $L^2(\torus)$ a Hilbert space.
 For each $n\in \Z^d$, we define $e_n: \R^d \rightarrow \C$ by:
 \begin{equation} \label{basis for Hilbert space}
 e_n(x):= e^{2\pi i \langle n, x\rangle}.
 \end{equation}
 This countable collection of exponentials turns out to form a complete orthonormal basis for $L^2(\torus)$. 
The orthogonality is the first step, which we prove next.   For the proof that the exponentials span $L^2(\torus)$ and are complete, 
we refer the reader to \cite{EinsiedlerWardBook}.

\bigskip
 
 \begin{thm}   [{\bf Orthogonality relations for the exponentials $e_n(x)$ on the torus}]
 \index{orthogonality relations for the exponentials $e_n(x)$}  
 \label{Orthogonality relations for the exponentials $e_n(x)$}
 
\begin{equation}  
 \int_{[0,1]^d}    e_n(x)    \overline{e_m(x)} dx = 
 \begin{cases}  
1    &      \mbox{if } n=m  \\ 
0  &        \mbox{if  not}.
\end{cases}
\end{equation}
\end{thm}
\begin{proof}
Because of the geometry of the cube, we can proceed in this case by separating the variables.   If $n \not= m$, then there is at least one index $k$ for which $n_k \not= m_k$.   We compute:
\begin{align*}
 \int_{[0,1]^d}    e_n(x)    \overline{e_m(x)} dx 
&=       \int_{[0,1]^d}    e^{2\pi i \langle n-m, x \rangle} dx \\
&=     \int_0^1  e^{2\pi i (n_k -m_k)x_k}  dx     \int_{[0,1]^{d-1}}    \prod_{j\not=k}  e^{2\pi i (n_j -m_j)x_j}  dx       \\
&=     \int_0^1  e^{2\pi i (n_k -m_k)x_k}  dx      \int_{[0,1]^{d-1}}    \prod_{j\not=k}  e^{2\pi i (n_j -m_j)x_j}  dx      \\
&=    \left( \frac{e^{2\pi i (n_k-m_k)}-1}{2\pi i (n_k-m_k)}    \right)    \int_{[0,1]^{d-1}}    \prod_{j\not=k}  e^{2\pi i (n_j -m_j)x_j}  dx  =0, 
\end{align*}
because $n_k -m_k$ is a nonzero integer.
\end{proof}
 
 Because $L^2(\torus)$ is also an inner product space, it still enjoys the Cauchy-Schwartz inequality.
 Intuitively, the space $L^2(\torus)$ should be a cozier little space than $L^1(\torus)$.   This intuition can be made more
  rigorous by the following Lemma, despite the fact that $L^2(\R^d) \not\subset L^1(\R^d)$.   More generally, given any compact and convex set
  $\P \subset \R^d$, and any $p\geq 1$, we  define 
  \begin{equation}
  L^p(\P):=      \{  f:\P \rightarrow \C \bigm |     \int_{\P} |f(x)|^p dx < \infty \}.
  \end{equation}

 \begin{lem}   \label{proper containment of L^2 in L^1 for torus}
We have the following proper containments:
 \begin{enumerate}[(a)]
 \item        \label{proper containment of spaces over the torus}
 $L^2(\torus) \subset L^1(\torus)$.
 \item       \label{lem:proper containment, second part} 
 In general, given any compact and convex set $\P\subset \R^d$, 
  $L^2(\P) \subseteq L^1(\P)$.
 \end{enumerate}
  \end{lem}
 \begin{proof}
 Given $f \in  L^2(\torus)$, we must show that $f \in  L^1(\torus)$.  Using the Cauchy-Schwartz inequality for 
 $ L^2(\torus)$, applied to $f$ and the constant function $h(x) \equiv 1$ on $\torus$, we have:
 \begin{align*}
\int_\torus |f(x)| dx &= \int_\torus |f(x) h(x) |dx \\
& \leq  
      \left(      \int_\torus      |f(x)|^2 dx  \right)^{\frac{1}{2}}  
      \left(      \int_\torus      |h(x)|^2 dx \right)^{\frac{1}{2}}  \\
&=  \left(      \int_\torus      |f(x)|^2 dx  \right)^{\frac{1}{2}},
\end{align*}
 so we see that $f$ is absolutely integrable over the torus $\torus$.   To show that the containment in part  \ref{proper containment of spaces over the torus} is proper, for $d=1$, we can consider the following function on
$[0, 1]$:
 \[
 f(x):= 
 \begin{cases}
 \frac{1}{\sqrt x} & \text{ if } x  \in (0, 1], \\
 0 &  \text{ if } x=0.
 \end{cases}
 \]
 So $\int_0^1 f(x) dx =  2 x^{\frac{1}{2}} \Big|_0^1=2$, but $\int_0^1 |f(x)|^2  dx =
 \int_0^1 \frac{1}{x} dx  = \infty$.  Hence $f \in L^1(\torus)$, but $f \notin \subset L^2(\torus)$.
 
 The general case of part \ref{proper containment of spaces over the torus} for arbitrary dimension follows easily from this example.    For part \ref{lem:proper containment, second part}, once we know that $L^2(\P)$ is a Hilbert space (\cite{RudinGreenBook}), it follows that it has a Cauchy-Schwartz inequality, so the same proof of part \ref{proper containment of spaces over the torus} works.
 \end{proof}


 \bigskip
 \subsection{Fourier series: fast and loose} 
 
 Let's see how we can expand (certain) functions in a Fourier series, as well as find a formula for their series coefficients,  in a foot-loose and carefree way - i.e. abandoning all rigor for the moment.
 
 Given that the sequence of exponential functions $\{e_n(x) \}_{n \in \Z^d}$ forms a basis
for the infinite dimensional vector space $V:=L^2(\torus)$, 
we know from Linear Algebra that any function $f \in V$ may be written in terms of this basis:
 \begin{equation} 
 f(x) = \sum_{n \in \Z^d}  a_n e_n(x).
 \end{equation}
 
How do we compute the Fourier coefficients $a_n$?  Let's go through the intuitive process here, ignoring convergence issues.  Well, again by Linear Algebra, we take the inner product of both sides with a fixed basis element $e_k(x)$:
 \begin{align*}
 \langle f(x), e_k(x) \rangle &= \langle \sum_{n \in \Z^d} a_n e_n(x), e_k(x) \rangle \\
 &=  \sum_{n \in \Z^d}  a_n \langle e_n(x), e_k(x) \rangle \\
  &=  \sum_{n \in \Z^d}  a_n \, \delta(n,k) \\
  &= a_k
 \end{align*}
 where we've used   the orthogonality relations, Theorem \ref{Orthogonality relations for the exponentials $e_n(x)$}  above,   in the third equality.  We also used the standard notation $\delta(n,k) := 0$ if $n\not=k$, and $\delta(n,k):=1$ if $n=k$.
 Therefore, it must be the case that 
 \begin{align*}
 a_k &=  \langle f(x), e_k(x) \rangle \\
 &:= \int_{[0,1]^d}   f(x) \overline{   e^{2\pi i \langle k, x \rangle} } dx \\
 &=  \int_{[0,1]^d}   f(x) e^{-2\pi i \langle k, x \rangle} dx,
 \end{align*}
 also called the
 {\bf Fourier coefficients} of $f$.

 
 \bigskip
 \subsection{Fourier series: slow and rigorous}
 
Let's record now the rigorous statements of the intuitive arguments that we constructed in the previous section.  We may think of a periodic function on $\R^d$ as a function belonging to $L^2(\torus)$.

 \begin{thm}[{\bf Fourier series for functions on $\torus$}]
  \label{Fourier series for periodic functions}
 \index{Fourier series for periodic functions} 
The set of exponentials 
\[
\{ e_n(x) \bigm | n \in \Z^d  \}
\]
 form a {\bf complete orthonormal basis} for $L^2(\torus)$.  
 Moreover, we have the following:
 \begin{enumerate}[(a)]
 \item Every function $g \in L^2(\torus)$ has a {\bf Fourier series} 
 \begin{equation} \label{The Fourier series}
 g(x) = \sum_{n\in \Z^d} c_n e^{2\pi i \langle n, x \rangle},
 \end{equation}
where the convergence in \eqref{The Fourier series} takes place in the $L^2$ norm on the torus $\torus$.
\item The {\bf Fourier coefficients} $c_n$ may be computed via the formula:
\begin{equation}  \label{Fourier coefficients}
 c_n = \int_{[0,1]^d}  g(t) e^{-2\pi i \langle n, t\rangle} dt,
\end{equation}
for all $n\in \Z^d$. 

\item \rm{(\bf The Parseval identity})
The function $g\in L^2(\torus)$ in \eqref{The Fourier series}  satisfies
\begin{equation} \label{True Parseval identity}
\int_{[0, 1]^d} | g(x)|^2 dx = \sum_{n \in \Z^d} |c_n|^2.
\end{equation}
\end{enumerate}
\hfill $\square$
\end{thm}
(For a proof, see \cite{EinsiedlerWardBook},  p. 96)
At the risk of overstating the obvious, we note that the equality in \eqref{True Parseval identity} is simply equality between real numbers. 
 We also note that the Fourier coefficients above are integrals over the unit cube
 $[0,1]^d$, and may also be thought of as $c_n = \langle g, e_n \rangle$, the projection of $g$ onto each basis element.  To summarize, we've encountered the following types of transforms so far:
\begin{equation} \label{both integrals}
\int_{[0,1]^d}  g(t) e^{-2\pi i \langle n, t\rangle} dt, \text{ and } \int_{\R^d}  g(t) e^{-2\pi i \langle n, t\rangle} dt.
\end{equation}
 To disambiguate, the first integral in \eqref{both integrals} arises from periodic functions on $\R^d$, and it appears as a Fourier coefficient in 
 Theorem \ref{Fourier series for periodic functions}.
 The second integral is our old friend the Fourier transform.  How are the two integrals related to each other? This is exactly the magic of the Poisson summation formula, Theorem \ref{Poisson.Summation}.

In the pretty proof of Poisson summation, we begin with a Fourier series of a periodized version of $f$, 
and end up showing that its Fourier coefficients, by a small miracle of nature, turn out to also be Fourier transforms of $f$. 
 
 A natural question is:  
\begin{question}
Which functions have a pointwise convergent Fourier series?
\end{question}
But this question turns out to be rather difficult, and many lifetimes have been devoted to related questions.
It is a fact of life that the Fourier series of an arbitrary continuous function on $\R^d$ may fail to converge uniformly, or even pointwise.  
However, there is some good news.    As it turns out, if we impose some smoothness conditions on $f$, then $f$ does have a Fourier series which converges pointwise, as we'll see next.

\bigskip
\section{Piecewise smooth functions have convergent Fourier series}

In this section we'll restrict attention to the real line.   We'd like to rigorously define the intuitive idea of a function that is almost continuous, in the sense of being continuous on an interval except for 
 finitely many finite jump discontinuities.

Precisely, given real numbers $a, b$, we define a function 
$f:[a, b] \rightarrow \C$ to be {\bf piecewise continuous} on $[a, b]$ if the following two conditions are met:
\begin{enumerate}
\item $f$ is continuous on $(a, b)$, except possibly on a finite set of points
 \[
 S:=\{x_1, \dots, x_N\}\subset [a, b]. 
 \]
\item  The left-hand and right-hand limits of $f$ exist at each of the points $x_k \in S$:
\begin{equation*}
f(x_k^+) := \lim_{\varepsilon \rightarrow 0 \atop \varepsilon >0} f(x_k+\varepsilon) \ \text{   exists, and }
f(x_k^-)  := \lim_{\varepsilon \rightarrow 0 \atop \varepsilon >0} f(x_k-\varepsilon) \ \text{ exists}.
\end{equation*}
\end{enumerate}
Furthermore, we define a function 
$f:[a, b] \rightarrow \C$ to be {\bf piecewise smooth} on $[a, b]$ if both $f$ and its derivative $f'$ are piecewise continuous on $[a, b]$.    We'll also say that a function $f:\R\rightarrow \C$ is piecewise continuous/smooth on $\R$ if it is piecewise continuous/smooth on every finite interval.
We have the following refined version of Theorem \ref{Fourier series for periodic functions}, on the real line.
 \begin{thm} \label{theorem:Fourier series convergence to the mean}
 Let $f:\R \rightarrow \C$ be a periodic function, with domain $[0 , 1]$, and piecewise smooth  on $\R$.  
  Then, for each $t \in \R$, we have
 \begin{equation}  \label{convergence of Fourier series to the mean}
   \lim_{N\rightarrow \infty}  \sum_{n= -N}^N  c_n  e^{2\pi i n t} =
   \frac{  f(t^+)  + f(t^-)  }{2},
 \end{equation}
 where $c_n:= \int_0^1 f(x) e^{-2\pi i x n} dx$ are the Fourier coefficients of $f$.
 
 \rm{(For a proof of Theorem \ref{theorem:Fourier series convergence to the mean} see \cite{Travaglini})}.
 \hfill $\square$
\end{thm}

We will come back to these {\bf partial Fourier sums}, \index{partial Fourier sums}
occurring in Theorem \ref{theorem:Fourier series convergence to the mean}, 
and defined by
\begin{equation} \label{partial sums}
S_N f(t):= \sum_{n= -N}^N  c_n e^{2\pi i n t}.
\end{equation}

There is also a natural and easy extension of Parseval's identity \eqref{True Parseval identity}.  
Given any two functions $f, g\in L^2(\torus)$, we've seen in \eqref{The Fourier series} 
that 
\[
f(x) = \sum_{n\in \Z^d} a_n \, e^{2\pi i \langle n, x \rangle}, \text{ and }
g(x) = \sum_{n\in \Z^d} b_n \, e^{2\pi i \langle n, x \rangle},
\]
 both converging in the $L^2(\torus)$ norm.
\begin{thm}
If $f, g, \in  L^2(\torus)$, then with the notation above we have
\[
\int_{\torus} f(x) \overline{g(x)} dx =  \sum_{n\in \Z^d} a_n \overline{b_n}.
\]
\hfill $\square$
\end{thm}



\subsection{The first periodic Bernoulli polynomial}

 To see a concrete instance of Theorem \ref{Fourier series for periodic functions}, we study the 
  function $P_1(x)$, which we've briefly encountered before, as the first periodic Bernoulli polynomial.  This function turns out to be so important that it deserves its own section here.   We recall its definition:
 \begin{equation}  \label{def of P_1 again}
 P_1(x):=  
 \begin{cases}  
\{ x \}   - \frac{1}{2}                 &\mbox{if }  x \notin \Z,   \\   
0                          & \mbox{if } x \in \Z.
\end{cases}
 \end{equation}
 
 It's easy to see that
$P_1 \in L^1(\mathbb T)$, so it has a Fourier series, by Theorem \ref{Fourier series for periodic functions}, part (a):
\begin{equation}
P_1(x) = \sum_{n \in \Z} c_n e^{2\pi i n x},
\end{equation}
and the equality here means equality in the $L^2(\mathbb T)$ norm.
\begin{figure}[htb]
\begin{center}
\includegraphics[totalheight=1.8in]{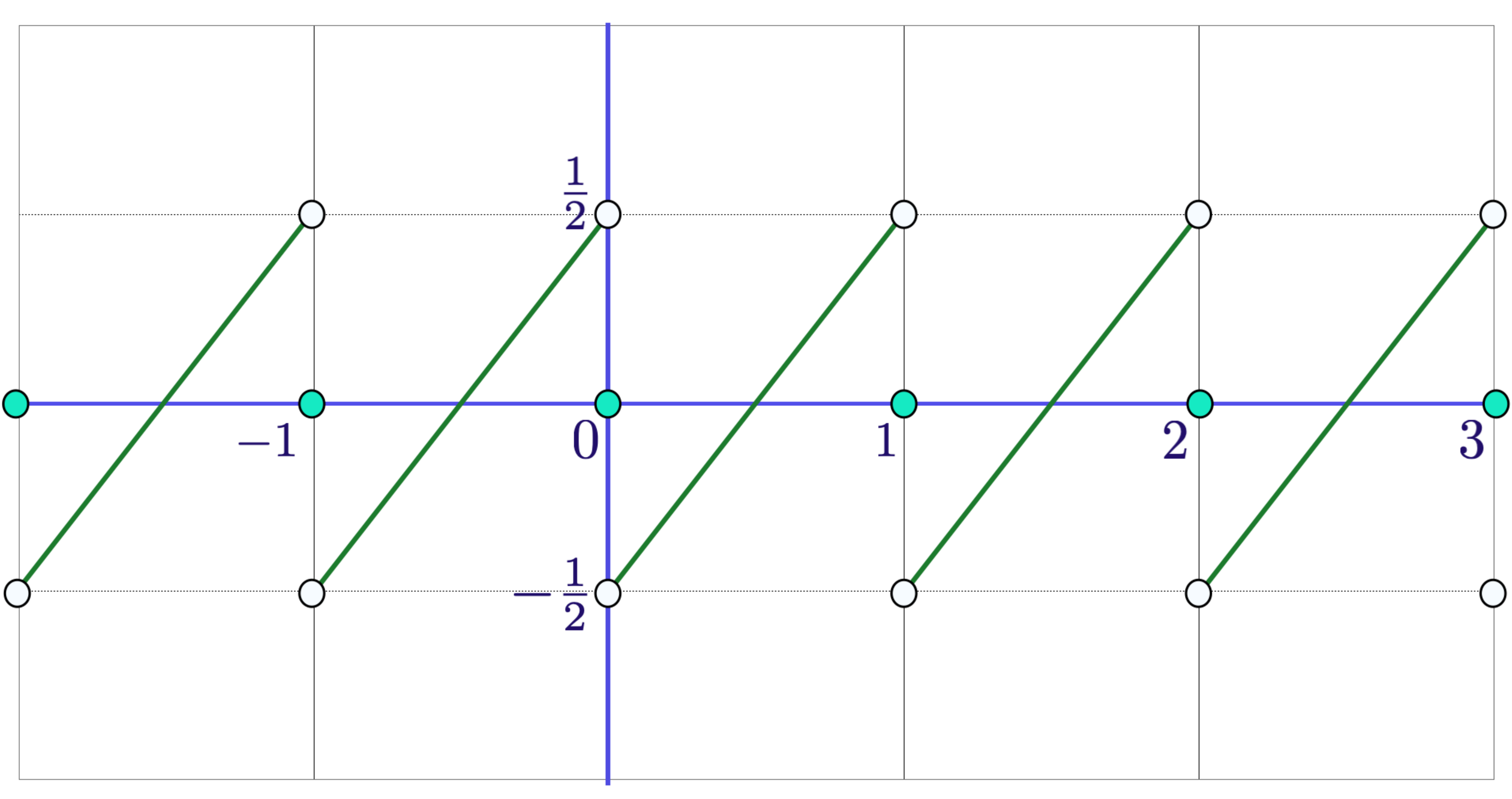}
\end{center}
\caption{The first periodic Bernoulli polynomial $P_1(x)$ }
\end{figure}
Let's compute the Fourier coefficients of $P_1$, according to 
Theorem \ref{Fourier series for periodic functions}, part (b).  We will use integration by parts:
\begin{align*}
c_n&= \int_0^1 \left( \{  x  \} -  \tfrac{1}{2} \right) e^{-2\pi i n x} dx 
=  \int_0^1 x    e^{-2\pi i n x} dx   -    \tfrac{1}{2} \int_0^1  e^{-2\pi i n x} dx \\
&=  x \frac{ e^{-2\pi i n x}  }{-2\pi i n }  \Big|_0^1 - \int_0^1  \frac{ e^{-2\pi i n x}  }{-2\pi i n } dx 
= \frac{ 1 }{-2\pi i n }  - 0  =  \frac{ 1 }{-2\pi i n },
\end{align*}
when $n \not=0$.  For $n=0$, we have $c_0 = \int_0^1 (x - \tfrac{1}{2})dx = 0$.
Hence we have the Fourier series 
\begin{equation}  \label{Fourier series of P_1 in the norm}
P_1(x) =\{ x \} -\frac{1}{2}=    - \frac{ 1 }{2\pi i  }  \sum_{n \in \Z  \atop n\not=0}  \frac{ 1 }{ n }  e^{2\pi i n x},
\end{equation}
where the latter equality means convergence in the $L^2(\torus)$ norm.  But we'd like pointwise convergence of the latter series!  In fact, this follows from Theorem 
\ref{theorem:Fourier series convergence to the mean}, as follows.

\begin{cor}  \label{example of theorem on pointwise convergence}
We have the pointwise convergent Fourier series
\begin{align}
 \lim_{N\rightarrow \infty}    -\frac{1}{2\pi i}    \sum_{-N \leq n \leq N \atop n\not=0}   \frac{1}{n} e^{2\pi i n x} 
 =   \{ x \} -\frac{1}{2},
\end{align}
valid for all $x \in \R$.
\end{cor}
\begin{proof}
First, we fix any $x \notin \Z$.  Theorem \ref{theorem:Fourier series convergence to the mean}
 allows us to conclude that we have pointwise convergent sums:
\begin{align}
 \lim_{N\rightarrow \infty}    -\frac{1}{2\pi i}    \sum_{-N \leq n \leq N \atop n\not=0}   \frac{1}{n} e^{2\pi i n x} 
 &=    \frac{  P_1(x^+)  + P_1(x^-)  }{2} \\
 &=   \{ x \} -\frac{1}{2},
\end{align}
For $x \in \Z$, we can also check that the equality 
\eqref{example of theorem on pointwise convergence} holds by observing that 
\[
\sum_{-N \leq n \leq N  \atop n\not=0}   \frac{1}{n} e^{2\pi i n x} =  
\sum_{-N \leq n \leq N  \atop n\not=0}  \frac{1}{n} = 0,
\]
while $\frac{  P_1(x^+)  + P_1(x^-)  }{2} =\tfrac{1}{2}\left(- \frac{1}{2}  + \frac{1}{2}\right)  =0 $ as well, which is consistent with the definition 
\eqref{def of P_1 again}
of $P_1(x)$ at the integers.
\end{proof}

Next, we can give a classical application of the Fourier series \eqref{Fourier series of P_1 in the norm}
using Parseval's identity \eqref{True Parseval identity}:
\[
\int_0^1 |P_1(u)|^2  du = \sum_{n \in \Z} |a_n|^2.
\]
Let's simplify both sides:
\begin{align*}
 \sum_{n \in \Z} |a_n|^2 &=  \frac{1}{4\pi^2} \sum_{n \in \Z-\{0\}} \frac{1}{n^2} =  
 \frac{1}{2\pi^2} \sum_{n \geq 1} \frac{1}{n^2},
 \end{align*}
 while
 \begin{align*}
 \int_0^1 |P_1(u)|^2  du = \int_0^1 \left( \{  x  \} - \frac{1}{2} \right)^2 dx 
 &=  
   \int_0^1 \left(x - \frac{1}{2}\right)^2 dx 
  = \frac{1}{12}.
\end{align*}
Therefore  
\[
\sum_{n \geq 1} \frac{1}{n^2} = \frac{\pi^2}{6},
\]
a number-theoretic identity that goes back to Euler. 
In a  similar manner one can evaluate the Riemann zeta function at all positive even integers, 
using the cotangent function
(Exercise \ref{Riemann zeta function, and Bernoulli numbers}).

\bigskip
Another natural question arises.
\begin{question}
What sort of functions $f:\torus \rightarrow \C$ are uniquely determined by all of their Fourier coefficients?
\end{question}
To describe a partial answer, we recall the space of all continuous functions on the torus:
\begin{equation}
C(\torus):= \{ f: \torus \rightarrow \C \bigm |   f \text{ is continuous on } \torus \}.
\end{equation}
\begin{thm}
Let $f \in C(\torus)$, and suppose that $\hat f(n) = 0$ for all $n \in \Z^d$.   Then $f(x)=0$, 
for all $ x \in [0, 1]^d$.

In particular, if $f, g \in C(\torus)$ and $\hat f(n) =\hat g(n)$ for all $n \in \Z^d$, 
then $f(x) = g(x)$ for all $x \in [0, 1]^d$.

\hfill $\square$ 
\end{thm}
In other words, a continuous function on the torus is uniquely determined by its Fourier coefficients  
(see \cite{EinsiedlerWardBook} for a proof).

\bigskip
\section{As $f$ gets smoother, $\hat f$ decays faster}
\label{As f gets smoother, the FT decays faster}

There is a very basic and important relationship between the level of smoothness of $f$, 
and the  speed with which  $\hat f$ tends to $0$ as $x \rightarrow \infty$.   To capture this relation very concretely, let's compute things on the real line, 
to see how the FT interacts with the derivative.
\begin{lem}   \label{FT interacts with derivative}
Let $f \in L^1(\R)$. 
\begin{enumerate}[(a)] 
\item \label{FT of derivative}
If $f$ is piecewise smooth, and also enjoys $f' \in L^1(\R)$, then:
\[
\widehat{ f' }(\xi) = (2\pi i)  \xi \hat f(\xi).
\]
\item  \label{explicit decay rate}
More generally, let $k\geq 0$, suppose that $f$ has $k$ derivatives, $f^{(k)}$ is piecewise smooth, 
 and that 
 we also have $f^{(k+1)} \in L^1(\R)$.
Then:
\[
\widehat{ f^{(k+1)} }(\xi) = (2\pi i \xi)^{k+1} \hat f(\xi).
\]
\item  \label{derivative of the FT}
Now we suppose that $x f(x) \in L^1(\R)$.  Then:
\[
\frac{d}{d\xi}  \F(f)(\xi) =  (-2\pi i) \, \F(x f(x))(\xi). 
\]
\end{enumerate}
\end{lem}
\begin{proof}
To prove part \ref{FT of derivative}, we notice that $\lim_{x\rightarrow \infty} f(x) = f(0) + \int_0^\infty f'(x) dx$, using the hypothesis $f' \in L^1(\R)$. 
Using the hypothesis $f \in L^1(\R)$, we know that the Riemann-Lebesgue Lemma \ref{Riemann--Lebesgue lemma} implies that $\lim_{x\rightarrow \infty} f(x)=0$.
Similarly, $\lim_{x\rightarrow -\infty} f(x) =0$.  Integration by parts now gives us:
\begin{align*}
\widehat{ f' }(\xi) &= \int_{\R}  f'(x) e^{-2\pi i x \xi} dx  
= f(x)e^{-2\pi i x \xi}\Big |_{-\infty}^\infty -  \int_{\R}  f(x)  (-2\pi i \xi) e^{-2\pi i x \xi} dx   \\
&= 2\pi i \xi \int_{\R}  f(x)  e^{-2\pi i x \xi} dx   := 2\pi i \xi \hat f(\xi).
\end{align*}
Part \ref{explicit decay rate}   follows from \ref{FT of derivative} by induction on $k$.  To prove part \ref{derivative of the FT}, we have:
\begin{align*}
\F(x f(x))(\xi) &:=  \int_{\R}  xf(x) e^{-2\pi i x \xi} dx = \frac{1}{-2\pi i} \int_{\R}  \frac{d}{d\xi} f(x) e^{-2\pi i x \xi} dx\\
&=-\frac{1}{2\pi i} \frac{d}{d\xi} \int_{\R}   f(x) e^{-2\pi i x \xi} dx 
=-\frac{1}{2\pi i} \frac{d}{d\xi} \hat f(\xi).
\end{align*}
\end{proof}

It follows from Theorem \ref{FT interacts with derivative}, part  \ref{explicit decay rate}, that we have an explicit decay rate for the Fourier coefficients of a periodic function $f$,
assuming that $f$ is sufficiently smooth.  To obtain the following Corollary, we can simply use the fact that $f^{(k+1)} \in L^1(\R)$ implies that $ \widehat{ f^{(k+1)}}$ is uniformly bounded: 
$\Big|  \frac{1}{(2\pi )^{k+1}}  \widehat{ f^{(k+1)} }(\xi) \Big| < C$, for a positive constant $C$.

\begin{cor} \label{cor: f smoother implies FT of F decays faster}
If $f$  has $k$ continuous derivatives, and  we also have $f^{(k+1)} \in L^1(\R)$,
then there is a constant $C>0$ such that:
\begin{equation}
| \hat f(\xi) | < C \frac{1}{  | \xi|^{k+1}}, 
\end{equation}
for all $\xi \not=0$.
\hfill $\square$
\end{cor}
In other words, we now understand
the dictum ``as $f$ gets smoother, $\hat f$ decays faster'' in a precise quantitative manner:  
if  $f$ has $k$ derivatives, then $\hat f$ decays faster than a polynomial of degree $k$.

    
\bigskip
\section{How fast do Fourier coefficients decay?}
    
In a manner completely analogous to the previous Section \ref{As f gets smoother, the FT decays faster}, 
we can repeat the important idea of integration by parts to see how fast
Fourier coefficients decay, and here we may expect even better results, because we will integrate over the compact unit cube (equivalently over $\torus$), rather than over the non-compact space $\R^d$.   We first work things out in dimension $1$, recalling that 
  the Fourier coefficients of $f$ are defined by $c_n:= \int_0^1 f(x) e^{-2\pi i n x} dx$, for all $n \in \Z$.   
For the sake of the reader, we  recall the space of functions $C^k(\mathbb T)$ from 
\ref{k derivatives on the torus}, which have $k$ continuous derivatives.  We also recall that  
 $f \in L^1(\mathbb T)$ means $\int_0^1f(x) dx$ is finite, and that $f(x+1) = f(x)$, 
 for all $x\in [0, 1]$.  Finally, we note that the same conclusion of the Riemann-Lebesgue lemma \ref{Riemann--Lebesgue lemma} also holds for
 functions $f \in L^1(\torus)$, with exactly the same proof that we gave in Lemma \ref{Riemann--Lebesgue lemma}.

\begin{thm} \label{decay of Fourier coefficients}
Let $f\in L^1(\mathbb T)$.
\begin{enumerate}[(a)]
\item  \label{decay of Fourier coefficients for C^1}
  If $f \in C^1(\mathbb T)$, then its Fourier coefficients satisfy
\begin{equation}
\lim_{|n|\rightarrow \infty}  |n c_n| = 0.
\end{equation}
In other words, $|c_n| = o\left( \frac{1}{n} \right)$.
\item  \label{decay of Fourier coefficients for C^k}
More generally, fix an integer $k\geq 1$.  If $f \in C^k(\mathbb T)$, then its Fourier coefficients satisfy
\begin{equation}
\lim_{|n|\rightarrow \infty}  |n^k c_n| = 0.
\end{equation}
In other words, $|c_n| = o\left( \frac{1}{n^k} \right)$.
\end{enumerate}
\end{thm}
\begin{proof}    
We compute the Fourier coefficients using integration by parts.  For each $n \not=0$, we have:
\begin{align*}
c_n &:= \int_0^1 f(x) e^{-2\pi i n x} dx 
=   \left[  f(x) \frac{ e^{-2\pi i n x} }{-2\pi i n}\right] \Big|_0^1 + \frac{ 1}{2\pi i n}\int_0^1 f'(x) e^{-2\pi i n x} dx \\
&=  \frac{   f(1)-f(0)  }{  -2\pi i n  }   + \frac{ 1}{2\pi i n}\int_0^1 f'(x) e^{-2\pi i n x} dx \\
&=  \frac{ 1}{2\pi i n}\int_0^1 f'(x) e^{-2\pi i n x} dx,
\end{align*}
using the periodicity of $f$.   Because $f'$ is continuous, the Riemann-Lebesgue lemma on $L^1(\mathbb T)$ gives us
$\lim_{|n|\rightarrow \infty} \int_0^1 f'(x) e^{-2\pi i n x} dx = 0$.  So we see that 
\[
|n c_n| \rightarrow 0, \  \text{ as } |n| \rightarrow \infty,
\]
completing part \ref{decay of Fourier coefficients for C^1}.   Part \ref{decay of Fourier coefficients for C^k} follows easily by induction on $k$, repeating  the same integration by parts computation
  above, exactly $k$ times. 
\end{proof}    
    We note that the same proof works with even weaker hypotheses in part \ref{decay of Fourier coefficients for C^k}.   Namely, given  an integer $k\geq 1$, all we require is that 
$f^{(j)}$ is continuous on $\mathbb T$, for $0\leq j < k$, and $f^{(k)} \in L^1(\mathbb T)$.

Let's see a concrete application of these ideas (see Note \ref{Note:GregKuperberg}).
\begin{thm} \label{Euler-Maclaurin type identity}
Suppose that $f\in C^k(\mathbb T)$, for a fixed integer $k\geq 1$.   Then:
\begin{equation}
\int_0^1 f(x) dx = \frac{1}{N} \sum_{m=0}^{N-1} f\left(\tfrac{m}{N}\right) + o\left(\tfrac{1}{N^{k}}\right),
\end{equation}
as $N\rightarrow \infty$.
\end{thm}
\begin{proof}
Because $f$ is periodic on $\R$, we follow  ``Hecke's dictum'';   namely, we first expand $f$ into its Fourier series, which is guaranteed by 
Theorem \ref{Fourier series for periodic functions}:
\[
f(x) = \sum_{n\in \Z} c_n e^{2\pi i n x}. 
\]
Since this Fourier series converges absolutely, we may interchange the finite sum with the series:
\begin{align*} \label{first step of finite sum approximation}
 \frac{1}{N} \sum_{m=0}^{N-1} f\left(\tfrac{m}{N}\right) &=  \frac{1}{N} \sum_{m=0}^{N-1} 
 \sum_{n\in \Z} c_n e^{2\pi i n \tfrac{m}{N}}     \\
 &=  \sum_{n\in \Z} c_n 
 \Big( \tfrac{1}{N} \sum_{m=0}^{N-1}  e^{2\pi i n \tfrac{m}{N}} \Big) \\
 &=\sum_{n\in \Z} c_{Nn}, 
\end{align*}
using Exercise \ref{DivisibilityUsingExponentials} (the harmonic detector for divisibility).
Next, we recall that the constant term is $c_0= \int_0^1 f(x) dx$, and we separate out this term from the latter series:
\begin{align*}
 \frac{1}{N} \sum_{m=0}^{N-1} f\left(\tfrac{m}{N}\right) 
  =  \int_0^1 f(x) dx +  \sum_{n\in \Z \atop {n \not=0 }   } c_{Nn}, 
\end{align*}
Now we can use the (little-o) rate of decay of the Fourier coefficients, given by 
Theorem \ref{decay of Fourier coefficients}, part \ref{decay of Fourier coefficients for C^k}, to write 
$|c_{Nn}| <  \frac{C}{(Nn)^k}$ for \emph{all} constants $C>0$.  We conclude that
\[
 \sum_{n\in \Z \atop {n \not=0 }   }   |c_{Nn}|  
<   C \sum_{n\in \Z \atop {n \not=0 }}    \frac{1}{  N^{k}  |n|^{k}} 
<   2 C\, \zeta(k)  \frac{1}{  N^{k}},  
\] 
for all constants $C>0$.  So as $N\rightarrow \infty$, the error term 
$\sum_{n\in \Z \atop {n \not=0 }   } c_{Nn}$ is $o\left( \frac{1}{  N^k}      \right)$, as claimed. 
\end{proof}

It is worth mentioning that although our proof of Theorem \ref{Euler-Maclaurin type identity}
does not cover the case $k=0$,  this case
 is also true because it represents the Riemann sum approximation to the integral. 


\section{The Schwartz space}     \label{nice functions} 
\index{Schwartz space}

We saw in Section \ref{As f gets smoother, the FT decays faster}  that a  function $f:\R^d \rightarrow \C$  in the space domain,  that has $k$ derivatives, corresponds to a function $\hat f$ in the Fourier transform domain.
If we `take this idea to the limit', so to speak,   What does that last adjective mean?  
Following the ideas of Laurent Schwartz, we can make rigorous sense of the words `rapidly decreasing', as follows.

We recall that our definition of a `nice function' was any function $f:\R^d \rightarrow \C$ for which the Poisson summation formula holds.  Here we give our first family of sufficient conditions for a function $f$ to be nice. 
A {\bf Schwartz function}  $f: \R \rightarrow \C$ is defined as any infinitely smooth function ($f \in C^\infty(\R)$)
 that satisfies the following
 growth condition:  
 \begin{equation}
 |x^a  \frac{d}{dx^k} f(x)| \text{ is bounded on } \R,
 \end{equation}
  for all integers $a, k \geq 0$.
In particular, a Schwartz function decreases faster than any polynomial function, as $|x|$ tends to infinity.

 \medskip
 \begin{example}
 \rm{
 The Gaussian function $G_t(x) := e^{-t ||x||^2}$ is a Schwartz function, for each fixed
 $t>0$.     To see this, we first consider $\R^1$, where we note that the $1$-dimensional Gaussian
 is a Schwartz function, as follows.  We observe that for all positive integers $k$, 
 $ \frac{d}{dx^k} G_t(x) = H_n(x) G_t(x)$, where $H_n(x)$ is a univariate polynomial in $x$ (which also depends on the parameter $t$, but we think of $t$ as a constant).
 Since $\lim_{x\rightarrow \infty}    \frac{x^a \, H_n(x)}{ e^{t ||x||^2}} =0$, for all positive integers $a$, we see that  $G_t(x)$ is a Schwartz function.
 Now we note that the product of Schwartz functions is again a Schwartz function;   hence the $d$-dimensional Gaussian, $G_t(x) := e^{-t ||x||^2} = \prod_{k=1}^d  e^{-t x_k^2} $, a product of $1$-dimensional Gaussians, is a Schwartz function. 
 
  Some might say the Gaussian is the quintessential Schwartz function, partly because it is also an eigenfunction of the Fourier transform, as we'll see below. 
 }
 \hfill $\square$
 \end{example}

 \medskip
 \begin{example}  \label{example: the abs value exponential}
 \rm{  
We define $f(x) := e^{-2\pi t |x|}$ on the real line, for a fixed $t >0$.  
To see  that $f$ is \emph{not} a Schwarz function, we merely have to observe that $f$ is not differentiable at $x=0$.  
To be a Schwartz function, $f$ would have to be infinitely differentiable everywhere on $\R$.

\begin{figure}[htb]
\begin{center}
\includegraphics[totalheight=2in]{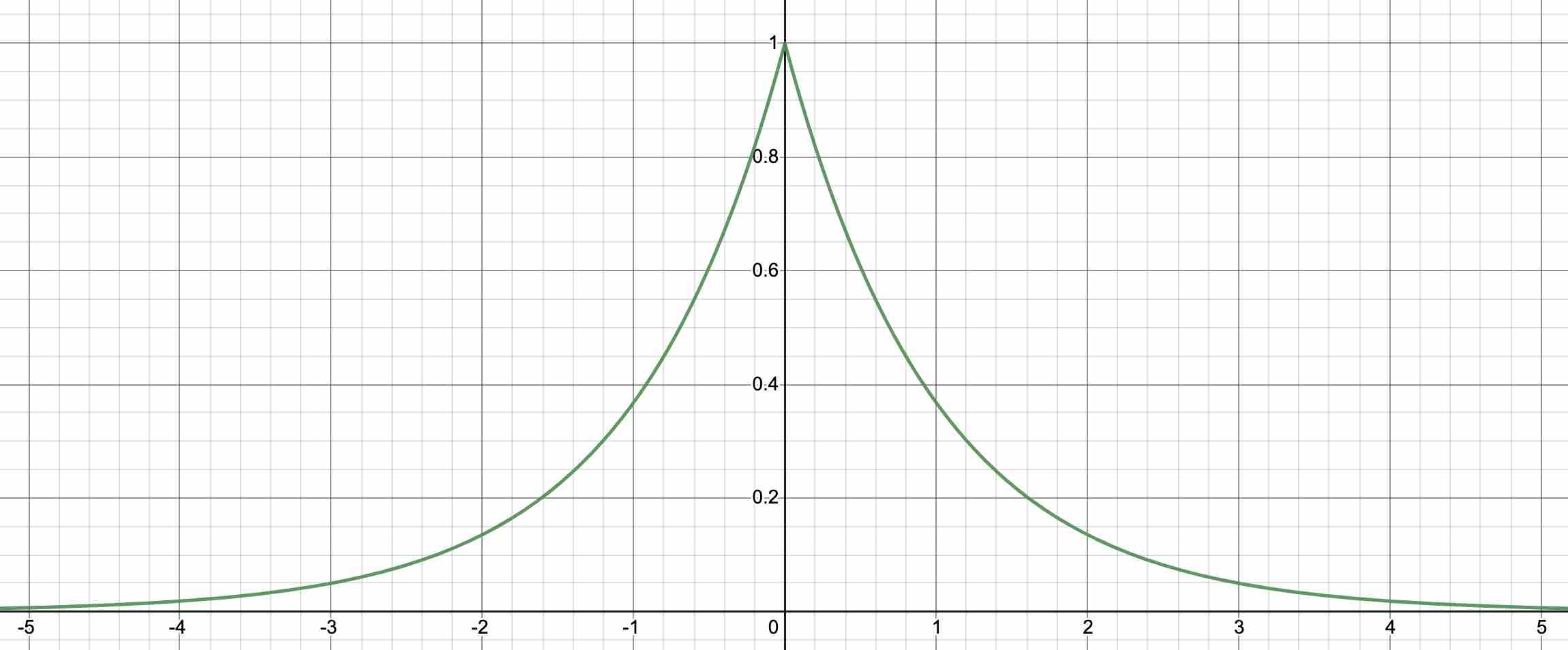}
\end{center}
\caption{The function  $e^{-|x|}$. }
\end{figure}

Interestingly, we can also see that $f$ is not a Schwartz function in another way - by computing its Fourier transform and observing that it is not rapidly decreasing:
\begin{align*}
\hat f(\xi)  &:= \int_\R e^{-2\pi t |x| -2\pi i x \xi }dx \\
&=    \int_{-\infty}^0   e^{2\pi t x -2\pi i x \xi } dx
+ \int_0^{+\infty}   e^{-2\pi t x -2\pi i x \xi } dx    \\
&= \int_{-\infty}^0   e^{2\pi  x(t - i \xi)    } dx
+ \int_0^{+\infty}   e^{-2\pi  x (t+ i \xi)  } dx    \\
&= \frac{  e^{2\pi  x(t - i \xi)}    }{2\pi (t-i\xi)}\Big|_{x=-\infty}^{x=0} + 
\frac{  e^{-2\pi  x(t + i \xi)}    }{-2\pi (t+i\xi)}\Big|_{x=0}^{x=\infty} \\
&= \frac{  1   }{2\pi (t-i\xi)} + 
\frac{  1   }{2\pi (t+i\xi)} \\
&= 
\frac{  t   }{\pi (t^2 + \xi^2)},
\end{align*}
valid for all $\xi \in \R$.   Because the Fourier transform
\begin{equation}   \label{FT of the abs value exponential}
\frac{  t   }{\pi (t^2 + \xi^2)}
\end{equation}
 is not a rapidly decreasing function, we have another proof that $f$ is not a Schwartz function.

This example is interesting in that $f$ is infinitely differentiable everywhere, except at one point, namely
 $x =0$.  Yet this local lack of smoothness - at only a single point -  is enough to cause a global change in decay for its Fourier transform. 
 }
 \hfill $\square$
 \end{example}

It is just as easy to define Schwartz functions on $\R^d$ as well. 
For any $k:= (k_1, \dots, k_d) \in \Z_{\geq 0}^d$, we can define the multivariable differential operator
\[
D_k := \frac{\partial}{\partial x_1^{k_1} \cdots  \partial x_d^{k_d}}. 
\]
\begin{example}
In $\R^1$, this is the usual $k$'th derivative, namely
$D_k f(x) := \frac{d}{dx^k}f(x)$.  
In $\R^2$,  for example, we have $D_{(1,7)} f(x) :=  \frac{\partial}{\partial x_1 \partial x_2^7}f(x)$.
\hfill $\square$
\end{example}

The {\bf order} of the differential operator
$D_k$ is by definition $|k|:= k_1+ \cdots + k_d$. 
To define spaces of differentiable functions, we call a function $f : \R^d \rightarrow \C$ 
a $C^m$-function if all partial derivatives
 $D_k f$ of order $|k| \leq m$ exists and are continuous.     We denote the collection of all such 
 $C^m$-functions on Euclidean space by  $C^m(\R^d)$.   
 When considering {\bf infinitely-differentiable functions on Euclidean space}, we denote this space by $C^\infty(\R^d)$.

So we see that in $\R^d$,  we can define {\bf Schwartz functions} \index{Schwartz function}
 similarly to our previous definition: they are infinitely differentiable
 functions $f:\R^d \rightarrow \C$ such that for all vectors $a, k \in \Z_{\geq 0}^d$
 we have:
  \begin{equation}
 |x^a  D_k f(x)| \text{ is bounded on } \R^d,
 \end{equation}
 where $x^a:= x_1^{a_1} \cdots x_d^{a_d}$ is the standard multi-index notation.
We also define the {\bf Schwartz space} $S(\R^d)$ to be set of all Schwartz functions $f:\R^d \rightarrow \C$.

 \begin{thm} \label{Schwartz goes to Schwartz}
 The Fourier transform maps the Schwartz space $S(\R^d)$ one-to-one, onto itself.
 (See Exercise \ref{Schwartz space convolution invariance})
  \end{thm}
 In fact, more is true:  the mapping $f \rightarrow \hat f$ from $S(\R^d)$ to itself is an isometry.   
 The proof of this 
 fact uses the Parseval relation below. 
  And now that we know the definition of rapid decay, we see that an obvious consequence of 
 Corollary \ref{cor: f smoother implies FT of F decays faster} 
 is the following:  
 \begin{equation}  \label{smooth implies FT is rapidly decreasing}
 \text{ If } f \text{ is infinitely smooth, then } \hat f  \text{ is rapidly decreasing}. 
\end{equation}
In fact, we can combine some of the ideas above to record another useful fact. 
\begin{lem} \label{useful Schwartz fact} 
Let $\phi:\R^d \rightarrow \C$ be  compactly supported and infinitely smooth.
Then 
\[
\phi \in \S(\R^d).
\] 
\end{lem}
 \begin{proof}
 Because $\phi$ is compactly supported, we know that $\hat \phi$ is infinitely smooth (differentiation under the integral).  Moreover, the assumption that $\phi$ is infinitely smooth implies that
  $\hat \phi$ is rapidly descreasing,  by \eqref{smooth implies FT is rapidly decreasing}.   So now we know that 
  $\hat \phi$ is both rapidly decreasing and infinitely smooth - i.e. a Schwartz function.   
  Applying Theorem \ref{Schwartz goes to Schwartz}, we see that its Fourier transform is also a Schwartz function.  Namely,  using Fourier inversion, we conclude that $\hat{\hat \phi}(-x) = \phi(x) \in S(\R^d)$. 
 \end{proof} 
 The functions satisfying the conditions of Lemma \ref{useful Schwartz fact}  are also called {\bf bump functions}.   
The curious reader might ask:  `are there any functions at all that satisfy the condition of Lemma \ref{useful Schwartz fact}'?
The answer is that there are many, though we are almost always interested in their properties, rather than their explicit form (but see Appendix \ref{Appendix:bump functions and inner products}).


 \section{Poisson Summation I} \label{Poisson Summation section}
\index{Poisson summation formula}

We  introduce the Poisson summation formula, one of the most useful tools in analytic number theory, and in discrete / combinatorial geometry.  This version of Poisson summation holds for Schwartz functions.  
 There are many different families of sufficient conditions that a function $f$ can satisfy,  in order for Poisson summation to be applicable to $f$.


\begin{figure}[htb]
\begin{center}
\includegraphics[totalheight=2.2in]{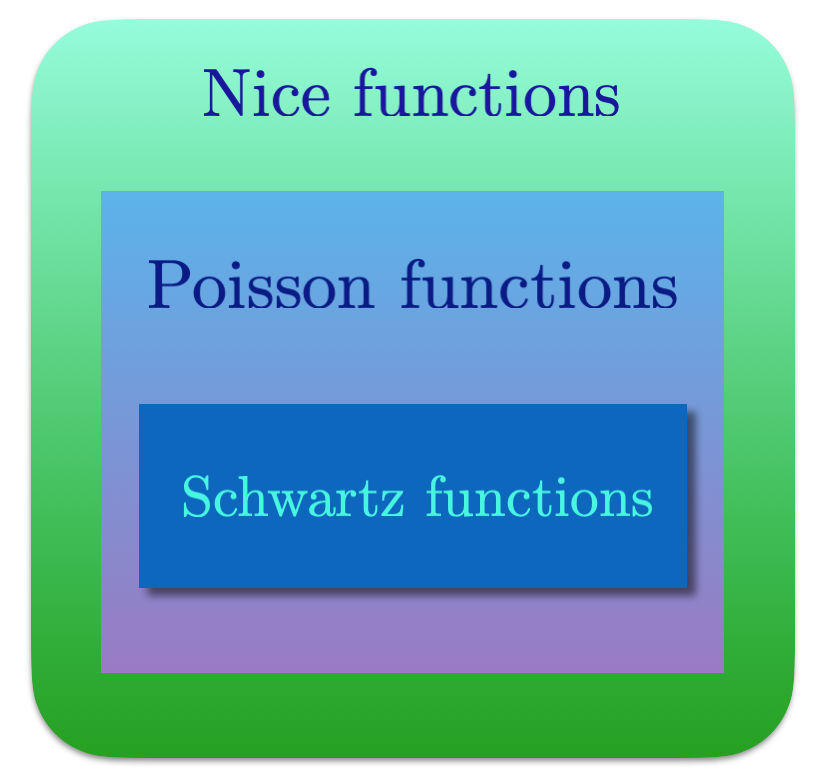}
\end{center}
\caption{Function spaces for Poisson summation}
\label{nice functions, containment}
\end{figure}


\bigskip
\begin{thm}[Poisson summation formula, I]
 \index{Poisson summation formula}   \index{Poisson summation formula} 
 \label{Poisson.Summation}
Given a    {\it Schwartz function}    $f~:~\R^d \rightarrow \C$, we have
\begin{equation}  \label{Poisson.summation1} 
\sum_{n \in \Z^d}  f(n+ x)  =  \sum_{\xi \in \Z^d} \hat f(\xi) e^{2\pi i \langle \xi, x \rangle},
\end{equation}
valid for all $x \in \R^d$.  
In particular, we have:
\begin{equation}  \label{Poisson.summation2}
\sum_{n \in \Z^d}  f(n)  =  \sum_{\xi \in \Z^d} \hat f(\xi).
\end{equation}
Both sides of \eqref{Poisson.summation1} converge absolutely, and are continuous functions on $\R^d$.
\end{thm}
\begin{proof}
If we let $F(x):= \sum_{n \in \Z^d}  f(n+ x)$, then we notice that $F$ is periodic on $\R^d$, with the cube $[0,1)^d$ as a fundamental domain.  The argument is easy: fix any $m \in \Z^d$. Then
$F(x + m) =  \sum_{n \in \Z^d}  f(n+ x + m) =  \sum_{k \in \Z^d}  f(x + k)$, because $\Z^d + m = \Z^d$. 
By Theorem \ref{Fourier series for periodic functions}, $F$ has a fourier series, so let's compute it:
\[
F(x) := \sum_{k \in \Z^d} a_k e^{2\pi i \langle k, x \rangle},
\]
where $a_k = \int_{[0,1)^d} F(u) e^{2\pi i \langle k, u \rangle}du$ for each fixed $k\in \Z^d$.
Let's see what happens if we massage $a_k$ a bit:
\begin{align}
a_k &:= \int_{[0,1)^d} F(u) e^{-2\pi i \langle k, u \rangle} du  \\
&=\int_{[0,1)^d} \sum_{n \in \Z^d}  f(n+ u) 
                    e^{-2\pi i \langle k, u \rangle} du  \\
&=\sum_{n \in \Z^d}   \int_{[0,1)^d}  f(n+ u)  \label{outersum}
                    e^{-2\pi i \langle k, u \rangle} du.
\end{align}
The interchange of summation and integral in the latter step is allowed by Theorem \ref{Application of dominated convergence}, which is an application of the dominated convergence theorem, 
 because the integrand satisfies
$ | f(n+ u)  e^{-2\pi i \langle k, u \rangle} |  =    | f(n+ u) | \in L^1(\R^d)$.  The latter absolute integrability of $f$ is
 due to the fact that $f$ is a Schwartz function. 

Now we fix an $n\in\Z^d$ in the outer sum of   \eqref{outersum}, and make the change of variable in the integral: 
$n+u := w$, so that $du = dw$.  
A critical step in this proof is the fact that as $u$ varies over the cube ${[0,1)^d}$, $w:= n+u$ varies 
over all of $\R^d$ because we have a tiling
\index{tiling}
 of Euclidean space by the unit cube:   ${[0,1)^d} + \Z^d = \R^d$.
We note that under this change of variable, 
$e^{-2\pi i \langle k, u \rangle} = e^{-2\pi i \langle k, w-n \rangle} = 
e^{-2\pi i \langle k, w \rangle}$, because $k, n \in \Z^d$ and hence 
$e^{2\pi i \langle k, n \rangle} =1$.
Therefore, we finally have:
\begin{align}\label{finally, the FT}
a_k = \sum_{n \in \Z^d}   \int_{n+ [0,1)^d}  f(w) 
                    e^{-2\pi i \langle k, w \rangle} dw
       = \int_{\R^d}  f(w) 
                    e^{-2\pi i \langle k, w \rangle} dw := \hat f(k),
\end{align}
so that $F(x) = \sum_{k \in \Z^d} a_k e^{2\pi i \langle k, x \rangle} 
= \sum_{k \in \Z^d} \hat f(k) e^{2\pi i \langle k, x \rangle}$.
\end{proof}

We define a function $f:\R^d \rightarrow \C$ to be a {\bf nice function} 
if both $f, \hat f \in L^1(\R^d)$, and if the Poisson summation formula
\begin{equation} \label{nice functions}
\sum_{n \in \Z^d}  f(n+ x)  =  \sum_{\xi \in \Z^d} \hat f(\xi) e^{2\pi i \langle \xi, x \rangle}
\end{equation}
holds for $f$ pointwise,  for each $x \in \R^d$.   In addition, we'll always assume absolute convergence of both sides of \eqref{nice functions}.

We will give various different sets of sufficient conditions for a
function $f$ to be nice.  
Figure \ref{nice functions, containment} suggests a simple containment relation between some of these function spaces, 
as we will easily prove.

There are a few things to notice about the classical, and pretty proof of Theorem \ref{Poisson.summation1}.  The first is that we began with any  square-integrable function $f$ defined on all of $\R^d$,
 and forced a periodization of it, which was by definition $F$.   This is known as the ``folding'' part of the proof. 
Then, at the end of the proof, there was the ``unfolding'' process, where we summed an integral over a lattice, and  because the cube tiles $\R^d$, the sum of the integrals transformed into a single integral over $\R^d$.  

The second thing we notice is that the integral $\int_{\R^d} f(x) e^{-2\pi i \langle x, \xi \rangle} dx$,
 which is by definition the Fourier transform of $f$,
appears quite naturally due to the tiling of $\R^d$ by the unit cube $[0,1)^d$.  
Hopefully there will now be no confusion as to the difference between the integral over the cube, and the integral over $\R^d$, both appearing together in this proof.

\section{Useful convergence lemmas, in preparation for Poisson summation II}

To prepare ourselves for Poisson's original summation formula, which we give in the next section, we will see here Poisson's hypotheses for the growth of $f$ and $\hat f$, together with the immediate convergence consequences they carry.

\begin{lem} \label{Poisson bound implies L^1}
Let $f:\mathbb R^d \rightarrow \C$ be a function that enjoys the bound 
\[
|f(x)|\leq\frac{C}{(1+||x||)^{d+\delta}}, 
\]
for all $x\in \mathbb R^d$, and for constants $C, \delta >0$ that are independent of $x$. 
Then $f\in L^1(\R^d)$.
\end{lem}
\begin{proof}
Consider the cube 
$Q_n:= [-n,n]^d$ and let $D_n:=Q_{n+1} - Q_n$ denote the set difference;  in other words, $D_n$ is the cubical shell between the cube $Q_n$ and the cube $Q_{n+1}$. 
We have $\R^d=\bigcup_{n\geq 0} D_n$, and $D_0=Q_1$.  Also, we note that on each shell $D_n$, $\frac{1}{\|x\|} \leq \frac{1}{n}$, so that: 
\begin{align}
    \int_{\R^d}|f(x)|dx
    &=\sum_{n\geq0}\int_{D_n}|f(x)|dx\\
    &=  \int_{D_0}|f(x)|dx  +  
         \sum_{n\geq 1}\int_{D_n}|f(x)|dx  \\
    &\leq\frac{C}{2^{d+\delta}}+\sum_{n\geq1}\int_{D_n}\frac{C}{(1+n)^{d+\delta}}dx \\
    &=\frac{C}{2^{d+\delta}}+\sum_{n\geq1}   \frac{C}{(1+n)^{d+\delta}} \int_{D_n} dx \\
    &=\frac{C}{2^{d+\delta}}+\sum_{n\geq1}\frac{C}{(1+n)^{d+\delta}}\left((2n+2)^{d}-(2n)^{d}\right)\\
    &=\frac{C}{2^{d+\delta}}+2^dC\sum_{n\geq1}\frac{1}{(1+n)^{d+\delta}}\left((n+1)^{d}-n^{d}\right)\\ \label{Big-O}
    &=\frac{C}{2^{d+\delta}}+\sum_{n\geq1}\frac{O(n^{d-1})}{(1+n)^{d+\delta}}\\
    &=\frac{C}{2^{d+\delta}}+\sum_{n\geq1}O\left(\frac{1}{n^{1+\delta}}\right)\\
    &=\frac{C}{2^{d+\delta}}+O\left(\sum_{n\geq1}\frac{1}{n^{1+\delta}}\right)<\infty,
\end{align}
where we've used the fact that the constant in the Big-O of equation \eqref{Big-O} is independent of $n$, so that we can move the series inside. 
\end{proof}

For the absolute summability of functions satisfying the same growth condition of the previous lemma, we have the following.
\begin{lem} 
\label{Poisson bound implies absolutely summable}
Let $f:\mathbb R^d \rightarrow \C$ be a function that enjoys the bound 
\[
|f(x)|\leq\frac{C}{(1+||x||)^{d+\delta}}, 
\]
for all $x\in \mathbb R^d$, and for constants $C, \delta >0$ that are independent of $x$. 
Then the series
\[
\sum_{k\in \Z^d} f(x+k) 
\]
converges uniformly and absolutely for all $x\in\R^d$.
\end{lem}
\begin{proof}
We will restrict attention to $x \in [0, 1)^d$, because the function $F(x):= \sum_{k\in \Z^d} f(k+x)$, if  convergent, forms 
a periodic function of $x \in \R^d$, with the unit cube $[0, 1)^d$ being a period.  We also note for all $x \in [0, 1)^d$,
we have the bound $\|x\| \leq \sqrt d$. 

We consider the tail of the series, for any given $N>0$:
\begin{align}
 |   \sum_{k \in \mathbb{Z}^{d}  \atop \|k\| > N}   f(k+x) | 
& \leq    \sum_{k \in \mathbb{Z}^{d}  \atop \|k\| > N}  \left | f(k+x) \right |  
  \leq C \sum_{k \in \mathbb{Z}^{d}  \atop \|k\| > N}      \frac{1}{(1+\|k+x\|)^{d+\delta}}  \\
& \leq   \sum_{k \in \mathbb{Z}^{d}  \atop \|k\| > N}     \frac{1}{  \left(1+\frac{\|k+x\|}{1+\sqrt d}  \right)^{d+\delta}} 
     \\ \label{second equality here}
&=\sum_{k \in \mathbb{Z}^{d}  \atop \|k\| > N}     \frac{\left(1+\sqrt{d}\right)^{d+\delta}}{\left(1+\sqrt{d}+\|k+x\|  \right)^{d+\delta}}
     \\  \label{third line here}
&\leq   C_{d, \delta}  \sum_{k \in \mathbb{Z}^{d}  \atop \|k\| > N}       \frac{1}{(1+\|k\|)^{d+\delta}}
     \\  \label{surface area of sphere approximation}
&=  C_{d, \delta}     \sum_{n \geq N}\frac{1}{(1+n)^{d+\delta}}  O \left(  n^{d-1} \right)    \\
&=\sum_{n \geq N}O\left(\frac{1}{n^{1+\delta}}  \right)\\
&=O\left(\sum_{n \geq N}\frac{1}{n^{1+\delta}}  \right) \rightarrow 0, \text{ as } N\rightarrow \infty,
\end{align}
and the last bound is independent of $x$.  In passing from \eqref{second equality here} to \eqref{third line here}, 
we used the estimate
$\|k + x \| \geq \|k \| - \| x\| \geq  \|k\|-\sqrt{d}$, and $C_{d, \delta}:= \left(1+\sqrt{d}\right)^{d+\delta}$.  The equality in \eqref{surface area of sphere approximation} is due to the fact that the number of integer points $k \in \Z^d$ that lie on a sphere of radius $n$ is $O\left(\text{surface area of } nS^{d-1}     \right) 
= O \left(  n^{d-1} \right)  $.  We've shown that the series 
$\sum_{k\in \Z^d}\left | f(k+x) \right |  $ converges uniformly on $\R^d$.
\end{proof}

We note that the only reason for having $(1+ \|x\|)^{d+ \delta}$ in the denominators of the bounds, instead of simply $\|x\|^{d+ \delta}$,
is to give simultaneously a bound at the origin, as well as any nonzero $x$.


\bigskip
\section{Poisson summation II, \'a la Poisson}

There are various different families of functions for which the adjective `nice'  applies, in \eqref{nice functions},
 and one of the simplest to understand  is the Schwartz class of functions. 
But there is a more general family of nice functions that is extremely useful, given by Poisson himself, as follows. 

 \begin{thm}[Poisson summation formula, II]   \index{Poisson summation formula}
 \label{nice2}
Suppose that for some positive constants $\delta$, $C$, and for all $x \in \R^d$, we have the bounds:
 \begin{align} \label{growth conditions for Poisson}
   |f(x)| < \frac{C}{     (1+\|x\|)^{d+\delta}    }
    \text{  \, and  \,  }     |\hat{f}(x)|  <         \frac{C}{     ( 1+\|x\|)^{d+\delta}   }.
\end{align}
Then we have the pointwise equality:
\begin{equation}  \label{Poisson summation, take 2}
\sum_{n \in \Z^d}  f(n+ x)  =  \sum_{\xi \in \Z^d} \hat f(\xi) e^{2\pi i \langle \xi, x \rangle},
\end{equation}               
for each $x\in \R^d$.  In addition,  both sides of \eqref{Poisson summation, take 2} 
converge absolutely, and are continuous functions on $\R^d$.          
\end{thm}
 \begin{proof}
Step $1$.  \  The growth conditions \eqref{growth conditions for Poisson} allow us to conclude that both $f, \hat f \in L^1(\R^d)$, by Lemma \ref{Poisson bound implies L^1}.   
This implies that both  $f, \hat f \in L^2(\R^d)$, by  the elementary Lemma \ref{both f and its FT in L^1 implies L^2}.
We also know that the Fourier transform of an $L^1$ function must be uniformly continuous on $\R^d$, and so both $f$ and $\hat f$ are uniformly continuous (Lemma \ref{uniform continuity}).  

Step $2$.  \ The hypothesis regarding the growth conditions  \eqref{growth conditions for Poisson} implies that the series defined by 
$F(x):= \sum_{n \in \Z^d}  f(n+ x)$ converges uniformly on $[0, 1]^d$, as we showed in Lemma \ref{Poisson bound implies absolutely summable}. 
It follows that this series must also converge in the $L^2$-norm on $[0, 1]^d$.  So  $F \in L^2(\torus)$, and it must therefore 
possess a Fourier series, which converges to it in the $L^2$-norm:
\begin{equation}
F(x) = \sum_{n \in \Z^d}  a_n \,  e^{2\pi i \langle n, x \rangle}.
\end{equation}

Step $3$. \  Next, we compute the Fourier coefficients $a_k$.  
This is almost the same step that already appeared in the proof of Theorem \ref{Poisson.Summation}, but we repeat it for completeness, and also because the interchange of sum and integral below is justified in a different way. 
\begin{align}
a_k &:= \int_{[0,1)^d} F(u) e^{-2\pi i \langle k, u \rangle} du  \\
&=\int_{[0,1)^d} \sum_{n \in \Z^d}  f(n+ u) 
                    e^{-2\pi i \langle k, u \rangle} du  \\
&=\sum_{n \in \Z^d}   \int_{[0,1)^d}  f(n+ u)  \label{outersum}
                    e^{-2\pi i \langle k, u \rangle} du.
\end{align}
The interchange of summation and integral in the latter step is allowed by the uniform convergence of the series
$ \sum_{n \in \Z^d}  f(n+ x)$.
We fix an $n\in\Z^d$ in the outer sum of   \eqref{outersum}, and make the change of variable in the integral: 
$n+u := w$. As $u$ varies over the cube ${[0,1)^d}$, $w:= n+u$ varies 
over all of $\R^d$ because the unit cube tiles the whole space:   
\[
{[0,1)^d} + \Z^d = \R^d.
\]
We also have 
$e^{-2\pi i \langle k, u \rangle} = e^{-2\pi i \langle k, w-n \rangle} = 
e^{-2\pi i \langle k, w \rangle}$, because $k, n \in \Z^d$ and hence 
$e^{2\pi i \langle k, n \rangle} =1$.
Finally:
$a_k = \sum_{n \in \Z^d}   \int_{n+ [0,1)^d}  f(w) 
                    e^{-2\pi i \langle k, w \rangle} dw
       = \int_{\R^d}  f(w) 
                    e^{-2\pi i \langle k, w \rangle} dw := \hat f(k)$.

Step $4$. \  Since each summand $f(n+x)$ is a continuous function of $x$, and since the convergence is uniform, the function $F(x)$ must also be continuous.
Finally, we'd like to pass from the convergence of the Fourier series in the $L^2$-norm,  to pointwise and uniform convergence.  For this task  we can use 
Lemma \ref{norm convergence plus absolute convergence implies equality}, 
assuming that we can show the absolute convergence of the Fourier series 
$\sum_{n \in \Z^d}  a_n \, e^{2\pi i \langle n, x \rangle} = 
\sum_{n \in \Z^d}  \hat f(n)  e^{2\pi i \langle n, x \rangle}$.
But this absolute convergence follows from the same Lemma \ref{Poisson bound implies absolutely summable}, with $f$ replaced by $\hat f$, because the same growth bounds \eqref{growth conditions for Poisson}
are also assumed for $\hat f$. 
To summarize this last step,  we know that $F$ is continuous, and the previous remarks allow us to use Lemma   \ref{norm convergence plus absolute convergence implies equality} to conclude that the Fourier series for $F$ converges pointwise and uniformly to $F(x)$. 
 \end{proof}
 
We call a function that enjoys the bounds  \eqref{growth conditions for Poisson} a {\bf Poisson function}, 
because Sim\'eon Denis Poisson  proved Theorem \ref{nice2} 
between $1823$ and   $1827$ \cite{Travaglini}.

Poisson's Theorem \ref{nice2} is a {\bf stronger} version of Poisson summation 
than  Theorem \ref{Poisson.Summation} above.  
 To justify this latter claim, we need to show that any Schwartz function also satisfies the growth conditions \eqref{growth conditions for Poisson}, but this is clear because Schwartz functions (and their transforms) decay faster than any polynomial, hence faster than the bounds given by \eqref{growth conditions for Poisson}.

We call the space of functions that satisfy the hypotheses of Theorem \ref{nice2}, the {\bf Poisson space} of functions, in honor of the mathematician that discovered this class.   As we've just seen, the suggestion of Figure \ref{nice functions, containment}  is correct, showing that the Schwartz space is contained in the Poisson space.
\begin{question}
Are there some natural necessary and sufficient conditions for Poisson summation?
\end{question}
This is an important open question. 
In other words, we may ask what are the inherent limitations of functions that satisfy Poisson summation?   Although there are well 
over $20$ different versions of sufficient conditions in the literature on Poisson summation, there are currently no known necessary and sufficient conditions for Poisson summation to hold.
It is natural to wonder what would happen if we only make the assumption that
\[
f\in L^1(\R^d)  \text{ and } \hat f\in L^1(\R^d)?
\]
 Is such an $f$ always a `nice'  function?
Sadly, the answer is ``no''  in general, and there is an important counterexample, 
by Yitzhak Katznelson (\cite{Katznelson}, Ch. VI, p. 143, Exercise 15).

\begin{figure}[htb]
\begin{center}
\includegraphics[totalheight=3.8in]{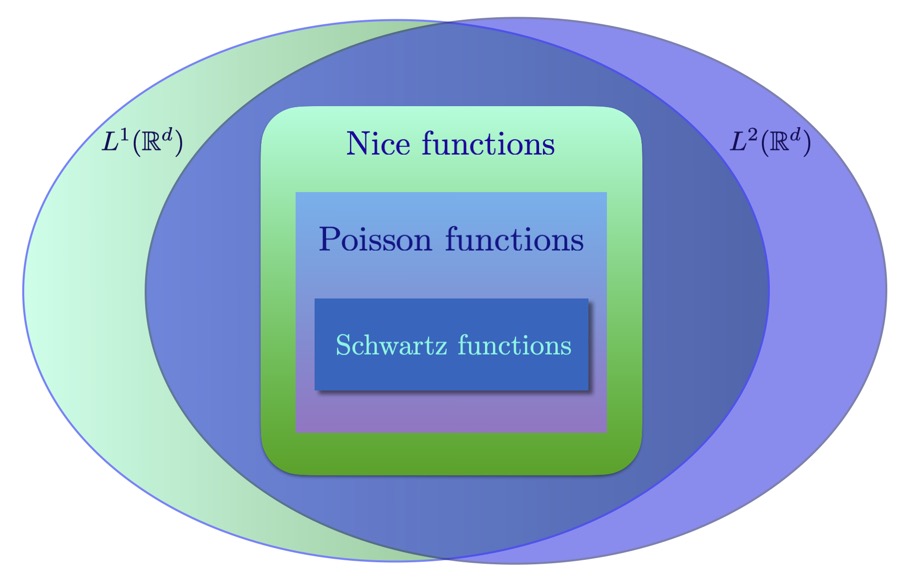}
\end{center}
\caption{A more detailed Venn diagram than Figure \ref{nice functions, containment},  for function spaces related to nice functions.}   
 \label{Refined nice functions}
\end{figure}

There are  many other families of nice functions in the literature, which include hypotheses such as  `functions of bounded variation', and `absolutely continuous' functions. We'll not delve into these other families here, but the reader may glance at 
Figure \ref{Refined nice functions} for a slightly more refined relationship between nice functions and the $L^1$ and $L^2$ spaces.
To justify the new containments that is suggested by Figure \ref{Refined nice functions}, we recall that a nice function $f$ was  defined  in \eqref{nice functions} to include the property that both $f, \hat f \in L^1(\R^d)$.   
By Lemma \ref{both f and its FT in L^1 implies L^2}, we know that therefore both
$f, \hat f \in L^2(\R^d)$ as well, so Figure \ref{Refined nice functions} is correct.


\bigskip
\section{An initial taste of general lattices, anticipating Chapter \ref{chapter.lattices}}

\begin{defi}\label{Def: first def of a lattice}
A   {\bf lattice}  \index{lattice}  is defined by the integer linear span of a fixed set of linearly independent vectors  $\{ v_1, \dots, v_m \} \subset \R^d$:
\begin{equation}\label{def.lattice}
\L :=  \left\{  n_1  v_1  + \cdots +  n_m v_m \in \R^d   \bigm |     \text{ all } n_j \in \Z    \right\}.
\end{equation}
\end{defi}
Although the {\bf integer lattice}  $\Z^d$ is the most common lattice, 
we often need to consider other types of lattices.
Any lattice $\L\subset \R^d$ can also be defined by:
\begin{equation}\label{alternate form of a lattice}
\L := \left\{   
     \begin{pmatrix} |  &  |  &  ...   & |  \\  
                        v_1  &  v_2  &  ...& v_m   \\  
                          |  &  |  & ... &  |  \\ 
        \end{pmatrix}
\begin{pmatrix}
n_1 \\
 \vdots \\
  n_m \\
\end{pmatrix}
\   \biggm |     \    
\begin{pmatrix}
n_1 \\
 \vdots \\
  n_m \\
\end{pmatrix}
   \in  \Z^m   
\right\}
 := M(\Z^m),
\end{equation}
where by definition, $M$ is the $d \times m$ matrix whose columns are the vectors $v_1, \dots, v_m$.  This set of basis vectors 
$\{ v_1, \dots, v_m\}$ is called a {\bf basis} \index{lattice basis}
for the lattice $\L$, and $m$ is called the {\bf rank} of the lattice $\L$.  
In this context, we also use the notation ${\rm rank}(\L) = m$.  
Any invertible matrix $M$ that appears in \eqref{alternate form of a lattice} is called 
a  {\bf basis matrix} \index{basis matrix}
for the lattice $\L$.
 Most of the time, we will be interested in {\bf full-rank}  \index{full rank lattice} 
 lattices, which means that $m=d$;  however, sometimes we will also be interested in lattices that have lower rank, and it is important to understand them.  
 The {\bf determinant} of a full-rank lattice  $\L := M(\Z^d)$ is defined by 
 \[
 \det \L := |\det M|,
 \]
 and we'll see in Chapter    that $\det \L$ is independent of the choice of basis matrix $M$.

 
\bigskip
\section{Poisson summation III, for general lattices}

We will use a slightly more general version of  the Poisson summation formula, which holds for any lattice, and which follows rather quickly from the Poisson summation formula above.  
We define a (full-rank) lattice $\L:= M(\Z^d) \subset \R^d$, the image of the integer lattice under an invertible linear transformation $M$.  
The {\bf dual lattice} \index{dual lattice} of $\L$ is defined by
$\L^* := M^{-T}(\Z^d)$, where $M^{-T}$ is the inverse transpose matrix of the real matrix $M$ (see Section \ref{dual lattice} for more on dual lattices).

As we've seen in Lemma \ref{FT under linear maps},  Fourier Transforms behave beautifully under compositions with any linear  transformation.
 We will use this fact again in the proof of the following extension of Poisson summation, which holds for 
 all lattices $\L$ and is quite standard.  We recall that a Poisson function $f$ by definition 
 satisfies the growth conditions  \eqref{growth conditions for Poisson}.

\begin{thm}  [Poisson summation formula, III]
\label{The Poisson Summation Formula, for lattices}
  \index{Poisson summation formula for lattices} 
Given a full-rank lattice $\L\subset \R^d$, and a Poisson function $f: \R^d \rightarrow \C$, we have
\begin{equation}  \label{PoissonSummationForLattices} 
\sum_{n \in \L} f(n+x) = \frac{1}{\det \L}  \sum_{m \in \L^*}  
                           \hat f(m)  e^{2\pi i \langle x, m \rangle},
\end{equation}
valid for all $x \in \R^d$.  In particular, we have
\begin{equation}  \label{Poisson.summation3}  \index{Poisson summation formula}
\sum_{n \in \L} f(n) = \frac{1}{\det \L}  \sum_{\xi \in \L^*}  
                           \hat f(\xi).
\end{equation}
Both sides of \eqref{PoissonSummationForLattices} converge absolutely and are continuous functions on $\R^d$.
\end{thm}
\begin{proof}
Any lattice (full-rank) may be written as $\L := M(\Z^d)$, so that $\det \L := |\det M|$.
Using the Poisson summation formula \eqref{Poisson.summation1}, with the change of variable $n = Mk$, 
with $k \in \Z^d$, we have:
\begin{align*} 
\sum_{n \in \L} f(n) &=  \sum_{k \in \Z^d} (f\circ M)(k)  \\
&=                              \sum_{\xi \in \Z^d} \widehat{(f\circ M)}(\xi)   \\
&=                    \frac{1}{|\det M|}     \sum_{\xi \in \Z^d}  \hat f\left(M^{-T}  \xi \right) \\
&=\frac{1}{\det \L} \sum_{m \in \L^*} \hat f(m).
\end{align*}
where in the third equality we used the elementary `Stretch' Lemma \ref{FT under linear maps}, and in the fourth equality we used the definition of the dual lattice $\L^*:= M^{-T} \Z^d$.
\end{proof}

As an afterthought, it turns out that the special case \eqref{Poisson.summation3} 
also easily implies the general case, namely \eqref{PoissonSummationForLattices} 
(Exercise \ref{going backwards in Poisson summation}). 

\begin{figure}[htb]
\begin{center}
\includegraphics[totalheight=3.4in]{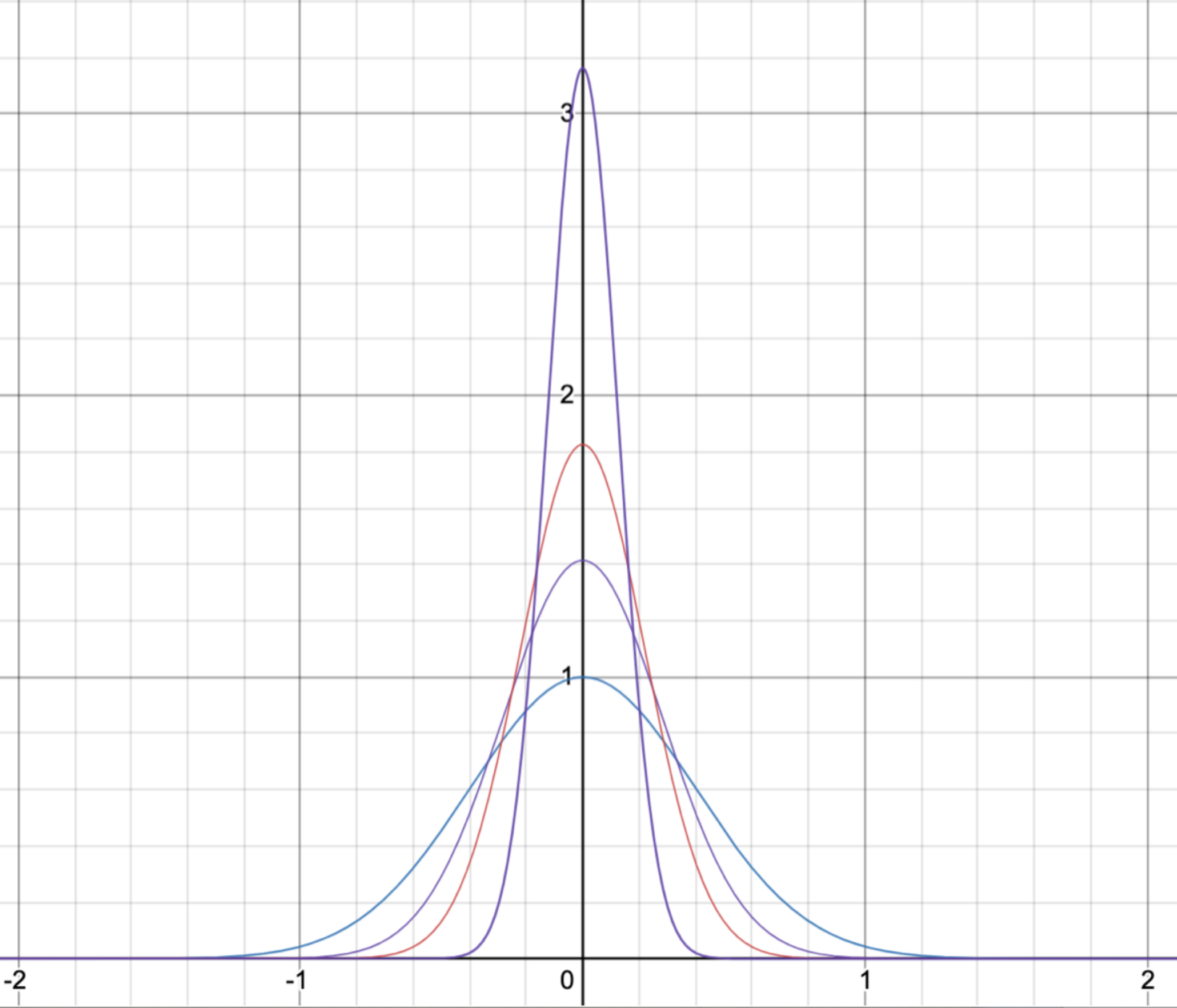}
\end{center}
\caption{The Gaussian family of functions $G_t(x)$ with $t = 1, t=.5,  t= .3, $ and $t=.1$ respectively. }   
 \label{pic of Gaussians}
\end{figure}

A traditional application of the Poisson summation formula is the quick derivation of 
 the functional equation of the theta function.  
We first define the Gaussian function by:
\begin{equation}  \index{Gaussian}
G_t (x) :=   t^{-\frac{d}{2}}    e^{ -\frac{\pi}{t}  || x ||^2 },
\end{equation}
for each fixed $t >0$, and for all $x \in \R^d$, as depicted in Figure \ref{pic of Gaussians}.

Two immediately interesting properties of the Gaussian are:
 \begin{equation}
\int_{\R^d}    G_{t} (x) dx = 1,
\end{equation}
for each $t>0$,  and 
\begin{equation}\label{transform of the Gaussian}
\hat G_t( m ) =  e^{ -\pi  t || m ||^2 },
\end{equation}
properties which are important in Statistics as well 
(Exercises \ref{Gaussian1} and \ref{Gaussian2}).  Each fixed $\varepsilon$ gives us one Gaussian function and intuitively, as $\varepsilon \rightarrow 0$, this sequence of Gaussians approaches the ``Dirac delta function'' at the origin, which is really known as a ``generalized function'', or ``distribution''  (Note \ref{Dirac delta}).

\begin{example} 
\rm{
The classical theta function \index{theta function}  (for the integer lattice) is defined by:
\begin{equation} \label{theta function for the integer lattice}
\theta(t) =  \sum_{n\in \Z^d}  e^{ -\pi  t || n ||^2 }.
\end{equation}
This function plays a major role in analytic number theory.  One of its first historical applications was carried out 
 by Riemann himself, who proved its functional equation (eq. \eqref{theta functional equation} below) and then applied a ``Mellin transform'' to it, to prove the functional equation of the Riemann zeta function $\zeta(s):= \sum_{n=1}^\infty \frac{1}{n^s}$. 
We claim that the theta function has the functional equation
\begin{equation} \label{theta functional equation}
\theta\left(  \frac{1}{t}  \right)   = t^{\frac{d}{2}} \theta(t),
\end{equation}
for all $t>0$.  This will follow immediately from the Poisson summation formula for Schwartz functions, namely \eqref{Poisson.summation2},  by using  \index{Poisson summation formula}
 $f(x):= G_t(x)$.   Using our knowledge of its FT, from
\eqref{transform of the Gaussian}, we have:
\begin{align*}
\sum_{n\in \Z^d}  G_t(n) &= \sum_{\xi \in \Z^d}  \hat G_t(\xi) \\
&= \sum_{\xi \in \Z^d}  e^{ -\pi  t || \xi ||^2 } := \theta(t).      
\end{align*}
Since by definition $\sum_{n\in \Z^d}  G_t(n) := 
t^{-\frac{d}{2}}   \sum_{n\in \Z^d} e^{ -\frac{\pi}{t}  ||n||^2 } 
:=  t^{-\frac{d}{2}} \theta\left(\frac{1}{t}\right)$,   \eqref{theta functional equation}   is proved.
}
  \hfill   $\square$
\end{example}



\bigskip
\section{The convolution operation}      \label{* is born}
\index{convolution} 
 For  $f,g \in L^1(\R^d)$,  their {\bf convolution} is defined by 
 \begin{equation} \label{def of convolution}
 (f * g)(x) = \int_{\R^d}   f(x-y) g(y) dy.
 \end{equation} 
But sometimes it is useful not to assume that we have absolutely integrable functions, and therefore 
we'll also use definition \eqref{def of convolution} to include any
functions $f, g$, for which the latter 
integral  still converges (see Examples \ref{ex heaviside},   \ref{ex ramp} below).   
It is possible to think intuitively of this analogue of multiplication as:  ``this is how waves like to multiply", via Lemma \ref{convolution theorem} \ref{convolution under FT}.    We have the following basic relations for the convolution operation.

 \begin{lem}   \label{convolution theorem}
For all $f, g, h  \in L^1(\R^d)$, we have:
 \begin{enumerate}[(a)]
 \item   $f*g \in L^1(\R^d)$.  
 	\label{part 1:convolution theorem}
 \item  $ \widehat{(f * g)}(\xi) = {\hat f}(\xi)  {\hat g}(\xi)$. 
 	 \label{convolution under FT}
 \item   $ f*g = g*f, \  \   f*(g*h)= (f*g)*h $, and $ \   f*(g+h) = f*g + f*h$.
 	\label{part 3:convolution theorem}
 \item   $\|f*g\|_1 \leq     \|f\|_1   \|g\|_1$.  
 	\label{part 4:convolution theorem}
\item  More generally, when $f\in L^p(\R^d), \ \ g \in L^1(\R^d)$, with $1\leq p <\infty$, then we have
$f*g \in L^p(\R^d)$ and
\[
 \|f*g\|_p \leq  \|f\|_p \|g\|_1.
\]
 \end{enumerate}
  \end{lem}
\begin{proof}
To prove part \ref{convolution under FT}, we use Fubini's Theorem 
(Theorem \ref{Fubini} in the Appendix):
\begin{align*}
\widehat{(f * g)}(\xi) &:= \int_{\R^d}  e^{-2\pi i \langle x, \xi \rangle}   \left(\int_{\R^d}f(x-y) g(y) dy \right)    dx \\
&= \int_{\R^d}  g(y)   e^{-2\pi i \langle y, \xi \rangle}  dy \int_{\R^d}f(x-y) e^{-2\pi i \langle x-y, \xi \rangle}  dx \\
&= \int_{\R^d}  g(y)   e^{-2\pi i \langle y, \xi \rangle}  dy \int_{\R^d}f(x) e^{-2\pi i \langle x, \xi \rangle}  dx \\
&:=  {\hat f}(\xi)  {\hat g}(\xi),
\end{align*}
where we've used the translation invariance of the measure, in the penultimate equality.  

To prove part \ref{part 4:convolution theorem}, we use Fubini's theorem again, and the triangle inequality for integrals: 
\begin{align*}
\| f * g\|_1&:= \int_{\R^d}    \left |    \int_{\R^d}f(x-y) g(y) dy   \right |  dx \\
&\leq    \int_{\R^d}     \int_{\R^d}    \left |    f(x-y) g(y) \right |  dy    dx \\
&=\int_{\R^d}     \int_{\R^d}    \left |    f(y) g(y) \right |  dy    dx  \\
&= \int_{\R^d}     \left |    f(y)  \right |  dy   \int_{\R^d}    \left |  g(y) \right |   dx  \\
&:=  \|f\|_1   \|g\|_1.
\end{align*}
For the proofs of the remaining parts,  we recommend Rudin's book \cite{RudinGreenBook}.
\end{proof}

 Lemma \ref{convolution theorem} \ref{convolution under FT}
  means that convolution of functions in the space domain corresponds to 
 the usual multiplication of functions in the frequency domain (and vice-versa).

 \begin{example}
 \rm{
 When  $\P := [-\frac{1}{2},  \frac{1}{2}]$, the convolution of $1_\P$ with itself is drawn in Figure \ref{pic of convolution of indicator}. We can already see that this convolution is a continuous function, 
 hence a little smoother than  the discontinuous function $1_\P$. 
 Using Lemma \ref{convolution theorem} we have
 \[
 \widehat{     (1_\P * 1_\P)     }(\xi) =    \hat{1}_\P(\xi)   \hat{1}_\P(\xi) = 
 \left( \frac{\sin(\pi \xi)}{\pi \xi} \right)^2.
 \]
We've used  equation \ref{ClassicalExample} in the last equality, for the Fourier transform of our interval $\P$ here.  Considering the graph in Figure \ref{pic of sinc2}, for the Fourier transform of the convolution $(1_\P * 1_\P)$, we see that this positive function is already much more tightly concentrated near the origin, as compared with
$\rm{sinc}(x):= \hat 1_\P(\xi)$.  We work out all of the details for this $1$-dimensional function, and generalize it,
 in Example \ref{hat function} below.

  \begin{figure}[htb]
\begin{center}
\includegraphics[totalheight=1.4in]{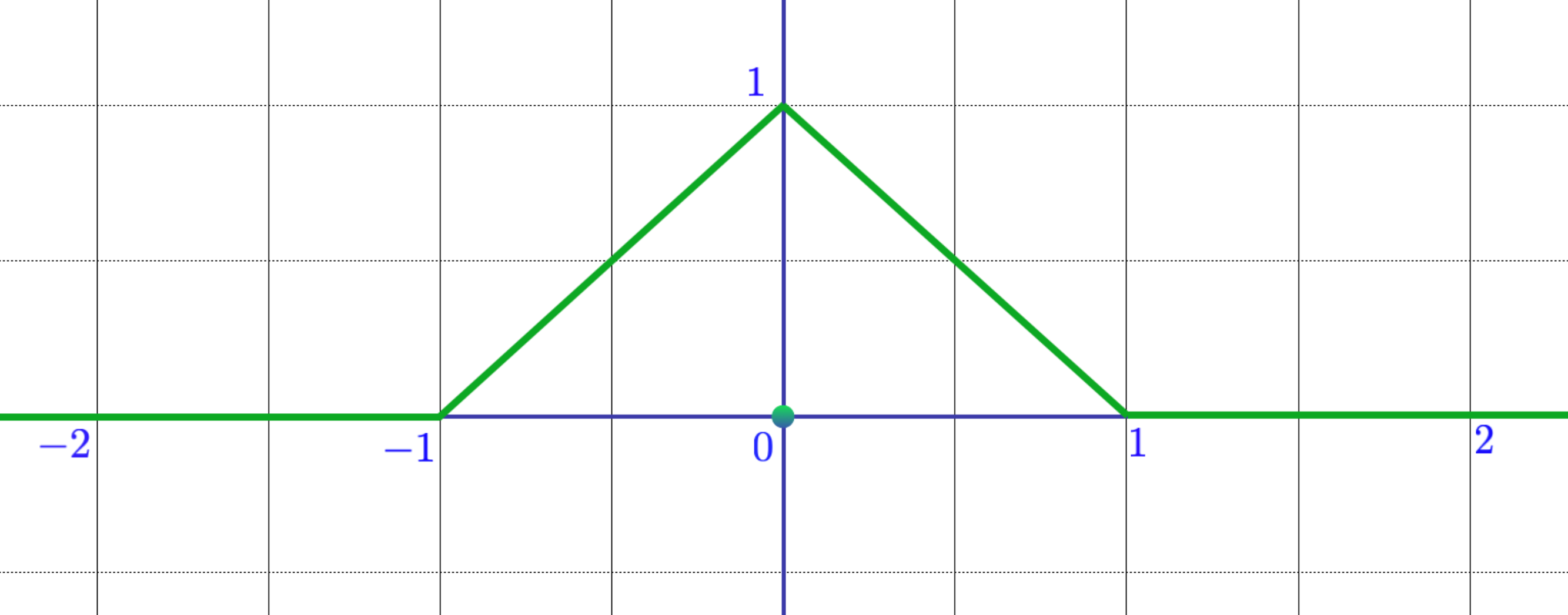}
\end{center}
\caption{The function $\left(  1_\P * 1_\P \right) (x)$, with $\P:= \left[ -\frac{1}{2}, \frac{1}{2} \right]$ }   
 \label{pic of convolution of indicator}
\end{figure}

 \begin{figure}[htb]
\begin{center}
\includegraphics[totalheight=1.5in]{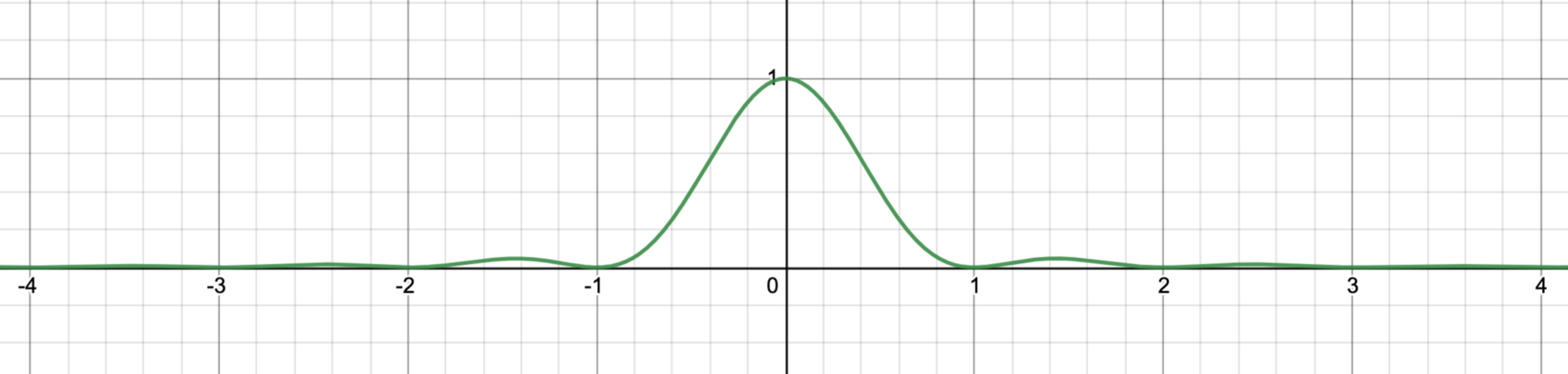}
\end{center}
\caption{The Fourier transform  $ \widehat{  \left(   1_\P * 1_\P   \right)   } (\xi)$, 
which is equal to the infinitely smooth, nonnegative function
 $\left( \frac{\sin(\pi \xi)}{\pi \xi} \right)^2  := \rm{sinc}^2(\xi)$. }   
 \label{pic of sinc2}
\end{figure}

}
 \hfill $\square$
 \end{example}
Another useful bit of intuition about convolutions  is that they are a kind of averaging process, and that 
 the convolution of two functions becomes smoother than either one of them.
 For our applications, when we consider the indicator function $1_\P(x)$ for a polytope $\P$,  
 then this function is not continuous on $\R^d$, so that the Poisson summation formula does not
 necessarily hold for it.  But if we consider the convolution of $1_\P(x)$ with a Gaussian, for example, 
 then we arrive at the $C^\infty$ function 
 \[
 (1_\P * G_t)(x), 
 \]
 for which the Poisson summation
  does hold.  In the sequel, we will use the latter convolved function in tandem with Poisson summation to study  ``solid angles".

\medskip
\begin{example} \label{convolution of general bodies}
\rm{
For any bounded measurable sets $K, L \subset  \R^d$, we have
\begin{align}
(1_K * 1_L)(y) &:= \int_{\R^d} 1_{K}(x) 1_{L}(y-x) dx \\
&:= \int_{\R^d} 1_{K}(x) 1_{-L+y}(x) dx \\
&= \int_{\R^d} 1_{K \cap (-L + y)}(x) dx \\
&= \int_{K \cap (-L + y)}     dx \\  \label{volume formula for convolution}
& = \vol\left( K \cap (-L + y) \right),
\end{align}
so that the convolution of indicator functions gives volumes, and this simple connection is one of the entry points of Fourier analysis into convex geometry.  
}
\hfill $\square$
\end{example}

\begin{figure}[htb]
\begin{center}
\includegraphics[totalheight=1.5in]{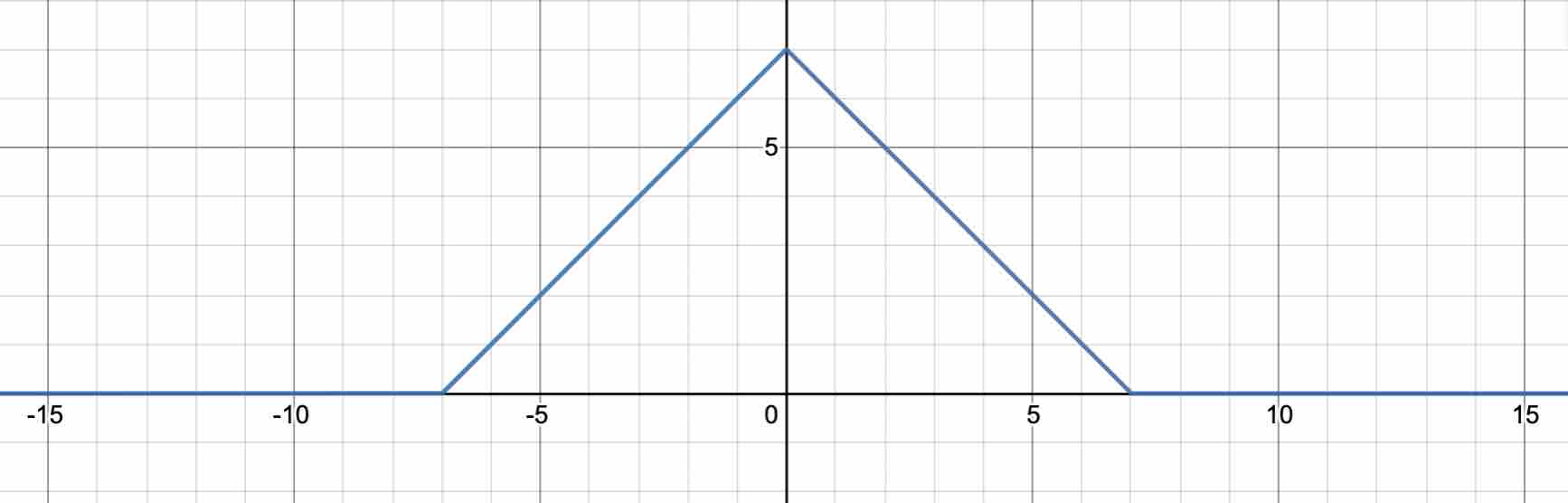}
\end{center}
\caption{The hat function $ 1_{\left[ -r, r \right]} * 1_{\left[ -r, r \right]}$ of Example \ref{hat function}, with $r = 3.5$.}   
 \label{pic of hat function, take 2}
\end{figure}

\medskip
\begin{example} \label{hat function}
\rm{
As a special case of Example \ref{convolution of general bodies},  consider the case
$K= L := \left[ -r, r \right] \subset \R$.  So we now know, 
by \eqref{volume formula for convolution}, that 
\begin{align}
g(x):=\left(1_{\left[ -r, r \right]} * 1_{\left[ -r, r \right]} \right)(x) = 
\vol\Big( \left[ -r, r \right] \cap (\left[ -r, r \right] + x) \Big),
\end{align}
making it clear that for $x \leq -2r$ and $x \geq 2r$, we have 
$\vol\Big( \left[ -r, r \right] \cap (\left[ -r, r \right] + x) \Big)=0$.  
Precisely,  when $x \in [-2r, 0]$, we have the function 
\[
g(x):= \vol\Big( \left[ -r, r \right] \cap (\left[ -r, r \right] + x) \Big) = | x-2r | = x+ 2r,
\] 
Finally, when $x \in [0, 2r]$, we have the function 
$g(x):= \vol\Big( \left[ -r, r \right] \cap (\left[ -r, r \right] + x) \Big) = | x-2r | = 2r - x $.  To summarize, we have
\[
g(x) = 
\begin{cases}
2r-|x| & \text{ if } x \in [-2r, 2r] \\
0  &  \text{ if not.}
\end{cases}
\]
Due to its shape, $g$ is sometimes called the {\bf hat function}, which is clearly a continuous function on $\R$. 
The hat function is extremely useful in many applications.  For example, we can use it to build up functions that are compactly supported on $\R$, and yet whose Fourier transform is 
\emph{strictly positive} on $\R$ - see Exercise \ref{positive FT over R}. 
}
\hfill $\square$
\end{example}

Given these examples, it is natural to wonder when the convolution is continuous:
\begin{question}\rm{[Rhetorical]} 
\label{convolution question}
Given any convex sets $A, B\subset \R^d$, is $1_A* 1_B(x)$ continuous for all $x \in \R^d$?
\end{question}
We can answer question \ref{convolution question} in the affirmative, in 
Exercises \ref{continuity of convolution for two L^1 bounded functions}
 and \ref{continuity of convolution for two L^2 functions}.

\begin{figure}[htb]
\begin{center}
\includegraphics[totalheight=1.4in]{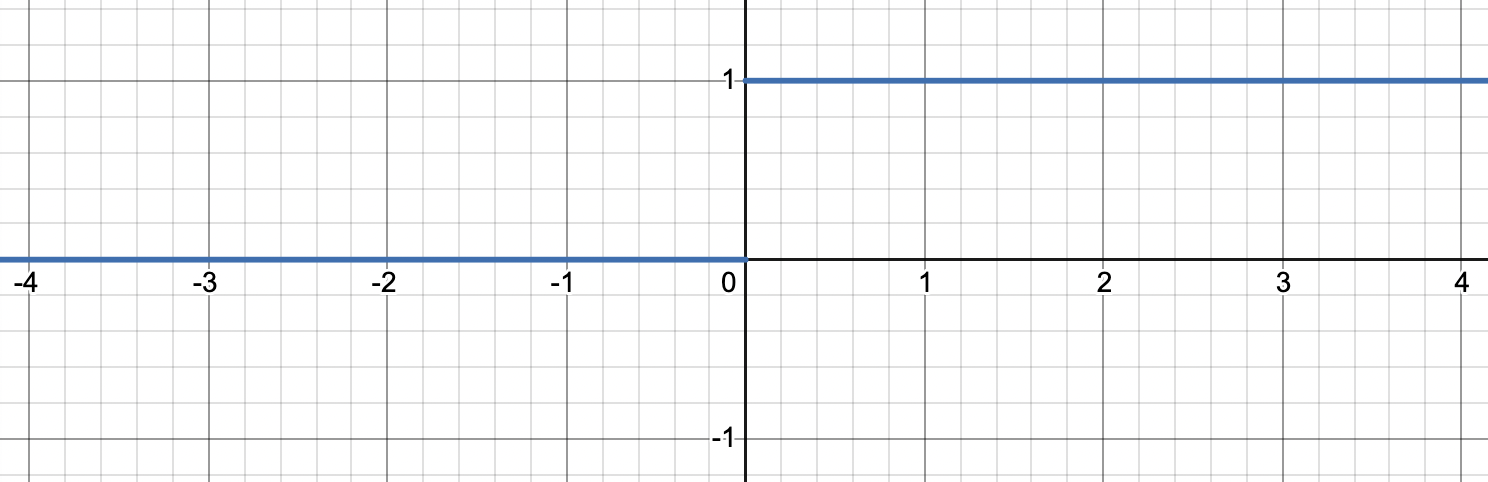}
\end{center}
\caption{The heaviside function $H_0(x)$ }   
 \label{Heaviside function}
\end{figure}

\bigskip
\begin{example} \label{ex heaviside}
\rm{
The {\bf Heaviside function} is defined by
\begin{equation}  \label{def of heaviside}
H_a(x):=   \begin{cases}  
1  &\mbox{if } x \geq a \\ 
0  & \mbox{if } x < a,
\end{cases}
\end{equation}
where $a$ is any fixed real number.   Although the Heaviside function is clearly not absolutely integrable over $\R$,
we may still use the same definition \eqref{def of convolution} for its convolution with a function  $f\in L^1(\R)$:
\begin{equation}\label{Heaviside convolution}
(f*H_0)(x):= \int_{\R} f(x-y) H_0(y)  dy=\int_{0}^\infty  f(x-y) dy=\int_{-\infty}^x  f(t) dt,
\end{equation}
a convergent integral.
}
\hfill $\square$
\end{example}

\begin{figure}[htb]
\begin{center}
\includegraphics[totalheight=2.4in]{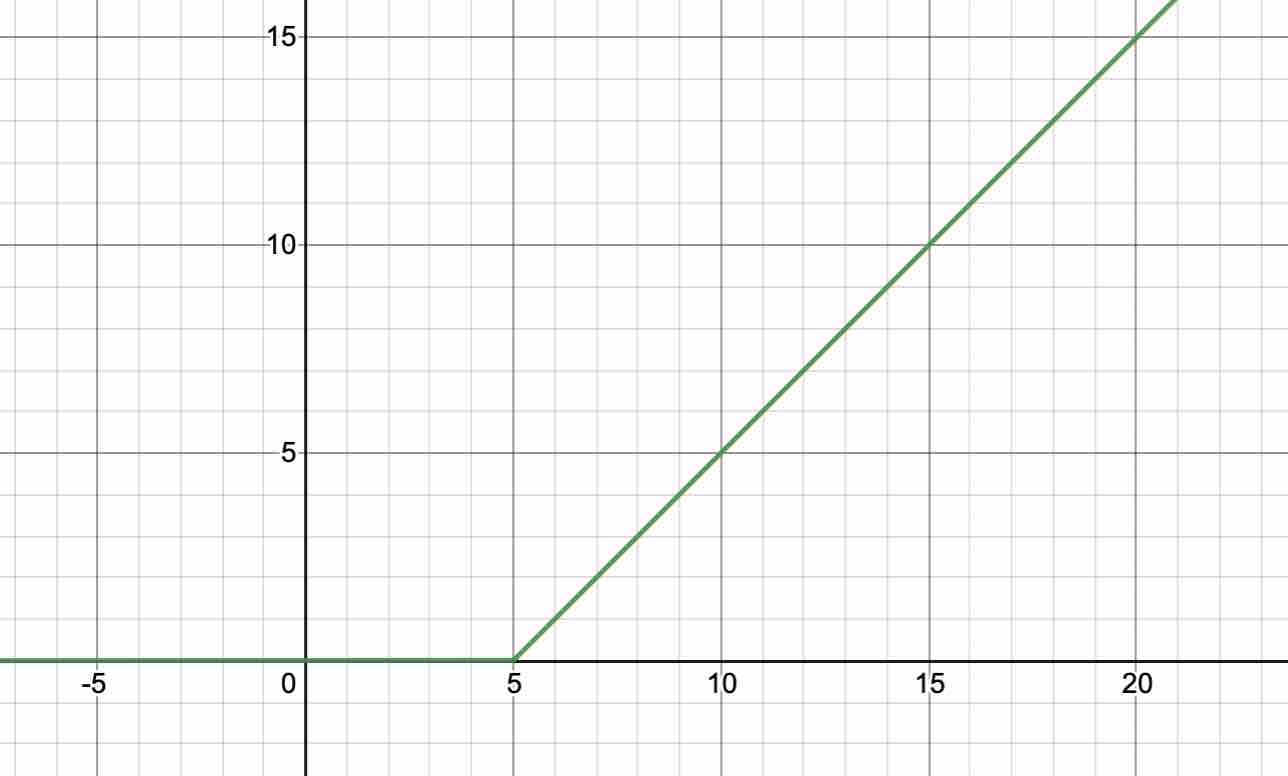}
\end{center}
\caption{The ramp function $r_5(x)$ }   
 \label{Ramp function}
\end{figure}

\bigskip
\begin{example} \label{ex ramp}
\rm{
The {\bf ramp function} is defined by 
\begin{equation} \label{def of ramp}
r_a(x):=   \begin{cases}  
x  &\mbox{if } x \geq a \\ 
0  & \mbox{if } x < a,
\end{cases}
\end{equation}
where $a$ is any fixed real number.  It is evident that we also have
$r_0(x) = \max\{ x, 0\}$.  It is also clear that $r'_a(x) = H_a(x)$.   The ramp function is ubiquitous in the analysis of machine learning algorithms, where it is called
the ReLu (Rectified Linear Unit) function.  
There is an elegant relationship between the ramp function and the Heaviside function:
\begin{equation} \label{claim for two heavisides}
H_0*H_0 = r_0,
\end{equation}
so we see that convolution makes sense here despite the fact that none of these functions are in $L^1(\R)$!
To check the latter claim \eqref{claim for two heavisides}, we use \eqref{Heaviside convolution} above:
\begin{align*}
H_0*H_0(x)&:=   \int_{-\infty}^x  H_0(t) dt =
 \begin{cases}  
 \int_0^x dx &\mbox{if } x \geq 0 \\ 
0  & \mbox{if } x<0
\end{cases} \\
&= 
\begin{cases}
x &\mbox{if } x \geq 0 \\ 
0  & \mbox{if } x<0
\end{cases} 
:= r_0(x).
\end{align*}
There is also a straightforward extension:  $H_a*H_b = r_{a+b}$ (Exercise \ref{Heaviside and ramp}). 
}
\hfill $\square$
\end{example}

\bigskip
\subsection{The support of a convolution}

Given two functions whose support is bounded, it's natural to wonder what the support of their convolution looks like. 
The very first observation is that if we have two closed, convex bodies $A, B \subset \R^d$, then:
\begin{equation}
\supp\left( 1_A* 1_B \right) = A + B,
\end{equation}
where the right-hand side uses the Minkowski sum of two sets (Exercise \ref{support of convolution}). 
There is a deeper result by Titchmarsh \cite{Titchmarsh2} in the case of $d=1$, and J. L. Lions \cite{Lions} in general dimension, that gives a very precise answer.
\begin{thm}[Titchmarsh and Lions]  \label{thm:Lions}
Let $f, g:\R^d \rightarrow \R$ have bounded support.   Then:
\begin{equation}
\conv\left(\supp (f*g) \right)= \conv\left( \supp f\right)  + \conv\left( \supp g\right), 
\end{equation}
\end{thm} 
where the right-hand side means we are taking the Minkowski sum of two convex bodies. 
\hfill $\square$

\medskip
The proof of this theorem is beyond the scope of this book,  although
 Exercise \ref{support of convolution}  gives another useful special case of Theorem \ref{thm:Lions}.

\bigskip
\section{More relations between $L^1(\R^d)$ and $L^2(\R^d)$}
\label{Section: more relations between L^1 and L^2}

Having seen convolutions, with various examples, we can now  return to the question:  
\begin{question}
What is the image of the space $L^1(\R^d)$ under the Fourier transform?
\end{question}

It seems that there is no known  `complete' answer to this open question yet;  however, an apparently lesser-known 
but elegant result, due to W. Rudin, is the following correspondence.
\begin{thm}[Rudin]   \label{RudinAmazingConvolutions}
\begin{equation}
f \in L^1(\R^d) \iff   \hat f = g*h, \text{ with } g, h \in L^2(\R^d).
\end{equation}
\hfill $\square$
\end{thm} 

In words, Theorem \ref{RudinAmazingConvolutions} tells us that 
the image of $L^1(\R^d)$ under the Fourier transform consists precisely of the set of convolutions 
$g*h$, where $g, h \in L^2(\R^d)$
(See  \cite{RudinGroups}, Theorem 1.6.3,  p.~27).  

Here is an outline of a proof for the easy direction: suppose that $g, h \in L^2(\R^d)$.   Because we want to find
 a solution in $f$, to the equation $\hat f = g*h$, it's natural to try  $f := \widehat{g*h} = \hat g \cdot \hat h$.  Let's try it, by defining
 \[
 f:= \hat g \cdot \hat h. 
 \]
 Because the Fourier transform acting on $L^2(\R^d)$ is an isometry,
 we have $\hat g, \hat h \in L^2(\R^d)$.  Also, the product of two $L^2$ functions in an $L^1$ function 
 (eq. \eqref{product of two L^2 functions is L^1}), so we conclude that 
 $f:= \hat g \cdot \hat h \in L^1(\R^d)$, as required.

This ongoing dance between the $L^1$ and $L^2$ spaces has more to offer. 
\begin{lem}  \label{both f and its FT in L^1 implies L^2} 
If $f \in L^1(\R^d)$ and $\hat f \in L^1(\R^d)$, then both $f, \hat f \in L^2(\R^d)$.  
\end{lem}
\begin{proof}
Because $\hat f \in L^1(\R^d)$, we know by the basic inequality \eqref{first bound for the FT} that $f$ must be bounded on $\R^d$:  $|f(x)| \leq M$ for some $M>0$.   We now compute:
\[
\int_{\R^d} |f(x)|^2 dx \leq  \int_{\R^d} M |f(x)| dx= M  \int_{\R^d}  |f(x)| dx < \infty,
\]
where the last inequality holds because $f \in L^1(\R^d)$ by assumption.   So $f \in L^2(\R^d)$.  
Precisely the same reasoning applies to $\hat f$, so that $ \hat f \in L^2(\R^d)$ as well.
\end{proof}

\medskip
Sometimes we are given a function $f \in L^2(\R^d)$, and we would like to know what extra properties $f$ needs to possess in order to place it in $L^1(\R^d)$.
\begin{lem}  \label{lem:from L^2 to L^1 and vice versa} 

\begin{enumerate}[(a)] 
\item  \label{part (a):L^2 to L^1}
Suppose that  $f$ vanishes outside a compact set $E \subset \R^d$.  Then:
\[
f  \in L^2(\R^d) \implies  f \in L^1(\R^d).
\] 
\item   \label{part (b):L^1 to L^2}
Suppose that $f$ is bounded on $\R^d$.   Then: 
\[
f  \in L^1(\R^d)  \implies f  \in L^2(\R^d).
\]
\end{enumerate}
\end{lem}
\begin{proof} To prove part \ref{part (a):L^2 to L^1}, we may use the Cauchy-Schwartz inequality:
\begin{align*}
\int_{\R^d} \left|f(x)\right|dx &= \int_{E} 1\cdot \left|f(x)\right|dx 
\leq    \left(   \int_E dx               \right)^{\frac{1}{2}}
          \left(   \int_E  |f(x)|^2 dx  \right)^{\frac{1}{2}}  \\
&= m(E)^{\frac{1}{2}}  \left(   \int_E  |f(x)|^2 dx  \right)^{\frac{1}{2}} < \infty.
\end{align*}
To prove part \ref{part (b):L^1 to L^2}, suppose $|f| < M$, for some bound $M>0$.   
Then $|f^2(x)| \leq M |f(x)|$ for all $x\in \R^d$, and we therefore have:
\[
\int_{\R^d} |f(x)|^2 dx \leq M \int_{\R^d} |f(x)| dx < \infty.
\]
\end{proof}

 \bigskip
 \subsection{How natural is the Fourier transform?}

 We close this section by thinking a bit about another natural question.  We've already seen in 
  Lemma \ref{convolution theorem} \ref{convolution under FT}
  that if $f, \hat f \in L^1(\R^d)$, then for each fixed $\xi \in \R^d$,  the map 
 \[
 \Phi_\xi:  f \rightarrow \hat f(\xi),
 \]
 is a complex homomorphism  from $L^1(\R^d)$ to $\C$. In other words, we already know that  $\Phi_\xi(f g) := \widehat{(fg)}(\xi) = 
 \hat f(\xi) \hat g(\xi) :=  \Phi_\xi(f)  \Phi_\xi(g)$.

Are there other linear transforms that act on $L^1(\R^d)$ as a homomorphism into $\C$?  It turns out there are not! 
 The Fourier transform is the unique homomorphism here, which means that it is very natural,
 and in this algebraic sense the Fourier transform is unavoidable. So we may as well befriend it.
\begin{thm} \label{Rudin's theorem about homomorphisms of L^1}
Suppose $\phi: L^1(\R^d)  \rightarrow \C$ is a nonzero complex homomorphism.
Then for each $f \in  L^1(\R^d)$, there exists a unique $t \in \R^d$ such that 
\[
\phi(f) = \hat f(t).
\]
\hfill $\square$
\end{thm}
The reader may consult Rudin's book \cite{RudinGreenBook}, Theorem 9.23, for a detailed proof. 
There is also a much more general version of Theorem \ref{Rudin's theorem about homomorphisms of L^1},  in the context of any locally compact abelian group, which Rudin proves in his book ``Fourier analysis on groups" (\cite{RudinGroups}, Theorem 1.2.2,  p.~7).

\bigskip
\section{The Dirichlet Kernel}
Using convolutions, we may now also go back to the partial sums of a Fourier series, which we have defined in \eqref{partial sums}
by  
  \begin{equation}
S_N f(t):= \sum_{n= -N}^N  \hat f(n) e^{2\pi i n t}.
\end{equation}
We compute:
\begin{align*}
S_N f(t)
&:= 
\sum_{n= -N}^N  \hat f(n) e^{2\pi i n t} =  
 \sum_{n= -N}^N  \int_0^1   f(x) e^{-2\pi i x n} dx \, e^{2\pi i n t} \\
&=
 \int_0^1   f(x)    \sum_{n= -N}^N   e^{2\pi i (t-x) n} dx\\
&:= (f*D_N)(t),
\end{align*} 
 where this convolution is defined on the $1$-Torus (the circle), and where we introduced the important definition
\begin{equation}
D_N(x):=  \sum_{n= -N}^N   e^{2\pi i x n},
\end{equation}
known as the {\bf Dirichlet kernel}.  
But look how naturally another convolution came up!  We've just proved the  following elementary Lemma.
\begin{lem}
If $f \in L^2(\mathbb T)$, then
\[
S_N f(t) = (f*D_N)(t),
\]
where this convolution is taken over $[0, 1]$.
\hfill $\square$
\end{lem}
It's therefore very natural to study the behavior of the Dirichlet kernel on its own.  In Exercise \ref{first Dirichlet kernel}, we showed that
the Dirichlet kernel has the closed form
\[
D_N(x) = \frac{\sin \left( \pi x(2N + 1) \right) }{\sin(\pi x)}.
\]
 \index{Dirichlet kernel}

\begin{figure}[htb]
\begin{center}
\includegraphics[totalheight=3.2in]{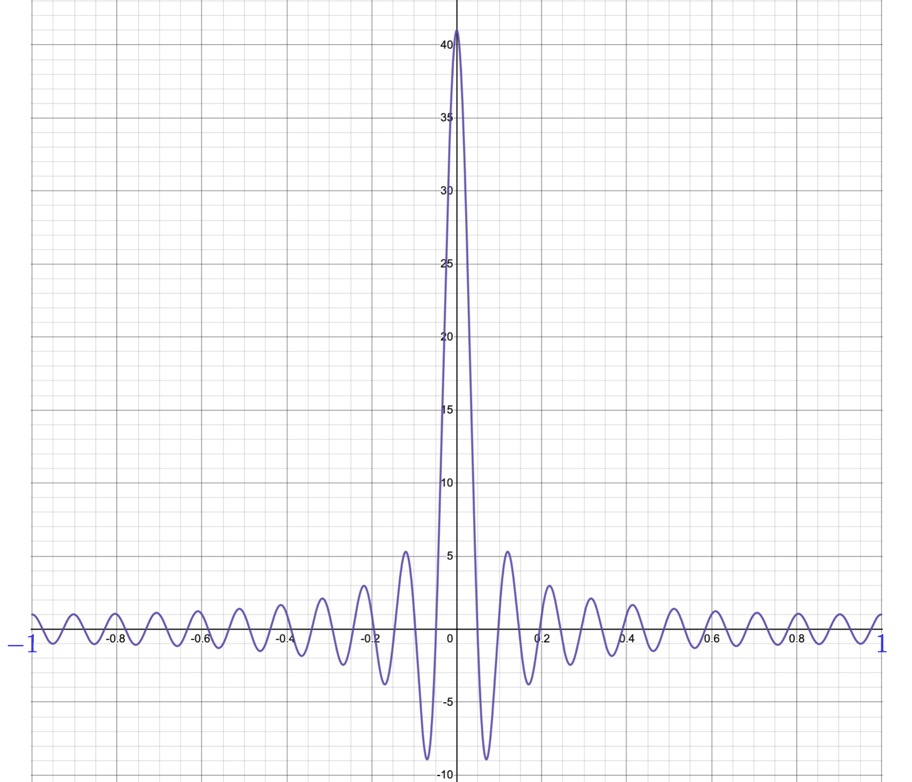}
\end{center}
\caption{The Dirichlet Kernel $D_{20}(x)$, restricted to the interval $[-1, 1]$ }   
 \label{pic of Dirichlet Kernel}
\end{figure}
 
It's clear from the definition of $D_N(x)$ that it is a periodic function of $x$, with period $1$, and if we restrict our attention to the interval $[-1, 1]$, then its graph appears in Figure \ref{pic of Dirichlet Kernel}.  It turns out the the $L^1$ norm of the Dirichlet kernel becomes unbounded as $N\rightarrow \infty$, and this phenomenon is responsible for a lot of results about pointwise divergence of Fourier series, a very delicate subject that is replete with technical subtleties.
There are even examples of continuous functions $f$ whose partial Fourier sums 
\index{partial Fourier sums}
$S_N f ( x )$  do not converge anywhere (\cite{Travaglini}, Theorem 4.19).
However, the Dirichlet kernel is also useful for proving pointwise
convergence theorems, such as the important 
Theorem \ref{theorem:Fourier series convergence to the mean}.


\section{The extension of the Fourier transform to $L^2$: \\
Plancherel}  
\index{Plancherel Theorem}  
So far we've worked with the Fourier transform that is defined only for functions that belong to $L^1(\R^d)$.   But sometimes we have a function that is \emph{not} in $L^1$, but 
we'd still like to study its transform.  Our prime example was 
$\hat 1_C$, a function that is not absolutely integrable for any bounded set $C$.

So how do we extend the Fourier transform to all of $L^2(\R^d)$?  Plancherel did it, and as it turns out there is a unique extension, because $L^1(\R^d) \cap L^2(\R^d)$ is dense in $L^2(\R^d)$.  
\begin{thm}[Plancherel, 1911]   
\label{thm:Plancherel}
Let $f \in L^2(\R^d)$.   Then there is a map 
\[
\Phi: L^2(\R^d) \rightarrow  L^2(\R^d)
\]
such that $\Phi(f)$ has the following properties: 
\begin{enumerate}[(a)]
\item If $f\in L^1$ as well, then we already have a proper definition of its 
Fourier transform, so we set  $\Phi(f):=\hat f$.
\item  (Plancherel's formula) $\| \Phi(f) \|_2 = \| f \|_2$. \label{Plancherel's magical formula}
\item $\Phi$ is a surjective Hilbert space isomorphism.  
\item  (Fourier inversion for $L^2$) We define:
\begin{align*}
A_R(\xi)&:= \int_{\|x\| < R} f(x) e^{-2\pi i \langle x, \xi \rangle} dx, \\
B_R(x)  &:= \int_{\|x\| < R} \Phi(f) (\xi) \, e^{2\pi i \langle x, \xi \rangle} d\xi.
\end{align*}
Then we have:
\begin{align}
&\lim_{R\rightarrow \infty} \|A_R(\xi) - \Phi(f) \|_2 = 0, \label{first L^2 inversion formula}\\
\text{ and }& \lim_{R\rightarrow \infty} \|B_R(x) - f \|_2 = 0.  
                                                                  \label{second L^2 inversion formula}
\end{align}
\end{enumerate}
This unique extension $\Phi(f)$ of the Fourier transform will henceforth be denoted by the same symbol: $\hat f$.
\hfill $\square$
\end{thm}
We refer the reader to Rudin's book \cite{RudinGreenBook} for a nice proof.
Equations \eqref{first L^2 inversion formula}
 and \eqref{second L^2 inversion formula} are the {\bf Fourier inversion formulas} for
$L^2$-functions.  We notice that `there is no free lunch' in following sense.  Although we were able to extend the Fourier transform to all of $L^2$, the convergence is \emph{not} pointwise convergence, but rather convergence in norm.  This sometimes causes some trouble, but it is part of life. 


At the risk of overstating the obvious, we note the good news that the equality in 
Plancherel's formula, which is part \ref{Plancherel's magical formula} of Plancherel's theorem above, is simply an equality between two real numbers.
Here is an interesting application for our focused study of indicator functions. 
\begin{cor}\label{FT of 1_Q is in L^2}
Given a bounded measurable set $Q\subset \R^d$, we have $\hat 1_Q  \in L^2(\R^d)$. 
Moreover:
\begin{equation}
\int_{\R^d} \left| \hat 1_Q(x) \right|^2 dx =  \vol Q.
\end{equation}
\end{cor}
\begin{proof} We know that $1_Q \in L^2(\R^d)$, so
by Plancherel's Theorem \ref{thm:Plancherel}, 
part \ref{Plancherel's magical formula}, we have
\[
\int_{\R^d} \left| \hat 1_Q(x) \right|^2 dx = \int_{\R^d} \left| 1_Q(x) \right|^2 dx 
= \int_Q dx =  \vol Q.
 \]
\end{proof}
Now let's consider a related function: $g(x):= (1_A * 1_B)(x)$, which we may equivalently
 rewrite as 
 \[
 (1_A * 1_B)(x) = \vol\left( Q \cap (-B + x) \right),
 \]
using \eqref{volume formula for convolution}.   
A natural question is whether or not $g$ is continuous.  
\begin{cor}
Let $A, B \subset \R^d$ be two bounded measurable sets. 
For $g(x):= (1_A * 1_B)(x)$, we have:
\begin{enumerate}[(a)]
\item $\hat g \in L^1(\R^d)$.
\item $g$ is continuous on $\R^d$. \label{second part of continuity of convolution}
\end{enumerate}
\end{cor}
\begin{proof}
 We know that trivially $1_Q  \in L^1(\R^d)$, and $1_{-Q}  \in L^1(\R^d)$, 
 and this implies that 
$\hat g := \F\left(   1_Q * 1_{-Q} \right) = \hat 1_Q \hat 1_{-Q}$, via Lemma \ref{convolution theorem} \ref{convolution under FT}. 
The Cauchy-Schwarz inequality now gives us:
\[
    \int_{\R^d}     |  \hat g(\xi)  |    \,   d \xi
      = \int_{\R^d}     |   \hat 1_Q(\xi)  |    |\hat 1_{-Q}(\xi) |    \,   d \xi
      \leq       \left(     \int_{\R^d}        |  \hat 1_Q  (\xi)  |^2        \,d \xi          \right)^{1/2}
                   \left(     \int_{\R^d}       |  \hat 1_{-Q}  (\xi)  |^2      \,d \xi         \right)^{1/2}
  < \infty,
\] 
where the finiteness of the last expression follows from 
Corollary \ref{FT of 1_Q is in L^2}.  A detailed proof of part 
\ref{second part of continuity of convolution}
is contained in Exercise  \ref{continuity of convolution for two L^2 functions}.
\end{proof}

\begin{cor}
\label{extended Plancherel identity}
For all $f, g \in L^2(\R^d)$, we have $\langle f, g\rangle =  \langle \hat f, \hat g  \rangle$.  
In other words:
\begin{equation}
         \int_{\R^d} f(x) \overline{g(x)} dx = \int_{\R^d} \hat f(x) \overline{\hat g(x)} dx.
\end{equation}
\hfill $\square$
\end{cor}
\begin{proof}
The elementary polarization identity (Exercise \ref{Plancherel extended}) 
tells us that
\[
\langle f, g \rangle = \frac{1}{2} \left( \| f\|^2 + \| g \|^2 - \|f+g\|^2\right).
\]
By Plancherel's Theorem \ref{thm:Plancherel}, 
part \ref{Plancherel's magical formula}, we know that  
$\| f\|^2 =  \| \hat f\|^2,  \| g \|^2 =  \| \hat g \|^2$, and $\| f+g \|^2 =  \| \hat f+  \hat g) \|^2$, so 
we have $\langle f, g \rangle = \langle \hat f,  \hat g \rangle $.
\end{proof}

\medskip
\begin{example}
\rm{
As we recall, the sinc function, defined by 
\begin{equation*}
\sinc(\xi):=  \begin{cases}  
\frac{\sin(\pi \xi)}{\pi \xi},     &\mbox{if } \xi \not= 0 \\ 
1  & \mbox{if } \xi= 0,
\end{cases}
\end{equation*}
plays an important role (in many fields), and was our very first example of the Fourier transform of a polytope: $\sinc = \hat 1_{\left[-\tfrac{1}{2}, \tfrac{1}{2} \right]}$.
Here we'll glimpse another aspect of the importance of sinc functions, as an application of Plancherel's theorem. 
 Let's prove that
\begin{equation}
\int_\R   \sinc(x-n) \sinc(x-m) dx = 
 \begin{cases}  
1    &      \mbox{if } n=m  \\ 
0  &        \mbox{if } n\not=m.
\end{cases}
\end{equation}
Although $\sinc \notin L^1(\R)$, we do have
$\sinc = \hat 1_{\left[-\tfrac{1}{2}, \tfrac{1}{2} \right]} \in L^2(\R)$, 
by Corollary \ref{FT of 1_Q is in L^2}.
Using Plancherel's theorem, we know that $\F(\sinc(x-n))(\xi)$ is well-defined as an 
$L^2(\R^d)$ function, and Corollary \ref{extended Plancherel identity} gives us:
\begin{align*}
\int_\R   \sinc(x-n) \, \sinc(x-m) dx 
&= \int_\R   \F(\sinc(x-n))(\xi)\,  \overline{  \F( \sinc(x-m) )(\xi) }  d\xi \\
&= \int_\R   1_\P(\xi) e^{2\pi i \xi n} \, 1_\P(\xi) \overline{  e^{2\pi i \xi m}  }d\xi\\
&=  \int_\P   e^{2\pi i \xi (n-m)}  d\xi \\
&= \delta(n, m),
\end{align*}
where  $\P:= [-\frac{1}{2}, \frac{1}{2}]$, and where we've used 
 the orthogonality of the exponentials over $\P$
 (Exercise \ref{orthogonality for exponentials}).  
 So we see that the collection of functions 
 \[
 \left\{  \sinc(x - n) \bigm |  n \in \Z   \right\}
 \]
  forms an orthonormal collection of functions in the Hilbert space $L^2([-\tfrac{1}{2},  \tfrac{1}{2}] )$, relative to its inner product.  It turns out that when we study
 Shannon's sampling theorem, these translated sinc functions are in fact a \emph{complete} orthonormal basis for the Hilbert subspace
  of $L^2(\R)$ that consists of `bandlimited functions' (see Theorem \ref{Shannon}).
 }
\hfill $\square$
\end{example}





\section{Approximate identities} 

It is a sad fact of life that there is no identity in $L^1(\R^d)$ for the convolution product - in other words, there is no function
$h \in L^1(\R^d)$ such that 
\begin{equation}\label{if there was an identity}
f*h = f
\end{equation}
 for all $f \in L^1(\R^d)$.

 Why is that?  Suppose there was such a function $h\in L^1(\R^d)$.  Then taking the Fourier transform of both sides of \eqref{if there was an identity}, we would also
have   
\begin{equation} \label{taking FT's...}
 \hat f  \ \hat h= \widehat{f*h}  = \hat f, 
\end{equation}
for all $f\in L^1(\R^d)$.  Picking an $f$ whose transform is nowhere zero, we can divide both sides of \eqref{taking FT's...}
by $\hat f$, to conclude that $\hat h \equiv 1$, the constant function.   But by the Riemann-Lebesgue Lemma \ref{Riemann--Lebesgue lemma}, we know that $\hat h$ must go to $0$ as $|x| \rightarrow \infty$, which is a contradiction. 
\index{Riemann-Lebesgue lemma}

Nevertheless, it is still interesting to think about what would happen if we were able to apply the
inverse Fourier transform to $\hat h$, formally applying the Fourier transform to the 
equation $\hat h = 1$ to get:
\begin{equation}
h(x) = \int_{\R^d} e^{2\pi i \langle x, \xi \rangle} dx,
\end{equation}
 an extremely interesting integral that unfortunately diverges.  In note \ref{Dirac delta}, we mention briefly that such observations became 
 critically important for the development of generalized functions that do play the role of the identity for convolutions, and much more.
 
Although there is no identity element for convolutions, it turns out that using sequences of functions we can get close!   Here is how we may do it, and as a consequence we will be able to rigorously apply the Poisson summation formula to a wider class of functions, including smoothed versions of the indicator function of a polytope.

Fix a function $\phi \in L^1(\R^d)$, such that  $\int_{\R^d} \phi(x) dx = 1$.    
 Beginning with any such function $\phi$, we construct an {\bf approximate identity}  by defining the sequence of functions 
\begin{equation}\label{approximate identity}
\phi_n(x):= n^d \phi(n x),
\end{equation}
for each $n = 1, 2, 3, \dots$.
  It's easy to check that we also have   $\int_{\R^d}\phi_n(x) dx  = 1$, for all $n\geq~1$ (Exercise \ref{total mass 1}). 
   So scaling $\phi$ by these $n$'s has the effect of 
squeezing $\phi$ so that it is becomes concentrated near the origin, while maintaining a total mass of $1$. 
Then intuitively a sequence of 
such $\phi_n$ functions approach the ``Dirac delta-function" at the origin (which is a distribution, not a function).

There are many families of functions that give an approximate identity.  
In practice, we will seldom have to specify exactly which sequence $\phi_n$ we pick, because we will merely use the existence of such a sequence to facilitate the use of Poisson summation. 
Returning now to the motivation of this section, we can recover the next-best-thing to an identity for the convolution product, as follows.  


\begin{thm}\label{approximate identity convolution}
Suppose we are given a  function $f \in L^1(\R^d)$, such that $p \in \R^d$ is a point of continuity for $f$.
Fix an approximate identity $\phi_n(x)$, and assume $f*\phi$ exists.   Then we have:
\begin{equation}\label{basic smoothing}
\lim_{n \rightarrow \infty}  \left(f * \phi_n \right)(p) = f(p).
\end{equation}
\end{thm}
\begin{proof} 
We begin by massaging the convolution product:
\begin{align*}
(\phi_n*f)(p) &:=  \int_{\R^d} \phi_n(x) f(p-x) dx \\
&=   \int_{\R^d} \phi_n(x) \Big(f(p-x) - f(p) + f(p) \Big) dx \\
&=          \int_{\R^d} \phi_n(x) \Big(f(p-x) - f(p) \Big) dx   + 
                 f(p) \int_{\R^d} \phi_n(x) dx       \\
&= f(p) +  \int_{\R^d} \phi_n(x) \Big(f(p-x) - f(p) \Big) dx,
\end{align*}
using the assumption that $\int_{\R^d} \phi_n(x) dx=1$.  
Using the definition of 
$\phi_n(x):= n^{d} \phi(n x)$, and making a change of variable 
$u= n x$ in the latter integral,  we have:
\begin{align*}
(\phi_n*f)(p) &:= f(p) +  \int_{\R^d} \phi(u) \Big(   f\left(p-   \frac{1}{n} u\right) - f(p) \Big) du.
\end{align*}
In the second part of the proof, we will show that as $n \rightarrow \infty$, the latter integral tends to zero.  
We will do this in two steps, first bounding the tails of the integral in a neighborhood of infinity, and then bounding the integral in a neighborhood of the origin.

Step $1$. \  Given any $\varepsilon >0$, we note that the latter integral converges, so the `tails are arbitrarily small'.  In other words, there exists an $r > 0$ such that
\[
\left|       \int_{\| u \| > r} \phi(u) \left(f\left(p-   \frac{1}{n} u\right) - f(p) \right) du \right|      < \varepsilon.
\]

Step $2$.  \  Now we want to bound $\int_{\| u \| < r} \phi(u) \left( f\left(p-\frac{1}{n} u\right)  - f(p) \right) du$.
We will use the fact that $ \int_{\R^d} | \phi(u) | du = M$, a constant.  Also, by continuity of $f$ at $p$, 
we can pick an $n$ sufficiently large, such that:
\[
\left| f\left(p-\frac{1}{n} u\right) - f(p) \right| < \frac{\varepsilon}{M},
\]
when $\|   \frac{1}{n} u   \| < r$.   Putting all of this together, and using the triangle inequality for integrals, 
we have the bound
\begin{align*}
\Big|     \int_{\| u \| < r} \phi(u) \left(f\left(p-\frac{1}{n} u\right)  - f(p) \right) du  \Big| 
     &\leq   \int_{\| u \| < r} |  \phi(u)   |  \left|  f\left(p-\frac{1}{n} u\right) - f(p) \right| du < \varepsilon.
\end{align*}
Therefore, as $n \rightarrow \infty $, we have $(\phi_n*f)(p) \longrightarrow f(p)$.
\end{proof}

We note that a point of discontinuity of $f$, Theorem \ref{approximate identity convolution} 
may be false even in dimension $1$, as the next example shows. 
\begin{example}
Let $f(x):= 1_{[0,1]}(x)$, which is discontinuous at $x=0$ and $x=1$.  We claim that 
for $p=1$, for example, we have 
\[
\lim_{n \rightarrow \infty}  (f * \phi_n)(p) = \frac{1}{2} f(p),
\]
so that the result of Theorem  \ref{approximate identity convolution}
 does not hold at this particular $p$, because $p$ lies on the boundary
of the $1$-dimensional polytope $[0,1]$.    When $p \in \interior([0,1])$, however, Theorem  
\ref{approximate identity convolution} does hold.
\hfill $\square$
\end{example}


\section{Poisson summation IV: a practical Poisson summation formula} 

In practice, we want to apply Poisson summation to indicator functions $1_\P$ of polytopes and general convex bodies.  
With this in mind, it's useful for us to have our own, home-cooked version of Poisson summation that is made for this culinary purpose. 

Throughout this section, we fix any compactly supported, nonnegative function $\varphi \in L^2(\R^d)$, with 
$\int_{\R^d} \varphi(x) dx = 1$, 
and we 
set $\varphi_\varepsilon(x) := \frac{1}{\varepsilon^d} \varphi   \left(   \frac{x}{\varepsilon}  \right)$, for each $\varepsilon > 0$.

\bigskip
\begin{thm}[Poisson summation formula IV]        \label{PracticalPoisson}
\index{Poisson summation formula}

Let $f(x) \in L^2(\R^{d})$ be a compactly supported function, and suppose that for each $x\in \R^d$, 
we have:
\begin{equation} \label{hypothesis, practical Poisson}
f(x)  =\lim_{\varepsilon\rightarrow0^{+}}      
\left(
\varphi_{\varepsilon}\ast f
\right)(x).
\end{equation}
Then the  following hold:

\begin{enumerate}[(a)]
\item
For each $\varepsilon>0$, we have absolute convergence:
$
\sum_{m\in\Z^d}      \left| \widehat{\varphi}\left(  \varepsilon
m\right)  \widehat{f}\left(  m\right)  \right| <+\infty.
$
\item           \label{practical Poisson summation, first version}
For all sufficiently small $\varepsilon >0$, and for each fixed $x\in \R^d$, we have the
pointwise equality:
\begin{equation}   
\sum_{n\in\mathbb{Z}^{d}}
\left(
\varphi_{\varepsilon}\ast f
\right)\left(  n+x\right)  
=
 \sum_{m\in\Z^d}             \widehat{\varphi}\left(
\varepsilon m\right)  \widehat{f}\left(  m\right)  
e^{2\pi i    \langle   m,  x  \rangle}.
\end{equation}
\item
\begin{equation}   \label{practical Poisson summation, second version}
\sum_{n\in\mathbb{Z}^{d}}f\left(  n+x\right)  
=
\lim_{\varepsilon \rightarrow 0}    
 \sum_{m\in\Z^d}             \widehat{\varphi}\left(
\varepsilon m\right)  \widehat{f}\left(  m\right)  
e^{2\pi i    \langle   m,  x  \rangle}.
\end{equation}
\end{enumerate}
\end{thm}
Because both $f$ and $\varphi_\varepsilon$ are compactly supported,
 the left-hand-sides of equations 
\eqref{practical Poisson summation, first version} and 
\eqref{practical Poisson summation, second version} 
are finite sums.  
\hfill $\square$

For a detailed proof of Theorem \ref{PracticalPoisson}, see \cite{BrandoliniColzaniTravagliniRobins1}.

An interesting aspect of this version of Poisson summation is that it can sometimes even apply to functions $f$
that are only piecewise continuous on $\R^d$, as long as \eqref{hypothesis, practical Poisson} holds.  Our prime example is of course
\[
f(x):= 1_\P(x),
\]
the indicator function of a polytope $\P$, and more generally $1_Q$ for a compact set $Q$ with reasonable behavior, such as a convex body. 
In Chapter \ref{Chapter.Minkowski}, we will use this version of Poisson summation, Theorem \ref{PracticalPoisson}, 
 to prove Theorem \ref{zero set of the FT of a polytope}.  

An interesting tool that gets used in the proof of Theorem \ref{PracticalPoisson} is a 
{\bf Plancherel-Polya type inequality}, as follows.
\begin{lem}   
Suppose that $ f \in L^1(\R^d), \hat f \in L^1(\R^d)$, and $f$ is compactly supported.  Then there exists a constant $c > 0$, depending on the support of $f$, such that
\begin{equation}\label{Plancherel-Polya inequality}
    \sum_{n \in \Z^d} 
   |  \hat f(n)  |  \leq c   \int_{\R^d}    | \hat f(\xi)  | d \xi.
\end{equation}
\end{lem}
\begin{proof} 
We define a new function $\psi$, which is infinitely smooth, and compactly supported, with $\psi(x) = 1$
for all $x$ in the support of $f$.  So we have $f(x) = \psi(x) f(x), \forall x \in \R^d$, and 
therefore $\hat f(\xi) =(\hat\psi * \hat f)(\xi)$ (using $\hat f \in   L^1(\R^d)$).  Because $\psi$ is smooth, we know that $\hat \psi$ is rapidly decreasing (by Corollary \ref{cor: f smoother implies FT of F decays faster}), 
and we have
\begin{align}
    \sum_{n \in \Z^d}    |  \hat f(n)  |  &=  
    \sum_{n \in \Z^d}  \left |   \int_{\R^d}  \hat \psi( n - \xi) \hat f(\xi) d\xi \right |  \\
    &\leq  \sum_{n \in \Z^d}   \int_{\R^d}  \left |   \hat \psi( n - \xi) \hat f(\xi) \right | d\xi   \\
    &= \int_{\R^d}  \sum_{n \in \Z^d} \left | \hat \psi( n - \xi)    \right |   \left |    \hat f(\xi)   \right |  d\xi \\
  & \leq  \sup_{\xi \in \R^d}   \left( \sum_{n \in \Z^d} \left |  \hat \psi( n - \xi)  \right |  \right)
      \int_{\R^d}  |  \hat f(\xi) |  d\xi   \\
      &\leq c  \int_{\R^d}  |  \hat f(\xi) |  d\xi.
\end{align}
The constant $c$ depends on $\psi$, and hence on the support of $f$.  To justify the last step, we note that 
$g(\xi):=  \sum_{n \in \Z^d} |   \hat \psi( n - \xi)  |$
 is a periodic function of $\xi$, with the unit cube $[0, 1]^d$ being a fundamental domain, so it suffices to show that $g$ is bounded on the unit cube.  
But due to the rapid decay of $\hat \psi$,  we may apply the Weierstrass $M$-test to conclude that 
the series $g$ is a uniformly convergent sum of continuous functions;   hence $g$ is itself a continuous function on a compact set (the cube), and  in fact achieves its maximum there. 
\end{proof}

The reader may consult \cite{SchmeisserSickel2000}, for example, for more information about related Plancherel-Polya type inequalities. 
In general, there are many functions $f \in L^1(\R)$ such that   $\sum_{n \in \Z}    |  \hat f(n)  | $ diverges, yet 
 $ \int_{\R}    | \hat f(\xi)  | d \xi$ converges, so that \eqref{Plancherel-Polya inequality} is false for these functions 
 (Exercise \ref{exercise:Plancherel-Polya type inequalities}).


\bigskip
\section{The Fourier transform of the ball} 
Whenever considering packing or tiling by a convex body $B$, we have repeatedly seen that taking the Fourier transform of the body, namely
$\hat 1_B$,  is very natural, especially from the perspective of Poisson summation.
It's also very natural to consider the FT of a ball in $\R^d$. 

To compute the Fourier transform of $1_{B(r)}$,  a very classical computation, we first define 
the \textbf{Bessel function}  \index{Bessel function}
$J_p$ of order $p$ (\cite{EpsteinBook}, page 147), which comes up naturally here:     
\begin{equation}\label{Bessel definition}
J_p(x) :=\left(\frac{x}{2} \right)^p \frac{1}{\Gamma\left(p + \frac12\right)\sqrt{\pi}}\int_0^\pi  e^{ix\cos\varphi} \sin^{2p}(\varphi) \, d\varphi,
\end{equation}
valid for $p>-\frac{1}{2}$, and all $x \in \R$.  When $p = n + \frac{1}{2}$ with $n \in \Z$, there are also the following relations with elementary trigonometric functions:   \index{Bessel function}
\begin{equation} 
J_{n + \frac{1}{2}}(x) = (-1)^n  \sqrt{   \frac{2}{\pi }    }
x^{n + \frac{1}{2}}
\left( \frac{1}{x}  \frac{d}{dx} \right)^n
\left( \frac{ \sin x }{x} \right).
\end{equation}
For example:
\begin{equation}
J_{\frac{1}{2}}(x) = \sqrt{ \frac{2}{\pi x} } \sin x,  \quad \text{ and }
J_{\frac{3}{2}}(x) =  \sqrt{ \frac{2}{\pi x} } \left( \frac{ \sin x}{x} - \cos x \right).
\end{equation}

We call a function $f:\R^d \rightarrow \C$  {\bf radial}   \index{radial function}
if it is invariant under all rotations of $\R^d$.  In other words, we have the definition
\[
f \text{ is radial } \iff  f \circ M = f,
\]
for all $M \in SO_d(\R)$, the orthogonal group.   Another way of describing a radial function is to say that the function $f$ is constant on each sphere that is centered at the origin, so that a radial function only depends on the norm of its input:  $f(x) = f(\|x\|)$, for all $x\in \R^d$.

A very useful fact in various applications of
Fourier analysis (in particular medical imaging)  is that the Fourier transform of a radial function is again a radial function (Exercise \ref{radial function transform}).

\bigskip
\begin{lem} \label{FT of the ball}
The Fourier transform of $B_d(r)$, the ball of radius $r$ in $\R^d$ centered at the origin, is
\[\hat{1}_{B_d(r)}(\xi) := \int_{B_d(r)} e^{-2\pi i     \langle \xi, x \rangle      } dx 
= \left(\frac{r}{\| \xi\|}\right)^{\frac{d}{2}}J_{\frac{d}{2}}\big(2\pi r\|\xi\|\big).\].   \index{Bessel function}
\end{lem}
\begin{proof}
Taking advantage of the inherent rotational symmetry of the ball, and also using the fact that the Fourier transform 
of a radial function is again radial (Exercise \ref{radial function transform}), we have:
\[
\hat{1}_{B_d(r)}(\xi) = \hat{1}_{B_d(r)}(0, \dots, 0, \|\xi\|),
\]
for all $\xi \in \R^d$.   With $r=1$ for the moment, we therefore have:
\[
\hat{1}_{B}(\xi) =
\int_{\| x\| \leq 1} e^{-2\pi i x_d\|\xi\|}\, dx_1\, \dotsc\, dx_d,
\]
Now we note that for each fixed $x_d$, the function being integrated is constant and the integration domain for the
variables $x_1, \dots, x_{d-1}$ is a $(d-1)$-dimensional ball of radius $(1-x_d^2)^{1/2}$. 
By equation  \eqref{volume of ball and sphere}, the
volume of this ball is $(1-x_d^2)^{\frac{d-1}{2}} \frac{     \pi^{\frac{d-1}{2}}    }{   \Gamma\left(\frac{d+1}{2} \right)   }$, we have
\[
\hat{1}_{B}(\xi)
= \frac{   \pi^{\frac{d-1}{2}}   }{    \Gamma(\frac{d+1}{2})    } 
\int_{-1}^1 
e^{-2\pi i x_d\|\xi\|}
(1-x_d^2)^{\frac{d-1}{2}}   \,dx_d
= \frac{ \pi^{  \frac{d}{2}         }}{\sqrt{\pi}\Gamma\left(\frac{d+1}{2}\right)} 
\int_0^\pi 
 \, e^{2\pi i\|\xi\|\cos\varphi}
 \sin^d\varphi  \,d\varphi.
\]
Using the definition~\eqref{Bessel definition} of the $J$-Bessel function, we get 
\[\hat{1}_{B}(\xi) = \| \xi\|^{-\frac{d}{2}  }J_{\frac{d}{2}  }\big(2\pi\|\xi\|\big),\]
and consequently 
\[\hat{1}_{B_d(r)}(\xi) = \left(\frac{r}{\| \xi\|}\right)^{\frac{d}{2}  }J_{\frac{d}{2}  }\big(2\pi r\|\xi\|\big).
\qedhere\]
\end{proof}

\bigskip
\begin{example}\label{ex:integral using Bessel functions}
\rm{
Using the $J$-Bessel functions, let's work out an explicit evaluation of the following interesting
 integrals, for all $p>0$:
\begin{equation} \label{identity of J-Bessel example}
\int_0^\pi  \sin^{2p}(\varphi) \, d\varphi = 
\sqrt \pi  \frac{    \Gamma\left(p + \frac12\right)    }{       \Gamma\left(p + 1  \right)}.
\end{equation}
Whenever we raise a negative real number to an arbitrary real exponent, some care has to be taken to avoid `branch problems' with the definition of exponentiation.  But over the latter domain of integration, we are considering the nonnegative function $\sin(\varphi)\geq 0$, so everything is copacetic.   We will use the following equivalent formulation for the $J_p$ Bessel function in terms of a hypergeometric series:
\begin{equation} \label{Bessel infinite series}
J_p(x) = \frac{x^p}{2^p} \sum_{k=0}^\infty (-1)^k \frac{x^{2k} }{2^{2k} k! \, \Gamma(p+k+1) }.
\end{equation}
(\cite{EpsteinBook}, p. 684).
  Using the definition of the Bessel function \eqref{Bessel definition}, we can rewrite it slightly:
\begin{equation}\label{fancy integral 1} \index{Bessel function}
\frac{J_p(x)}{x^p}   2^p \sqrt{\pi}   \,  \Gamma\left(p + \frac12\right)
=
\int_0^\pi  e^{ix\cos\varphi} \sin^{2p}(\varphi) \, d\varphi.
\end{equation}
Taking the limit as $x\rightarrow 0$, we can safely move this limit inside the integral in \eqref{fancy integral 1}
because we are integrating a differentiable function over a compact interval:
\begin{align*}
\int_0^\pi  \sin^{2p}(\varphi) \, d\varphi 
&= \lim_{x\rightarrow 0}
\frac{J_p(x)}{x^p}   2^p \sqrt{\pi}  \,  \Gamma\left(p + \frac12\right).
\end{align*}
So if we knew the asymptotic limit $ \lim_{x\rightarrow 0} \frac{J_p(x)}{x^p}$, we'd be in business. 
From \eqref{Bessel infinite series}, we may divide both sides by $x^p$, and then take the limit as $x\rightarrow 0$ to obtain the constant term of the remaining series, giving us
\[
 \lim_{x\rightarrow 0} \frac{J_p(x)}{x^p} =  \frac{1}{2^p \Gamma(p+1)}.
 \]
Altogether, we have
\begin{align*}
\int_0^\pi  \sin^{2p}(\varphi) \, d\varphi 
&= \lim_{x\rightarrow 0}
\frac{J_p(x)}{x^p}   2^p \sqrt{\pi}  \,  \Gamma\left(p + \frac12\right) \\
&= \frac{1}{2^p \Gamma(p+1)}     2^p \sqrt{\pi}  \,  \Gamma\left(p + \frac12\right) \\
&=  \sqrt \pi  \frac{    \Gamma\left(p + \frac12\right)    }{       \Gamma\left(p + 1  \right)},
\end{align*}
valid for all $p>0$.
}
In the special case that $p$ is a positive integer, the latter identity can of course be written in terms of a ratio of factorials (Exercise \ref{the integral of the example in terms of factorials}). 
\hfill $\square$
\end{example}


\bigskip


\bigskip
\section{Uncertainty principles}    
\label{section:uncertainy principles}

\begin{quote}
Uncertainty is the only certainty there is, and knowing how to live with insecurity is the only security. 

-- John Allen Paulos
\end{quote}

Perhaps the most basic type of an \emph{uncertainy principle} is the fact that if a function $f$ is compactly supported,
then its Fourier transform $\hat f$ cannot be compactly supported - Theorem \ref{basic uncertainty principle} below.   
Similar impossible constraints, placed simultaneously  on both $f$ and $\hat f$, 
have become known as {\bf uncertainty principles}.  Perhaps the most famous of these, originating in quantum mechanics, is Heisenberg's 
discovery, as follows. 

\begin{thm}[Heisenberg uncertainty principle]
 \index{uncertainty principle, Heisenberg}
Let $f\in L^2(\R^d)$, with the normalization assumption that $\int_{\R^d} |f(x)|^2 dx=1$.   Then:
\begin{equation}
\int_{\R^d} \|x\|^2 |f(x)|^2 dx  \int_{\R^d} \|x\|^2 |\hat f(x)|^2 dx  \geq \frac{1}{16 \pi^2},
\end{equation}
with equality holding if and only if $f$ is equal to a Gaussian. 

\rm{(For a proof see \cite{OsgoodBook}, or \cite{DymMcKean}.   )}
\hfill $\square$
\end{thm}

\bigskip
\begin{thm}[Hardy uncertainty principle] 
\label{Hardy uncertainty principle}  \index{uncertainty principle, Hardy}
Let $f\in L^1(\R^d)$ be a function that enjoys the property that 
\begin{equation*}
   |f(x)|   \leq A e^{-\pi c x^2 }    \text{  and    }   \      |\hat f(\xi) |    \leq  B e^{-\pi \xi^2/c}, 
\end{equation*}
 for all $x, \xi\in \R^d$,  and for some constants $A, B, c >0$.      
 
 Then $f(x)$ is a scalar multiple of the Gaussian $e^{-\pi c x^2}$.

\rm{(For a proof see  \cite{Hardy.uncertainty})}
 \hfill $\square$
\end{thm}

\bigskip
But perhaps the most ``elementary'' kind of uncertainty principle is the following basic fact, which is useful to keep in mind.
\begin{thm}   \label{basic uncertainty principle} 
Let $f\in L^1(\R^d)$ be a function that is supported on a compact set in $\R^d$. 
Then $\hat f$ is not supported on any compact set in $\R^d$.

\rm{(For a proof see  \cite{EpsteinBook})}
 \hfill $\square$
\end{thm}


\section*{Notes}
\begin{enumerate}[(a)]
\item  \label{Fourier books}
There are some wonderful introductory books that develop Fourier analysis from first principles, such as the books by
Stein and Shakarchi \cite{SteinShakarchi} and Giancarlo Travaglini \cite{Travaglini}.
The reader is also encouraged to read more advanced but fundamental introductions to Fourier analysis, in particular  the  books by Mark Pinsky \cite{MarkPinsky},  Edward Charles Titchmarsh \cite{Titchmarsh},  Antoni Zygmund  \cite{Zygmund}, 
 Einsiedler and Ward \cite{EinsiedlerWardBook}, Dym and McKean \cite{DymMcKean},
and of course the classic: Stein and Weiss \cite{SteinWeiss}.   In addition, the book \cite{Terras} by Audrey Terras is a good introduction to Fourier analysis 
on finite groups, with applications.   
Another excellent introduction to Fourier analysis, which is more informal and focuses 
on various applications, is given by Brad Osgood \cite{OsgoodBook}.

\item  There are some ``elementary''  techniques that we will use, from the calculus of a complex variable, but which require essentially no previous knowledge in this field.   In particular, suppose we have two analytic functions
 $f:\C \rightarrow \C$ and $g:\C \rightarrow \C$, such that $f(z_k) = g(z_k)$ for a convergent sequence of complex numbers $z_k \rightarrow L$, where $L$ is any fixed complex number.    Then $f(z) = g(z)$ for all 
 $z \in \C$.  
 
 The same conclusion is true even if the hypothesis is relaxed to the assumption that
 both $f$ and $g$ are meromorphic functions, as long as the sequence and its limit stay away from the poles of $f$ and $g$.
 
\item  \label{Dirac delta} 
The ``Dirac delta function"  is part of the theory of ``generalized functions'' and may be intuitively defined by the full sequence
of Gaussians $G_t (x) :=   t^{-\frac{d}{2}}    e^{ -\frac{\pi}{t}  || x ||^2 }$, taken over all $t>0$.   The observation that there is no identity for the convolution product 
on $\R^d$ is a clear motivation for a theory of generalized functions, beginning with the Dirac delta function.
Another intuitive way of ``defining''  the Dirac delta function is: 
\[
\delta_0(x) := 
 \begin{cases}  
\infty     &      \mbox{if }  x=0 \\ 
0  &        \mbox{if  not},
\end{cases}  
\]
even though this is not a function.  But in the sense of distributions (i.e. generalized functions), we have
$\lim_{\rightarrow 0} G_t(x) = \delta_0(x)$.

More rigorously, the $\delta$-function belongs to a theory of distributions  that was developed  by Laurent Schwartz
\index{Schwartz, Laurent}
 in the 1950's and by S.L. Sobolev in 1936, where we can think of generalized functions as linear functionals on the space of all bump functions on $\R^d$ (see the book by Lighthill \cite{Lighthill} for a nice introduction to generalized functions).  

Such generalized functions were originally used by the Physicist Paul Dirac in 1920, before the rigorous mathematical theory was even created for it,  in order to better understand quantum mechanics.  \index{Dirac, Paul}

\item  \label{Note:GregKuperberg} I'd like to thank Greg Kuperberg for very helpful comments, and in particular for introducing 
me to the statement of Theorem \ref{Euler-Maclaurin type identity}, for which we still cannot find a published reference.

\item
It is sometimes interesting to derive analogues between norms in $\R^d$ and norms in an infinite dimensional function space.  
Among the many norm relations in $\R^d$, we mention one elementary but interesting relation:
\[
\| x \|_1 \leq \sqrt{n} \  \| x \|_2, 
\]
for all vectors $ x \in \R^d$, where $\|x\|_1:= |x_1|+\cdots + |x_d|$, and $\|x\|_2:= \sqrt{x_1^2+\cdots + x_d^2}$.
 (see Exercise \ref{elementary norm relations} for more practice with related norm relations).
 At this point the curious reader might wonder  ``are there any other inner products on $\R^d$,   besides the usual inner product 
 $\langle x, y \rangle:= \sum_{k=1}^d x_k y_k$?"   A classification of all inner products that exist on $\R^d$ is given in Exercise \ref{All norms on Euclidean space}.

\item 
Of great practical importance, and historical significance, a  {\bf bump function} is defined as any infinitely smooth function 
on $\R^d$, which is compactly supported.  In other words, a bump function enjoys the following properties:
    \begin{itemize}
       \item   $\phi$ has compact support on $\R^d$.
       \item   $\phi \in C^\infty(\R^d)$.
    \end{itemize}
Bump functions are also called {\bf test functions}, and if we consider the set of all bump functions
on $\R^d$, under addition, we get a vector space $V$, whose dual vector space is called the space of {\bf distributions on $\R^d$}.

\item   Theorem \ref{nice2}, originally appearing in Poisson's work, also appear in Stein and Weiss' book \cite{SteinWeiss};  here we 
gave a slightly different exposition.

\item The cotangent function, appearing in some of the exercises below,  is the unique {\it meromorphic function} that has a simple pole at every integer, with residue 1 (up to multiplication by an entire function with the same residues).   The cotangent function also forms an entry point for Eisenstein series in number theory, through the corresponding partial fraction expansion of its derivatives. 

\item  A deeper exploration into projections and sections of the unit cube in $\R^d$ can be found in
 ``The cube - a window to convex and discrete geometry'', by Chuangming~Zong \cite{Zong.book}. 
In \cite{KoldobskyBook}, Alexander Koldobsky gives  a thorough introduction to sections of convex bodies, intersection bodies, and the Busemann-Petty problem. 

\item  There are numerous other identities throughout mathematics that are equivalent to special cases of Poisson summation, such as the Euler-MacLaurin summation formula, the Abel-Plana 
 formula, and the Approximate sampling formula of signal analysis
(see \cite{Butzer.etal} for a nice treatment of such equivalences for functions of $1$ real variable, and functions of $1$ complex variable).

\item There is an important and fascinating result of Cordoba \cite{Cordoba1} which states the following.  Let 
$A:= \{ x_k\}_{k\in \Z},  B:= \{ y_k\}_{k\in \Z}$ be two discrete sets in $\R^d$.  Suppose that for all 
Schwartz functions $f$, we have 
\[
\sum_{k\in \Z} f(x_k) = \sum_{k\in \Z} \hat f(y_k).
\]
 Then both of the sequences $A$ and $B$ must be lattices in $\R^d$, and $B= A'$, its dual lattice.

\item Finally, it's worth mentioning that the term `Harmonic analysis'   is simply a more general theory, extending the notion of Fourier analysis to other groups, besides Euclidean space.  
\end{enumerate}

\bigskip

\section*{Exercises}
\addcontentsline{toc}{section}{Exercises}
\markright{Exercises}


\begin{quote}                         
``In theory, there is no difference between theory and practice; 

 \ \ but in practice, there is!''  \ \ \ \ -- Walter J. Savitch  \\
 \end{quote}

\medskip
\begin{prob} $\clubsuit$  \label{elementary norm relations}
\rm{
On $\R^d$ the $L^2$-norm is defined by 
$\|x\|_2:= \sqrt{ x_1^2 + \cdots + x_d^2}$,  the $L^1$-norm is defined by
$\|x\|_1:=  |x_1| + \cdots + |x_d|$, and the $L^\infty$-norm is defined by 
$\|x\|_\infty:= \max\{ |x_1|, \dots, |x_d| \}$. 

Prove the following four norm relations:
\[
\|x\|_\infty        \leq       \| x \|_2      \leq            \| x \|_1   \leq     \sqrt{d}  \,   \| x \|_2    \leq   d \,  \|x\|_\infty,  
\]    
for all $x \in \R^d$.
}
\end{prob}


\medskip
\begin{prob}  \label{exercise:hyperbolic cosine and sine} 
\rm{
We know that the functions $u(t) := \cos t = \frac{e^{it} + e^{-it}}{2}$ and $v(t) := \sin t = 
 \frac{e^{it} - e^{-it}}{2i}$ are natural, partly because they parametrize the unit circle: $u^2 + v^2 = 1$.   Here we see that there are other similarly natural functions, parametrizing the hyperbola.

\begin{enumerate}[(a)]
\item  Show that the following functions parametrize the hyperbola $u^2 - v^2 = 1$:
\[
u(t) := \frac{e^t + e^{-t}}{2}, \ \ \ v(t) := \frac{e^t - e^{-t}}{2}.
\]
(This is the reason that the function $\cosh t:= \frac{e^t + e^{-t}}{2}$ is called the hyperbolic cosine, and the function
$\sinh t := \frac{e^t - e^{-t}}{2}$ is called the hyperbolic sine)
\item  The hyperbolic cotangent is defined as $\coth t := \frac{ \cosh t }{ \sinh t} = 
\frac{ e^t + e^{-t}}{e^t - e^{-t}}$.   Using Bernoulli numbers, show that $t \coth t$ has the  Taylor series:
\[
t \coth t = \sum_{n=0}^\infty  \frac{2^{2n}}{(2n)!} B_{2n} t^{2n}.
\]
\end{enumerate}
}
\end{prob}

\medskip  
\begin{prob}  $\clubsuit$  
\label{compute FT for exponential of abs value}
\rm{
Prove that:
\[
\frac{t}{\pi} \sum_{n \in \Z} \frac{1}{n^2 + t^2} = \sum_{m \in \Z} e^{-2\pi t |m|}.
\]  

Hint.  \ Think of Poisson summation, applied to the function $f(x):= e^{-2\pi t |x|}$.
}
\end{prob}

\medskip
\begin{prob} $\clubsuit$  \label{Riemann zeta function, and Bernoulli numbers}
\rm{
Here we evaluate the Riemann zeta function at the positive even integers.
\begin{enumerate}[(a)]
\item  Show that 
\[
 \sum_{n\in \Z}      e^{-2 \pi t |n|}   =
 \frac{    1 + e^{-2\pi t}    }{     1-e^{-2\pi t}     } := \coth(\pi t),
 \]
for all $t>0$.  
 
 \bigskip

\item Show that the cotangent function has the following well-known partial fraction expansion:
\[
\pi \cot(\pi x) = \frac{1}{x} + 2x \sum_{n=1}^\infty \frac{1}{ x^2 - n^2},
\]
valid for any $x \in \R - \Z$.

\item  Let $0 < t < 1$.    Show that 
\[
\frac{t}{\pi} \sum_{n\in \Z}  \frac{1}{n^2 + t^2} =  \frac{1}{\pi t}  +
 \frac{2}{\pi} \sum_{m=1}^\infty   (-1)^{m+1}  \zeta(2m) \  t^{2m-1},
\]
where $\zeta(s):= \sum_{n=1}^\infty \frac{1}{n^s}$ is the Riemann zeta function, initially defined by the latter series, which is valid for  all $s \in \C$ with $Re(s) >1$.

\item     Here we show that we may quickly evaluate the Riemann zeta function at all even integers, as follows. We recall the definition of the Bernoulli numbers, namely:
\[  
\frac{z}{e^z - 1}  =  1 - \frac{z}{2}  + \sum_{m \geq 1}  \frac{B_{2m}}{2m!}  z^{2m}.
\]

Prove that for all $m \geq 1$,
\[
 \zeta(2m) = \frac{(-1)^{m+1} }{2} \frac{ (2\pi)^{2m}}{  (2m)!} B_{2m}.
\]
Thus, for example, using the first $3$ Bernoulli numbers, we have: $\zeta(2) = \frac{\pi^2}{6}$, $\zeta(4) = \frac{\pi^4}{90}$, and $\zeta(6) = \frac{\pi^6}{945}$.
\end{enumerate}
}
\end{prob}


\medskip
\begin{prob}   \label{Chebyshev polys}
\rm{
 For each $n\geq 1$, let $T_n(x) = \cos(nx)$.
For example, $T_2(x) = \cos(2x) =2 \cos^2(x) - 1$, so $T_2(x) = 2u^2 -1$, a polynomial in $u:= \cos x$.
\begin{enumerate}[(a)]
\item     Show that for all $n\geq 1$, $T_n(x)$ is a polynomial in $\cos x$.  
\item    Can you write $x^n + \frac{1}{x^n}$ as a polynomial in the variable $x + \frac{1}{x}$?  
Would your answer be related to the polynomial $T_n(x)$?  What's the relationship in general? 
For example, $x^2 + \frac{1}{x^2} = \Big( x + \frac{1}{x}\Big)^2 - 2$.   
\end{enumerate}
}
\end{prob}
 
Notes.  The polynomials $T_n(x)$ are very important in applied fields such as approximation theory, and optimization,
 because they have many useful extremal properties.  
They are called Chebyshev polynomials. \index{Chebyshev polynomials}

\medskip
\begin{prob} \label{sec - its own Fourier transform}
\rm{
The hyperbolic secant is defined by 
\[
{\rm sech}(\pi x) := \frac{2}{e^{\pi x} + e^{-\pi x}}, \text{ for }  x \in \R.
\]
\begin{enumerate}[(a)]
\item \label{eigenfunction of FT}
Show that  ${\rm sech}(\pi x)$  is its own Fourier transform:
\[
\F({\rm sech})(\xi) = {\rm sech}(\xi), 
\]
for all $\xi \in \R$.
\item  \label{bounded above by Gaussian}
Show that  ${\rm sech}(\pi x)$ can never be bounded above by any Gaussian, in the precise sense that the following claim  is {\bf impossible}:  there exists a constant  $c>0$ such that for all $x \in \R$ we have:
\[
{\rm sech}(\pi x) \leq e^{-cx^2}.
\]
\end{enumerate}

Notes.  For part  \ref{eigenfunction of FT}, the easiest path is through the use of
basic complex analysis (but the reader might find a different path!). For part \ref{bounded above by Gaussian}, it may be helpful to look at Hardy's uncertainty principle, Theorem \ref{Hardy uncertainty principle}.   We can also conclude from Hardy's uncertainty principle that any eigenfunction $f$ of the Fourier transform cannot be bounded above by a Gaussian, aside from the case that $f$ is itself a Gaussian.
}
\end{prob}

\medskip
\begin{prob} 
\rm{
Using the previous exercise, conclude that 
\[
\int_{\R}  \frac{1}{e^{\pi x} + e^{-\pi x}}  dx = \frac{1}{2}.
\] 
}
\end{prob}


\medskip
\begin{prob}    $\clubsuit$
\rm{
Prove that
\[
\int_0^1 \left(\{ax\} - \frac{1}{2} \right)\left( \{bx\} - \frac{1}{2} \right)dx 
= \frac{ \rm{gcd}^2(a, b) }{12ab},
\]
for all positive integers $a, b$.  As always, $\{x\}$ is the fractional part of $x$.

Notes.  This integral is called a {\bf Franel integral}, and there is a substantial literature about related integrals.
In 1924, J\'er\^ome Franel related this integral to the Riemann hypothesis, and to Farey fractions. 
}
\end{prob}

\medskip
\begin{prob} 
\rm{
Given an even function $f \in L^1(R^d)$,  such that $ \hat f \in L^1(R^d)$ as well,
we clearly always have $f(x) :=  \tfrac{1}{2} \left(  f(x) + \hat f(x)  \right)    +   \tfrac{1}{2} \left(  f(x) - \hat f(x) \right)$.  Show that:
\begin{enumerate}[(a)]  
\item  The function $g(x):= \tfrac{1}{2} \left(  f(x) + \hat f(x)  \right)$ is an eigenfunction of the Fourier transform acting on $L^1(\R^d)$, with 
eigenvalue equal to $+1$.  
\item  Similarly, the function  
$h(x):=  \tfrac{1}{2} \left(  f(x) - \hat f(x) \right)$ is an eigenfunction of the Fourier transform acting on $L^1(\R^d)$, with 
eigenvalue equal to $-1$. 
\end{enumerate} 
}
\end{prob}

\medskip
\begin{prob}    $\clubsuit$ \label{Schwartz space convolution invariance}
\rm{
Let $f: \R \rightarrow \C$ belong to the Schwarz class of functions on $\R$, denoted by $S(\R)$. 
Show that $\hat f \in S(\R)$ as well.
}
\end{prob}

\medskip
\begin{prob} 
\rm{
Here we define $f(x):=  1+\sin (2 \pi x)$, for all $x \in \R$.  We note that $f$ is a periodic function of $x\in \R$, with period $1$, so it may be considered as a function on the torus $\mathbb T$.
\begin{enumerate}[(a)]  
\item Using Theorem \ref{Euler-Maclaurin type identity}, find the little-o asymptotics (with $N\rightarrow \infty$) for the finite sum defined above.
\item \label{second part of example for our Euler-Maclaurin-type theorem}
Show directly that $ \sum_{m=0}^{N-1} f(\frac{m}{2N}) 
= N+\frac{N}{2^{N-1}}$. 
\item From part \ref{second part of example for our Euler-Maclaurin-type theorem}, conclude
 (independently of Theorem  \ref{Euler-Maclaurin type identity}) that we get the same little-o asymptotics that 
  Theorem  \ref{Euler-Maclaurin type identity} predicts.

Notes.  For part \ref{second part of example for our Euler-Maclaurin-type theorem},  you might
 begin with the polynomial identity 
 \[
 1+ z+ z^2 + \cdots + z^{N-1} = \prod_{r=1}^{N-1}   (z - e^{\frac{2\pi i r}{N}}   ).
 \]
\end{enumerate}
}
\end{prob}

\medskip
\begin{prob}   $\clubsuit$ \label{All norms on Euclidean space}
\rm{
Here we answer the very natural question ``What are the other inner products on $\R^d$, besides the usual inner product $\langle x, y \rangle:= \sum_{k=1}^d x_k y_k$ ?"  
Here we show that all inner products are related to each other via positive definite matrices, as follows.  
We recall from Linear Algebra that a symmetric 
matrix is called positive definite if all of its eigenvalues are positive. 
Prove that the following two conditions are equivalent:
\begin{enumerate}
 \item  $  \langle x, y   \rangle$  is an inner product on $\R^d$. 
\item   $  \langle x, y  \rangle  := x^T M y$, for some positive definite matrix $M$.
\end{enumerate}
}
\end{prob}

\bigskip
\begin{prob}   
For any positive real numbers $a < b <  c < d$, define 
\[
f(x) := 1_{[a, b]}(x) + 1_{[c, d]}(x).
\]
Can you find $a,b,c,d$ such that $\hat f(\xi)$ is nonzero for all $\xi \in \R$? 
\end{prob}

\medskip
\begin{prob}  $\clubsuit$
Show that for $f, \hat f \in L^1(\R^d)$, the only eigenvalues of the linear operator 
\[
f \rightarrow \hat f
\]
 are $\{ 1, -1, i, -i \}$, and find  functions in  $ L^1(\R^d)$ that achieve each of these eigenvalues.
\end{prob}

\medskip
\begin{prob}    $\clubsuit$  \label{going backwards in Poisson summation}
Show that the special case of Poisson summation, \ref{Poisson.summation3}, implies the general case, Theorem \ref{The Poisson Summation Formula, for lattices}.  \index{Poisson summation formula}
\end{prob}


\medskip
\begin{prob}     \label{Gaussian1}   $\clubsuit$
We recall the definition of the Gaussian:  for each fixed $\varepsilon >0$, and for all $x \in \R^d$, they are defined by 
\begin{equation}    \index{Gaussian}
G_{\varepsilon} (x) :=  \frac{1}{\varepsilon^{\frac{d}{2}}}  e^{ -\frac{\pi}{\varepsilon}  || x ||^2 }.
\end{equation}
Show that: 
  \[
\int_{\R^d}    G_{\varepsilon} (x) dx = 1.
\]
\end{prob}

\medskip
 \begin{prob}   \label{Gaussian2}  $\clubsuit$ (Hard-ish)
Show that, for all $m \in \R^d$, the Fourier transform of the Gaussian
$G_{\varepsilon}(x)$  is:
\[
\hat G_\varepsilon( m ) =  e^{ -\pi  \varepsilon || m ||^2 }.
\]
Conclude that for each fixed $n\in \R^d$
\[
\F\left(
\frac{1}{\varepsilon^{\frac{d}{2}}}  e^{ -\frac{\pi}{\varepsilon}  || x+n ||^2 }
\right)(\xi)
=
e^{ -\pi  \varepsilon || \xi ||^2 } e^{2\pi i \langle \xi, n \rangle}.
\]
\end{prob} 

\medskip
\begin{prob} 
We define the translation operator $T_h: L^2(\R^d) \rightarrow L^2(\R^d)$ by
$(T_h f)(x) := f( x - h)$, for any fixed $h>0$.   Show that convolution commutes with the translation operator, as follows:
\[
T_h(f*g) = (T_h f)*g = f*(T_h g). 
\]
\end{prob} 
Notes.   Using standard Linear algebra terminology, this is called translational equivariance.

\medskip
\begin{prob} \label{continuity of convolution for two L^1 bounded functions}
Prove that  if  $f, g \in L^1(\R^d)$ are bounded functions,  then $f*g$ is continuous on $\R^d$.
\end{prob}

\medskip
\begin{prob} $\clubsuit$
 \label{continuity of convolution for two L^2 functions}
\rm{
Suppose that we fix any $f, g \in L^2(\R^d)$.  Here we carry the reader through a detailed 
proof that
$f*g$ is always continuous on $\R^d$, by using a mixture of convergence in $L^2(\R^d)$ and pointwise convergence.
\begin{enumerate}[(a)]
\item  \label{first part of continuity, L^2 convolution}
For any sequence of functions $f_n \in L^2(\R^d)$ with the property that
          $\lim_{n\rightarrow \infty} f_n = f$ in $L^2(\R^d)$, show that 
    \begin{equation}
    \lim_{n \rightarrow \infty} (f_n * g)(x)= (f*g)(x), 
    \end{equation}
    for each $x \in \R^d$.
\item  \label{second part of continuity, L^2 convolution}
Define the translation operator $T_h: L^2(\R^d) \rightarrow L^2(\R^d)$ by
$(T_h f)(x) := f( x - h)$, for a fixed $h>0$.  Show that 
\[
\lim_{h \rightarrow 0} T_h f = f,  \ \text{ in } L^2(\R^d).
\]
\item \label{third part of continuity, L^2 convolution}
Show that   $\lim_{h\rightarrow 0} (T_h f * g)(x) = (f*g)(x)$, for each fixed $x \in \R^d$.
\item \label{final part of continuity, L^2 convolution}
Conclude that  $f*g$ is continuous on $\R^d$.
\end{enumerate}
}
\end{prob}

Notes.  It follows from either 
Exercise \ref{continuity of convolution for two L^1 bounded functions} or 
Exercise \ref{continuity of convolution for two L^2 functions}
that if $A, B \subset \R^d$ are convex bodies (or finite unions of convex bodies), then
$(1_A*1_B)(x)= \vol\left( A \cap (-B + x) \right)$ is a continuous function of $x \in  \R^d$.

\medskip
\begin{prob}
\rm{
Show by example that $f, g \in L^1(\R^d)$ does not necessarily imply that 
$f g \in L^1(\R^d)$ (here $f g $ is the usual product of functions).
}
\end{prob}

\medskip
 \begin{prob}   \label{Plancherel extended}  $\clubsuit$
\rm{
For all  $f, g \in L^2(\R^d)$, prove  that:
\begin{enumerate}[(a)]
\item 
\[
\langle f, g \rangle = \frac{1}{2} \left( \| f\|^2 + \| g \|^2 - \|f+g\|^2\right).
\]
\item  
\[
\langle f, g\rangle =  \langle \hat f, \hat g  \rangle.
\]
\end{enumerate}
}
\end{prob}

\medskip
 \begin{prob}   \label{total mass 1}  $\clubsuit$
\rm{
Given any approximate identity sequence $\phi_\varepsilon$, as defined in   \eqref{approximate identity},
 show that  for each~$\varepsilon~>~0$,
 \[
 \int_{\R^d}\phi_\varepsilon(x) dx  = 1.
 \]
 }
\end{prob}

\medskip
 \begin{prob}   \label{continuity and equality a.e. implies equality}  $\clubsuit$
\rm{
Let $E\subset \R^d$ be any set, and suppose we have two continuous functions 
$f, g:E \rightarrow \C$.

If  $f = g$ for almost every $x\in E$, prove that $f=g$ for all $x\in E$.    
}
\end{prob}

\medskip
 \begin{prob}  \label{strictly less than the FT at zero}
  $\clubsuit$
\rm{
Under some positivity assumptions for $f$, the Fourier transform of $f$ achieves its 
\emph{unique
 maximum} at the origin.  More precisely, we have the following.
\begin{enumerate}[(a)]
\item Suppose $f\in L^1(\R^d)$, and $f(x) >0$ for all $x\in \R^d$.  Prove that
\[
|\hat f(\xi)| < \hat f(0), 
\]
for all $\xi \in \R^d$.

\item Now let $\P\subset \R^d$ be a $d$-dimensional convex set.  Prove that
for all $\xi\not=0$, we have
\begin{equation}
\hat 1_\P(\xi) <  \vol \P.
\end{equation}
\end{enumerate}
}
\end{prob}

\medskip
 \begin{prob} 
\rm{
Show that the ramp function, defined in \eqref{def of ramp}, also has the representation:
\begin{equation}
r_0(x) = \frac{ x + |x| }{2},
\end{equation}
for all $x\in \R$.  

Notes. \ Some books, particularly in approximation theory, use the notation $r_0(x) := x_+$.
}
\end{prob}

\medskip
 \begin{prob}   
 \label{positive FT over R}   $\clubsuit$
\rm{
 Here we show how to construct compactly supported functions $f:\R\rightarrow \C$
whose Fourier transform is {\bf strictly positive} on all of $\R$.
Fix any two incommensurable real numbers $r, s$ (meaning that $\frac{r}{s} \notin \Q$),
and define
\[
f:= 1_{[-r, r]}*1_{[-r, r]}  + 1_{[-s, s]}*1_{[-s, s]},
\]
which is a sum of two hat functions, as depicted in Figure \ref{A sum of two hat functions}. 
Prove that for all $\xi \in \R$, we have $ \hat f(\xi) >0$.
}

\begin{figure}[htb]
\begin{center}
\includegraphics[totalheight=2in]{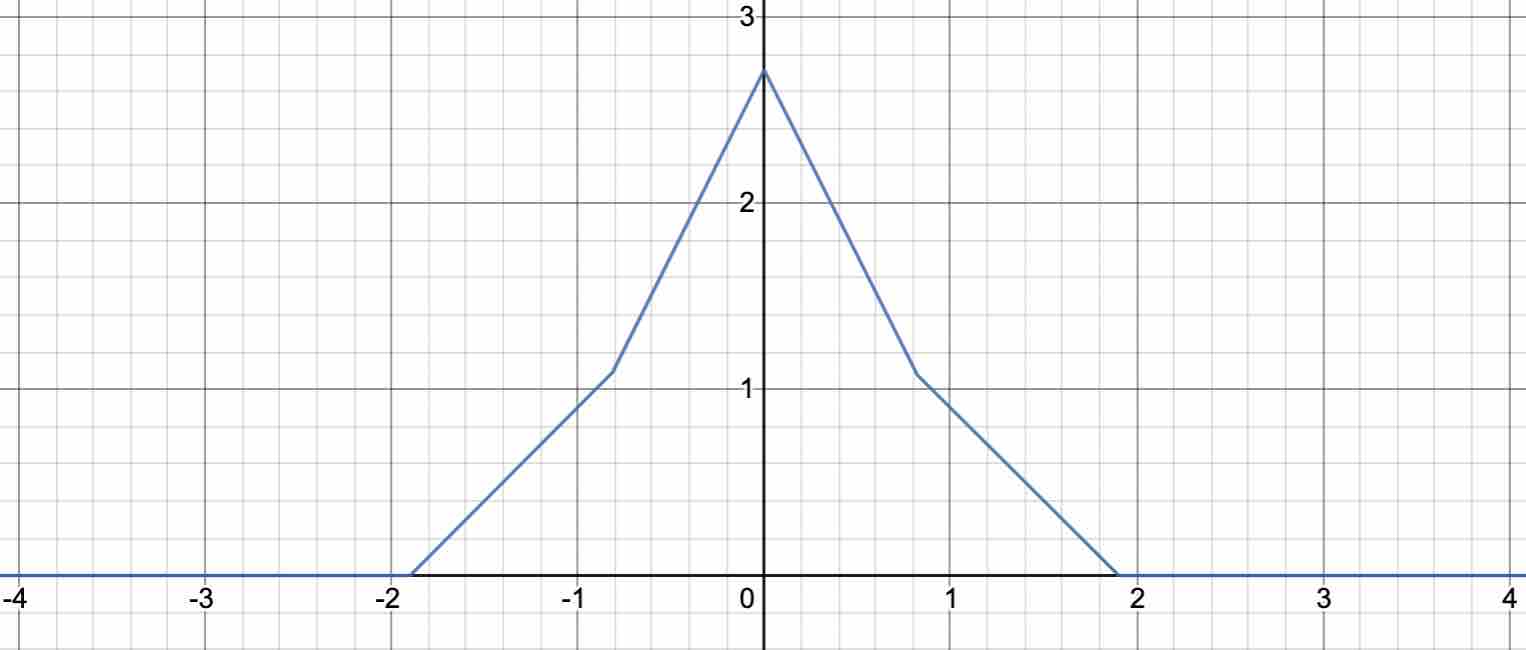}
\end{center}
\caption{The function $f$ of Exercise \ref{positive FT over R}, a sum of two hat functions, 
with $s = \sqrt{\frac{2}{3}}$, and  $r= 1.9$}   
 \label{A sum of two hat functions}
\end{figure}

\begin{figure}[htb]
\begin{center}
\includegraphics[totalheight=4.5in]{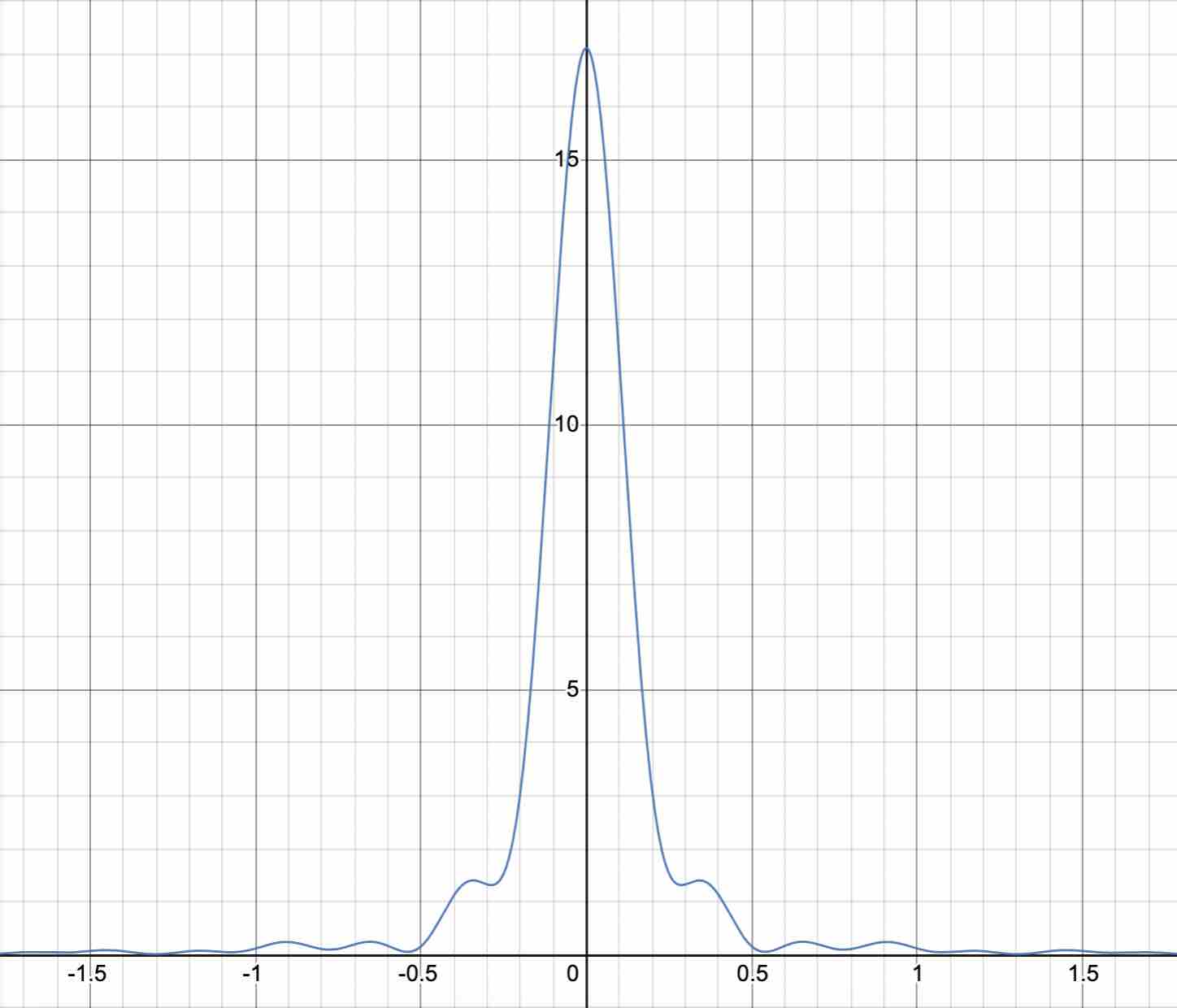}
\end{center}
\caption{The \emph{strictly positive} Fourier transform $\hat f(\xi)$ of Exercise \ref{positive FT over R}, with the two incommensurable numbers $s = \sqrt{\frac{2}{3}}$, and  $r= 1.9$}   
 \label{strictly positive FT on R}
\end{figure}

Notes.   This construction can be extended to higher dimensions, once we know more about the Fourier transforms of balls in $\R^d$ - see Exercise \ref{positive FT over R^d}. 
\end{prob}


\medskip
 \begin{prob}  
 \label{Heaviside and ramp}
   $\clubsuit$
Show that for all $a, b \in \R$, we have:
\[
H_a*H_b = r_{a+b},
\] 
where $H_a$ is the heaviside function of \eqref{def of heaviside}, and $r_a$ is the ramp function of \eqref{def of ramp}.
\end{prob}

\medskip
 \begin{prob}  
 \label{exercise:Plancherel-Polya type inequalities}    $\clubsuit$
\rm{
 Here we show that the absolute convergence of a series, and the absolute convergence of the corresponding integral, are independent of each other.
 \begin{enumerate} [(a)]
 \item
Find a function $f :\R \rightarrow \C$ such that   $\sum_{n \in \Z}    |  \hat f(n)  | $ diverges, yet 
 $ \int_{\R}    | \hat f(\xi)  | d \xi$ converges.   
  \item  On the other hand,  find a function $f :\R \rightarrow \C$ such that   
  $\sum_{n \in \Z}    |  \hat f(n)  | $  converges, yet 
  $ \int_{\R}    | \hat f(\xi)  | d \xi$ diverges.   
 \end{enumerate}
 }
 \end{prob}
 Notes.   Exercise \ref{exercise:Plancherel-Polya type inequalities} 
 shows that there the Plancherel-Polya inequality holds only for a special class of functions.

 \medskip
\begin{prob} 
We recall that $C(\R^d)$ is the function space consisting of all continuous functions on $\R^d$.  Show that
\begin{enumerate}[(a)]
\item $C(\R^d)\cap L^1(\R^d) \not\subset C(\R^d)\cap L^2(\R^d)$.
\item $C(\R^d)\cap L^2(\R^d) \not\subset C(\R^d)\cap L^1(\R^d)$.
\end{enumerate}
\end{prob}

\medskip
\begin{prob} \label{tricky application of Poisson summation}
\rm{
Here is a slightly different version of Poisson summation, which is easy to prove.   If $g:\R^d\rightarrow \C$ is infinitely smooth, and compactly supported, prove that 
\[
\sum_{n \in \Z^d}  \hat g(n)  = \sum_{n \in \Z^d}  g(n),
\]
and of course the right-hand side is a finite sum.
}
\end{prob}


\chapter{
\blue{Classical geometry of numbers \\
Part I: \, Minkowski meets Siegel}  
}
\label{Chapter.Minkowski}
\begin{wrapfigure}{R}{0.49\textwidth}
\centering
\includegraphics[width=0.28\textwidth]{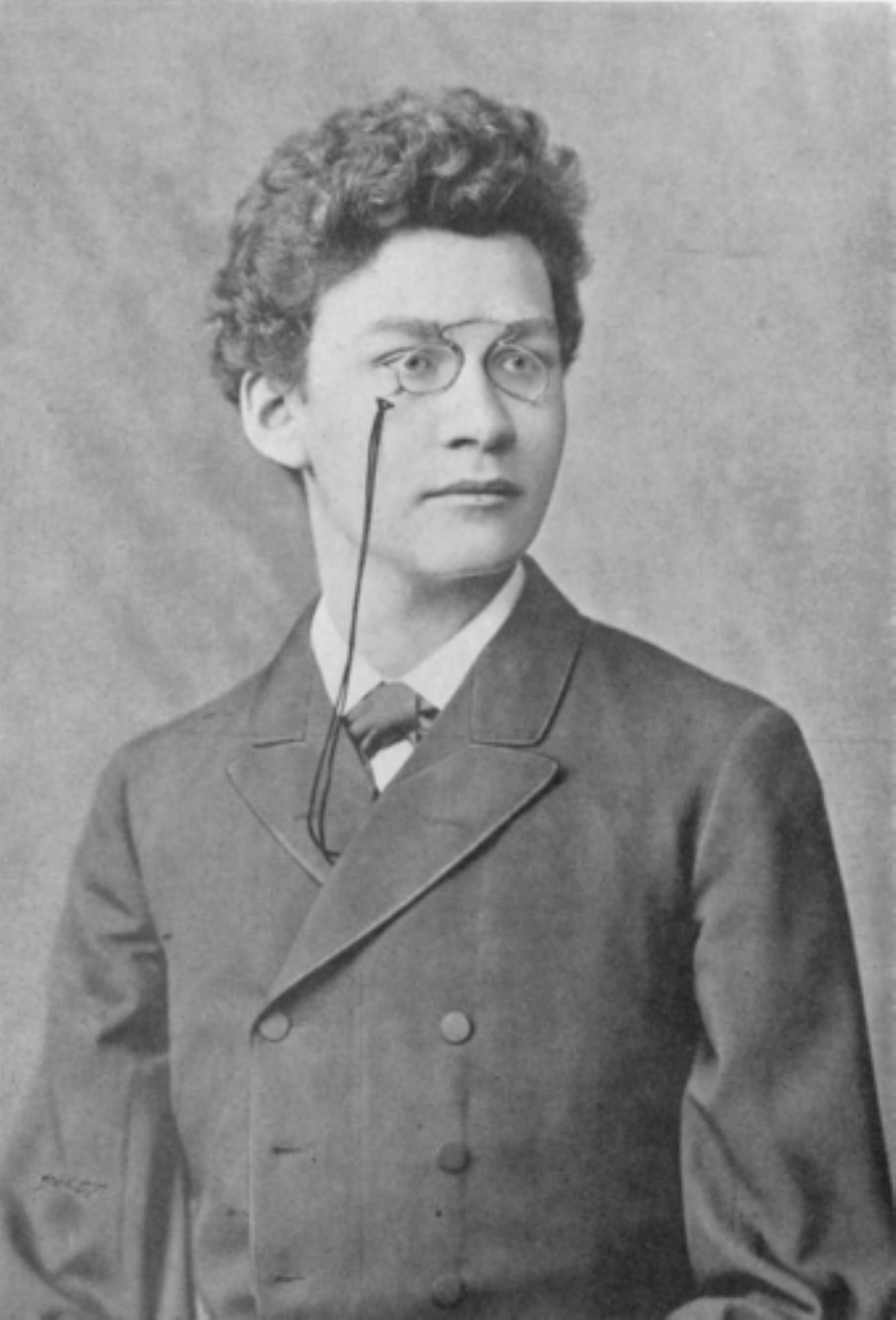}
\caption{Hermann Minkowski}
\end{wrapfigure}

\label{Geometry of numbers}
 \index{Siegel's formula} \index{Minkowski}
\begin{quote}                         
    ``Henceforth space by itself, and time by itself, are doomed to fade away into mere shadows, and only a kind of union of the two will preserve an independent reality.''
  -- Hermann Minkowski    \index{Minkowski, Hermann}
\end{quote}

\section{Intuition}
To see a wonderful and fun application of Poisson summation, we give a relatively easy extension of Minkowski's first theorem, in the Geometry of Numbers.   
Minkowski's theorem gives the existence of an integer point inside  symmetric 
bodies in $\R^d$, once we know their volume is sufficiently large.   

We'll explore, and prove,  a more powerful identity which 
is now a classical result of Carl Ludwig Siegel (Theorem \ref{Siegel}), yielding an identity between Fourier transforms of convex bodies and their volume.  Our proof of this identity of Siegel uses Poisson summation, applied to the convolution of an indicator function with itself.  \index{Poisson summation formula}

The geometry of numbers is an incredibly beautiful field, and too vast to encompass in just one chapter (see note  \ref{new books, geometry of numbers}).  This chapter, as well as chapters
\ref{Chapter.geometry of numbers II} and \ref{Chapter.geometry of numbers III}, which treat the classical geometry of numbers, together give just a taste of a giant and thriving forest.

\bigskip
\section{Minkowski's first convex body Theorem}


\begin{wrapfigure}{L}{0.55\textwidth}
\centering
\includegraphics[width=0.5\textwidth]{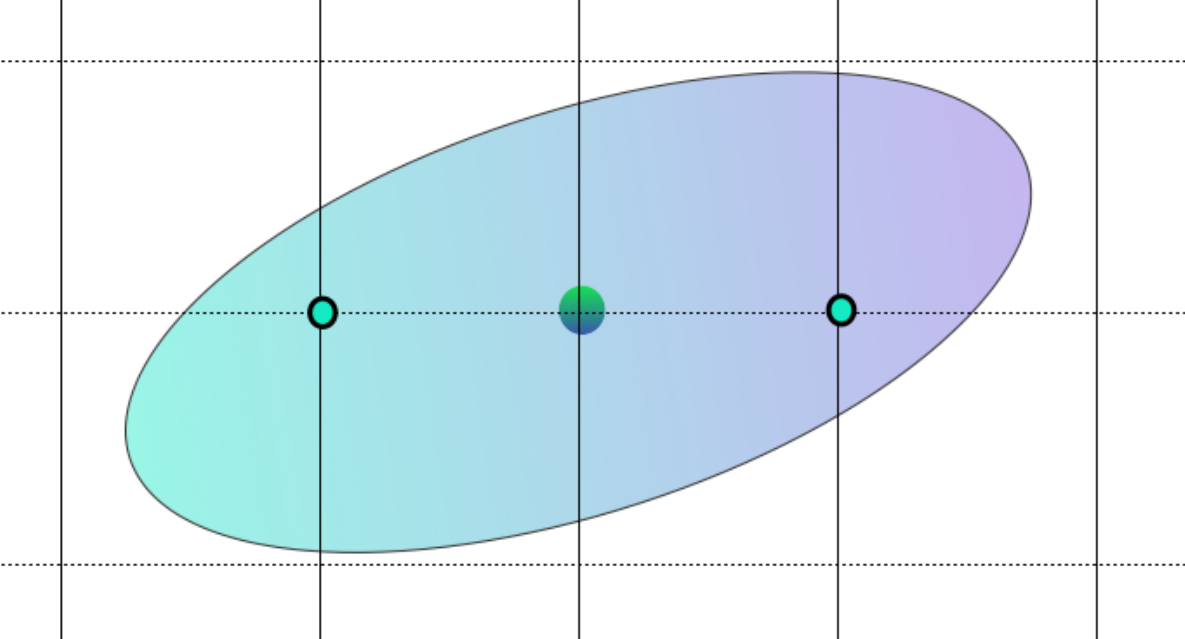}
\caption{A convex, centrally symmetric body in $\R^2$, with area bigger than $4$, containing two nonzero integer points.}    \label{convex body}
\end{wrapfigure}
Minkowski initiated the field that we call today `the geometry of numbers', around 1890.   To begin, we define a {\bf body} $\P$ in $\mathbb R^d$ as a compact set.  In other words,  $\P$ is a bounded, closed set.  Most of the time, it is useful to work with convex bodies that enjoy the following symmetry. 
 We call a body $\P$  {\bf centrally symmetric}, also called {\bf symmetric about the origin},  if  for all $\x \in \R^d$ we have
 \begin{equation} \label{definition of symmetric body}
 \x \in \P  \iff  -\x \in \P.
 \end{equation}

\bigskip 

A body $\P$ is called {\bf symmetric} if some translation of $\P$ is symmetric about the origin.  
For example, the ball 
$\{ x\in \R^d \mid \| x\| \leq 1\}$ is centrally symmetric, and the translated ball 
$
\{ x\in \R^d \mid \| x- w\| \leq 1\}
$
 is symmetric, but not centrally symmetric.   An initial, motivating question in the geometry of numbers is:
 
 \begin{question}\label{Rhetorical question, centrally symmetric}
 {\rm[Rhetorical]}  How large does a convex body $\P$ have to be in order to contain a nonzero integer point?
 \end{question}
 
\begin{figure}[htb]
\begin{center}
\includegraphics[totalheight=2.5in]{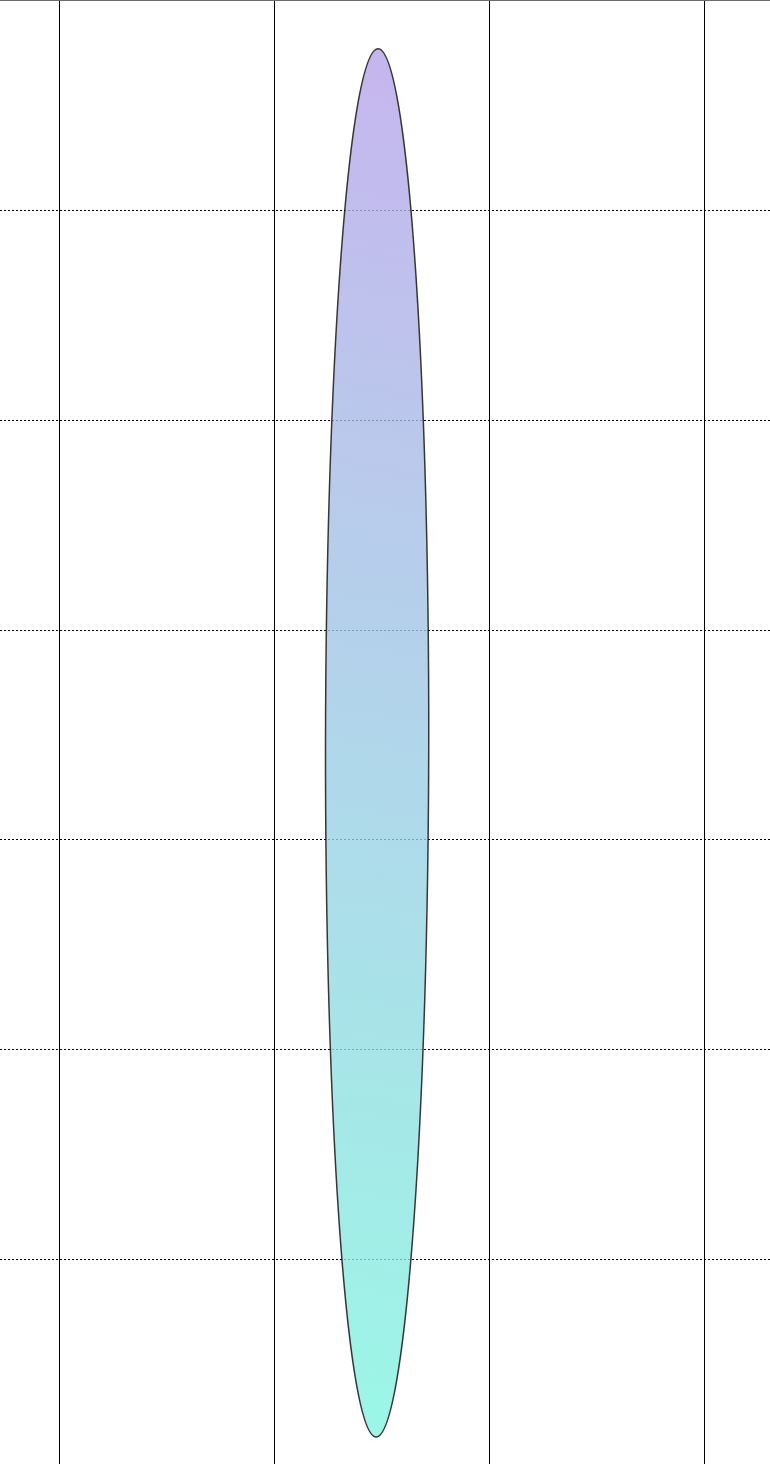}
\end{center}
\caption{A convex symmetric body in $\R^2$, which is not centered at the origin, may be constructed with arbitrarily large volume and simultaneously with 
no integer points.}    \label{nonexample}
\end{figure}

 If we are not careful, then Figure \ref{nonexample}, for example, shows that $\P$ can be as large as we like, and yet never contain an integer point.  So without further hypotheses, there are no positive answers to Question \ref{Rhetorical question, centrally symmetric}.  Therefore, it is natural to assume that our body $\P$ is positioned in a `nice' way relative to the integer lattice, and centrally symmetry is a natural assumption in this respect.
 
\bigskip
\begin{thm}[Minkowski's first convex body Theorem for $\Z^d$] 
\label{Minkowski's convex body theorem, for Z^d}
Let $K$ be a $d$-dimensional convex body in $\mathbb R^d$, symmetric about the origin.
\begin{equation}
\text{ If }   \vol K > 2^d,  \text{  then }   K
\text{  must contain a nonzero integer point in its interior}. 
\end{equation}
\hfill $\square$
\end{thm}

Sometimes this classical and very useful result of Minkowski is stated in its contrapositive form: 
Let $K \subset \R^d$ be any convex body, symmetric about the origin.   
\begin{equation}
\text{  If the only integer point
in the interior of }  K  \text{ is the origin, then } \vol K \leq 2^d.
\end{equation}

It is natural -  and straightforward - to extend this result to any lattice $\L:= M(\Z^d)$, by simply applying the linear transformation $M$ to both the integer lattice, and to the convex body $K$. The conclusion is the following, which is the version that we will prove as a consequence of Siegel's Theorem \ref{Siegel}.

\begin{thm}[Minkowski's first convex body Theorem for a lattice $\L$]  
\label{Minkowski convex body Theorem for L}
Let $K$ be a $d$-dimensional convex body in $\mathbb R^d$, symmetric about the origin, and let
$\L$ be a (full rank) lattice in $\R^d$.
\begin{equation}  \label{Minkowski 2}
\text{ If }   \vol K > 2^d (\det \L),  \text{  then }   K
\text{  must contain a nonzero point of } \L \text{ in its interior}. 
\end{equation}
\end{thm}
\begin{proof}
The proof appears below - see \hyperlink{first proof of Minkowski}{\rm{``first proof of Minkowski''}}.    
We also give a second proof in Section \ref{section:Blichfeldt}, using Blichfeldt's methods - see
\hyperlink{second proof of Minkowski}{\rm{``second proof of Minkowski''}}.
\end{proof}
These very important initial results of Minkowski \cite{Minkowski} have found applications in algebraic number theory, diophantine analysis, combinatorial optimization, and other fields.
In the next section we show that Minkowski's result \eqref{Minkowski 2} follows as a special case of Siegel's formula.

\begin{figure}[htb]
 \begin{center}
\includegraphics[totalheight=1.8in]{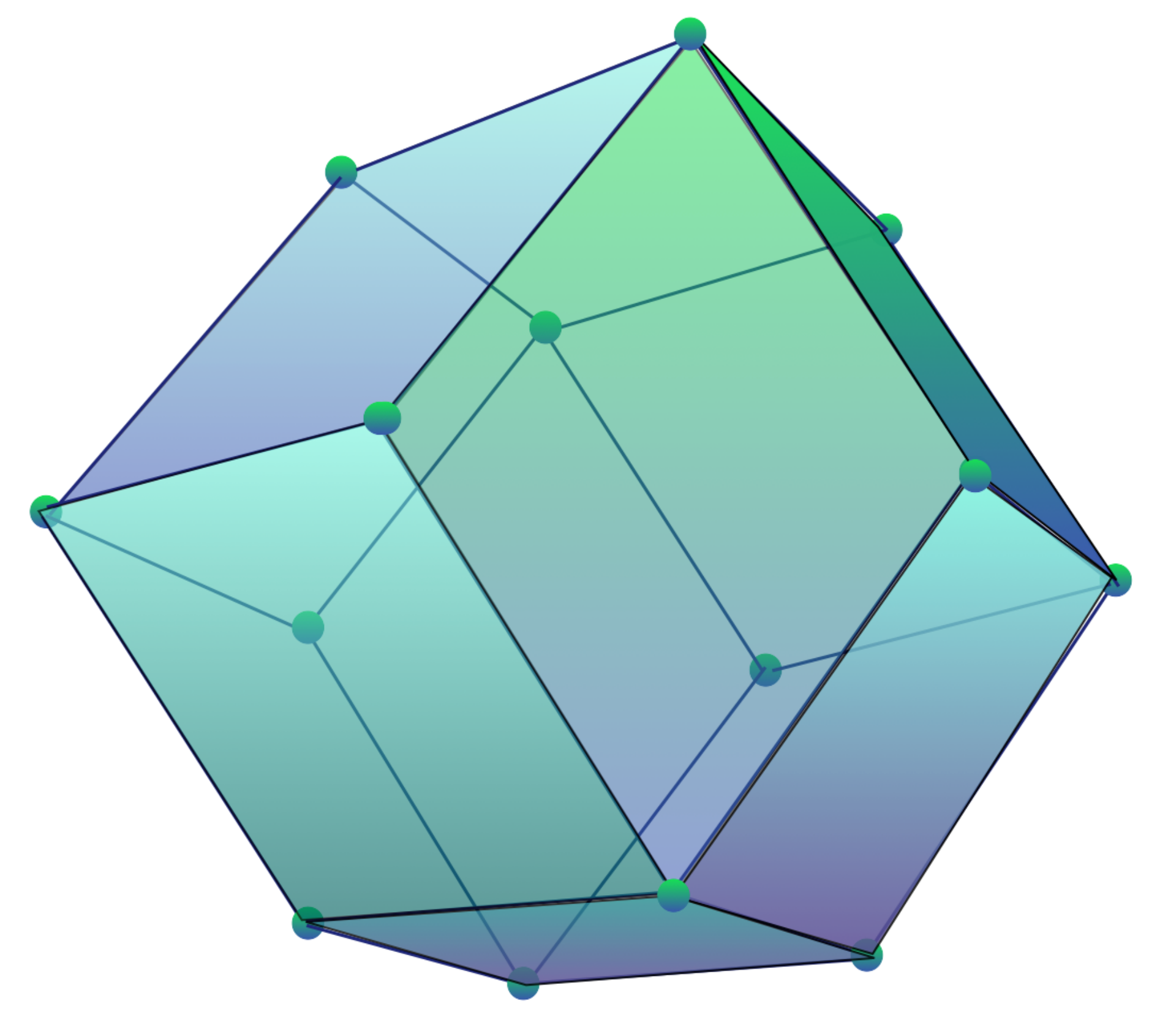}
 \end{center}
\caption{The Rhombic dodecahedron, a $3$-dimensional symmetric polytope that tiles $\R^3$ by translations, and is another extremal body for Minkowski's convex body Theorem. }
\label{The Rhombic dodecahedron}
\end{figure}

\section{Siegel's extension of Minkowski:  \\
a Fourier transform identity for convex bodies}

\begin{quote}                         
    ``Behind every inequality there is an equality - find it.''
    
  -- Basil Gordon  \index{Gordon, Basil}
\end{quote}

Minkowski's Theorem \ref{Minkowski convex body Theorem for L}  suggests that behind his inequality 
$2^d > \vol K$ there may hide an interesting equality:
\[
2^d = \vol K + \text{(some positive error term)}.
\]
A natural and motivating question is ``what form does this positive error term take?"   Siegel found it, and as we will soon see, it naturally leads us
to the Fourier transform of $K$.
First, an important construction in the geometry of numbers is the  {\bf Minkowski sum} of convex bodies. \index{Minkowski sum}
 Given two convex bodies $K, L \subset \R^d$, their Minkowski sum is defined by
\[
K + L := \{ x + y \mid x \in K, y \in L\}.
\]
Another related construction, appearing in some of the results below, is 
\[
K- L := \{  x-y \mid x \in K, y \in L\},
\]
the Minkowski difference of $K$ and $L$.  Of course, we also have $K-L = K + (-L)$, where 
$-L:=\{ -x \mid x \in \L\}$.
A very useful special case is the gadget known as the
 {\bf Minkowski symmetrized  body} of $K$, \index{symmetrized body}
defined by 
\begin{equation}
\frac{1}{2} K - \frac{1}{2} K,
\end{equation}
and often also called the {\bf difference body} of $\frac{1}{2}K$.  
Given any set $K\subset \R^d$, the difference body $K-K$ is centrally symmetric.  To see this,  suppose $x\in K-K$, so we may write
$x= y-z$, with $y, z \in K$.  Then $-x = z-y \in K-K$. 

In addition,  we have the fortuitous and easy fact that a  convex set  $K\subset \R^d$ is centrally symmetric if and only if we have the equality
\begin{equation} \label{centrally symmetric set}
\frac{1}{2} K - \frac{1}{2} K=K.
\end{equation}
 (Exercise \ref{c.s. C equals its symmetrized body}).  Now suppose we are given two convex bodies $K, L\subset \R^d$.  Then the resulting bodies $K+L$, $K-L$ turn out to also be convex (Exercise \ref{convexity of K-K}).  
Another important geometric notion is the dilation of a convex body by a positive real number~$t$:
 \[
 tK := \{ t x \mid x\in K\}, 
 \]
The most basic version of Siegel's theorem is the following identity, which assumes that a convex body $K$ is symmetric about the origin. 
\bigskip
\begin{thm} [Siegel]  \label{Siegel}  \index{Siegel's formula}
Let $K$ be any $d$-dimensional convex body in $\mathbb R^d$, symmetric about the origin, and suppose that the only integer point in the interior of $K$ is the origin.  Then 
\begin{align} \label{Siegel, version 1}
2^d
&= \vol K +  \frac{4^d}{  \vol K    }   \sum_{\xi \in \Z^d - \{0\}} \left| \hat 1_{\frac{1}{2} K}(\xi) \right|^2.
\end{align}
\hfill $\square$
\end{thm}
We now prove the following extension of Siegel's Theorem   \eqref{Siegel},  namely \eqref{Siegel.formula} below, which applies to bodies that are not necessarily convex, nor necessarily symmetric about the origin.  
Our proof of Theorem \ref{Siegel for general lattices} below
 consists of yet another application of Poisson summation. \index{Poisson summation formula}
It turns out that if $K$ is any convex body, then $f:= 1_{\{\frac{1}{2} K\}}*1_{\{-\frac{1}{2} K\}}$ is a nice function (Exercise \ref{convolution of indicators is a nice function}), 
in the sense that Poisson summation  \eqref{nice functions} holds for $f$.
So Theorem \ref{Siegel} is a consequence  of the following extension to bodies that are 
not necessarily convex or symmetric.

\medskip
\begin{thm}[Siegel's formula, for a body $K$, and a lattice $\L$]   
\label{Siegel for general lattices}
\index{Siegel's formula}
Let $K\subset \R^d$ be a body (compact set)  for which
the convolution $1_{\frac{1}{2} K}*1_{-\frac{1}{2} K}$ is a nice function.  
If the only integer point in the interior of the difference body 
$ \frac{1}{2}K -  \frac{1}{2}K$ is the origin, then 
\begin{equation}\label{Siegel.formula}
2^d
= \vol K +  \frac{4^d}{  \vol K    }   \sum_{\xi \in \Z^d - \{0\}} \left| \hat 1_{\frac{1}{2} K}(\xi) \right|^2.
\end{equation}

More generally, if we replace the lattice $\Z^d$ by any full-rank lattice $\L$, and assume that the only lattice point of $\L$ in the interior of  $ \frac{1}{2}K -  \frac{1}{2}K$ is the origin, then we have:
\begin{equation}\label{Siegel formula 2}
2^d \det \L
= \vol K +  \frac{4^d}{  \vol K    }   \sum_{\xi \in \L^* - \{0\}} \left| \hat 1_{\frac{1}{2} K}(\xi) \right|^2.
\end{equation}
\end{thm}

\begin{proof} 
We start with the function
\begin{equation}
f(x):= \left(   1_{\frac{1}{2} K}*1_{-\frac{1}{2} K} \right) (x),
\end{equation}
which is continuous on $\mathbb R^d$, and we plug $f$ into Poisson summation
 \eqref{Poisson.summation2}:
\index{Poisson summation formula}
\begin{align}
\sum_{n \in \Z^d}  f(n) =  \sum_{\xi \in \Z^d}  \hat f(\xi).  
\end{align} 
We first compute the left-hand-side of Poisson summation, using the definition of $f$:
\begin{align}
\sum_{n \in \Z^d}  f(n)  &= \sum_{n \in \Z^d}  \int_{\R^d}    1_{\frac{1}{2} K}(y) 1_{-\frac{1}{2} K}(n - y) dy \\
&= \sum_{n \in \Z^d}  \int_{\R^d}    1_{\frac{1}{2} \rm{int} K}(y) 1_{-\frac{1}{2} \rm{int} K}(n - y) dy,
\end{align}
where the last step follows from the fact that the integral does not distinguish between a convex set or its closure.   
Now we follow the definition of containment: $y \in \frac{1}{2} K$ and $n - y \in -\frac{1}{2} K$ imply that the integer point $n \in  \frac{1}{2} K -\frac{1}{2} K$.  But by hypothesis $ \frac{1}{2} K -\frac{1}{2} K$ contains the origin as its {\em only} interior integer point, so the left-hand-side of the 
Poisson summation formula contains only one term, namely the $n=0$ term:
\begin{align}
\sum_{n \in \Z^d}  f(n)  &= \sum_{n \in \Z^d}  \int_{\R^d}    1_{\frac{1}{2} K}(y) 1_{-\frac{1}{2} K}(n - y) dy \\
&= \int_{\R^d}    1_{\frac{1}{2} K}(y) 1_{-\frac{1}{2} K}(- y) dy \\
&= \int_{\R^d}    1_{\frac{1}{2} K}(y)  dy \\
&= \vol \left( {\frac{1}{2} K} \right) = \frac{\vol K}{2^d}.
\end{align}
\noindent
On the other hand, the right-hand-side of Poisson summation gives us:
\begin{align}
\sum_{\xi \in \Z^d}  \hat f(\xi)  
&=  \sum_{\xi \in \Z^d}    {\hat 1}_{\frac{1}{2} K}(\xi)  {\hat 1}_{-\frac{1}{2} K}(\xi)     \\
&=  \sum_{\xi \in \Z^d}  \int_{\frac{1}{2} K} e^{2\pi i \langle \xi, x \rangle} dx  
\int_{-\frac{1}{2} K} e^{2\pi i \langle \xi, x \rangle} dx \\
&=  \sum_{\xi \in \Z^d}  \int_{\frac{1}{2} K} e^{2\pi i \langle \xi, x \rangle} dx  
\int_{\frac{1}{2} K} e^{2\pi i \langle -\xi, x \rangle} dx \\
&=  \sum_{\xi \in \Z^d}  \int_{\frac{1}{2} K} e^{2\pi i \langle \xi, x \rangle} dx  
\ \overline{
\int_{\frac{1}{2} K} e^{2\pi i \langle \xi, x \rangle} dx
} \\   \label{pulling out the zero term}
&=  \sum_{\xi \in \Z^d} \left| \hat 1_{\frac{1}{2} K}(\xi) \right|^2 \\
&=  \left| \hat 1_{\frac{1}{2} K}(0) \right|^2 + \sum_{\xi \in \Z^d - \{0\}} \left| \hat 1_{\frac{1}{2} K}(\xi) \right|^2 \\
&= \frac{\vol^2 K}{4^d} + \sum_{\xi \in \Z^d - \{0\}} \left| \hat 1_{\frac{1}{2} K}(\xi) \right|^2,
\end{align}
where we have pulled out the $\xi=0$ term from the series \eqref{pulling out the zero term}. 
 So we've arrived at
\begin{align*}
\frac{\vol K}{2^d}
&= \frac{\vol^2 K}{4^d} + \sum_{\xi \in \Z^d - \{0\}} \left| \hat 1_{\frac{1}{2} K}(\xi) \right|^2,
\end{align*}
yielding the required identity:
\begin{align*}
2^d  &=      \vol K +   \frac{4^d}{\vol K}\sum_{\xi \in \Z^d - \{0\}} \left| \hat 1_{\frac{1}{2} K}(\xi) \right|^2.
\end{align*}

Finally, to prove the stated extension to all lattices $\L$, we use the slightly more general 
form of Poisson summation, Theorem \ref{The Poisson Summation Formula, for lattices}, 
valid for any lattice $\L$:
\begin{align}
\sum_{n \in \L}  f(n) =  \frac{1}{\det \L}\sum_{\xi \in \L^*}  \hat f(\xi).  
\end{align} 
All the steps of the proof above are identical, except for the factor of $\frac{1}{\det \L}$, so that we arrive at
the required identity of Siegel for arbitrary lattices:
\begin{align}\label{Siegel, take 2}
\frac{\vol K}{2^d}
&= \frac{\vol^2 K}{4^d \det \L} + \frac{1}{\det \L}\sum_{\xi \in \L^* - \{0\}} \left| \hat 1_{\frac{1}{2} K}(\xi) \right|^2.
\end{align} 
\end{proof}


The proof of Minkowski's convex body Theorem for lattices, namely Theorem 
\ref{Minkowski convex body Theorem for L} above, now follows immediately.

\begin{proof}[Proof of Theorem \ref{Minkowski convex body Theorem for L}]
 \rm{[Minkowski's convex body Theorem for a lattice} $\L$]
 \hypertarget{first proof of Minkowski}
 Applying Siegel's Theorem \ref{Siegel for general lattices} to the centrally symmetric body $K$, 
 we see that the lattice sum on the right-hand-side of  identity \eqref{Siegel.formula}
  contains only non-negative terms.  It follows that we immediately get the analogue of Minkowski's result for a given cenetrally symmetric body $K$ and a lattice $\L$, in its contrapositive form: 
\begin{align}   
&\text{If the only lattice point of $\L$ in the interior of $K$  is the origin, }  \\
&\text{then } 2^d \det \L \geq \vol K.  
\end{align}
\end{proof}
In fact, we can easily extend 
Minkowski's Theorem \ref{Minkowski convex body Theorem for L},  
using the same ideas of the latter proof, by using 
Siegel's Theorem \ref{Siegel for general lattices} 
so that it applies to non-symmetric bodies as well (but there's a small `catch' - see Exercise \ref{Extending Minkowski to nonconvex bodies}).

Enrico Bombieri \cite{Bombieri} found an extension of Siegel's formula  \ref{Siegel for general lattices}, allowing the body to contain any number of lattice points.  In recent work, Michel Martins and S. R.\cite{MartinsRobins} found an extension of Bombieri's results, using the cross-covariogram of two bodies.  


\bigskip
\section{Tiling and multi-tiling Euclidean space   by translations of polytopes}
\index{tiling}
First, we give a `spectral' equivalence for the tiling of Euclidean space by a single polytope, using only translations by a lattice.  It will turn out that the case of equality in Minkowski's convex body Theorem is characterized precisely by the polytopes that tile $\R^d$ by translations.  These bodies are called extremal bodies. 

More generally, we would like to also consider the notion of multi-tiling, as follows.   We say that a polytope
$\P$ {\bf $k$-tiles $\R^d$ by using a set of translations $\L$} if for some integer $k$, we have
\begin{equation}
\sum_{n \in \L} 1_{\P + n}(x) = k,
\end{equation}
for all $x \in \R^d$, except those points $x$ that lie on the boundary of $\P$ or its translates under $\L$ (and of course these exceptions form a set of measure $0$ in $\R^d$).   In other words, 
$\P$ is a $k$-tiling body if almost every $x \in \R^d$ is covered by exactly $k$ translates of $\P$.

Other synonyms for $k$-tilings in the literature are {\bf multi-tilings}  of
$\R^d$, or {\bf tiling  at level $k$}.
When $\L$ is a lattice, we will say that such a $k$-tiling is {\bf periodic}. 
A common research theme is to search for tilings which are not necessarily periodic, but this is a difficult problem in general.
The classical notion of tiling, such that there are no overlaps between the interiors of any two tiles, corresponds here to the case $k=1$.  We have the following wonderful dictionary between multi-tiling $\R^d$ by translations of a convex body $\P$, and a certain vanishing 
property of the Fourier transform of $\P$, due to Mihalis Kolountzakis (\cite{Kolountzakis1}, \cite{Kolountzakis2}). 

\bigskip
\begin{thm} [Kolountzakis]
\label{zero set of the FT of a polytope}
Suppose that $\P\subset \R^d$ is a  compact set with positive $d$-dimensional volume.  The following two properties are equivalent:
\begin{enumerate}[(a)]
\item  $\P$ $k$-tiles $\R^d$ by translations with a lattice $\L$.
\item
\[
\hat 1_\P(\xi) = 0,
\]
 for all nonzero $\xi \in \L^*$, the dual lattice. 
\end{enumerate}
Either of these conditions also implies that $k = \frac{\vol \P}{ \det \L}$, an integer.
\end{thm}
\begin{proof}
We begin with the definition of multi-tiling, so that by assumption
\begin{equation} \label{by definition of k-tiling}
\sum_{n \in \L} 1_{\P + n}(x) = k,
\end{equation}
for all $x \in \R^d$ except those points $x$ that lie on the boundary of $\P$ or its translates under $\L$ (and of course these exceptions form a set of measure $0$ in $\R^d$).  A trivial but useful 
observation is that 
\[
1_{\P + n}(x) =1  \iff  1_{\P}(x-n) =1, 
\]
so we can rewrite the defining identity \eqref{by definition of k-tiling}  as $\sum_{n \in \L} 1_{\P}(x-n) = k$.  Now we notice that the left-hand-side is a periodic function of $x$, namely 
\[
F(x) := \sum_{n \in \L} 1_{\P}(x-n)
\]
 is periodic in $x$ with $\L$ as its set of periods.   This is easy to see: if we let $l \in \L$, then $F(x + l) = \sum_{n \in \L} 1_{\P}(x+ l -n) = \sum_{m \in \L} 1_{\P}(x+ m) =  F(x)$, because the lattice $\L$ is invariant under a  translation by any vector that belongs to it. 

The following `intuitive proof' would in fact be rigorous if we were allowed to use  `generalized functions', 
but since we do not use them in this book, we label this part of the proof as `intuitive', and we then give a rigorous proof, using functions rather than generalized functions.

[{\bf Intuitive proof}]  \  By Theorem \ref{Fourier series for periodic functions}, we may expand $F$ into its Fourier series, because it is a periodic function on $\R^d$.  Now by Poisson summation, namely 
Theorem \ref{The Poisson Summation Formula, for lattices}, we know that its Fourier coefficients are the following:
\begin{equation}\label{Poisson version of k-tiling}
 \sum_{m \in \L} 1_{\P}(x+m) = \frac{1}{\det \L} \sum_{\xi \in \L^*} \hat 1_{\P}(\xi) e^{2\pi i \langle \xi, x \rangle},
\end{equation}
If we now make the assumption that $\hat 1_{\P}(\xi)=0$ for all nonzero $\xi \in \L^*$, then by \eqref{Poisson version of k-tiling} this assumption is equivalent to 
\[
\sum_{m \in \L} 1_{\P}(x+m) = \frac{\hat 1_\P(0)}{\det \L} =  \frac{\vol \P}{\det \L}.
\]
This relation means that we have a $k$-tiling, where $k:= \frac{\vol \P}{\det \L}$. 
Now we replace the intuitive portion of the proof with a rigorous proof.

[{\bf Rigorous proof}] \  In order to apply Poisson summation, \index{Poisson summation formula}
 it is technically necessary to replace $1_P(x)$ by 
a smoothed version of it, in \eqref{Poisson version of k-tiling}.   Because this process is so common and useful in applications, this proof is instructive.  
We pick an approximate identity $\phi_n$, which is also compactly supported and continuous.  Applying the Poisson summation formula of Theorem \ref{PracticalPoisson} to the smoothed function $1_P*\phi_n$, we get:
\begin{align}\label{rigorous limit for multi-tiling}
 \sum_{m \in \L} \left( 1_\P*\phi_n   \right)(x+m)  &=  \frac{1}{\det \L} \sum_{\xi \in \L^*} \hat 1_{\P}(\xi)  \hat \phi_n(\xi) e^{2\pi i \langle \xi, x \rangle}.
\end{align}
Using the fact that the convolution of two compactly supported functions is itself compactly supported,  we see  that $1_\P*\phi_n $ is again compactly supported.  Thus the sum on the LHS of \eqref{rigorous limit for multi-tiling} is a finite sum.  
Performing a separate computation, we take the limit as $n\rightarrow \infty$ inside this finite sum, and using 
Theorem \ref{approximate identity convolution} (due to the continuity of $1_\P*\phi_n$), we obtain
\[
\lim_{n\rightarrow \infty}  \sum_{m \in \L} \left( 1_\P*\phi_n   \right)(x+m) 
=  \sum_{m \in \L}  \lim_{n\rightarrow \infty}    \left( 1_\P*\phi_n   \right)(x+m) 
=\sum_{m \in \L}  1_\P(x+m).
\]
Now using our Poisson summation IV (Theorem \ref{PracticalPoisson},  part \ref{practical Poisson summation, first version}), we have
\begin{align} \label{smoothed-out RHS}
\sum_{m \in \L}  1_\P(x+m) =    \frac{1}{\det \L} 
     \sum_{\xi \in \L^*} \hat 1_{\P}(\xi)  \hat \phi_n(\xi) 
     \, e^{2\pi i \langle \xi, x \rangle}.
\end{align}
 for all sufficiently large values of $n$.  Separating  the term $\xi=0$ on the RHS of this Poisson summation formula, we have:
\begin{align}
\sum_{m \in \L}  1_\P(x+m) &=   \frac{\hat 1_\P(0)}{\det \L} 
   +   \sum_{\xi \in \L^*-\{0\}} \hat 1_{\P}(\xi)  \hat \phi_n(\xi) 
     \, e^{2\pi i \langle \xi, x \rangle}\\   \label{equivalent condition for tiling}
   &=  
   \frac{\vol \P}{\det \L} 
   +   \sum_{\xi \in \L^*-\{0\}} \hat 1_{\P}(\xi)  \hat \phi_n(\xi) 
     \, e^{2\pi i \langle \xi, x \rangle}.
\end{align}
Now, $\hat 1_{\P}(\xi) =0$ for all $\xi \in \L^*-\{0\}$ in  
\eqref{equivalent condition for tiling} will hold
\begin{align*}
&\iff
\sum_{m \in \L}  1_\P(x+m) =   \frac{\vol \P}{\det \L}, 
\end{align*}
an equivalent condition which we may write as $\sum_{m \in \L}  1_\P(x+m) = k$, where necessarily
$k:=   \frac{\vol \P}{\det \L}$.   The condition $\sum_{m \in \L}  1_\P(x+m) = k$ means that $\P$ $k$-tiles
$\R^d$ by translations with the lattice $\L$, and also implies that $k$ must be an integer. 
\end{proof}

\begin{figure}[htb]
 \begin{center}
\includegraphics[totalheight=2.7in]{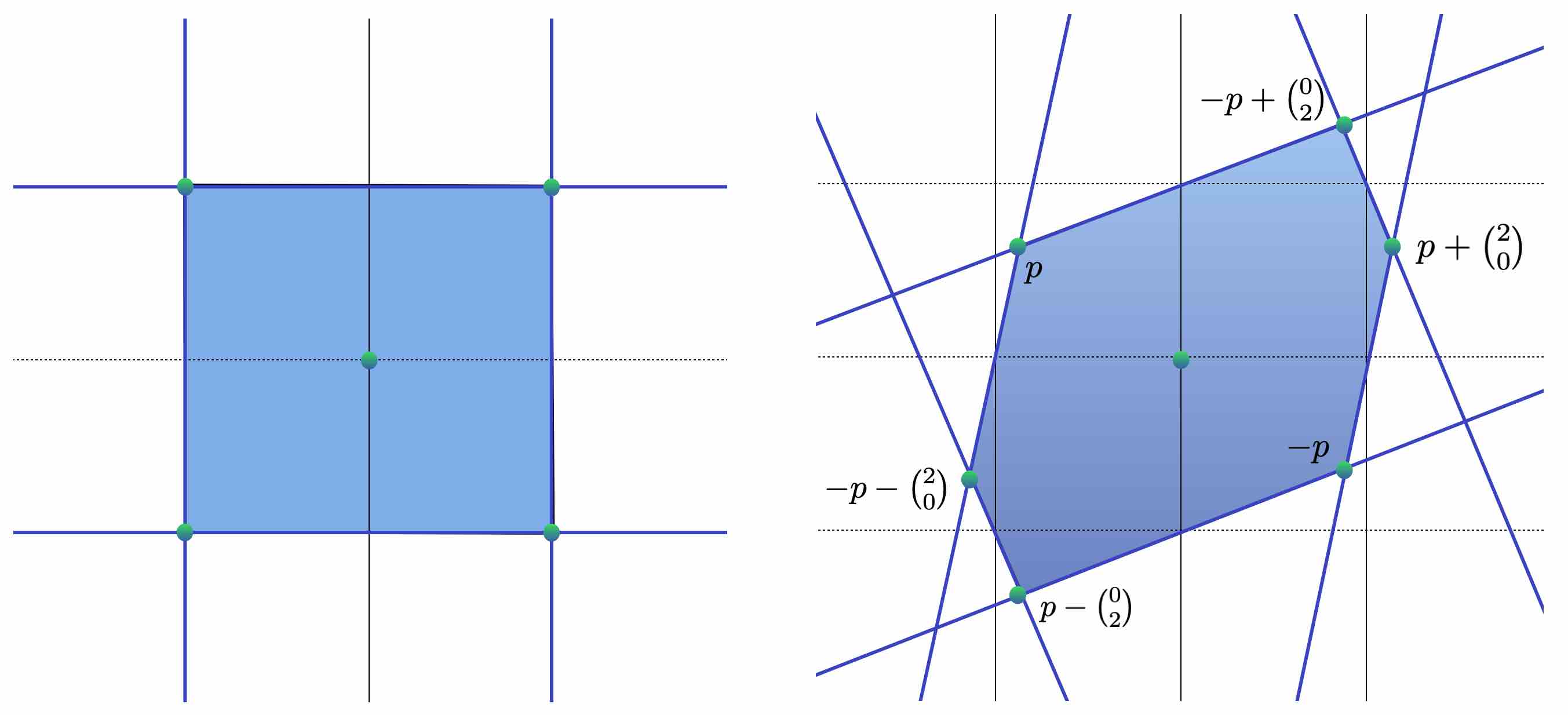}
 \end{center}
\caption{The square on the left is an extremal body in $\R^2$, relative to the integer lattice, and has area $4$ as expected.   More generally, there is a continuous ($2$-parameter) family
of extremal hexagons, each having area $4$, with no integer points in their interior.  This continuous family is parametrized by the point $p\in \R^2$ in the figure on the right.}
 \label{Extremal body in R^2}
\end{figure}


In 1905, Minkowski gave necessary conditions for a polytope $\P$ to tile $\R^d$ by translations.  Later, Venkov and independently McMullen found sufficient conditions as well, culminating in the following fundamental result.

\begin{thm}[Minkowski-Venkov-McMullen] \label{Minkowski-Venkov-McMullen}
A polytope $\P$ tiles $\R^d$ by translations if and only if the following $3$ conditions hold:
\begin{enumerate}
\item $\P$ is a symmetric polytope.
\item The facets of $\P$ are symmetric polytopes.
\item  Fix any face $F\subset \P$ of codimension $2$, and project $\P$ onto the $2$-dimensional plane that is orthogonal to the $(d-2)$-dimensional affine span of $F$.   Then this projection is either a parallelogram, or a centrally symmetric hexagon.
\end{enumerate}
\end{thm}


\bigskip
\section{Extremal bodies}  \label{extremal bodies}
\index{extremal body}

An {\bf extremal body}, relative to a lattice $\L$,  is a convex symmetric body $K$ which contains exactly one lattice point of $\L$ in its interior, and such that
\[
\vol K = 2^d(\det \L).
\]
In other words, an extremal body satisfies the hypotheses of 
Minkowski's inequality (Theorem \ref{Minkowski convex body Theorem for L}), and attains the equality case.  

If we just look at Siegel's equation  \eqref{Siegel.formula} a bit more closely, we quickly get a nice corollary that 
arises by combining Theorem \ref{zero set of the FT of a polytope} and Siegel's Theorem \ref{Siegel}.
Namely, equality occurs in Minkowski's  convex body theorem if and only if 
$K$ tiles \index{tiling} $\R^d$ by translations.  
Let's prove this.
\begin{thm}[Extremal bodies] \label{thm:extremal bodies}
Let $K$ be any convex, centrally symmetric subset of $\R^d$, and fix a full-rank lattice $\L \subset \R^d$. 
Suppose that the only point of $\L$ in the interior of $K$ is the origin.  
Then:
\begin{quote}
$2^d \det \L= \vol K \  \iff \ \frac{1}{2}K$ tiles $\R^d$ by translations with the lattice $\L$.
\end{quote}
\end{thm}
\begin{proof}
By Siegel's formula \eqref{Siegel formula 2}, we have
\begin{equation}
2^d \det \L
= \vol K +  \frac{4^d}{  \vol K    }   \sum_{\xi \in \L^* - \{0\}} \left| \hat 1_{\frac{1}{2} K}(\xi) \right|^2.
\end{equation}
Therefore, the assumption $2^d \det \L = \vol K$ holds $\iff$ 
\begin{equation}
0 =   \frac{4^d}{  \vol K    }   \sum_{\xi \in \L^* - \{0\}} \left| \hat 1_{\frac{1}{2} K}(\xi) \right|^2,
\end{equation}
 $\iff$ all of the non-negative summands $ \hat 1_{\frac{1}{2} K}(\xi) =0$, for all nonzero $\xi \in \L^*$.   Now we would like to 
 use Theorem \ref{zero set of the FT of a polytope} to show the required tiling equivalence, namely that 
$\frac{1}{2} K$ tiles $\R^d$ by translations with the lattice $\L$.  We have already verified condition (a) of Theorem \ref{zero set of the FT of a polytope}, applied to the body $\frac{1}{2} K$, namely that
  $\hat 1_{\frac{1}{2} K}(\xi) =0$, for all nonzero $\xi \in \L^*$.  
  
  To verify condition (b) of Theorem 
   \ref{zero set of the FT of a polytope}, we notice that because $\vol\left(  \frac{1}{2} K  \right)= \frac{1}{2^d}\vol K$, it follows that 
   $2^d \det \L= \vol K$ is equivalent to 
   $1 = \frac{ \vol \left( \frac{1}{2} K \right)}{ \det \L}$, so that we may apply Theorem 
   \ref{zero set of the FT of a polytope} with $\P:= \frac{1}{2} K$, and with the multiplicity $k:=1$. 
\end{proof}

\medskip
There is an extension of Theorem  \ref{Minkowski-Venkov-McMullen} (the Minkowski-Venkov-McMullen result) to multi-tilings, as follows. 
\begin{thm}\cite{GravinShiryaevRobins} \label{k-tiling theorem, GravinShiryaevRobins}
If a polytope $\P$ multi-tiles $\R^d$ by translations with a discrete set of vectors, then 
\begin{enumerate}[(a)]
\item $\P$ is a symmetric polytope.
\item The facets of $\P$ are symmetric polytopes.
\end{enumerate}
\end{thm}

In the case that $\P\subset \R^d$ is a rational polytope, meaning that all the vertices of $\P$ have rational coordinates, the latter two necessary conditions for multi-tiling become sufficient conditions as well \cite{GravinShiryaevRobins}.




\begin{figure}[htb]
 \begin{center}
\includegraphics[totalheight=2in]{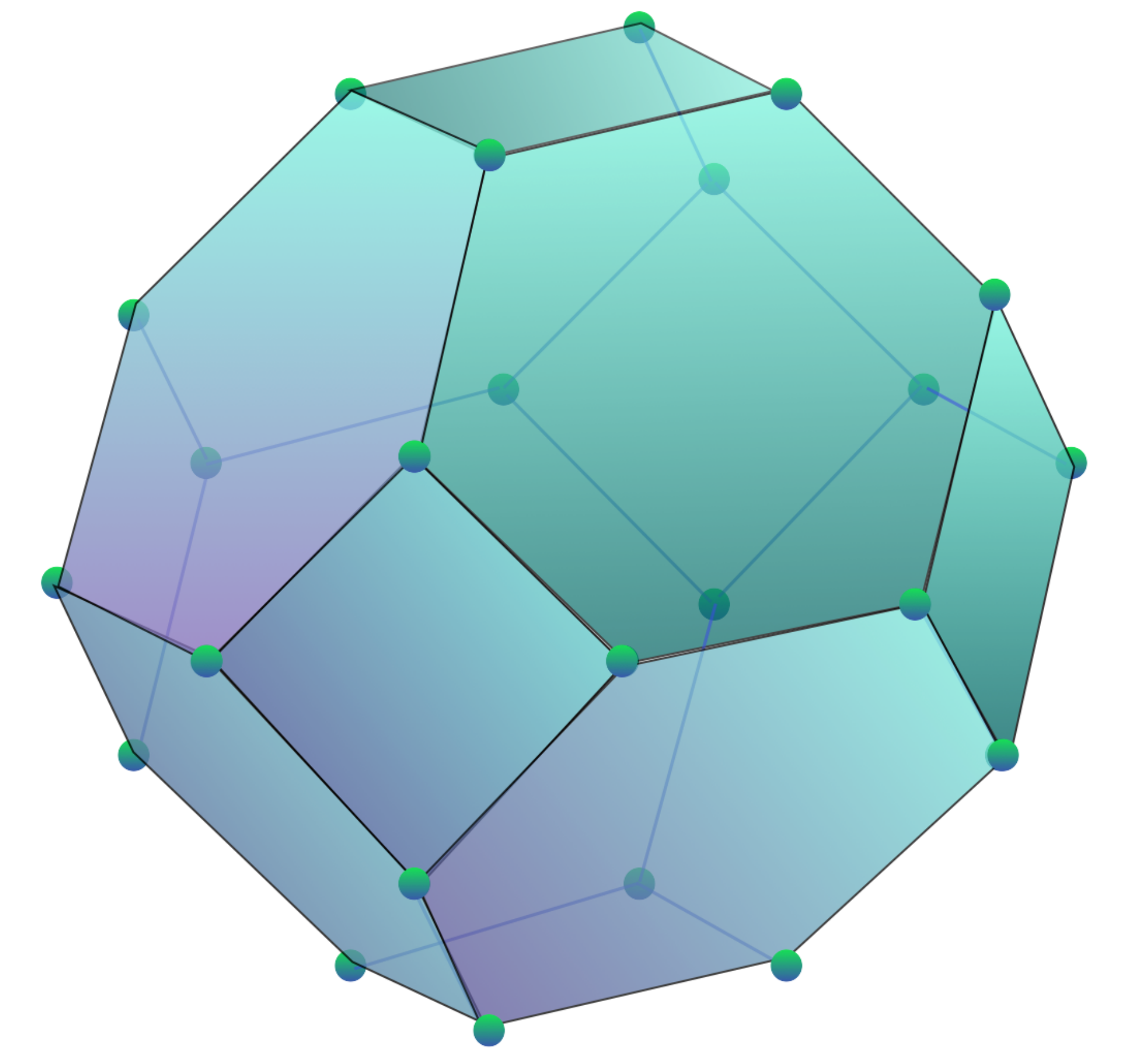}
 \end{center}
\caption{The truncated Octahedron, one of the $3$-dimensional polytopes that tiles $\R^3$ by translations. }
\label{TruncatedOctahedron}
\end{figure}

 \begin{question}[Rhetorical]
 \label{question about equality of polytopes, given equality of their transforms on a lattice}
Is it possible to find two distinct polytopes $\P, Q$ such that
\begin{equation}\label{distinct polytopes, yet their FT's agree on a lattice}
\hat 1_\P(\xi) = \hat 1_{Q}(\xi) \text{ for all } \xi \in \Z^d ?
\end{equation}
\end{question}

\begin{example}
{\rm
We finish this section by answering 
Question \ref{question about equality of polytopes, given equality of their transforms on a lattice}. 
Take any two distinct 
extremal bodies, relative to the lattice $\Z^d$, say $\P, Q\subset \R^d$.  By definition of an extremal body, we know that their volumes must be both equal to $2^d$, so $\hat 1_\P(0) =\vol \P =2^d = \vol Q = \hat 1_{Q}(0)$. 
Moreover, by Theorem \ref{zero set of the FT of a polytope}, we also have:
\[
\hat 1_\P(\xi) =0= \hat 1_{Q}(\xi) \text{ for all } \xi \in \Z^d \setminus \{0\}.
\]
\hfill $\square$
}
\end{example}


 \section{Zonotopes, and centrally symmetric polytopes} 
 \label{Centrally symmetric polytopes}
 \index{centrally symmetric polytope}
   
It's both fun and instructive to begin by seeing how very simple Fourier methods can give us deeper insight into the geometry of symmetric polytopes.   The reader may glance at the definitions above, in 
\eqref{definition of symmetric body}.

\bigskip
\begin{example} 
\rm{
Consider the cross-polytope $\Diamond \subset \R^3$, 
\index{cross-polytope}
defined in Chapter \ref{Chapter.Examples}.  This is a centrally symmetric polytope, but each of its facets is {\em not} a symmetric polytope, because its facets are triangles. 
}
\hfill $\square$
\end{example}

If \emph{all} of the $k$-dimensional faces of a polytope $\P$ are symmetric, for each $1\leq k \leq d$, 
then $\P$ is called a {\bf zonotope}. 
 \index{zonotope}
  Zonotopes form an extremely important class of polytopes, and have various equivalent formulations.
\begin{lem}
A polytope $\P \subset \R^d$ is a zonotope $\iff$  $\P$ has one of the following properties. 
\begin{enumerate}[(a)]
\item  $\P$ is a projection of some $n$-dimensional cube.
\item   $\P$ is the Minkowski sum of a finite number of line segments.
\item Every face of $\P$ is symmetric.
\end{enumerate}
\hfill $\square$
\end{lem}
A projection here means any affine transformation of $\P$, where the rank of the associated matrix may be less than $d$.

 Zonotopes 
 have been very useful in the study of tilings (\cite{Ziegler}, \cite{BeckRobins}).  
 \index{tiling}
 For instance, in dimension $3$, the only  polytopes that tile $\R^3$ by translations with a lattice are zonotopes, and there is a list of $5$ of them (up to an isomorphism of their face posets), called the {\bf Fedorov solids}, and drawn in Figure \ref{Fedorov solids}  (also see our Note \ref{Fedorov Note} below).  \index{Fedorov solids}
 
 By definition, any zonotope  is a symmetric polytope, but the converse is not true; for example, 
the cross-polytope 
\index{cross-polytope}
 is symmetric, but it has triangular faces, which are not symmetric, so the crosspolytope is not a zonotope.

\begin{figure}[htb]
 \begin{center}
\includegraphics[totalheight=2in]{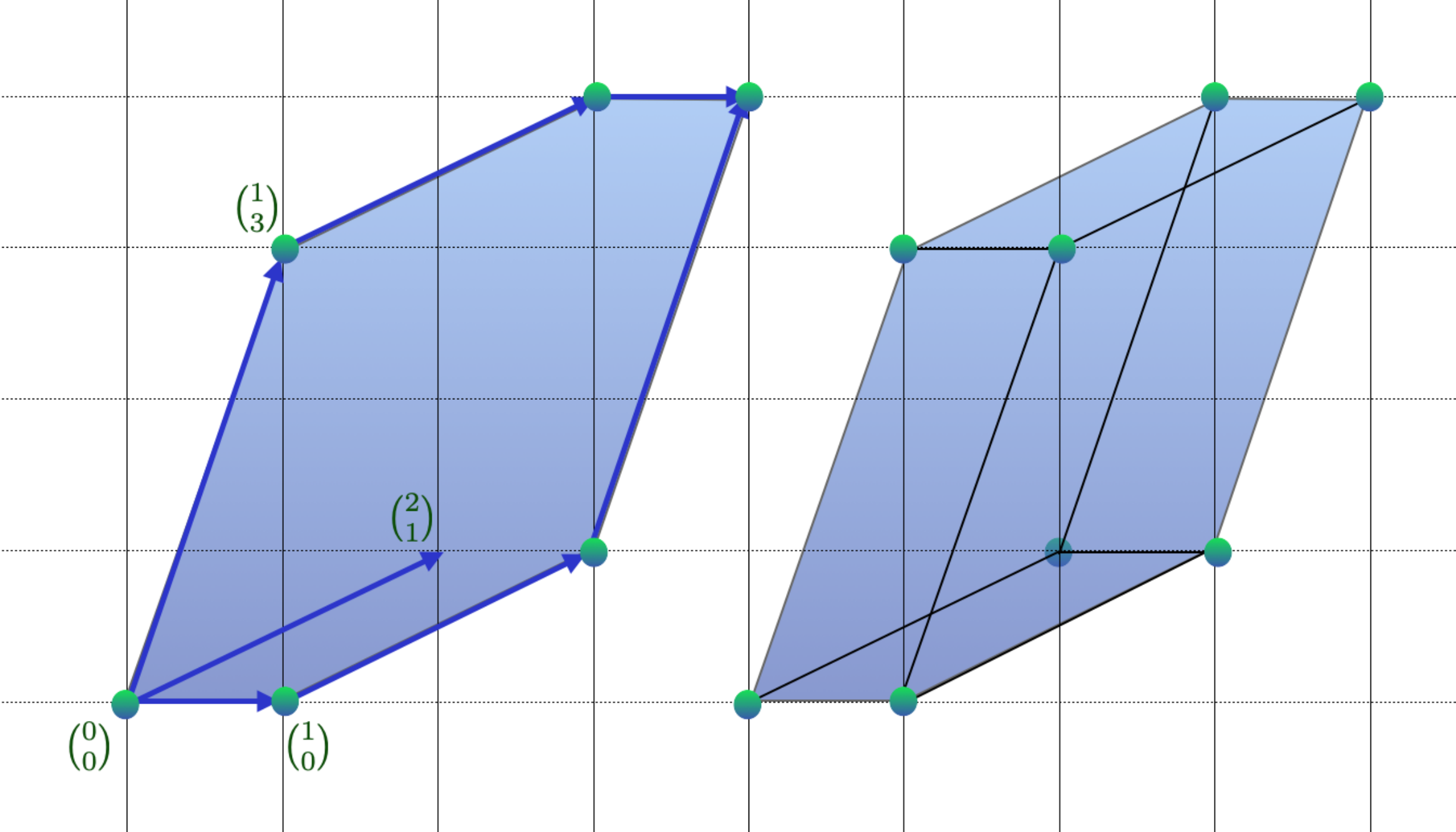}
 \end{center}
\caption{The Minkowski sum of $3$ line segments in the plane, forming a $2$-dimensional zonotope, as described in Example \ref{Minkowski sum of 3 line segments}. }
\label{a zonotope}
\end{figure}

\bigskip
\begin{example}\label{Minkowski sum of 3 line segments}
\rm{
Consider the following $3$ line segments in $\R^2$:  
$\conv\{ \icol{0\\0}, \icol{1 \\0} \}, \conv\{ \icol{0\\0}, \icol{2 \\1} \}$, and 
$\conv\{ \icol{0\\0}, \icol{1 \\ 3} \}$.
The Minkowski sum of these three line segments, by definition a zonotope in $\R^2$, is 
the symmetric hexagon whose vertices are $\icol{0\\0}, \icol{1 \\0},  \icol{2 \\1}, \icol{3 \\3}, \icol{3 \\1}, \icol{4 \\3}$.  
Notice that once we graph it, in Figure \ref{a zonotope}, the graph is hinting to us
 that this body is a projection of a $3$-dimensional cube, and indeed this turns
out to be always true for Minkowski sums of line segments.
}
\hfill $\square$
\end{example}
 \index{Minkowski sum}

\bigskip

\begin{example} \label{truncated octahedron}
\rm{
A particular embedding of the truncated octahedron $\P$, drawn in Figure \ref{TruncatedOctahedron},
 is given by the convex hull of the set of $24$ vertices defined by all permutations of $(0, \pm 1, \pm 2)$.
 We note that this set of vertices can also be thought of as the orbit of just the one point $(0, 1, 2)\in \R^3$ under the hyperoctahedral group (see 
 \cite{BillChen} for more on the hyperoctahedral group).  It turns out that this truncated octahedron $\P$ tiles $\R^3$ by translations with a lattice 
(Exercise  \ref{tiling using the truncated octrahedron}).
}
\hfill $\square$
\end{example}

 \begin{figure}[htb]
 \begin{center}
\includegraphics[totalheight=2in]{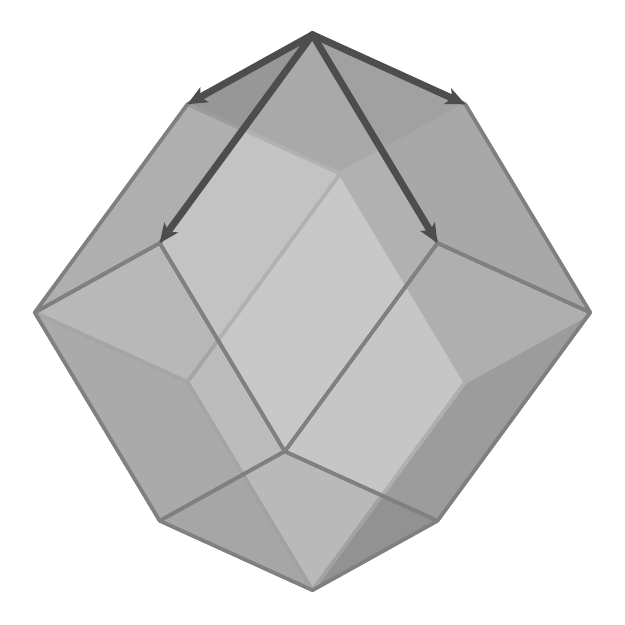}
 \end{center}
\caption{A $3$-dimensional zonotope, called the rhombic dodecahedron, showing in bold its $4$ line segments whose Minkowski sum generate the object. }
\label{first 3d zonotope pic} 
\end{figure}
\index{Minkowski sum}

Next, we show that it's possible to detect whether any compact set $S$ is centrally symmetric,
 by just observing whether its Fourier transform is real-valued.
 \begin{lem}\label{symmetric iff FT is real}
Suppose that  $S \subset \R^d$ is a compact set.   Then:
\[
S \text{ is symmetric about the origin } \iff   \  \hat 1_{S}(\xi) \in \R, 
\]
for all  $\xi \in \R^d$.
\end{lem}
\begin{proof}
 Suppose that $S$ is symmetric about the origin, meaning that $S = -S$.  Then we have:
 \begin{align}
 \overline{ \hat 1_S(\xi)} := \overline{  \int_{S}   e^{2\pi i \langle \xi, x \rangle} dx} 
 &=  \int_{S}   e^{-2\pi i \langle \xi, x \rangle} dx \\
 &=  \int_{-S}   e^{2\pi i \langle \xi, x \rangle} dx \\
 &=  \int_{S}   e^{2\pi i \langle \xi, x \rangle} dx :=   \hat 1_S(\xi), \\
 \end{align}
 showing that the complex conjugate of $\hat 1_S$ is itself, hence that it is real-valued.  

Conversely, suppose that $\hat 1_{S}(\xi) \in \R$, for all $\xi \in \R^d$.  
We must show that $S=-S$.   
We first compute:
\begin{align}
\hat 1_{-S}(\xi) := \int_{-S}   e^{2\pi i \langle \xi, x \rangle} dx     
&= \int_{S}   e^{-2\pi i \langle \xi, x \rangle} dx \\
 &=  \overline{      \int_{S}   e^{2\pi i \langle \xi, y \rangle} dy         }        \\
 &:= \overline{    \hat 1_S(\xi)           }   \\
 &=   \hat 1_S(\xi),
 \end{align}
 for all $\xi \in \R^d$, where we have used the assumption that  $ \hat 1_S(\xi)$ is real-valued in the last equality.  But Theorem \ref{the FT of a convex set determines the set} tells us that in this case:
$\hat 1_{S}(\xi) = \hat 1_{-S}(\xi) \text{ for all } \xi \in \R^d   \   \iff  \  S = -S$.

 \end{proof}

 \begin{figure}[htb]
 \begin{center}
\includegraphics[totalheight=4in]{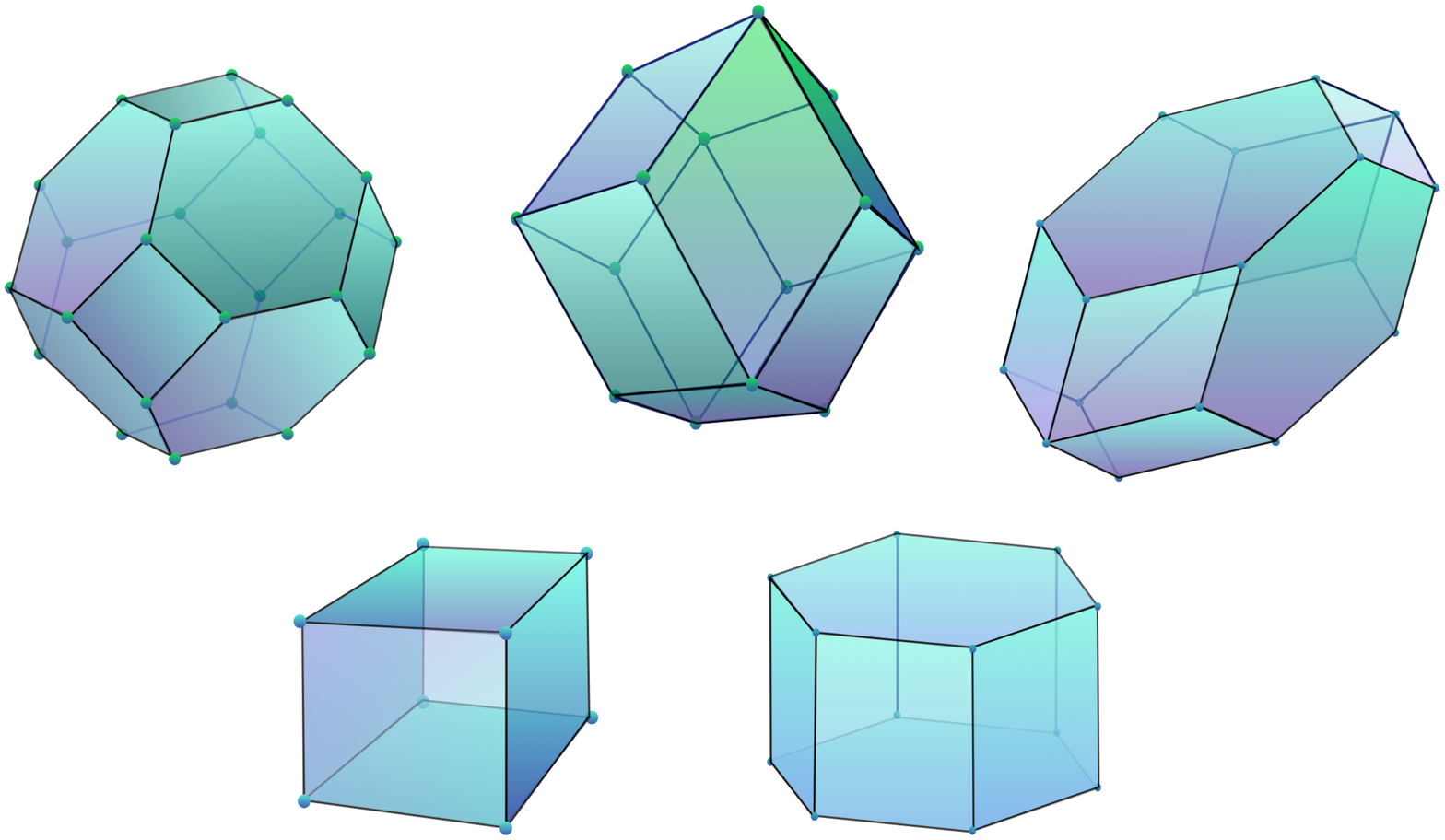}
 \end{center}
\caption{The Fedorov solids, the only $3$-dimensional polytopes that tile $\R^3$ by translations.  All $5$ of them are zonotopes, and they are also extremal bodies for Minkowski's convex body theorem. 
 The top three, from left to right, are:  the Truncated octahedron, the Rhombic dodecahedron,  and the Hexarhombic dodecahedron.  The bottom two are the cube and the hexagonal prism.  }
\label{Fedorov solids}  \index{Fedorov solids}
\end{figure}

 \bigskip
 \begin{example} \rm{
 The interval $\P:= [-\frac{1}{2}, \frac{1}{2}]$ is a symmetric polytope, and indeed we can see that 
 its Fourier transform
 $\hat 1_\P(\xi)$ is real-valued, namely we have $\hat 1_\P(\xi) =  \sinc(\xi)$, as we saw in equation 
 \eqref{SincFunction}.}
 \hfill $\square$
 \end{example}

 \bigskip
  \begin{example} \rm{
 The cross-polytope $\Diamond_2$ is a symmetric polytope, and as we verified in dimension $2$,
  equation \eqref {Fourier transform of 2d crosspolytope},
  its Fourier transform
 $1_{\Diamond_2}(\xi)$ is real-valued. }
 \hfill $\square$
 \end{example}


Alexandrov \cite{Alexandrov}, and independently Shephard \cite{ShephardSymmetricPolytopes}, proved the following remarkable fact.
\begin{thm}[Alexandrov and Shephard]   \label{cs1}   \index{Alexandrov, A. D. }   \index{Shephard}
\label{Alexandrov-Shepard thm}
Let $P$ be any real, $d$-dimensional polytope, with $d \geq 3$.  If all of the facets of  $P$ 
are symmetric, then $P$ is symmetric.
\hfill $\square$
\end{thm}

\begin{example}
\rm{
The converse to the latter result is clearly false, as demonstrated by the cross-polytope in dimension $d > 2$:  it is centrally symmetric, but its facets are not symmetric because they are simplices and we know that 
no simplex (of dimension $\geq 2$) is symmetric
 (Exercise \ref{no simplex is symmetric}).  
 }  
\hfill $\square$
\end{example}

\begin{wrapfigure}{R}{0.49\textwidth}
\centering
\includegraphics[width=0.20\textwidth]{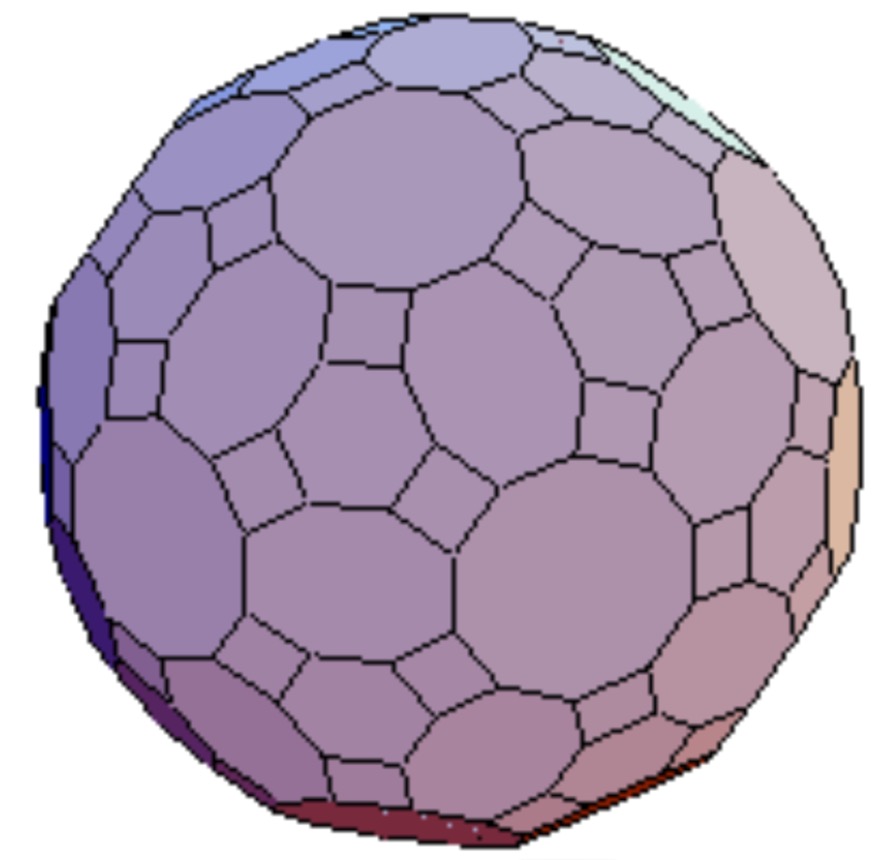}
\caption{A $3$-dimensional  zonotope that does not tile $\R^3$ by translations. }
\label{complex zonotope}
\end{wrapfigure}

Suppose we consider $3$-dimensional polytopes $\P$, and ask which ones enjoy the property that all of their $2$-dimensional faces are symmetric?  Because $1$-dimensional faces are always symmetric, and because Theorem \ref{Alexandrov-Shepard thm} tells us that $\P$ itself must also be symmetric,  the answer is that $\P$ must be a zonotope - in other words all of its faces are symmetric. 

Moving up to $4$-dimensional polytopes, our curiosity might take the next step: 
which $4$-dimensional polytopes enjoy the property that all of their $3$-dimensional faces are symmetric?  Must they also be zonotopes?  The $24$-cell is a good counterexample, because it has triangular $2$-dimensional faces, and hence is not a zonotope.   On the other hand, the $24$-cell  tiles $\R^4$ by translations with a lattice (it is the Voronoi cell of the D$4$ lattice), and therefore by Theorem \ref{Minkowski-Venkov-McMullen} its $3$-dimensional facets must be symmetric. 

What if we ask which $4$-dimensional polytopes enjoy the property that all of their $2$-dimensional faces are symmetric? 
Peter McMullen \cite{McMullen4}
discovered the wonderful conclusion that all of their faces must be symmetric - in other words they must be zonotopes - and that much more is true. 

\begin{thm}[McMullen]   \label{McMullen extension to Alexandrov}
Let $P$ be any real, $d$-dimensional polytope, with $d \geq 4$.
Fix any positive integer $k$ with
 $2 \leq k \leq d-2$. 
 
If  the $k$-dimensional faces of $\P$ are symmetric, then   $\P$ is a zonotope.
\hfill $\square$
\end{thm}

\bigskip

One might wonder what happens if we `discretize the volume' of a symmetric body $K$, by counting integer points, and then ask for an analogue of Minkowski Theorem \ref{Minkowski's convex body theorem, for Z^d}.
 In fact,
Minkowski already had a result about this too
(and he had so many beautiful ideas that it's hard to put them all in one place!).  We give Minkowski's  own elegant and short proof. 

\begin{thm}[Minkowski, 1910]  \label{Minkowski's 3^d theorem}
Let $K\subset \R^d$ be any $d$-dimensional, convex, centrally symmetric set.
If the only integer point in the interior of $K$ is the origin, then 
\begin{equation} \label{Minkowski, $3^d$}
\left | K \cap \Z^d  \right |  \leq 3^d.
\end{equation}
\end{thm}
\begin{proof}
We define the map $\phi: \Z^d \rightarrow \left(   \Z/3\Z  \right)^d$, by reducing each coordinate modulo $3$. 
Now we claim that when restricted to the set $K\cap \Z^d$, our map $\phi$ is $1-1$.   The statement of the theorem follows directly from this claim.  So let $x,y \in K\cap \Z^d$,  and suppose $\phi(x) = \phi(y)$.  Then, by definition of the map $\phi$, we have 
\begin{equation}  \label{z point}
n:= \frac{1}{3}(x-y) \in \Z^d,  
\end{equation}
Now we define $C$ to be the {\bf interior} of the convex hull of $x, -y$, and $0$.   Because $K$ is symmetric, and $x, y\in K$, we know that  $-y \in K$ as well, so that $C \subset \text{int}(K)$.   Now using the convexity of $C$, we also see that  $n \in C$, because $n$ is a non-trivial convex linear combination of $0, x, -y$. 

Therefore  $n \in \text{int}(K)$ as well.  Altogether, $n \in \text{int}(K) \cap \Z^d = \{0\}$, which forces $n =0$.   
Hence $x-y=0$. 
\end{proof}

Theorem \ref{Minkowski's 3^d theorem} is often called { \bf Minkowski's $3^d$ theorem}. \index{Minkowski's $3^d$ theorem}  
An  immediate and natural question is:  which bodies account for the `equality case'?  One direction is easy to see:  if $K$ is the integer cube $[-1, 1]^d$, then it is clear that $K$ is symmetric about the origin, and the only integer point in its interior is the origin.  
In addition, $\vol K = 2^d$, and $K$ contains precisely $3^d$ integer points.   It is a bit surprising, perhaps, that only in 2012 was it proved that this integer cube is the only case of equality in  Minkowski's $3^d$ theorem \cite{DraismaMcAllisterNill}.


\bigskip
\section{Sums of two squares, via Minkowski's theorem}
In $1625$, Albert Girard appears to have been the first to observe (without proof) that if we have 
a prime $p\equiv 1 \pmod 4$, then $p= a^2+b^2$ for some positive integers $a, b$, and that up to order such a representation is unique.  It's easy to see that if a prime $p\equiv 3 \pmod 4$ then it cannot be written as a sum of two integer squares, because every square mod $4$ is congruent to either $0$ or $1$ $\pmod 4$, and hence a sum of two integer squares must be congruent to either $0, 1$, or $2 \pmod 4$.

Fermat popularized this result, which now bears his name, although Fermat
did not provide a proof of this statement.  The first recorded proof was discovered by Euler, in $1752$, and employs the idea of infinite descent. 
Here we'll give a proof of this result by appealing to Minkowski's convex body theorem, namely Theorem \ref{Minkowski convex body Theorem for L}.

To warm-up, the reader may want to solve the following elementary and classical number-theory problem: given a prime $p\equiv 1 \pmod 4$, prove that there exists an integer $m$ such that $-1 \equiv m^2 \pmod p$ 
(Exercise \ref{-1 a quadratic residue mod p}).

\smallskip
\begin{thm}[Sum of two squares]
Let $p$ be an odd prime.   Then $p= a^2 + b^2$ is solvable in integers $a, b \iff p\equiv 1 \pmod 4$.
\end{thm}
\begin{proof}
We know from above that there exists an integer $k$ such that $-1 \equiv k^2 \pmod 4$, and we've also seen that the prime $p$ must satisfy $p\equiv 1 \pmod 4$.  We define the lattice $\L:= M(\Z^2)$, with 
$M:= 
 \big(\begin{smallmatrix}
\ 1 & 0  \\
k & p
\end{smallmatrix}
\big)$,
a lattice which manifestly has determinant $p$.  For our convex body, we'll pick the ball 
$B:= \{ x \in \R^2 \mid \| x\| \leq \sqrt{2p} \}$.  The volume of $B$ is $2p\pi$, and we can now check that the hypotheses of Minkowski's Theorem \ref{Minkowski convex body Theorem for L}
are satisfied:  
\[
 \vol B > 2^2 \det M \iff   (2p)\pi > 4p,
 \]
 which is true.   Hence there exists an integer 
point $\icol{a\\b} \in \L$ in the interior of $B$.   But any point in $\icol{a\\b} \in \L$ must satisfy
 \[
 \icol{a\\b} = 
 \big(\begin{smallmatrix}
\ 1 & 0  \\
k & p
\end{smallmatrix}
\big) \icol{m\\n} =  \icol{m\\mk+np},
\]
for some $m, n \in \Z$.  We now have:
\[
a^2 + b^2 = m^2 + (mk +np)^2 \equiv m^2(1+k^2) \equiv 0 \pmod p,
\]
so $p \mid a^2+b^2$.  Finally, we'll use the fact that $\icol{a\\b}$ is also in the interior of the body $B$, giving us $a^2 + b^2 < 2p$.  Together with $p \mid a^2+b^2$, we arrive at $p =a^2+b^2$.
\end{proof}
This proof shows one small aspect of Minkowski's powerful geometry of numbers, using simple ideas in geometry to conclude nontrivial number-theoretic facts.   Minkowski's Theorem \ref{Minkowski convex body Theorem for L} can also be used to prove Lagrange's theorem, namely that every integer may be written as a sum of $4$ squares (for a proof see \cite{Hardy.Wright.book}, for example).


\bigskip
\section{The volume of the ball, and of the sphere} 
\label{Volume of the ball, the Gamma function}
\index{Volume of the ball} \index{Gamma function}

The most symmetric of all convex bodies is the ball, and here we will explicitly compute the 
volumes of $d$-dimensional balls and the volumes of $(d-1)$-dimensional spheres.  
For these very classical computations, we need the Gamma function:
\begin{equation}
\Gamma(x):=  \int_0^\infty  e^{-t} t^{x-1} dt,
\end{equation}
valid for all $x>0$.    The Gamma function $\Gamma(x)$ interpolates smoothly between the integer values of the factorial function $n!$, in the following sense.
\begin{lem} \label{Gamma properties}
Fix $x>0$.  Then
\begin{enumerate}[(a)]
\item  $\Gamma(x+1) = x \Gamma(x)$.
\item $\Gamma(n+1) = n!$, for all nonnegative integers $n$.
\item   $\Gamma\Big( \frac{1}{2} \Big) = \sqrt \pi$.
\item $\Gamma$ extends to  an infinitely smooth function on the complex plane, except at $0$  and at the
         negative integers, where it has simple poles. 
\end{enumerate}
\end{lem}
The verifications of parts (a), (b), and (c) are good exercises (Exercise \ref{prove Gamma properties}), and we don't want to deprive the reader of that pleasure.  Part (d) requires some knowledge of complex analysis, but we include the statement here for general knowledge. 

What is the volume of the unit ball $B:= \left\{ x \in \R^d  \mid  \|x\| \leq 1 \right\}$?  And what about the volume of the unit sphere $S^{d-1} := \left\{  x\in \R^d \mid \| x \| \right\}= 1 \|$?

\begin{lem} \label{lem:volume of ball and sphere}
For the unit ball $B$, and unit sphere $S^{d-1}$, we have:
\begin{equation} \label{volume of ball and sphere}
\vol B  = \frac{ \pi^{\frac{d}{2}} }{\Gamma\left(\frac{d}{2} +1\right)}, \text{ and  } 
\vol \left(S^{d-1}\right) = \frac{2 \pi^{\frac{d}{2}} }{ \Gamma\left(\frac{d}{2}\right)}.
\end{equation}
\end{lem}
\begin{proof}
We let $\kappa_{d-1}:= \vol(S^{d-1})$ denote the surface area of the unit sphere 
$S^{d-1} \subset \R^d$.   We use polar coordinates in $\R^d$, meaning that we may write each $x \in \R^d$ in the form $x = (r, \theta)$, where $r>0$ and $\theta \in  S^{d-1}$.   Thus $\|x\| = r$, and we also have the calculus fact that $dx = r^{d-1} dr d\theta$.  

 Returning to our Gaussians $e^{-\pi \|x\|^2}$, we may recompute their integrals using polar coordinates in $\R^d$:   
\begin{align*}
1=\int_{\R^d}   e^{-\pi \|x\|^2} dx &= \int_{S^{d-1} } 
             \int_0^{\infty}    e^{-\pi r^2}  r^{d-1} dr \, d\theta \\
&=\kappa_{d-1} \int_0^{\infty}    e^{-\pi r^2}  r^{d-1}  dr \\
&= \kappa_{d-1}   \frac{1}{   2\pi^{\frac{d}{2}}      }
 \int_0^{\infty}    e^{- t}  t^{\frac{d}{2} -1}  dt,
\end{align*}
where we've used $t:= \pi r^2$, implying that 
$r^{d-1}  dr  = r^{d-2}  r dr =
\Big(\frac{t}{\pi} \Big)^{\frac{d-2}{2}} \frac{dt}{2\pi} $.  Recognizing the latter integral as 
$\Gamma\left(\frac{d}{2}\right)$, we find that 
$1= \frac{ \kappa_{d-1} }  { 2\pi^{\frac{d}{2}}  }  \Gamma\left(\frac{d}{2}\right)$, as desired. 
 
For the volume of the unit ball $B$, we have:
\[
\vol B = \int_0^1  \kappa_{d-1} r^{d-1}  dr  =   \frac{\kappa_{d-1}}{d}
=  \frac{\pi^{\frac{d}{2}}}{\frac{d}{2}   \Gamma\left(\frac{d}{2}\right)  }
=  \frac{\pi^{\frac{d}{2}}}{  \Gamma\left(\frac{d}{2} +1\right)  }.
\]
\end{proof}
It is an easy fact, but worth mentioning,   that we may also 
rewrite the formulas \eqref{volume of ball and sphere} in terms of ratios of factorials by using the recursive properties of the $\Gamma$ function (Exercise \ref{explicit volume for ball and sphere}).
While we are at it, let's dilate the unit ball by $r>0$, and recall our definition of the ball of radius $r$:
\[
B_d(r):= \left\{ x \in \R^d  \mid  \|x\| \leq r \right\}.
\]
We know that for any $d$-dimensional body $K$, we have $\vol(rK) = r^d \vol K$, so we also get the volumes of the ball of radius $r$, and the sphere of radius $r$:
\begin{equation} \label{dilated volumes of balls and spheres}
\vol  B_d(r) = \frac{ \pi^{\frac{d}{2}} }{\Gamma\left(\frac{d}{2} +1\right)} r^d , \text{ and  } 
\vol \left(r S^{d-1}\right) = \frac{2 \pi^{\frac{d}{2}} }{ \Gamma\left(\frac{d}{2}\right)} r^{d-1}.
\end{equation}
Intuitively,  the derivative of the volume is the surface area, and now we can confirm this intuition:
\[
\frac{d}{dr} \vol B_d(r) = \frac{ d \pi^{\frac{d}{2}} }{\Gamma\left(\frac{d}{2} +1\right)} r^{d-1}
=\frac{ 2\frac{d}{2}  \pi^{\frac{d}{2}} }{\frac{d}{2} \Gamma\left(\frac{d}{2} \right)} r^{d-1}
=\frac{ 2 \pi^{\frac{d}{2}} }{\Gamma\left(\frac{d}{2}\right)} r^{d-1} =\vol \left(r S^{d-1}\right).
\]

\bigskip
\section{Classical geometric inequalities}

 It turns out that the volume of the difference body 
  \index{symmetrized body}
$\frac{1}{2} K - \frac{1}{2} K$, which appeared quite naturally in some of the proofs above, can be related in a rather precise manner to the volume of $K$ itself.
 The consequence is the following {\bf Rogers-Shephard inequality}:
\begin{equation}  \label{Rogers-Shephard inequality}
\vol K     \leq   \vol \left( \frac{1}{2} K - \frac{1}{2} K \right)     \leq       {2d \choose d} \vol K,
\end{equation}
where equality on the left holds $\iff$ $K$ is a symmetric body, and equality on the right 
holds $\iff$ $K$ is a simplex (see \cite{RogersShephard}, and Cassels \cite{CasselsBook}).  
There is also an extension of the Rogers-Shephard inequality to two distinct convex bodies 
$K, L\subset \R^d$:
\begin{equation}  
   \vol \left( K - L \right)     \vol \left( K \cap L \right)     \leq       {2d \choose d} \vol K \vol L.
\end{equation}
(\cite{RogersShephard} and  \cite{Gutierrez.Jimenez.Villa}).
A quick way of proving \eqref{Rogers-Shephard inequality} is by using the ubiquitous
{\bf Brunn-Minkowski inequality}.  To set it up, two sets $A, B\subset \R^d$  are called {\bf homothetic}
 \index{homothetic} if $A = \lambda B + v$, for some 
 fixed $v\in \R^d$, and some $\lambda >0$ (or either $A$ or $B$ consist of just one point).

\begin{thm}[Brunn-Minkowski inequality]. \index{Brunn-Minkowski inequality}
If $K$ and $L$ are convex subsets of $\R^d$, then 
\begin{equation}
\vol(K+L)^\frac{1}{d} \geq \vol(K)^\frac{1}{d} + \vol(L)^\frac{1}{d},
\end{equation}
with equality if and only if $K$ and $L$ lie in parallel hyperplanes or  are homothetic to each other.
\hfill $\square$
\end{thm}
 (see \cite{Schneider.book}, section $7.1$, for a proof and a thorough introduction to this inequality)


\bigskip
\section{Minkowski's theorems on linear forms}
\index{Minkowski's theorem on linear forms}

There is a quick and wonderful application of Minkowski's first theorem
to products of linear forms.  
\bigskip
\begin{thm} [Minkowski - homogeneous linear forms]
\label{Minkowski's thm on linear forms, I}
For each $1\leq i \leq d$, let 
 \[
 L_i(x):= a_{i, 1}x_1 + \cdots + a_{i, d}x_d
 \]
  be  linear forms with real coefficients $a_{i, j}$, and suppose that the matrix $A$
  formed by these coefficients $a_{i, j}$ is invertible.  Suppose further that there exists positive numbers $\lambda_1, \dots, \lambda_d$ with the property that $\lambda_1 \lambda_2  \dots, \lambda_d \geq |\det A|$.
  
 Then  there exists a nonzero integer vector $n\in \Z^d$ such that 
 \begin{equation}
 \left| L_1(n) \right | \leq \lambda_1, \cdots,  \left| L_d(n) \right | \leq \lambda_d.
 \end{equation}
\end{thm}
\begin{proof}
We define the body
\begin{equation}
\P:= \{ x \in \R^d \mid |L_k(x) | \leq \lambda_k, \text{ for each } 1\leq k \leq d\},
\end{equation}
which is a centrally-symmetric parallelepiped.  To compute $\vol \P$, we 
note that the image of $\P$ under the linear transformation $A$ is 
$Q:= A(\P) = \{ x \in \R^d \mid |x_k | \leq \lambda_k, 1\leq k \leq d\}$, which clearly has volume 
$\vol Q = 2^d \lambda_1 \cdots \lambda_d$.  Therefore 
\[
\vol \P = \vol A^{-1} Q = \frac{1}{\det A} \vol Q = 
\frac{1}{\det A}  2^d \lambda_1 \cdots \lambda_d \geq 2^d,
\]
the last inequality holding by assumption.  By Minkowski's Theorem \ref{Minkowski convex body Theorem for L}, $\P$ contains a nonzero integer point, and we're done.
\end{proof}

\begin{cor}[Minkowski - product theorem for homogeneous linear forms]
\label{Minkowski's thm on linear forms, II}
For each $1\leq i \leq d$, let 
 \begin{equation}
 L_i(x):= a_{i, 1}x_1 + \cdots + a_{i, d}x_d
 \end{equation}
  be  linear forms with real coefficients $a_{i, j}$, and suppose that the matrix $A$
  formed by these coefficients $a_{i, j}$ is invertible.   If $d>1$, then there exists
a nonzero integer vector $n\in \Z^d$ such that 
 \begin{equation} \label{inequality for product of linear forms}
 \left| L_1(n) L_2(n) \cdots L_d(n) \right | \leq |\det A|.
 \end{equation}
 \end{cor}
\begin{proof}
We can simply use Theorem \ref{Minkowski's thm on linear forms, I} with
$\lambda_1 = \lambda_2 = \cdots = \lambda_d := |\det A|^\frac{1}{d}$.  Since
$\lambda_1 \cdots \lambda_d = |\det A| $, the conclusion of Theorem  \ref{Minkowski's thm on linear forms, I} gives us the existence of a nonzero integer point $n$ that satisfies the following:
\[
 \left| L_1(n) \right | \left | L_2(n) \right | \cdots \left | L_d(n) \right | 
 \leq \lambda_1 \cdots \lambda_d = |\det A|.
\] 
\end{proof}

It's worth mentioning that there are various ways to strengthen 
Corollary \ref{Minkowski's thm on linear forms, II}.  For example, it is possible to replace 
the inequality in \eqref{inequality for product of linear forms} by a strict inequality
 \cite{HlawkaSchoissengeierTaschner}.


\bigskip
\section{Poisson summation as the trace of a compact linear operator}    
\label{section:Poisson summation as a trace}
\index{Compact linear operators and Poisson summation}

Now that we've seen a few applications of Poisson summation (and there will more throughout the book), it's natural to wonder if there is something a little deeper going on here.  In this brief section we carry the reader through a more general context for Poisson summation, as the trace of a certain linear operator. The modern context for this extension is called the spectral theory of compact operators.  For more about the spectral theory of noncompact operators as well, the reader is invited to peruse Audrey Terras' book \cite{Terras.Harmonic1}.

Suppose we are given a compact set $\P\subset \R^d$ of positive $d$-dimensional volume, and a
continuous function $K(x, y): \P \times \P \rightarrow \C$.  Then we can define a corresponding operator $T_K: L^2(\P)\rightarrow L^2(\P)$ by
\begin{equation}
T_K(f)(x) := \int_{\P} K(x, y) f(y) dy.
\end{equation}
The function $K(x, y)$ is called a {\bf kernel}. 
The operator $T_K$ is clearly linear, and indeed $T_K(\alpha f+\beta g)~=~\alpha~T_K(f)~+~\beta~T_K(g)$ follows from the 
linearity of the integral.
We call  $T_K$ a {\bf positive operator} if
$\langle T_K(f) , f \rangle > 0$ for all nonzero functions $f$.  Finally, the kernel (as well as the operator) is called
 {\bf self-adjoint} if $K(x, y) = \overline{ K(y, x) }$ for all $x, y \in \P$.   A standard fact is that all of the eigenvalues of $T_K$ are real. By the spectral theorem for compact, self-adjoint linear operators (see \cite{EinsiedlerWardBook}), we know that $T_K$ has an orthonormal basis of eigenvectors 
 $\{v_1, v_2,  v_3, \dots \}$, which correspond  to its nonzero eigenvalues 
$\{     \lambda_1, \lambda_2,  \lambda_3,   \dots \}$.

James Mercer proved the following useful theorem \cite{Mercer}.
\begin{thm}[Mercer, 1909]
Suppose that $T_K$ is a positive, self-adjoint operator on a compact set $\P\subset\R^d$. 
Then:
\begin{equation}
K(x, y) = \sum_{n=1}^\infty \lambda_n  v_n(x) \overline{ v_n(y)},
\end{equation}
and the series converges absolutely and uniformly.
\hfill $\square$
\end{thm}

The {\bf trace} of the linear operator $T_K$ is defined by 
${ \rm Trace}(T_K):=\sum_{n=1}^\infty \lambda_n$.  If $T_K$ satisfies the hypotheses of Mercer's theorem, then we have also have the following immediate Corollary:
\begin{equation} \label{Cor of Mercer}
 \int_{\P} K(x, x) dx =  \int_{\P}  \sum_{n=1}^\infty \lambda_n  v_n(x) \overline{ v_n(x)} dx
= \sum_{n=1}^\infty \lambda_n  \int_{\P} \left|v_n(x)\right|^2 dx =  \sum_{n=1}^\infty \lambda_n, 
\end{equation}
which is the trace of $T_K$. 

So what does all of this abstraction
 have to do with Poisson summation, the reader might ask? Well, let's pick
 $\P:= \torus$, the $d$-dimensional torus, and let's fix a Schwartz function $f:\R^d \rightarrow \C$.  We may now consider the linear operator defined by
\begin{equation}
L_f(g)(x) := (f*g)(x) := \int_{\R^d} f(x-y) g(y) dy,
\end{equation}
for all $x \in \R^d$, and for all $ g \in L^2(\torus)$.  
The interplay between the torus and $\R^d$ is intended here, and in fact we have:
\begin{align*}
L_f(g)(x) &:=  \int_{\R^d} f(x-y) g(y) dy 
= \sum_{n \in \Z^d} \int_{\left( \R^d/\Z^d \right)  -n} f(x-y) g(y) dy \\
&= \sum_{n \in \Z^d} \int_\torus  f(x-y+n) g(y) dy \\
&=\int_\torus   \left( \sum_{n \in \Z^d}  f(x-y+n) \right) g(y) dy \\
&:=\int_\torus  K(x, y) g(y) dy,
\end{align*}
where we've defined our kernel $K(x, y):=  \sum\limits_{n \in \Z^d}  f(x-y+n)$ in the last equality above.   What are the eigenfunctions of $L_f$? We claim that they are precisely the exponentials $e_k(x):= e^{2\pi i \langle x, k\rangle}$, indexed by $k \in \Z^d$!  We can compute:
\begin{equation}
L_f( e_k)(x):= \int_{\R^d} f(y) e^{2\pi i \langle x-y, k\rangle}dy 
 = e^{2\pi i \langle x, k\rangle}     \int_{\R^d} f(y) e^{-2\pi i \langle y, k\rangle} dy
 =\hat f(k) e_k(x),
\end{equation}
proving that each function $e_k(x)$ is an eigenfunction of $L_f$, with eigenvalue $\hat f(k)$. Using the completeness of this set of orthonormal exponentials $\{e_k(x) \mid k \in \Z^d\}$ in the Hilbert space $L^2(\torus)$, it's also possible to show that these are all of the eigenfunctions.   So we see that the
 trace of $L_f$ equals 
 \begin{equation}\label{first trace formula}
{ \rm Trace}(L_f)  := \sum_{n=1}^\infty \lambda_n = \sum_{k\in \Z^d} \hat f(k).
\end{equation}
On the other hand, if we assume that $L_f$ is a self-adjoint positive operator (for this particular $f$), then
\eqref{Cor of Mercer} tells us that the trace may also be computed in another way:
\begin{equation}\label{second trace formula}
{ \rm Trace}(L_f) =  \int_{\torus} K(x, x) dx :=    \int_\torus    \sum_{n \in \Z^d}  f(n) dx
=   \sum_{n \in \Z^d}  f(n)  \int_\torus dx = \sum_{n \in \Z^d}  f(n).
\end{equation}
So we've arrived at the Poisson summation formula (for Schwartz functions)
\[
\sum_{n \in \Z^d}  f(n)  = \sum_{k\in \Z^d} \hat f(k)
\]
 by computing the trace of the linear operator $L_f$, acting on the Hilbert space $L^2(\torus)$.





\bigskip

\section*{Notes}
\begin{enumerate}[(a)]

\item Siegel's original proof of Theorem  \ref{Siegel}   used Parseval's identity, but the ``Fourier-spirit'' of the two proofs is similar.

\item Minkowski's book \cite{Minkowski} in $1896$ was the first treatise to develop the threads between convex geometry, Diophantine approximation, and the theory of quadratic forms.  This book marked the birth of the geometry of numbers.

\item In Exercise \ref{equivalent statements for unimodular triangles} below, we see three equivalent conditions for a $2$-simplex to be unimodular.   In higher dimensions, a $d$-simplex will not satisfy all
three conditions, and hence this exercise shows one important `breaking point' between $2$-dimensional and $3$-dimensional discrete geometry.

\item There is a very important tool in number theory, called the Selberg trace formula, which extends Poisson summation to hyperbolic space.  See, for example, Audrey Terras' book \cite{Terras.Harmonic1}. 

\item The Poisson summation formula also extends to all locally compact abelian groups, and this field has a vast  literature - see, for example \cite{Terras.Harmonic1}.

\item  \label{new books, geometry of numbers}
There are a growing number of interesting books on the geometry of numbers.   An excellent
encyclopedic text  is  Gruber and Lekkerkerker's \cite{GruberLekkerkerker} ``Geometry of Numbers".
Another encyclopedic reference is Peter Gruber's own book \cite{GruberBook}. 

 Two other excellent and classic introductions are Siegel's book \cite{SiegelBook}, and Cassels' book \cite{CasselsBook}.  An expository introduction to some of the elements of the Geometry of numbers, at a level that is even appropriate for high school students, is given by Olds, Lax, and Davidoff  \cite{OldsBook}.   
For upcoming books, the reader may also consult Martin Henk's lecture notes `Introduction to geometry of numbers' \cite{Henk3}, and 
the book by Lenny Fukshansky and Stephan Ramon Garcia, `Geometry of Numbers'  \cite{FukshanskyBook}. 

\item   \label{Brunn-Minkowski}  The Brunn-Minkowski inequality is fundamental to many branches of mathematics, including the geometry of numbers.   A wonderful and encyclopedic treatment of the 
Brunn-Minkowski inequality, with its many interconnections, appears in Rolf Schneider's book ``The Brunn-Minkowski theory" \cite{Schneider.book}.

\item  \label{Fedorov Note}
The Fedorov solids are depicted, and explained via the modern ideas of Conway and Sloane, in an excellent
expository article by David Austin \cite{DavidAustin}. 
For a view into  the life and work of Evgraf Stepanovich Fedorov, as well as
a fascinating account of how Fedorov himself thought about the $5$ parallelohedra, 
 the reader may 
consult the article by Marjorie Senechal and R. V. Galiulin  \cite{SenechalGaliulin}.  The authors of
 \cite{SenechalGaliulin} also discuss the original book of Fedorov, called
\emph{An Introduction to the Theory of Figures}, published in 1885, which is now considered a pinnacle of modern crystallography. 
Fedorov later became one of the great crystallographers of his time.   

In $\R^4$, it is known that there are $52$ different combinatorial types of $4$-dimensional parallelohedra.   In $\R^5$, the complete classification of all the combinatorial types of 
$5$-dimensional paralellohedra was completed in 2016 \cite{Dutour et al.parallelohedra}, where the authors found $110, 244$ of them.

\item The field of multi-tiling is still rapdily growing.  One of the first important papers in this field was by Mihalis Koloutzakis \cite{Kolountzakis1}, who related the multi-tiling problem to a famous technique known as the idempotent theorem, and thereby proved
that if we have a multi-tiling in $\R^2$ with any discrete set of translations, then we also have a multi-tiling with a finite union of lattices. 
A  recent advance is an equivalence between multi-tiling and 
certain Hadwiger-type invariants, given by Nir Lev and Bochen Liu \cite{LevLiu}.  Here the authors show as well that
for a generalized polytope $\P \subset \R^d$   (not necessarily convex or connected), if $\P$ is spectral, then $\P$ 
is equidecomposable by translations to a cube of equal volume.  

Another natural question in multi-tiling, which is still open, is the following:
\begin{question} \label{multi-tiling - what is the discrete set of translations}
Suppose that $\P$ multi-tiles with a discrete set of translations $D$.  Do we really need the set $D$ of translates of $\P$ to be a very complicated discrete set,
or is it true that just a  finite union of lattices suffices?
Even better, perhaps one lattice always suffices?
\end{question}
   In this direction, Liu proved recently that if we assume that $\P$ multi-tiles with a finite union of lattices, 
   then $\P$ also multi-tiles with a single lattice \cite{Liu}.  This is big step in the direction of answering Question 
   \ref{multi-tiling - what is the discrete set of translations} in general.  An earlier, and smaller step, was taken in \cite{GravinKolountzakisRobinsShiryaev}, where the authors answered part of Question \ref{multi-tiling - what is the discrete set of translations} in $\R^3$, reducing the search from an arbitrary discrete set of translations, to translations by a finite union of lattices.  Taken together, the latter two steps imply that in $\R^3$ (and in $\R^2$), any multi-tiling with a discrete set of translations also occurs with just a one lattice.  
   
   
In a different direction, the work of Gennadiy Averkov  \cite{Averkov} analyzes the equality cases for an extension of Minkowski's theorem, relating those extremal bodies to multi-tilers. 
 In \cite{YangZong}, Qi Yang and Chuanming Zong show that
 the smallest $k$ for which we can obtain a nontrivial $k$-tiling in $\R^2$ is $k=5$, and the authors characterize
  those $5$-tiling bodies, showing in particular  that if a convex polygon is a $5$-tiler, then it must be either an octagon, or a decagon. 
  In \cite{Zong2022}, Zong and his collaborators continue the latter research to show that the smallest $k$ for which we can obtain a nontrivial $k$-tiling in $\R^3$ is $k=5$.  These investigations naturally lead to the general question:
  
\begin{question}\label{smallest k}
   In $\R^d$ (for $d \geq 4$), what is the smallest integer $k$ such that there exists a $d$-dimensional polytope 
   $\P$ that $k$-tiles  $\R^d$ (nontrivially)  by translations?
\end{question}

\item \label{Fuglede conjecture}
We say that a body $\P$ (any compact subset of $\R^d$) 
 is `spectral' if the function space $L^2(\P)$ possesses an orthonormal, complete basis of exponentials. 
There is a fascinating and vast literature about such  spectral bodies, relating them to tiling, 
\index{tiling}
and multi-tiling problems. 
One of the most interesting and natural questions in this direction is the following conjecture, by Bent Fuglede \cite{Fuglede74}. 

The Fuglede conjecture asks whether the following is true.
 \begin{question}\label{Fuglede}
$\P$ tiles $\R^d$ by translations $\iff$ $\P$ is spectral?
\end{question}
Terry Tao disproved the Fuglede conjecture for some nonconvex bodies.  Indeed,  in 2003 
Alex Iosevich, Nets Katz, and Terry Tao \cite{IosevichKatzTao} proved that 
the Fuglede conjecture is true for all convex domains in $\R^2$. 
In 2021, this conjecture was proved for all convex domains (which must necessarily be polytopes by
an additional simple argument), in the fundamental work of  Nir Lev and M\'at\'e Matolcsi \cite{LevMatolcsi}.

In a related direction, Sigrid Grepstad and Nir Lev \cite{GrepstadLev} showed that for any bounded, measurable subset $S\subset \R^d$, if $S$ multi-tiles by translations with a discrete set, then
 $S$ has a Riesz basis of exponentials.

 \item  We have seen  
  that the zero set of the Fourier transform of a polytope is very important, in that Theorem    \ref{zero set of the FT of a polytope}
  gave us a necessary and sufficient condition for multi-tiling.   But the zero set of the FT also gives more information, and an interesting application of the information content in the zero set is the Pompeiu problem. \index{Pompeiu problem}
 The Pompeiu problem is an ancient problem (defined in 1929 by Pompeiu) that asks the following:   
 which bodies $\P \in \R^d$ are uniquely characterized by the collection of their integrals over $\P$, and over all rigid motions of $\P$? 
 An equivalent formulation is the following.
 \begin{question} \label{Pompeiu conjecture}
    Given a body $\P$ with nonempty interior, does there exist a nonzero continuous function $f$ that allows for the 
the vanishing of all of the integrals
\begin{equation}\label{Pompeiu question}
    \int_{M(\P)} f(x) dx = 0,
\end{equation}
taken over all rigid motions $M$, including translations? 
 \end{question}
A body $\P \subset \R^d$, for which the answer to the question above is affirmative, is said to have the Pompeiu property.  

Even for convex bodies $\P$, it is still an open problem in general dimension whether 
$\P$ has the Pompeiu property.
It is known, by the work of Brown, Schreiber, and Taylor  \cite{BrownSchreiberTaylor} that $\P$ has the Pompeiu property $\iff$ 
the collection of Fourier transforms  $\hat 1_{\sigma(\P)}(z)$, taken over all rigid motions 
$\sigma$ of $\R^d$,  
have  a common zero $z$.  It was also known that all polytopes have the Pompeiu property.  
Recently, in \cite{FabricioSinai1}, Fabricio Machado and SR showed that the zero set of the FT does not contain (almost all) circles whose center is the origin, and as a consequence we get a simple new proof that all polytopes have the `Pompeiu property'.   
\end{enumerate}

\newpage
\section*{Exercises}
\addcontentsline{toc}{section}{Exercises}
\markright{Exercises}

\begin{quote}    
``Every problem has a creative solution''.

--  Folklore
 \end{quote}

\begin{quote}    
``Every problem has a  solution that is simple, neat, and wrong''.  

--  Mark Twain
 \end{quote}
\medskip

\medskip
\begin{prob}\label{easy first problem}
\rm{
Suppose that in $\R^2$, we are given a symmetric, convex body $K$ of area $4$, which contains only the origin. 
Prove that $K$ must tile $\R^2$ by translations.
}
\end{prob}

\medskip
\begin{prob}   \label{convexity of K-K}  $\clubsuit$
\rm{
Given $d$-dimensional compact, convex  sets $K, L \subset \R^d$, prove that $K+L$ is convex, and that $K-L$ is convex.
}
\end{prob}

\medskip
\begin{prob}  \label{practice with Minkowski sums}
\rm{
Given $d$-dimensional compact, convex  sets $A, B \subset \R^d$, prove that:
\[
A \cap B \subseteq \tfrac{1}{2} A + \tfrac{1}{2} B \subseteq \conv\left( A \cup B \right),
\]
and show that equality holds in either of the two containments $\iff A=B$. 
}
\end{prob}

\medskip
\begin{prob}  \label{more practice with Minkowski sums}
\label{Convexity: A+A = 2A}
\rm{
It is easy to see that essentially by definition, $A \subset \R^d$ is convex $\iff  A+A \subset 2A$.
\begin{enumerate}[(a)]
\item   Given any convex subset $A \subset \R^d$, prove that: 
\[
A + A = 2A.
\]
\item Find a counter-example to show that the converse is false;  in other words, it is false 
that $A+A=2A \implies A$ is convex.

\end{enumerate}
}
\end{prob}

\medskip
\begin{prob} \label{c.s. C equals its symmetrized body}
$\clubsuit$
\rm{
Suppose initially  that  $C \subset \R^d$ is any set.  
\begin{enumerate}[(a)]
\item  Show that 
\begin{equation}
  \frac{1}{2}C - \frac{1}{2}C = C
    \  \implies \text{ $C$ is centrally symmetric}. 
\end{equation}
\item
Show that 
\begin{equation} \label{equivalence for a convex set}
C  \text{ is centrally symmetric and convex }  \iff  \frac{1}{2} C - \frac{1}{2} C = C.
\end{equation}
\item  Find an example of a centrally symmetric set $C$ that is not convex, and satisfies  
\[
 \frac{1}{2} C - \frac{1}{2} C \not= C. 
 \]
\end{enumerate}
}
\end{prob}

\medskip
\begin{prob} \label{convexity interacts with Minkowski sums}
$\clubsuit$
\rm{
Given any convex sets $A, B \subset \R^d$, show that
\[
\conv(A+B) = \conv A + \conv B.
\]
}
\end{prob}

\medskip
\begin{prob} $\clubsuit$ \label{equivalent statements for unimodular triangles}
\rm{
Suppose we have a triangle $\Delta$ whose vertices $v_1, v_2, v_3$ are integer points.  Prove that the following 
 properties are equivalent:
\begin{enumerate}[(a)]
\item  $\Delta$ has no other integer points inside or on its boundary (besides its vertices).  
\item  $Area(\Delta) = \frac{1}{2}$.
\item  $\Delta$ is a unimodular triangle, which in this case means that  $v_3 - v_1$ and $v_2- v_1$ form a basis for $\Z^2$.
\end{enumerate} 
(Hint: You might begin by ``doubling'' the triangle to form a parallelogram.)
}
\end{prob}

\medskip
\begin{prob} \label{exercise:unimodular simplex}
\rm{
Show that in $\R^d$, an integer simplex $\Delta$ is unimodular $\iff \vol \Delta =  \frac{1}{d!}$.
}
\end{prob}

\medskip
\begin{prob}  \label{non-unimodular but empty simplex}
\rm{
In $\R^3$, find an integer simplex $\Delta$ that has no other integer points inside or on its boundary (other than its vertices of course),  but
such that $\Delta$  is not a unimodular simplex.  
}
\end{prob}

\medskip
\begin{prob}\label{support of convolution}  $\clubsuit$
\rm{
Recalling the definition of the support of a function $f$ from \eqref{def of support}, show that:
 \begin{enumerate}[(a)]
 \item
Suppose that we are given two closed, convex bodies $A, B \subset \R^d$.   Show that
\[
\supp ( 1_A * 1_B) =  A + B,
\]
where the addition is the Minkowski addition of sets.
\index{Minkowski sum}
\item
More generally, if two functions $f, g:\R^d \rightarrow \C$ are compactly supported,
show that
\[
 \supp(f*g)\subseteq  \closure\left(   \supp(f) +  \supp(g)  \right),
\]
the closure of the Minkowski sum of their individual supports.
\end{enumerate}
}
\end{prob}

Notes.  For a vast generalization, see Theorem \ref{thm:Lions}.

\medskip
\begin{prob}
\rm{
Let $K\subset \R^d$ be a convex body, and let its $(d-1)$-dimensional boundary be denoted by $B$.  Show that we have the equality of Minkowski sums:
\[
B + B = K + K. 
\]
}
\end{prob}

\medskip
\begin{prob}  \label{FT of a polytope is not Schwartz}
\rm{
Prove that for any polytope $\P$, $\hat 1_{\P}$ is  not a Schwartz function.
}
\end{prob}

\medskip
\begin{prob} \label{convolution of indicators is a nice function}
$\clubsuit$
\rm{
(hard-ish)  Show that if $K$ is any convex body, then $1_K*1_{-K}$ is a nice function, in the sense of \eqref{nice functions}.
In other words, show that the Poisson summation formula holds for the function $f(x):= \left( 1_K*1_{-K} \right)(x)$.

Hint.  Use the Parseval identity, valid for functions $f \in L^2(\R^d)$.  For this particular exercise, feel free to use the results of all of the later sections (though in general we refrain from such a `look ahead').
}
\end{prob}


\medskip
\begin{prob}   \label{Cantor set} \index{Cantor set}
\rm{
We first define the following sets recursively:
 \[
 C_0 := [0, 1],  \   C_1 := [0, \tfrac{1}{3}] \cup [\tfrac{2}{3}, 1], \dots ,
 C_n:= \tfrac{1}{3} C_{n-1} \cup  \left\{ \tfrac{1}{3} C_{n-1} + \tfrac{2}{3} \right\},
 \] 
and now the classical {\bf Cantor set} is defined by their infinite intersection:
\[
\mathcal C:= \cap_{n=0}^\infty C_n.
\]
It is a standard fact (which you may assume here) that the Cantor set $\mathcal C$ is compact, uncountable, and has measure $0$. Despite these facts, show that its difference body satisfies the somewhat surprising identity:
\[
\mathcal C-\mathcal C = [-1, 1].
\]
}
\end{prob} 
Notes.   There is a nice article \cite{Kraft} about such difference sets, written for undergraduates.


\medskip
\begin{prob}  
\rm{
Show that any regular hexagon in the plane cannot tile by translations with the integer lattice $\Z^2$. 
}
\end{prob}

\medskip
\begin{prob}   \label{tiling using the truncated octrahedron}
\rm{
Show that  the truncated octahedron,  defined in Example  \ref{truncated octahedron},  tiles $\R^3$ by using only 
 translations with a lattice.  Which lattice can you use for this tiling? 
 }
\end{prob} 


\medskip
\begin{prob} \label{an application of Cauchy-Schwartz 1}
\rm{
Define $f(x):= a \sin x + b \cos x$, for constants $a,b\in \R$.
Show that the maximum value of $f$ is $ \sqrt{ a^ 2 + b ^2 } $, and occurs
when $\tan x = \frac{ a}{b}$. 
}
\end{prob}


\medskip
\begin{prob}
\rm{
Find an example of a symmetric polygon $\P \subset \R^2$ that multi-tiles (nontrivially) with multiplicity $k = 5$.

Notes.  A trivial multi-tiling for $\P$ is by definition a multi-tiling that uses $\P$, with some multiplicity $k>1$, but such that
 there also exists a $1$-tiling (classical) using the same $\P$ (but perhaps using a different lattice).
 }
\end{prob} 


\medskip
\begin{prob}  \label{elementary condition on the covariogram}
\rm{
Let $K\subset \R^d$ be a convex and centrally symmetric set.   Show that 
\begin{equation*}
\frac{1}{2}K \cap \left(     \frac{1}{2}K +   n    \right) \not=  \phi \iff n \in K.
\end{equation*}
}
\end{prob}


\medskip
\begin{prob}  \label{explicit volume for ball and sphere}
\rm{
Using Lemma \ref{lem:volume of ball and sphere}, 
 show that for the unit ball $B$ and unit sphere $S^{d-1}$ in $\R^d$, we have:
\begin{enumerate}[(a)]
\item
\[
\vol S^{d-1}  =  
\begin{cases}
 \frac{ \left(2 \pi\right)^{\frac{d}{2}} }{  2\cdot 4 \cdot 6   \cdots (d-2)          }, &  \text{if } d  \text{ is even},   \\
   \frac{ 2 \left(2 \pi\right)^{\frac{d-1}{2}} }{  1\cdot 3 \cdot 5   \cdots (d-2)          },                          &   \text{if }  d \text{ is odd}.
 \end{cases}
\]
\item
\[
\vol B =  
\begin{cases}
 \frac{ \left(2 \pi\right)^{\frac{d}{2}} }{  2\cdot 4 \cdot 6   \cdots d          }, &  \text{if } d  \text{ is even},   \\
   \frac{ 2 \left(2 \pi\right)^{\frac{d-1}{2}} }{  1\cdot 3 \cdot 5   \cdots d        },                          &   \text{if }  d \text{ is odd}.
 \end{cases}
 \]
\end{enumerate}
}
\end{prob}

\medskip
\begin{prob} $\clubsuit$
\label{-1 a quadratic residue mod p}
\rm{Suppose we are given a prime $p\equiv 1 \pmod 4$.  Prove that there exists an integer $m$ such that $-1 \equiv m^2 \pmod p$.
}

(Hint: you can assume ``Euler's little theorem'':  
\[
a^{\phi(n)} \equiv 1 \pmod n, \text{  for all coprime integers } a, n,
\]
where $\phi(n)$ is the Euler $\phi$-function.)
\end{prob}


\medskip
\begin{prob} \label{Extending Minkowski to nonconvex bodies}
$\clubsuit$
\rm{
Here we use Siegel's theorem \ref{Siegel for general lattices} to give the following extension of 
Minkowski's classical Theorem  \ref{Minkowski convex body Theorem for L}, but 
for bodies $K$ that are not necessarily symmetric, nor necessarily convex.  

Namely, let $K$ be any bounded, measurable subset of $\R^d$, with positive $d$-dimensional volume.
 Let $B:= \frac{1}{2}K - \frac{1}{2}K$ be the symmetrized  body of $K$ (hence $B$ is a centrally symmetric set containing the origin). 
Let $\L$ be a (full rank) lattice in $\R^d$.  Prove the following statement:
\begin{equation*} 
\text{ If }   \vol K > 2^d (\det \L),  \text{  then }   B
\text{  must contain a nonzero point of } \L \text{ in its interior}. 
\end{equation*}

Notes.  We note that the positive conclusion of the existence of a nonzero integer point holds only for the symmetrized body $B$, 
with no guarantees for any integer points in $K$.  
}
\end{prob} 


\medskip
\begin{prob}  $\clubsuit$  \label{prove Gamma properties}
\rm{
Prove the elementary properties of the $\Gamma$ function (parts (a), (b), and (c)), in Lemma \ref{Gamma properties}.
}
 \end{prob}


\begin{prob}    \label{Mobius meets Poisson} 
{\rm
Using the M\"obius $\mu$-function, defined in problem \ref{SumOfPrimitiveRootsOfUnity}, prove the following inversion formula for infinite series.
\[
\text{If } g(x):= \sum_{n=1}^\infty f(mx), \text{ then } f(x) = \sum_{m=1}^\infty \mu(m) g(mx).
\]
To make everything completely rigorous, can you formulate a sufficient convergence criterion for $f$ in order to make the latter statement true?
}
\end{prob}


\medskip
\begin{prob}  
\rm{
We recall that a polytope $\P\subset \R^d$ is symmetric if there exists a vector $v \in \R^d$ such that
$\P - v$ is symmetric about the origin.  Prove that the following are equivalent:
\begin{enumerate}[(a)]
\item $\P$ is symmetric.
\item There exists a vector $v\in \R^d$ such that for all $\xi \in \R^d$ we have:
\[
\hat 1_\P(\xi) \, e^{2\pi i \langle \xi, v\rangle} \in \R.
\]
\end{enumerate}
}
 \end{prob}




\chapter{
\blue{
An introduction to Euclidean lattices}
}
 \label{chapter.lattices}
\index{lattice}  \index{integer lattice}

\begin{quote}   
``Lattices quantify the idea of periodic structures.''                  

-- Anonymous
\end{quote}

\begin{quote}     
``Less is more........more or less.''
   
 -- Ludwig Mies van der Rohe
\end{quote}

\begin{figure}[htb]
 \begin{center}
\includegraphics[totalheight=2.6in]{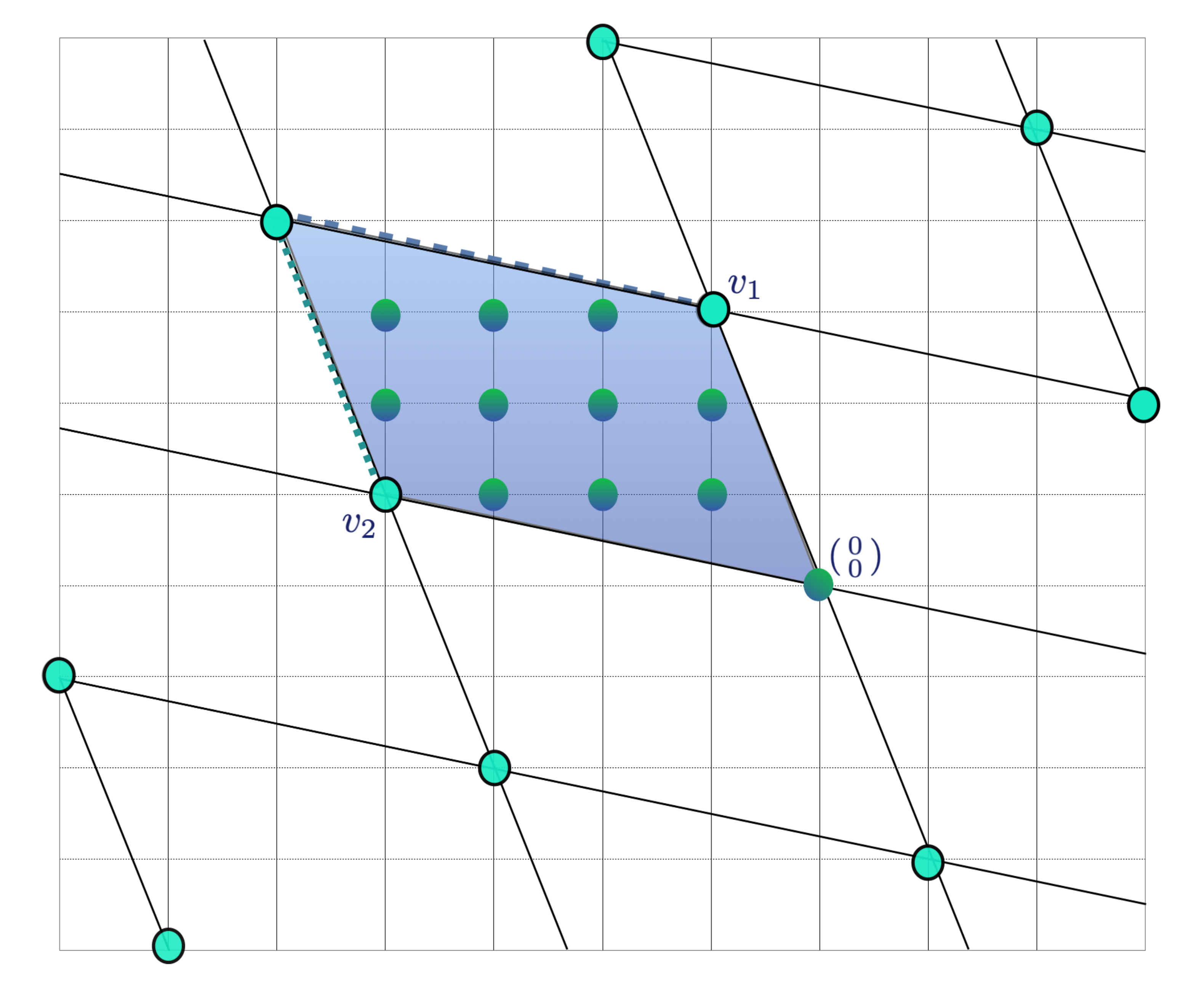}
\end{center}
\caption{A fundamental parallelepiped (half-open), for a lattice $\L$, generated by the vectors $v_1$ and $v_2$.} \label{parallelepiped1}
\end{figure}

\begin{wrapfigure}{R}{0.49\textwidth}
\centering
\includegraphics[width=0.31\textwidth]{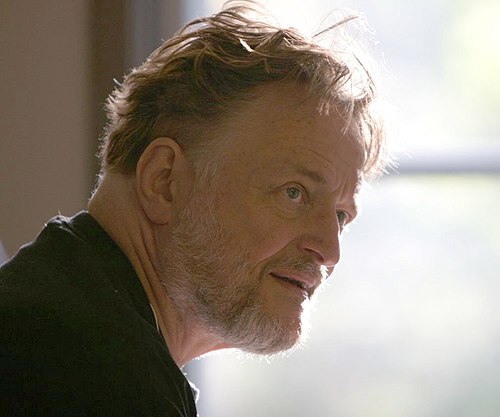}
\caption{John Conway}
\label{JOHN CONWAY PIC}
\end{wrapfigure}

\bigskip \bigskip
\section{Intuition}

We introduce Euclidean lattices here, 
which may be thought of intuitively as regularly-spaced points in 
$\R^d$, with some  hidden number-theoretic structure.   Another intuitive way to think of lattices is that they are one of the most natural ways to {\bf discretize Euclidean space}. 
A lattice in $\R^d$ is also  the most natural extension of an infinite set of equally-spaced points on the real line.     In the real-world, lattices come up very naturally when we study crystals, for example.

Perhaps it's not surprising that number theory comes in through 
study of the integer lattice $\Z^d$ and more general lattices, as they are a $d$-dimensional extension of the integers $\Z$.   Moreover, whenever we 
study almost any periodic behavior, lattices naturally come up, essentially from the 
definition of {\bf periodicity} in Euclidean space. \index{periodicity}   John Conway and Neil Sloane have an encyclopedic treatment of lattices \cite{ConwaySloan.book}, and have contributed greatly to the development of the subject.  And of course, wherever there are lattices, there are also Fourier series, as we saw in Chapter  \ref{Fourier analysis basics}.

\section{Introduction to lattices}
For the sake of completeness, we recall again the definition of a general lattice, and in this chapter we'll understand lattices in more detail. 
\begin{defi}
A   {\bf lattice}  \index{lattice}  is defined by the integer linear span of a fixed set of linearly independent vectors  $\{ v_1, \dots, v_m \} \subset \R^d$:
\begin{equation}\label{def.lattice}
\L :=  \left\{  n_1  v_1  + \cdots +  n_m v_m \in \R^d   \bigm |     \text{ all } n_j \in \Z    \right\}.
\end{equation}
\end{defi}
The first thing we might notice is that, by Definition \ref{def.lattice}, a lattice may also be written as follows:
\begin{equation}
\L := \left\{   
     \begin{pmatrix} |  &  |  &  ...   & |  \\  
                        v_1  &  v_2  &  ...& v_m   \\  
                          |  &  |  & ... &  |  \\ 
        \end{pmatrix}
\begin{pmatrix}
n_1 \\
 \vdots \\
  n_m \\
\end{pmatrix}
\   \biggm |     \    
\begin{pmatrix}
n_1 \\
 \vdots \\
  n_m \\
\end{pmatrix}
   \in  \Z^m   
\right\}
 := M(\Z^m),
\end{equation}
where by definition, $M$ is the $d \times m$ matrix whose columns are the vectors $v_1, \dots, v_m$.  This set of basis vectors 
$\{ v_1, \dots, v_m\}$ is called a {\bf basis} \index{lattice basis}
for the lattice $\L$, and $m$ is called the {\bf rank} of the lattice $\L$.  
In this context, we also use the notation ${\rm rank}(\L) = m$.

We will call $M$ a {\bf basis matrix} \index{basis matrix}
for the lattice $\L$.
 But there are always infinitely many other bases for $\L$ as well, and Lemma \ref{changing basis matrices} below shows how they are related to each other.

Most of the time, we will be interested in {\bf full-rank}  \index{full rank lattice} 
 lattices, which means that $m=d$;  however, sometimes we will also be interested in lattices that have lower rank, and it is important to understand them.  
 The {\bf determinant} of a full-rank lattice  $\L := M(\Z^d)$ is defined by 
 \[
 \det \L := |\det M|.  
 \]
 The determinant of a lattice measures how \emph{coarseness} of the lattice - the larger the determinant, the coarser the lattice.

 It is easy (and necessary) to prove that our definition of $\det \L$ is independent of the choice of basis matrix $M$, which is the content of Lemma
 \ref{changing basis matrices} below. 
 To better understand lattices, we need the {\bf unimodular group}, which we write as 
$GL_d(\Z)$,    \index{unimodular matrix}   under matrix multiplication:
\begin{equation}  \label{Definition of the unimodular group}
GL_d(\Z) := 
\left\{  M \bigm |        M \text{ is a }
 d \times d \text{ integer matrix,  with} \   |\det M| = 1 
 \right\}.
\end{equation}
\index{unimodular group}
\noindent
The elements of $\rm{GL_d}(\Z)$ are called {\bf unimodular matrices}. \index{unimodular matrix}
By definition, this group of matrices includes both the identity $I$ and the negative identity $-I$.
The easy fact that $GL_d(\Z)$ really is a group, under matrix multiplication, is a standard and easy fact \cite{MorrisNewman}. 
\begin{example}
 \rm{
Some typical elements of $\rm{GL_2}(\Z)$ are 
\[
S = \big(\begin{smallmatrix}
\ 0 & 1  \\
-1 & 0
\end{smallmatrix}
\big), 
T:=  \big(\begin{smallmatrix}
1 & 1  \\
1 & 0
\end{smallmatrix}
\big), 
 -I := \big(\begin{smallmatrix}
-1 &  \ 0  \\
\   0 & -1
\end{smallmatrix}
\big),
\text{ and \  }
 \big(\begin{smallmatrix}
1 &  n  \\
  0 & 1
\end{smallmatrix}
\big),
\]
where $n\in \Z$.  Interestingly, there is still no complete understanding of all of the subgroups of 
$GL_2(\Z)$ (see Newman's book \cite{MorrisNewman}).
}
\hfill $\square$
\end{example}
Now we suppose a lattice $\L$ is  defined by two different basis matrices:  $\L  = M_1(\Z^d)$ and
$\L  = M_2(\Z^d)$.   Is there a nice relationship between $M_1$ and $M_2$?   

\begin{lem}\label{changing basis matrices}
If a full-rank lattice $\L \subset \R^d$ is defined by two different basis matrices $M_1$, and $M_2$,
then 
\[
M_1 = M_2 U,
\]
where $U \in \rm{GL_d}(\Z)$,  a unimodular matrix. 
In particular, $\det \L$  is  independent of the 
choice of basis matrix $M$.
\end{lem}
\begin{proof}
By hypothesis, we know that the columns of $M_1$, say $v_1, \dots, v_d$,  form a basis of $\L$, and that the columns of $M_2$, say $w_1, \dots, w_d$,  also form a basis of $\L$.   So we can begin by writing each fixed basis vector $v_j$ in terms of all the basis vectors $w_k$:  
\[
v_j = \sum_{k=1}^d c_{j,k} w_k, 
\]
for each $j = 1, \dots, d$, and for some $c_{j,k} \in \Z$.   We may collect all $d$ of these identities into matrix form:
\[
M_1 = M_2 C,
\]
where $C$ is the  integer matrix whose entries are defined by the integer coefficients $c_{j,k}$ above. 
Conversely, we may also write each basis vector $w_j$ in terms of the basis vectors $v_k$:
$w_j = \sum_{k=1}^d d_{j,k} v_k$, for some $d_{j,k}\in\Z$, getting another matrix identity:
\[
M_2 = M_1 D.
\]
Altogether we have 
\[
M_1 = M_2 C = (M_1 D) C,
\]
and since $M_1^{-1}$ exists by assumption, we get  $DC= I$, the identity matrix.  Taking determinants, we see that 
\[
| \det D | |  \det C | = 1,
\]
and since both $C$ and $D$ are integer matrices, they must belong to $\rm{GL_d}(\Z)$, by definition.
Finally, because, because a unimodular matrix $U$ has $|\det U|=1$, we see that any two basis $M_1, M_2$ matrices satisfy 
$|\det M_1 | = |\det M_2 |$.
\end{proof}

\begin{lem} \label{Automorphisms of lattices} 
The group of one-to-one, onto, linear transformations from $\Z^d$ to itself 
is equal to
the unimodular group $GL_d(\Z)$.  \index{unimodular group}
\hfill $\square$
\end{lem}
Try to prove this yourself, and for much more about the delicate internal structure
of $GL_d(\Z)$, even for $d=2$, see Morris Newman's book \cite{MorrisNewman}.
 
 \begin{example}
 \rm{
 In $\R^1$, we have the integer lattice $\Z$, but we also have lattices of the form $r\Z$, for any real number $r$.  
 It's easy to show that any lattice in $\R^1$ is of this latter type (Exercise \ref{lattices in R^1}).
 For example, if $r = \sqrt 2$, then all integer multiples of $\sqrt 2$ form a $1$-dimensional lattice.  }
 \hfill $\square$
 \end{example}

\begin{figure}[htb]
 \begin{center}
\includegraphics[totalheight=2.8in]{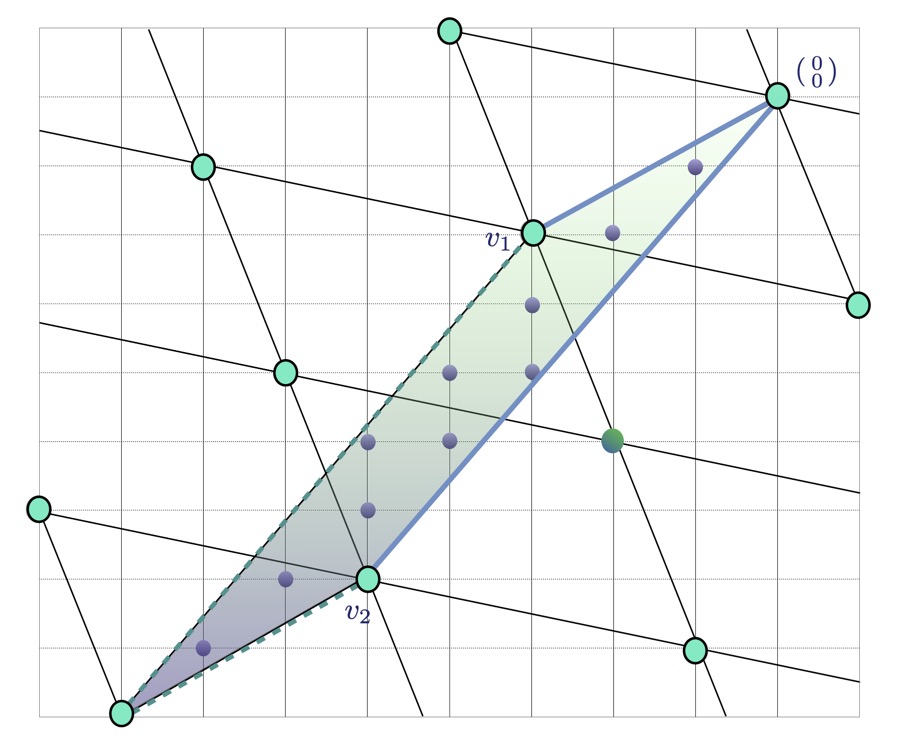}
\end{center}
\caption{A second fundamental parallelepiped for the same lattice $\L$ as in 
Figure \ref{parallelepiped1}}  \label{parallelepiped2}
\end{figure}

 \begin{example}
  \rm{
 In $\R^2$, consider the lattice $\L$ generated by the two integer vectors 
 $v_1:=\icol{-1\\3}$ 
 and $v_2:= \icol{-4\\   1}$, drawn in Figure  \ref{parallelepiped1}. 
 A different choice of basis for the same lattice $\L$ is $\{  \icol{-3\\-2}, \icol{-8\\ -9}   \}$, drawn in 
 Figure \ref{parallelepiped2}.   We note that $\det \L = 11$, and indeed the areas of both half-open parallelepipeds equals $11$.     }
 \hfill $\square$
 \end{example}

\medskip
 A {\bf fundamental parallelepiped} \index{fundamental parallelepiped}
 for a lattice $\L$ with basis $\{  v_1, \dots, v_m   \}$   is:
 \begin{equation} \label{def:half-open parallelepiped}
 \Pi:= \left\{ \lambda_1 v_1 + \cdots + \lambda_m v_m  \bigm |     \text{  all }   0 \leq  \lambda_k < 1   \right\},
 \end{equation}
 also known as a {\bf half-open parallelepiped}. 
Any lattice $\L$ has infinitely many fundamental parallelepipeds and 
it is a nice fact of life that they are all images of one another by the unimodular group (Exercise \ref{fundamental domains}) .

We have the pleasant property
that $\Pi$  tiles $\R^d$ by translations with vectors
from $\L$, and with no overlaps.   Let's make this intuition more precise, in the following lemma.
We recall that for any real $\alpha$, 
\[
\lfloor \alpha \rfloor \text{ is the greatest integer not exceeding } \alpha, 
\]
and 
\[
\{ \alpha\}:= \alpha - \lfloor \alpha \rfloor 
\]
 is called the {\bf fractional part} \index{fractional part} 
 of $\alpha$.   Clearly $0 \leq \{ \alpha\} < 1$.

\begin{lem}  \label{tiling by translates of Pi}
Suppose we are given a full rank lattice $\L\subset \R^d$, and a fundamental parallelepiped $\Pi$ for $\L$,  
as in Definition \eqref{def:half-open parallelepiped}.
   Then any $x\in \R^d$ may be written uniquely as
\[
x = n + y
\] 
where $n \in \L$, and $y \in \Pi$.  Consequently, $\Pi$ tiles $\R^d$ by translations with $\L$.
\end{lem}
 \begin{proof}
We know that $\Pi$ is formed by a basis for the lattice $\L$, and we can label the basis elements by
 $v_1, \dots, v_d$.   These $d$ vectors also form a basis for $\R^d$,  
 so in particular any $x\in \R^d$ may be written as
 \[
 x = \sum_{j=1}^d \alpha_j v_j.
 \]
Writing each $\alpha_j:= \lfloor \alpha_j \rfloor + \{ \alpha_j \}$, we have
 \[
 x = \sum_{j=1}^d \lfloor \alpha_j \rfloor  v_j   + \sum_{j=1}^d \{ \alpha_j \} v_j := n + y,
 \]
 where we've defined $n :=  \sum_{j=1}^d \lfloor \alpha_j \rfloor  v_j$, and 
$y:= \sum_{j=1}^d \{ \alpha_j \} v_j$.   Since  $\lfloor \alpha_j \rfloor \in \Z$, we see that
$n \in \L$.    Since  $0\leq \{ \alpha_j \} < 1$, we see that $y \in \Pi$. 

To prove uniqueness, suppose we are given $x:= n_1 + y_1 = n_2 + y_2$, where $n_1, n_2 \in \L$ and $y_1, y_2 \in \Pi$.   So by definition 
$y_1 = \sum_{j=1}^d \{ \alpha_{j, 1} \} v_j$ and $y_2=\sum_{j=1}^d \{ \alpha_{j, 2} \} v_j$.
Then $y_1 - y_2 = n_2 - n_1 \in \L$, which means that $ \alpha_{j, 1} - \alpha_{j, 2} \in \Z$.  But  $0 \leq \alpha_{j, 1}<1$ and $0 \leq \alpha_{j, 2}<1$  implies that $ \alpha_{j, 1} - \alpha_{j, 2}=0$.   Therefore  $y_1 = y_2$, and so $n_1 = n_2$. 
 \end{proof}
 
 \bigskip
 
 It follows from the uniqueness statement of Lemma \ref{tiling by translates of Pi}, for example, that the origin is the unique lattice point of $\L$
 that lies in any fixed fundamental parallelepiped of $\L$.

 How do we define the determinant of a ``lower dimensional'' lattice?  Well, let's begin with
 a lattice $\L\subset \R^d$ of rank $r\leq d$.  
 We can observe how the squared lengths of vectors in $\L$ behave with respect to a given basis of $\L$:
 \begin{equation} 
 \| x \|^2 = \left\langle  \sum_{j=1}^r c_j v_j, \,  \sum_{k=1}^r c_k v_k \right\rangle =
 \sum_{1\leq j, k \leq r} c_j c_k \langle v_j, \,   v_k \rangle := c^T M^T M  c,
 \end{equation}
 where $M^TM$ is an $r\times r$ matrix whose columns are  basis vectors of $\L$. 
With this as motivation, we define:
\begin{equation}\label{def. of sublattice determinant}
\det \L := \sqrt{ M^T M},
\end{equation}
called the {\bf determinant of the lattice} $\L$.  \index{determinant of a general lattice}
 This definition coincides, as it turns out,
 with the Lebesgue measure of any fundamental parallelepiped of $\L$ (Exercise \ref{equivalence between determinants of a sublattice}).   
 
 We may sometimes also use the following ubiquitous inequality of Hadamard, which gives
 a bound on the determinant of any invertible matrix, and hence on the determinant of a lattice.
Hadamard's inequality can be intuitively visualized:  if we keep all the lengths of the sides of a parallelepiped constant,  and consider all possible parallelepipeds $\P$ 
with these fixed side lengths, then the volume of $\P$ is maximized exactly when $\P$ is rectangular.

\begin{thm}[Hadamard's inequality]  \label{Hadamard inequality}  \index{Hadamard's inequality}
Given a non-singular matrix $M$, over the reals, whose column vectors are 
$v_1, \dots, v_d$, we have:
\[
|\det M|    \leq     \| v_1 \|    \|v_2\| \cdots \|v_d\|,
\]
\end{thm}
with equality if and only if all of the $v_k$'s are pairwise orthogonal.
\begin{proof}
We use the following matrix decomposition from Linear Algebra:  $M = QR$, where $Q$ is an orthogonal
matrix, $R:= [r_{i,j}]$ is an upper-triangular matrix, and $r_{kk} > 0$ (this decomposition is a well-known consequence of the Gram-Schmidt process applied to the columns of M).  So now we know that $|\det Q| = 1$, and $\det R = \prod_{k=1}^d r_{kk}$, and it follows that
\[
|\det M| = |\det Q \det R| = \det R = \prod_{k=1}^d r_{kk}.
\]
Let's label the columns of $Q$ by $Q_k$, and the columns of $R$ by $R_k$.  We now consider the matrix
$M^T M = R^T Q^T Q R = R^T R$.  Comparing the diagonal elements on both sides of $M^T M = R^T R$,
we see that $\| Q_k\|^2 = \| R_K \|^2$.  But we also have $\| R_K \|^2 \geq r_{kk}^2$, so that
$\|Q_k \| \geq r_{kk}$.
Altogether we have
\begin{equation}  \label{product formula}
|\det M| = \prod_{k=1}^d    r_{kk}  \leq  \prod_{k=1}^d   \|Q_k \|.
\end{equation}
The case of equality occurs if and only if  $\| R_K \|^2 = r_{kk}^2$ for all $1\leq k\leq d$, and this case of equality would mean
that $R$ is a diagonal matrix.  Thus, we have  equality in inequality \eqref{product formula} if and only if
$M^T M = R^T R$ is a diagonal matrix,  which means that the columns of $M$ are mutually orthogonal.
\end{proof}

 
 \bigskip
 \section{Sublattices}

\begin{figure}[htb]
 \begin{center}
\includegraphics[totalheight=3.3in]{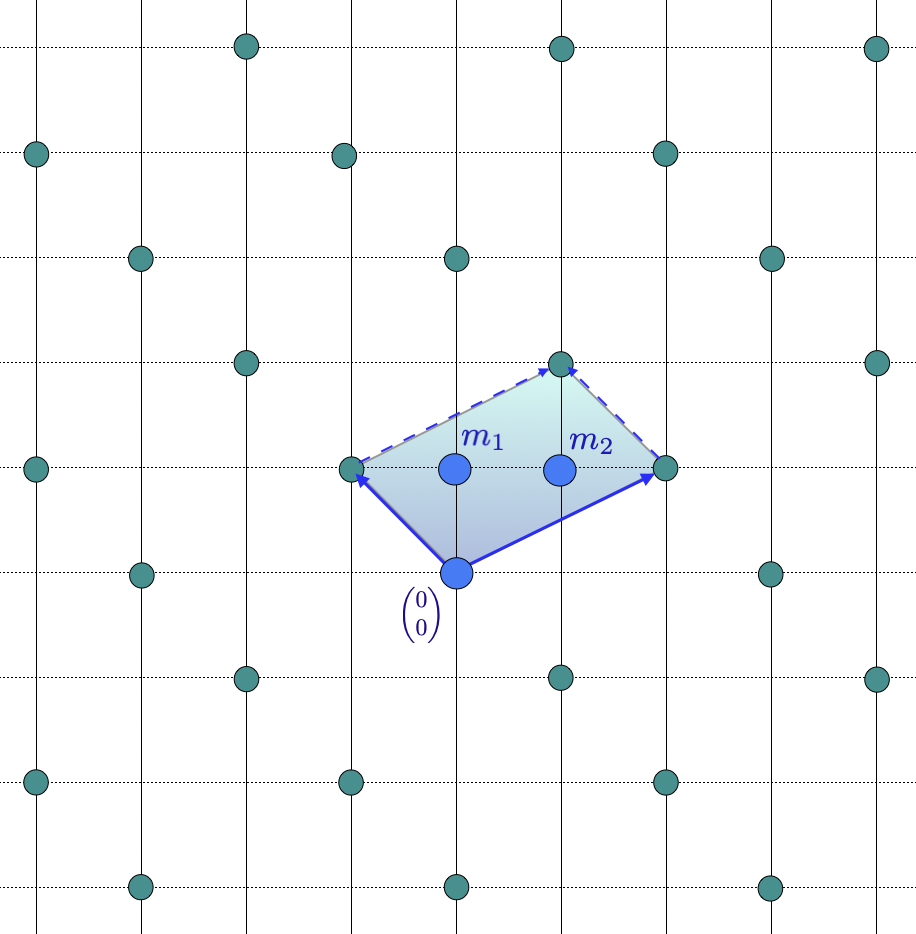}
\end{center}
\caption{The sublattice \green{$\L_0\subset \Z^2$} of Example \ref{Ex: sublattices 1}, drawn with the thickened green lattice points}  
\label{sublattice2}
\end{figure}

 Given two lattices $\L_0\subset \R^d$, and $\L \subset \R^d$, such that 
 $\L_0 \subseteq \L$, we say that {\bf $\L_0$ is a sublattice of $\L$}.    \index{sublattice}  
Sublattices that have the same rank are rather interesting, and extremely useful in applications.  So we'll usually focus on sublattices $\L_0 \subseteq \L$ such that $\rm{rank}(\L) = \rm{rank}(\L_0)$.   In this context, we sometimes call $\L_0$ a {\bf coarser lattice}, 
and $\L$ a {\bf finer lattice}.
 Given a sublattice $\L_0$ of $\L$, both of the same rank, a crucial idea is to  think of all of the translates 
 of  $\L_0$ by an element of the finer lattice $\L$:
 \begin{equation}
 \L / \L _0:= \left\{    \L_0 + m  \bigm |  m \in \L            \right\}.
 \end{equation}
 Each such translate $\L _0+ m$ is called a {\bf coset}  \index{coset}
 of $\L_0$ in $\L$.  The collection $\L / \L_0$ of all of these cosets is called a 
 {\bf quotient group}, and as we'll see shortly this is a very interesting finite set.   \index{quotient group}
 
\begin{figure}[htb]
\begin{center}
\includegraphics[totalheight=3.1in]{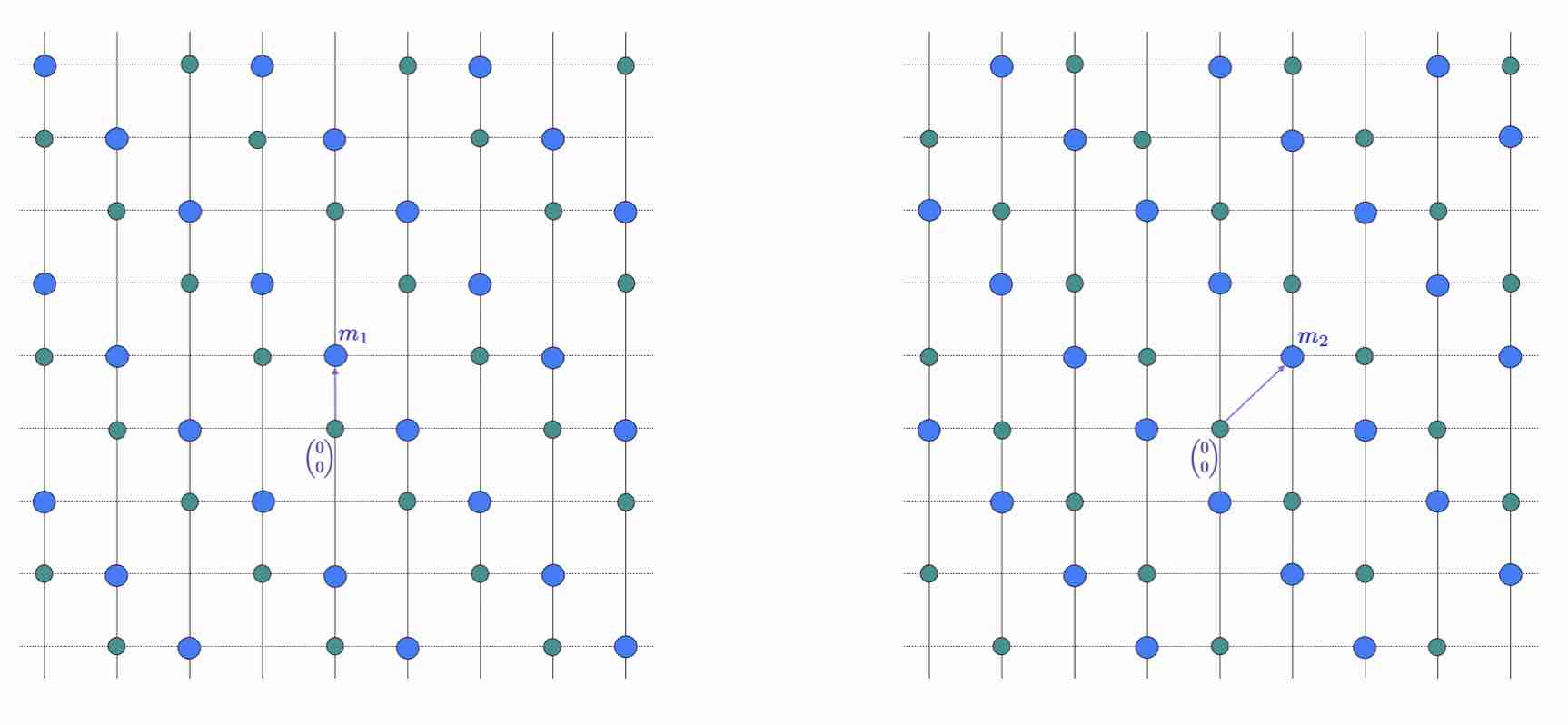}
\end{center}
\caption{Left: the thickened blue points represent the nontrivial coset \blue{$\L_0 + m_1$} 
of the sublattice $\L_0$.
Right: the thickened blue points represent the nontrivial coset \blue{$\L_0 + m_2$}
 of the sublattice $\L_0$.
(see Example \ref{Ex: sublattices 1}).
}  
\label{sublattice3}
\end{figure}

 \begin{example}\label{Ex: sublattices 1}
 \rm{
   Figure \ref{sublattice2} shows a  sublattice \green{$\L_0$} of the integer lattice $\Z^2$, with a fundamental parallelepiped $\Pi$ that is generated by the two vectors $\icol{-1 \\  \ 1}$ and $\icol{2 \\ 1}$.      
Here $\rm{area}(\Pi)=3$, and there are exactly $3$ cosets of $\L_0$ in $\Z^2$ (not a coincidence, as we'll see in Theorem \ref{sublattice index} below).  These $3$ cosets are: the trivial coset $\L_0$, and the two nontrivial cosets $\L_0 + m_1$ and $\L_0 +m_2$, drawn in Figure \ref{sublattice3} with thickened blue points.
 }
 \hfill $\square$
 \end{example}

\begin{figure}[htb]
 \begin{center}
\includegraphics[totalheight=2.2in]{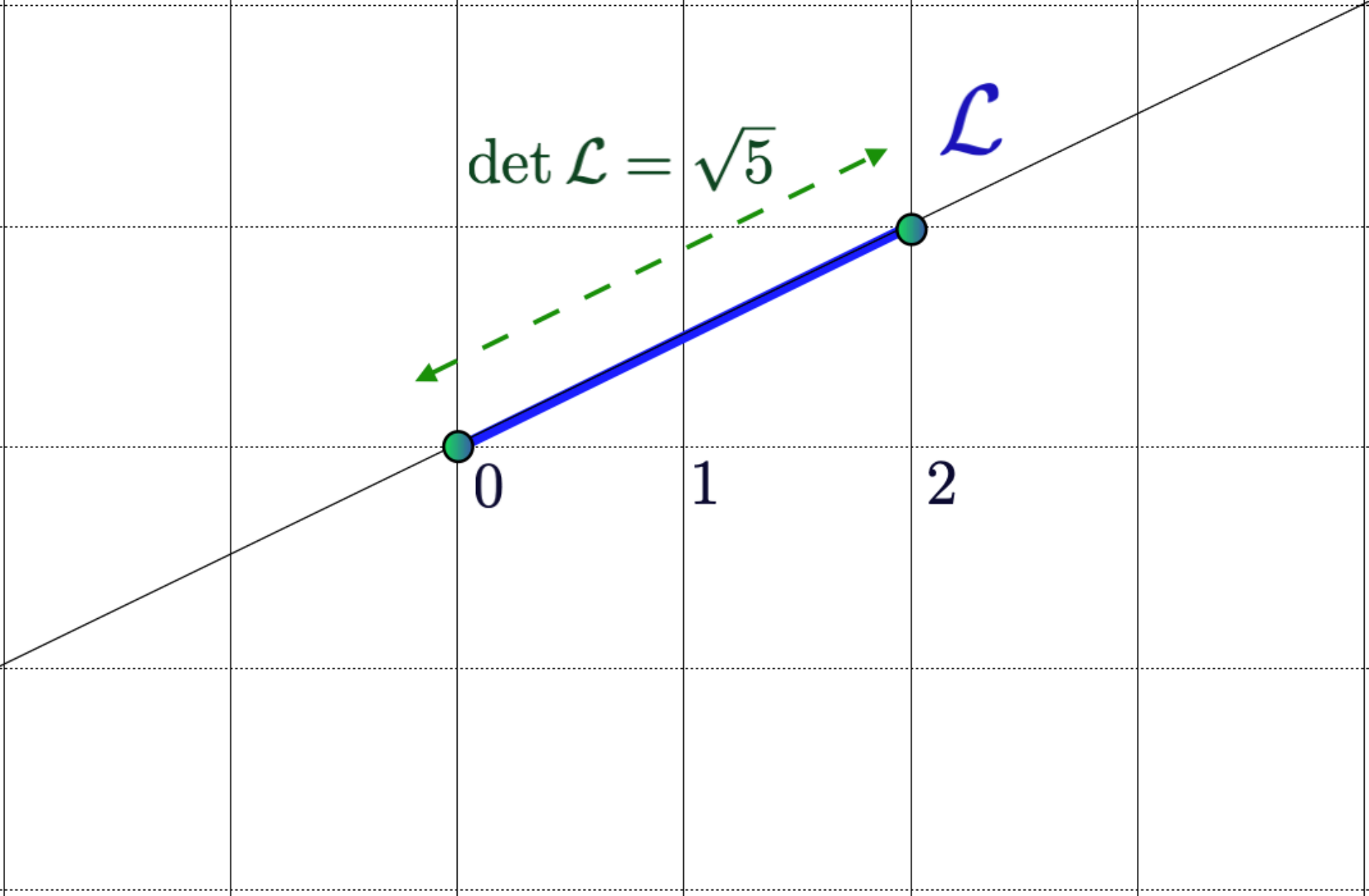}
\end{center}
\caption{
A sublattice $\L \subset \Z^2$ of rank $1$, which has just one basis vector. 
Here $\L$ has a $1$-dimensional fundamental parallelepiped, showing that 
$\det \L = \sqrt{v^T v} = \sqrt 5$, consistent with Definition \ref{def. of sublattice determinant}.
}  
\label{sublattice, rank 1}
\end{figure}

As an example of a lower-dimensional sublattice, Figure \ref{sublattice, rank 1} shows a rank $1$ sublattice of the integer lattice $\Z^2$, together with its determinant.

To better understand sublattices and some of their many subtleties, it's useful to first understand how many fundamental parallelpipeds of a lattice $\L$ are contained in a large ball, asymptotically.   Here we follow the geometric approach taken by Barvinok (\cite{Barvinok.A.course.in.convexity}, p. 287).

\begin{figure}[htb]
 \begin{center}
\includegraphics[totalheight=5in]{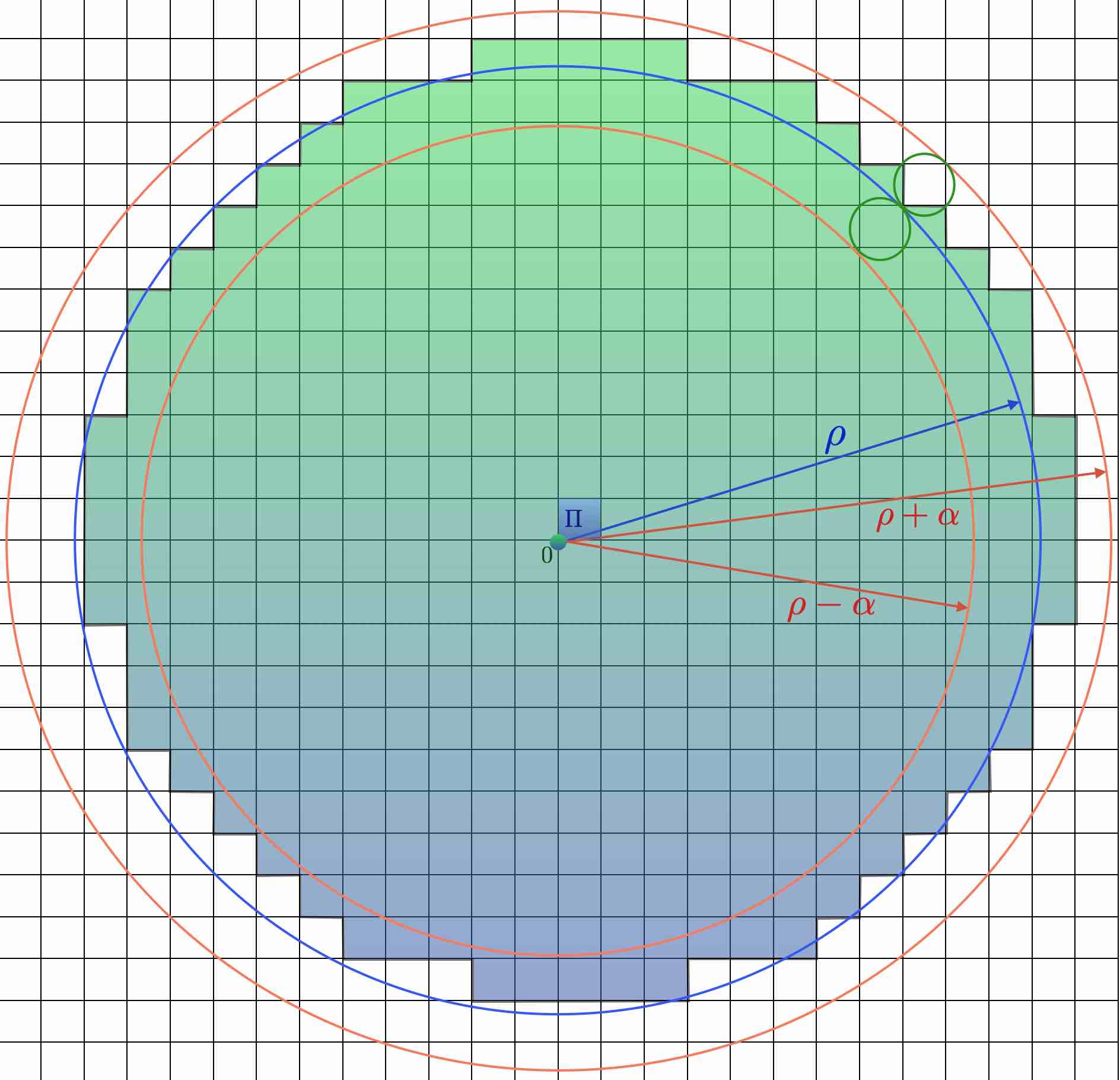}
\end{center}
\caption{Here we use the integer lattice $\Z^2$ to illustrate the proof of 
Theorem \ref{Prepping for main result on lattices}.  The shaded squares represent the number of integer points in a ball of radius $\rho$, after identifying each integer point with its northeast square. }  
\label{packing squared into a large circle}
\end{figure}

 \begin{thm} \label{Prepping for main result on lattices}
 Let $\L \subseteq \R^d$ be a lattice, and let $B_\rho:= \{  x \in \R^d \mid    \| x\| \leq \rho \} \subset \R^d$ be the ball of radius $\rho >0$. 
 Then:
 \begin{enumerate}[(a)]
 \item    \label{first part of prepping for lattices}
 \[
  \lim_{\rho \rightarrow \infty} \frac{ \left| \L \cap B_\rho \right|  }{  \vol B_\rho  } =  \frac{1}{\det \L}.
  \]
\item     \label{second part of prepping for lattices}
In general, for any $x \in \R^d$ we have:
 \[
  \lim_{\rho \rightarrow \infty} \frac{ \left|(\L+x) \cap B_\rho \right|  }{  \vol B_\rho  } =  \frac{1}{\det \L}.
  \]
\end{enumerate}
\end{thm}
 Figure \ref{packing squared into a large circle} may be helpful to the reader, while digesting the proof.

\begin{proof} 
We let $\Pi$ be a fundamental parallelepiped of $\L$.  Considering the set of all lattice points $n \in \L$ that are contained in the ball $B_\rho$ of radius $\rho$, we may use each of these points to translate $\Pi$:
\[
A_\rho:= \bigcup_{n \in \L\cap B_\rho} (\Pi + n).
\]
If we associate to each such lattice point $n\in \L$ the unique translate of $\Pi$ that lies to its northeast direction, then we have the collection of fundamental parallelepipeds
that are drawn with the shaded green squares in Figure \ref{packing squared into a large circle}.
By Lemma \ref{tiling by translates of Pi}, we know that the lattice translates of $\Pi$ tile $\R^d$, so we have
\[
\vol A_\rho =  \left |  \L\cap B_\rho  \right| \vol \Pi.
\]
Because $\Pi$ is bounded, it is contained in some ball $B_{\alpha}$, with radius $\alpha >0$.   For the construction of the proof, 
we'll think of a `band' of diameter  $4\alpha$ placed around the perimeter of $B_\rho$, where the boundary of this band consists of the two large orange circles in Figure \ref{packing squared into a large circle}.
While it is true that some portion of the translated copies of $\Pi$ in $\bigcup_{n \in \L\cap B_\rho} (\Pi + n)$
leak out of the ball $B_\rho$,  we can nevertheless cover them as well by taking the Minkowski sum of $B_\rho$ with $B_{2\alpha}$, obtaining $B_{r+2\alpha}$.  Therefore $A_\rho \subset B_{\rho+2\alpha}$.

To see the inclusion $B_{\rho-2\alpha} \subset A_\rho$, we note that by using 
Lemma  \ref{tiling by translates of Pi} again, we may conclude that each point of $B_{\rho-2\alpha}$ is contained in some translate $\Pi+n$, with $n \in B_{\rho}\cap \L$.   Putting everything together, we therefore have:
\[
\vol B_{\rho-2\alpha} \leq  \vol A_\rho = \left |  \L\cap B_\rho  \right| \vol \Pi \leq \vol B_{\rho+2\alpha},
\]
which we'll rewrite as 
\begin{equation} \label{end of first part}
\frac{
\vol B_{\rho-2\alpha}
}
{
\vol B_{\rho}
}
\leq   
  \frac{
  \left    |  \L\cap B_\rho  \right|       
  }
  {
  \vol B_{\rho}
  }
  \vol \Pi 
 \leq 
\frac{
  \vol B_{\rho+2\alpha}
  }
  {
  \vol B_{\rho}
  }.
\end{equation}
So if we prove that 
$\lim_{\rho \rightarrow \infty} \frac{\vol B_{\rho-2\alpha}}{\vol B_{\rho}} = 
\lim_{\rho \rightarrow \infty} \frac{\vol B_{\rho+2\alpha}}{\vol B_{\rho}} =1$, then by \eqref{end of first part}
we will have proved  part \ref{first part of prepping for lattices}.  Recalling from
 \eqref{dilated volumes of balls and spheres} that $\vol B_{\rho} = 
 \frac{ \pi^{\frac{d}{2}} }{\Gamma\left(\frac{d}{2} +1\right)} \rho^d$, we can finish the computation:
\begin{equation} \label{little limit of 1}
\lim_{\rho \rightarrow \infty} \frac{\vol B_{\rho-2\alpha}}{\vol B_{\rho}} =
\lim_{\rho \rightarrow \infty} \frac{(\rho-2\alpha)^d}{\rho^d} = 1,
\end{equation}
and similarly $\lim_{\rho \rightarrow \infty} \frac{\vol B_{\rho+2\alpha}}{\vol B_{\rho}}=1$.

To prove part \ref{second part of prepping for lattices}, we first show that the following set inclusions hold:
\begin{equation}\label{set inclusions for prepping}
x +\left( \L \cap B_{\rho - \|x\|}\right)   \  \  \subset \ \   \left(  \L+x \right) \cap B_\rho 
\  \  \subset \ \ 
x +   \left(\L  \cap  B_{\rho+\|x\|}\right).
\end{equation}
To prove the left-hand inclusion, let $y =   x+z$, with $z \in  \L \cap B_{\rho - \|x\|} $.  Then
$\| y \| \leq \|x\| + \|z\| \leq  \|x\| + (\rho - \|x\|) = \rho$, giving us $y \in B_\rho$.  Also, $z \in \L \implies y \in \L+x$, which together with $y \in B_\rho$ proves the first inclusion.   To prove the second inclusion in \eqref{set inclusions for prepping}, 
let $y \in  \left(  \L+x \right) \cap B_\rho $, so that by assumption
$y = x + n$, with $n \in \L$ and $\|y\| \leq \rho$.  It remains to show that $\|n\| \leq \rho + \|x\|$, but this follows from $\rho \geq \|y\| \geq \|n \| - \|x\|$.

From \eqref{set inclusions for prepping}, we have
\[
\left | 
  \L  \cap B_{\rho - \|x\|}
\right|
\leq 
\left|   \left(  \L +x\right) \cap B_\rho  \right|
\leq
\left|  \L \cap  \left( B_{\rho+\|x\|}  \right) \right|,
\]
which we'll rewrite as
\[
\frac{ \left |   \L \cap B_{\rho - \|x\|}  \right| }{\vol B_{\rho- \|x\|}}
\frac{\vol B_{\rho- \|x\|} }{  \vol B_\rho}
\leq 
\frac{\left|   \left(  \L+x \right) \cap B_\rho  \right|}{\vol B_\rho}
\leq
\frac{ \left |   \L \cap B_{\rho + \|x\|}  \right| }{\vol B_{\rho+ \|x\|}}
\frac{\vol B_{\rho+ \|x\|} }{  \vol B_\rho}.
\]
From part \ref{first part of prepping for lattices}, we  know that
$\lim_{\rho \rightarrow \infty} \frac{ \left |   \L \cap B_{\rho - \|x\|}  \right| }{\vol B_{\rho- \|x\|}} = 
\frac{1}{\det \L}
= \lim_{\rho \rightarrow \infty} \frac{ \left |   \L \cap B_{\rho + \|x\|}  \right| }{\vol B_{\rho+ \|x\|}}$, and from \eqref{little limit of 1} we know that 
$\lim_{\rho \rightarrow \infty}  \frac{\vol B_{\rho- \|x\|} }{  \vol B_\rho} = 1 = 
\lim_{\rho \rightarrow \infty}  \frac{\vol B_{\rho+ \|x\|} }{  \vol B_\rho}$, finishing the proof.

\end{proof}

\bigskip
 \begin{thm} \label{sublattice index}
 Let $\L_0 \subseteq \L$ be any two lattices of the same rank, so by definition $\L_0$ is a sublattice of $\L$.  
 Let $\Pi$ be any fundamental parallelepiped for
 $\L_0$.    Then 
 \begin{enumerate}[(a)]
 \item  
 \label{first part of sublattice thm}
 $\Pi \cap \L$ contains each coset representative of $\L / \L_0$ exactly once. 
 \item     
 \label{second part of sublattice thm}
  $\frac{ \det \L_0 }{ \det \L  }$ is a positive integer, and  is equal to the number of 
cosets of $\L_0$ in $\L$.  In other words, we have a finite abelian group $\L / \L_0$, whose size is
\[
\left| \L / \L_0 \right| =  \frac{ \det \L_0 }{ \det \L  }.
 \]
\item    
\label{third part of sublattice thm}
  Consequently, $\left|\Pi \cap \L \right| =\frac{ \det \L_0 }{ \det \L  }$.
\end{enumerate}
\end{thm} 
\begin{proof}
To prove part \ref{first part of sublattice thm}, we define $f:  \Pi \cap \L \rightarrow \L / \L_0$ by $f(n)=\L_0+n$.
We must show that $f$ is bijective.   
 To show $f$ is surjective, suppose that we are given any coset $\L_0+n$. By Lemma \ref{tiling by translates of Pi}, we know that $n = l_0 + x$, where $l_0 \in \L_0$ and 
 $x\in \Pi$.  Now $x = n - l_0$, and both $n, l_0 \in \L$ ($l_0 \in \L_0 \subset \L$), implying that $x \in \L$.  This proves surjectivity, because $x \in  \Pi \cap \L$ and 
 $f(x):= \L_0 + x = \L_0 +n -  l_0 = \L_0 + n$.
 
For the injectivity of $f$, suppose that $f(n_1) = f(n_2)$, where $n_1, n_2 \in  \Pi \cap \L$.  Then 
 $\L_0+n_1 = \L_0+n_2$, so that $n_1 - n_2 \in \L_0$.  But the only element in $\Pi$ that lies in $\L_0$ is the origin.  Therefore $n_1 -n_2 = 0$. 
 
To prove part \ref{second part of sublattice thm}, which is more interesting, 
we begin by letting $C:= \Pi \cap \L$, a finite set  of coset representatives for  $\L / \L_0$.  
So by definition we have $\L = \bigcup_{x\in C} \left( \L_0 + x \right)$.   Intersecting both sides of the latter identity with a ball $B_\rho$ of radius $\rho$, it follows from the disjointness of the latter union that
\begin{equation}\label{almost finishing main sublattice theorem}
\left |  \L \cap B_\rho \right | = \sum_{x\in C}  \left |  \left( \L_0 + x \right) \cap B_\rho  \right |.
\end{equation}

From Theorem \ref{Prepping for main result on lattices}, we know that 
 \[
  \lim_{\rho \rightarrow \infty} \frac{ \left| \L \cap B_\rho \right|  }{  \vol B_\rho  } 
  =  \frac{1}{\det \L},  \text{ and } 
  \lim_{\rho \rightarrow \infty} \frac{ \left|(\L_0 + x) \cap B_\rho \right|  }{  \vol B_\rho  }
  =  \frac{1}{\det \L_0},
    \]
for any $x \in \R^d$.  Dividing both sides of \eqref{almost finishing main sublattice theorem} by $\vol B_\rho$ and letting $\rho \rightarrow \infty$, we get:
\begin{equation}
\frac{1}{\det \L} =  \frac{1}{\det \L_0} \sum_{x\in C}  1 
=  \frac{1}{\det \L_0}  \left | \Pi \cap \L\right |,
\end{equation}
which finished the proofs of both \ref{second part of sublattice thm} and 
\ref{third part of sublattice thm}.
\end{proof}

\bigskip
 \begin{example}
 {\rm
 Let $\L:= \Z^d$, and $\L_0:= 2\Z^d$, the sublattice consisting of vectors all of whose coordinates are even integers.   So $\L_0 \subset \L$, and the quotient group $\L/ \L_0$ consists of the cosets 
 $\left\{   2\Z^d + n \bigm |   n \in \Z^d \right\}$.    It is (almost) obvious that the number of elements of the latter set is exactly $2^d$, and this observation is also a special case of  
 Theorem \ref{sublattice index}:
 \begin{equation}
  \Z^d / 2\Z^d =  \frac{\det 2\Z^d}{ \det  \Z^d} =2^d.
   \end{equation}

 We may also think of this quotient group $ \Z^d / 2\Z^d $ as the discrete unit cube, namely 
 $\left\{   0, 1  \right\}^d$, a common object in theoretical computer science, for example.
 \hfill $\square$
 }
 \end{example}


\bigskip
\section{Discrete subgroups - \\ an alternate definition of a  lattice}
\label{alternate def of a lattice - subgroups}

The goal here is to give another useful way to define a lattice.   The reader does not need any background in group theory, because the ideas here are self-contained, given some background in basic linear algebra.

\smallskip
\begin{defi} \label{discrete subgroup}
\begin{enumerate}[(a)]  
 We define a {\bf discrete subgroup}  \index{discrete subgroup}
 of $\R^d$ as a set $S \subset \R^d$, together with the operation of vector addition between all of its elements, which enjoys the following two properties. 
\item   {\bf  [The subgroup property]}   If    $x, y \in S$,  then $ x-y \in S$.   \\
          \label{discrete subgroup.first part}  
\item    {\bf  [The discrete property]}   There exists a positive real number $\delta >0$, 
such that \\
 the distance between any two distinct points of   $S$  
is at least  $\delta$.  \\    \label{discrete subgroup.second part}
\end{enumerate}
\end{defi}

In particular, it follows from Definition \ref{discrete subgroup} \ref{discrete subgroup.first part}  
that the zero vector must be in $S$, because for any $x \in S$, it must be the case that $x - x \in S$.   
The distance function that we alluded to in Definition \ref{discrete subgroup} \ref{discrete subgroup.second part}  
 is the usual Euclidean distance function, which we denote here by 
 \[
 \|x-y\|_2:= \sqrt{ \sum_{k=1}^d (x_k - y_k)^2}.  
 \]

\begin{example}
\rm{
The lattice $\Z^d$ is a discrete subgroup of $\R^d$.    In dimension $1$, the lattice
$r\Z$ is a discrete subgroup of $\R$, for any fixed $r>0$.   Can we think of discrete subgroups that are not lattices?   The answer is given by Lemma \ref{discrete subgroup equivalence} below.
}
\hfill $\square$
\end{example}

The magic here is the following very useful way of going back and forth between this new notion of a discrete subgroup of $\R^d$,  and our Definition \ref{def.lattice} of a lattice.  The idea of using this alternate Definition \ref{discrete subgroup}, as opposed to our previous Definition \ref{def.lattice} of a lattice,  is that it gives us a {\bf basis-free} \index{basis-free} way of discovering and proving facts about lattices. 

\bigskip
\begin{lem}\label{discrete subgroup equivalence}
$\L \subset \R^d$ is a lattice $\iff$ $\L$ is a discrete subgroup of $\R^d$.
\hfill $\square$
\end{lem}
(For a proof see \cite{GruberBook}).

\bigskip
\begin{example}
\rm{
Given any two lattices $\L_1, \L_2 \subset \R^d$, let's show that $S := \L_1 \cap \L_2$ is also a lattice.   
First, any lattice contains the zero vector, and it may be the case that their intersection consists of only the zero vector.   For any vectors $x, y \in S$, we also have $x,y \in \L_1$, and $x, y \in \L_2$, hence 
by the subgroup property of $\L_1 $ and of $\L_2$, we know that both $x-y \in \L_1$, and $x-y \in \L_2$.  In other words,   $x-y \in \L_1 \cap \L_2:= S$.   To see why the discrete property of 
Definition  \ref{discrete subgroup} holds here, we just notice that since $x-y \in \L_1$, we already know that
$| x-y | > \delta_1$, for some $\delta_1>0$;  similarly, because  $x-y \in \L_2$, we know that
$| x-y | > \delta_2$ for some $\delta_2>0$.    So we let $\delta:= \min(\delta_1, \delta_2\}$, and we have shown that $S$ is a discrete subgroup of $\R^d$.  By Lemma \ref{discrete subgroup equivalence}, we see that $S$ is a lattice.  

If we had used Definition \ref{def.lattice} of a lattice to show that $S$ is indeed a lattice, it would require us to work with bases, and this proof would be longer and less transparent.
}
\hfill $\square$
\end{example}

\bigskip

\begin{example} \label{A_d example}
 \rm{
Consider the following discrete set of points in $\R^d$:
 \[
A_{d-1}:= \left\{  x \in \Z^d  \bigm |  \sum_{k=1}^d x_k =0     \right\},
\]
for any $d\geq 2$, as depicted in Figure \ref{A_d}.    Is $A_d$ a lattice?   Using the definition  
\ref{def.lattice} of a lattice, it is not obvious that $A_d$ is a lattice,  because we would have to exhibit a basis, but
it turns out that the following set of vectors may be shown to be a basis:  $ \left\{e_2 - e_1, e_3 - e_1, \cdots e_d - e_1  \right\}$, and hence $A_d$ is a sublattice
of $\Z^d$, of rank $d-1$ (Exercise  \ref{basis for A_d}).

\begin{figure}[htb]
 \begin{center}
\includegraphics[totalheight=3.8in]{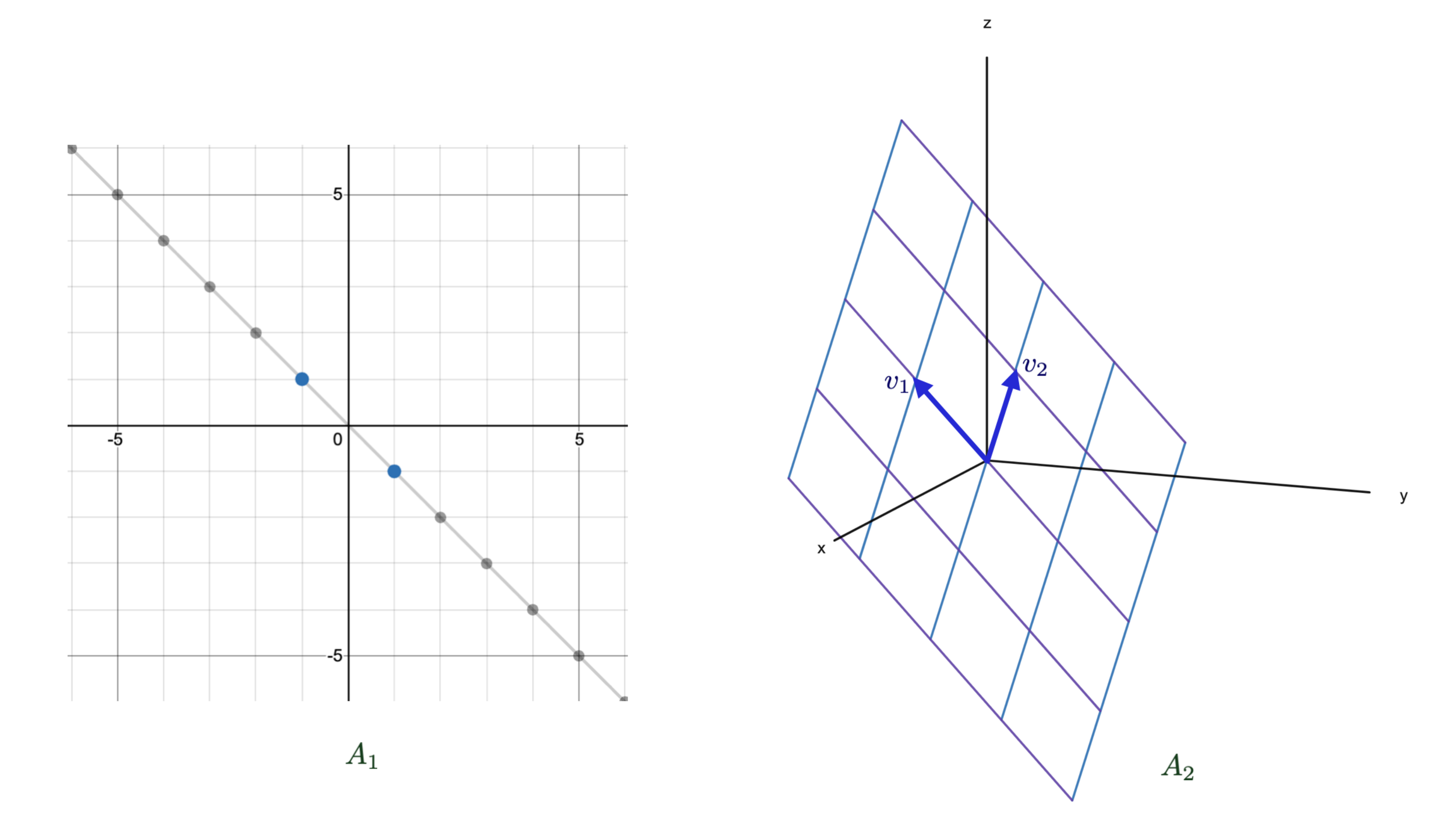}
\end{center}
\caption{The lattice $A_1$,  and 
the lattice $A_2$, with basis $\left\{ v_1, v_2  \right\}$}     \label{A_d}
\end{figure}

Just for fun, we will use Lemma    \ref{discrete subgroup equivalence} to show that
$A_d$ is indeed a lattice.   To verify the subgroup property of 
Definition \ref{discrete subgroup} \ref{discrete subgroup.first part} 
suppose that $x, y \in A_d$.  Then by definition we have
$\sum_{k=1}^d x_k =0$ and  $ \sum_{k=1}^d y_k =0$.  So $\sum_{k=1}^d (x_k - y_k)=0$, implying that $x-y \in A_d$.

To verify the discrete property of 
Definition \ref{discrete subgroup} \ref{discrete subgroup.second part}  
suppose we are given two distinct points $x, y \in A_d$.   We can first compute their ``cab metric'' distance function, in other
words the $L^1$-norm defined by
\[
 \| x-y \|_1:=  |x_1 - y_1| + \cdots + |x_d - y_d|,
\]
By assumption, there is at least one coordinate where $x$ and $y$ differ, say the $k$'th coordinate.   Then 
$    \| x-y \|_1  :=  |x_1 - y_1| + \cdots + |x_d - y_d| \geq 1$, because all of the coordinates are integers, and $x_k \not= y_k$ by assumption.  Since the $L^1$-norm and the $L^2$-norm
 are only off
by $\sqrt d$ (by Exercise \ref{elementary norm relations}), we  have:
\[
\sqrt{d}  \| x-y \|_2   \geq \| x-y \|_1    \geq 1, 
\]
so the property \ref{discrete subgroup} \ref{discrete subgroup.second part}   
is satisfied with $\delta := \frac{1}{\sqrt{d}}$, and we've shown that $A_d$ is a lattice.  
\hfill $\square$
}
\end{example}

We note that the lattices $A_d$ defined in Example \ref{A_d example}  are very important in many fields of Mathematics, including Lie algebras (root lattices),  Combinatorial geometry, and Number theory.

\section{Lattices defined by congruences}

In this section we develop some of the theory in a concrete manner.   A classic example of a lattice defined by an auxiliary algebraic construction is the following.  
Suppose we are given a constant integer vector $(c_1, \dots, c_d) \in \Z^d$, where we further assume 
that $\gcd(c_1, \dots, c_d) = 1$.   Let
\begin{equation}\label{Lattice from congruences}
C := \left\{  x \in \Z^d   \bigm |   c_1 x_1 + \cdots + c_d x_d \equiv 0 \mod N    \right\},
\end{equation}
where $N$ is a fixed positive integer.

Is $C$ a lattice?  
Indeed, we can see that $C$ is a lattice by first checking 
Definition \ref{discrete subgroup} \ref{discrete subgroup.first part}.  For any $x, y \in C$, we have
$c_1 x_1 + \cdots + c_d x_d \equiv 0 \mod N$ and $c_1 y_1 + \cdots + c_d y_d \equiv 0 \mod N$.
Subtracting these two congruences gives us  
$c_1 (x_1-y_1) + \cdots + c_d (x_d-y_d) \equiv 0 \mod N$, so that $x-y \in C$.  The verification 
of Definition \ref{discrete subgroup} \ref{discrete subgroup.second part} if left to the reader, and its logic is similar to Example  \ref{A_d example}.

There is even a simple formula for the volume of a fundamental parallelepiped for $C$:
\begin{equation}
 \det C = N,
 \end{equation}
as we prove below, in Lemma \ref{lemma:lattice defined by congruence}.
But perhaps we can solve an easier problem first.  Consider the {\bf discrete hyperplane} defined by:
\[
H:= \left\{  x \in \Z^d  \bigm |   c_1 x_1 + \cdots + c_d x_d =0   \right\},
\]
Is $H$ a lattice?  We claim that $H$ itself is indeed a sublattice of $\Z^d$, and has rank $d-1$.   
Since this verification is quite similar to the arguments above, we leave this as Exercise \ref{hyperplane lattice}.

The fundamental  parallelepiped (which is $(d-1)$-dimensional) of $H$  also has a wonderful formula, as follows.  First, we recall
a general fact (from Calculus/analytic geometry) about hyperplanes, 
namely that the distance $\delta$ between any two parallel hyperplanes \\  
$c_1 x_1 + \cdots + c_d x_d =  k_1$ and 
$c_1 x_1 + \cdots + c_d x_d =  k_2$ is given by 
\begin{equation} \label{distance between two hyperplanes}
\delta  = \frac{ |k_1 - k_2|}{\sqrt{ c_1^2 + \cdots + c_d^2}}.
\end{equation}
(see Exercise  \ref{distance between hyperplanes})

\medskip
\begin{lem}   \label{wonderful hyperplane determinant formula}
For any latttice defined by a discrete hyperplane \\
$H:= \left\{  x \in \Z^d  \bigm |   c_1 x_1 + \cdots + c_d x_d =0   \right\}$,  with 
$\gcd(c_1, \dots, c_d) = 1$, we have:
\begin{equation}  \label{tricky volume of sublattice}
\det H = \sqrt{ c_1^2 + \cdots + c_d^2}.
\end{equation}
\end{lem}
\begin{proof}
We first fix a basis $\{v_1, \dots, v_{d-1}\}$ for the $(d-1)$-dimensional sublattice defined by 
$H:= \left\{  x \in \Z^d  \mid   c_1 x_1 + \cdots + c_d x_d =0   \right\}$.   We  adjoin to this basis one new vector, namely any integer vector $w$ that translates $H$ to its `hyperplane companion'
$H + w$, which we define by 
\[
H+w:= ~\left\{  x \in \Z^d  \bigm |   c_1 x_1 + \cdots + c_d x_d =1   \right\}.   
\]
It's easy - and fun - to see that there are no integer points strictly between these two hyperplanes $H$ and $H+w$ (Exercise \ref{tiling the integer lattice with hyperplanes}), and so
the parallelepiped $\P$ formed by the edge vectors $v_1, \dots, v_{d-1}, w$ is a fundamental domain for $\Z^d$, hence has volume $1$.   

 On the other hand, we may also calculate the volume of $\P$ by multiplying the volume of its base times its height, using \eqref{distance between two hyperplanes}:
  \begin{align}
1= \vol \P &= (\text{volume of the base of } \P)(\text{height of } \P) \\ 
&= (\det H)\cdot \delta \\
&= (\det H) \frac{1}{\sqrt{c_1^2 + \cdots + c_d^2}},
\end{align}
and so   $\det H = \sqrt{ c_1^2 + \cdots + c_d^2}$.
\end{proof}

It follows directly from the definition \ref{Lattice from congruences} of $C$ that we may write the lattice $C$ as a countable, disjoint union of translates of $H$:
\begin{equation}
C := \left\{  x \in \Z^d   \bigm |    c_1 x_1 + \cdots + c_d x_d = kN, \text{ where } k = 1, 2, 3, \dots   \right\}.
\end{equation}

To be concrete, let's work out some examples.

\begin{example}
\rm{
Using  Lemma \ref{wonderful hyperplane determinant formula},  we can easily compute the determinant of the $A_d$ lattice from Example \ref{A_d example}: 
\[
\det A_d = \sqrt{1 + 1 + \cdots + 1} =  \sqrt d.  
\]
\hfill $\square$
}
\end{example}

\medskip

\begin{example}\label{Congruence lattice 2d}
\rm{
As in Figure \ref{lattice 2d}, 
 consider the set of all integer points $(m, n) \in \R^2$ that satisfy
\[
 2m + 3n \equiv 0 \mod 4. 
 \]
   In this case the related hyperplane is the line $2x+3y = 0$, and the solutions to the latter congruence may be thought of as a union of discrete lines:
\[
C = \left\{  {x\choose y} \in \Z^2   \bigm |   2x+3y = 4k, \text{ and } k \in \Z \right\}.
\]
}
\begin{figure}[htb]
 \begin{center}
\includegraphics[totalheight=3.7in]{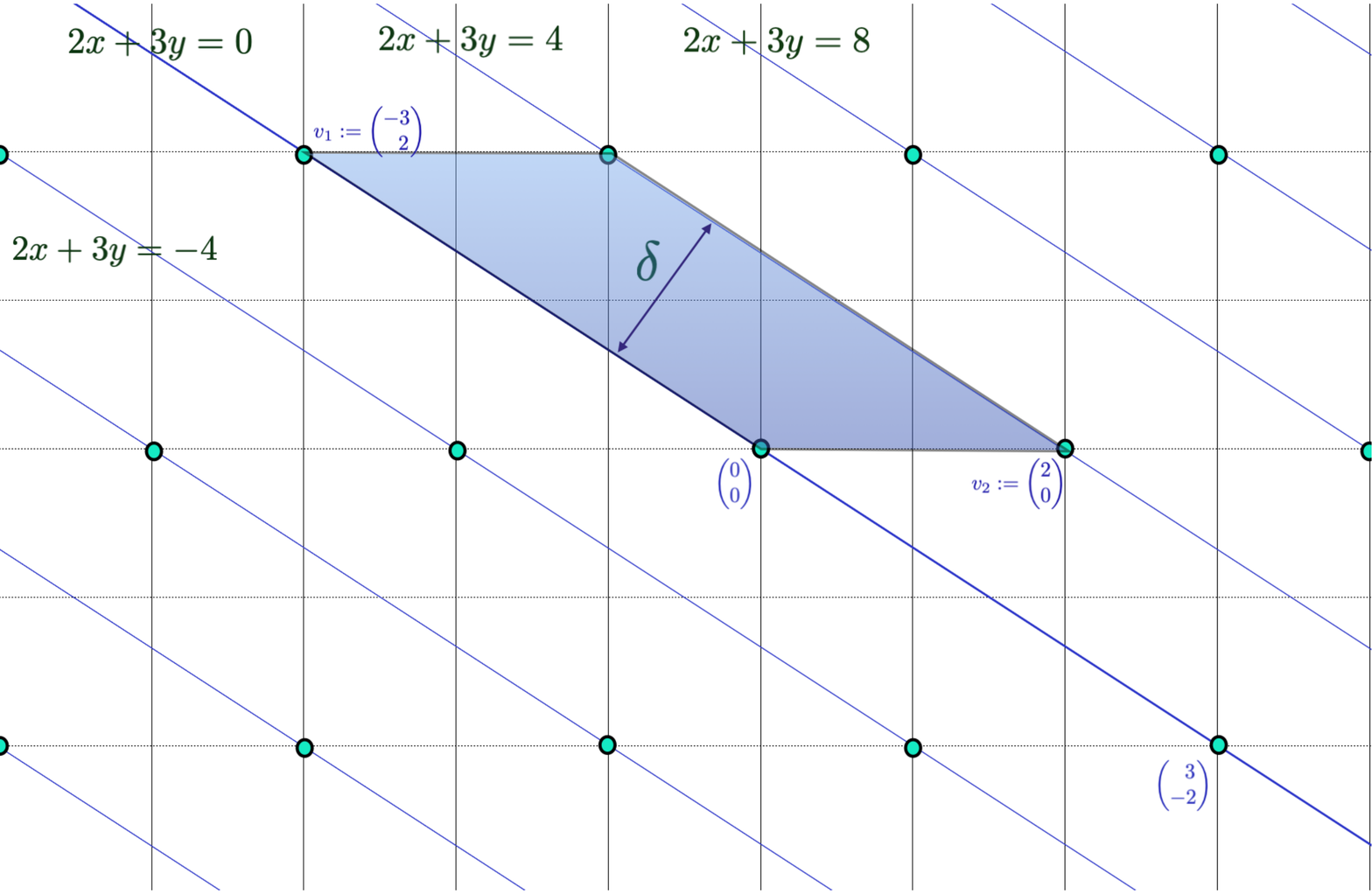}
\end{center}
\caption{The lattice of Example \ref{Congruence lattice 2d}}  \label{lattice 2d}
\end{figure}
\rm{
In other words, our lattice $C$, a special case of \eqref{Lattice from congruences}, can in this case be  visualized 
in Figure \ref{lattice 2d} as a disjoint union of discrete lines.  If we denote the distance between any two of these 
adjacent discrete lines  by $\delta$, then using \eqref{distance between two hyperplanes} we have:
\begin{equation}
\delta = \frac{4}{\sqrt{3^2 + 2^2}}.
\end{equation}
Finally, the determinant of our lattice $C$ here is the area of the 
shaded parallelepiped:
\begin{equation}
\det C =  \delta \sqrt{3^2 + 2^2} = 4.
\end{equation}
}
\hfill $\square$
\end{example}
Eager to prove the volume relation  $ \det C = N$, 
we can use the ideas of Example \ref{Congruence lattice 2d} as a springboard for this generalization.
Indeed, Example \ref{Congruence lattice 2d} and the proof of Lemma  
\ref{wonderful hyperplane determinant formula} both
suggest that we should compute the volume of a fundamental parallelepiped $\P$, for  the lattice $C$ (as opposed to the lattice $\Z^d$),  
by using a fundamental domain for its base, and then by multiplying  its volume by the height of $\P$.

 \begin{lem}\label{lemma:lattice defined by congruence}
Given a constant integer vector $(c_1, \dots, c_d) \in \Z^d$, with $\gcd(c_1, \dots, c_d) = 1$, 
let 
\begin{equation}
C := \left\{  x \in \Z^d   \bigm |   c_1 x_1 + \cdots + c_d x_d \equiv 0 \mod N    \right\},
\end{equation}
where $N$ is a fixed positive integer.    Then $C$ is a lattice, and 
\[
\det C = N. 
 \]
 \end{lem}
 \begin{proof}
We  fix a basis $\{v_1, \dots, v_{d-1}\}$ for the $(d-1)$-dimensional sublattice defined by 
$H:= \left\{  x \in \Z^d  \mid   c_1 x_1 + \cdots + c_d x_d =0   \right\}$, and we adjoin to this basis one new vector, namely any integer vector $w$ that translates $H$ to its nearest discrete hyperplane companion 
\[
H+w:= \left\{  x \in \Z^d  \bigm |   c_1 x_1 + \cdots + c_d x_d =N   \right\}. 
\]
  Together, the set of vectors 
$ \{  v_1, \dots, v_{d-1}, w   \} $
form the edge vectors of a fundamental parallelepiped $\P$ for the lattice $C$, whose hight $\delta$ is the distance between these two parallel
 hyperplanes $H$ and $H+w$.  Using \eqref{distance between two hyperplanes}, we can
may calculate the volume  of $\P$ (which is by definition equal to $\det C$)  by multiplying the volume of its `base' times its `height':
\begin{align} 
\det C &= (\text{volume of the base of } \P)(\text{height of } \P) = (\det H)\delta \\    \label{proof of det C = N}
&= (\det H) \frac{N}{\sqrt{c_1^2 + \cdots + c_d^2}} = N,
\end{align}
using the fact that $\det H = \sqrt{c_1^2 + \cdots + c_d^2}$ from 
Lemma \ref{wonderful hyperplane determinant formula}.
\end{proof}


\bigskip
\section{The Gram matrix}

There is another very natural matrix that we may use to study lattices, which we can motivate as follows.   
 Suppose we are given any basis for a lattice $\L\subset \R^d$, say $\beta:= \{ v_1, \dots, v_r \}$, 
where $1 \leq r \leq d$.   By definition $\L = M(\Z^d)$, and ${\rm rank}(\L) = r$,  
where the columns of $M$ are defined by the basis vectors from $\beta$, and so $M$ is a $d\times r$ matrix. 
We can therefore represent any $x\in \R^d$ uniquely in terms of the basis $\beta$ like this:
 \begin{equation} \label{writing a vector in terms of a basis}
 x =  c_1 v_1 + \cdots + c_r v_r, 
 \end{equation}
 and the squared length of $x$ is:
 \begin{equation} \label{Gram matrix positive semidefinite}
 \| x \|^2 = \left\langle  \sum_{j=1}^r c_j v_j, \,  \sum_{k=1}^r c_k v_k \right\rangle =
 \sum_{1\leq j, k \leq r} c_j c_k \langle v_j, \,   v_k \rangle := c^T M^T M  c,
 \end{equation}
 where $c:= (c_1, \dots, c_r)^T$ is the coefficient vector defined by
 \eqref{writing a vector in terms of a basis}.
 
   It's therefore very natural to focus on the matrix $M^T M$, whose entries are the inner products
   $\langle v_j, v_k \rangle$ of all the basis vectors of the lattice $\L$, so 
 we define 
 \[
 G:= M^TM, 
 \]
 a {\bf Gram matrix} for $\L$.  It's clear from the computation above in 
 \eqref{Gram matrix positive semidefinite} that $G$ is positive definite.  Although
  $G$ does depend on which basis of $\L$ we choose, it is an elementary fact that $\det G$ is independent of the basis of $\L$.

Because we are always feeling the urge to learn more Linear Algebra, we would like to see
why any real symmetric matrix $B$ is the Gram matrix of some set of vectors.  To see this,  we apply the
  Spectral Theorem:   $B = P D P^T$, for some orthogonal matrix $P$ and a diagonal matrix $D$ with nonnegative diagonal elements.  So we can write $B = (P \sqrt D) (P \sqrt D)^T := M^T M$, where we defined the matrix $M:= (P \sqrt D)^T$, so that the columns of $M$ are the vectors whose corresponding dot products form the symmetric matrix $B$, and now $B$ is a Gram matrix.

  To review some more linear algebra, suppose we are given a real symmetric matrix $A$.
We recall that such a matrix is called {\bf positive definite} if in addition we have the positivity condition
\[
x^T A x > 0,
\]
for all  $x \in \R^d$.  Equivalently, all of the eigenvalues of $A$ are positive.  The reason is easy: 
$Ax = \lambda x$ for a non-zero vector $x \in \R^d$ 
implies that 
\[
x^T A x := \langle x, Ax \rangle = \langle x, \lambda x \rangle = \lambda \| x \|^2,
\]
so that  $x^T A x >0 $ if and only if $\lambda > 0$.     In the sequel, if we only require a symmetric matrix
$A$ that  enjoys the property $x^T A x \geq 0$ for all
$x\in \R^d$, then we call such a matrix {\bf positive semidefinite}.

Also, for a full-rank lattice $\L:= M(\Z^d)$, we see that $B:= M^T M$ will be positive definite if and only if $M$ is invertible, so that the columns of $M$ are a basis of $\L$.  Since a positive definite matrix is symmetric by definition, we've proved:
  
\begin{lem}  Suppose we are given a real symmetric matrix $B$.   Then:
\begin{enumerate}[(a)]
\item  $B$ is positive definite if and only if it is the Gram matrix of a full-rank lattice.
\item $B$ is positive semidefinite if and only if it is the Gram matrix of some set of vectors. 
\end{enumerate}
\hfill $\square$
\end{lem}
  
 What about reconstructing a lattice, knowing only one of its Gram matrices?  This is almost possible to accomplish, up to an orthogonal transformation, as follows.

 \begin{lem}\label{reconstructing a lattice basis from its Gram matrix}
 Suppose that $G$ is an invertible matrix, whose spectral decomposition is 
 \[
 G = P D P^T.  
 \]
 Then 
 \begin{equation} 
 G = X^T X \quad  \iff    \quad  X = Q \sqrt{D} P^T, 
 \end{equation}
for some orthogonal matrix $Q$.
 \end{lem}
\begin{proof}
The assumption $G = X^T X$ guarantees that $G$ is symmetric and has positive eigenvalues, so by the Spectral Theorem we have:
\[
G = P D P^T,
\] 
where $D$ is a diagonal matrix consisting of the positive eigenvalues of $G$, and $P$ is an orthogonal matrix consisting of eigenvectors of $G$.   Setting $X^T X = P D P^T$, we must have 
\begin{equation} \label{technical orthogonal identity}
I = X^{-T} PDP^T X^{-1} = (X^{-T} P \sqrt D) (X^{-T} P \sqrt D)^T, 
\end{equation}
where we define $\sqrt D$ to be the diagonal matrix whose diagonal elements are the positive square roots of the eigenvalues of $G$.  From \ref{technical orthogonal identity}, it follows that 
$X^{-T} P \sqrt D$ is an orthogonal matrix, let's call it $Q^{-T}$.  Finally, 
$X^{-T} P \sqrt D = Q^{-T}$ implies that $X = Q \sqrt D P^T$.
\end{proof} 
 
 So Lemma \ref{reconstructing a lattice basis from its Gram matrix} allows us to reconstruct a lattice $\L$, up to an orthogonal transformation, by only knowing one of its Gram matrices.



\bigskip 
\section{Dual lattices}  \index{dual lattice} \label{dual lattice}

Every  lattice $\L:= M(\Z^d)$ has a {\bf dual lattice}, which we have already encountered in the Poisson summation formula
for arbitrary lattices.   The dual lattice of a full-rank lattice $\L$ was defined by:
\begin{equation}\label{first definition of dual lattice}
\L ^*= M^{-T}(\Z^d).
\end{equation}
But there is another way to define the dual lattice of a lattice $\L \subset \R^d$ (of any rank), which is coordinate-free:
\begin{equation}\label{second definition of dual lattice}
\L^* := \left\{   x \in \R^d \mid \langle x, n \rangle \in \Z, \text{ for all }  n \in \L      \right\}.
\end{equation}

\begin{lem} 
The two definitions above, \eqref{first definition of dual lattice} and \eqref{second definition of dual lattice}, are equivalent.
\begin{proof}
We let $A := \L^*:= M^{-T}(\Z^d)$, and
 $B:=  \left\{   x \in \R^d \bigm |     \langle x, n \rangle \in \Z, \text{ for all }  n \in \L      \right\}$.
We first fix any $x \in A$.   To show $x \in B$, we fix any
 $n\in \L$, and we now have to verify that $\langle x, n \rangle \in \Z$.  
By assumption, $x = M^{-T}m$ for some $m \in \Z^d$, and $n= M k$, for some $k\in\Z^d$.  Therefore 
\[
\langle x, n \rangle =\langle     M^{-T}m, n \rangle =  \langle   m, M^{-1}n \rangle 
=           \langle   m, k \rangle   \in \Z, 
\]
because both $m,k \in \Z^d$.    So we have $A \subset B$.
For the other direction, suppose that $y\in B$, so by definition 
\begin{equation*}\label{def. of B}
\langle y, n \rangle \in \Z,   \text{  for all } n \in \L.
\end{equation*}
We need to show that $y = M^{-T} k$ for some $k\in \Z^d$, which is equivalent to $M^T y \in \Z^d$.   
Noticing that the $k$'th element of $M^T y$ is  
$\langle n, y \rangle$ with $n$ belonging to a basis of $\L$, we are done, by \eqref{def. of B}.
Therefore $A=B$.
\end{proof}
\end{lem}

\bigskip
\begin{example}  
\rm{
Let $\L := r \Z^d$, the integer lattice dilated by a positive real number $r$.  It's dual lattice is $\L^* = \frac{1}{r} \L$, because a basis for $\L$ is $M := rI$, implying that a basis matrix for $\L^*$ is $M^{-T} = \frac{1}{r} I$.
We also notice that $\det \L = r^d$, while $\det \L^* = \frac{1}{r^d}$.
}
\hfill $\square$
\end{example}

A fundamental relation between a full-rank lattice and its dual follows immediately from Definition 
\ref{first definition of dual lattice}:  $\det(\L^*) := \det(M^{-T}) = \frac{1}{\det M}= \frac{1}{\det \L}$, 
which we record as:
\begin{equation} \label{FundamentalDualRelation}
(\det \L )(\det \L^*) = 1.
\end{equation}
If we consider any integer sublattice of $\Z^d$, say $\L \subset \Z^d$, together with its dual lattice $\L^*$ in the same space, some interesting relations unfold between them.   Let's consider an example. 
 
\bigskip
\begin{example} \label{dual lattice example}
 \rm{
In $\R^2$, let 
$\L := \left\{m   \icol{1\\1}    + n  \icol{1\\  4} \bigm |    m,n \in \Z   \right\}$, 
a lattice with $\det \L = 3$ that is depicted in Figure \ref{DualLattice} by the larger green balls. 
Its dual lattice is
\[
\L^*:= \left\{
\frac{1}{3} \left(   a\icol{\  4\\-1} +b    \icol{ \ \,  -1\\  \  \  \ 1}   \right)    \bigm |    a, b \in \Z
\right\}, 
\]
whose determinant equals $\frac{1}{3}$, and is depicted in Figure 
 \ref{DualLattice} by the smaller orange balls.   So $\L$ is a  coarser lattice than $\L^*$.  That is,
 $\L^* \supseteq \L$.
}

\begin{figure}[htb]
 \begin{center}
\includegraphics[totalheight=4.3in]{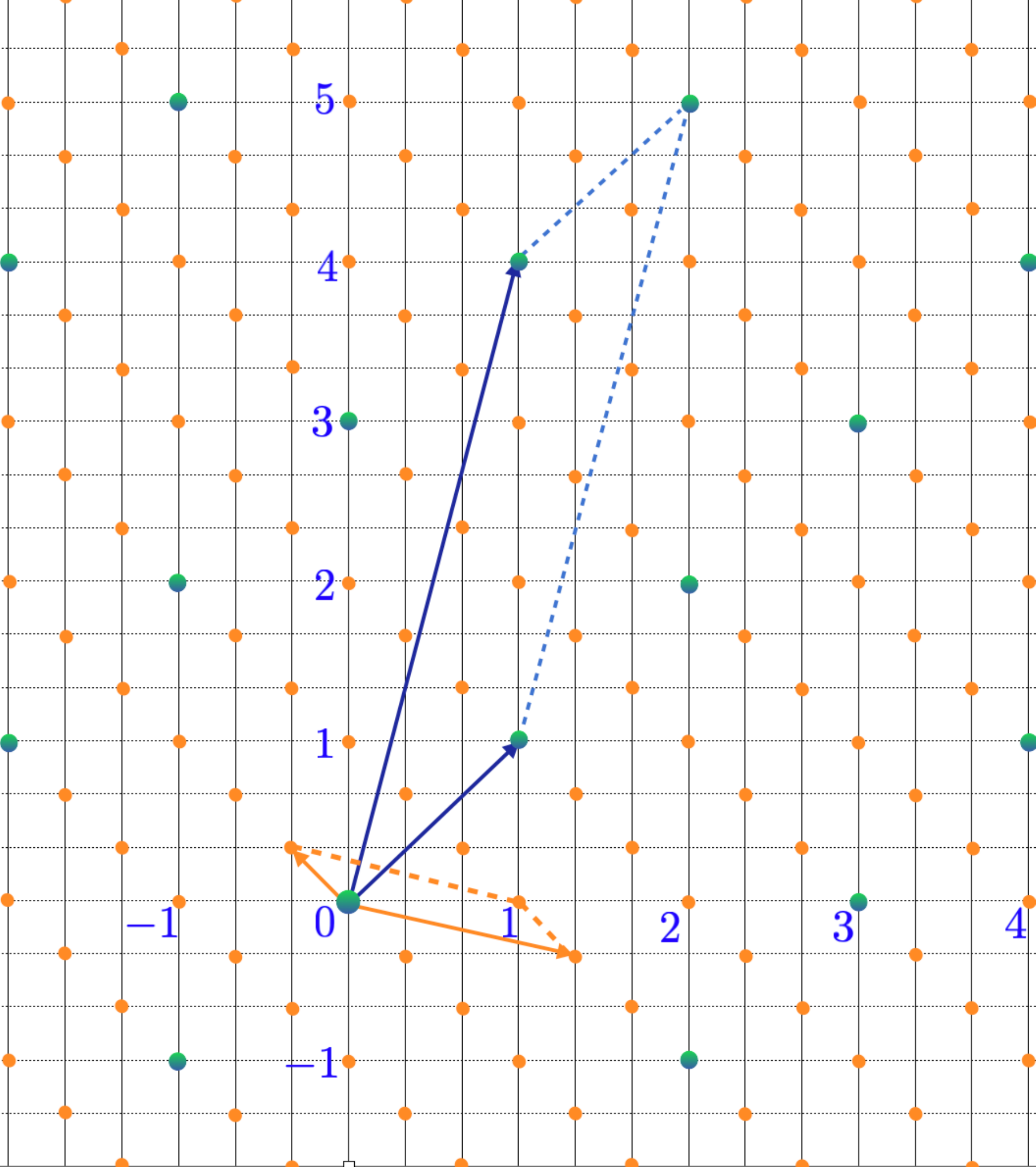}
\end{center}
\caption{The lattice $\L$ of Example \ref{dual lattice example}, depicted by the green points.  Its dual lattice  $\L^*$ is depicted by the orange points.  Here $\L^* \supset \L$ and is therefore 
a finer lattice.}
\label{DualLattice}
\end{figure}

We can verify the relation \eqref{FundamentalDualRelation} here: $\det \L^* = \frac{1}{3} = \frac{1}{\det \L}$.    
We may notice here that
$\L^*/\L$ forms a finite group of order $9 = (\det \L)^2$, which is equal to the number of cosets of
the coarser lattice $\L$ in the finer lattice $\L^*$.
\hfill $\square$
\end{example}

\begin{question}[Rhetorical]
When is it true that $\L^* \supseteq \L$?  In other words, for which lattices is 
the dual lattice a refinement of the original lattice?
\end{question}
To study this phenomenon, a lattice  $\L \subset \R^d$  is called an {\bf integral lattice} 
if 
\[
\langle x, y \rangle \in \Z \text{ for all } x, y \in \L.   
\]
It follows directly from our definition of an integral lattice that  
\[
\L^* \supseteq \L \iff \L \text{ is an integral lattice},
\]
and in this case we have a finite abelian group $\L^* / \L$, called the {\bf discriminant group}.
\begin{lem}\label{lem:discriminant lattice}
For a full-rank integral lattice $\L$, we have:
\begin{equation}
\left | \L^* / \L \right | = (\det \L)^2.
\end{equation}
\end{lem}
\begin{proof}
we recall Theorem \ref{sublattice index}:
\[
\left | \L^* / \L \right | = \frac{ \det \L}{ \det \L^* }
= \frac{ \det \L}{ (1/\det \L) } = (\det \L)^2.
\]
\end{proof}
\begin{example}
{\rm
Clearly, any integer sublattice $\L \subseteq \Z^d$ is also an integral lattice.  But there are others, as the next example shows.
}
\hfill $\square$
\end{example}
\begin{example}\label{orthogonal transformation of Z^d}
{\rm
To see different kinds of integral lattices, we can take any orthogonal linear transformation of
$\Z^2$.  Let's fix an angle $0< \theta < \frac{\pi}{2}$ with $\cos \theta$ and $\sin \theta$ irrationals and linearly independent over the rationals.  We define:
\[
M:= \big(\begin{smallmatrix}
\cos \theta & \ \   -\sin \theta  \\
\sin \theta & \ \ \  \cos \theta
\end{smallmatrix}
\big),
\ \ \ 
\L:= \left\{ 
M\icol{a\\b} \mid \icol{a\\b}\in \Z^2
\right \}.
\]
For any two lattice vectors $u, v \in \L$, we have 
$\langle u, v \rangle = \langle M \icol{a\\b}, M \icol{c\\d} \rangle = 
\langle M^T M \icol{a\\b},  \icol{c\\d} \rangle = 
\langle \icol{a\\b},  \icol{c\\d} \rangle \in \Z$, so that 
our lattice $\L$ is also an integral lattice.  We notice, though,  that in this example $\L$  has no nonzero integer vectors at all !   

Let's compute the dual lattice here:  $\L^*$ is given by the matrix 
$M^{-T} = M$.  In other words, here we have  $\L^* = \L$.  Is  this a coincidence?
}
\hfill $\square$
\end{example}

\begin{example} \label{dual lattice example 2}
 \rm{
Continuing with Example \ref{dual lattice example}, we have an integral lattice $\L\subset \R^2$, so that $\L^* \supseteq \L$.  Here $\det \L=3$, $\det \L^* = \frac{1}{3}$, confirming the claim from
Lemma \ref{lem:discriminant lattice}, namely that 
\[
\left | \L^* / \L \right | = \frac{3}{\left(\tfrac{1}{3}\right)} = 9 = (\det \L)^2.
\]
}
\hfill $\square$
\end{example}

Next, we call $\L \subset \R^d$ a {\bf unimodular lattice} if $\det \L = 1$.   The collection of all  unimodular lattices is quite important in number theory 
and we'll see it again later, in  Siegel's mean value theorem \ref{thm:Siegel mean values}.  
We say that a lattice $\L$ is {\bf self dual} if $\L^* = \L$.   Chasing these elementary
 ideas around, the following observation is immediate (Exercise \ref{proof of equivalence of a self-dual lattice}).
\begin{lem}\label{equivalence of a self-dual lattice}
The following are equivalent:
\begin{enumerate}[(a)]
\item $\L$ is self-dual.
\item $\L$ is an integral unimodular lattice.
\end{enumerate}
\hfill $\square$
\end{lem}

\begin{example}
{\rm
Continuing with Example \ref{orthogonal transformation of Z^d}, we saw that the lattice defined by
\[
\L:= \left\{ 
M\icol{a\\b} \mid \icol{a\\b}\in \Z^2
\right \},
\]
with $M:= \big(\begin{smallmatrix}
\cos \theta & \ \   -\sin \theta  \\
\sin \theta & \ \ \  \cos \theta
\end{smallmatrix}
\big)$, was self-dual.  According to Lemma \ref{equivalence of a self-dual lattice}, $\L$ 
should therefore also be an integral, unimodular lattice - and indeed it is. 
}
\hfill $\square$
\end{example}



\bigskip
\section{Some important lattices}

Throughout this section, we'll fix the special vector  
\begin{equation}\label{vector of all 1/2}
w:= 
\begin{pmatrix}
 \tfrac{1}{2} \\
 \\
\tfrac{1}{2} \\
\\
\vdots \\
\\
\tfrac{1}{2} \\
\end{pmatrix}
\in \R^d.
\end{equation}

\begin{example}\label{D_n lattices}
 \rm{
The $D_n$ lattice is defined by 
\[
D_n:= \left\{  x\in \Z^n \bigm |      \sum_{k=1}^n x_k \equiv  0 \mod 2   \right\}. 
\]

\begin{figure}[htb]
 \begin{center}
\includegraphics[totalheight=3.2in]{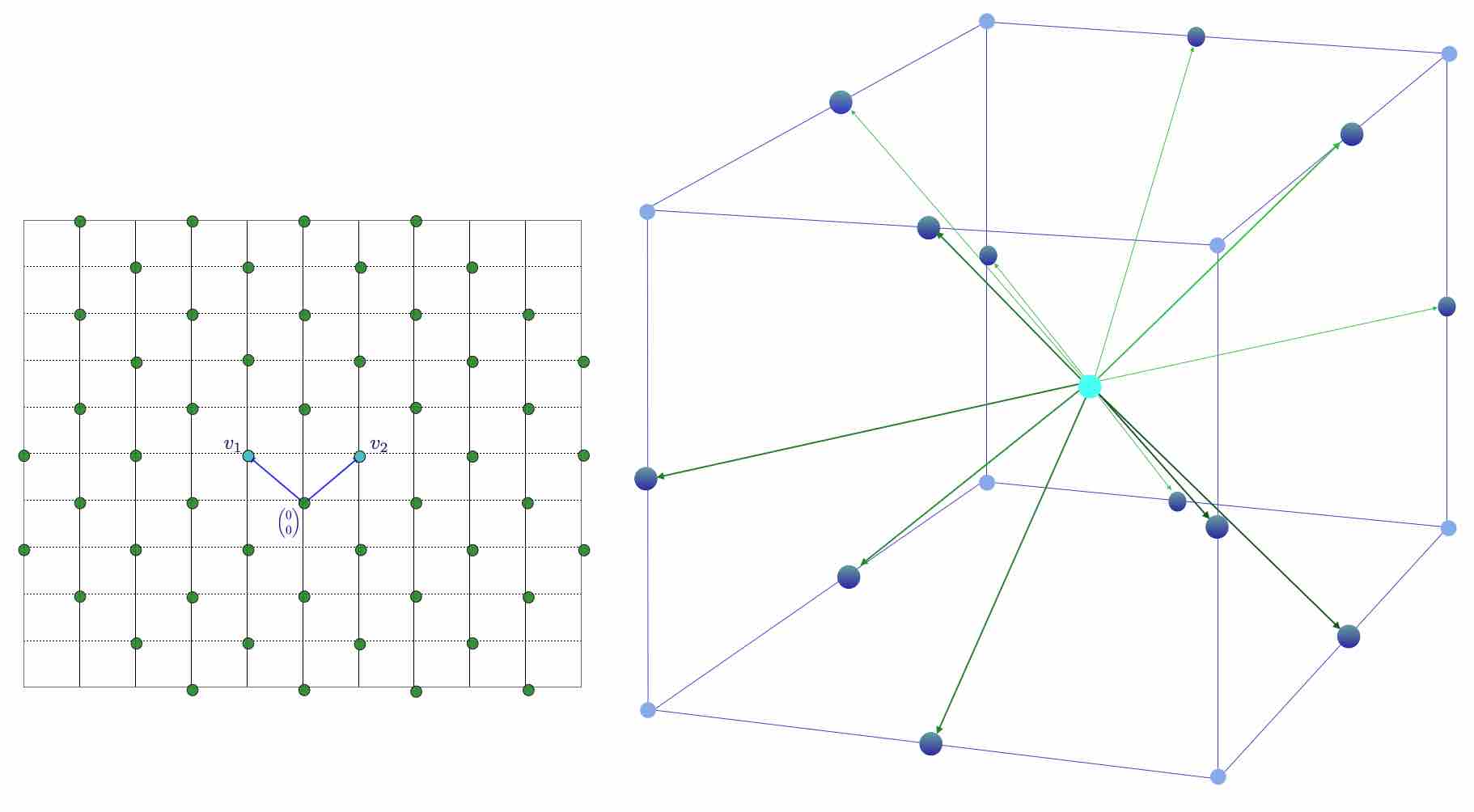}
\end{center}
\caption{Left: the $D_2$ lattice.  Right:  the $12$ shortest nonzero vectors of the $D_3$ lattice, inscribed in the cube $[-1, 1]^3$.}  
\label{The $D_2$ and $D_3$ lattices}
\end{figure}

In $\R^4$, the $D_4$ lattice turns out to be a fascinating object of study.
The Voronoi cell  $\text{Vor}_0(D_4)$ is called the {\bf $24$-cell}, \index{$24$-cell} 
and is depicted in Figure \ref{24-cell}.  It is a $4$-dimensional polytope with some incredible properties - for example, it is one of the few polytopes that is self-dual.  It is also an example of a polytope $\P$ in the lowest possible dimension $d$ (namely $d=4$) such that $\P$ tiles $\R^d$ by translations, 
and yet $\P$ is not a zonotope. 

By Lemma \ref{lemma:lattice defined by congruence}, 
we see that $\det D_n = 2$.   
The lattice $D_n$ is often called the ``checkerboard''  lattice, because $\det D_n = 2$ means there are exactly two cosets
in $\Z^d / D_n$.   Finally, the dual lattice $D_n^*$  is equal to the lattice
\[
\\Z^d  \cup \left( \Z^d +  w \right),
\]
which we leave for the pleasure of the reader (Exercise \ref{dual of D_n}).  
}
\hfill $\square$
\end{example}

\begin{example}\label{def:E_8 lattice}
\rm{
The $E_8$ lattice is defined by 
\[
E_8:= D_8 \cup (D_8 + w),
\]
with $w$ defined in \eqref{vector of all 1/2}.
It's a nice exercise to show that the latter definition in fact gives us a lattice (Exercise \ref{E_8 is a lattice}).
 It turns out that $E_8$  is a self-dual lattice (Exercise \ref{E_8 is self-dual}). $E_8$ is also an even, unimodular lattice. Moreover, it also turns out
 that $E_8$ gives the optimal solution to the sphere packing problem in $\R^8$ (See Chapter \ref{Sphere packings}).  This lattice has amazing symmetries, and
is important in the physics of string theory, as well as data transmission.
}
\hfill $\square$
\end{example}

\begin{example}
\rm{
We define the hyperplane $H:= \{ x \in \Z^8 \mid  x_1+ \cdots + x_8 = 0\}$ in $\R^8$. 
Then we have the lattice
\[
E_7:= E_8 \cap H,
\]
which has rank $7$.
We also fix $V:= \{ x \in \R^8 \mid x_2+ x_3 + x_4 + x_5 + x_6 + x_7 = x_1 + x_8 = 0 \} \subset \R^8$, a vector subspace of dimension $6$, and we define 
\[
E_6:= E_8 \cap V,
\]
a lattice of rank $6$. Both $E_6$ and $E_7$ are, almost by definition, sublattices of $E_8$.
}
\hfill $\square$
\end{example}

\begin{example}
\rm{
Last but not least is the famous Leech lattice in $\R^{24}$.  It is the unique even unimodular lattice in $\R^{24}$ that has no vectors of length
$\sqrt 2$ (see \cite{CasselsBook} for a proof, among other constructions for the Leech lattice).   
There are many constructions of the Leech lattice, none of which are trivial, and one of which involves the important Golay binary code \cite{ConwaySloan.book}.  Coding theory, which is a discrete version of sphere packing, is a fascinating and important topic for another day. 
}
\hfill $\square$
\end{example}



\section{The Hermite normal form}

If a lattice satisfies $\L \subset \Z^d$, we'll call it an integer sublattice.  
We may recall that any lattice $\L \subset \R^d$ has infinitely many bases, so it may seem impossible at first to associate a single matrix with a given lattice.   However, there is an elegant way to do this, as follows. 

\begin{example}\label{first ex. of HNF}
{\rm
Suppose we are given a lattice $\L$ as the integral span of the vectors 
\[
v_1:= \icol{3\\1}, v_2:= \icol{-2\\    \   2},
\]
which clearly has determinant $8$. 
Then any integer linear combinations of $v_1$ and $v_2$ is still in $\L$.  In particular, mimicking Gaussian elimination, we place $v_1$ and $v_2$ as rows of a matrix, and row-reduce over the integers:
\begin{align*}
\begin{pmatrix}
\  \ 3 & \  1 \\
-2 &  \ 2 
\end{pmatrix}
\rightarrow 
\begin{pmatrix}
3 & \ 1 \\
1 &  \ 3 
\end{pmatrix}
\rightarrow 
\begin{pmatrix}
0 & \ -8 \\
1 &  \  \  \ 3 
\end{pmatrix}
\rightarrow 
\begin{pmatrix}
1 & \  \ 3 \\
0 &   -8 
\end{pmatrix}
\rightarrow 
\begin{pmatrix}
1 & 3 \\
0 &   8 
\end{pmatrix},
\end{align*}
where at each step we performed row operations (over $\Z$) that did not change the lattice.  Hence we have a reduced basis for $\L$, consisting of $\icol{1\\3}$ and $\icol{0\\8}$. 

\begin{figure}[htb]
 \begin{center}
\includegraphics[totalheight=2.2in]{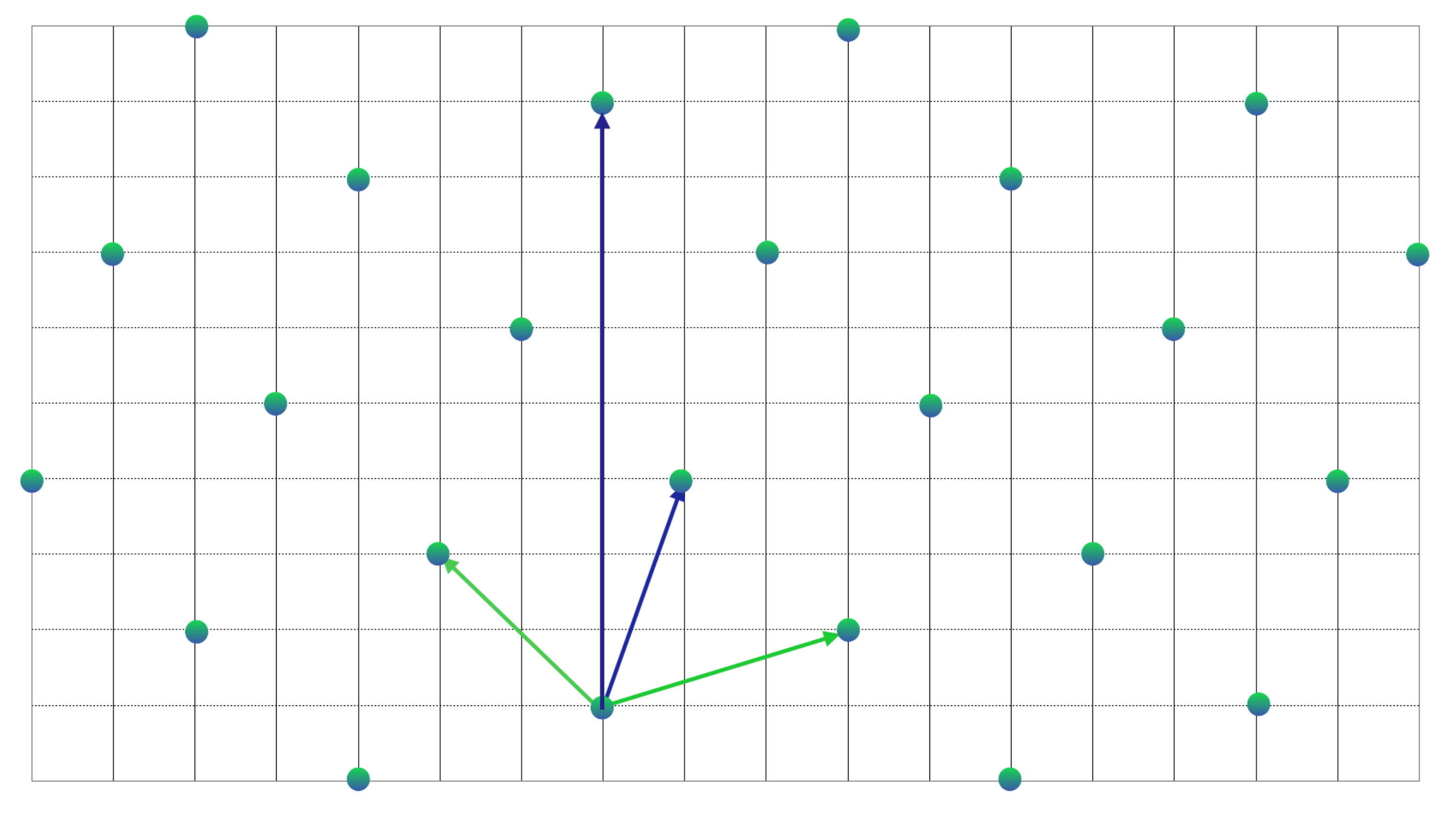}
\end{center}
\caption{
The lattice $\L$ of Example \ref{first ex. of HNF}, depicted by the bold green points, and showing the 
original basis $\{ v_1, v_2\}$ of $\L$,  and the Hermite-reduced basis of $\L$.  Here it's geometrically clear that both are bases for the same lattice $\L$.
 }  
\label{HNF.pic}
\end{figure}

We notice that the resulting matrix is upper-triangular, with positive integers on the diagonal, nonnegative integers elsewhere, and in each column the diagonal element is the largest element in that column.  

There is another way to interpret the matrix reductions above, by using unimodular matrices, as follows.  The first reduction step can be accomplished by the multiplication on the left by a unimodular matrix:
\[
\begin{pmatrix}
1 & \  0 \\
1 &  \ 1 
\end{pmatrix}
\begin{pmatrix}
\ 3 & \  1 \\
-2 &  \ 2 
\end{pmatrix}
=
\begin{pmatrix}
3 & \ 1 \\
1 &  \ 3 
\end{pmatrix}
\]
Similarly, each step in the reduction process can be interpreted by multiplying on the left by some new unimodular matrix, so that at the end of the process we have a product of unimodular matrices times 
our original matrix $\begin{pmatrix}
\ 3 & \  1 \\
-2 &  \ 2 
\end{pmatrix}$.  Because a product of unimodular matrices is yet another unimodular matrix, we can see that we arrived at a reduction of the form:
\[
U
\begin{pmatrix}
\ 3 & \  1 \\
-2 &  \ 2 
\end{pmatrix}
=
\begin{pmatrix}
1 & \ 3 \\
0 &   8 
\end{pmatrix},
\]
where $U$ is a unimodular matrix.
\hfill $\square$
}
\end{example}
The point of Example \ref{first ex. of HNF} is that  a similar matrix reduction persists for all integer lattices, culminating in the following result, which just hinges on the fact that $\Z$ has a division algorithm.

\begin{thm} \label{theorem.HNF}
Given an invertible integer $d\times d$ matrix $M$,  there exists a unimodular matrix $U$ with 
$UM = H$, such that $H$ satisfies the following conditions:
\begin{enumerate}
\item  $[H]_{i, j} = 0$ if $i>j$.
\item  $[H]_{i, i} > 0$, for each $1\leq i \leq d$.
\item  $0 \leq [H]_{i, j}  < [H]_{i, i} $, for each $i >j $. \label{third property}
\end{enumerate}
Property \ref{third property} tells us that each diagonal element $[H]_{i, i}$ in the $i$'th column of $H$ is the largest element in the $i$'th column. 

Moreover, the matrix $H$ is the unique integer matrix that satisfies the above conditions.   
\hfill $\square$
\end{thm}
 The matrix $H$ in Theorem \ref{theorem.HNF} is called the {\bf Hermite normal form} of  $M$.
To associate a unique matrix to a given integral full-rank lattice $\L \subset \R^d$, we first choose
 any basis of $\L$, and we then construct a $d\times d$ integer matrix $M$ whose rows are the basis vectors that we chose.   We then apply Theorem \ref{theorem.HNF} to $M$, arriving at an integer matrix $H$ whose rows are another basis of $\L$, called the {\bf Hermite-reduced basis}. 

\begin{cor}
There is a one-to-one correspondence between full-rank integer sublattices in $\R^d$ and 
integer $d \times d$ matrices  in their Hermite Normal Form.
\hfill $\square$
\end{cor}

\begin{example}
Given any $2$-dimensional lattice $\L\subset \Z^2$, with a basis matrix $M$, we can use the Hermite-normal form of $M$ to get the following basis for $\L: \{ \icol{1\\0},  \icol{p\\q} \}$, for some 
nonnegative integers $p, q$ (Exercise \ref{ }).
\hfill $\square$
\end{example}



\section{The Voronoi cell of a lattice}

The {\bf Voronoi cell} of  a lattice $\L$, at the origin, is defined by
\begin{equation}
\text{Vor}_0(\L) := \left\{ x \in \R^d \bigm |     \|x\|  \leq   \|x - v\|, \ \text{ for all  } v \in \L    \right\}.
\end{equation}
In other words, the Voronoi cell  $\text{Vor}_0(\L)$ of a lattice $\L$ is the set of all point in space that are closer to the origin than to any other lattice point in $\L$.   Because the origin  wins the battle of minimizing this particular distance function,  it is also possible to  construct the Voronoi cell by using half-spaces.  Namely, for each $v\in \L$, we define the half-space 
 \[
 H_v:= \left\{ x \in \R^d \bigm |   \langle x, v \rangle  \leq   \tfrac{1}{2} \|v\|    \right\},
 \]
 and we observe that the Voronoi cell may also be given by
 \[
{\rm Vor}_0(\L)  = \bigcap_{v\in \L- \{0\}} H_v,
 \]

\begin{figure}[htb]
 \begin{center}
\includegraphics[totalheight=4.2in]{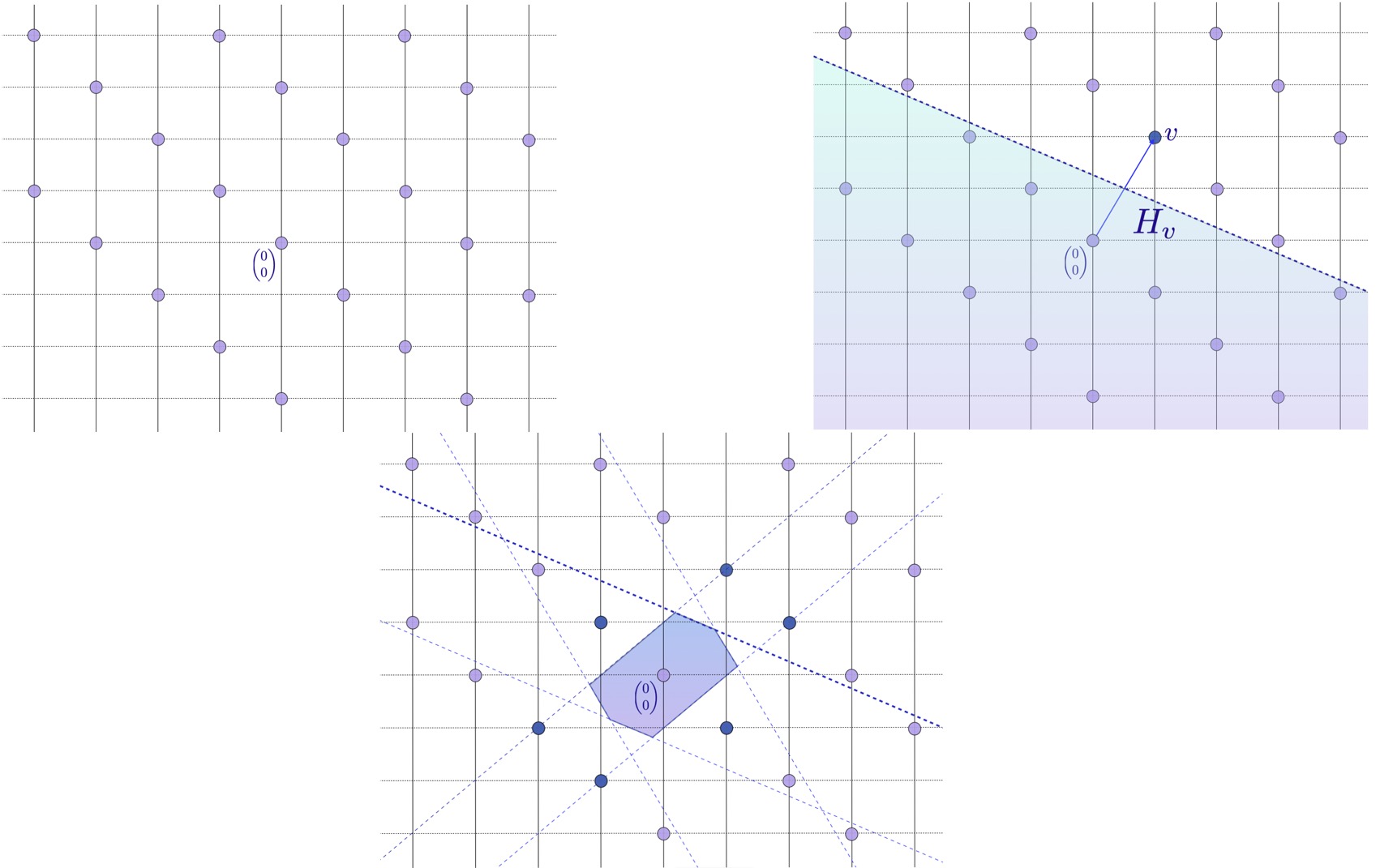}
\end{center}
\caption{Top left:  a sublattice $\L$ of $\Z^2$, of index $3$. Top right: $v \in \L$ is one of the $6$ relevant vectors, with its corresponding half-plane $H_v$, helping to define the Voronoi cell at the origin.
Bottom: The Voronoi cell $\text{Vor}_0(\L)$, \index{Voronoi cell} a symmetric hexagon of area $3$, with its $6$ relevant (heavy blue) lattice points of $\L$.
 }  
\label{VoronoiConstruction}
\end{figure}
 
as drawn in Figure \ref{VoronoiConstruction}.
It is easy to observe that the Voronoi cell of a lattice is symmetric about the origin, convex, and compact (Exercise \ref{facts about Vor cell}).  So we may expect that Minkowski's theorems apply to $\text{Vor}_0(\L)$, as we see in the proof of Lemma \ref{basic Voronoi lemma} below.
It's also useful to define an analogous Voronoi cell located at each lattice point $m \in \L$:
\begin{equation}
\text{Vor}_m(\L) := \left\{ x \in \R^d \bigm |     \|x - m \|  \leq   \|x - v\|, \ \text{ for all  } v \in \L    \right\}.
\end{equation}
A moment's thought (but this is good practice - Exercise \ref{translating the Voronoi cell around})
reveals that a translation of the Voronoi cell at the origin is exactly the Voronoi cell at another lattice point of $\L$, namely:   
\begin{equation} \label{translated Voronoi cells}
\text{Vor}_0(\L) + m =  \text{Vor}_m(\L).
\end{equation}
\begin{lem} \label{basic Voronoi lemma}
Given a full-rank lattice $\L\subset \R^d$, whose Voronoi cell at the origin is $K$, we have:
\begin{enumerate}[(a)]
\item   $K$ tiles $\R^d$ by translations with $\L$. \label{part 1 of Voronoi}
\item
$
\vol (K)  = \det \L.     \label{part 2 of Voronoi}
$
\end{enumerate}
\end{lem}
\begin{proof}
Part  \ref{part 1 of Voronoi} follows from  the observation that
any $x \in \R^d$, there exists a lattice point $m \in \L$ that is at least as close to $x$ as it is to any other lattice point of $\L$. 
In other words, $  \|x - m\|  \leq   \|x - v\|, \forall v\in \L$, and so 
$x \in \text{Vor}_m(\L)$.
From \eqref{translated Voronoi cells} we see that $x$ is covered by the translate $\text{Vor}_0(\L) + m$. It's also clear that as $n$ varies over $\L$, all of the interiors of the translates $\text{Vor}_0(\L) + n$ are disjoint, so that 
$K:= \text{Vor}_0(\L)$ tiles $\R^d$ by translations with $\L$.
To prove part  \ref{part 2 of Voronoi},  we let $B:= 2K$.   
By Theorem \ref{thm:extremal bodies} (regarding extremal bodies), we know that 
$\frac{1}{2}B = K$ tiles $\R^d$ with the lattice $\L$ if and only if  $ \vol(B) = 2^d \det \L$.   Since 
\ref{part 1 of Voronoi} tells us that
$K=\frac{1}{2}B$ tiles with the lattice $\L$, 
 we see that
$\vol K = \vol\Big(\frac{1}{2}B\Big) = \frac{1}{2^d} \vol B = \det \L$.
\end{proof}
The proof above shows that the Voronoi cell  of $\L$ is also an extremal body for $\L$, according to Theorem \ref{thm:extremal bodies}.

A fascinating open problem is the Voronoi conjecture, named after 
the Ukrainian mathematician Georgy Voronoi, who formulated it in 1908.   
Two polytopes $\P, Q$ are called {\bf affinely equivalent} \index{affinely equivalent} if 
$\P = M(Q) + v$, where $M \in GL_d(\R)$, and $v\in \R^d$. 

\medskip
\begin{question}[The Voronoi conjecture]   \index{Voronoi conjecture}
Does a polytope $\P$  tile $\R^d$ by translations $\iff \P$ is the Voronoi cell of some lattice $\L$, or $\P$ is affinely equivalent 
to such a Voronoi cell?
\end{question}
The main difficulty in the Voronoi conjecture appears to be
 the apriori search among all of the (infinitely many) possible affinely equivalent images of such a Voronoi cell.

\medskip
\begin{example}
\rm{
For the lattice $A_n \subset \R^{n+1}$ defined in Example \ref{A_d example}, its Voronoi cell turns out to have beautiful and important properties:  $\text{Vor}(A_2) \subset \R^3$ is a hexagon, 
$\text{Vor}(A_3) \subset \R^4$ is a truncated octahedron (one of the Fedorov solids), and so on (see Conway and Sloane \cite{ConwaySloan.book}).
}
\end{example}

\begin{figure}[htb]
 \begin{center}
\includegraphics[totalheight=2in]{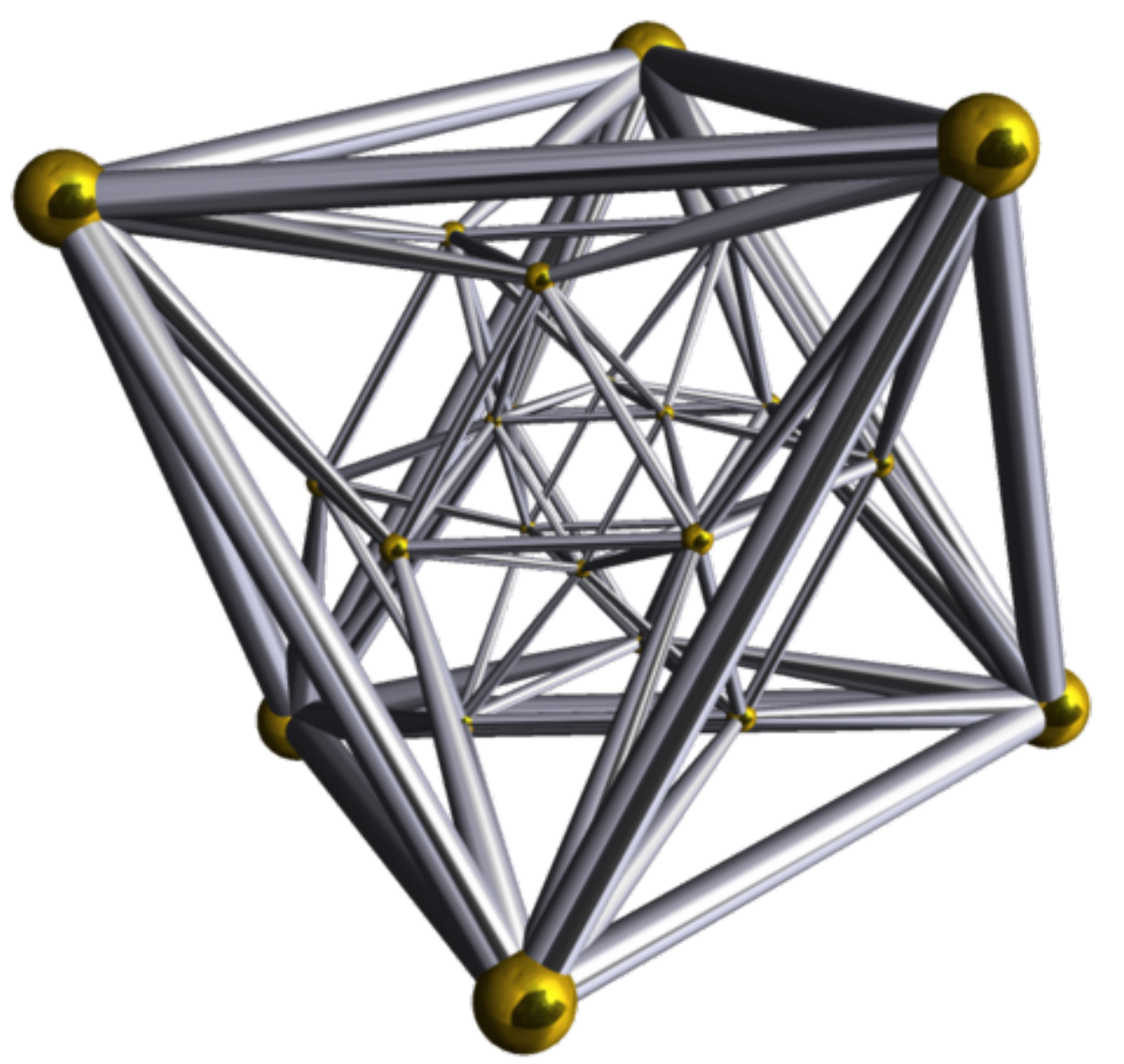}
\end{center}
\caption{The Voronoi cell of the $D_4$ lattice in $\R^4$, known as the $24$-cell.}  \label{24-cell}
\end{figure}

\bigskip
\section{Characters of lattices}\index{characters of a lattice}
\label{sec:characters of a lattice}

For each lattice point  $n \in \L$, we associate a function called a {\bf character} of $\L$, which we 
define by:
\begin{equation}\label{def: character of a lattice}
\chi_n(x) := e^{2\pi i \langle n, x \rangle},
\end{equation}
for all $x \in \R^d$.   
If we multiply these characters together by defining $\chi_n \chi_m:= \chi_{n+m}$, then
\[
G_\L:=\{ \chi_n \mid n \in \L\}
\]
 forms a group, under multiplication of functions, 
called  the {\bf group of characters} of $\L$.  
To see that this multiplication makes sense, we can compute:
\begin{equation}\label{Phi a homo}
(\chi_n \chi_m) (x) :=
e^{2\pi i \langle n, x \rangle} e^{2\pi i \langle m, x \rangle}
=e^{2\pi i \langle n+m, x \rangle} 
:= \chi_{n+m}(x).
\end{equation}
But much more is true.
\begin{thm}
\[
G_\L \simeq \L,
\]
an  isomorphism of groups.
\end{thm}
\begin{proof}
We consider the natural map $\Phi: \L \rightarrow  G_\L$ defined by $\Phi(n):= \chi_n$.  We'll show that
$\Phi$ is an isomorphism, so it is necessary to prove that $\phi$ is a bijective homomorphism, by definition.  By \eqref{Phi a homo} above, 
the computation 
\[
\Phi(n+m) = \chi_{n+m}= \chi_n \chi_m = \Phi(n) \Phi(m)
\]
already shows that $\Phi$ is a  homomorphism.  The surjectivity of $\Phi$ is clear from the definition of $\Phi$.  The more interesting direction is to show that $\Phi$ is injective.  First,
 $\Phi(n) = \Phi(m) \implies  e^{2\pi i \langle n, x \rangle}  =  e^{2\pi i \langle m, x \rangle}$ for all 
 $x \in \R^d \ \iff   e^{2\pi i \langle n-m, x \rangle} =1$ 
for all $x\in \R^d \  \iff$ 
\begin{equation} 
\label{using inner product is an integer, in Phi}
 \langle n-m, x \rangle \in \Z,   
 \end{equation}
for all  $x \in \R^d$.
Suppose to the contrary that $n-m \not=0$.  Consider the open ball 
\[
B:=\left \{ x \in \R^d \mid \| x\| < \frac{1}{\|n-m\|} \right \}, 
\]
and pick any $x \in B$.
The Cauchy-Schwartz inequality gives us:
\[
\left | \langle n-m, x \rangle \right | \leq \|n-m\| \|x\| < \|n-m\|  \frac{1}{\|n-m\|} =1,
\]
so that by \eqref{using inner product is an integer, in Phi} we now have $ \langle n-m, x \rangle=0$ for all $x \in B$.  But this implies $n-m=0$, a contradiction.
\end{proof}

Intuitively, one of the huge benefits of group characters is that by using the magic of just two-dimensional complex numbers, we can study high-dimensional lattices.
\begin{example}
For the integer lattice $\Z^d$, its group of characters comprises the following functions:
\[
\chi_n(x):= e^{2\pi i \langle n, x  \rangle},
\]
for each $n \in \Z^d$.  
\hfill $\square$
\end{example}

Now we allow ourselves the luxury of being slightly more general and free to think about any finite group.  Although we focused thus far on discrete (sub)groups in $\R^d$
defined in \eqref{discrete subgroup}, 
the reader may consult \cite{Herstein} for the definition of any group.

\begin{lem}\label{lem: orthogonality of characters of a lattice}
Let $G$ be any finite group, and $\chi: G\rightarrow \C\setminus \{0\}$ a nontrivial homomorphism of $G$.
\end{lem}
\begin{enumerate}[(a)]
\item \label{part 1 of lemma: vanishing group character sum}
We have:
\[
\sum_{g \in G} \chi(g) =0.
\]
\item \label{part 2 of lemma: vanishing group character sum}
For any two distinct homomorphisms $\chi, \psi:G\rightarrow \C\setminus \{0\}$, we have:
\[
\sum_{g \in G} \chi(g)\overline{\psi(g)}  =0.
\]
\end{enumerate}
\begin{proof}
To prove part \ref{part 1 of lemma: vanishing group character sum}, 
we first note that because $\chi$ is nontrivial, there exists a nonzero element 
$g_0 \in G$ such that $\chi(g_0) \not=1$.
We have:
\[
\chi(g_0)\sum_{g \in G} \chi(g)
=\sum_{g \in G}   \chi(g_0)  \chi(g)
=\sum_{g \in G}   \chi(g_0 g)
=\sum_{g \in G} \chi(g),
\]
where the last step follows from the fact that multiplication by $g_0$ permutes all the elements of $G$ (Exercise \ref{permuting the elements of a group}).  
So we now have $(\chi(g_0)-1)\sum_{g \in G} \chi(g)=0$,
 and because $\chi(g_0)\not=1$, 
 we conclude that $\sum_{g \in G} \chi(g)=0$.

To prove part \ref{part 2 of lemma: vanishing group character sum}, we define $\phi:= \chi \psi^{-1}$, which is another character of $G$, where we have 
$\phi(x):= \chi(x) \psi^{-1}(x) =  \chi(x) \overline{\psi(x)}$, for all $x\in G$.  Moreover, by assumption $\phi$ is not the trivial character ($\phi$ and $\psi$ are distinct), so that part \ref{part 1 of lemma: vanishing group character sum} applies to the character $\phi$ and we're done.
\end{proof}

For some applications, it's useful to somehow transfer the problem of summing a function over a sublattice (or superlattice) of $\L$, to the problem of summing essentially the same function over $\L$.
The following is the  classical orthogonality relation for a finite abelian group
$\Z^d/ M\Z^d$, but we prefer to phrase it in terms of the lattice $\L:= M(\Z^d)$, for future applications to lattices.
\begin{cor}[Orthogonality relations for characters of a lattice]
\label{Orthogonality relations for characters of a lattice}
Let $\L\subset \Z^d$ be a full-rank integer sublattice, so we may write $\L:= M(\Z^d)$, with $M$ an invertible integer matrix.  Then we have:
\begin{equation}
\label{indicator function of a lattice, in terms of characters, I}
\frac{1}{|\det M|} \sum_{g \in \Z^d/ M\Z^d}   e^{2\pi i \langle M^{-T} g, m \rangle}
=
\begin{cases}  
1    &      \mbox{if }  m \in \L\\ 
0  &        \mbox{if } m \notin \L,
\end{cases}
\end{equation}
for all $m \in \Z^d$.
\end{cor}
\begin{proof}
If $m\in \L$, then $m = Mn$ for some $n\in \Z^d$.  Therefore 
$e^{2\pi i \langle M^{-T} g, m \rangle} = e^{2\pi i \langle  g, M^{-1} Mn \rangle}
= e^{2\pi i \langle  g, n \rangle}=1$, giving us 
\[
\sum_{g \in \Z^d/ M\Z^d}   e^{2\pi i \langle M^{-T} g, m \rangle} 
=\sum_{g \in \Z^d/ M\Z^d}  1
=\left | \Z^d/ M\Z^d \right | = |\det M|,
\]
which proves the first part. 
On the other hand, if $m\notin \L$, then $m\rightarrow e^{2\pi i \langle M^{-T} g, m \rangle}$ is a nontrivial homomorphism of the finite group $G:= \Z^d/ M\Z^d$, so that we have the required vanishing by part 
 \ref{part 1 of lemma: vanishing group character sum} of
  Lemma \ref{lem: orthogonality of characters of a lattice}.
\end{proof}
It's also very useful to think of Corollary \ref{Orthogonality relations for characters of a lattice} in the following way.  The right-hand side of \eqref{indicator function of a lattice, in terms of characters, I}
is by definition $1_\L$, the indicator function of the lattice.  So we have the alternate form:
\begin{equation}\label{indicator function of a lattice in terms of characters, II}
1_\L(m)
=\frac{1}{|\det M|} \sum_{g \in \Z^d/ M\Z^d}   e^{2\pi i \langle M^{-T} g, m \rangle},
\end{equation}
for all $m \in \Z^d$.

\begin{thm}
\label{character identity for the lattice sum of a function}
Let $\L\subset \Z^d$ be a full-rank integer sublattice, so we may write $\L:= M(\Z^d)$, with $M$ an invertible integer matrix. 
Given an absolutely summable function $f: \Z^d \rightarrow \C$, we have:
\begin{equation}
\sum_{n\in \L} f(n) = \frac{1}{\det M}
\sum_{g \in \Z^d/ M\Z^d} \ 
 \sum_{n \in \Z^d} f(n) e^{2\pi i \langle M^{-T} g, n \rangle}.
\end{equation}
\end{thm}
\begin{proof}
\begin{align*}
\sum_{n\in \L} f(n) = \sum_{n\in \Z^d} 1_\L(n) f(n) 
&=\frac{1}{|\det M|} \sum_{n\in \Z^d}  
\ \sum_{g \in \Z^d/ M\Z^d}   e^{2\pi i \langle M^{-T} g, n \rangle} f(n)  \\
&=\frac{1}{|\det M|}  \sum_{g \in \Z^d/ M\Z^d}  
\ \sum_{n\in \Z^d} e^{2\pi i \langle M^{-T} g, n \rangle} f(n),
\end{align*}
where we used \eqref{indicator function of a lattice in terms of characters, II} in the
penultimate equality.
\end{proof}
\begin{example}
{\rm
Suppose we consider the arithmetic progression \[
\L:= \{ n\in \Z \mid n \equiv 0 \pmod 7\} := 7\Z, 
\]
which is of course a $1$-dimensional integer sublattice of $\Z$.  Here the finite group is 
$G = \Z/7\Z$, so that $|\det M | = 7$.  Here we see that for any function 
$f:\Z \rightarrow \C$ such that $\sum_{n\in \Z}|f(n)| < \infty$, 
Theorem \ref{character identity for the lattice sum of a function} gives us 
\begin{equation*}
\sum_{n \equiv 0 \pmod 7} f(n) = 
\frac{1}{7}
\sum_{g \in \Z/ 7\Z} \ 
 \sum_{n \in \Z} f(n) e^{\frac{2\pi i  gn}{7}}.
\end{equation*}
}
\hfill $\square$
\end{example}

\bigskip

\section*{Notes}
\begin{enumerate}[(a)]
\item The important families of lattices $A_n, B_n,  C_n, D_n\subset \R^n$ are 
called root lattices  (in all dimensions $n \geq 1$).  These lattice arise naturally in the classification of Lie Algebras, the combinatorics of Weyl chambers, and representation theory.  We've only glimpsed $A_n$ and $D_n$ in this chapter.
The curious reader may consult  Conway and Sloane's book \cite{ConwaySloan.book} for a lot more detail, which also gives more information about the $5$ sporadic lattices 
$E_6, E_7, E_8, F_4, G_2$, as well as many properties of all the root lattices.   Here, the
 index always signifies the dimension of the lattice.

\item The special lattice $D_4$ is currently thought to be the correct candidate for the densest sphere packings in dimension $4$ (see Chapter \ref{Sphere packings}). 

\item A lattice $\L\subset \R^d$ is called an {\bf even} lattice if $\langle x, x \rangle \in 2\Z$, 
for all $x \in \L$, and $\langle x, y \rangle \in \Z$, for all $x, y \in \L$. 
 It is a fact that the special lattice $E_8$ is the only even, unimodular lattice in $\R^8$.
A slightly deeper fact is that the only dimensions $d$ for which there exists an even, unimodular lattice are $d \equiv 0 \text{ mod } 8$.  This fact is closely tied to the theta function of such a lattice  
 \cite{ConwaySloan.book}.

\item The classic book of Martinet \cite{Martinet} develops many algebraic connections between lattices,  
 semi-simple algebras, root systems, quaternions, and quadratic forms.

\item The discriminant group $\L^*/\L$ arises naturally in the classification of lattices, but it also arises naturally in many different fields.  For example,  in the theory of chip-firing (\cite{CarolineKlivans}, Theorem 4.6.6),  it is shown that for a finite graph
$G$ the discriminant groups of the cut and flow lattices of $G$ are isomorphic. 

\item Theorem \ref{sublattice index} is usually proved using the Hermite-normal form of an integer matrix.  Here we chose this geometric route party because of its intrinsic beauty, and partly because its philosophy matches the discrete geometric path of this book.

\item We'll see an interesting application of dual lattices in Section \ref{The covering radius of a lattice, and its packing radius}, where we can `transfer the complexity' of computing a natural covering parameter of a lattice to computing 
a natural packing parameter of its dual lattice.

\end{enumerate}

\section*{Exercises}
\addcontentsline{toc}{section}{Exercises}
\markright{Exercises}

\begin{quote}    
``The only way to learn mathematics is to do mathematics"

--  Paul Halmos
 \end{quote}
\medskip

\medskip
\begin{prob}   \label{Dual of the integer lattice}
We recall that a lattice $\L$ is called self dual if $\L^* = \L$.   Prove that for any lattice $\L \subset \R^d$, we have 
$(\L^*)^* = \L$.
\end{prob}

\medskip
\begin{prob} \label{E_8 is a lattice}
Show that $E_8$, defined in Example \ref{def:E_8 lattice}, is in fact a lattice.
\end{prob}

\medskip
\begin{prob}\label{E_8 is self-dual}
Show that $E_8$ is self-dual: $E_8^* = E_8$. 
\end{prob}

\medskip
\begin{prob} $\clubsuit$ 
 \label{distance between hyperplanes}
Show that the distance $\delta$ between any two parallel hyperplanes,  described by 
 $c_1 x_1 + \cdots + c_d x_d =  k_1$ and 
$c_1 x_1 + \cdots + c_d x_d =  k_2$,  is equal to:
\begin{equation*}
\delta  = \frac{ |k_1 - k_2|}{\sqrt{ c_1^2 + \cdots + c_d^2}}.
\end{equation*}
\end{prob}

\medskip
\begin{prob}   
Suppose we are given a full-rank  sublattice of the integer lattice: $\L \subset \Z^d$.   
Prove that there is point of $\L$ on the $x$-axis.
\end{prob}

\medskip
\begin{prob}     \label{lattices in R^1}  $\clubsuit$
Let $\L$ be a lattice in $\R^1$.  Show that $\L = r\Z$ for some real number $r$.
\end{prob}

\medskip
\begin{prob} \label{matrix form for det of dual lattice}
Suppose we are given a rank $k$ lattice $\L\subset \R^d$, with  $1\leq k \leq d$. 
If $M$ is a basis matrix for $\L$, then prove that  the matrix $ M(M^TM)^{-1}$ gives a basis for the dual lattice $\L^*$.
\end{prob}

\medskip
\begin{prob} 
Show that for any two lattices $L, M\subset \R^d$,  we have $L\subseteq M \iff M^* \subseteq L^*$.
\end{prob}

\medskip
\begin{prob} \label{dual of D_n}
Prove that we have the following description for the dual lattice of $D_n$:
\[
D_n^* = \Z^d \cup   \left(  \Z^d  +  \left(  \tfrac{1}{2}, \cdots, \tfrac{1}{2} \right)^T \right).
\]
\end{prob}

\medskip
\begin{prob}     \label{Eisenstein lattice} 
The {\bf hexagonal lattice} is the $2$-dimensional lattice defined by
 \[
 \L := \{ m + n \tau \mid  m,n \in \Z\}, \text{ where }  \tau:= e^{2\pi i/3}.
 \]   
 Prove that $\det \L = \frac{\sqrt 3}{2}$, and give a description of the dual lattice to the hexagonal lattice. 
 \end{prob}


\medskip
\begin{prob} [hard]   \label{minimal lattice in R^2} 
Show that the  hexagonal lattice attains the minimal value for Hermite's constant
in $\R^2$, namely $\gamma_2^2 = \frac{2}{\sqrt{3}}$. 
 \end{prob}

\medskip
\begin{prob}  \label{special basis in R^2} 
Let $\L \subset \R^2$ be any rank $2$ lattice.  Show that there exists a basis $\beta:= \{ v, w \}$ of $\L$ such that
the angle $\theta_\beta$   between $v$ and $w$  satisfies
\[
\frac{\pi}{3}  \leq \theta_\beta \leq \frac{\pi}{2}.
\]
 \end{prob}

\medskip
\begin{prob}  \label{Hadamard's inequality, exercise}
Suppose that $M$ is a $d\times d$ matrix, all of whose $d^2$ elements are bounded by $B$. 
Show that $|\det M| \leq B_d d^{\frac{d}{2}} $.
 \end{prob} 

(Hint:  consider Hadamard's inequality \ref{Hadamard inequality})

Notes.  It follows from this exercise that if all of the elements of $M$ are $\pm 1$, then
$|\det M| \leq  d^{\frac{d}{2}} $.  If it's further true that all of the rows of $M$ are pairwise orthogonal, then $M$ is called a Hadamard matrix.  So we see from this exercise that $M$ is a Hadamard matrix $\iff$  $|\det M| =  d^{\frac{d}{2}}$.
Hadamard matrices are important in combinatorics.
It is known that if $d > 2$, then Hadamard matrices can only possibly exist when $4 \mid d$.  But for each $d = 4m$, it is not known whether a 
$d\times d$  Hadamard matrix exists, except for very small cases.

\medskip
\begin{prob} $\clubsuit$  \label{basis for A_d}
Show that the following set of vectors is a basis for $A_d$:  
\[
 \left\{e_2 - e_1, \ e_3 - e_1, \cdots , \ e_d - e_1  \right\},
\]
 where the $e_j$ are the standard basis vectors.  Hence $A_d$ is a rank-$(d-1)$
  sublattice of $\Z^d$. 
\end{prob}

\medskip
\begin{prob}
$\clubsuit$
\label{proof of equivalence of a self-dual lattice}
Prove Lemma \ref{equivalence of a self-dual lattice}, namely that
the following are equivalent:
\begin{enumerate}[(a)]
\item $\L$ is self-dual.
\item $\L$ is an integral unimodular lattice.
\end{enumerate}
\end{prob}

\medskip
\begin{prob}  $\clubsuit$      \label{Orthogonality.For.Characters.Of.A.sublattice}
\index{orthogonality relations for lattices} 
Here we prove the {\bf orthogonality relations for characters of a lattice $\L$}.  We will do it for any sublattice $\L \subset \Z^d$.  Let $D$ be a fundamental parallelepiped for  $\L$.
Using the notation in Exercise \ref{character group}, prove that for any two characters $\chi_a, \chi_b \in G_\L$,
we have:
\begin{equation}
\frac{1}{\det \L} \sum_{n \in D \cap \Z^d } \chi_a(n) \overline{\chi_b(n)} = 
\begin{cases}  
1    &      \mbox{if } \chi_a = \chi_b \\ 
0  &        \mbox{if not}.
\end{cases}
\end{equation}
\end{prob}


\medskip
\begin{prob}  $\clubsuit$  \label{fundamental domains}
Prove that any two fundamental parallelepipeds (as defined in the text) of $\L$, say $D_1$ and $D_2$, must be related to each other by an element of the unimodular group:   
\[
D_1 = M(D_2), 
\]
for some $M \in GL_d(\Z)$.
\end{prob}

\medskip
\begin{prob}    \label{number of integer sublattices of index n, R^2}
Let $f(n)$ be the number of distinct integer sublattices of index $n$ in $\Z^2$. 
We recall from elementary number theory the function $\sigma(n) := \sum_{d | n} d$, the sum of the divisors of $n$ (including $d=n$ itself). 
Show that
\[
f(n) = \sigma(n).
\]
\end{prob}

\medskip
\begin{prob}  $\clubsuit$  \label{equivalence between determinants of a sublattice}
Given a sublattice $\L\subset \R^d$ of rank $r$, show that our definition of its determinant, 
namely $\det \L :=  \sqrt{M^T M}$, conincides with the Lebesgue measure of any of its fundamental parallelepipeds.  

(Here $M$ is a $d\times r$ matrix whose columns are basis vectors of $\L$)
\end{prob}

\medskip
\begin{prob} 
Show that a set of vectors $v_1, \dots, v_m \in \R^d$, where $1\leq m \leq d$, are linearly independent
$\iff$  their Gram matrix is nonsingular.  
\end{prob}

\medskip
\begin{prob}
Prove that for any given lattice $\L \subset \R^2$, any two(nonzero) shortest  linearly independent 
vectors for $\L$ generate the lattice  $\L$.

Notes.  \ As a reminder, the first two shortest nonzero vectors of $\L$ may have equal length. 
We note that in dimensions $d \geq 5$, such a claim is false in general, as problem \ref{counterexample in RË5} below shows. 
\end{prob}

\medskip
\begin{prob} \label{counterexample in RË5}
Find a lattice $\L \subset \R^5$ such that any set of five shortest nonzero vectors of $\L$ do not generate $\L$. 
\end{prob}

\medskip
\begin{prob}\label{Hermite normal form: basis for 2-dim'l lattices}
Given any $2$-dimensional lattice $\L:= M(\Z^2)$, use the Hermite-normal form of $M$ 
to prove that $ \{ \icol{1\\0},  \icol{p\\q} \}$ is a basis for $\L$, for some 
nonnegative integers $p, q$.
\end{prob}

\medskip
\begin{prob}  $\clubsuit$  \label{hyperplane lattice}
 Consider the {\bf discrete hyperplane} defined by:
\[
H:= \left\{  x \in \Z^d  \bigm |   c_1 x_1 + \cdots + c_d x_d =0   \right\},
\]
Show that $H$ is a sublattice of $\Z^d$, and has rank $d-1$.
\end{prob}

\medskip
\begin{prob}  $\clubsuit$  \label{tiling the integer lattice with hyperplanes}
\index{discrete hyperplane}
\rm{
Suppose we are given a discrete hyperplane $H$, as in Exercise \ref{hyperplane lattice}.
 \begin{enumerate}[(a)]
 \item Prove there exists a vector $x\in \R^d$ such that 
\[
\{ H+ kx \bigm |  k \in \Z \} = \Z^d.
\]
\item  Prove that there are no integer points strictly between $H$ and $H + x$.
\end{enumerate}
}
\end{prob}

Notes.     
You may assume Bezout's identity,  \index{Bezout's identity} which states that
 if $\gcd(c_1, \dots, c_d)=1$
then there exists 
an integer vector $(m_1, \dots, m_d)$ such that $c_1 m_1 + \cdots + c_d m_d = 1$.
This exercise shows that we can always tile the integer lattice $\Z^d$ with discrete
translates of a discrete hyperplane $H$.

\medskip
\begin{prob}\label{Ellipsoid problem}
\rm{
Here we give the details for \eqref{ellipsoid}, the definition of an ellipsoid in $\R^d$.  Starting over again,
we fix an orthonormal basis $\{ b_1, \dots, b_d\}$ for $\R^d$, and 
we define the following matrix:
\[
M := \begin{pmatrix} |  &  |  &  ...   & |  \\  
                        c_1 b_1  & c_2 b_2 &  ...& c_d b_d   \\  
                          |  &  |  & ... &  |  \\ 
        \end{pmatrix},
\]
where the $c_k$'s are positive scalars. 
We  now apply the linear transformation $M$ to the unit sphere 
$S^{d-1}:= \{ x \in \R^d \mid \| x \|^2  = 1\}$ in 
$\R^d$, and we recall what this entails. We define the $\text{Ellipsoid}_M:=  
M(S^{d-1})$, a $(d-1)$-dimensional object.  In the spirit of review,  we recall that by
 definition
$M(S^{d-1}) :=  \{ u \in \R^d \mid u = Mx,  x \in S^{d-1}  \}$.   
\begin{enumerate}[(a)]
\item  Show that 
\begin{equation} \label{equation of ellipsoid}
\text{Ellipsoid}_M =
\left\{    x\in \R^d  \bigm |       \sum_{j=1}^d \frac{{\langle x, b_j\rangle}^2}{c_j^2} =1  \right\}.
\end{equation}
\item  We recall that  the unit ball in $\R^d$ is defined by 
$B := \left\{ x \in \R^d \bigm |   \| x \|^2  \leq 1\right\}$.
Show that for the open ellipsoid body $E$ (a $d$-dimensional object), as defined in \eqref{open ellipsoid}, we have the
$d$-dimensional volume formula:
\[
\vol(E) = \vol B \prod_{j=1}^d c_j.
\]
\end{enumerate}
}
\end{prob}

\medskip
\begin{prob} 
\rm{
We will use the equation \eqref{equation of ellipsoid} definition of an ellipsoid.
We can extend the previous exercise in the following way.  Let $A$ be  {\bf any} $d \times d$ real matrix,
and look at the action of $A$ on the unit sphere
 $S^{d-1} \subset \R^d$. Show:

(a)  If  $\rm{rank}(A) = d$, then $A(S^{d-1})$ is a $d$-dimensional ellipsoid, defined by an equation of the form
\eqref{equation of ellipsoid}.

(b)  If $\rm{rank}(A) := r < d$, then $A(S^{d-1})$ is an $r$-dimensional ellipsoid.
}
\end{prob}

\medskip
\begin{prob}   \label{square root of a matrix}
\rm{
Suppose that $A$ is a positive definite, real matrix. 
Solve for (i.e. characterize) all matrices $X$ that are the `square roots' of $A$:
\[
A = X^2.
\] 
}
\end{prob}

\medskip
\begin{prob} 
\rm{
Suppose that a certain $2$-dimensional lattice $\L$ has a Gram matrix 
\[
G := \begin{pmatrix}
\ 2 & -1 \\
 -1 &  \  2 \end{pmatrix} .  
\]
Reconstruct $\L$ (i.e. find a basis for $\L$), up to an orthogonal transformation.
}
\end{prob}

\medskip
\begin{prob}  
\rm{
Find a $2$ by $2$ matrix $M$ that enjoys one of the properties of a positive semidefinite
matrix, namely that  $x^T M x \geq 0$, for all $x\in \R^2$, but such that $M$ is not symmetric.  
}
\end{prob}

\medskip
\begin{prob}  
\rm{
Given any $3$ lattices with $L_1 \subseteq L_2 \subseteq L_3$, show that their indices are multiplicative in the following sense:
\[
\left | L_3/L_1\right | =   \left | L_3/L_2\right |  \left | L_2/L_1\right |.
\] 
}
\end{prob}

\medskip
\begin{prob}  
\rm{
To count the number of sublattices of a fixed index, let's 
define $N_d(k)$ to be the number of integer sublattices of $\Z^d$ that have a fixed index $k$, for any given positive integer $k$.
\begin{enumerate}[(a)]
\item Prove that $N_d(2) = 2^d -1$.
\item Can you find a formula for $N_d(k)$, at least in the case that $k$ is prime?
\end{enumerate} 
}
\end{prob}

Notes.  \ Here it may be useful to think about the Hermite-normal form.  See \cite{Zong.lattices} for a recent study of $N_d(k)$, and also of the number of sublattices of fixed index that are equivalent under the unimodular group $GL_d(\Z)$.

\medskip
\begin{prob}  \label{exercise:2by2 positive definite matrix}
\rm{
Suppose we are given a real $2$ by $2$ matrix $A$.   Prove that
\[
A \text{ is positive definite } \iff  \text{ both } \rm{trace}(A) >0 \text{ and } \det A > 0.
\] 
}
\end{prob}

\medskip
\begin{prob} 
\label{Ex:smith normal form}
\rm{(not trivial) 
Here we prove the  existence of the (geometric) {\bf Smith normal form} of a lattice.
\index{Smith normal form}
Namely, let $\L_0\subset \L$ be a 
sublattice of $\L$.   Then  there exists a basis $\{ v_1, \dots, v_d \}$ of $\L$, 
 and positive integers $k_1,... ,k_d$ such that:
 \begin{enumerate}[(a)]
 \item $\{ k_1 v_1, \cdots, k_d v_d \}$ is a basis for $\L_0$.
 \item  $k_j \mid k_{j+1}$ for $j=1, \dots, d-1$. 
 \end{enumerate}
}
\end{prob}

\medskip
\begin{prob} 
\label{discriminant group in terms of Smith Normal Form}
\rm{(assuming some background about finite abelian groups) 
Here we'll assume the existence (and notation) of the {\bf Smith normal form} from Exercise \ref{Ex:smith normal form}.  Prove  that for a full-rank lattice $\L\subset \R^d$ with $\L \subset \L^*$, we have the following explicit form for its discriminant group:
\[
\L^*/\L  =   \Z/k_1 \Z \times \cdots \times \Z/k_d \Z.
\]
}
\end{prob}

\medskip
\begin{prob}  
\label{structure of Discriminant group of D_n}
\rm{
Prove that for the lattice $D_n$, its discriminant group $D_n^*/ D_n$ has the following structure:
\begin{enumerate}[(a)]
\item $D_n = \Z/2\Z \times \Z / 2\Z \iff $ $n$ is even.
\item $D_n = \Z/4\Z  \iff $ $n$ is odd.
\end{enumerate} 
}
\end{prob}

\medskip
\begin{prob}  
\label{permuting the elements of a group}
\rm{
For the readers who may not be familiar with groups in general, see \cite{Herstein} for the definition of a group. 
Given any group $G$ (not necessarily finite), and any element $g \in G$, prove that
\[
gG = G.
\]
}
\end{prob}


\chapter{
\blue{
Classical geometry of numbers \\ 
Part II: \, Blichfeldt's theorems}    
}  
\label{Chapter.geometry of numbers II}
\index{Blichfeldt}

\begin{quote}   
``Simplicity is the ultimate sophistication.''                  

-- Leonardo Da Vinci
\end{quote}

\begin{wrapfigure}{R}{0.47\textwidth}
\centering
\includegraphics[width=0.23\textwidth]{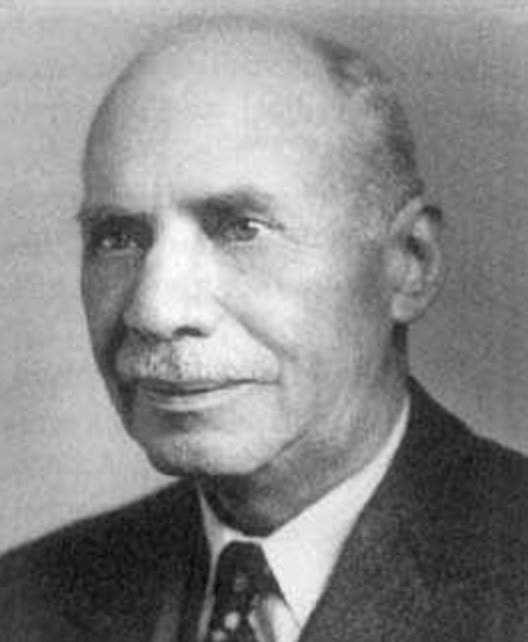}
\caption{Hans Blichfeldt}
\end{wrapfigure}

\bigskip
\section{Intuition}

\begin{figure}[htb]
 \begin{center}
\includegraphics[totalheight=3.3in]{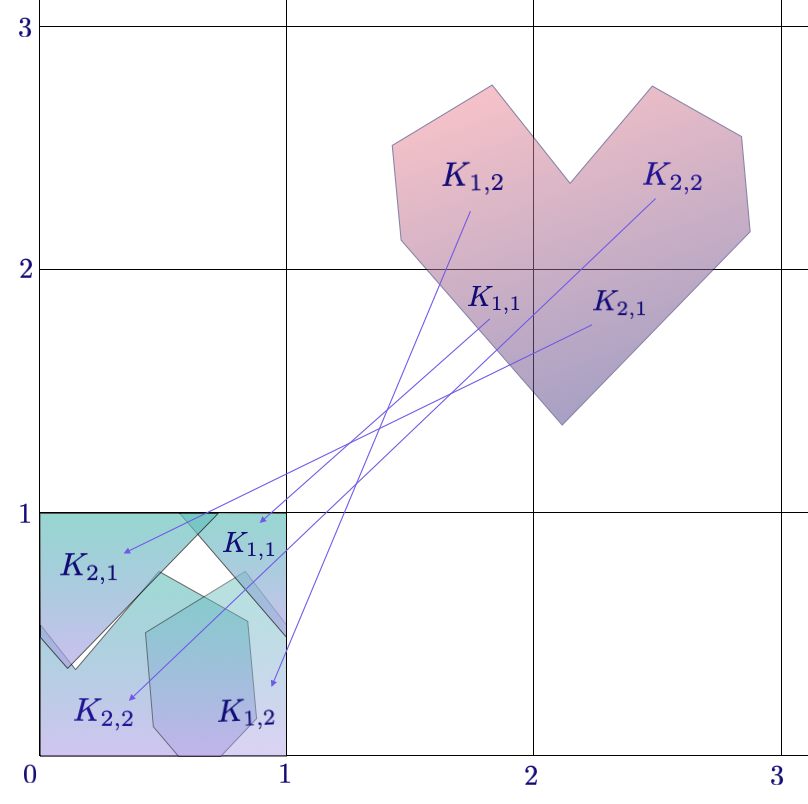}
\end{center}
\caption{An outline of Blichfeldt's elementary argument, using a `broken heart'.  This result has enjoyed a lot of applications.}  
\label{Blichfeldt2}
\end{figure}

There is a beautifully simple and powerful idea, in the geometry of numbers, due to Hans Blichfeldt, who discovered it in 1914.  Here is a simple illustration of it - suppose we have a body $K\subset \R^2$, whose area is bigger than $1$.   Now, obviously $K$ intersects each little integer square $[m, m+1]\times[n, n+1]$, in some little region $K_{m,n}$, as in Figure \ref{Blichfeldt2}.  After translating all of these little regions to the unit square $[0,1]^2$, it must be the case that there exists 
a point $p$ in the interior of 
$[0,1]^2$ that is covered by at least $2$ integer translates of the little regions $K_{m,n}$.  Thinking through it over a fresh cup of coffee, we conclude that there
are (at least) two points $x, y \in K$ that enjoy the property  $x-y\in \Z^2$, as the overlapping regions in the unit square of Figure \ref{Blichfeldt2} suggest.  This elegant conclusion is part of a stronger phenomenon, which was thoroughly developed by Blichfeldt, and which we now study.
Corollary \ref{Blichfeldt's lemma} below gives another proof of this same fact, but with more generality, including any full-rank lattice in dimension $d$.

\bigskip
\section{Blichfeldt's Theorem}
\label{section:Blichfeldt}
In this section we study a classical and powerful result of Blichfeldt, which will also give an alternate proof of Minkowski's first theorem.  We begin with some elementary functional analysis (Theorem \ref{thm:Blichfeldt1}), and then we develop some straightforward number-theoretic  consequences, such a classical {\bf pigeon-hole} geometric principle (Corollary \ref{Blichfeldt's lemma}), both of which have found many applications. 

Robert Remak \cite{Remak1} further extended Blichfeldt's work (in a beautiful and elementary way), and 
we feel that it's very useful to follow his more general route. Namely, there is a dance between counting and computing volumes, and more generally proceeding from the function-analytic approach to the derivation of very concrete combinatorial consequences.  
\begin{thm}[Remak, 1927]\index{Remak, Robert}
\label{thm:Blichfeldt1}
Let $f \in L^1( \R^d )$ be a nonnegative function, and let $\L \subset \R^d$ be
a full-rank lattice.  
\begin{enumerate}[(a)]
\item \label{part a of Blichfeldt thm 1}
There exists a point $y \in \R^d$ such that 
\begin{equation}
(\det \L) \sum_{n \in \L} f(y+n) \geq  \int_{\R^d} f(x) dx.
\end{equation}
\item \label{part b of Blichfeldt thm 1}
On the other hand, there also exists a point $z \in \R^d$ such that
\begin{equation}
(\det \L) \sum_{n \in \L} f(z+n) \leq  \int_{\R^d} f(x) dx.
\end{equation}
\end{enumerate}
\end{thm}
\begin{proof}
We fix a basis $\{ v_1, \dots, v_d\}$ for the lattice $\L$, and we consider its fundamental parallelepiped $\Pi:=
\left\{ \lambda_1 v_1 + \cdots + \lambda_m v_m  \bigm |     \text{  all }   0 \leq  \lambda_k < 1   \right\}$. 
By Lemma \ref{tiling by translates of Pi}, each $x \in \R^d$ can be written uniquely as $x = v + n$, with $v \in \Pi, n \in \L$.  We therefore have:
\begin{equation}\label{unfolding for Blichfeldt}
\int_{\R^d} f(x) dx = \sum_{n \in \L} \int_{\Pi} f(v+n) dv 
=\int_{\Pi} \left(   \sum_{n \in \L}  f(v+n)  \right) dv.
\end{equation}
We may of course assume that $f$ is not the zero function.
Because $f$ is nonnegative by the hypothesis of the theorem, we have $\int_{\R^d} f(x) dx>0$, so there exists a positive constant $c>0$ such that 
\begin{equation} \label{assumption for Blichfeldt}
\int_{\R^d} f(x) dx = c \det \L.  
\end{equation}
To prove part \ref{part a of Blichfeldt thm 1}, suppose it was true that $\sum_{n \in \L}  f(v+n) < c$  for all $v \in \Pi$.
Then  using \eqref{unfolding for Blichfeldt} we would obtain:
\[
\int_{\R^d} f(x) dx =\int_{\Pi} \left(   \sum_{n \in \L}  f(v+n)  \right) dv < c \int_{\Pi} dv
= c \det \L,
\]
contradicting \eqref{assumption for Blichfeldt}.  Therefore there exists at least one point
 $y \in \Pi$ such that 
 \[
 \sum_{n \in \L}  f(y+n) \geq c := \frac{1}{\det \L} \int_{\R^d} f(x) dx.
 \]
 To prove part \ref{part b of Blichfeldt thm 1}, suppose it was true that $\sum_{n \in \L}  f(v+n) > c$  for all $v \in \Pi$.  Then using \eqref{unfolding for Blichfeldt}, we would obtain:
\[
\int_{\R^d} f(x) dx =\int_{\Pi} \left(   \sum_{n \in \L}  f(v+n)  \right) dv > c \int_{\Pi} dv
= c \det \L,
\]
again contradicting \eqref{assumption for Blichfeldt}.  Therefore there exists at least one point
 $z \in \Pi$ such that 
 \[
 \sum_{n \in \L}  f(z+n) \leq c := \frac{1}{\det \L} \int_{\R^d} f(x) dx.
 \]

\end{proof}
Next, we apply Theorem \ref{thm:Blichfeldt1} to the indicator function of any 
set $S\subset \R^d$, arriving at the classical and useful result, known as ``Blichfeldt's lemma'' \cite{Blichfeldt2}.

\begin{cor}[Blichfeldt's lemma, 1914] \label{Blichfeldt's lemma}
Let $\L\subset \R^d$ be a full-rank lattice, and let $S$ be any subset of $\R^d$, 
whose volume $\vol(S)$ is also allowed to be $\infty$.   We fix any positive integer $m$.
If we have
  \[
  \vol(S) > m \det \L,
  \]
   then there exist $m+1$ distinct points $p_1, \dots, p_{m+1}\in S$ such that their pairwise differences $p_i - p_j$ are all lattice points of $\L$.
\end{cor}
\begin{proof}
Using the function $1_S$ in Remak's Theorem \ref{thm:Blichfeldt1}, we know there exists 
some $y \in \R^d$ such that 
\begin{equation}
(\det \L) \sum_{n \in \L} 1_S(y+n) \geq  \int_{\R^d} 1_S(x) dx = \vol(S) > m \det \L,
\end{equation}
where the second inequality above just follows by assumption.  So we arrive at
\begin{equation}\label{penultimate inequality of Cor of Blichfeldt}
 \sum_{n \in \L} 1_S(y+n) >  m,
\end{equation}
which implies that $\sum_{n \in \L} 1_S(y+n) \geq   m+1$, since the left-hand-side of 
\eqref{penultimate inequality of Cor of Blichfeldt} is an integer.   But the latter inequality
 means that there are at least $m+1$ distinct points $n_k\in \L$ such that $y+n_j \in S$, which is the desired conclusion (with $p_j:= y+n_j\in S$).
\end{proof}
Even the case $m=1$ of Corollary \ref{Blichfeldt's lemma} is very useful, 
and we record it separately, as it is one of the best known results in the geometry of numbers.
\begin{cor}[Case $m=1$ of Blichfeldt's lemma] 
\label{m=1 Blichfeldt's lemma} 
Suppose we are given a full-rank lattice 
 $\L \subset \R^d$, and any  subset $S \subset \R^d$  such that 
 $\vol(S) > \det \L$.   Then there exist (at least two) distinct points $a, b \in S$  such that $a-b \in \L$.
\end{cor}
Now we can give another proof of Minkowski's first theorem 
(Theorem \ref{Minkowski convex body Theorem for L}), as a quick consequence of Corollary \ref{m=1 Blichfeldt's lemma}.
\begin{proof}\rm{ [Second proof of Minkowski's first theorem]}  \label{second proof of Minkowski}
 \hypertarget{second proof of Minkowski} We define $K:= \frac{1}{2}B$, so we have:
 \[
 \vol K = \frac{1}{2^d}\vol B > \det \L,
 \]
 the latter inequality holding by the assumption of 
 Theorem \ref{Minkowski convex body Theorem for L}.  Since $\vol K > \det \L$, 
 Corollary \ref{m=1 Blichfeldt's lemma} tells us that there exist \emph{distinct}
  points $a, b \in K$ such that $a-b \in \L$.  If we show that $n:= a-b$ is also in $B$, we're done. 
To this end, we notice that  $2a \in B$ and $2b \in B$.  Since B is centrally symmetric, we also 
have $-2b \in B$, so that
\[
n = \tfrac{1}{2}(2a) + \tfrac{1}{2}(-2b) \in \tfrac{1}{2}B + \tfrac{1}{2}B = B,
\]
where the latter equality holds because $B$ is convex.
\end{proof}

There is another useful consequence of Remak's Theorem \ref{thm:Blichfeldt1}, originally due to Van der Corput.
\begin{cor}[Van der Corput] 
\label{cor:upper and lower bounds using Remak's thm}
Suppose we are given a full-rank lattice 
 $\L \subset \R^d$, and any  subset $S \subset \R^d$.  Then there exists vectors $v_1, v_2 \in \R^d$ 
such that 
\begin{equation} 
\label{lattice point enumerator bounds volume above and below}
\left | \L \cap (S+v_2) \right| \leq \frac{\vol S}{\det \L}  \leq \left | \L \cap (S+v_1) \right|.
\end{equation}
\end{cor}
\begin{proof}
To prove the right-hand inequality in 
\eqref{lattice point enumerator bounds volume above and below}, we apply Remak's Theorem \ref{thm:Blichfeldt1}, part \ref{part a of Blichfeldt thm 1}, to the indicator function $1_S$, which gives us the existence of a $v_1\in \R^d$ such that
\begin{equation}
 (\det \L) \sum_{n \in \L} 1_S(-v_1+n) \geq  \int_{\R^d} 1_S(x) dx 
:= \vol S.
\end{equation}
But $ (\det \L) \sum_{n \in \L} 1_S(-v_1+n) = (\det \L) \sum_{n \in \L} 1_{S+v_1}(n) = 
(\det \L)  \left | \L \cap (S+v_1) \right|$, and we're done.  The left-hand side of inequality 
\eqref{lattice point enumerator bounds volume above and below} is proved in exactly the same manner, this time applying Remak's Theorem \ref{thm:Blichfeldt1}, part \ref{part b of Blichfeldt thm 1}, to the indicator function $1_S$.
\end{proof}

We note that in practice, the upper bound in 
Corollary \ref{cor:upper and lower bounds using Remak's thm} has found more applications.


\bigskip
\section{Van der Corput's inequality for convex bodies}


There is a natural extension of Minkowski's first theorem 
(Theorem \ref{Minkowski convex body Theorem for L}) to convex, centrally symmetric sets
that contain any number of integer points, known as Van der Corput's inequality \cite{Van.der.Corput.1}.

\begin{cor}[Van der Corput's inequality, 1936] \label{Van der Corput}
\index{Van der Corput}
Let $\L\subset \R^d$ be a full-rank lattice, and let $K$ be a centrally symmetric, convex $d$-dimensional
set in $\R^d$, 
whose volume is also allowed to be $\infty$.   We fix any positive integer $m$.
\begin{enumerate}[(a)]
\item If $ \vol(K) > m 2^d \det \L$, then  $\left | K\cap \L \right | \geq 2m+1$.  \label{first part of Van der Corput}
\item If  $ \vol(K) = m 2^d \det \L$, and we also assume that $K$ is compact, then 
$\left | K\cap \L \right | \geq 2m+1$.  \label{second part of Van der Corput}
\end{enumerate}
\end{cor}
\begin{proof}
To prove (a), let's apply Blichfeldt's lemma \ref{Blichfeldt's lemma} to the set $\tfrac{1}{2}K$, whose volume equals 
$\tfrac{1}{2^d}\vol K$.  So there exist $m+1$ distinct points 
$\tfrac{1}{2} p_1, \dots, \tfrac{1}{2} p_{m+1}\in \tfrac{1}{2}K$ with the property that all of their pairwise differences 
 $\tfrac{1}{2} p_i - \tfrac{1}{2} p_j$ are distinct, and $\tfrac{1}{2} p_i - \tfrac{1}{2} p_j \in \L$.  We define an ordering on $\R^d$ by saying that for any two points $x, y \in \R^d$, $x > y$ if the first coordinate of $x$ is larger than the first coordinate of $y$.   Without loss of generality we assume that $p_1 > p_2 > \cdots > p_{m+1}$.  
 Defining $q_k:= \tfrac{1}{2} p_k - \tfrac{1}{2} p_{m+1}$, we've already seen that $q_k \in \L$.  With the ordering defined above, we also have $q_k > q_{k+1}$, and in particular the points $0, \pm q_1, \dots, \pm q_m$ are all distinct.

So it suffices to show that the $2m+1$ distinct points $0, \pm q_1, \dots, \pm q_m$ all belong to $K$.  But
\[
q_k:= \tfrac{1}{2} p_k - \tfrac{1}{2} p_1 \in  \tfrac{1}{2}K - \tfrac{1}{2} K = K, 
\]
where we used the convexity and central symmetry of $K$ in the last equality (recalling Exercise  \ref{c.s. C equals its symmetrized body}). 
We leave part \ref{second part of Van der Corput} as Exercise \ref{Van der corput exercise 1}.
\end{proof}

\begin{example}
\rm{
We consider the long and thin box described by 
\[
K:=\{ x\in \R^d \mid  |x_k| < 1, \text{ for } k = 1, 2, \dots, d-1, \text{ and } |x_d| < m \},
\]
 for any fixed positive integer $m$. 
 Here $K$ contains precisely the integer points 
  $(0, 0, \dots, \pm k)$, for $k \in \{ 0, 1, \dots, m-1\}$.  To summarize, $K$  contains exactly $2m-1$ integer points.  Using the fact that $\vol K = m 2^d$, we now see that
  Corollary \ref{Van der Corput}, 
  part \ref{first part of Van der Corput}
  is sharp. 
 }
 \hfill $\square$
\end{example}
It's clear that the case $m=1$ of  
Corollary \ref{Van der Corput},  part \ref{first part of Van der Corput}
 is Minkowski's Theorem \ref{Minkowski convex body Theorem for L}; indeed when $m=1$, the hypothesis $ \vol(K) >  2^d \det \L$ tells us, via Minkowski's Theorem, that 
  $K$ must contain at least one nonzero lattice point $p\in \L$.  But by the central symmetry of $K$, we know that $-p\in \L$ as well, so that we have $K\cap \L \supset  \{-p, 0, p\}$.    This is of course equivalent to the conclusion of Corollary \ref{Van der Corput}, 
  part \ref{first part of Van der Corput} for $m=1$.

Sometimes it's useful to state Van der Corput's inequality 
(Corollary \ref{Van der Corput}, part \ref{first part of Van der Corput})
 in its contrapositive form, using the (trivial) fact that the number of interior lattice points in a centrally-symmetric body 
 $K$ is always an odd integer:
\begin{equation}\label{contrapositive to Van der Corput}
\text{If } \left | \interior K\cap \L \right | \leq 2m-1, \text{then } \vol K \leq m 2^d \det \L.
\end{equation}
Interestingly, $85$  years passed since the paper of 
Van der Corput \cite{Van.der.Corput.1}, before the equality cases of Corollary \ref{Van der Corput} were completely classified in Averkov's recent work \cite{Averkov}.

  \medskip
\begin{question}[Rhetorical] \label{Van der Corput for non cs bodies}
What about finding  a ``Van der Corput''-type inequality for bodies $K$ that are not necessarily centrally symmetric?
\end{question}

There is an ``easy-fix'' that gives us a positive answer to Question \ref{Van der Corput for non cs bodies}.
We notice that the only time we used central symmetry in the proof of Corollary \ref{Van der Corput} 
was at the very end of the proof.   So we get a more general conclusion, for the body 
$\tfrac{1}{2}K - \tfrac{1}{2}K$ (instead of $K$), with precisely the same proof of 
Corollary \ref{Van der Corput}.  

\begin{cor}\label{general Van der Corput} \index{Van der Corput's inequality}
Let $\L\subset \R^d$ be a full-rank lattice, and let $K$ be a convex $d$-dimensional
set in $\R^d$, 
whose volume is also allowed to be $\infty$.   We fix any positive integer $m$.
\begin{enumerate}[(a)]
\item If $ \vol(K) > m 2^d \det \L$, then  $\left | \left(\tfrac{1}{2}K - \tfrac{1}{2}K\right)\cap \L \right | \geq 2m+1$.  \label{first part of general Van der Corput} 
\item If  $ \vol(K) = m 2^d \det \L$, and we also assume that $K$ is compact, then 
\[
\left | \left(\tfrac{1}{2}K - \tfrac{1}{2}K\right)\cap \L \right | \geq 2m+1.  
\]  \label{second part of general Van der Corput}
\hfill $\square$
\end{enumerate}
\end{cor} 

\begin{example}
\rm{
Let's consider all convex integer polygons $K \subset \R^2$ with exactly one integer point 
in the interior of 
$\tfrac{1}{2}K - \tfrac{1}{2}K$, which we'll assume to be the origin.

To begin, it's again useful to phrase 
Corollary \ref{general Van der Corput}  part \ref{first part of general Van der Corput} , in its contrapositive form:
\begin{equation}\label{contrapositive to general Van der Corput}
\text{If } \left | \interior \left(\tfrac{1}{2}K - \tfrac{1}{2}K\right)\cap \Z^2 \right | \leq 2m-1, \text{then } 
\vol K \leq m 2^d := 4m.
\end{equation}
If  $m=1$, then by our assumption in \eqref{contrapositive to general Van der Corput}
there is precisely $1$ interior integer point belonging to the interior of $\tfrac{1}{2}K - \tfrac{1}{2}K$.  The conclusion of \eqref{contrapositive to general Van der Corput}
 is that $\vol K \leq 4$. 
 }
 \hfill $\square$
\end{example}

\bigskip

\section*{Notes} \label{Notes.chapter.Geometry of Number II}
\begin{enumerate}[(a)]
\item
We mention another result of Blichfeldt, which goes in the other direction to the previous theorems, giving us a lower bound on the volume by assuming it contains enough integer points.
\begin{thm}[Blichfeldt, 1921] \index{Blichfeldt}
Suppose that $K\subset \R^d$ is a $d$-dimensional convex body that contains at least $d$ linearly independent integer points (possibly on its boundary).  Then:
\[
\vol K \geq \frac{1}{d!} \left( \left|K\cap \Z^d \right| - d \right).
\]
\hfill $\square$
\end{thm}
Blichfeldt's latter bound is best-possible, in the sense that equality is achieved, for example,
 by the following countable collection of integer simplices in each dimension:
\begin{equation}
\Delta_k:= \conv\{ 0, k \, e_1, e_2, e_3, \dots, e_d \},
\end{equation}
defined for each positive integer $k$.  A moment's thought gives 
$\vol \Delta_k= \frac{k}{d!}$,  as well as  $| \Delta_k \cap \Z^d | = d+k$.  For more information, 
see \cite{HenkHenzWills}.
\end{enumerate}


\bigskip
\section*{Exercises}
\addcontentsline{toc}{section}{Exercises}
\markright{Exercises}

\begin{quote}
``Math is dirty, if it is done right.'' 

-- Günter Ziegler
\end{quote}

\medskip
\begin{prob}$ \clubsuit$
\rm{
Let $K\subset \R^d$ be a convex set of finite volume (but not necessarily bounded). 
\begin{enumerate}[(a)]
\item Prove that if  $\vol K > 1$, then $K$ must contain an integer point of $\Z^d$.  
\item Prove that if  $\vol K > m$, for any positive integer $m$, then $K$ must contain 
at least $m$ distinct integer points of $\Z^d$.
\end{enumerate}
}
\end{prob}

\medskip
\begin{prob} \label{one iteration of Minkowski symmetrization}
\rm{
Suppose we have a convex, compact set  $K\subset \R^d$ (but $K$ is not necessarily centrally symmetric). 
 We define $Q:= \tfrac{1}{2}K - \tfrac{1}{2}K$.  We already know the (trivial)
  fact that $Q$ is centrally symmetric. 
  
Prove that $\tfrac{1}{2}Q - \tfrac{1}{2}Q = Q$.
}
\end{prob}

\medskip
\begin{prob}
\rm{
Suppose we're given a centrallly-symmetric hexagon $H$ in the plane, with $\vol H = 8$.
Prove that $H$ contains at least $5$ integer points (some of which might lie on its boundary as well). 
}
\end{prob}

\medskip
\begin{prob} \label{elementary inequality, concavity 1}
Prove the inequality
\[
(x+y)^r \geq x^r + y^r,
\]
valid for all $r>1$ and $x, y > 0$. 
\end{prob}

\medskip
\begin{prob}  
\label{convergence of basic lattice sum}
With the usual norm $\|n\|:= \sqrt{n_1^2 + \cdots + n_d^2}$, prove that  if $r\in \R$, then
\[
\sum_{n\in \Z^d} \frac{1}{ \| n \|^r} < \infty \iff r > d.
\]
\end{prob}

\medskip
\begin{prob}\label{Arithmetic-geometric inequality application}
\rm{ 
Given positive numbers $a_1, \dots a_d$, with the property that $\prod_{k=1}^d a_k = 1$, prove that 
\[
(1+a_1)(1+a_2) \cdots (1+a_d) \geq 2^d.
\]
}
\end{prob}

\medskip
\begin{prob} \label{Van der corput exercise 1}
 $\clubsuit$
\rm{
Prove part (b) of Van der Corput's Theorem \ref{Van der Corput}. Namely, we are given a full-rank lattice 
 $\L\subset \R^d$  and a compact, convex, centrally symmetric $d$-dimensional set $K\subset \R^d$, 
together with any positive integer $m$.  Prove that if $ \vol(K) = m 2^d \det \L$, then
\[
\left | K\cap \L \right | \geq 2m+1. 
\]
}
\end{prob}

\medskip
\begin{prob} 
\rm{
Suppose $K\subset \R^d$ is a convex body.  
If there is exactly one integer point in the interior of
 $Q:= \tfrac{1}{2}K - \tfrac{1}{2}K$, must there also exist at least one integer point in the interior of $K$?
}
\end{prob}

\medskip
\begin{prob} 
\rm{
Suppose we are given $n$ bounded sets $S_1, \dots, S_n\subset \R^d$, and $n$ positive numbers $c_1, \dots, c_n$.  Prove that there exists a single vector $y\in \R^d$ such that:
\[
c_1\left| (S_1 + y) \cap \Z^d \right| +\cdots + c_1\left| (S_1 + y) \cap \Z^d \right| \geq
c_1(\vol S_1)+ \cdots + c_n(\vol S_n).
\]
}
\end{prob}
Hint. \ Apply Remak's Theorem \ref{thm:Blichfeldt1} to an appropriate linear combination of indicator functions.


\chapter{\blue{The Fourier transform of a polytope via its vertex description: \\  The Brion theorems} }
\label{chapter.Brion}

\begin{quote}
``See in {\bf nature} the cylinder, the sphere, the cone.''

-- Paul C\'ezanne 
\end{quote}

\begin{figure}[htb]
 \begin{center}
\includegraphics[totalheight=2in]{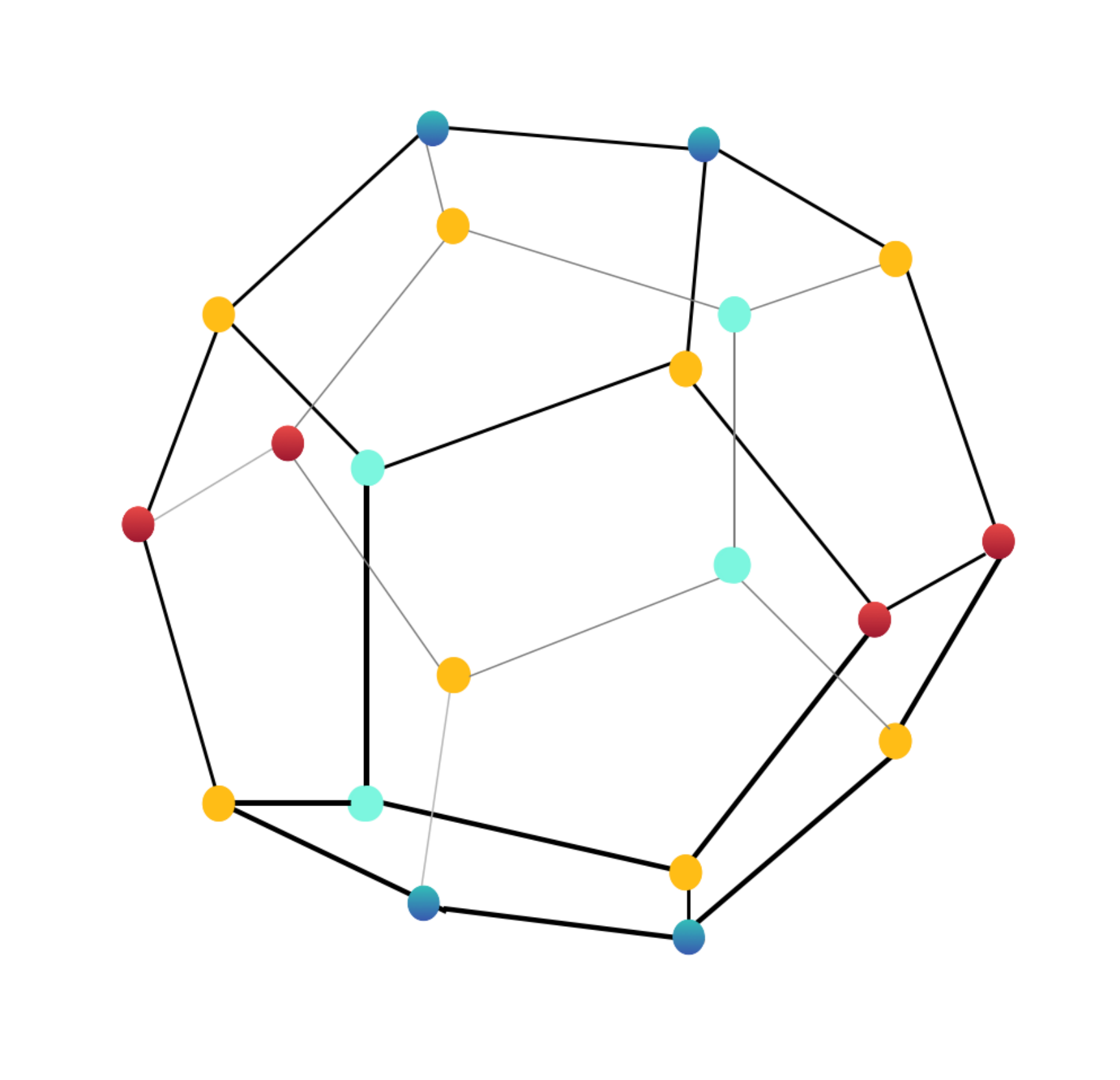}
\end{center}
\caption{The Dodecahedron in $\R^3$, an example of a simple polytope.  
In Exercise \ref{FT of a Dodecahedron},  we compute its Fourier-Laplace 
transform by using Theorem \ref{brion, continuous form} below.}  \label{Dodecahedron} 
\end{figure}

\bigskip
\section{Intuition}
Here we introduce the basic tools for computing precise expressions for the Fourier transform of a polytope.    To compute transforms here, we assume that we
 are given the vertices of a polytope $\P$ , together with the local geometric information at each vertex of $\P$, namely its neighboring vertices
  in $\P \subset \R^d$.  It turns out that computing the Fourier-Laplace transform of the tangent cone at each vertex of $\P$ completely characterizes the Fourier transform of $\P$. 

One of the basic results here, called the discrete version of Brion's Theorem (\ref{brion, discrete form}), 
may be viewed as an extension of the finite geometric sum in dimension $1$,  to sums in integer cones, in dimension $d$.  
Some basic families of polytopes are introduced, including simple polytopes and their polars, which are simplicial polytopes.  These families of polytopes play an important role in the development of Fourier analysis on polytopes.
\begin{figure}[htb]
 \begin{center}
\includegraphics[totalheight=2.6in]{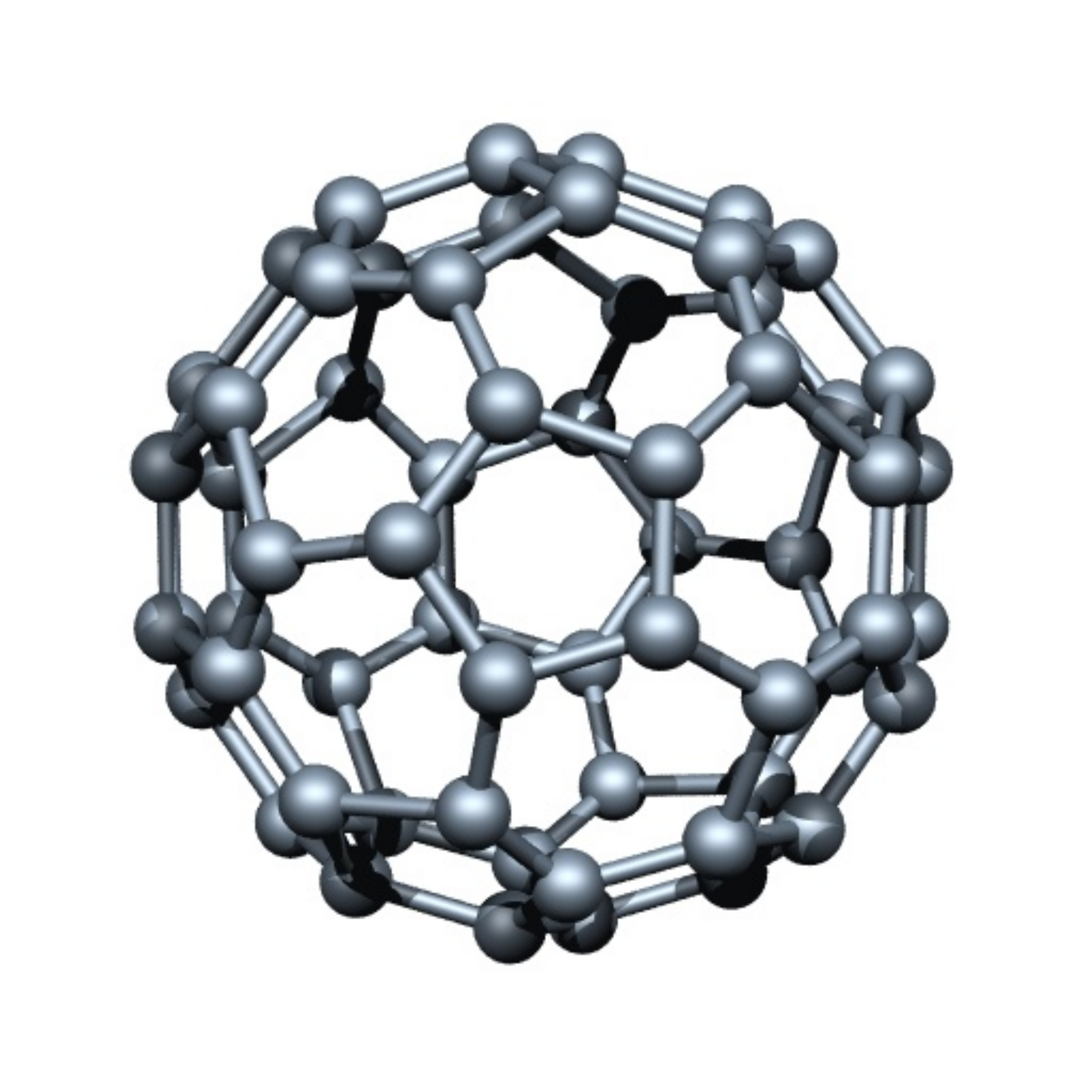}
\end{center}
\caption{The $C_{60}$ Carbon molecule, also known as a buckeyball,  is another example of a simple polytope.   The nickname ``buckeyball' came from Buckminster Fuller, who 
used this molecule as a model for many other tensegrity structures. 
\index{simple polytope}
(the graphic is used with permission from
Nanografi, at https://phys.org/news/2015-07-scientists-advance-tunable-carbon-capture-materials.html) }  \label{C60} 
\end{figure}

\bigskip
\section{Cones, simple polytopes, and simplicial polytopes}

One of the most important  concepts in combinatorial geometry is the definition of a {\bf cone $\K \subset \R^d$, with an apex $v$}, defined by;
\begin{equation}  \label{def of a cone}
\K:= \left\{ v+  \sum_{k=1}^N  \lambda_k w_k \mid   \lambda _k \geq 0   \right \}.
\end{equation}
The {\bf edge vectors}  of $\K$ are those vectors among the $w_1, \dots, w_N$ (not necessarily all of them) which belong to the boundary 
$\partial \K$ of $\K$.
A fun exercise is to show that the following two conditions are equivalent:
 \begin{enumerate}[(a)]
\item  A cone $\K$ has an apex at the origin.
\item  $\K$ is a cone that enjoys the property $\lambda \K = \K$, for all $\lambda >0$. 
\end{enumerate}
(Exercise \ref{cone equivalence}).

We note that according to definition  \eqref{def of a cone}, an apex need not be unique - in Figure \ref{Cones, pointed and unpointed}, the cone on the left has a unique apex, while the cone on the right has infinitely many apices.
If the vectors $w_1, \dots, w_N$ span a $k$-dimensional subspace of $\R^d$, we say that the 
cone $\K$ has {{\bf dimension $k$}.
When a $k$-dimensional cone $\K \subset \R^d$ has exactly $k$ linearly independent edge vectors $w_1, \dots w_k \in \R^d$, 
 we call such a cone a {\bf simplicial cone}.

  \begin{figure}[htb]
 \begin{center}
\includegraphics[totalheight=2.3in]{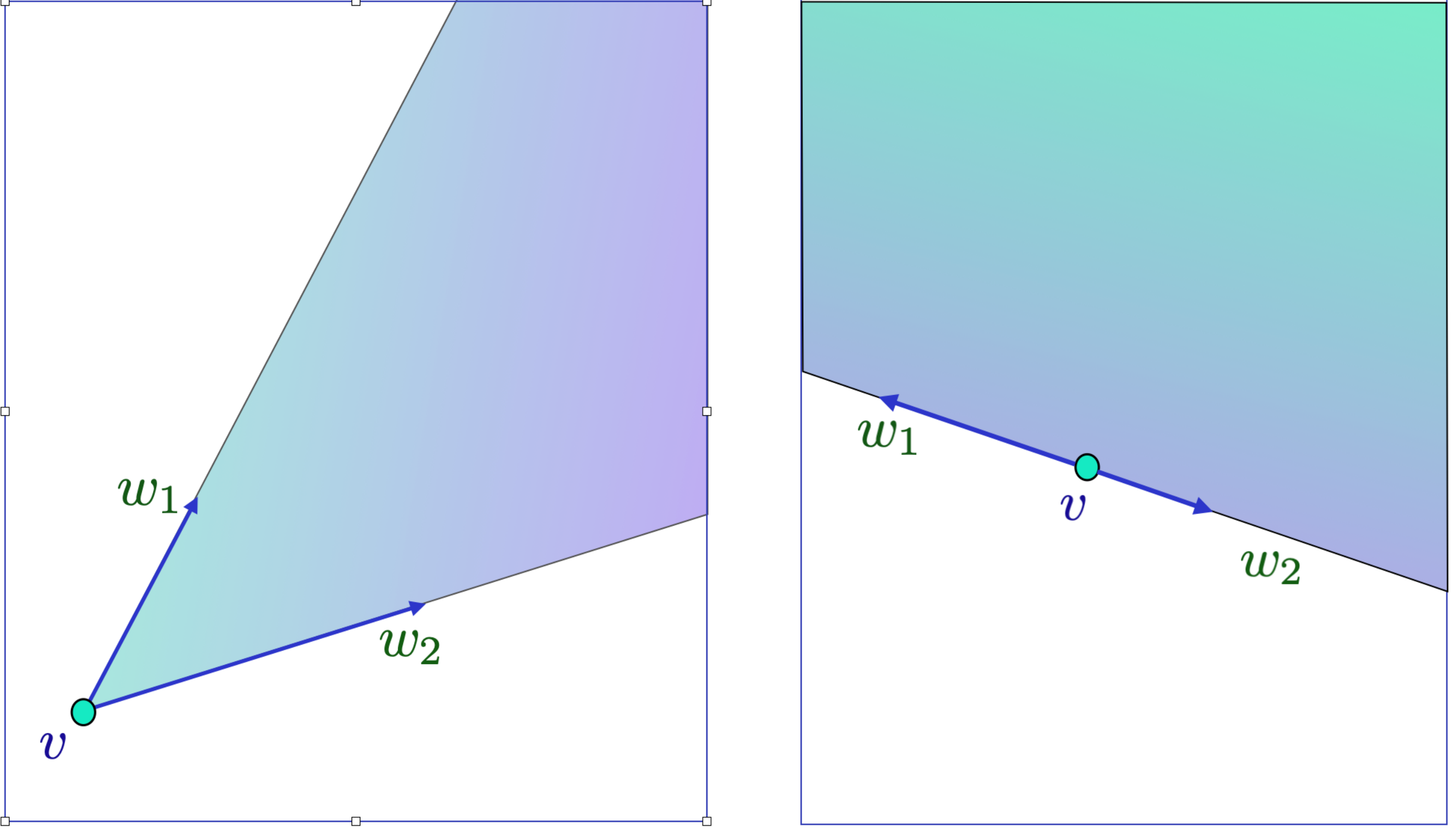}
\end{center}
\caption{The cone on the left is pointed, and has edges $w_1, w_2$.   The cone on the right, with edges  $w_1, w_2$, is not pointed, and in this case it is also a half-space. 
\index{cones} }  \label{Cones, pointed and unpointed} 
\end{figure}

 A {\bf pointed cone}
 \index{cone, pointed}
  is a cone $\K\subset \R^d$ with apex $v$, which enjoys the further property that 
  there exists a hyperplane $H$ with  $H\cap \K = v$.  
  The following $4$ conditions give equivalent characterizations of a pointed cone~$\K$:
 \begin{enumerate}[(a)]
 \item  $\K$ has a unique apex.
\item  There does not exist a vector $u\in \R^d$
 such that $\K + u = \K$.
\item The translated cone $C:= \K-v$, with apex at the origin, enjoys  
$C \cap (-C) = \{0\}$.
 \item  $\K$ does not contain an entire line. \label{non-pointed cone contains a line}
\end{enumerate}
(Exercise \ref{pointed cone equivalence}).  

  We note that every cone has an apex, it's just that the apex may not be unique, for example when $\K$ is a half-space.  All cones are unbounded regions, by definition, so some care will have to be taken when integrating over them.  On the other hand, they are  `almost linear', because for a cone with apex at the origin, we have
\[
x, y \in \K  \ \implies  \ x+y \in \K.
\]
This closure property, which does not exist for polytopes,  makes cones extremely helpful in the analysis of polytopes (for example, Section \ref{section.Brianchon-Gram}).

 An  $n$-dimensional polytope $\P\subset \R^d$  is called a {\bf simplicial polytope} if every facet of $\P$  is a simplex. Equivalently:
 \begin{enumerate}[(a)]
\item Each facet  of $\P$ has  exactly $n$ vertices.
\item Each  $k$-dimensional face of $\P$ has  exactly $k+1$ vertices, for $0\leq k \leq n-1$.
\end{enumerate}
 It is a fun exercise to show that any simplicial cone is always a pointed cone (Exercise \ref{simplicial implies pointed}), but the converse is clearly false.  

By contrast with the notion of a simplicial polytope, we have the following `polar' family of polytopes. 

An $n$-dimensional polytope $\P\subset \R^d$  is called a  {\bf simple} polytope if every vertex is contained in exactly $n$ edges of $\P$.
Equivalently:
 \begin{enumerate}[(a)]
\item Each vertex of $\P$ is contained in exactly $n$ of its facets.
\item  Each $k$-dimensional face of $\P$ is contained in exactly $d-k$ facets, for all $k \geq 0$.
\end{enumerate}




\medskip
\begin{example}
\rm{
Any $d$-dimensional simplex $\Delta$ is a simple polytope.   In fact, any $k$-dimensional
face of the simplex $\Delta$ is also a simplex, and hence a simple polytope of lower dimension.

The $3$-dimensional dodecahedron, in Figure \ref{Dodecahedron}, is also a simple polytope. Its edge graph, which  is always a planar graph for a convex polytope, in this case consists of $20$ vertices, $30$ edges, and $12$ faces.
}
\hfill $\square$
\end{example}

\bigskip
\begin{example}  
\rm{
A  $d$-dimensional simplex also happens to be a simplicial polytope.  The $3$-dimensional icosahedron is a simplicial polytope. 
}
\hfill $\square$
\end{example}

It is a nice exercise to show that the only polytopes which are both simple and simplicial are either
simplices, or $2$-dimensional polygons  (Exercise \ref{simplicial AND simple}).

\begin{example}
\rm{
The $d$-dimensional cube $[0,1]^d$ is a simple polytope.  Its polar polytope, which is the cross-polytope 
$\Diamond$ (see \eqref{cross polytope}), is a simplicial polytope. 
}
\hfill $\square$
\end{example}

One might ask:  are the facets of a simple polytope necessarily simplicial polytopes? 
Again, an example helps here.  
\begin{example}
\rm{
The $120$-cell is a  $4$-dimensional polytope whose  $3$-dimensional boundary  is composed of $120$ dodecahedra \cite{SchleimerSegerman}.
The $120$-cell is a simple polytope,  but because all of its facets are dodecahedra, it does not have any simplicial facets. 
}
\hfill $\square$
\end{example}
 
As becomes apparent after comparing the notion of a simple polytope with that of a simplicial polytope, these two types of polytopes are indeed polar to each other, in the sense of polarity that we've already encountered in definition \eqref{polar polytope, definition}

\begin{lem}
$\P\subset \R^d$ is a simple polytope   $\iff$  $\P^o$ is a simplicial polytope.
\end{lem}
(see Gr\"unbaum \cite{Grunbaum} for a thorough study of this polarity).  
This polarity between simple and simplicial polytopes suggests a stronger connection between our geometric structures thus far, and the combinatorics inherent in the partially ordered set of faces of $\P$.  Indeed, Gr\"unbaum put it elegantly:
 \begin{quote}
 ``In my opinion, the most satisfying way to approach the definition of polyhedra
is to distinguish between the combinatorial structure of a polyhedron,
and the geometric realizations of this combinatorial structure.'' \cite{Grunbaum2}
\end{quote}


 \bigskip
\section{Tangent cones, and the Fourier transform of a simple polytope}

An important step for us is to work with the Fourier-Laplace transform of a cone, and then build some 
 theorems that allow us to simplify many geometric computations, by using the frequency domain on the Fourier transform side.  
 

 We may define the {\bf tangent cone} \index{tangent cone}  of each face  $\F \subset \P$ as follows:
\begin{equation}\label{tangentcone}
\K_{\F} =  \left\{ q+ \lambda(p-q)      \mid  q\in \F,  p\in \P, \lambda \in \R_{\geq 0}   \right\}.
\end{equation}
We note that in general $\K_{\F}$ does not necessarily contain the origin. 
The tangent cone is also known as the {\bf cone of feasible directions}.  Intuitively,  we can imagine standing at the point 
$q\in \F$, and looking in the direction of 
all points that belong to $P$.  Then we take the union of all of these directions.

\begin{figure}[htb]
 \begin{center}
\includegraphics[totalheight=3in]{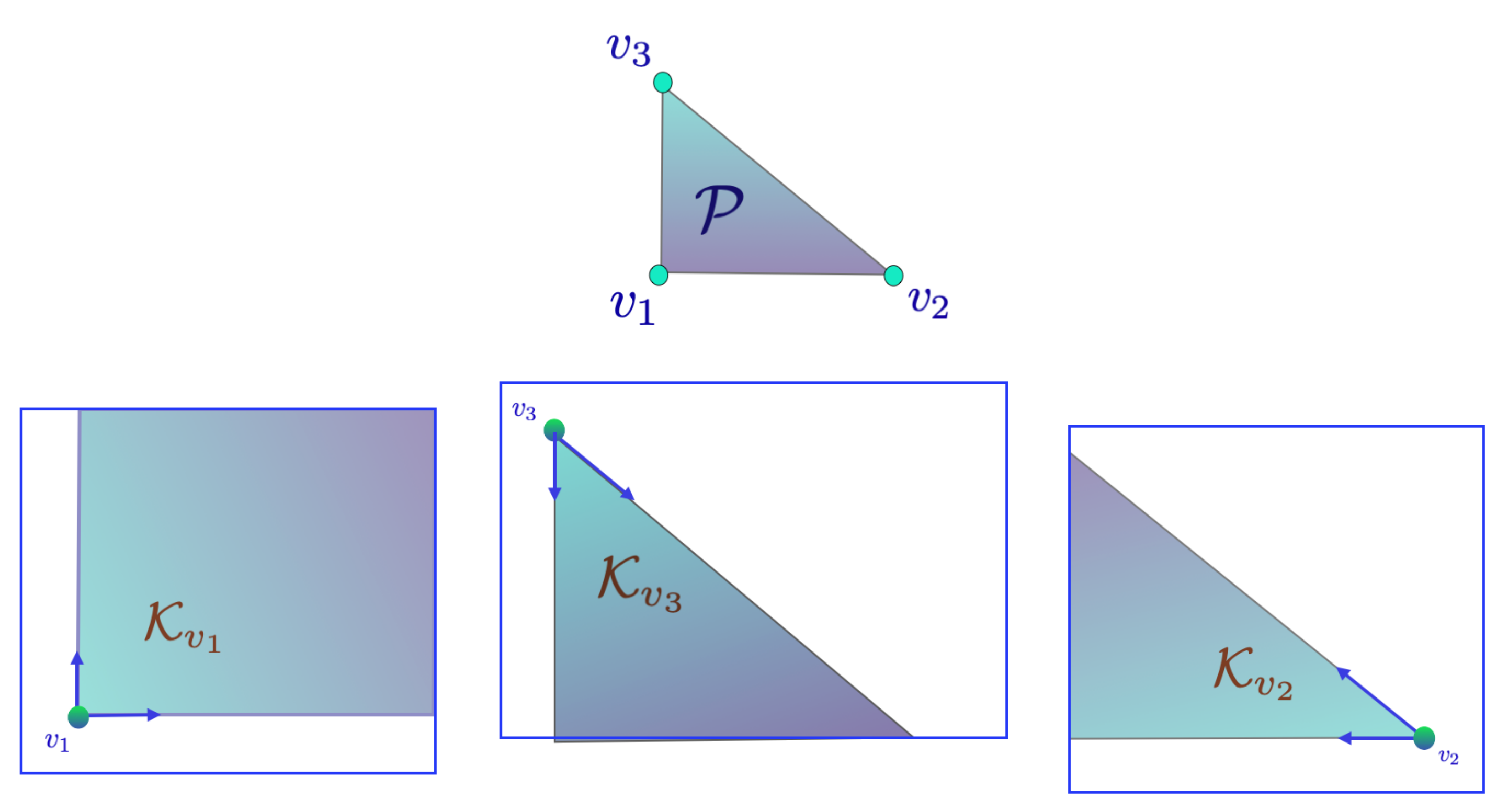}
\end{center}
\caption{The triangle $\P$ has three vertex tangent cones:  $\K_{v_1}, \K_{v_2}, \K_{v_3}$.  The picture is meant to signify that
these cones are, of course, unbounded.
\index{tangent cones} }  
\label{Tangent Cones1} 
\end{figure}

In the case that the face $F$ is a vertex of $\P$, we call this tangent cone a {\bf vertex tangent cone}.
The vertex tangent cone $\K_v$, which is a cone with apex $v$, may also be 
 generated by the edge vectors $v_k - v$, where $ [v_k, v]$ is an edge of $\P$:
  \begin{equation} \label{main definition of tangent cone}
\K_v  = \{  v+ \sum_{k=1}^N  \lambda_k (v_k-v)   \mid \text{ all }  \lambda_k  \geq 0, \text{ and the }
v_k  \text{ are the neighboring vertices of } v\},
  \end{equation}
 a construction we will often use in practice. 
 
 The tangent cone of an edge of a $3$-dimensional convex polytope is an infinite wedge containing the whole line passing through that edge, while the tangent cone of a vertex (for a convex polytope) never contains a whole line (Exercise \ref{Exercise.tangent cone of a vertex}).  
For non-convex polytopes, there are many competing definition for the vertices, and not all of them agree.
 One definition for the vertices of non-convex polytopes appears in \cite{BaranyAkopyanRobins}, using Fourier transforms of cones.  But in this chapter we focus mainly on convex polytopes.


\begin{example}
 \rm{
For the unit cube $\square := [0,1]^d$, the tangent cone at the vertex $v=0$ is
\[
\K_0 = \left\{  \lambda_1 {\bf e_1}+ \lambda_2 {\bf e_2} +  
                      \lambda_3 {\bf e_3}    + \cdots + \lambda_d {\bf e_d}  \mid \lambda_k \geq 0
                      \right\},
\]
which also happens to be the {\bf positive orthant} $\R^d_{\geq 0}$.
On the other hand, the tangent cone of $\square$ at the vertex $v =(1,0, \dots, 0)$ is:
\[
\K_v = v + \left\{   \lambda_1 (-{\bf e_1})+ \lambda_2 {\bf e_2} +  
                      \lambda_3 {\bf e_3}    + \cdots + \lambda_d {\bf e_d}  \mid \lambda_k \geq 0\right\} ,
\]
where ${\bf e_j}$ is the standard unit vector along the $j$'th axis. 
}
\hfill $\square$
\end{example}

 \begin{example} 
 \rm{ To relate some of these definitions, consider a $d$-dimensional simplex $\Delta\subset \R^d$. 
 Located at each of its vertices $v \in \Delta$, we have a tangent cone $K_v$, as in \eqref{main definition of tangent cone}, and here $K_v$ is a simplicial cone.   The simplex $\Delta$ is both a simple polytope and a simplicial polytope.    }
 \hfill $\square$
 \end{example}

 
\bigskip
\section{The Brianchon-Gram identity}   \label{section.Brianchon-Gram}
 \bigskip

The following combinatorial identity, called the Brianchon-Gram identity, may be thought of as a geometric inclusion-exclusion principle.  This identity is quite general, holding true for any convex polytope, simple or not.    For a proof of the following result see, for example, \cite{BarvinokEhrhartbook} or  \cite{BeckRobins}.

 \begin{thm}[Brianchon-Gram identity]\label{Brianchon}
 \index{Brianchon-Gram identity}
Let $\P$ be any convex polytope.  Then
\begin{equation}\label{BG}
 1_\P = \sum_{\F \subseteq \P} (-1)^{dim \F} 1_{\K_F},
\end{equation}
where the sum takes place over all faces of $\P$, including $\P$ itself.  
 \hfill $\square$
 \end{thm}

 \medskip

It turns out that the Brianchon-Gram relations \eqref{BG} can be shown to be equivalent (in the sense that one easily implies the other) 
to the {\bf Euler-Poincare relation} \index{Euler-Poincare relation}  (Exercise \ref{Euler equivalent to Brianchon-Gram})
 for the face-numbers \index{face-numbers} 
of a convex polytope, which says that 
\begin{equation}\label{Euler}
f_0 - f_1 + f_2 - \cdots + (-1)^{d-1} f_{d-1} +  (-1)^{d} f_{d}= 1.  
\end{equation}
Here $f_k$ is the number of faces of $\P$ of dimension $k$.

\medskip
\begin{example} 
\rm{
If we let $\P$ be a $2$-dimensional polygon (including its interior of course) with $V$ vertices, then if must also have $V$ edges, and exactly $1$ face,  so that  \eqref{Euler} tells us that $V - V + 1 = 1$, which is 
not very enlightening, but true.
}
\hfill $\square$
\end{example}

\medskip
\begin{example} 
\rm{
If we let $\P$ be a $3$-dimensional polytope with $V$ vertices, $E$ edge, and $F$ facets, then \eqref{Euler} tells us that $f_0 - f_1 + f_2 - f_3 = 1$, which means that $V - E + F - 1 = 1$.  So we've retrieved Euler's well known formula
\[
V-E+F=2
\] for the Euler characteristic of $3$-dimensional polytopes.
}
\hfill $\square$
\end{example}

\medskip
To gain some facility with the Euler characteristic, we consider if it is possible to construct a polytope in $\R^3$ all of whose facets are hexagons (which are not necessarily regular).  
We claim that this is impossible. 
\begin{lem}
There can be no convex polytope $\P\subset \R^3$ with only hexagonal facets. 
\end{lem}
\begin{proof}
Suppose to the contrary that all the facets of $\P$ are hexagons (not necessarily regular).  By the assumption that $\P$ is a polytope, we know that each edge of $\P$ bounds exactly two facets.   To relate the facets to the edges, consider that 
each facet contains exactly $6$ edges, giving us $6F=2E$.    Combining this latter identity with Euler's formula, we obtain
$V-E+F= V-2F$.

Now let's relate the facets to the vertices.  Each vertex meets at least three facets, and each hexagonal facet contains exactly six vertices.   From the perspective of the facets towards the vertices, we get
$6F \geq 3V$, so that $V \leq 2F$.  Putting things together, we arrive at
\[
2= V-E+F= V-2F  \leq 0, 
\]
a contradiction.
\end{proof}


\section{Brion's formula for the Fourier transform
\\  of a simple polytope}

 Brion \index{Brion} proved the following extremely useful result, Theorem \ref{brion, continuous form},
 concerning the Fourier-Laplace transform
 of a {\it simple polytope}  $\P$.  To describe the result, we consider each vertex $v$ of $\P$, and we fix the $d$ 
 edge vectors $w_1(v), \dots, w_d(v)$ that emanate from $v$.  We recall that the nonnegative real span of the edge vectors $w_k(v)$ generate the vertex 
 tangent cone $\K_v$, and that  these
 edge vectors  are not necessarily required to be unit vectors.   Placing these edge vectors 
 as columns of a matrix $M_v$, we define 
 \[
 \det \K_v := | \det M_v |,
 \]
 the absolute value of the determinant of the ensuing matrix.  

\medskip
\begin{thm}[{\bf Brion's theorem - the continuous form, 1988}] 
\label{brion, continuous form}
\index{Brion's theorem - the continuous form}
Let $\P \subset \R^d$ be a  simple, $d$-dimensional real polytope.  
Then
\begin{equation}\label{transform formula for a simple polytope}
  \int_\P e^{-2\pi i \langle u,  \xi  \rangle} \, du  =  
 \left(     \frac{1}{2\pi i}    \right)^d        \sum_{ v \text{ {\rm a vertex of }} \P } 
  \frac{ e^{-2\pi i \langle v,  \xi   \rangle} \det \K_v}
  { \prod_{ k=1 }^d \langle w_k(v), \xi \rangle } 
\end{equation}
for all $\xi \in \R^d$ such that the denominators on the right-hand side do not vanish.
\hfill $\square$
\end{thm}

Brion's Theorem \ref{brion, continuous form} is one of the cornerstones of Fourier transforms of polytopes.  
We note that the determinant $\det \K_v$ clearly depends on our choice of edge vectors $w_1, \dots, w_d$ for the cone $\K_v$, but it is straightforward (and interesting for applications)
 that the quotient $ \frac{  \det \K_v}{ \prod_{ k=1 }^d \langle w_k(v), \xi \rangle }$ does not depend on the choice of edge vectors (Exercise \ref{independent of edge vectors}).

 This new proof of Brion's theorem uses some of the Fourier techniques that we've developed so far. 
Because we promised a friendly approach, we first give a short outline of the relatively simple ideas of the proof.

 
{\bf Step $1$}. \ We begin with the Brianchon-Gram identity (a standard first step) involving the indicator functions of all of
the tangent cones of $\P$.

\medskip
{\bf Step $2$}. \ We now multiply both sides of the Brianchon-Gram identity \eqref{BG} with the function
 $e^{2\pi i \langle x, \xi \rangle - \varepsilon \| x \|^2}$, where we fix an $\varepsilon >0$, and then we will 
 integrate over all $x \in \mathbb R^d$.  Using these integrals, due to the damped Gaussians for each fixed $\varepsilon >0$, we are able to keep the {\it same domain of convergence} for all of our ensuing functions.
 
\medskip
{\bf Step $3$}. \ Now we let $\varepsilon \rightarrow 0$ and prove that the limit of each integral gives us something meaningful. 
Using integration by parts, we prove  that for any vertex tangent cone $\K$ 
 the corresponding integral
$\int_{\K}  e^{-2\pi i \langle x, \xi \rangle - \varepsilon \| x \|^2} dx$
  converges,  as $\varepsilon \rightarrow 0$, to the desired exponential-rational function.   
  In an analogous but easier manner, we will also prove that the corresponding integral over a 
  non-pointed cone (which includes all faces of positive dimension)  converges to zero, completing the proof.

\medskip \noindent 
In many of the traditional proofs of Theorem \ref{brion, continuous form}, the relevant Fourier-Laplace integrals over the vertex tangent cones have disjoint domains of convergence, lending the feeling that something magical is going on with the disjoint domains of convergence.  Getting around this problem by defining functions that have the same domain of convergence (throughout the proof) was exactly the motivation for this proof.
We favor a slightly longer but clearer expositional proof over a shorter, more obscure proof.   The reader familiar with some physics might notice that this proof idea resembles simulated annealing with a Gaussian.

We also note that throughout the proof we will work over $\xi \in \R^d$, and we don't require any analytic continuation.
Onto the rigorous details of the proof.   First, a technical but crucial Lemma.
 \begin{lem}\label{IntegByParts}
Let $\K_v$ be a $d$-dim'l simplicial pointed cone, with apex $v$, and edge vectors $w_1, \dots, w_d \in \R^d$.  
Then
\begin{equation}\label{LimitDim.d}
 \lim_{\varepsilon \rightarrow 0}   
  \int_{\K_v}   e^{-2\pi i \langle x, \xi \rangle - \varepsilon \|x\|^2} dx
  =
  \ \left( \frac{1}{2\pi i} \right)^d  \frac{ e^{-2\pi i \langle v,  \xi \rangle}  \det \K_v  }
  { \prod_{ k=1 }^d \langle w_k(v) , \xi \rangle },
 \end{equation}
 for all $\xi \in \R^d$ such that $ \prod_{ k=1 }^d \langle w_k(v) , \xi \rangle \not=0$.
 \end{lem}
\begin{proof}
We begin by noticing that we may prove the conclusion in the case that $v=0$, the origin, and for simplicity write $\K_v := \K$ in this case. 
First we make a change of variables, mapping the simplicial cone $\K$ to the nonnegative orthant $\R^d_{\geq 0}$ by the matrix $M^{-1}$, where $M$ is the $d$ by $d$ matrix whose columns are precisely the vectors $w_k$.  Thus, in the integral of \eqref{LimitDim.d}, we let $x:= My$, 
with $y \in \R_{\geq 0}^d$, so that $dx = \left| \det M \right|  dy$.  Recalling that by definition $\det \K =| \det M |$, we have
 \begin{equation}
  \int_{\K}   e^{-2\pi i \langle x, \xi \rangle - \varepsilon ||x||^2} dx
  =
 \left| \det \K \right|
   \int_{\R_{\geq 0}^d}   e^{-2\pi i \langle My, \xi \rangle - \varepsilon ||M y||^2} dy.
 \end{equation}

It is sufficient to therefore show the following limiting identity: 
 \begin{equation}\label{SimplerLimit}
 \lim_{\varepsilon \rightarrow 0}   
  \int_{\R_{\geq 0}^d}   e^{-2\pi i \langle My, \xi \rangle - \varepsilon ||My||^2} dy
  =
  \ \left( \frac{1}{2\pi i} \right)^d  \frac{  1}{ \prod_{ k=1 }^d \langle w_k(v) , \xi \rangle }. 
 \end{equation}

 To see things very clearly, we first prove the $d=1$ case.  Here we must show that 
\begin{equation}\label{LimitDim.1}
 \lim_{\varepsilon \rightarrow 0}   
  \int_{0}^\infty   e^{-2\pi i x \xi  - \varepsilon x^2} dx
  =
  \frac{1}{2\pi i \xi},
 \end{equation}
for all $\xi \in \R-\{0\}$, and we see that even this $1$-dimensional case is interesting.  We proceed with integration by parts by letting
$dv:= e^{-2\pi i x \xi}dx$ and $u:=     e^{ - \varepsilon x^2}$, to get
\begin{align}   
  \int_{0}^\infty   e^{-2\pi i x \xi  - \varepsilon x^2} dx 
  &=    e^{ - \varepsilon x^2} \frac{e^{-2\pi i x \xi}}{-2\pi i \xi} \Big |_{x=0}^{x=+\infty} - 
  \int_0^\infty  \frac{e^{-2\pi i x \xi}}{-2\pi i \xi}  (-2\varepsilon x) e^{ - \varepsilon x^2} dx \\
  &= \frac{1}{2\pi i \xi} - \frac{\varepsilon}{\pi i \xi} 
   \int_{0}^\infty   x e^{-2\pi i x \xi  - \varepsilon x^2} dx \\  \label{last nasty integral} 
   &=  \frac{1}{2\pi i \xi} - \frac{1 }{\pi i \xi} 
   \int_{0}^{\infty}   e^{-2\pi i  \frac{u}{\sqrt{\varepsilon}}   \xi  }  u e^{-u^2} du 
\end{align}
where we've used the substitution $u:= \sqrt{\varepsilon}  x$ in the last equality \eqref{last nasty integral}.   
We now notice that
\[
\lim_{\varepsilon \rightarrow 0}   \int_{0}^{\infty}   e^{-2\pi i  \frac{u}{\sqrt{\varepsilon}}   \xi  }  u e^{-u^2} du
=\lim_{\epsilon \rightarrow 0}   \hat g\Big(\frac{\xi}{\sqrt\epsilon}\Big),
\]
where $g(u):=u e^{-u^2}1_{[0, +\infty]}(u)$ is an absolutely integrable function. 
Luckily, we know by the Riemann--Lebesgue lemma \ref{Riemann--Lebesgue lemma} 
\index{Riemann-Lebesgue lemma}  
that 
\[
\lim_{|w| \rightarrow \infty}   \hat g(w) =0, 
\]
and so we arrive at the desired limit \eqref{LimitDim.1}.

We now proceed with the general case, which just uses the $1$-dimensional idea above several times.  To prove \eqref{SimplerLimit}, we first fix the variables $y_2, \dots, y_d$ and perform integration by parts on $y_1$ first.   Thus, we let 
\begin{align}
dv_1 &:= e^{-2\pi i \langle My, \xi \rangle} dy_1=  
e^{-2\pi i \langle y, M^t \xi \rangle} dy_1= e^{-2\pi i \Big( y_1  \langle w_1, \xi \rangle + \cdots + 
y_d  \langle w_d, \xi \rangle \Big)}dy_1,
\end{align}
 thought of as a function of only $y_1$.   
 Carrying out the integration in the variable $y_1$, we have
$v_1 = e^{-2\pi i \langle y, M^t \xi \rangle} / \left(- 2\pi i \langle w_1, \xi \rangle \right)$.   
We let $u_1:=  e^{- \varepsilon ||My||^2}  $, also thought of as a function of $y_1$ alone. 
 We have $du_1 = -\varepsilon L(y) e^{- \varepsilon ||My||^2} dy_1$, where $L(y)$ is a real polynomial in $y$, whose coefficients come from the entries of $M$.
Integrating by parts in the variable $y_1$ now gives us
\begin{align}\label{SimplerLimitProof}  
 & \int_{\R_{\geq 0}^d}   e^{-2\pi i \langle My, \xi \rangle - \varepsilon ||My||^2} dy
  =
 \int_{\R_{\geq 0}^{d-1}} dy_2 \cdots dy_d  \left[   u_1 v_1 \Big |_0^\infty - \int_0^\infty v_1 du_1     \right]  \\
 & =  \int_{\R_{\geq 0}^{d-1}} dy_2 \cdots dy_d  \left[   
 \frac{ e^{-2\pi i \langle y, M^t \xi \rangle - \varepsilon ||M y||^2} }{  -2\pi i \langle w_1, \xi \rangle  }
  \Big |_{y_1=0}^{y_1 = \infty}  +   \frac{\varepsilon}{-2\pi i  \langle w_1, \xi \rangle}    
  \int_0^\infty   L(y) e^{-2\pi i \langle y, M^t \xi \rangle - \varepsilon ||My||^2} dy_1
   \right]  \\
  & =  \int_{\R_{\geq 0}^{d-1}} 
 \frac{ e^{2\pi i \langle t, M^t \xi \rangle - \varepsilon ||M t||^2}  }{  2\pi i \langle w_1, \xi \rangle  }dt 
  -  \frac{\varepsilon}{2\pi i  \langle w_1, \xi \rangle}    
  \int_{\R_{\geq 0}^{d}}  L(y) e^{-2\pi i \langle y, M^t \xi \rangle - \varepsilon || M y ||^2} dy \\
& =   \frac{1}{  2\pi i \langle w_1, \xi \rangle  }   \int_{\R_{\geq 0}^{d-1}} 
 e^{-2\pi i \langle t, M^t \xi \rangle - \varepsilon || M t ||^2} dt 
  -  \frac{\varepsilon}{2\pi i  \langle w_1, \xi \rangle}    
  \int_{\R_{\geq 0}^{d}}  L(y) e^{-2\pi i \langle y, M^t \xi \rangle - \varepsilon || M y ||^2} dy,
\end{align}
where we've used $t:= (y_2, \dots, y_d)$ in the $3$'rd equality.   We repeat exactly the same process of integration by parts as in \eqref{last nasty integral}, one variable at a time.  We observe that after $d$ iterations we get a sum of $d$ terms, where the first term does not contain any $\varepsilon$ factors, while
 all the other terms do contain $\varepsilon$ factors in the exponents.  Therefore, when we complete the $d$-many integration by parts iteratively, and 
finally let $\varepsilon$ tend to zero, only the leading term remains, namely 
$\left( \frac{-1}{2\pi i} \right)^d  \frac{  1}{ \prod_{ k=1 }^d \langle w_k,  \xi \rangle } $.  We've shown that \eqref{SimplerLimit} is true.
\end{proof}

\bigskip
\begin{proof}(of Theorem \ref{brion, continuous form})
We begin with the Brianchon Gram identity:
\begin{equation}\label{BG2}
1_\P = \sum_{\F \subseteq \P} (-1)^{dim \F} 1_{K_F}.
\end{equation}
We fix any $\xi \in \R^d$, and any $\varepsilon > 0$.  Multiplying both sides of \eqref{BG2} by 
$e^{-2\pi i \langle x, \xi \rangle - \varepsilon \|x\|^2}$, and integrate over all $x \in \R^d$, we have: 
\begin{equation}
\int_{\R^d} 1_\P(x)  e^{-2\pi i \langle x, \xi \rangle - \varepsilon \|x\|^2} dx
= \sum_{\F \subseteq \P} (-1)^{dim \F} \int_{\R^d} 1_{K_F}(x) 
 e^{-2\pi i \langle x, \xi \rangle - \varepsilon \|x\|^2} dx.
\end{equation}

\noindent
Equivalently, 
\begin{equation}\label{IntegratingBrianchon}
\int_{\P}  e^{-2\pi i \langle x, \xi \rangle - \varepsilon \|x\|^2} dx
= \sum_{\F \subseteq \P} (-1)^{dim \F} 
\int_{\K_F} 
 e^{-2\pi i \langle x, \xi \rangle - \varepsilon \|x\|^2} dx.
\end{equation}

For each fixed $\varepsilon > 0$, all integrands in \eqref{IntegratingBrianchon} 
are Schwartz functions, and so all of the integrals in the latter identity now converge absolutely (and rapidly).
We identify two types of tangent cones that may occur on the right-hand side of \eqref{IntegratingBrianchon}, for each face $\F \subseteq \P$. 

{\bf Case $1$}.   When $\F = v$, a vertex, we have the vertex tangent cone $\K_v$:  these are the tangent cones that exist for each vertex of $\P$.
It is a standard fact that all of these vertex tangent cones are pointed cones.  By hypothesis, all of our vertex tangent cones are simplicial cones, so letting $\varepsilon \rightarrow 0$ and calling on 
Lemma \ref{IntegByParts}, we obtain the required limit for 
$\int_{\K_v}      e^{2\pi i \langle x, \xi \rangle - \varepsilon \|x\|^2} dx$.

\smallskip
{\bf Case $2$}.   When $\F$ is not a vertex, we have the tangent cone $\K_\F$, and it is a standard fact
that in this case  $\K_\F$ always contains a line.  
Another standard fact in the land of polytopes is that each  tangent cone in this case 
may be written as  $\K_{\F} = \R^k \oplus \K_p$, the direct sum of a copy of Euclidean space with a pointed cone $\K_p$ for any point $p \in \F$.   As a side-note, it is also true that $\dim \F = k$.

We  would like to show that for all faces $\F$ that are not vertices of $\P$, the associated integrals tend to $0$:
\[
\int_{\K_F}    e^{-2\pi i \langle x, \xi \rangle - \varepsilon \|x\|^2} dx \rightarrow 0,
\]
as $\varepsilon \rightarrow 0$.  Indeed,
\begin{align}
\int_{\K_F} 
 e^{-2\pi i \langle x, \xi \rangle - \varepsilon \|x\|^2} dx & =  
 \int_{  \R^k \oplus \K_p  }    e^{-2\pi i \langle x, \xi \rangle - \varepsilon \|x\|^2} dx \\
 & =  \int_{  \R^k   }    e^{-2\pi i \langle x, \xi \rangle - \varepsilon \|x\|^2} dx
  \int_{  \K_p  }    e^{-2\pi i \langle x, \xi \rangle - \varepsilon \|x\|^2} dx. \label{product}
\end{align}
The integral $ \int_{  \R^k   }    e^{-2\pi i \langle x, \xi \rangle - \varepsilon \|x\|^2} dx$ 
is precisely the usual Fourier transform of a Gaussian, which is known to be the Gaussian 
$G_{\varepsilon}(x):= \varepsilon^{-k/2} e^{-\frac{\pi}{\varepsilon} \|x\|^2}$ by Exercise \ref{Gaussian2}. 
It is apparent that for any fixed nonzero value of $x\in \R^k$, we have 
$\lim_{\varepsilon \rightarrow 0}G_{\varepsilon}(x)=0$.   
Finally, by Lemma \ref{IntegByParts} again, the limit
 $\lim_{\varepsilon \rightarrow 0} \int_{  \K_p  }    
 e^{-2\pi i \langle x, \xi \rangle - \varepsilon \|x\|^2} dx$ is finite, because $\K_p$ is another pointed cone.
Therefore the product of the integrals in
  \eqref{product}          tends to zero, completing the proof.
\end{proof}

\section{The Fourier transform of any real polytope}

 Brion's theorem, which holds for simple polytopes,  is particularly useful  whenever we are given a polytope in terms of its local data at the vertices -  including the edge vectors for each vertex tangent cone.  We can then easily write down the Fourier transform of a simple polytope, by Theorem  \ref{brion, continuous form}. 
What happens, though, for non-simple polytopes? 
There is the following natural extension of Brion's Theorem \ref{brion, continuous form}  to all real polytopes, which is now easy to prove.

\medskip
\begin{thm}[{\bf Fourier transform of any real polytope}]
 \label{brion2}
Let $\P \subset \R^d$ be any $d$-dimensional  polytope.   Then:
\begin{equation} \label{the main deal!}
  \int_\P e^{-2\pi i \langle u,  \xi  \rangle} \, du  =
 \sum_{v \in V}   
\frac{e^{-2\pi i \langle v, \xi \rangle} }{(2\pi i)^d}     
\sum_{j=1}^{M(v)}   \frac{\det \K_j(v)  }{\prod_{k=1}^d  \langle  w_{j, k}(v),   \xi     \rangle},
\end{equation}
for all $\xi \in \R^d$ such that none of the denominators vanish:
$ \prod_{k=1}^d  \langle  w_{j, k}(v),   \xi     \rangle  \not=0$.  At each vertex vertex $v \in \P$, the vertex tangent cone $\K_v $
is triangulated into simplicial cones, using the notation $\K_v = \K_1(v)\cup \dots \cup \K_{M(v)}(v)$.
\end{thm}
\begin{proof}
The proof here is identical in almost every aspect to the proof of Theorem \ref{brion, continuous form}, except for {\bf Case} $1$ of its proof, above.   By contrast with the proof above of {\bf Case} $1$, here our vertex tangent cones $\K_v$ need not be simplicial.   However, we may triangulate each vertex tangent cone $\K_v$ into simplicial cones $\K_1(v)$, \dots $\K_{M(v)}(v)$, so that we have the disjoint union
$\K_v = \K_1(v)\cup \dots \cup \K_{M(v)}(v)$.   Therefore 
\begin{align*}
 \lim_{\varepsilon \rightarrow 0}   
       \int_{\K_v}   
        e^{-2\pi i \langle x, \xi \rangle - \varepsilon \|x\|^2} dx
  &= \lim_{\varepsilon \rightarrow 0} 
   \sum_{j=1}^{M(v)}     
       \int_{  \K_{j,v}   }  
       e^{-2\pi i \langle x, \xi \rangle - \varepsilon \|x\|^2} dx \\
 &=  \sum_{j=1}^{M(v)}  
   \lim_{\varepsilon \rightarrow 0}    
       \int_{  \K_{j,v}   }  
       e^{-2\pi i \langle x, \xi \rangle - \varepsilon \|x\|^2} dx \\
  &=  \ \left( \frac{-1}{2\pi i} \right)^d  
  \sum_{j=1}^{M(v)}  
  \frac{ e^{-2\pi i \langle v,  \xi \rangle}  \det \K_{j}(v) }
  { \prod_{ k=1 }^d \langle w_{j, k}(v) , \xi \rangle },
 \end{align*}
 where we've used Lemma \ref{IntegByParts} in the last equality, owing to the fact that all of the cones
 $\K_j(v)$ are simplicial.
 The calculation above is valid 
  for each $\xi \in \R^d$ such that $ \prod_{ k=1 }^d \langle w_{j, k}(v) , \xi \rangle \not=0$ for all vertices $v$ and all $j = 1, \dots, M(v)$.
\end{proof}

The nonvanishing condition 
$ \prod_{k=1}^d  \langle  w_{j, k}(v),   \xi     \rangle  \not=0$ may be restated more combinatorially as follows.  Let $\mathcal H$ be the finite union of hyperplanes, where each hyperplane is defined by 
\begin{equation}\label{def:the hyperplane arrangement defined by zero set}
\mathcal H:=  \{ \xi \in \R^d \mid   \langle  w_{j, k}(v),   \xi     \rangle =0\}.
\end{equation}
 In other words, $\mathcal H$ is the union of all hyperplanes that are orthogonal to any edge of $\P$.  So the only restriction in \eqref{the main deal!} is that $\xi \notin \mathcal H$.  But again we emphasize that these `singularities' are removable singularities because after extending both sides to all $\xi \in \C^d$,  the left-hand side of  \eqref{the main deal!} is an entire function of $\xi \in \C^d$.


\bigskip
\section{Fourier-Laplace transforms of cones}
\label{Fourier Laplace transforms of cones}

What about the Fourier transform of a cone?  Well, if we naively try to use the same integrand over a cone, the integral will diverge. 
  But there is a way to fix this divergence by replacing the real vector 
$\xi \in \R^d$ by a complex vector $z \in \C^d$.
Let's consider what would happen if we formally replace the variable $\xi \in \R^d$ by a
complex vector $z := x+iy \in \C^d$, to obtain the transform:
\[
1_\P(z):= \int_\P   e^{-2\pi i \langle u,  z\rangle} \, du.
\]
Our inner product $\langle u, z \rangle := u_1 z_1 + \cdots + u_d z_d$ is always 
the usual inner product on $\R^d$, defined
 without using the Hermitian inner product here.  In other words, we simply use the usual  inner product on $\R^d$, and then formally substitute complex numbers $z_k$ into it.
 This means, by definition, that
\begin{align}
\int_\P e^{-2\pi i \langle u,  z\rangle} \, du   &:= \int_\P e^{-2\pi i \langle u,  x+iy \rangle} \\
&:= \int_\P e^{-2\pi i \langle u,  x\rangle} e^{2\pi \langle u,  y\rangle} \,  du,
\end{align}
so that we have an extra useful real factor of $e^{2\pi \langle u,  y\rangle}$ that makes the integral 
converge quite rapidly over unbounded domains, provided that 
$ \langle u,  y \rangle < 0$.
If we set $y=0$, then it's clear that we retrieve the usual Fourier transform of $\P$, while if we set $x=0$, we get a new integral, which we call the {\bf Laplace transform} of $\P$.  Finally, the {\bf Fourier-Laplace transform} \index{Fourier-Laplace transform} 
of $\P$ is defined by:
\[
\hat 1_\P(z) := \int_\P e^{-2\pi i \langle u,  z\rangle} \, du
\]
valid for any $z \in \C^d$ for which the integral converges.

One clear reason for the use and flexibility of the full Fourier-Laplace transform (as opposed to just the Fourier transform) 
is the fact that for a cone $\K$,
its usual Fourier transform diverges.  But if we allow a complex variable $z\in \C^d$, then the integral does converge on a restricted domain.   Namely, the Fourier-Laplace transform of 
a cone $\K$ is defined by:
\[
\hat 1_\K(z) := \int_\K e^{-2\pi i \langle u,  z\rangle} \, du,
\]
for a certain set of $z\in \C^d$, but we can easily understand its precise domain of convergence. 
For an arbitrary cone $\K \subset \R^d$, we define its {\bf dual cone} 
\index{dual cone} by:
\begin{equation}\label{dual cone definition}
\K^* := \{  y \in \R^d \mid \langle y, u \rangle < 0 \text{ for all } u\in \K  \},
\end{equation}
which is an open cone.  As one might expect, there is an easy duality: 
$\K_1 \subset \K_2  \iff  \K_2^* \subset \K_1^*$ (Exercise  \ref{duality of dual cone}).



\bigskip
\begin{example}
\rm{
Given the $1$-dimensional cone $\K_0 := \R_{\geq 0}$, we compute its Fourier-Laplace transform:
\begin{align*}  
  \int_{\K_0} e^{-2\pi i u z} \, du =         \int_0^\infty e^{-2\pi i u z} \, du =           &=   \frac{1}{-2\pi i z}  e^{-2\pi i u (x+iy)}\Big|_{u=0}^{u=\infty} \\
 &=  \frac{1}{-2\pi i z}  e^{-2\pi i ux}   e^{2\pi uy}\Big|_{u=0}^{u=\infty} \\
  &=  \frac{1}{-2\pi i z}  (0-1)
  = \frac{1}{2\pi i}    \frac{1 }{    z   },
\end{align*} 
valid for all $z:= x + iy\in \C$ such that $y < 0$.  We note that for such a fixed complex 
$z$, $| e^{-2\pi i u z} | = e^{2\pi  u y }$ is a rapidly decreasing function of $u\in \R_{>0}$, 
because $y<0$.
}
\hfill $\square$
\end{example}

\begin{figure}[htb]
 \begin{center}
\includegraphics[totalheight=2.5in]{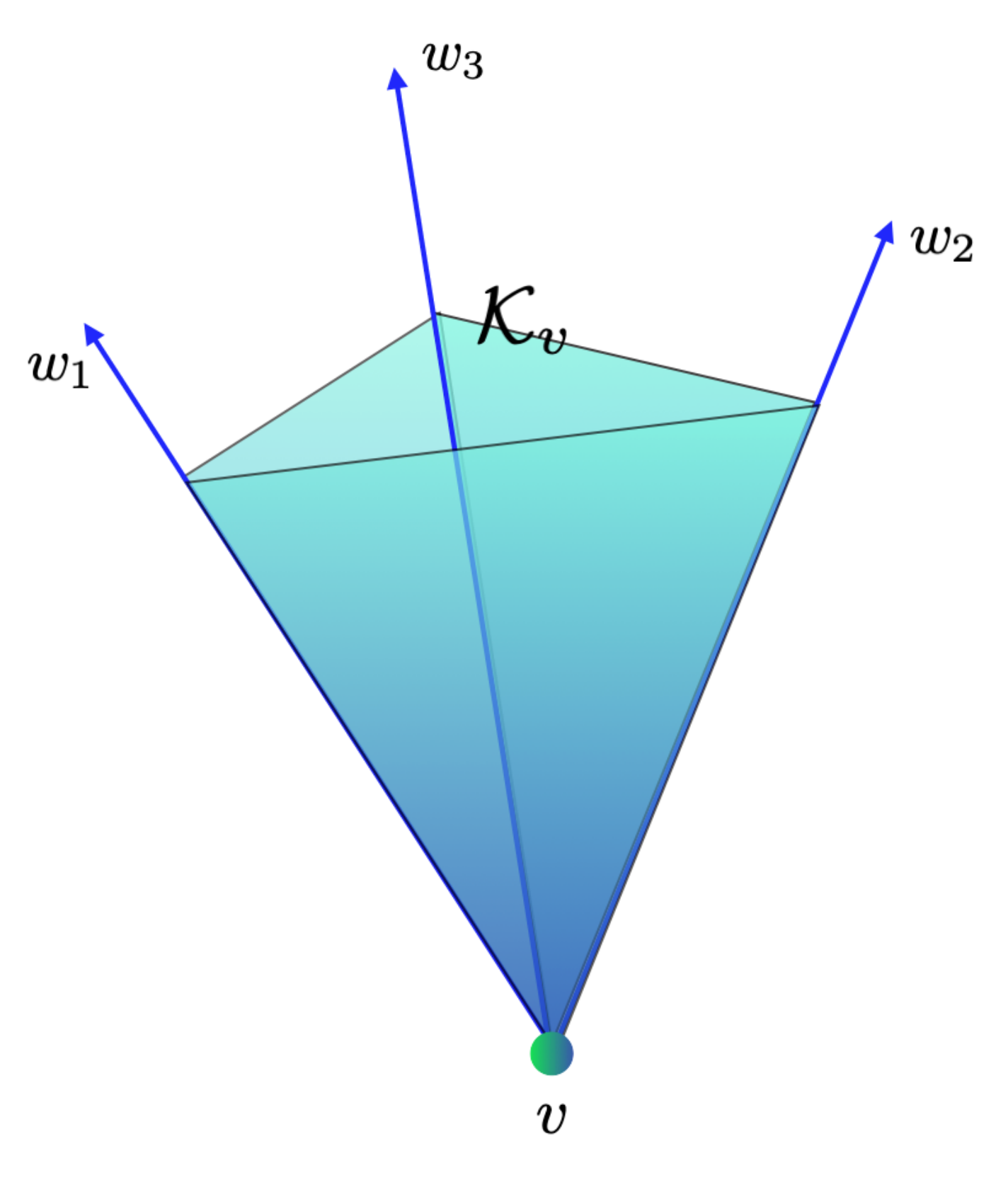}
\end{center}
\caption{A simplicial, pointed cone in $\R^3$, with apex $v$ and edge vectors $w_1, w_2, w_3$}  \label{Simplicial cone} 
\end{figure}

Now let's work out the Fourier-Laplace transform of a $d$-dimensional cone whose apex is the origin.
 \begin{lem} \label{F-L transform of a simplicial cone, apex at o}
 Let $\K \subset \R^d$ be a simplicial, $d$-dimensional cone, with apex at the origin.  
 If the edges of $\K$ are labelled $w_1, \dots, w_d$, then 
 \[
 \hat 1_K(z) :=  \int_\K e^{-2\pi i \langle u,  z\rangle} \, du  =   \frac{1}{(2\pi i)^d}    \frac{\det \K }{\prod_{k=1}^d  \langle  w_k, z     \rangle}.
  \]
Furthermore, the {\bf domain of convergence} for the latter integral is naturally associated with the dual cone, \index{dual cone} and it is given by:
 \[
 \{ z:= x + iy \in \C^d \mid   \  y \in \K^*  \}.
 \]
 \end{lem}
 \begin{proof}
We first compute the Fourier-Laplace transform of the positive orthant $\K_0 := \R_{\geq 0}^d$, with a complex vector $z = x + i y \in \C^d$:
\begin{align}   \label{transform of a cone}
\hat 1_{\K_0}(z) &:= \int_{\K_0}  e^{-2\pi i \langle z, u \rangle} du \\
&=   \int_{\R_{\geq 0}}  e^{-2\pi i z_1 u_1} du_1  \cdots  \int_{\R_{\geq 0}}  e^{-2\pi i z_d u_d} d u_d \\
&=  \prod_{k=1}^d    \frac{ 0- 1}{-2\pi i z_k}  
=  \left( \frac{1}{2\pi i}\right)^d    \frac{1}{ z_1 z_2 \cdots z_d}.   \label{trick2}
\end{align} 
Next, the positive orthant $\K_0$ may be mapped to the cone $\K$ by a linear transformation.  Namely, we may use the matrix $M$ whose columns are defined to be the edges of $\K$, so that by definition $\K = M(\K_0)$.   Using this mapping, we have:
\begin{align*}  
\hat 1_{\K}(z) &:= \int_{\K}  e^{-2\pi i \langle z, u \rangle} du \\
&=  |\det M|   \int_{\K_0}  e^{-2\pi i \langle z, M t \rangle} dt   \\
&=  |\det M|    \int_{\K_0}  e^{-2\pi i \langle M^T z,  t \rangle} dt   \\
&=     \left( \frac{1}{2\pi i}\right)^d    \frac{|\det M| }{\prod_{k=1}^d  \langle  w_k, z     \rangle}.
\end{align*} 
where in the second equality we've made the substitution $u = Mt$, with $t\in \K_0, u \in \K$, and 
$du = |\det M| dt$.   In the final equality, we used equation \eqref{trick2} above, noting that the $k$'th 
element of the vector $M^Tz$ is $\langle  w_k, z     \rangle$, and we note that by definition $|\det M| = \det \K$.

For the domain of convergence of the integral, we observe that 
\[
e^{-2\pi i \langle u, z\rangle} = e^{-2\pi i \langle u, x+iy\rangle} = 
e^{-2\pi i \langle u, x\rangle} e^{2\pi  \langle u, y\rangle},
\]
and because $ \left|   e^{-2\pi i \langle u, x\rangle}  \right| = 1$, the integral 
$\int_\K e^{-2\pi i \langle u, z\rangle} du$
converges   $\iff   \langle u, y \rangle < 0$ for all $u\in \K$.  But by definition of the dual cone, this means that $y \in \K^*$.
\end{proof}

\bigskip
\begin{example}
\rm{
Given the $2$-dimensional cone 
$\K := \{ 
\lambda_1 \big(\begin{smallmatrix}
1  \\
5 \\
\end{smallmatrix}
\big) +    
\lambda_2 \big(\begin{smallmatrix}
-3  \\
\ 2 \\
\end{smallmatrix}
\big)  \mid  \lambda_1, \lambda_2 \in \R_{\geq 0}  \}$, we compute its Fourier-Laplace transform, and find its domain of convergence.  By Lemma \ref{F-L transform of a simplicial cone, apex at o},
\begin{align*}  
 \hat 1_{\K}(z):=  \int_\K e^{-2\pi i \langle u,  z \rangle} \, du  &= 
      \frac{1}{(2\pi i)^2}    \frac{17}{(z_1 + 5z_2)(-3z_1 + 2z_2)},
\end{align*} 
valid for all $z =  \icol{z_1\\z_2} := x + iy$ such that $ y \in \K^*$.  Here the dual cone is 
given here by \\
$\K^* = \interior\{ 
\lambda_1 \big(\begin{smallmatrix}
\ 5  \\
-1 \\
\end{smallmatrix}
\big) +    
\lambda_1 \big(\begin{smallmatrix}
-2  \\
-3 \\
\end{smallmatrix}
\big)  \mid  \lambda_1, \lambda_2 \in \R_{\geq 0}  \}$.
}
\hfill $\square$
\end{example}

To compute the Fourier-Laplace transform of a simplicial cone $\K$ whose apex is $v \in \R^d$, we may first compute the transform of the translated 
 cone $\K_0:= \K - v$, whose apex is at the origin, using the previous lemma.
We can then use the fact that the Fourier transform behaves in a simple way under translations, namely
\[
\hat 1_{K + v}(z) = e^{2\pi i \langle z, v \rangle} \hat 1_K(z), 
\]
to obtain the following result (Exercise \ref{translating a cone}).

 \bigskip
\begin{cor} \label{transform of a translated cone}
Let $\K_v \subset \R^d$ be a simplicial $d$-dimensional cone, whose apex is $v \in \R^d$.  Then
\begin{equation} \label{IMPORTANT cone transform}
{\hat 1}_{\K_v}(z) :=   \int_{\K_v}   e^{-2\pi i \langle u,  z\rangle} \, du  =   \frac{1}{(2\pi i)^d}   
 \frac{ e^{-2\pi i \langle v, z \rangle}  \det \K_v  }{\prod_{k=1}^d  \langle  w_k, z     \rangle},
\end{equation}
a rational-exponential function.
More generally, for any $d$-dimensional cone $\K_v\subset \R^d$ with apex $v$, we can 
 always triangulate $\K_v$ into $M(v)$ simplicial subcones $\K_j(v)$ \cite{DRS}, and apply the previous result to each simplicial 
 subcone, obtaining:
\begin{equation} \label{general cone transform}
{\hat 1}_{\K_v}(z) :=   \int_{\K_v}   e^{-2\pi i \langle u,  z\rangle} \, du  =   
\frac{e^{-2\pi i \langle v, z \rangle} }{(2\pi i)^d}     
\sum_{j=1}^{M(v)}   \frac{\det \K_j(v)  }{\prod_{k=1}^d  \langle  w_{j, k}(v),   z     \rangle},
\end{equation}
a rational-exponential function.  
\hfill $\square$
\end{cor}

For a non-simple polytope, the question of computing efficiently the Fourier-Laplace transforms of all of its tangent cones becomes unwieldy, as far as we know (this problem is related to the $P \not= NP$ problem).  In fact, even computing the volume of a polytope is already known to be NP-hard in general, and the volume is just the Fourier transform evaluated at one point:  $\vol \P = 1_\P(0)$.

\medskip
 \begin{example}
 \rm{
 Let's work out a $2$-dim'l example of Brion's Theorem \ref{brion, continuous form}, using Fourier-Laplace transforms of  tangent cones.  
We will find the rational-exponential function for the Fourier-Laplace transform of the triangle $\Delta$, whose vertices are defined by 
 $v_1:= \icol{0\\0}$,  $v_2:= \icol{a\\0}$, and $v_3:= \icol{0\\b}$, with $a>0, b>0$.  
 
 First, the tangent cone at the vertex $v_1:= \icol{0\\0} $ is simply the nonnegative orthant in this case, with edge vectors $w_1 = \icol{1\\0}$ and $w_2 = \icol{0\\1}$.  Its determinant, given these two edge vectors, is equal to $1$.  Its Fourier-Laplace transform is 
 \begin{equation}
  \int_{\K_{v_1}} e^{-2\pi i  \langle x,  z\rangle} \, dx = 
   \  \frac{1}{(2\pi i)^2} \,  \frac{1}{z_1 z_2},
 \end{equation}
 and note that here we must have both $\Im(z_1)>0$ and $\Im( z_2)>0$ in order to make the integral converge.  Here we use the standard notation $\Im(z)$ is the imaginary part of $z$.

The second tangent cone at vertex $v_2$ has edges $w_1 = \icol{-a\\   \  b}$ and 
$w_2 = \icol{ \   0\\    -b}$ (recall that we don't have to normalize the edge vectors at all).  Its determinant has absolute value equal to $ab$, and its 
Fourier-Laplace transform is

 \begin{equation}
  \int_{\K_{v_2}} e^{-2\pi i  \langle x,  z\rangle} \, dx = 
    \left( \frac{1}{2\pi i} \right)^2  \frac{(ab) e^{-2\pi i a z_1} }{(-a z_1 + b z_2)(-a z_1)},
\end{equation}
 and here the integral converges only for those $z$ for which $\Im( -az_1 + bz_2) >0$ and $\Im( -a z_1 ) >0$.
 
Finally, the third tangent cone at vertex $v_3$ has edges $w_1 = \icol{ \ a\\  -b}$ and 
$w_2 =\icol{\ 0\\ -b}$. Its determinant has absolute value equal to $ab$, 
and its Fourier-Laplace transform is

 \begin{equation}
  \int_{\K_{v_3}} e^{-2\pi i  \langle x,  z\rangle} \, dx = 
    \left( \frac{1}{2\pi i} \right)^2  \frac{(ab) e^{-2\pi i b z_2} }{(a z_1 - b z_2)(-b z_2)}.
\end{equation}
 and here the integral converges only for those $z$ for which $\Im( az_1 - bz_2) >0$ and $\Im( -b z_2 ) >0$.

 We can again see quite explicitly the disjoint domains of convergence in this example, so that there is not even one value of $z \in \C^2$ for which all three Fourier-Laplace transforms of all the tangent cones converge simultaneously.  Despite this apparent shortcoming, Brion's identity \eqref{brion, continuous form} still tells us that we may somehow still add these local contributions of the integrals at the vertices combine to give us a formula for the Fourier-Laplace transform of the triangle:
 \begin{equation}
\hat 1_{\Delta}(z) := \int_{\Delta} e^{-2\pi i  \langle x,  z\rangle} dx = 
    \left( \frac{1}{2\pi i} \right)^2 
     \left(
       \frac{1}{z_1 z_2}
    +  \frac{-b\ e^{-2\pi i a z_1} }{(-a z_1 + b z_2) z_1}
    + \frac{ -a \ e^{-2\pi i b z_2} }{(a z_1 - b z_2) z_2}
     \right),
\end{equation} 
which is \emph{now} magically valid for all generic $(z_1, z_2) \in \C^2$; in other words, it is now valid for all
$(z_1, z_2) \in \C^2$ except those values which make the denominators vanish.
}
\hfill $\square$
 \end{example}


\section{The Fourier transform of a polygon} \index{polygon, Fourier transform}

Here we give an efficient formula for the Fourier transform of any polygon, namely 
Corollary~\ref{FT of a polygon!}.
Let's begin with a simple and natural question:  what is the Fourier transform of a hexagon?

\begin{figure}[htb]
\begin{center}
\includegraphics[totalheight=3in]{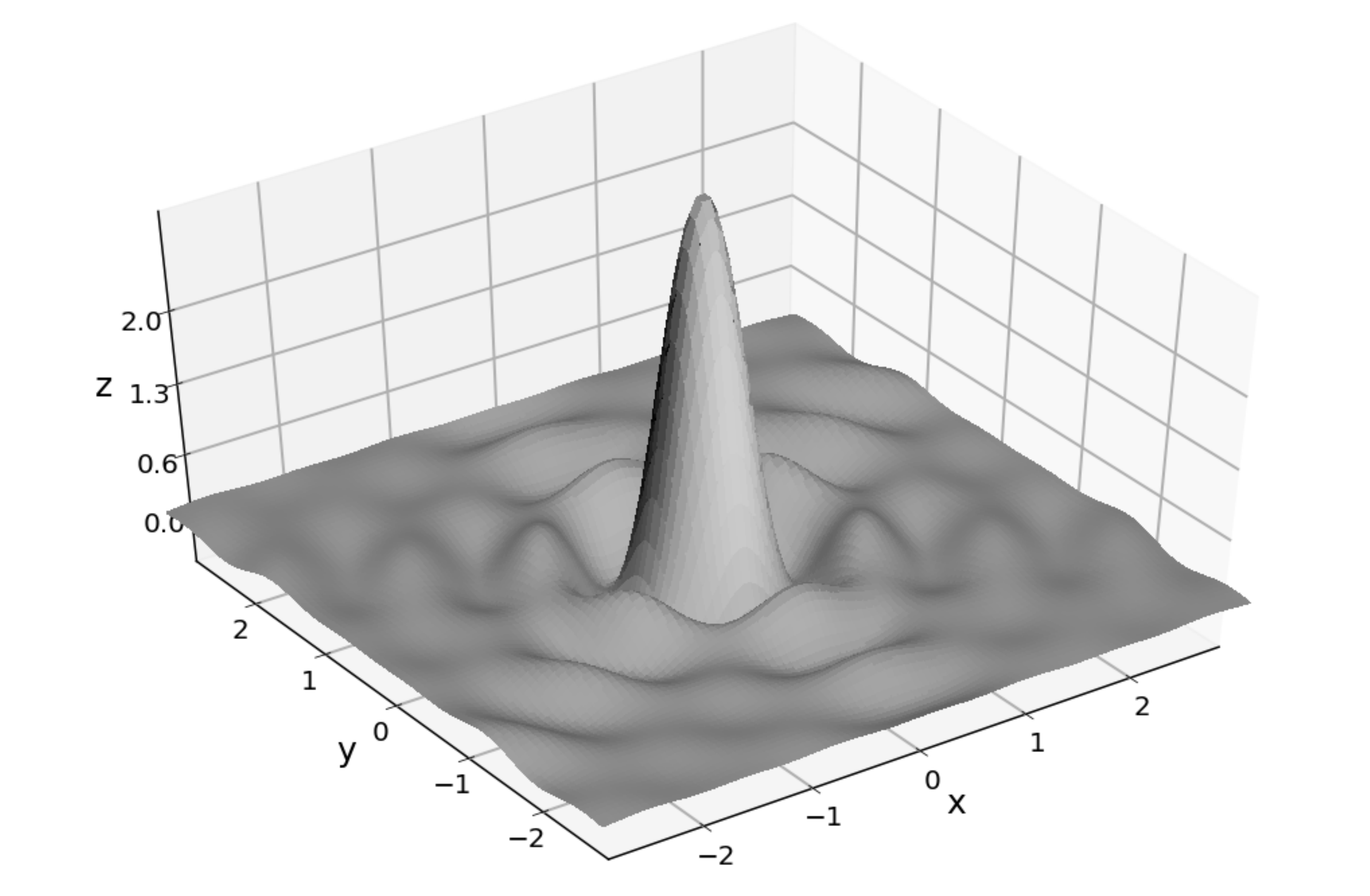}
\end{center}
\caption{A graph of the Fourier transform $\hat 1_H(x,y)$ of the symmetric hexagon $H$
 in equation \eqref{eq:sinc-hexagonal}.  The graph suggests that the largest peak occurs at the origin, which indeed is always the case, by Exercise \ref{strictly less than the FT at zero}.   }    
\label{HexagonPic}
\end{figure}

\bigskip
 \begin{example}\label{FT of a symmetric hexagon}
 \rm{
 Suppose we have a hexagon $H$ that is symmetric about the origin.  We know that its Fourier transform is real-valued, by Lemma \ref{symmetric iff FT is real}, so that here it makes sense to form a 
$3$-dimensional graph of the points
 $(x, y, \hat 1_H(x,y))$, as in Figure \ref{HexagonPic}.  
 

To be concrete, let's define a (parametrized) hexagon $H$ with the following vertices: 
\[
v_1 = \Big(\frac{2c}{\sqrt{3}}, 0\Big),\ \ 
v_2 = \Big(\frac{c}{\sqrt{3}}, c\Big),\ \ 
v_3 = \Big(\frac{-c}{\sqrt{3}}, c\Big),\ \ 
v_4 = -v_1,\ \ v_5 = -v_2,\ \ v_6 = -v_3,
\]
for each fixed parameter $c>0$.
Just for fun, our hexagon is scaled so that it has an inscribed circle of radius $c$, which may be useful in future applications.

We first use Brion's Theorem  \ref{brion, continuous form} to compute the Fourier Transforms of the $6$ vertex tangent cones of $H$.
For $v_1$, the two rays defining $K_{v_1}$ are $w_1 := v_2 - v_1 = (-\frac{c}{\sqrt{3}}, c)$ and $w_2 := v_6-v_1 = (-\frac{c}{\sqrt{3}}, -c)$, so the Fourier Transform of $K_{v_1}$ is:
\[
\hat 1_{K_{v_1}}(z) 
= \frac{e^{-2\pi i \frac{2c}{\sqrt{3}}z_1}}{(-2\pi i)^2} \frac{\frac{2c^2}{\sqrt{3}} } {(-\frac{c}{\sqrt{3}}z_1 + c z_2)(-\frac{c}{\sqrt{3}}z_1 - c z_2)} 
= \frac{2\sqrt 3 }{(2\pi )^2} \frac{ e^{ -\frac{4\pi i  c }{\sqrt{3}} z_1}}
{(-z_1 + \sqrt 3 z_2)(z_1 + \sqrt 3 z_2)}.
\]

For $v_2$, the two rays are $w_1 := v_3 - v_2 = (-\frac{2c}{\sqrt{3}}, 0)$ and $w_2 := v_1 - v_2 = (\frac{c}{\sqrt{3}}, -c)$, giving us:
\[
  \hat 1_{K_{v_2}}(z) 
= \frac{e^{-2\pi i (\frac{c}{\sqrt{3}}z_1+cz_2)}}{(-2\pi i)^2} \frac{\frac{2c^2}{\sqrt{3}} } {\frac{-2c}{\sqrt{3}}z_1(\frac{c}{\sqrt{3}}z_1 - c z_2)} 
= \frac{\sqrt 3}{(2\pi )^2} \frac{ e^{-2\pi c i(\frac{1}{\sqrt{3}} z_1+ z_2)}}
{z_1( z_1 - \sqrt 3 z_2)}.
\]

For $v_3$, the two rays are $w_1 := v_4 - v_3 = (-\frac{c}{\sqrt{3}}, -c)$ and $w_2 := v_2 - v_3 = (\frac{2c}{\sqrt{3}}, 0)$, giving us:
\[
\hat 1_{K_{v_3}}(z) 
= \frac{e^{-2\pi i (-\frac{c}{\sqrt{3}}z_1+ cz_2)}}{(-2\pi i)^2} \frac{\frac{c}{\sqrt{3}} } {(-\frac{c}{\sqrt{3}}z_1 - c z_2) \frac{2c}{\sqrt{3}}z_1}
= \frac{\sqrt 3}{(2\pi )^2} \frac{ e^{-2\pi c i(-\frac{1}{\sqrt{3}}z_1+ z_2)}}
{z_1(z_1  +\sqrt 3z_2)}.
\]

By the inherent symmetry of our hexagon $H$, the computations for the other tangent cones are just $\hat 1_{K_{-v}}(z) =   1_{K_{v}}(-z)$, so we have:
\begin{equation}\label{eq:sinc-hexagonal}
\begin{aligned}
&\hat 1_H(z_1, z_2) := \int_{H}e^{-2\pi i \langle \xi, z \rangle} d\xi \\
&=   \hat 1_{K_{v_1}}(z) +  \hat 1_{K_{v_1}}(-z)     
+  \hat 1_{K_{v_2}}(z) +   \hat 1_{K_{v_2}}(-z)    
+  \hat 1_{K_{v_3}}(z) +   \hat 1_{K_{v_3}}(-z)    \\
&=
\frac{\sqrt{3}}{2\pi^2}\left(
\frac{2 \cos(\frac{4\pi c}{\sqrt{3}} z_1)}{(-z_1 + \sqrt 3 z_2)(z_1 + \sqrt 3 z_2)} 
+ 
\frac{\cos\big(\frac{2\pi c}{\sqrt{3}} z_1 + 2\pi c z_2 \big)}{z_1(z_1 - \sqrt 3 z_2)}
+ 
\frac{\cos\big(\frac{2\pi c}{\sqrt{3}} z_1 - 2\pi c z_2 \big)}{z_1(z_1 + \sqrt 3 z_2)} \right).
\end{aligned}
\end{equation}
}
\hfill $\square$
\end{example}

More generally, suppose we are given the vertices of a polygon $\P\subset \R^2$, so that
\[
\P := \conv\{ v_1, \dots, v_N\},
\]
the convex hull of its vertices.  Brion's Theorem  \ref{brion, continuous form} again gives us
a closed form in terms of the Fourier transforms of its vertex tangent cones. To this end, we first compute the FT of each of its vertex tangent cones:
\[
\hat 1_{\K_{v_k}}(\xi) =  \frac{ e^{-2\pi i \langle v,  \xi   \rangle} \det \K_{v_k}}
  {  \langle v_{k+1}-v_k, \xi \rangle  \langle v_{k-1}-v_k, \xi \rangle },
\]
where $\det \K_{v_k} = 
 \begin{pmatrix} |  &    |  \\  
                      v_{k+1}-v_k &    v_{k-1}-v_k     \\  
                          |  &    |  \\ 
        \end{pmatrix} $ is the invertible $2\times 2$ real matrix whose columns are the edge vectors that are incident with the vertex $v_k$.  For any real convex polygon, 
we may order its vertices  in a counter-clockwise orientation $v_1, \dots, v_N$, with the definition $v_{N+1}:= v_1$.

With this notation, we have proved the following expression for the Fourier transform of a polygon, as a direct consequence of Brion's Theorem  \ref{brion, continuous form}.

\begin{cor}\label{FT of a polygon!}
Given any convex polygon $\P\subset \R^2$, its Fourier transform has the formula:
\[
   \int_\P e^{-2\pi i \langle u,  \xi  \rangle} \, du  =
   -\frac{1}{4\pi^2}           
  \sum_{ k=1}^N
     \frac{ e^{-2\pi i \langle v,  \xi   \rangle} \det \K_{v_k}}
  {  \langle v_{k+1}-v_k, \xi \rangle  \langle v_{k-1}-v_k, \xi \rangle },
\]
for all $\xi \in \R^2$ such that $\xi$ is not orthogonal to any edge of $\P$.
\hfill $\square$
\end{cor}


\bigskip
\section{Each polytope has its moments} \index{moments}
\index{simple polytope}
\index{volume}

The following somewhat surprising formula for the volume of a simple polytope gives us a very rapid algorithm
for computing volumes of simple polytopes.   We note that it is an NP-hard problem \cite{Barany} 
 to compute volumes of general polytopes, without fixing the dimension. 
Nevertheless, there are various other families of polytopes whose volumes possess tractable algorithms.

\bigskip
\begin{thm}[Lawrence \cite{LawrenceVolume}]
   \label{volume of a simple polytope}
Suppose $\P \subset \R^d$ is a simple,  $d$-dimensional polytope.
For a vertex tangent cone $\K_v$ of $\P$, fix a set of edges of the cone, say 
$w_1(v), w_2(v), \dots, w_d(v) \in \R^d$. 
Then    \index{volume of a simple polytope}
\begin{equation}    \label{formula for the volume of simple polytope}
\vol \P = \frac{ (-1)^d }{ d! } \sum_{ v \text{ {\rm a vertex of }} \P } 
\frac{ {\langle   v, z \rangle}^{d} \det \K_v }
{ \prod_{ k=1 }^d \langle w_k(v), z \rangle}
\end{equation}
for all $z\in \C^d$ such that $z$ does not belong to the 
finite union of hyperplanes that are orthogonal to any edge of $\P$. 

More generally, for any integer $k \ge 0$, we have the {\bf moment formulas}:
\index{moment formulas}

\begin{equation}
\int_\P {\langle   x, z \rangle}^{k} dx = \frac{ (-1)^d k! }{ (k +d)! } \sum_{ v \text{ {\rm a vertex of }} \P } 
\frac{ {\langle   v, z \rangle}^{k+d} 
\det \K_v    }{ \prod_{ m=1 }^d \langle w_m(v), z \rangle} \, .
\end{equation}
\end{thm}
\bigskip
\begin{proof}
We begin with Brion's identity \eqref{transform formula for a simple polytope}, and we substitute $z := t z_0$ for a fixed
complex vector  $z_0\in \C^d$, and any positive real value of $t$:
\begin{equation*}\label{transform formula for a simple polytope 2}
  \int_\P e^{-2\pi i \langle u,  z_0\rangle t} \, du  =  
 \left(     \frac{1}{2\pi i}    \right)^d        \sum_{ v \text{ {\rm a vertex of }} \P } 
  \frac{ e^{-2\pi i \langle v,  z_0\rangle t} \det \K_v}
  { t^d \prod_{ m=1 }^d \langle w_m(v), z_0 \rangle }.
\end{equation*}
Now we expand both sides in their Taylor series about $t=0$.   The left-hand-side becomes:
\begin{align*}
  \int_\P \sum_{k=0}^\infty  \frac{1}{k!} \left(-2\pi i \langle u,  z_0\rangle t \right)^k    \, du  
  &=  
 \left(     \frac{1}{2\pi i}    \right)^d        \sum_{ v \text{ {\rm a vertex of }} \P } 
  \frac{ \sum_{j=0}^\infty  \frac{1}{j!} \left(-2\pi i \langle v,  z_0\rangle t \right)^j  \det \K_v}
  { t^d \prod_{ m=1 }^d \langle w_m(v), z_0 \rangle }   
  \end{align*}
 Integrating term-by-term on the left-hand-side, we get:
  \begin{align}\label{main identity for moments}
 \sum_{k=0}^\infty  \frac{t^k}{k!} (-2\pi i )^k   \int_\P   \langle u,  z_0\rangle^k   \, du  
  &=  
    \left(     \frac{1}{2\pi i}    \right)^d       \sum_{ v \text{ {\rm a vertex of }} \P } 
 \frac{ \det \K_v}{\prod_{ m=1 }^d \langle w_m(v), z_0 \rangle}
   \sum_{j=0}^\infty  \frac{t^{j-d}}{j!} (-2\pi i )^j  {\langle v,  z_0\rangle}^j.  
  \end{align}
Comparing the coefficients of $t^k$ on both sides, we have:
\begin{equation*}
 \frac{(-2\pi i )^k}{k!}   \int_\P   \langle u,  z_0\rangle^k   \, du  = 
  \left( \frac{1}{2\pi i}  \right)^d    \sum_{ v \text{ {\rm a vertex of }} \P } 
   \frac{ \det \K_v}{\prod_{ m=1 }^d \langle w_m(v), z_0 \rangle}
   \frac{1}{(k+d)!} (-2\pi i )^{k+d}  {\langle v,  z_0 \rangle}^{k+d}.
\end{equation*}
Simplifying, we arrive at the moment formulas, for each $k\geq 0$:
\begin{equation*}
   \int_\P   \langle u,  z_0\rangle^k   \, du  =  (-1)^d \frac{k!}{(k+d)!} 
  \sum_{ v \text{ {\rm a vertex of }} \P } 
   \frac{   {\langle v,  z_0 \rangle}^{k+d}  \det \K_v    }{\prod_{ m=1 }^d \langle w_m(v), z_0 \rangle }.
\end{equation*}
In particular, when $k=0$, we get the volume formula \eqref{formula for the volume of simple polytope}.
\end{proof}

The following interesting identities are also consequences of the proof above, and were discovered by Brion and Vergne \cite{brionvergne}.
\begin{cor}
\label{additional Brion-Vergne type identities}
Suppose $\P \subset \R^d$ is a simple,  $d$-dimensional polytope.  For each $0\leq j \leq d-1$, we have:
\begin{equation}
 \sum_{ v \text{ {\rm a vertex of }} \P } 
   \frac{   {\langle v,  z \rangle}^{j}  \det \K_v    }{\prod_{ m=1 }^d \langle w_m(v), z \rangle }=0,
\end{equation}
for all $z\in \C^d$  such that $z$ does not belong to the 
finite union of hyperplanes that are orthogonal to any edge of $\P$. 
\end{cor}
\begin{proof}
We may go back to \eqref{main identity for moments}, and stare at that Laurent series in $t$. We notice 
that the singular part in $t$, which contains exactly the terms with $t^j$ with $j= 0, -1, \dots, -(d-1)$, must vanish
because the left-hand side of that identity does not contain any singular terms in $t$.
\end{proof}

\bigskip
\section{The zero set of the Fourier transform}

Now we know enough to derive some new results, regarding the real zero set of the Fourier transform: $Z_\R(\P):= \{ x\in \R^d \mid \hat 1_\P(\xi) = 0 \}$.

\begin{cor}
\label{THE NEW RESULT, vanishing of the FT}
Let $\P\subset \R^d$ be a $d$-dimensional integer polytope $\P \subset \R^d$.  Then:
\begin{equation}
\hat 1_\P(\xi) = 0,
\end{equation}
 for each integer point $\xi  \in \Z^d$
that does not belong to the 
finite union of hyperplanes orthogonal to an edge of $\P$,
\end{cor}
\begin{proof}
Any integer polytope may be triangulated into integer simplices (not necessarily unimodular simplices), and we'll call such a collection of simplices $T$.
We consider any of these integer simplices, say $\Delta$.  By Brion's Theorem \ref{brion, continuous form}, 
we have
\begin{equation}
\hat 1_\Delta(\xi) =  
 \left(     \frac{1}{2\pi i}    \right)^d        \sum_{ v \text{ {\rm a vertex of }} \P } 
  \frac{ e^{-2\pi i \langle v,  \xi   \rangle} \det \K_v}
  { \prod_{ k=1 }^d \langle w_k(v), \xi \rangle } 
\end{equation}
for all $\xi \in \R^d$ such that the denominators on the right-hand side do not vanish. Here we've used the fact that $\Delta$ is a simple polytope.
In particular, for an integer point $\xi \in \Z^d$ (which does not belong to the 
finite union of hyperplanes that are orthogonal to any edge of $\P$), we have $\langle \xi, v \rangle \in\Z$, using the assumption that the vertices $v\in \P$ are integer points. Consequently, we have
\begin{align}
\hat 1_\Delta(\xi) &=  
 \left(     \frac{1}{2\pi i}    \right)^d        \sum_{ v \text{ {\rm a vertex of }} \P } 
  \frac{  \det \K_v}
  { \prod_{ k=1 }^d \langle w_k(v), \xi \rangle } 
= 0,
\end{align}
by Corollary \ref{additional Brion-Vergne type identities} (the $j=0$ case).  Summing all of the Fourier transforms of the simplices that belong to our triangulation (and ignoring their boundaries because the FT vanishes there), we arrive at 
\begin{equation}
\hat 1_\P(\xi) = \sum_{\Delta \in T} \hat 1_\Delta(\xi) = 0.
\end{equation}
\end{proof}

\begin{figure}[htb]
\begin{center}
\includegraphics[totalheight=2.4in]{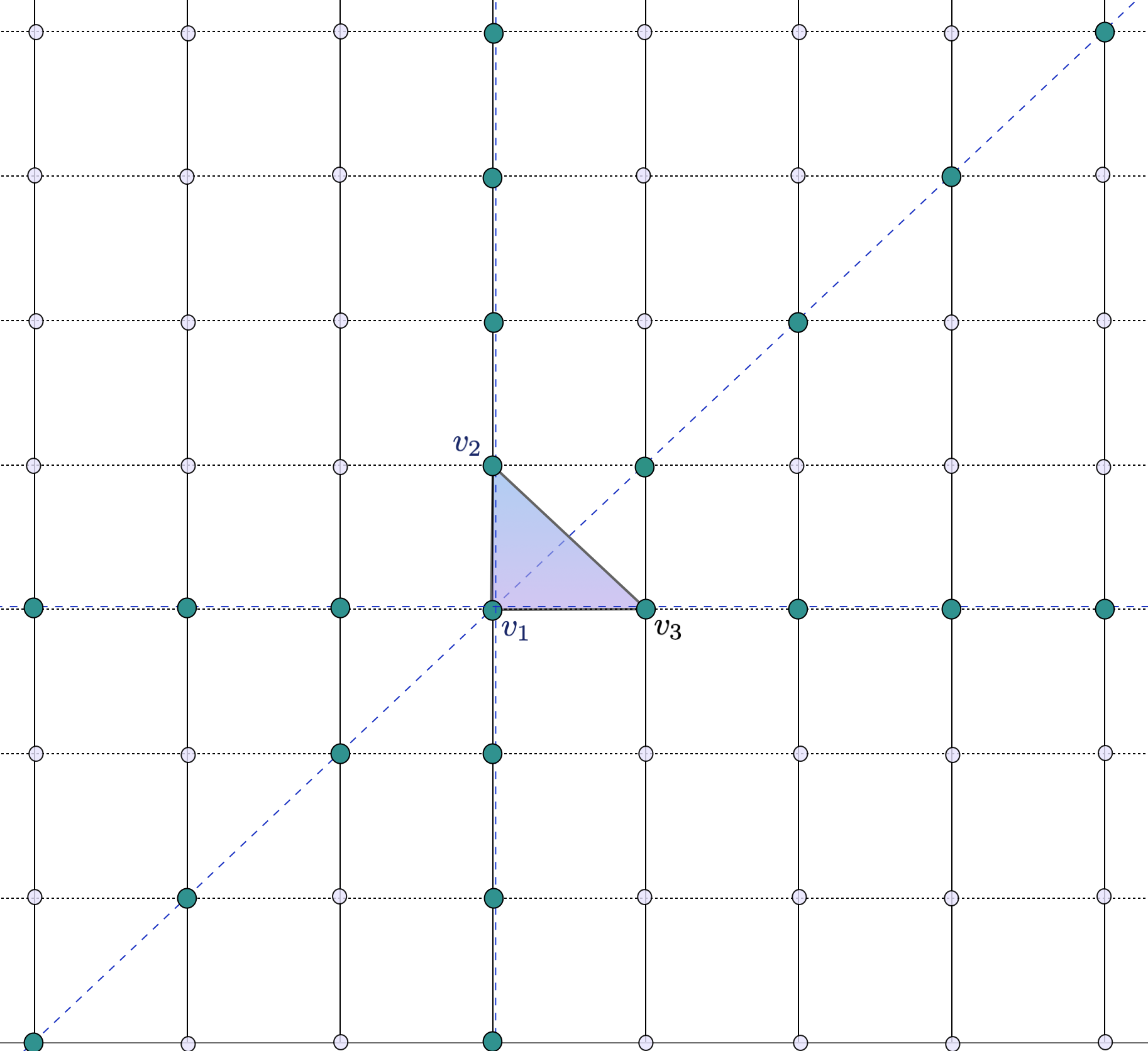}
\end{center}
\caption{
Here the green integer points depict the special integer frequencies of the Fourier transform 
$\hat 1_{\Delta}(\xi)$, where $\Delta$ is the standard triangle.  
These special integer frequency vectors are orthogonal to the three sides of $\Delta$.  
For all other nonzero integer points (the generic integer frequencies), 
$\hat 1_{\Delta}(\xi) =0$, according to Corollary \ref{THE NEW RESULT, vanishing of the FT}.
}    
\label{Special Frequencies}
\end{figure}

Given a polytope $\P \subset \R^d$, we call a vector $\xi \in \R^d$ a {\bf generic frequency} (relative to $\P$) \index{generic frequencies}
 if $\xi$ is not orthogonal to any edge of $\P$ (and hence not orthogonal to any other face of $\P$).  
 All other $\xi \in \R^d$ are orthogonal to some edge of $\P$, and are called {\bf special frequencies} (see also Section \ref{sec: generic vs. special frequencies}).

 We recall 
 the hyperplane arrangement defined by the finite collection of hyperplanes orthogonal to any edge of 
 $\P$:
\begin{align}\label{the hyperplane arrangement H}
\mathcal H &:= \{  \xi \in \R^d \mid \langle \xi, E \rangle = 0, \text{  for any edge $E$ of $\P$}   \} \\
&=\{\text{special frequencies}\},
\end{align}
which came up naturally in the general formula for the Fourier transform of a polytope (Theorem \ref{brion2}).
It's clear from the definitions above that the special frequencies are $\xi \in \mathcal H$, and the generic frequencies are 
$\xi \notin \mathcal H$. 

Hence Corollary \ref{THE NEW RESULT, vanishing of the FT} may be restated as follows. 
For an integer polytope $\P\subset \R^d$, we have:
\begin{equation}
\hat 1_\P(\xi) = 0,
\end{equation}
for all generic frequencies $\xi \in \Z^d$. In other words, $\{\text{generic integer frequencies}\}  \subset Z_\R(\P)$.
It's natural to wonder if the latter vanishing of the transform is sufficient to identify a polytope among the collection of all convex bodies, as follows. 
\begin{conjecture}\label{conj: recovering a polytope from knowledge of its null set}
Suppose we know that $\P$ is a convex body in $\R^d$.  Suppose further that we are given the data:
 \[
 \hat 1_\P(\xi) = 0, \text{ for all } \xi\in \Z^d \setminus \mathcal H, 
 \]
 where $\mathcal H$ is some finite collection of hyperplanes passing through the origin.  
Then: 
\begin{enumerate}[(a)]
\item $\P$ is a polytope.
\item Moreover, $\P$ is an integer polytope, and $\mathcal H$ is precisely the collection of hyperplanes that are orthogonal 
to all of the edges of $\P$.
\end{enumerate}
\end{conjecture}
Although Conjecture \ref{conj: recovering a polytope from knowledge of its null set} appears here for the first time, it highlights the importance of the zero set of the Fourier transform.

Kobayashi \cite{Kobayashi1} asked the following question.
\begin{question}\label{Kobayashi question}
Does the zero set $Z_\C(\P):= \{ \zeta \in \C^d \mid \hat 1_\P(\zeta) = 0 \}$
 determine the convex body $\P$, among all convex bodies, up to translations?
\end{question}
We've already seen, in Theorem \ref{zero set of the FT of a polytope} (Kolountzakis' vanishing criterion), that 
if we only assume 
that $\left\{ \Z^d \setminus \{0\}\right\}  \subset Z_\R(\P)$, then even this very sparse assumption on the zero set is already equivalent 
to $\P$  multi-tiling Euclidean space. 

We finish this section by reinterpreting Brion's
Theorem  \ref{brion2}, using the  meromorphic continuing the real vector $\xi$, and the hyperplane arrangement $\mathcal H$ of \eqref{the hyperplane arrangement H}. 
We may extend the Fourier transform $ \hat  1_{K_v}(z)$ of a rational cone  to all of $\C^d$,  using the bold-face notation  $ \hat {\bf 1}_{K_v}(z)$,
 by using the fact that it is a rational function in several variables:
\begin{equation} \label{meromorphic continuation identity}
  \hat {\bf 1}_{K_v}(z) := 
\frac{e^{-2\pi i \langle v, z \rangle} }{(2\pi i)^d}     
\sum_{j=1}^{M(v)}   \frac{\det \K_j(v)  }{\prod_{k=1}^d  \langle  w_{j, k}(v),   z    \rangle},
\end{equation}
for all $z \notin \mathcal H$.
With this notation we may rewrite Brion's Theorem \ref{brion2}, for any real polytope $\P$,  as follows:
\begin{equation}
  \int_\P e^{-2\pi i \langle u,  z\rangle} \, du  =  
\sum_{v\in V}     \hat {\bf 1}_{K_v}(z),
\end{equation}
valid for all $z \notin \mathcal H$.


\bigskip
\section*{Notes} \label{Notes.chapter.Brion}
 \begin{enumerate}[(a)]
 \item There is a lot more literature about the zero set of the Fourier transform of a convex body.  For more information, we refer the reader to \cite{Bianchi1}, \cite{Kobayashi1}, \cite{Kobayashi2}
 \cite{Kolountzakis2}.
 
\item  There is a large literature devoted to triangulations of cones, polytopes, and general point-sets, 
and the reader is invited to consult the excellent and encyclopedic book on triangulations, by 
Jes\'us de Loera, J\"org Rambau, and Francisco Santos \cite{DRS}.

\item  The notion of a {\bf random polytope} has a large literature as well, and although we do not go into this topic here, one classic survey paper is by Imre B\'ar\'any \cite{Barany}.

\item  The attempt to extend Ehrhart theory to non-rational polytopes, whose vertices have some irrational coordinates, is ongoing.   
The pioneering papers of Burton Randol \cite{Randol1} \cite{Randol2} extended integer point counting to algebraic polytopes, meaning that their vertices are allowed to have coordinates that are algebraic numbers.  Recently, a growing number of papers are considering all real dilates of a rational polytope, which is still rather close to the Ehrhart theory of rational polytopes. 
 
 In this direction, it is natural to ask how much more of the geometry of a given polytope $\P$ can be captured by
  counting integer points in all of its positive real dilates.    Suppose we translate a 
$d$-dimensional integer polytope $\P \subset \R^d$ by an integer vector $n \in \Z^d$.  
The standard Ehrhart theory gives us an invariance principle, namely the equality of the Ehrhart polynomials for $\P$ and $\P + n$:
\[
L_{\P+n}(t) = L_\P(t), 
\]
for all  \emph{integer} dilates $t>0$.  

However, when we allow $t$ to be a positive \emph{real} number, then it is in general {\bf false} that 
\[
L_{\P+n}(t) = L_\P(t)  \text{  for all }  t > 0. 
\]
 In fact, these two Ehrhart functions are so different in general, that by the breakthrough of Tiago Royer \cite{Tiago1}, 
  it's even possible to uniquely reconstruct the polytope $\P$ if we know all the counting quasi-polynomials  
$L_{\P+n}(t)$, for all integer translates $n \in \Z^d$.   In other words, the work of \cite{Tiago1} shows that for two 
rational polytopes $\P, Q \subset \R^d$, we have:
\[
L_{\P+n}(t) = L_{Q+n}(t) \text{ for all } n\in\Z^d \text{ and all } t>0 \iff \P=Q.
\]
 It is rather astounding that just by counting integer points in sufficiently many translates of $\P$, we may completely reconstruct the whole polytope $\P$ uniquely.    Royer further demonstrated  \cite{Tiago2} that such an idea also works if we replace a polytope by any symmetric convex body.   It is now natural to try to prove the following extended question.
\begin{question}
\rm{Suppose we are given polytopes $\P, Q \subset \R^d$.   Can we always find a finite subset 
$S\subset \Z^d$ (which may depend on $\P$ and Q)
 such that 
 \[
 L_{\P+n}(t) = L_{Q+n}(t)  \text{ for all }  n \in S, \text{ and all }    t >0 \  \iff  \ \P = Q? 
 \]
 }
\end{question}
\end{enumerate}

\newpage

\section*{Exercises}
\addcontentsline{toc}{section}{Exercises}
\markright{Exercises}

\begin{quote}
``It is better to solve one problem five different ways, than to solve five problems one way.''

-- George P\'olya
\end{quote}

\medskip
\begin{prob}  $\clubsuit$ \label{independent of edge vectors}
\rm{
Although $\det \K_v$ depends on the choice of the length of each edge of $\K_v$, show that the ratio  
$ \frac{ |\det \K_v|  }{\prod_{k=1}^d  \langle  w_k(v), z  \rangle}$
 remains invariant if we replace each edge $w_k(v)$ of a simplicial cone by a constant
positive multiple of it, say $\alpha_k w_k(v)$ with $\alpha_k>0$.

(Here $z$ is any generic complex vector, meaning that  $\langle  w_k(v), z  \rangle \not=0$ for all $1\leq k\leq d$ ).
}
\end{prob}

\medskip
\begin{prob}
Consider the regular hexagon $\P \subset \R^2$, whose vertices are the $6$'th roots of unity.
 \begin{enumerate}[(a)]
\item  Compute the area of $\P$ using Theorem   \ref{volume of a simple polytope}. 
\item Compute all of the moments of $\P$, as in Theorem    \ref{volume of a simple polytope}. 
\end{enumerate}
\end{prob}

\medskip
\begin{prob}
Compute the Fourier transform of the triangle $\Delta$ whose vertices are given by
\[
(1, 0), (0, 1), (-c, -c),
\]
where $c>0$.
\end{prob}

\medskip
\begin{prob} \label{translating a cone}
 $\clubsuit$ 
Prove Corollary \ref{transform of a translated cone} for a simplicial cone  $\K_v$, whose apex is $v$, by 
translating a cone whose vertex is at the origin, to get:
\[
{\hat 1}_{\K_v}(z) :=   \int_{\K_v}   e^{-2\pi i \langle u,  z\rangle} \, du  =   \frac{1}{(2\pi i)^d}   
 \frac{ e^{-2\pi i \langle v, z \rangle}  \det \K_v  }{\prod_{k=1}^d  \langle  w_k, z     \rangle}.
\]
\end{prob}

\medskip
\begin{prob} 
Using some of the idea in Lemma \ref{LimitDim.d}, prove the following:
\begin{enumerate}[(a)]
\item  For all nonzero $\alpha \in \R$,
\[
\lim_{\varepsilon \rightarrow 0} \int_0^\infty \cos(\alpha x) \, e^{-\varepsilon |x|^2} dx = 0.
\]
\item  For all nonzero $\alpha \in \R$,
\[
\lim_{\varepsilon \rightarrow 0}    \int_0^\infty \sin(\alpha x) \, e^{-\varepsilon |x|^2} dx = \frac{1}{\alpha}.
\]
\end{enumerate}
\end{prob}

\medskip
\begin{prob} \label{Pyramid over a square}
Consider the following $3$-dimensional polytope $\P$, whose vertices are as follows:
\[
\left\{  (0, 0, 0),   \   (1, 0, 0), \  (0, 1, 0), \  (1,  1,  0),  \  (0, 0, 1) \right\}.  
\]
``a pyramid over a square".   Compute its Fourier-Laplace transform $\hat 1_\P(z)$.
\end{prob}

 \begin{figure}[htb]
\begin{center}
\includegraphics[totalheight=2in]{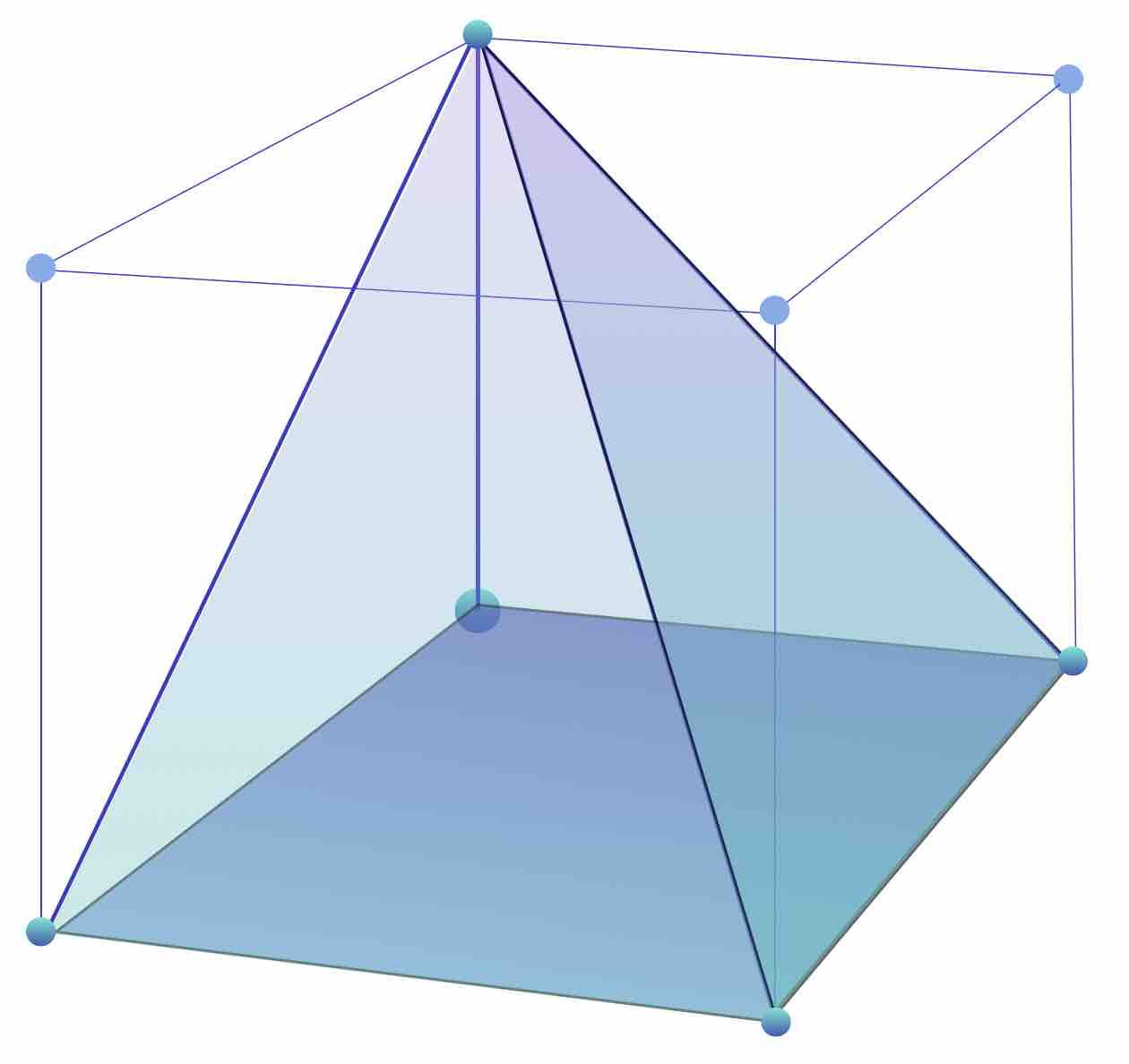}
\end{center}
\caption{The pyramid over a square, in Exercise \ref{Pyramid over a square}   } 
\label{pyramid over a square}
\end{figure}

\medskip
\begin{prob}
\rm{
We recall that the $3$-dimensional cross-polytope (also called an octahedron)
was defined by $\Diamond:=\left\{ \left( x_1, x_2, x_3 \right) \in \R^d \mid
 \, \left| x_1 \right| + \left| x_2 \right|  + \left| x_3 \right| \leq 1 \right\}$. 
Compute the Fourier-Laplace transform of $\Diamond$ by using 
Theorem \ref{brion2}.
\index{cross-polytope}

Notes. \ Here not all of the tangent cones are simplicial cones, so 
we may triangulate each vertex tangent cones into simplicial cones, or you may try your own methods.
}
\end{prob}

\medskip
\begin{prob}[hard-ish]   \label{FT of a Dodecahedron}
\rm{
Here we will find the Fourier transform of a dodecahedron $\P$, centered at the origin. 
Suppose we fix the following $20$ vertices of $\P$:
\[
\left\{       \left(\pm 1,\  \pm 1, \ \pm 1 \right),   \   
\left(0, \ \pm \phi, \ \pm \frac{1}{\phi} \right),  \   \left(\pm \frac{1}{\phi},  \ 0, \  \pm \phi \right),    \ 
  \left( \pm \phi,  \ \pm \frac{1}{\phi}, \  0 \right)
\right\},
\]
where $\phi:= \frac{1+\sqrt{5}}{2}$.  It turns out that $\P$ is a simple polytope.  Compute its Fourier-Laplace transform using Theorem \ref{brion, continuous form}.

Notes.    All of the vertices of $\P$ given here can easily be seen to lie on a sphere $S$ of radius $\sqrt{3}$, and 
this is a regular embedding of the dodecahedron.  It is also true (though a more difficult fact) that these $20$ points maximize the volume of any polytope whose $20$ vertices lie on the surface of this sphere $S$. 
}
\end{prob}

\begin{figure}[htb]
  \begin{center}
\includegraphics[totalheight=2.3in]{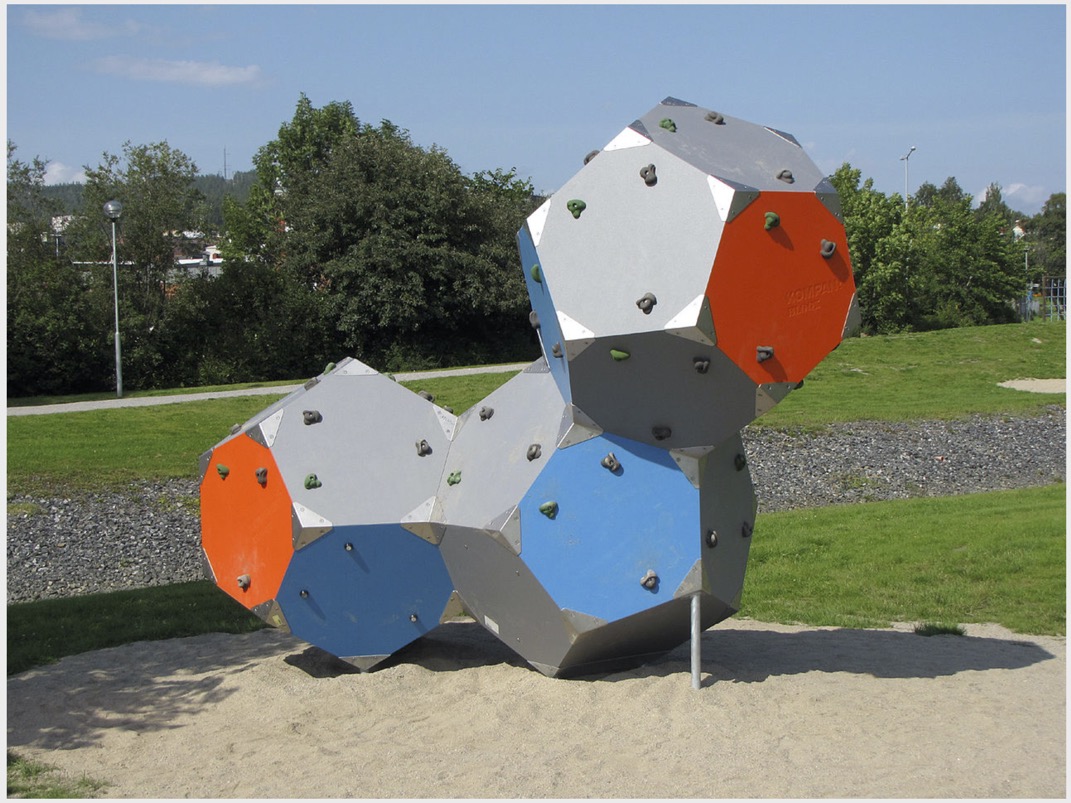}
   \end{center}
\caption{A climbing wall in Sweden, made up of Dodecahedrons, showing one of their real-life applications}  \label{Dodecahedron}
\end{figure}

\medskip
\begin{prob}
Define the $3$-dimensional polytope 
$\P := \rm{ conv }\{ (0,0,0), (1,0,0), (0,1,0), (0,0, 1), (a,b,c)  \}$, where we fix real the positive real numbers $a,b,c$.   Compute $\hat 1_\P(z)$, by computing the Fourier-Laplace transforms of its tangent cones.  

(Note.  Here, not all of the tangent cones are simplicial cones).
\end{prob}

\medskip
\begin{prob} \label{Pyramid over a cube}
This exercise extends Exercise \ref{Pyramid over a square} to $\R^d$, as follows.
Consider the $d$-dimensional polytope $\P$, called a ``pyramid over a cube", defined by the convex
hull of the unit cube $[0,1]^{d-1} \subset \R^{d-1}$, with the point $(0, 0, \dots, 0, 1) \in \R^{d}$.
Compute its Fourier-Laplace transform $\hat 1_\P(z)$. 
\end{prob}

\medskip
\begin{prob}  $\clubsuit$     \label{cone equivalence}
Show the following two conditions are equivalent:
 \begin{enumerate}[(a)]
\item  A cone $\K$ has an apex at the origin.
\item  $\K$ is a cone that enjoys the property $\lambda \K = \K$, for all $\lambda >0$. 
\end{enumerate}
\end{prob}

\medskip
\begin{prob}  $\clubsuit$ \label{simplicial implies pointed}
Suppose we are given a $d$-dimensional simplicial cone $\K \subset \R^d$ (so be definition $\K$ has exactly $d$ edges).
Show that $\K$ must be pointed. 
\end{prob}

\medskip
\begin{prob}  $\clubsuit$  \label{Exercise.tangent cone of a vertex}
Show that for any polytope $\P\subset \R^d$, a vertex tangent cone $\K_v$ never contains a whole line. 
\end{prob}

\medskip
\begin{prob}  $\clubsuit$     \label{pointed cone equivalence}
\index{cone, pointed}
\rm{
 Show that if $\K$ is a cone with an apex $v$ (not necessarily a unique apex), the
  following conditions are equivalent:
 \begin{enumerate}[(a)]
 \item  $\K$ is a pointed cone.
\item There exists a hyperplane $H$ such that  $H\cap \K = v$.  
\item The translated cone $C:= \K-v$, with apex at the origin, enjoys  $C \cap (-C) = \{0\}$.
 \item  $\K$ has a unique apex.
 \item  $\K$ does not contain an entire line. 
\end{enumerate}
}
\end{prob}

\medskip
\begin{prob} $\clubsuit$ \label{simplicial AND simple}
Show that the only polytopes that are both simple and simplicial are either
simplices, or $2$-dimensional polygons.
\end{prob}

\blue{
For problems   \ref{elementary properties of polarity} - \ref{dual after translation},  
we recall - for the sake of disambiguation with the polar set below - that for any cone $\K\subset \R^d$, its dual cone $\K^*$ 
was defined (recalling \eqref{dual cone definition}) by 
\[
\K^* := \{  y \in \R^d \mid \langle y, u \rangle < 0 \text{ for all } u\in \K  \}.
\]
}

\medskip
\begin{prob} $\clubsuit$ \label{duality of dual cone}
Show that if we have reverse inclusions for dual cones.  Namely:
\[
\K_1 \subset \K_2  \iff  \K_2^* \subset \K_1^*.
\]
\end{prob}

\medskip
\begin{prob}  \label{dual cones and Minkowski sums}  
\index{Minkowski sum}
Show that if we take the Minkowski sum $K_1 + K_2$ of two cones $\K_1, \K_2 \subset \R^d$, then 
polarity interacts with Minkowski sums in the following pleasant way:
\[
\left(\K_1 +\K_2\right)^* = \K_1^* \cap \K_2^*.
\]
\end{prob}

\blue{
For problems   \ref{elementary properties of polarity} - \ref{dual after translation},  
given any set $S\subset \R^d$, we define its {\bf polar set}  by
\[
S^o := \left \{y \in \R^d \mid \langle y, z \rangle \leq 1 \text{ for all } z \in S \right \},
\]
which may sometimes be unbounded.  Note that this definition is consistent with our previous definition
of the polar polytope in \eqref{polar polytope, definition}.  We also note here the distinction between a polar set
and the dual cone.  Throughout, we've defined duality only for cones, to disambiguate between the two notions.
}

\medskip
\begin{prob}      \label{elementary properties of polarity}
\rm{
Here are some elementary properties of polarity, applies to general sets.
 \begin{enumerate}[(a)]
 \item   If $A \subset B \subset \R^d$, show that $B^o \subset A^o$.
 \item  For $A \subset \R^d$, show that $A \subset \left(A^o\right)^o$.
\item If $A_1, \dots, A_k \subset \R^d$, show that $\left(\cup_{j=1}^k A_j \right)^o = \cap_{j=1}^k A_j^o$.
\item For $A\subset \R^d$, we have $A = A^o \iff A=B_r$, a ball of radius $r$, centered at the origin.
\end{enumerate}
}
\end{prob}

\medskip
\begin{prob} $\clubsuit$ \label{dual after translation}
 For any fixed translation vector $v\in \R^d$, prove that 
\begin{equation}
\left( S + v \right)^o = \left\{  \frac{1}{1+\langle v, y \rangle} \, y \mid y \in S^o \right\}.
\end{equation}
\end{prob}

\medskip
\begin{prob}  \label{polytope from pentagons}
Suppose we try to construct a polytope $\P \subset \R^3$ all of whose facets are pentagons (not necessarily regular). 
Show that
$
 F\geq 12,
 $
 where $F$ is the number of facets of $\P$.
\end{prob}

\medskip
\begin{prob}  \label{Euler equivalent to Brianchon-Gram}  $\clubsuit$ 
\rm{
\begin{enumerate}[(a)]
\item  Show that the Brianchon-Gram relations \eqref{BG} imply
 the Euler-Poincare relation for the face-numbers  of a convex polytope $\P$:
\begin{equation}
f_0 - f_1 + f_2 - \cdots + (-1)^{d-1} f_{d-1} +  (-1)^{d} f_{d}= 1,
\end{equation}
where  $f_k$ is the number of faces of $\P$ of dimension $k$.
\item \label{b} (hard) \  Conversely, given a $d$-dimensional polytope $\P\subset \R^d$, show that the Euler-Poincare relation above implies the Brianchon-Gram relations:
 \[
 1_\P(x) = \sum_{\F \subset \P} (-1)^{dim \F} 1_{\K_F}(x),
 \]
  for all $x\in \R^d$.
\end{enumerate}

Notes.  Interestingly, even though the above two conditions are equivalent, condition \ref{b} is often more useful in practice, because we have a free variable $x$, over which we may sum or integrate.
}
\end{prob}

\medskip
\begin{prob}  
Find a $2$ dimensional integer polygon $\P$ such that, for any integer point $n\in \Z^2$ there exists $t>0$ with
\[
L_{\P+n}(t) \not= L_\P(t).
\]
\end{prob}

Notes.  \ When $t$ is restricted to be a positive integer, it is of course true that 
$L_{\P+n}(t) = L_\P(t)$.  The point here is that when working with all positive dilates, the differences between integer polytopes becomes more pronounced.



 \chapter{\blue{What is an angle in higher dimensions?}}
 \label{Angle polynomial}
 \index{angle polynomial} \index{solid angle}

\begin{quote} 
``Everyone else would climb a peak by looking for a path somewhere in the mountain. Nash would climb another mountain altogether and from that distant peak would shine a searchlight back onto the first peak.''

-- Donald Newman
\end{quote} 
 
 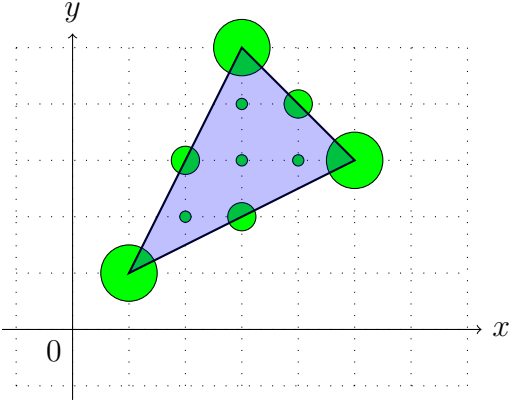
\begin{figure}[!h] 
		\centering
		\begin{tikzpicture}[scale=.75]
			\draw (0,0) node[below left] {$0$};
			\draw[loosely dotted] (-1,-1) grid (7,5);
			\draw[->] (-1.25,0) -- (7.25,0) node[right] {$x$};
			\draw[->] (0,-1.25) -- (0,5.25) node[above] {$y$};
			\draw[fill = green] (3,4) circle (.1cm);
			\draw[fill = green] (4,3) circle (.1cm);
			\draw[fill = green] (3,3) circle (.1cm);
			\draw[fill = green] (3,2) circle (.25cm);
			\draw[fill = green] (2,2) circle (.1cm);
			\draw[fill = green] (2,3) circle (.25cm);
			\draw[fill = green] (4,4) circle (.25cm);
			\draw[fill = green] (5,3) circle (.5cm);
			\draw[fill = green] (3,5) circle (.5cm);
			\draw[fill = green] (1,1) circle (.5cm);
			\draw[thick] (1,1) -- (5,3) -- (3,5) -- cycle;
			\filldraw[nearly transparent, blue] (1,1) -- (5,3) -- (3,5) -- cycle;
			
		\end{tikzpicture}  
		\caption{A discrete volume of the triange $\P$, called the angle polynomial of $\P$.  
		Here we sum   local angle weights, relative to $\P$, at all integer points.} 
		\label{triangle solid angle sum}
	\end{figure}

\section{Intuition}

There are infinitely many ways to discretize the classical notion of volume, and here we offer a second path, using `local solid angles'.   Given a rational polytope $\P$, we will place small spheres at all integer points in $\Z^d$, and compute the proportion of the local intersection of each small sphere with $\P$.  This discrete, finite sum, gives us a new method of discretizing the volume of a polytope, and it turns out to be a more symmetric way of doing so.   To go forward, we first discuss how to extend the usual notion of `angle' to higher dimensions, and then use Poisson summation again to pursue the fine detail of this new discrete volume.

\bigskip
\section{Defining an angle in higher dimensions}  \label{Chapter.solid.angles}

The question of how an angle in two dimensions extends to higher dimensions is a basic one
in discrete geometry.   A natural way to extend the notion of an angle is to consider
a cone $\K \subset \R^d$, place a sphere centered at the apex of $\K$, and then compute the proportion 
of the sphere that intersects $\K$.    This intuition is captured more rigorously by the following integral:
\begin{equation} \label{integral def. of solid angle}
\omega_\K = \int_\K e^{-\pi \| x \|^2} dx.
\end{equation}
called the  {\bf solid angle of the cone} $\K$.  \index{solid angle}
The literature has other synonyms for solid angles, arising in different fields, 
 including the {\bf volumetric moduli} \cite{GourionSeeger},  
 \index{volumetric moduli}
 and the {\bf volume of a spherical polytope} 
 \index{volume of a spherical polytope} 
 \cite{BeckRobins}, \cite{DesarioRobins}, \cite{RicardoNhatSinai}.

\begin{figure}[htb]
 \begin{center}
\includegraphics[totalheight=2.3in]{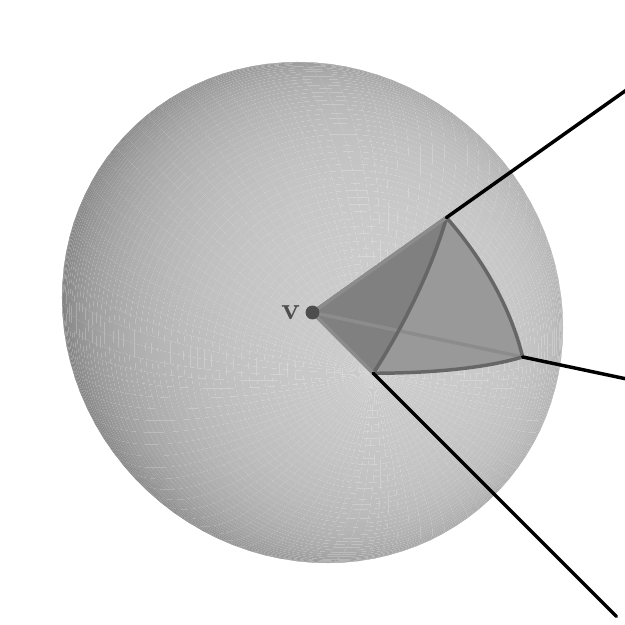}
\end{center}
\caption{A solid angle in $\R^3$ - note the equivalence with the area of the geodesic triangle on the sphere.} 
\end{figure}

We can easily show that the latter definition of a solid angle
is equivalent to the volume of a spherical polytope, 
\index{volume of a spherical polytope}
 using polar coordinates in $\R^d$, as follows.    We denote the 
unit sphere  by
$S^{d-1}:= \{  x\in \R^d \mid \| x\| = 1 \}$.  Then using the fact that the Gaussians give a probability distribution, namely $\int_{\R^d} e^{-\pi ||x||^2} dx = 1$  (which we know by Exercise \ref{Gaussian1}), we have 
\begin{align}
        \omega_\K &= \frac{\int_\K e^{-\pi\|x\|^2}dx}{\int_{\R^d} e^{-\pi\|x\|^2}dx} \label{second equality}
            \ = \ \frac{\int_0^{\infty} e^{-\pi r^2} r^{d-1} dr \int_{S^{d-1} \cap \K} d\theta}{\int_0^{\infty} 
            e^{-\pi r^2} r^{d-1} dr \int_{S^{d-1}} d\theta} \\
           &= \ \frac{\int_{S^{d-1} \cap \K} d\theta}{\int_{S^{d-1}} d\theta} \\ \label{normalized spherical volume}
            &= \ \frac{\vol_{d-1} \left(\K \cap S^{d-1}\right)}{\vol_{d-1} \left( S^{d-1} \right)},
\end{align}
where $\vol_{d-1}$ denotes the volume measure on the surface of the $(d-1)$-dimensional sphere $S^{d-1}$.  
We may think of \eqref{normalized spherical volume}  as the {\bf normalized volume} of a spherical polytope defined by the intersection of the cone $\K$
with the unit sphere.    Thus for any cone $\K\subset \R^d$, we have
\[
0 \leq \omega_\K \leq 1.
\]
 We used polar coordinates in the second equality \eqref{second equality} above:
$x= (r, \theta)$, with $r \geq 0,  \ \theta \in \S^{d-1}$.  
The Jacobian in the change of variables is $dx = r^{d-1} dr d\theta$.

We note that when $\K= \R^d$, so that here the cone is \emph{all} of Euclidean space, 
the integral \eqref{integral def. of solid angle}
 becomes
\[
\int_{\R^d} e^{-\pi ||x||^2} dx = 1, 
\]
by Exercise \ref{Gaussian1}.   This computation confirms that we do indeed have the proper normalization 
with $\omega_\K = 1$ if and only if  $\K = \R^d$.

\begin{example} 
\rm{
 If $\K\subset \R^d$ is a half-space, then
$\omega_\K = \frac{1}{2}$.   If $\K:= \R_{\geq 0}^d$, the positive orthant, then 
\begin{align*}
\omega_\K &=  \int_{\R^d_{\ge 0}} e^{-\pi ||x||^2} dx =  
\left(  \int_{\R_{\ge 0}} e^{-\pi u^2} du  \right)^d  =  \frac{1}{2^d}.
\end{align*}
So in the plane, the positive quadrant takes up $\frac{1}{4}$ of the whole plane. In $\R^3$, the positive octant takes up $\frac{1}{8}$ of the whole space, etc.
\hfill $\square$
}
\end{example}

We might wonder: ``Do we really need to use Gaussians to define these solid angles?'' The clear answer is ``no'', as the following example shows.   But one reason to favor Gaussians over other radially symmetric functions is that they behave beautifully under convolutions and Fourier transforms, as we'll see later in Lemma \ref{basic connection for solid angle}.
\begin{example}
\rm{
Let $\K\subset \R^d$ be a $d$-dimensional polyhedral cone, and fix $s > \frac{d}{2}$. Then we have:
\[
\omega_K = 
\frac{ 
\pi^{\tfrac{d}{2}} \Gamma(s) 
}
{
\Gamma(s - \tfrac{d}{2})
}
 \int_K \frac{dx}{  \left( 1+ \|x\|^2 \right)^s}.
\]
The reader may enjoy proving this from scratch, in Exercise \ref{another integral for the solid angle},
(or see Proposition 2.3.2 in \cite{GervasioSantos}).
}
\hfill $\square$
\end{example}

\bigskip
\section{Local solid angles for a polytope, and Gaussian smoothing}

Here we want to define solid angles relative to a fixed polytope.   
So given any polytope $\P \subset \R^d$, we fix any point $x \in \R^d$ and define a local solid angle relative to $\P$ as follows. 
 The normalized {\bf solid angle} \index{solid angle} fraction that a  
$d$-dimensional polytope $\P$ subtends  at any point $x \in \R^d$ is defined  by 
\begin{equation}\label{def. of solid angle}
\om_\P(x)=\lim_{\eps \to 0}     \frac{\vol(S^{d-1}(x,\eps) \cap \P)}{\vol\left(S^{d-1}(x,\eps)\right)}.
\end{equation}

Here, $\omega_{\P}(x)$ measures the fraction of a small $(d-1)$-dimensional sphere 
$S^{d-1}(x,\varepsilon)$
centered at $x$, that intersects the polytope $\P$.  
 We will use the standard notation for the interior of a convex body, namely $\interior(\P)$, and for the boundary of a convex body, namely $\partial \P$.  As a side-note, we mention  that balls and spheres can be used interchangeably in this definition, meaning that the fractional weight given by  \eqref{def. of solid angle}
 is the same using either method (see Exercise  \ref{two definitions for a solid angle}).

It follows from the definition of a solid angle that $0 \leq     \om_\P(x) \leq 1$,  for all $x \in \R^d$, and that
\begin{equation}\label{def: solid angle at x}
\om_\P(x) =
  \begin{cases}
      1 & \text{if }       x \in \interior(\P)        \\
     0  & \text{if }        x  \notin \P.
    \end{cases}  
 \end{equation}
But when $x \in \partial \P$, we have $0< \om_\P(x) <1$.   For example, 	if $x$ lies on a 
codimension-two face of $\P$,  then  $\om_\P(x)$ is the fractional dihedral angle 
subtended by $\P$ at $x$.	 
	
To define one type of discrete volume for any polytope $\P\subset \R^d$,  we fix a 
positive integer $t$, and define the finite sum 
\begin{equation} \label{anglesum1}
A_\P(t) :=  \sum_{n\in \Z^d} \om_{t\P}(n),
\end{equation}
where $t\P:=\{tx \mid x \in \P\}$ is the $t$'th dilation of the polytope $\P$.   
In other words, $A_\P(1)$ is by definition the discrete volume for $\P$ which is 
obtained by placing at each integer point $n\in \Z^d$ the weight
$\om_{\P}(n)$, and summing all of the weights over all $n \in \Z^d$.

\begin{example}
\rm{
In Figure \ref{triangle solid angle sum}, the solid angle sum of the polygon $\P$ is
 \[
 A_\P(1) = \theta_1 + \theta_2 +  \theta_3 + 3\left( \tfrac{1}{2} \right)+ 4 = 6.   
 \]
 Here  the $\theta_j$'s are the three angles at the vertices of 
$\P$.
}
\hfill $\square$
\end{example}

Using purely combinatorial methods,  Macdonald  showed  that for any integer polytope $\P$, 
and for {\bf positive integer values}  of $t$,
 \begin{equation} \label{solidanglesum2}
A_\P(t) = (\vol \P)  t^d + a_{d-2} t^{d-2} + a_{d-4} t^{d-4} + \cdots +   
    \begin{cases}
      a_1 t & \text{if } d \text{ is odd},\\
      a_2 t^2  & \text{if } d  \text{ is even}.
    \end{cases}
\end{equation} 
We will call $A_\P(t)$ the {\bf angle-polynomial} of $\P$, \index{angle polynomial} 
for integer polytopes $\P$ and positive
integer dilations $t$.  However, when these restrictions are lifted, the sum still captures crucial
geometric information of $\P$, and we will simply call it the (solid) angle-sum of $\P$.

We define the {\bf heat kernel}, \index{heat kernel} 
for each fixed positive $\varepsilon$, by
\begin{equation}
G_{\varepsilon}(x) := \varepsilon^{-\frac{d}{2}}  e^{-\frac{\pi}{\varepsilon}  \| x \|^2},
\end{equation}
for all $x \in \R^d$.  By Exercises \ref{Gaussian1} and \ref{Gaussian2}, we know that
$\int_{\R^d} G_{\varepsilon}(x)dx = 1$ for each fixed $\varepsilon$, and that 
\begin{equation}  \label{Fourier transform of the Gaussian, formal}
\hat G_{\varepsilon}(\xi) = e^{-\varepsilon \pi \| \xi \|^2}.
\end{equation}

The convolution of the indicator function $1_\P$ by the heat kernel $G_{\varepsilon}$ will be called the 
{\bf Gaussian smoothing} \index{Gaussian smoothing}  
of $1_\P$:
\begin{align} \label{Gaussian smoothing}
(1_\P * G_{\varepsilon})(x) &:=  \int_{\R^d}   1_\P(y) G_{\varepsilon} (x-y) dy 
					      = \int_{\P}               G_{\varepsilon} (y-x) dy \\
			&=    \varepsilon^{-\frac{d}{2}}     \int_{\P}   e^{-\frac{\pi }{\varepsilon}   \| y-x \|^2} dy,
\end{align}
a $C^{\infty}$ function of $x\in \R^d$, and in fact a Schwartz function 
(Exercise  \ref{convolution of the indicator function with a Gaussian}).  

The following Lemma provides a first
crucial link between the discrete geometry of a local solid angle and the convolution of $1_\P$ with a 
Gaussian-based approximate identity.

\bigskip
\begin{lem} \label{basic connection for solid angle}
		Let $\P$ be a full-dimensional polytope in $\R^d$. Then for each point $x\in\R^d$, we have
\begin{equation}
 \lim_{\varepsilon \rightarrow  0} (1_\P * G_{\varepsilon})(x) = \omega_P(x). 
 \end{equation}
\end{lem}
	\begin{proof}
		We have
		\begin{align*}
		(1_\P * G_{\varepsilon})(x)  &=   \int_{\P} G_{\varepsilon} (y-x) dy \\
			           &= \int_{u \in P- x} G_{\varepsilon} (u) du \\
			           &= \varepsilon^{\frac{d}{2}}   \int_{\frac{1}{\sqrt{\eps}}(P- x)} G_1(v) dv.
		\end{align*}
		
In the calculation above, we make use of the evenness of $G_{\varepsilon}$ in the second equality. The substitution $v = u/\sqrt{\varepsilon}$ was used in the third equality.  Following those substitutions, 
we change the domain of integration from $P$ to the translated body $P- x$, and then to the dilation of $P-x$ by the factor $\frac{1}{\sqrt{\varepsilon}}$. 
		
Finally, when $\varepsilon$ approaches $0$, 
	$\frac{1}{\sqrt{\varepsilon}}(P- x)$ tends to a cone $K$ with apex at the origin, subtended by $P- x$.
	This  cone $K$ is in fact a translation of the tangent cone  \index{tangent cone}
 of $P$, at $x$.  We therefore arrive at
\[ 
\lim_{\varepsilon \to 0} (1_\P * G_{\varepsilon}) (x) = \int_{K} G_1(v)dv = \om_K(0) = \om_P(x). 
\]
\end{proof}

Putting things together,  the definition \ref{anglesum1} and 
Lemma \ref{basic connection for solid angle}
above tell us that
\begin{equation} 
       A_\P(t) = \sum_{n\in\Z^d}   \omega_{tP}(x) = 
        \sum_{n\in\Z^d}  \lim_{\varepsilon \rightarrow 0}  (1_{t\P} * G_{\varepsilon})(n). 
\end{equation}
We would like to interchange a limit with an infinite sum over a lattice, so that we may use Poisson summation, and although this is subtle in general, it's possible to carry out here, because the summands are rapidly decreasing. 
\begin{lem} \label{interchanging limit and sum, solid angle sum}
Let $\P$ be a full-dimensional polytope in $\R^d$. 
Then
\begin{equation} 
			\sum_{n\in\Z^d}  \lim_{\varepsilon \rightarrow 0}  (1_{t\P} * G_{\varepsilon})(n)
			= \lim_{\varepsilon \rightarrow 0} \sum_{n\in\Z^d} (1_{t\P} * G_{\varepsilon})(n). 
\end{equation}
\hfill $\square$
\end{lem}
(For a proof of Lemma \ref{interchanging limit and sum, solid angle sum} see \cite{RicardoNhatSinai}).

Next, we apply the Poisson summation formula to the Schwartz function \\
$f(x) := (1_\P * G_{\varepsilon})(x)$:
\begin{align} \label{Gaussian smoothing for Angle sum}
A_P(t) &=   \lim_{\varepsilon \rightarrow 0} \sum_{n\in\Z^d} (1_{t\P} * G_{\varepsilon})(n) \\
&= \lim_{\varepsilon \rightarrow 0} \sum_{\xi  \in\Z^d} \hat 1_{t\P}(\xi) \hat G_{\varepsilon}(\xi) \\
&= \lim_{\varepsilon \rightarrow 0} \sum_{\xi  \in\Z^d} \hat 1_{t\P}(\xi)  \ e^{-\varepsilon \pi \| \xi \|^2} \\
&= t^d \lim_{\varepsilon \rightarrow 0} \sum_{\xi  \in\Z^d} \hat 1_{\P}( t\xi)  
     \ e^{-\varepsilon \pi \| \xi \|^2} \\
&=  t^d \ \hat 1_{\P}( 0 )  +  \lim_{\varepsilon \rightarrow 0} t^d \sum_{\xi  \in\Z^d-\{0\}} \hat 1_{\P}( t\xi)  
     \ e^{-\varepsilon \pi \| \xi \|^2} \\
&=  t^d (\vol \P)  +  \lim_{\varepsilon \rightarrow 0}t^d \sum_{\xi  \in\Z^d-\{0\}} \hat 1_{\P}( t\xi)  
     \ e^{-\varepsilon \pi \| \xi \|^2},     
 \label{last line of Poisson summation}
\end{align}
where we used the fact that Fourier transforms interact nicely with dilations of the domain:
\[
\hat 1_{t\P}(\xi) = \int_{t\P} e^{-2\pi i \langle \xi, x \rangle} dx = 
 t^d  \int_{\P} e^{-2\pi i \langle \xi, ty \rangle} dy   = t^d  \int_{\P} e^{-2\pi i \langle t \xi, y \rangle} dy=
  t^d \hat 1_{\P}(t\xi).
\]
We also used the simple change of variable $x =t y$, with $y \in \P$, implying that $dx = t^d dy$,
as well as the Fourier transform formula for the heat kernel \eqref{Fourier transform of the Gaussian, formal}.

So far, we've proved the following.
\begin{lem}\label{half-way to the angle polynomial, via Poisson summation}
Given a real polytope $\P\subset \R^d$, its angle polynomial has the expression:
\begin{equation} \label{phase 2 of angle polynomial}
	A_\P(t) = 
	 t^d (\vol \P)  + t^d  \lim_{\varepsilon \rightarrow 0} 
	  \sum_{n\in\Z^d-\{0\}} ( \hat 1_{\P}(t\xi) * G_{\varepsilon})(n).
\end{equation}
\hfill $\square$
\end{lem}

Lemma \ref{half-way to the angle polynomial, via Poisson summation}
suggests  a polynomial-like behavior for the angle polynomial $A_\P(t)$.  Once we learn how to compute the Fourier transform of a polytope using Stokes' formula, in  
Chapter~\ref{Stokes' formula and transforms}, we will be able to continue this computation in \eqref{phase 2 of angle polynomial}, and many more geometric facts will unfold.

\bigskip
\section{$1$-dimensional polytopes}
\label{sec:1-dim polytopes}

Although this toy case is straightforward, we'll still encounter some interesting formulae and ideas. 
We may use our knowledge of the Fourier transform of a $1$-dimensional polytope $\P$, in
the right-hand-side of \eqref{phase 2 of angle polynomial}, namely a closed  interval $P:= [a, b]$.


Let's compute the angle polynomial of the $1$-dimensional polytope $\P:= [a,b]$, with
$a,~b~\in~\R$.  We will use our knowledge of the $1$-dimensional Fourier transform of an interval,
from Exercise \ref{transform.of.interval.a.to.b}, to compute:
\begin{align} 
A_P(t) &=  (b-a) t + \lim_{\varepsilon \rightarrow 0} \sum_{\xi  \in\Z-\{0\}} \hat 1_{\P}( t\xi)  \ e^{-\varepsilon \pi \xi^2} \\
&= (b-a) t  + \lim_{\varepsilon \rightarrow 0} \sum_{\xi  \in\Z-\{0\}}
\left(       \frac{e^{-2\pi i t\xi b}    -  e^{-2\pi i t\xi a}    }{-2\pi i \xi}       \right)
   e^{-\varepsilon \pi \xi^2}  \\  \label{strange limit1}
&= (b-a) t  
- \lim_{\varepsilon \rightarrow 0} \sum_{\xi  \in\Z-\{0\}}
      \frac{e^{-2\pi i tb\xi    -\varepsilon \pi \xi^2}}{2\pi i \xi}     
+  \lim_{\varepsilon \rightarrow 0} \sum_{\xi  \in\Z-\{0\}}
          \frac{e^{-2\pi i ta\xi   -\varepsilon \pi \xi^2}}{2\pi i \xi}  
\end{align}           
Throughout the latter computation, all series converge absolutely (and quite rapidly) due to the existence of the Gaussian damping factor~$e^{-\varepsilon \pi \xi^2}$.

Let's see what happens when we specialize the vertices $a$ or $b$ - perhaps we can solve for one of these new limits? 
So we set $a= 0, b \in \R \setminus \Z, t \in \R_{>0}$.    In this special case, 
one of the two series in \eqref{strange limit1} becomes:
\[
\sum_{\xi  \in\Z-\{0\}}
          \frac{e^{-2\pi i ta\xi } e^{  -\varepsilon \pi \xi^2}}{-2\pi i \xi} = 
  \sum_{\xi  \in\Z-\{0\}}
          \frac{e^{  -\varepsilon \pi \xi^2}}{-2\pi i \xi} = 0,
\]
because the summand is an odd function of $\xi$.  Since the solid angle at an integer vertex of an interval equals $\tfrac{1}{2}$, we already know by direct computation that in this case
\[
A_{[0, b]}(t) =
\begin{cases}
\frac{1}{2} +  \lfloor bt  \rfloor   & \text{if } bt \notin \Z  \\
 \lfloor bt  \rfloor    &  \text{if } bt \in \Z
\end{cases}
=
\frac{1}{2} - \frac{1}{2}1_{\Z}(bt)  +  \lfloor bt  \rfloor,
\]
for all $t>0$.  Here we've used a handy definition for the indicator function of the integers:
\[
1_{\Z}(x):= 
\begin{cases}
1 & \text{if } x \in \Z  \\
0 &  \text{if } x \notin \Z.
\end{cases}
\]
Solving \eqref{strange limit1}  for the other limit, we get:
\[        
\frac{1}{2}- \frac{1}{2}1_{\Z}(bt)  +  \lfloor bt  \rfloor 
=   b t  + 
\lim_{\varepsilon \rightarrow 0} \sum_{\xi  \in\Z-\{0\}}
\left(       \frac{e^{-2\pi i t\xi b}     }{-2\pi i \xi}       \right)
   e^{-\varepsilon \pi \xi^2} 
\]
After relabelling $bt := x\in \R$, we've just proved the following.

\begin{lem} \label{rigorous approach for P_1(x)}
For any $x\in \R$, we have
\[
\frac{1}{2\pi i} \lim_{\varepsilon \rightarrow0}  \sum_{\xi \in \Z - \{0\}}
 \frac{e^{-2\pi i x\xi    -\varepsilon \pi \xi^2}}{\xi}  
 =   x -  \lfloor x  \rfloor  - \frac{1}{2} +\frac{1}{2} 1_{\Z}(x).
 \]
 \hfill $\square$
\end{lem}
Now we can bootstrap our information from Lemma \ref{rigorous approach for P_1(x)} 
 by plugging its result back into equation \eqref{strange limit1}:
\begin{align} 
A_P(t)  &= (b-a) t  
- \left( bt - \lfloor bt \rfloor - \frac{1}{2} +\frac{1}{2} 1_{\Z}(bt) \right)
+  \left( at - \lfloor at \rfloor - \frac{1}{2} +\frac{1}{2} 1_{\Z}(at) \right) \\  
\label{recovering 1-dim'l solid angle poly}
&=  \lfloor bt \rfloor    -\frac{1}{2} 1_{\Z}(bt)  - \lfloor at \rfloor +\frac{1}{2} 1_{\Z}(at),
\end{align}            
and we've arrived at the angle polynomial for any $1$-dimensional polytope $\P:= [a, b]$, where $a, b \in \R$.  Of course, \eqref{recovering 1-dim'l solid angle poly}
 is  easy to check directly from the definition of the angle polynomial for an interval, but note that we also recovered a non-trivial limit in the process, namely Lemma \ref{rigorous approach for P_1(x)}.

\bigskip
\section{Pick's formula and Nosarzewska's inequality}

A polygon $\P$ is called an {\bf integer polygon} if all of its  vertices belong to the integer lattice $\Z^2$.
There is a wonderful relationship, discovered by George Pick in 1899, between the area of $\P$, 
and the number of integer points contained in $\P$ and on its boundary.
\begin{thm}[Pick's formula, 1899]
Let $\P$ be an integer polygon.   Then 
\begin{equation}\label{Pick's formula}
{\rm Area} \, \P  =  I + \frac{1}{2} B -1,
\end{equation}
where $I$ is the number of interior integer points in $\P$, and B is the number of boundary integer points in $\P$. 
\hfill $\square$
\end{thm}
There is an equivalent formulation of Pick's theorem in terms of local solid angle weights at each integer point. 
\begin{thm}[Pick's formula, reformulated with angle weights]
\label{Pick's formula, solid angle version}
Let $\P$ be an integer polygon.   Then 
\[
\sum_{n \in \Z^2}  \omega_\P(n) = {\rm Area} \, \P ,
\]
where  $\omega_\P(n)$ is the $2$-dimensional angle defined in \eqref{def: solid angle at x}.
\hfill $\square$
\end{thm}

\begin{figure}[!h]
		\centering
		\begin{tikzpicture}[scale=0.45]
			\draw (0,0) node[below left] {$0$};
			\draw[loosely dotted] (-1,-1) grid (7,5);
			\draw[->] (-1.25,0) -- (7.25,0) node[right] {$x$};
			\draw[->] (0,-1.25) -- (0,5.25) node[above] {$y$};
			\draw[fill = green] (3,4) circle (.1cm);
			\draw[fill = green] (4,3) circle (.1cm);
			\draw[fill = green] (3,3) circle (.1cm);
			\draw[fill = green] (3,2) circle (.25cm);
			\draw[fill = green] (2,2) circle (.1cm);
			\draw[fill = green] (2,3) circle (.25cm);
			\draw[fill = green] (4,4) circle (.25cm);
			\draw[fill = green] (5,3) circle (.5cm);
			\draw[fill = green] (3,5) circle (.5cm);
			\draw[fill = green] (1,1) circle (.5cm);
			\draw[thick] (1,1) -- (5,3) -- (3,5) -- cycle;
			\filldraw[nearly transparent, blue] (1,1) -- (5,3) -- (3,5) -- cycle;
			\draw (3,-2) node {$P_1$};
			
			\draw (9,2) node[scale = 2] {$\cup$};
			
			\draw (0+12,0) node[below left] {$0$};
			\draw[loosely dotted] (-1+12,-1) grid (7+12,5);
			\draw[->] (-1.25+12,0) -- (7.25+12,0) node[right] {$x$};
			\draw[->] (0+12,-1.25) -- (0+12,5.25) node[above] {$y$};
			\draw[fill = green] (5+12,2) circle (.1cm);
			\draw[fill = green] (4+12,2) circle (.1cm);
			\draw[fill = green] (3+12,2) circle (.25cm);
			\draw[fill = green] (5+12,3) circle (.5cm);
			\draw[fill = green] (6+12,2) circle (.5cm);
			\draw[fill = green] (1+12,1) circle (.5cm);
			\draw[thick] (1+12,1) -- (6+12,2) -- (5+12,3) -- cycle;
			\filldraw[semitransparent, blue] (1+12,1) -- (5+12,3) -- (6+12,2) -- cycle;
			\draw (3+12,-2) node {$P_2$};

			\draw (9+12,2) node[scale = 2] {$=$};

			\draw (0+24,0) node[below left] {$0$};
			\draw[loosely dotted] (-1+24,-1) grid (7+24,5);
			\draw[->] (-1.25+24,0) -- (7.25+24,0) node[right] {$x$};
			\draw[->] (0+24,-1.25) -- (0+24,5.25) node[above] {$y$};
			\draw[fill = green] (3+24,4) circle (.1cm);
			\draw[fill = green] (4+24,3) circle (.1cm);
			\draw[fill = green] (3+24,3) circle (.1cm);
			\draw[fill = green] (5+24,2) circle (.1cm);
			\draw[fill = green] (4+24,2) circle (.1cm);
			\draw[fill = green] (3+24,2) circle (.25cm);
			\draw[fill = green] (2+24,2) circle (.1cm);
			\draw[fill = green] (2+24,3) circle (.25cm);
			\draw[fill = green] (4+24,4) circle (.25cm);
			\draw[fill = green] (5+24,3) circle (.5cm);
			\draw[fill = green] (6+24,2) circle (.5cm);
			\draw[fill = green] (3+24,5) circle (.5cm);
			\draw[fill = green] (1+24,1) circle (.5cm);
			\draw[thick] (1+24,1) -- (6+24,2) -- (3+24,5) -- cycle;
			\draw[thick] (1+24,1) -- (5+24,3);
			\filldraw[nearly transparent, blue] (1+24,1) -- (5+24,3) -- (3+24,5) -- cycle;
			\filldraw[semitransparent, blue] (1+24,1) -- (5+24,3) -- (6+24,2) -- cycle;
			\draw (3+24,-2) node {$P_1 \cup P_2$};
			
		\end{tikzpicture}
		\caption{Additive property of the angle polynomial}
	\end{figure}
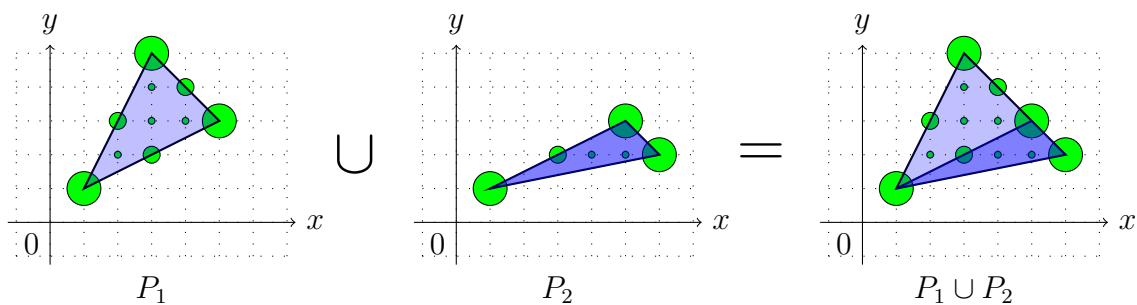

Pick's formula, here formulated as Theorem \ref{Pick's formula, solid angle version}, possesses a 
natural extension to higher dimensions.  Namely, in Theorem \ref{cs.facets} of 
Chapter \ref{Stokes' formula and transforms}, 
we extend Pick's formula to higher dimensions, with a detailed proof that invokes Stokes' theorem.

If we want to work with more general convex regions than polygons, 
 there is a related and  beautiful inequality, discovered in $1948$ by Maria Nosarzewska \cite{Nosarzewska},
  for any $2$-dimensional body.
\begin{thm}[Nosarzewska]\label{thm:Nosarzewska's inequality}
\index{Nosarzewska's inequality}
For a convex body $\P\subset \R^2$,  whose perimeter has length $S$, 
we have:
\begin{equation}\label{Nosarzewska inequality}
{\rm Area}\, \P  - \frac{1}{2} S <  \left|\P\cap \Z^2\right|  \leq 
 {\rm Area}\, \P +\frac{1}{2} S+1.
\end{equation}
\hfill $\square$
\end{thm}
\begin{proof}
To prove the upper bound in \eqref{Nosarzewska inequality}, we'll work with the 
 convex hull $\P_0$ of the interior integer points of $\P$.  Because $\P_0$ is an integer polygon, we may apply Pick's theorem to it.  We let 
 $I_0$ be the number of interior integer points of $\P_0$, $B_0$ be the number of boundary integer points of $\P_0$, and $S_0$ be the perimeter of $\P_0$, so that:
\begin{align*}
 {\rm Area}\, \P +\frac{1}{2} S +1   &\geq   {\rm Area}\, \P_0 +\frac{1}{2} S_0 +1 \\
 &\geq {\rm Area}\, \P_0 +\frac{1}{2} B_0 +1 \\
 & = I_0 + B_0 
 := \left|\P_0\cap \Z^2\right| = \left|\P\cap \Z^2\right|.
\end{align*}
The second inequality above uses the fact that each integer line segment in the plane has length at least $1$, so that the perimeter of an integer polygon must be greater than or equal to the number of integer points on it.  The equality  ${\rm Area}\, \P_0 +\frac{1}{2} B_0 +1 = I_0 + B_0$ above is
 true by Pick's formula \eqref{Pick's formula}.   
For the lower bound, we refer the reader to \cite{Nosarzewska}.
\end{proof}

For any convex body $\P\subset \R^2$, 
Nosarzewska's inequality (Theorem \ref{thm:Nosarzewska's inequality}) is a refinement of Jarnik's inequality 
$\left | \rm{Area}\, \P -  \left|\P\cap \Z^2\right|  \right | < S$.
In $1972$, Bokowski, Hadwiger, and Wills \cite{BokowskiHadwigerWills} 
extended the lower bound in \eqref{Nosarzewska inequality}  to all higher dimensions:
\[
\left | \P \cap \Z^d \right | > \vol \P - \frac{1}{2} S(\P),
\]
where $S(\P)$ is the surface area of the convex body $\P\subset \R^d$.   We might wonder if 
the upper bound of Nosarzewska's inequality \eqref{Nosarzewska inequality} also extends to higher dimensions directly. But a simple counter-example is an $\varepsilon$-neighborhood of a line segment along the first coordinate axis, for example, containing a fixed number of integer points, but whose surface area and volume are both arbitrarily small.  

To discover a result that circumvents the latter counter-example, one might look for an extra assumption on the linear independence of integer points contained in $\P$.  Such a result was given
 by Henk and Wills \cite{HenkWills}, as follows. 
\begin{thm}
Let $\P\subset \R^3$ be a body that contains 3 linearly independent integer points.   Then:
\begin{equation}
\left | \P \cap \Z^d \right | < \vol \P + 2 S(\P).
\end{equation}
\hfill $\square$
\end{thm}

\bigskip
\section{The Gram relations for solid angles}
\index{Gram relations}

\begin{question}[Rhetorical]
\label{sum of angles of a triangle, question}
When we were kids, we learned that the sum of the angles of a triangle equals $\pi$ radians. How does this theorem extend to higher dimensional polytopes?
\end{question} 

We describe the extension here, mainly due to Gram (but has a colorful history).   First, for each face
$F$ of a polytope $\P \subset \R^d$, we define the {\bf solid angle of} $F$, 
\index{solid angle of a face}
as follows.  Fix any
$x_0 \in \interior F$, and let 
\[
\omega_F := \omega_P(x_0). 
\]
We notice  that this definition is independent of $x_0$, as long as we restrict $x_0$ to the relative interior of $F$. 
\begin{example}
If $\P$ is the $d$-dimensional cube $[0, 1]^d$, then each of its facets $F$ has
$\omega_F = \frac{1}{2}$.   Moreover, it is a fact that for the cube, a face $F$ of dimension $k$ has a solid angle 
\[
\omega_F =\frac{1}{2^{d-k}},
\] 
for each $0 \leq k \leq d-1$ (Exercise \ref{solid angle of a face of the cube}).  In particular 
a vertex $v$ of this cube, having dimension $0$,  has solid angle $\omega_v = \frac{1}{2^d}$.
\hfill $\square$
\end{example}

Luckily, Question \ref{sum of angles of a triangle, question} has a beautifully simple answer, as follows. 

\begin{thm}[Gram relations] \label{Gram relations}  \index{Gram relations}
Given any $d$-dimensional polytope $P\subset \R^d$, we have
\[
\sum_{F\subset \P} (-1)^{\dim F} \omega_F = 0.
\]
\end{thm}
\hfill $\square$

(For a proof of Theorem \ref{Gram relations}, see \cite{BeckRobins}, for example).

\begin{example}
\rm{
Let's see what the Gram relations tell us in the case of a triangle $\Delta$.   
For each edge $E$ of $\Delta$, placing a small sphere at a point in the interior of $E$ means half of it is inside $\Delta$ and half of it is outside of $\Delta$, so that $\omega_E = \frac{1}{2}$.   Next, each vertex of $\Delta$ has a solid angle equal to the usual (normalized) angle $\theta(v)$ at that vertex.  
Finally $\Delta$ itself has a solid angle of $1$, because picking a point $p$ in the interior of $\Delta$, and placing a small sphere centered at $p$, the whole sphere will be contained in $\Delta$.  Putting it all together, the Gram relations read:
\begin{align*}
0 &= \sum_{F\subset \Delta} (-1)^{\dim F} \omega_F \\
&=    (-1)^0 (\theta(v_1) + \theta(v_2) +\theta(v_3) ) 
+ (-1)^1 \left(\frac{1}{2} + \frac{1}{2} +\frac{1}{2}\right) 
+ (-1)^2 \cdot 1 \\
&= \theta(v_1) + \theta(v_2) +\theta(v_3) -\frac{1}{2},
\end{align*}
which looks familiar!   We've retrieved our elementary-school knowledge, namely that the three
angles of a triangle sum to $\pi$ radians.  So the Gram relations really are an extension of this fact.
}
\hfill $\square$
\end{example}

What about $\R^3$?    
\begin{example}
\rm{
Let's see what hidden secrets lie behind the Gram relations for the standard simplex 
$\Delta \subset \R^3$.  
\index{standard simplex}
At the origin $v_0 = 0$, the tangent cone is the positive orthant, 
so that $\omega(v_0) = \frac{1}{8}$.  The other $3$ vertices all ``look alike'', in the sense that 
their tangent cones are all isometric, and hence have the same solid angle $\omega_v$.  
What about the edges?  
In general, it's  a fact that the solid angle of an edge equals the dihedral angle between the planes of
its two bounding facets (Exercise \ref{dihedral angle=solid angle}).  There are two types of edges here, as in the figure.  
For an edge $E$ which lies on the boundary of the skew facet,  we have the dihedral angle
$\cos \phi = \left\langle \frac{1}{\sqrt 3}\icol{ 1\\1\\1}, \icol{ 0\\1\\0} \right\rangle = \frac{1}{\sqrt 3}$, so that
$\omega_E = \phi =  \cos^{-1}\frac{1}{\sqrt 3}$.  It's straightforward that for the other type of edge,
 each of those $3$ edges has a solid angle of $\frac{1}{4}$.   Putting it all together, we see that
 \begin{align*}
0 &= \sum_{F\subset \Delta} (-1)^{\dim F} \omega_F \\
&=    (-1)^0 \left(  \frac{1}{8} + 3\omega_v   \right) 
+ (-1)^1 \left( 3 \frac{1}{4}+   3 \cos^{-1}\frac{1}{\sqrt 3}   \right)
+ (-1)^2  \frac{1}{2} \cdot 4  
+ (-1)^3 \cdot 1.
\end{align*}
Solving for $\omega_v$, we get $\omega_v =  \cos^{-1}\frac{1}{\sqrt 3} -\frac{1}{8}$.
So we were able to compute the solid angle of at a vertex of $\Delta$ in $\R^3$, using the
Gram relations, together with a bit of symmetry.
}
\hfill $\square$
\end{example}

Related to the topics above is the fact that the angle polynomial possesses the following fascinating
 functional equation 
(For a proof of Theorem \ref{Angle polynomial functional equation}, and an extension of it, see \cite{DesarioRobins}).

\begin{thm}[Functional equation for the angle polynomial] \label{Angle polynomial functional equation}  
\index{angle polynomial: functional equation}
Given a $d$-dimensional \\
rational polytope $\P\subset \R^d$, we may extend the domain of $A_{\P}(t)$ to all of $\R$ by using
the expression \ref{phase 2 of angle polynomial}.
It follows that 
\[
A_{\P}(-t) = A_{\P}(t),
\]
for all $t \in \R$.
\hfill $\square$
\end{thm}

\bigskip
\section{Bounds for solid angles} 

Throughout this section we're given a $d$-dimensional, simplicial, pointed cone  $\K$, with apex at the origin, and edge vectors $w_1, \dots, w_d \in \R^d$. 
 Let $M \in GL_{\R}(d)$ be the matrix whose columns are the edge vectors $w_j$.   We'll use the observation  
 that $M$ maps the positive orthant $\R^d_{\geq 0}$ bijectively onto $\K$.


 Gourion and Seeger \cite{GourionSeeger} gave some interesting bounds  
 for the solid angle $\omega_K$, in terms of the singular value decomposition of $M$ 
 (Theorem \ref{upper and lower bounds for the solid angle} below). For the 
  linear algebra definitions and applications of singular values, the reader may consult \cite{BisgardLinearAlgebraBook}.  First, some of the many basic and easy facts that make singular values useful are the following upper and lower bounds on linear transformations.
 
\begin{lem}\label{basic bounds for singular values}
Let $A \in \R^{m \times n}$, and $x \in \R^n$.  Then
\begin{equation}\label{eq: upper and lower bounds using singular values}
\sigma_{min}  \| x \| \leq \| A x \| \leq \sigma_{max} \| x \|,
\end{equation}
where $\sigma_{min}$ and $\sigma_{min}$ are the smallest and largest singular values of $A$, respectively. 
\hfill $\square$
\end{lem}
(See also \cite{Golub} for a proof of Lemma \ref{basic bounds for singular values})

\begin{thm}[Gourion and Seeger, 2010]
\label{upper and lower bounds for the solid angle}
Let $\sigma_{min}$  be the smallest singular value of $M$, and let $\sigma_{max}$ be the largest 
singular value of $M$.  Then:
\begin{equation}
\frac{| \det M | }{(2 \sigma_{max} )^d} \leq \omega_K 
\leq  \frac{| \det M | }{(2 \sigma_{min} )^d}.
\end{equation}
\end{thm}
\begin{proof}
\begin{align}
\omega_K &:= \int_K e^{-\pi \| x\|^2} dx = |\det M| \int_{\R^d_{\geq 0}} e^{-\pi \| Mu \|^2} du \\
&\leq |\det M| \int_{\R^d_{\geq 0}} e^{-\pi \sigma_{min}^2 \| u \|^2} du \\ 
&= \frac{|\det M|}{\sigma_{min}^d} \int_{\R^d_{\geq 0}} e^{-\pi  \| u \|^2} du \\ 
&= \frac{|\det M|}{2^d \sigma_{min}^d},
\end{align}
as claimed.  The inequality above  followed from the lower bound given by 
Lemma \ref{basic bounds for singular values}: $\| M u \| \geq  \sigma_{min} \|u\|$. 
The lower bound is proved in exactly the same manner, using the upper bound of 
Lemma \ref{basic bounds for singular values}: $\| M u \| \leq  \sigma_{max} \|u\|$.
\end{proof}

\begin{example}
{\rm
Let's consider the $3$-dimensional simplicial cone 
\[
K:= 
 \{  \lambda_1  
\left( \begin{smallmatrix}
1   \\ 0 \\ 0 
\end{smallmatrix}
\right) + 
\lambda_2 
\left( \begin{smallmatrix}
1   \\ 1 \\ 0 
\end{smallmatrix}
\right) 
+ 
\lambda_3 
\left( \begin{smallmatrix}
1   \\ 1 \\ 1 
\end{smallmatrix}
\right) \mid \lambda_1, \lambda_2, \lambda_3 \geq 0 \}.
\]
It's not difficult to show that its solid angle is $\omega_K = \frac{1}{48}$, 
by tesellating all of
 $\R^d$ with isometric images of $K$ (Exercise \ref{a unimodular cone}).
Computing (brute-force) the minimum and maximum singular values  for the matrix $M$ whose columns are the edge vectors of $K$, and substituting them into 
Theorem \ref{upper and lower bounds for the solid angle},  we get:
\[
0.01101 \leq  \omega_K \leq 0.73135.
\]
The latter lower bound gets much closer here to the true value  $\omega_K = \frac{1}{48} \approx .02083$.
}
\hfill $\square$
\end{example}


\bigskip
\section{The classical Euler-Maclaurin summation formula}

Here we show yet another application of Poisson summation, which has found great applications in number theory and numerical analysis: the classical Euler-Maclaurin (EM) summation formula. 

\begin{thm}[Euler-Maclaurin summation I]
Suppose $f:\R \rightarrow \C$ is infinitely smooth, compactly supported on $[a, b]$, and 
$f, \hat f \in L^1(\R)$. 
Then we have:
\begin{equation}
\sum_{a\leq n\leq b} f(n) - \int_{a}^b f(x) dx = f(b) P_1(b) - f(z) P_1(a) - \int_a^b f'(x) P_1(x) dx.
\end{equation}
\end{thm}
\begin{proof}
Applying Poisson summation to $f$, we have:
\begin{align*}
\sum_{n\in \Z} f(n) 
&= \sum_{\xi \in \Z} \hat f(\xi)  
= \hat f(0) +  \sum_{\xi \in \Z\setminus \{0\}} \hat f(\xi)  \\ 
&=  \int_{\R} f(x) dx+  \sum_{\xi \in \Z\setminus \{0\}} \hat f(\xi).
\end{align*}
Because $f$ is compactly supported on $[a, b]$, the latter equality becomes:
\begin{align}\label{first EM equality}
\sum_{a\leq n\leq b} f(n) -  \int_a^b f(x) dx
=  \sum_{\xi \in \Z\setminus \{0\}} \hat f(\xi).
\end{align}
Now we use integration by parts:
\begin{align}\label{second EM equality}
\hat f(\xi) &:= \int_\R f(x) e^{-2\pi i x \xi} dx = \int_a^b f(x) e^{-2\pi i x \xi} dx \\
&= f(b) \frac{e^{-2\pi i \xi b}}{-2\pi i \xi} - f(a) \frac{e^{-2\pi i \xi a}}{-2\pi i \xi}
-\int_a^b f'(x)   \frac{e^{-2\pi i \xi x}}{-2\pi i \xi} dx,
\end{align}
We'd like to plug the latter formula into \eqref{first EM equality}, but we'll do it carefully,
as follows: 
\begin{align}\label{third EM equality}
\sum_{a\leq n\leq b} f(n) - \int_a^b f(x) dx
&=   
 -  f(b) \lim_{N\rightarrow \infty} \sum_{-N\leq \xi \leq N \atop \xi \not=0}
          \frac{e^{-2\pi i \xi b}}{2\pi i \xi} 
 + f(a) \lim_{N\rightarrow \infty}  \sum_{-N\leq \xi \leq N \atop \xi \not=0} 
          \frac{e^{2\pi i \xi a}}{2\pi i \xi}\\
          \label{penultimate EM formula for proof}
&-\lim_{N\rightarrow \infty}   \int_a^b f'(x) \sum_{-N\leq \xi \leq N \atop \xi \not=0} 
                 \frac{e^{-2\pi i \xi x}}{2\pi i \xi} dx.
\end{align}
Because $f$ is infinitely smooth, both $\hat f$ and $\hat f'$ are rapidly decreasing. 
By Corollary \ref{example of theorem on pointwise convergence}, we know that 
\[
-\lim_{N\rightarrow \infty}   \sum_{-N\leq \xi \leq N \atop \xi \not=0} 
                 \frac{e^{-2\pi i \xi x}}{2\pi i \xi}  = \{ x\} - \frac{1}{2} := P_1(x),
\]
for $x \notin \Z$.  Therefore the required identity now follows from \eqref{penultimate EM formula for proof}.
\end{proof}
The hypothesis in this initial version of Euler-Maclaurin summation formula may be weakened considerably, but this will suffice for now.

Our proof above (originally due to G. H. Hardy) is not the easiest proof, but rather opens
up a path to higher dimensions.  Such an endeavor, in arbitrary dimension,  entails a long and winding road, so here we'll content ourselves with only a taste of it.  
One of the first applications of EM summation to number theory was the asymptotic approximation for the
tail of the Riemann zeta function:
\[
\sum_{n>x} \frac{1}{n^s} = O(x^{1-s}),
\]
for $s>1$ (\cite{ApostolBook}, Theorem 3.2).

\bigskip
\section{Further topics}

There is a fascinating conjecture related to bounding the smallest solid angle of any simplex, 
which was stated in \cite{Karasev.et.al}.
\begin{conjecture}
\label{conjecture: smallest solid angle of a simplex}
Any $d$-dimensional simplex $\P \subset \R^d$ has a solid angle not greater than the
solid angle of the $d$-dimensional regular simplex.
\end{conjecture}
For dimension $2$, the conjecture is trivial, but this problem becomes highly non-trivial in dimensions $d\geq 3$.  It is known to be true (though non-trivial) in dimensions $d =3 $ and $d=4$, as shown by Akopyan and Karasev \cite{AkopyanKarasev}. 
For further details related to Conjecture \ref{conjecture: smallest solid angle of a simplex} see \cite{Karasev.et.al}. In dimensions $d \geq 5$, 
Conjecture  \ref{conjecture: smallest solid angle of a simplex} is still open.

Having seen two different types of discrete volumes, we might wonder:
\begin{question}[Rhetorical]
How do we define more general discrete volumes of convex bodies?
\end{question}
Well, given any convex, compact set 
$K\subset \R^d$, and a bounded function
 $f:\R^d \rightarrow \C$, we consider the finite sum 
 \begin{equation}\label{discrete volume sum definition}
\vol_f(K):=  \sum_{n \in K\cap \Z^d} f(n).
 \end{equation}
 If $\vol_f(K)$ enjoys the property that 
\begin{equation}\label{discrete volume definition}
\lim_{t\rightarrow +\infty} \frac{1}{t^d} \vol_f(tK) = \vol K,
\end{equation}
then we call $\vol_f(K)$ a {\bf discrete volume} of $K$ \index{discrete volume}. In this chapter we considered $f_1(x):= \omega_K(x)$, the solid angle of $K$, at  each $x\in \R^d$.  In Chapter \ref{Ehrhart theory}, we'll consider Ehrhart's theory of discretized volumes, given by the constant function $f_2(x) = 1$, for all $x \in \R^d$.  Both $f_1$ and $f_2$ satisfy our definition 
\eqref{discrete volume definition} of discrete volume.

A moment's thought reveals that this more general definition given by 
\eqref{discrete volume definition} is equivalent to the following.
\begin{lem}
The function $f:\R^d \rightarrow \C$ defines a discrete volume in the sense of 
\eqref{discrete volume sum definition} and 
\eqref{discrete volume definition} if and only if:
\[
\int_{K} f(x) dx = \vol K.
\]
\end{lem}
\begin{proof}
By definition, the Riemann sum approximation to the integral gives us:
\begin{align*}
\int_K f(x) dx &:= \lim_{t\rightarrow +\infty} \frac{1}{t^d}  \sum_{n \in K\cap \frac{1}{t} \Z^d} f(n) 
= \lim_{t\rightarrow +\infty} \frac{1}{t^d}  \sum_{n \in tK\cap \Z^d} f(n) 
:=\lim_{t\rightarrow +\infty} \frac{1}{t^d}  \vol_f(tK) \\
&= \vol K,
\end{align*}
where the last equality holds if and only if $f$ gives a discrete volume, by 
definition \eqref{discrete volume definition}.
\end{proof}

\bigskip \bigskip

\section*{Notes} \label{Notes.chapter.angle.polynomial}
\begin{enumerate}[(a)]
\item  Let's compare and contrast the  two notions of discrete volumes that we have encountered so far. 
For a given rational polytope $\P$, we notice that the Ehrhart quasi-polynomial $L_\P(t)$ is invariant when 
we map $\P$ to any of its unimodular images.  That is, any rational polytope in the whole orbit of the unimodular group
 $\rm{SL}_d(\Z)(\P)$ has the same discrete volume $L_\P(t)$.  This is false for the second discrete volume $A_\P(t)$ - it is not invariant under the modular group 
(Exercise \ref{counterexamples for angle polynomials}).   But   $A_\P(t)$      is invariant under the
 large finite group of the isometries of $\R^d$ that preserve the integer lattice (known as the hyperoctahedral group). 
 
 So we see that $A_\P(t)$ is more sensitive to the particular embedding of $\P$ in space, because it is dependent upon a metric.   It is reasonable to expect that it can distinguish between ``more'' rational polytopes, but such a question remains to be formalized.  
 
 The angle polynomial  also has the advantage of being a much more symmetric polynomial, with half as many coefficients that occur
  in the Ehrhart polynomial of integer polytopes.  
  
  However,  $L_\P(t)$ has its advantages as well  -  to compute a \emph{local summand} for the angle polynomial $A_\P(t):= \sum_{n\in\Z^d}   \omega_{tP}(x) $ requires finding the volume of a local spherical polytope, while to compute a \emph{local summand} for the Ehrhart polynomial 
  $L_\P(t):= \sum_{n\in t \P \cap \Z^d} 1$ is quite easy: it is equal to $1$.
   
But as we have seen, computing the full global sum for $A_\P(t)$ turns out to have its simplifications.

\item There are natural ways to associate probabilities with solid angles -  see for example the work of Klain and Feldman  \cite{KlainFeldman}.

\item    Nhat Le Quang developed a thorough analysis of solid angle sums in $\R^2$, for all rational polygons, in his $2010$ undergraduate dissertation \cite{Nhat}.

\item  The recent work of Gerv\'asio  \cite{GervasioSantos}  gives an online implementation for the calculation of solid angles in any dimension, with open source code.  In fact, the thesis \cite{GervasioSantos}  contains extensive numerical computations that give empirical distributions for the bounds of Theorem \ref{upper and lower bounds for the solid angle} on random integer cones.
 
 \item  In \cite{RicardoNhatSinai}, there is an explicit description for some of the coefficients of the solid angle polynomial 
 $A_\P(t)$ of a $d$-dimensional polytope, for all positive real dilations $t>0$.   Indeed, the approach in \cite{RicardoNhatSinai} uses  the Fourier analytic landscape.

\item
There is also a characterization of $k$-tiling  $\R^d$ by using 
solid angle sums  \cite[Theorem 6.1]{GravinShiryaevRobins}, as follows.
\begin{lem}[Gravin, Robins, and Shiryaev]
\label{Gravin, Robins, Shiryaev solid angle k-tiling}
A polytope $P$ $k$-tiles $\R^d$ by integer translations if and only if
\[
\sum_{n\in \Z^d} \omega_{P + v}(n) = k,
\]
for every $v \in \R^d$.
\hfill $\square$
\end{lem}

\end{enumerate}


\bigskip
\section*{Exercises}
\addcontentsline{toc}{section}{Exercises}
\markright{Exercises}

\begin{quote}
  ``I haven't failed, I have just successfully found $10,000$ ways that won't work.''

-- Thomas Edison
\end{quote}

\medskip
\begin{prob}\label{a unimodular cone}
\rm{
Let $
\K =  \{  \lambda_1  
\left( \begin{smallmatrix}
1   \\ 0 \\ 0 
\end{smallmatrix}
\right) + 
\lambda_2 
\left( \begin{smallmatrix}
1   \\ 1 \\ 0 
\end{smallmatrix}
\right) 
+ 
\lambda_3 
\left( \begin{smallmatrix}
1   \\ 1 \\ 1 
\end{smallmatrix}
\right) \mid \lambda_1, \lambda_2, \lambda_3 \geq 0 \},
$
a simplicial cone.  Show that the solid angle of $\K$ is $\omega_\K = \frac{1}{48}$.
}
\end{prob}

\medskip
\begin{prob}
{\rm
We recall the $2$-dimensional cross-polytope
$
\Diamond:=\left\{ \left( x_1, x_2 \right) \in \R^2 \mid
 \, \left| x_1 \right| + \left| x_2 \right|  \leq 1 \right\}.
$
Find, from first principles, the angle quasi-polynomial for the rational polygon
$\P:= \frac{1}{3}\Diamond$, 
for all integer dilations of $\P$. 
}
\end{prob}

\medskip
\begin{prob}
{\rm
We recall that the $3$-dimensional cross-polytope was defined by 
\[
\Diamond:=\left\{ \left( x_1, x_2, x_3 \right) \in \R^3 \mid
 \, \left| x_1 \right| + \left| x_2 \right| + \left| x_3 \right| \leq 1 \right\}.
\]
Compute the angle polynomial of $A_{\Diamond}(t)$.
}
\end{prob}

 \medskip
\begin{prob}
{\rm
We recall that the $d$-dimensional cross-polytope \index{cross-polytope}
was defined by 
\[
\Diamond:=\left\{ \left( x_1, x_2, \dots, x_d \right) \in \R^d \mid
 \, \left| x_1 \right| + \left| x_2 \right| + \cdots + \left| x_d \right| \leq 1 \right\}.
\]
Compute the angle polynomial of $A_{\Diamond}(t)$.
}
\end{prob}

\medskip
\begin{prob}
{\rm
Let $\P$ be an integer zonotope.    Prove that the angle polynomial of $\P$ is
\[
A_{\P}(t) = (\vol \P)t^d,
\]
valid for all positive integers $t$.
}
\end{prob}
 Notes.  \ Although at this point in our development this problem may be challenging, 
once the reader uses Stokes' theorem  (Chapter \ref{Stokes' formula and transforms}) this problem will become quite easy.
 
\medskip
\begin{prob}  
{\rm
Using \eqref{recovering 1-dim'l solid angle poly}, find the angle quasi-polynomial $A_{\P}(t)$
for the $1$-dimensional polytope $\P:= [\frac{1}{2}, \frac{2}{3}]$.
}
\end{prob}

\medskip
\begin{prob}
{\rm
Generalizing the previous exercise, using \eqref{recovering 1-dim'l solid angle poly}, compute the angle quasi-polynomial
$A_{\P}(t)$ for any rational $1$-dimensional polytope $\P:=  [\frac{a}{c}, \frac{b}{d}]$. 
}
\end{prob}

\medskip
\begin{prob}  
{\rm
Define the rational triangle $\Delta$ whose vertices are $(0, 0), (1, \frac{N-1}{N}), (N, 0)$, where $N \geq 2$ is a fixed integer.
Find the angle quasi-polynomial $A_\Delta(t)$. 
}
 \end{prob}

\medskip
\begin{prob} $\clubsuit$  \label{another integral for the solid angle}
{\rm
Let $\K\subset \R^d$ be a $d$-dimensional polyhedral cone, and fix $s > \frac{d}{2}$. 
Prove that the solid angle $\omega_K$ has the alternate expression:
\[
\omega_K = 
\frac{ 
\pi^{\tfrac{d}{2}} \Gamma(s) 
}
{
\Gamma(s - \tfrac{d}{2})
}
 \int_K \frac{dx}{  \left( 1+ \|x\|^2 \right)^s}.
\]
}
\end{prob}

\medskip
\begin{prob} $\clubsuit$
\label{general integral for solid angle}
{\rm
Given a non-negative function $f:\R^d \rightarrow \R_{\geq 0}$, prove that the following are equivalent:
\begin{enumerate}[(a)]
\item (radially symmetric probability distribution) \ 
$f$ is radially symmetric, and 
 \[
 \int_{\R^d}  f(x) dx =1.
 \]
\item (solid angle integral) \ $\int_K  f(x) dx = \omega_K$
for all $d$-dimensional  polyhedral cones $\K\subset \R^d$.
\end{enumerate}
}
\end{prob}

\medskip
\begin{prob} $\clubsuit$
\label{counterexamples for angle polynomials}
{\rm
For each dimension $d$, find an example of an integer polytope $\P \subset \R^d$ and a unimodular matrix 
$U \in \rm{GL}_d(\Z)$, such that  the angle quasi-polynomials
$A_{\P}(t)$ and $A_{U(\P)}(t)$ are not equal to each other for all $t \in \Z_{>0}$.
}
\end{prob}

\medskip
\begin{prob} $\clubsuit$ 
\label{solid angle of a face of the cube}
{\rm
For the cube $\square:= [0, 1]^d$, show that any  face $F\subset \square$ that has dimension $k$ has the solid angle 
$\omega_F = \frac{1}{2^{d-k}}$.
}
\end{prob}

\medskip  
\begin{prob} $\clubsuit$
 \label{dihedral angle=solid angle}
{\rm
Show that the solid angle $\omega_E$ of an edge E ($1$-dimensional face) of a polytope 
equals the dihedral angle between the hyperplanes defined by
its two bounding facets.  (Hint: use the unit normal vectors for both of the bounding facets)
}
\end{prob}

\medskip
\begin{prob} 
{\rm
Using the Gram relations, namely Theorem \ref{Gram relations},  compute the solid angle at any vertex of the following regular tetrahedron:
\[
T:= \conv\Big\{ \icol{1\\0\\0} \icol{0\\1\\0},   \icol{0\\0\\1},     \icol{1\\1\\1}    \Big\}.
\]
}
\end{prob}

\medskip
\begin{prob} 
{\rm
Can you find a convex body $\P\subset \R^2$ that achieves the equality case in the upper bound of 
Nosarzewska's inequality \eqref{Nosarzewska inequality}?
 }
\end{prob}


 \chapter{\blue{The discrete Brion theorem:  Poisson summation strikes again} }
 \label{chapter:Discrete Brion}
 \index{Poisson summation} \index{discrete Brion theorem}

\begin{quote}   
``Everything you've learned in school as `obvious'   becomes less and less obvious as you begin to study the universe. 
 For example, there are no solids in the universe. 
 There's not even a suggestion of a solid. 
 There are no absolute continuums. 
 There are no surfaces. 
 There are no straight lines.''

-- Buckminster Fuller 
\end{quote}

\bigskip
\red{(Under construction)}

\section{Intuition}
As we saw in Theorem \ref{brion, continuous form}, there exists a wonderful way to decompose the Fourier transform of a polytope in terms of the Fourier-Laplace transforms of its vertex tangent cones.
We can now ask:
\begin{question}{\rm [Rhetorical]
\label{question:discrete Brion}
Is there a natural way to {\bf discretize} the continuous identity \eqref{transform formula for a simple polytope} of Brion, 
for the Fourier transform of a polytope?}
\end{question}

Another basic question we could ask is:
\begin{question} {\rm [Rhetorical]
How does the finite geometric sum in dimension $1$ extend to dimension $d$? }
\end{question}
As we'll see, these two questions are intertwined, and one answers the other.  
One useful way to make sense of Question \ref{question:discrete Brion} is to replace integrals with sums over the integer lattice:
\begin{equation}
 \int_\P e^{-2\pi i \langle u,  z  \rangle} \, du      \longrightarrow
  \sum_{n\in \Z^d} e^{2\pi i \langle z,  n  \rangle}.
\end{equation}
Such a descretization will lead us to a discrete version of Brion's Theorem, namely Theorem \ref{brion, discrete form} below.  Although the discrete Brion theorems of this chapter have several applications, for us the main application will be the enumeration of lattice points in polytopes, which is the Ehrhart theory of Chapter \ref{Ehrhart theory}. 

 
\bigskip
\section{Discretizing the Fourier-Laplace transform of a cone}  

We may also replace the integer lattice by any lattice $\L$, and the ensuing function is very similar.  
But since this is only a cosmetic change of variable, we can simplify life and work with the integer lattice. 
To this discrete end, we define the {\bf integer point transform} 
\index{integer point transform}   of a rational polytope $\P$ by
\[
\sigma_\P(z) :=  \sum_{n \in \P \cap \Z^d}    e^{2\pi i \langle n,  z\rangle},
\]
a discretization of the Fourier transform of $\P$. 

We may also think of the discretized sum $  \sum_{n\in \Z^d} e^{2\pi i \langle z,  n  \rangle}$
more combinatorially by making the change of variable $q_1:= e^{2\pi i z_1}, \dots, q_d:= e^{2\pi i z_d}$, so that 
we have $q_1^{n_1} q_2^{n_2} \cdots q_d^{n_d} =  e^{2\pi i n_1 z_1 + \cdots +  2\pi i n_d z_d} := e^{2\pi i \langle n, z \rangle}$. 
with this notation in mind, we define the {\bf multinomial notation} for a monomial in several variables: 
\[
q^n:= q_1^{n_1} q_2^{n_2} \cdots q_d^{n_d}.
\]
We will therefore sometimes use the equivalent definition 
\[
\sigma_\P(q):=  \sum_{n \in \P \cap \Z^d}    q^n.
\]

We similarly define the {\bf integer point transform of a rational cone} $\K_v$ by the series
\begin{equation}\label{def of integer point transform of a cone} 
\sigma_{\K_v}(z) :=    \sum_{n \in \K_v \cap \Z^d}    e^{2\pi i  \langle n,  z\rangle}.
\end{equation}
But even in dimension $1$ things can  get interesting, so let's see an example. 

\begin{example}[Finite geometric sums]
\rm{
Consider the $1$-dimensional polytope $\P := [a,b]$, where $a, b\in \Z$.   The problem is to compute the finite
geometric series:
\begin{align*}
\sum_{n \in \P \cap \Z}    e^{ 2\pi i n z}   &=  \sum_{a \leq n \leq b}    q^n,
\end{align*}
where we've set $q:= e^{2\pi i z}$.  Of course, we already know that it possesses a 
 `closed form' of the type:
\begin{align}  \label{verifying1}
 \sum_{a \leq n \leq b}    q^n &= \frac{q^{b+1} -q^{a} }{q - 1} \\
 &=  \frac{q^{b+1}}{q-1}     -  \frac{ q^{a} }{q - 1},  \label{verifying2}
\end{align}
because we already recognize this formula for a {\bf finite geometric sum}.    
On the other hand, anticipating the discrete form of Brion's theorem below, we first compute the discrete 
sum corresponding to the vertex tangent cone at the vertex $a$, namely $\sum_{a \leq n}    q^n$:
\begin{equation}
\label{cone identity1}
q^a+ q^{a+1} + \cdots = \frac{q^a}{1-q}.
\end{equation}
Now we compute the the sum corresponding to the vertex tangent cone at vertex $b$, namely  $\sum_{n \leq b}    q^n$:
\begin{equation}  \label{cone identity2}
q^b+ q^{b-1} + \cdots = \frac{q^b}{1-q^{-1}} =
\frac{q^{b+1}}{q-1} . 
\end{equation}
 Summing these two contributions, one from each vertex tangent cone, we get:
\begin{align*}
\frac{q^a}{1-q} + \frac{q^{b+1}}{q-1} =  \sum_{a \leq n \leq b}    q^n,
\end{align*}
by the finite geometric sum identity, thereby 
verifying Theorem \ref{brion, discrete form} for this example.
This example shows that Brion's Theorem \ref{brion, discrete form}  (the discrete version) may be thought of as a $d$-dimensional extension of the finite geometric sum.             

But something is still very wrong here - namely, identity  \eqref{cone identity1} converges for $|q|< 1$, 
while identity \eqref{cone identity2}   converges only for $|q|>1$, so there is not even one value of $q$ for which the required identity \eqref{verifying2}  is true.  So how can we make sense of these completely {\bf disjoint domains of convergence} ?!
}
\hfill $\square$
\end{example}

\bigskip
To resolve these conundrums, the very useful result of 
Michel Brion \cite{Brion} comes to the rescue.   
Our proof of Theorem \ref{brion, discrete form} discretizes the continuous form of Brion's Theorem \ref{brion, continuous form}, 
using the Poisson summation formula, to arrive at a discrete form of Brion's Theorem.

First, we need a slightly technical but easy Lemma.
\begin{lem}\label{technical lemma1}
Let $\K_v$ be a rational cone, with apex at $v$.  We pick any compactly supported and infinitely smooth
approximate identity $\phi_\varepsilon$.   Then:
\begin{align}  \label{lem:technical claim for limiting series}
     \lim_{\varepsilon \rightarrow 0}   
     \sum_{n \in \Z^d } 
        \left(       1_{\interior \K_v}(x)    e^{2\pi i \langle x,  z\rangle}*\phi_\varepsilon      \right)(n)
        =    \sum_{n \in \Z^d \cap \interior \K_v }      e^{2\pi i \langle n,  z\rangle} := \sigma_{ \interior \K_v}(z).
\end{align}
\end{lem}
\begin{proof}
We first note that by our assumptions on  $\phi_\varepsilon$, it lies in the Schwartz space  $S(\R^d)$, by 
 Lemma \ref{useful Schwartz fact}.
So $\phi_\varepsilon$ is rapidly decreasing.  Using the Weierstrass $M$-test, 
we see that the series 
$\sum_{n \in \Z^d } 
        \left(       1_{\interior \K_v}(x)    e^{2\pi i \langle x,  z\rangle}*\phi_\varepsilon      \right)(n)$
 converges uniformly in $\epsilon$, and because the summands are continuous functions of $\epsilon$, so is the whole series.  So we may take the limit as $\epsilon \rightarrow 0$ inside the series.   
Finally, using  Lemma \ref{approximate identity convolution}, and the continuity of the function 
$1_{\interior \K_v}(x)    e^{2\pi i \langle x,  z\rangle}*\phi_\varepsilon$ at all $x \in \R^d$, 
we have 
 $\lim_{\varepsilon \rightarrow 0}   
        \left(       1_{\interior \K_v}(x)    e^{2\pi i \langle x,  z\rangle}*\phi_\varepsilon      \right)(n)
        =  1_{\interior \K_v}(n)    e^{2\pi i \langle n,  z\rangle}$, from which 
        \eqref{lem:technical claim for limiting series} follows. 
\end{proof}


It turns out that the continuous form of Brion's theorem, namely Theorem \ref{brion, continuous form}, 
 can be used to prove the discrete form of Brion's theorem, namely Theorem \ref{brion, discrete form} below.  

\bigskip
\begin{thm}[{\bf Brion's theorem - the discrete form, 1988}]   \label{brion, discrete form}
\index{Brion's theorem - the discrete form}
Let $\P \subset \R^d$ be a  rational, $d$-dimensional  polytope, and let $N$ be the number of vertices of $\P$.
For each vertex $v$ of $\P$, we consider the open vertex tangent cone $\interior \K_v$  of $ \interior \P$, the interior of $\P$.  
Then
\begin{equation} \label{Discrete formula, Brion's theorem}
\sigma_{ \interior  \P}(z) =  \sigma_{ \interior  \K_{v_1}}(z) + \cdots  + \sigma_{ \interior  \K_{v_N}}(z).
\end{equation}
for all $z \in \C^d - S$, where $S$ is the hyperplane arrangement defined by the (removable) singularities of all of the transforms $\hat 1_{\K_{v_j}}(z)$.
\end{thm}
\begin{proof}      We will use the continuous version of Brion, namely Theorem \ref{brion, continuous form}, together with 
the Poisson summation formula, to deduce the discrete version here.   In a sense, the Poisson summation formula allows us to discretize the integrals.   \index{Poisson summation formula}

Step $1$.  [{\bf Intuition - fast and loose}] \  To begin, in order to motivate the rigorous proof that follows, we will use Poisson summation on a function $1_{\P}(n) e^{2\pi i \langle n,  z\rangle}$
 that ``doesn't have the right"  
to be used in Poisson summation, because $\hat 1_{\P} \notin L^1(\R^d)$ .  But this first step brings the intuition to the  foreground. 
   Then, in Step $2$,  we  will literally ``smooth'' out the lack of rigor in Step 1, by smoothing $1_\P$ with an approximate identity. 
\begin{align*}
\sum_{n \in \P \cap \Z^d} e^{2\pi i \langle n,  z\rangle}     &:= \sum_{n \in \Z^d}  1_{\P}(n) e^{2\pi i \langle n,  z\rangle}  \\
&=  \sum_{\xi \in \Z^d}  \hat 1_\P(z+ \xi)  \\
&=  \sum_{\xi \in \Z^d} \left(      \hat 1_{K_{v_1}}  (z+ \xi) + \cdots +     \hat 1_{K_{v_1}}  (z+ \xi)   \right)  \\
&=  \sum_{\xi \in \Z^d}   \hat 1_{K_{v_1}}  (z+ \xi) + \cdots +    \sum_{\xi \in \Z^d}  \hat 1_{K_{v_N}}  (z+ \xi)    \\
&=  \sum_{n \in \Z^d}   1_{K_{v_1}}  (n)  e^{2\pi i \langle n,  z\rangle}     + \cdots +   
        \sum_{n\in \Z^d}  1_{K_{v_N}}  (n)    e^{2\pi i \langle n,  z\rangle}   \\
     &:=      \sum_{n \in \Z^d \cap K_{v_1} }    e^{2\pi i \langle n,  z\rangle}  + \cdots +    
                \sum_{n \in \Z^d \cap K_{v_N} }   e^{2\pi i \langle n,  z\rangle},    
\end{align*}
where we have used the Poisson summation formula in the second and fifth equalities.
The third equality uses Brion's Theorem \ref{brion, continuous form} for the Fourier transform of $\P$.

Step $2$    [{\bf Rigorous proof}].  \ To make Step $1$ rigorous, we pick any compactly supported approximate identity 
$\phi_\varepsilon$, and form a smoothed 
version of the function in step $1$.  Namely we let 
\[
f_\varepsilon(x):=  (1_{\P}(x)   e^{2\pi i \langle x,  z\rangle})*\phi_\varepsilon(x),
\]
so that now we are allowed to apply Poisson summation \index{Poisson summation formula}
 to $f_\varepsilon$, because our choice of a smooth and compactly supported $\phi_\varepsilon$ implies that $f_\varepsilon$ is a Schwartz function.
Recalling Theorem \ref{approximate identity convolution}, we know that at a point $x\in \R^d$ of continuity of 
$1_{\P}(x)   e^{2\pi i \langle x,  z\rangle}$, we have
\[
\lim_{\varepsilon \rightarrow 0} f_\varepsilon(x) = 1_{\P}(x)   e^{2\pi i \langle x,  z\rangle}.   
\]
To proceed further,
it is therefore natural to consider points $x \in \interior \P$, the interior of $\P$, because $1_\P$ is continuous there, while it is not continuous on the boundary of $\P$.   To recap, we have so far the equalities 
\begin{align*}
\sum_{n \in \interior \P \cap \Z^d} e^{2\pi i \langle n,  z\rangle}   &:=   
                    \sum_{n \in \Z^d}     1_{ \interior \P}(x)   e^{2\pi i \langle x,  z\rangle} \\
 &= \sum_{n \in  \interior \P \cap \Z^d} \lim_{\varepsilon \rightarrow 0} f_\varepsilon(n)  \\
 &= \lim_{\varepsilon \rightarrow 0} \sum_{n \in \interior \P \cap \Z^d} f_\varepsilon(n),
 \end{align*}
 where we've used the fact that $f_\varepsilon$ is compactly supported, because it is the convolution of two compactly supported functions.  So the exchange above, of the sum with the limit,  is trivial because the sum is finite. 
With this in mind, the Poisson summation formula,  applied to the Schwarz function $f_\varepsilon$,
gives us:
\begin{align*}
&\sum_{n \in \interior \P \cap \Z^d} e^{2\pi i \langle n,  z\rangle}   = 
\lim_{\varepsilon \rightarrow 0}  \sum_{n \in  \interior  \P \cap \Z^d} f_\varepsilon(n)    
 = \lim_{\varepsilon \rightarrow 0}   \sum_{n \in \Z^d}  \left(1_{ \interior  \P} \  
               e^{2\pi i \langle x,  z\rangle})*\phi_\varepsilon\right) (n)  \\
 &=\lim_{\varepsilon \rightarrow 0}   \sum_{n \in \Z^d}  \F{ \big(   (1_{ \interior \P} \  
             e^{2\pi i \langle x,  z\rangle})*\phi_\varepsilon    \big)   }(\xi)  \\
 &= \lim_{\varepsilon \rightarrow 0}   \sum_{\xi \in \Z^d}  \hat 1_{ \interior \P}(z+ \xi)     \hat  \phi_\varepsilon(\xi)  \\
&=  \lim_{\varepsilon \rightarrow 0}  
                 \sum_{\xi \in \Z^d} \left(      \hat 1_{ \interior  \K_{v_1}}  (z+ \xi) + \cdots +     
                  \hat 1_{ \interior  \K_{v_1}}  (z+ \xi)   \right)
                   \hat  \phi_\varepsilon(\xi)  \\
&= \lim_{\varepsilon \rightarrow 0}  
               \sum_{\xi \in \Z^d}    \F{ \big(   (1_{ \interior   \K_{v_1}} \  e^{2\pi i \langle x,  z\rangle})*\phi_\varepsilon    \big)   }(\xi)       
                     + \cdots +    
     \lim_{\varepsilon \rightarrow 0}   \sum_{\xi \in \Z^d}  
                 \F{ \big(   (1_{ \interior    \K_{v_N}} \  e^{2\pi i \langle x,  z\rangle})*\phi_\varepsilon    \big)   }(\xi)  \\
&= \lim_{\varepsilon \rightarrow 0}  
               \sum_{n \in \Z^d}      (1_{ \interior  \K_{v_1}} \  e^{2\pi i \langle x,  z\rangle})*\phi_\varepsilon   (n)       
                     + \cdots +    
          \lim_{\varepsilon \rightarrow 0}       
                   \sum_{\xi \in \Z^d}     (1_{ \interior  \K_{v_N}} \  e^{2\pi i \langle x,  z\rangle})*\phi_\varepsilon   (n) \\
&=  \sigma_{ \interior  \K_{v_1}}(z) + \cdots  + \sigma_{ \interior  \K_{v_N}}(z),
\end{align*}
We've applied Theorem \ref{approximate identity convolution} to $f(n) := 1_{\interior \K_v}(n)$, 
for each $n \in \interior \K_v$, because $f$ is continuous at all such points.
The conclusion of Theorem  \ref{approximate identity convolution} is that 
\[
\lim_{\varepsilon \rightarrow 0}  
       \Big(  (1_{ \interior  \K_{v_1}} \  e^{2\pi i \langle x,  z\rangle})*\phi_\varepsilon   \Big)(n) = 
         1_{ \interior  \K_{v_1}}(n) \  e^{2\pi i \langle n,  z\rangle},
\]
and by Lemma \ref{technical lemma1} the last equality,  in the long string of equalities above,  is justified. 
\end{proof}

\bigskip
\begin{example} \label{example:standard triangle integer point transform}
\rm{
We can now recompute the integer point transform of  the {\bf standard triangle} in the plane, using Brion's Theorem \ref{brion, discrete form}.   
Namely,  for the standard triangle
\[
\Delta:= \conv(   \icol{0\\0},   \icol{1\\0},  \icol{0\\1}),
\]
as depicted in Figure \ref{standard triangle}, we find $\sigma_{\Delta}(z)$.
\noindent
\begin{figure}[htb]
\begin{center}
\includegraphics[totalheight=3.7in]{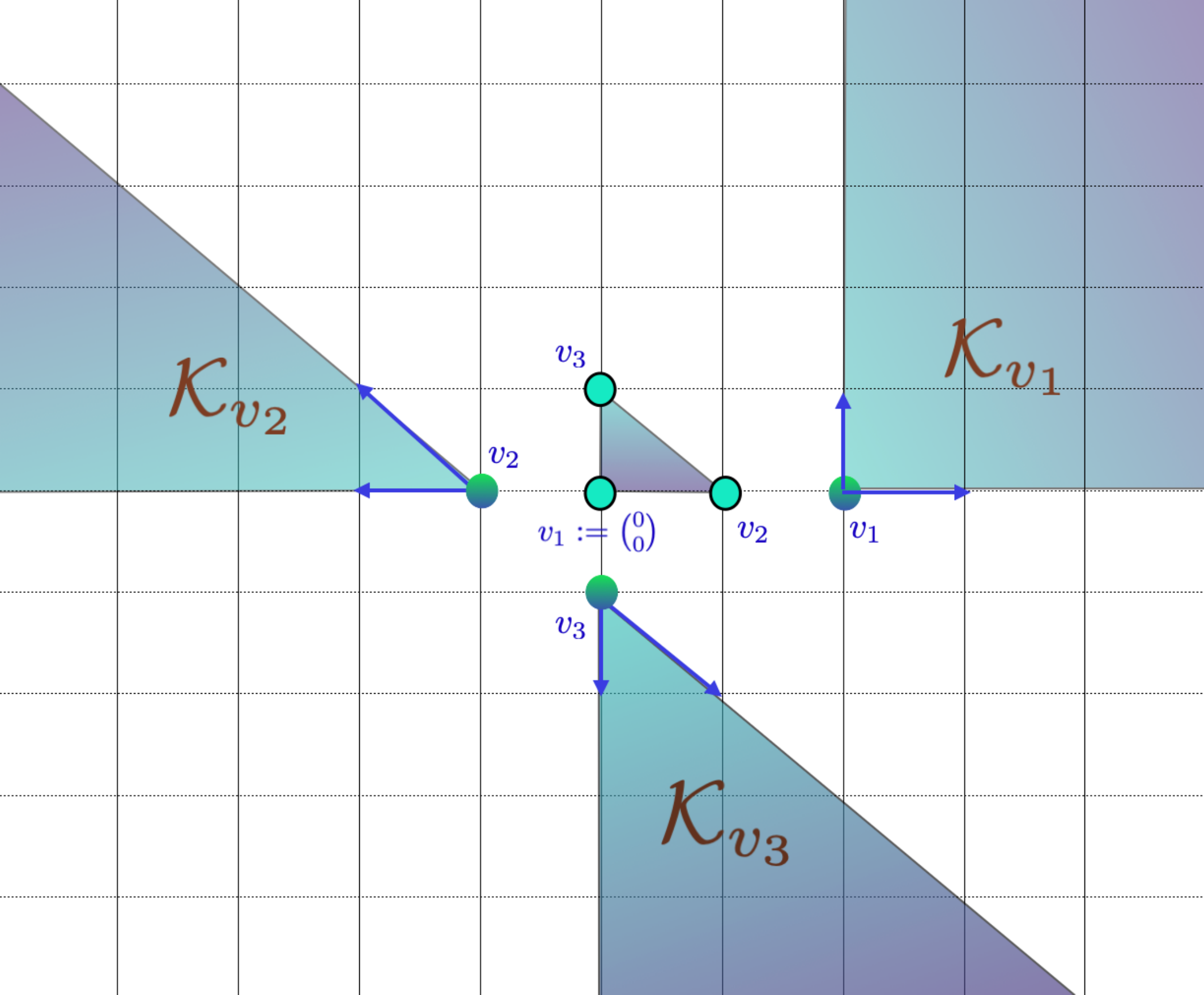}
\end{center}
\caption{The standard triangle, with its vertex tangent cones}    \label{standard triangle}
\end{figure}
By definition, the integer point transform of its vertex tangent cone $\K_{v_1}$ is
\begin{align*}
\sigma_{\K_{v_1}}(z) &:=  \sum_{n \in \K_{v_1} \cap \Z^d}    e^{\langle n, z\rangle}
= 
\sum_{n_1 \geq 0, n_2 \geq 0}    
e^{ \langle   n_1 \icol{1\\0}   +     n_2\icol{0\\1}  ,  
z\rangle} \\
& =   \sum_{n_1 \geq 0}     e^{  n_1 z_1}    \sum_{n_2 \geq 0}     e^{  n_2 z_2}  \\
& = \frac{1}{(1- e^{z_1} )(1-  e^{z_2})}.
\end{align*}

For the vertex tangent cone $\K_{v_2}$, we have
\begin{align*}
\sigma_{\K_{v_2}}(z) &:=  \sum_{n \in \K_{v_2} \cap \Z^d}    e^{\langle n, z\rangle}
= 
\sum_{n_1 \geq 0, n_2 \geq 0}    
e^{\langle \icol{1\\0}   +   n_1 \icol{-1\\ \ 0}   +     n_2\icol{-1\\ \ 1}  ,  
z\rangle} \\
& =   e^{z_1}  \sum_{n_1 \geq 0}     e^{  n_1 (-z_1)}    
\sum_{n_2 \geq 0}     e^{  n_2 (-z_1+z_2)}  \\
& = \frac{e^{z_1}}{(1- e^{-z_1} )(1-  e^{-z_1+z_2})}.
\end{align*}

Finally, for the vertex tangent cone $\K_{v_3}$, we have
\begin{align*}
\sigma_{\K_{v_3}}(z) &:= \sum_{n_1 \geq 0, n_2 \geq 0}    
e^{\langle \icol{0\\1}   +    n_1 \icol{ \ 0\\-1}   +     n_2\icol{\ 1\\-1},   z \rangle} \\
& =   e^{z_2}  \sum_{n_1 \geq 0}     e^{  n_1 (-z_2)}    
\sum_{n_2 \geq 0}     e^{  n_2 (z_1-z_2)}  \\
& = \frac{e^{z_2}}{(1- e^{-z_2} )(1-  e^{z_1-z_2})}.
\end{align*}
Altogether, using \ref{simplified discrete Brion identity} we have  
\begin{align}
\sigma_{\P}(z) &= \sigma_{\K_{v_1}}(z) + \sigma_{\K_{v_2}}(z) +\sigma_{\K_{v_3}}(z) \\
&= \frac{1}{(1- e^{z_1} )(1-  e^{z_2})} + \frac{e^{z_1}}{(1- e^{-z_1} )(1-  e^{-z_1+z_2})}
+\frac{e^{z_2}}{(1- e^{-z_2} )(1-  e^{z_1-z_2})}.  \label{last line}
\end{align}
}
\hfill $\square$
\end{example}

 \noindent
\begin{figure}[htb]
\begin{center}
\includegraphics[totalheight=3.3in]{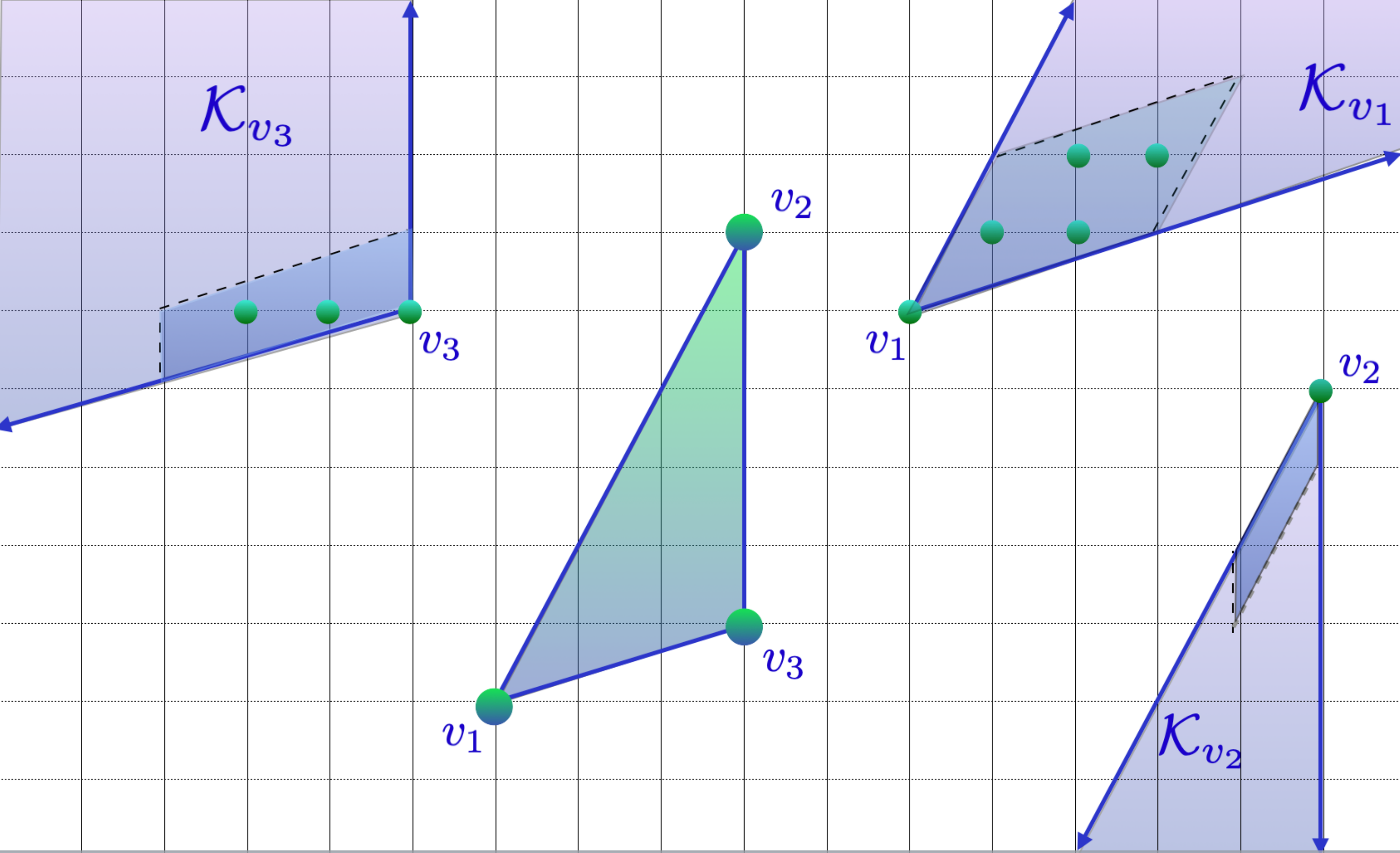}
\end{center}
\caption{A triangle with vertices $v_1, v_2, v_3$, and its vertex tangent cones}   \label{triangle}
\end{figure}


\bigskip
 \section{Examples, examples, examples}

 \bigskip
\begin{example}\label{triangle1}
Here we will compute the integer point transform of the triangle $\Delta$ defined by the convex hull of the points $\icol{0\\0}, \icol{3\\1}, \icol{3\\6}$, as shown in Figure \ref{triangle}.
We first compute the integer point transforms of all of its tangent cones.  For the vertex $v_1$, we already computed the integer point transform of its tangent cone in the previous example.  

For the vertex $v_2$, we notice that its vertex tangent cone is a unimodular cone, because
$| \det
\big(\begin{smallmatrix}
0 & -1  \\
-1 & -2
\end{smallmatrix}
\big)  | = 1$.    Its integer point transform is:
\begin{align*}
\sigma_{\K_{v_2}}(z) &:=  \sum_{n \in \K_{v_2} \cap \Z^d}    e^{\langle n, \  z\rangle}
= 
\sum_{n_1 \geq 0, n_2 \geq 0}    
e^{ \langle 
\icol{3\\6}   +   n_1 \icol{0\\-1}   +     n_2\icol{-1\\-2}  ,  
z\rangle} \\
& = e^{ 3z_1 + 6z_2)}   \sum_{n_1 \geq 0, n_2 \geq 0}    
e^{  n_1 (-z_2)}  e^{ n_2(-z_1 -2z_2)}    \\
& = \frac{ e^{3z_1 + 6z_2} }
{  
(1- e^{-z_2} )(1-  e^{ -z_1 -2z_2} )
}.
\end{align*}
Equivalently, using the notation from Example \ref{integer point transform of cone}     above, 
\[
\sigma_{\K_{v_2}}(z) := \sum_{n \in \K_{v_2} \cap \Z^d}   q^n  
=\frac{ q_1^3 q_2^6}
{
 (1-  q_2^{-1} )  (1-  q_1^{-1} q_2^{-2}).
}
\]

For vertex $v_3$, the computation is similar to vertex tangent cone $\K_{v_1}$, and
we have:
\begin{align*}
\sigma_{\K_{v_3}}(z) &:=  \sum_{n \in \K_{v_3} \cap \Z^d}    e^{\langle n, \  z\rangle}
=  \sum_{n_1 \geq 0, n_2 \geq 0}    e^{\langle 
\icol{3\\1}   +   n_1 \icol{-3\\-1}   +     n_2\icol{0\\1},  \ 
z\rangle} \\
&=   e^{3z_1 + z_2}    \sum_{n_1 \geq 0, n_2 \geq 0}    e^{(-3z_1-z_2) n_1} 
e^{2\pi i (z_2) n_2}   \\
&= e^{3z_1 + z_2}    \frac{ 1+  e^{-z_1}    +     e^{-2z_1}     }
{  
(1- e^{ 3z_1 + z_2} )(1-  e^{ z_2} )
} \\
&=    \frac{ e^{3z_1 + z_2}  +  e^{2z_1 + z_2}    +     e^{z_1 + z_2}     }
{  
(1- e^{3z_1 + z_2} )(1-  e^{z_2} )
} \\
&=   \frac{ q_1^3 q_2  +   q_1^2 q_2  +   q_1q_2    }
{  
(1-  q_1^{-3} q_2^{-1}  )(1-   q_2)  
}.  
\end{align*}
Finally, putting all of the three vertex tangent cone contributions together, 
Theorem \ref{brion, discrete form} gives us:
\begin{align*}
\sigma_{\Delta}(z) &= \sigma_{\K_{v_1}}(z) + \sigma_{\K_{v_2}}(z) + \sigma_{\K_{v_3}}(z) \\
&= \frac{   1+ q_1 q_2 +    {q_1}^2 q_2  + {q_1}^2 {q_2}^2 +    q_1^3 q_2^2  
}
{  (1- q_1^3  q_2 )  (1-  q_1 q_2^2)
}  
+ 
\frac{ q_1^3 q_2^6}
{
 (1-  q_2^{-1} )  (1-  q_1^{-1} q_2^{-2})
}
+
 \frac{ q_1^3 q_2  +   q_1^2 q_2  +   q_1q_2    }
{  
(1-  q_1^{-3} q_2^{-1}  )(1-   q_2)  
}.  
\end{align*}
\hfill $\square$
\end{example}

\bigskip
 \begin{example}\label{integer point transform of cone}
 \rm{
We work out the integer point transform $\sigma_\K(z)$ of the cone
\[
\K := \{ 
\lambda_1 \big(\begin{smallmatrix}
3  \\
1 \\
\end{smallmatrix}
\big) +    
\lambda_2 \big(\begin{smallmatrix}
1  \\
2 \\
\end{smallmatrix}
\big)  \mid  \lambda_1, \lambda_2 \in \R_{\geq 0}  \}, 
\]
Drawn in the figures below.   We note that here $\det \K = 5$, and that there are indeed $5$ integer points in
$D$, its half-open fundamental parallelepiped. 

\noindent
\begin{figure}[htb]
\begin{center}
\includegraphics[totalheight=3in]{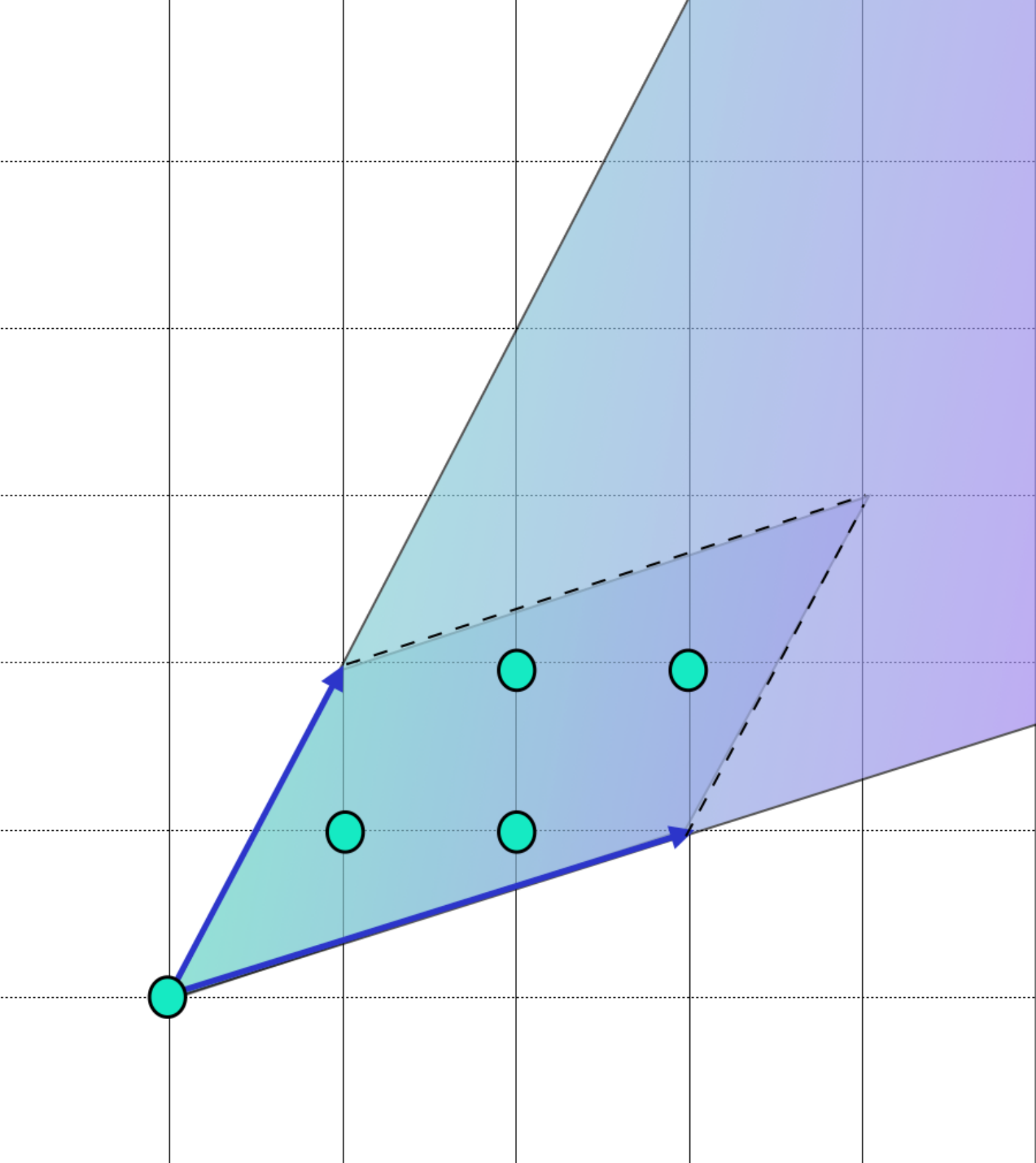}
\end{center}
\caption{The $5$ integer points in a fundamental parallelepiped $D$ of the cone $\K$.}    \label{cone2}
\end{figure}

\noindent
\begin{figure}[htb]
\begin{center}
\includegraphics[totalheight=2.5in]{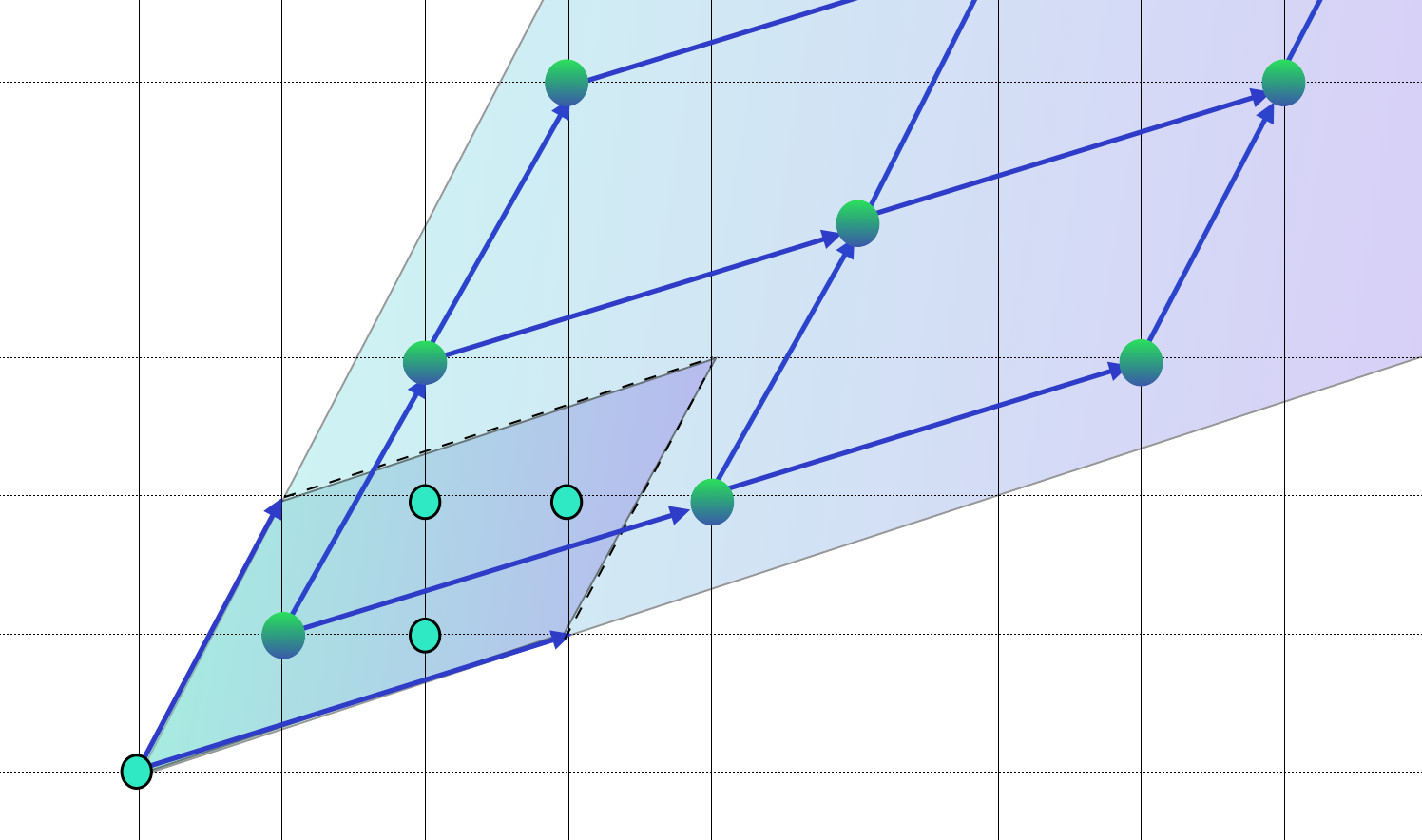}
\end{center}
\caption{The point ${1\choose 1}$  in $D$, with its images in $\K$ under translations by the edge vectors of $\K$.}   
 \label{cone3}
\end{figure}
We may `divide and conquer' the integer point transform $\sigma_\K(z)$, 
by breaking it up into $5$ infinite series,
one for each integer point in $D$, as follows:
\[
\sigma_\K(z) :=  \sum_{n \in \K\cap \Z^d}    e^{\langle n,  z\rangle}
:= \sum_{\icol{0\\0}} + \sum_{\icol{1\\1}} + \sum_{\icol{2\\1}} + \sum_{\icol{2\\2}} + \sum_{\icol{3\\2}},
\]
where 
\begin{align*}
\sum_{\icol{1\\1}}   &:=  \sum_{n_1 \geq 0, n_2 \geq 0}    
e^{\langle \big(\begin{smallmatrix}
1  \\
1 \\
\end{smallmatrix}  \big)   +
n_1 \big(\begin{smallmatrix}
3  \\
1 \\
\end{smallmatrix}
\big) +    
n_2 \big(\begin{smallmatrix}
1  \\
2 \\
\end{smallmatrix}
\big)  ,  z     \rangle} \\
&=
e^{ \langle  \big(\begin{smallmatrix}
1  \\
1 \\
\end{smallmatrix}  \big) , z \rangle}
\sum_{n_1 \geq 0, n_2 \geq 0}    
e^{ \langle  
n_1 \big(\begin{smallmatrix}
3  \\
1 \\
\end{smallmatrix}
\big) +    
n_2 \big(\begin{smallmatrix}
1  \\
2 \\
\end{smallmatrix}
\big)  ,  z\rangle} \\
&=
e^{ \langle  \big(\begin{smallmatrix}
1  \\
1 \\
\end{smallmatrix}  \big) , z \rangle}
\sum_{n_1 \geq 0}    
e^{n_1       \langle  
\big(\begin{smallmatrix}
3  \\
1 \\
\end{smallmatrix}
\big) , z\rangle}
\sum_{n_2 \geq 0}    
e^{n_2      \langle  
\big(\begin{smallmatrix}
1 \\
2 \\
\end{smallmatrix}
\big) , z\rangle}   \\
&=
\frac{e^{ z_1 + z_2}
}
{ (1-e^{3z_1 + z_2})(1-e^{z_1 + 2z_2})
},
\end{align*}
and similarly we have
\[
\sum_{\icol{2\\1}}= 
\frac{e^{2z_1 + z_2}
}
{ (1-e^{ 3z_1 + z_2})(1-e^{z_1 + 2z_2})
},
\]
\[
\sum_{\icol{2\\2}} = 
\frac{e^{ 2z_1 + 2z_2}
}
{ (1-e^{ 3z_1 + z_2})(1-e^{z_1 + 2z_2})
},
\]
\[
\sum_{\icol{3\\2}} = 
\frac{e^{ 3z_1 + 2z_2}
}
{ (1-e^{ 3z_1 + z_2})(1-e^{z_1 + 2z_2})
},
\]
and finally 
\[
\sum_{\icol{0\\0}}  = 
\frac{1}
{ (1-e^{3z_1 + z_2})(1-e^{z_1 + 2z_2})
}.
\]
To summarize, we have the following expression:
\[
\sum_{n \in \K\cap \Z^d}    e^{ \langle n,  z\rangle}
=\frac{   1+ e^{z_1 + z_2} + e^{ 2z_1 + z_2)} + e^{2z_1 + 2z_2} + e^{3z_1 + 2z_2 }
}
{  (1-e^{3z_1 + z_2})(1-e^{z_1 + 2z_2})
}.
\]
Equivalently, using our multinomial notation $q_j:= e^{ z_j}$, we have
\[
\sum_{n \in \K\cap \Z^d}   q^n  
=\frac{   1+ q_1 q_2 +    {q_1}^2 q_2  + {q_1}^2 {q_2}^2 +    q_1^3 q_2^2  
}
{  (1- q_1^3  q_2 )  (1-  q_1 q_2^2)
}.
\]
}
\hfill $\square$
\end{example}


\section{Integer point transforms of rational cones are rational functions}

The Examples  \ref{triangle1} and \ref{integer point transform of cone}  above suggest a general pattern, namely that integer point transforms are always rational functions, and that their numerators are polynomials that encode the integer points inside a fundamental parallelepiped $\Pi$
that sits at the vertex of each vertex tangent cone.   The proof of this general fact will be fairly easy - we only need to put several geometric series together, as in Figure \ref{cone3}.   Now that we've seen some examples, we can prove things in general. 

First, given any $d$-dimensional simplicial rational cone $\K\subset \R^d$, with integer edge vectors 
 $w_1, \dots, w_d \in \Z^d$, and apex $v\in \R^d$, we define the {\bf fundamental parallelepiped} of $\K$ by:
 \begin{equation}
 \Pi := \{   \lambda_1 w_1 + \cdots + \lambda_d w_d \mid 
 \text{ all }  0 \leq \lambda_j < 1 \},
 \end{equation}
a half-open, integer parallelepiped. 
 In the same way that we've encoded integer points in polytopes using 
$\sigma_\P(z)$,  we can encode the integer points in $\Pi$ by defining
\[
\sigma_\Pi(z) := \sum_{n \in \Z^d \cap \Pi} e^{\langle z, n \rangle}.
\]

For a rational simplicial cone $K_v$, it turns out that its integer point transform
\[
\sigma_{K_v}(z):= \sum_{n \in \K_v \cap \Z^d}    e^{ \langle n,  z \rangle}
\]
has a pretty structure theorem - it is a rational function of the variables $e^{z_1}, \dots, e^{z_d}$, as follows. 

\begin{thm}  \label{closed form for integer point transform of a cone}
Given a  $d$-dimensional simplicial cone $\K_v \subset \R^d$,  with apex $v \in \R^d$, and with 
$d$ linearly independent {\bf integer} edge vectors  
$w_1, w_2, \dots, w_d \in \Z^d$.   Then:  
\begin{equation}
\sigma_{K_v}(z) 
= 
\frac{   \sigma_{ \Pi + v}(z)     }
{       \prod_{k=1}^d  \left(      1 - e^{ \langle  w_k   , z \rangle}   \right)          }.
\end{equation}
\end{thm}
\begin{proof}
We claim that we can parametrize all of the integer points in the cone $\K_v$ precisely by
\begin{equation}\label{claim:the cone integer points} 
\K_v \cap \Z^d = \{  p + m_1 w_1 + \cdots + m_d w_d   \mid   p \in (\Pi+ v)\cap \Z^d, \text{ and all }  m_j  \in \Z_{\geq 0}\}.
\end{equation}
To prove \eqref{claim:the cone integer points}, we begin by writing each $m\in  \K_v \cap \Z^d$, by definition of the cone $\K_v$,
 as follows:
\[
m = v+ \lambda_1 w_1 + \cdots + \lambda_d w_d,
\]
 with the $\lambda_k \geq 0$.  This representation of $m$ is unique,
because $w_1, \dots, w_d$ is a basis for $\R^d$.   Now we use the fact that each
$\lambda_k = \lfloor \lambda_k \rfloor + \{ \lambda_k \}$, where $\{x\}$ is the fractional part of $x$:
\begin{align*}
m &= v+  \Big(  \{ \lambda_1 \} w_1 + \cdots + \{ \lambda_d \} w_d \Big)  +  
\lfloor \lambda_1 \rfloor w_1 + \cdots + \lfloor \lambda_d \rfloor w_d \\
&:=p  +  \lfloor \lambda_1 \rfloor w_1 + \cdots + \lfloor \lambda_d \rfloor w_d,
 \end{align*}
where we've defined $p:= v+  \Big(  \{ \lambda_1 \} w_1 + \cdots + \{ \lambda_d \} w_d \Big)$.  We now notice that
$p \in v+ \Pi$, and in fact $p \in \Z^d$, because $m, w_1, \dots, w_d  \in \Z^d$.

Since $\Pi$ tiles the cone $\K_v$ precisely by the translation vectors $w_1, \dots, w_d$, we  see 
that the set of all integer points in $\K_v$ is precisely 
the disjoint union of the sets 
\begin{equation} \label{typical integer points in the cone}
\{ p + k_1 w_1 + \dots + k_d w_d \mid  \, k_1, \dots, k_d \in \Z_{\geq 0}   \} 
\end{equation}
 (which we may think of as `multidimensional arithmetic progressions') , as $p$ varies over the integer points of $\Pi$ .
Finally, we expand each denominator in the following rational function as a geometric series to get:
\[
\frac{   \sigma_{ \Pi + v}(q)     }
{       \prod_{j=1}^d  \left(    1 - q^{ w_j}     \right)      } = 
\left(
\sum_{p \in (\Pi + v) \cap \Z^d}  q^p
\right)
\left(
\sum_{k_1 \geq 0} q^{k_1 w_1} 
\right)
\cdots
\left(
\sum_{k_d \geq 0} q^{k_d w_d} 
\right).
\]
Multiplying out all of these geometric series together, we see that the exponents look precisely like
the points in \eqref{typical integer points in the cone}. 
\end{proof}



\bigskip
\section*{Notes} \label{Notes.chapter.Brion}
 \begin{enumerate}[(a)]
\item    In the development of our text so far, we've observed 
that the discrete version of Brion's theorem 
(Theorem \ref{brion, discrete form}) followed from the continuous version of Brion's theorem 
(Theorem \ref{brion, continuous form}).   The tool we used in order to discretize Theorem \ref{brion, continuous form}  was the Poisson summation formula. 
 By contrast, the ideas in our previous book \cite{BeckRobins} developed in exactly the opposite direction: in that context we first proved the discrete Brion theorem, and then derived the continuous version from it. 
\end{enumerate}

\bigskip

\section*{Exercises}
\addcontentsline{toc}{section}{Exercises}
\markright{Exercises}

\medskip
\begin{prob}  \label{unimodular cone, integer point transform}
{\rm
Suppose that $\P\subset \R^d$ is a unimodular polytope, with vertex set $V$.   Using Theorem \ref{brion, discrete form}, show that its integer point transform is equal to:
\begin{equation}
\sigma_\P(z) = \sum_{v \in V} 
\frac{e^{\langle v, z \rangle} }{\prod_{k=1}^d \left( 1 - e^{\langle w_k, z\rangle}  \right) }.
\end{equation}
}
\end{prob}

\medskip
\begin{prob}  
{\rm
Fix a positive integer $m > 1$, and let $\P$ be the $2$-dimensional triangle 
whose vertices are given by 
$(0, 0),  (0, 1)$, and $(m, 0)$.   First compute the integer point transforms
$\sigma_{K_v}(z)$ for its three vertex tangent cones, and then compute the integer point transform $\sigma_\P(z)$. 
}
\end{prob}

\medskip
\begin{prob}  
{\rm
Given a positive integer $N$, find the integer point transform $\sigma_P(z)$ for the 
$2$-dimensional cone whose edge vectors are 
\[
w_1:= \icol{ 1 \\ 0}, \ w_2:=  \icol{1 \\ N}.
\]
}
\end{prob}

\medskip
\begin{prob}  
{\rm
More generally, given any coprime positive integer $p, q$, find the integer point 
transform $\sigma_P(z)$ for the 
$2$-dimensional cone whose edge vectors are 
\[
w_1:= \icol{ 1 \\ 0}, \ w_2:=  \icol{p \\ q}.
\]
}
\end{prob}

\medskip
\begin{prob}  \label{bound for integer point transform}
{\rm
Let $\P \subset \R^d$ be the $d$-dimensional polytope.   
\begin{enumerate}[(a)]
\item Prove the following inequality for the integer point transform:
\[
\left |  \sigma_\P(x) \right | \leq  \left |\Z^d \cap \P \right |,
\]
for all $x \in \R^d$. 
\item Is it true that for all $x \in [0, 1)^d \setminus \{0\}$, we have
\[
\left |  \sigma_\P(x) \right |  < \left |\Z^d \cap \P \right | ?
\]
\end{enumerate}
}
\end{prob}

\medskip
\begin{prob} 
{\rm 
Let $\P$ be the $3$-dimensional simplex whose vertices are given by 
$(0, 0, 0),  (1, 1, 0), (1, 0, 1)$, and $(0, 1, 1)$.   Compute the integer point transforms
of its four vertex tangent cones $\sigma_{K_v}(z)$, and then compute the
 integer point transform of $\sigma_\P(z)$. 
 }
\end{prob}

\medskip
\begin{prob}  
{\rm
Suppose we are given a $2$-dimensional simplicial integer cone $\K \subset \R^2$, together with its dual cone $\K^*$. 
Is there a simple relationship between
the integer point transforms $\sigma_{\K}(z)$ and $\sigma_{\K^*}(z)$ in dimension $2$?
}
\end{prob}

Notes.  For this problem, it's worth thinking about the relationship between the edge vectors of the
 fundamental parallelepipeds for $\K$ and $\K^*$.


 \chapter{\blue{Counting integer points in polytopes -  the Ehrhart theory} }
 \label{Ehrhart theory}
 \index{Ehrhart theory}
 
 \begin{quote}                         
 ``How wonderful that we have met with a paradox.  Now we have some hope of making progress. ''

  -- Niels Bohr \index{Niels Bohr}
   \end{quote}

 \begin{wrapfigure}{R}{0.58\textwidth}
\centering
\includegraphics[width=0.62\textwidth]{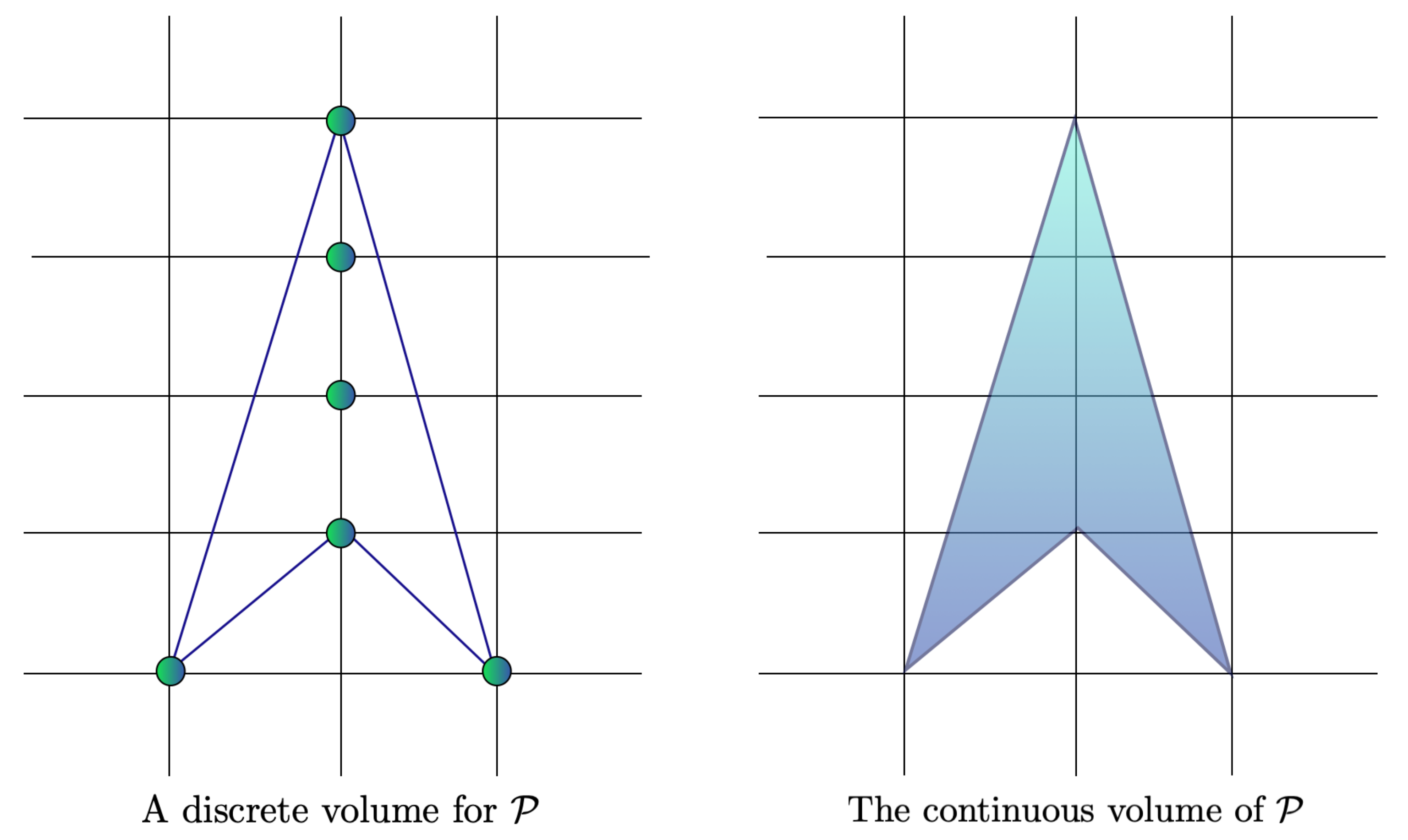}
\end{wrapfigure}

\red{(Under construction)}

 \section{Intuition}
 A basic question in discrete geometry is ``how do we discretize volume?" 

 One method of discretizing the volume of $\P$ is to count  the number of integer points in $\P$.
   Even in $\R^2$, this question may be highly non-trivial, depending on the arithmetic properties of the vertices of $\P$.   Ehrhart first considered integer dilations of a fixed, integer polytope $\P$, and 
studied the {\bf integer point enumerator}:
 \begin{equation} \label{combinatorial discrete volume} 
 L_{\P}(t):= | \Z^d \cap t\P |,
 \end{equation}
 where $t\P$ is the $t$'th dilate of $\P$, and $t$ is a positive integer.  Ehrhart showed
 that  $L_{\P}(t)$ is a polynomial in the positive integer parameter $t$,  known as the 
 {\bf Ehrhart polynomial} of $\P$.

 Viewed from the lens of Fourier analysis, Ehrhart polynomials may be computed by `averaging'
 the Fourier transform of a polytope over the full integer lattice:
 \begin{align}\label{Ehrhart is a discrete approximation to the volume}
 L_{\P}(t) &:= | \Z^d \cap t\P | = \sum_{n\in \Z^d} 1_{t\P}(n) = \sum_{\xi \in \Z^d}  \hat 1_{t\P}(\xi) \\
 &= \hat 1_{t\P}(0) + \sum_{\xi \in \Z^d \setminus \{0\}}  \hat 1_{t\P}(\xi)\\
 &=(\vol \P) t^d + \sum_{\xi \in \Z^d \setminus \{0\}}  \hat 1_{t\P}(\xi),
 \end{align}
where we've used Poisson summation in the third equality.  But because we may not use 
indicator functions directly in  Poisson summation, \index{Poisson summation formula}
 some care is required and the process of smoothing may be applied to $1_\P$. 
 
\noindent
\begin{figure}[htb]
\begin{center}
\includegraphics[totalheight=2.5in]{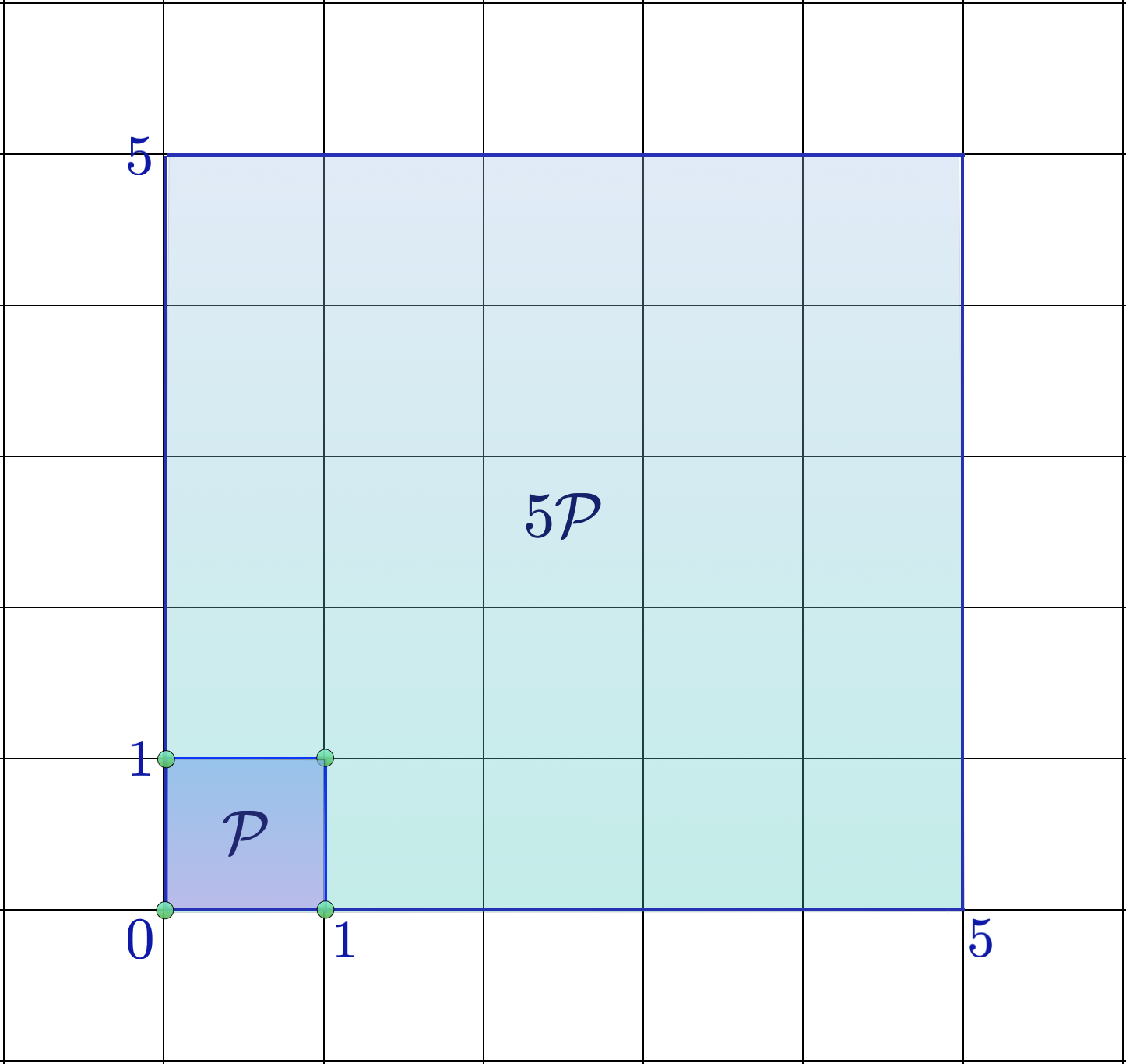}
\end{center}
\caption{Here the polytope $\P$ is the unit square, and we also have its $5$'th dilate $5\P$.}    \label{Ehrhart1}
\end{figure}

 As we've just seen in \eqref{Ehrhart is a discrete approximation to the volume}, the integer point enumerator $ | \Z^d \cap t\P |$, expanded using Poisson summation, has the primary term 
 $(\vol \P) t^d$. 
We recall the definition of the volume of $\P$, and of the Riemann integral:
\begin{align*}
\vol \P &:= \int_{\P} dx = \lim_{t \rightarrow \infty}   \frac{1}{t^d}  \sum_{n \in \P \cap \frac{1}{t} \Z^d} 1 = 
\lim_{t \rightarrow \infty} \frac{1}{t^d} \left | \P \cap \frac{1}{t} \Z^d \right |   \\
&= \lim_{t \rightarrow \infty} \frac{1}{t^d} \left | t\P \cap \Z^d \right |.
\end{align*}

More generally, given a function $f:\R^d \rightarrow \C$, we may sum the values of $f$ 
 at all integer points and observe how close this sum gets to the integral of $f$ over $\P$.  This approach is known as Euler-Maclaurin summation \index{Euler-Maclaurin summation}
 over polytopes, and is a current and exciting topic of a growing literature  (see Note \ref{EM summation note} below). 
  In this chapter we'll  also compare the above combinatorial method of discretizing volume, namely 
 \eqref{combinatorial discrete volume}, 
to our previous discrete volumes of Chapter~\ref{Angle polynomial}, which used solid angles.

 \bigskip
 \section{Computing integer points in polytopes via  the discrete Brion Theorem}
 
 \begin{example}\label{Ehrhart:the square}
 \rm{
 Probably the simplest example in $\R^2$ is the unit square $\P:= [0, 1]^2$.   
 As Figure \ref{Ehrhart1} suggests,
 the $t$-dilate $t\P$ here contains $(t+1)^2 = t^2 + 2t + 1$ points of the integer lattice $\Z^2$.   Here it was easy to conclude that
 $L_\P(t)$ was a polynomial function of $t \in \Z_{>0}$, but by a small miracle of nature 
a similar phenomenon occurs for \emph{all integer polytopes}  in $\R^d$.

\hfill $\square$
 }
 \end{example}
If all of the vertices of $\P$ have integer coordinates, we call $\P$ an {\bf integer polytope}.  
On the other hand, if all of the vertices of a polytope $\P$ have rational coordinates,  we call $\P$ a {\bf rational polytope}.
 \index{rational polytope}

 Let $\P \subset \R^d$ be a rational, $d$-dimensional  polytope, and let $N$ be the number of its vertices.
For each vertex $v$ of $\P$, we consider the vertex tangent cone $\K_v$  of $ \P$.   Once we dilate $\P$ by $t$, each
vertex $v$ of $\P$ gets dilated to become $tv$, and so each of the vertex tangent cones $\K_v$ of $\P$ simply get shifted to the 
corresponding vertex tangent cone $\K_{tv}$ of $t\P$.   Using the discrete Brion theorem (Theorem \ref{brion, discrete form}), we have
\begin{equation} \label{simplified discrete Brion identity}
  \sum_{n \in t\P \cap \Z^d}    e^{ \langle n,  z\rangle}  =        
   \sum_{n \in \K_{tv_1} \cap \Z^d}    e^{ \langle n,  z\rangle} + \cdots +
     \sum_{n \in \K_{tv_N} \cap \Z^d}    e^{\langle n,  z\rangle},
\end{equation}
for all $z \in \C^d - S$, where $S$ is the hyperplane arrangement defined by the (removable) singularities of all of the transforms $\hat 1_{\K_{v_j}}(z)$.
To simplify notation, we have absorbed the constant $-2\pi i$ into the complex vector $z$ by replacing
$z$ by $-\frac{1}{2\pi i } z$.
We recall that we rewrote \eqref{simplified discrete Brion identity} by using the notation:
\begin{equation}\label{discrete Brion, with rational functions}
\sigma_\P(z) = \sigma_{\K_{v_1}}(z)  + \cdots +  \sigma_{\K_{v_N}}(z).
\end{equation}
And now we notice that when $z = 0$,  the left-hand-side gives us precisely 
\[
  \sum_{n \in t\P \cap \Z^d}   1 := | \Z^d \cap t\P |, 
\]
which is good news - it is the Ehrhart polynomial  $L_{\P}(t)$, by definition.   The bad news is that  $z=0$ is a singularity of the right-hand-side of
 \eqref{discrete Brion, with rational functions}.   But  then again, there is  still more good news - we already saw in the previous chapter
  that it is a removable singularity.  So we may let  $z\rightarrow 0$, and discover what happens.

\bigskip
\begin{example} \label{Ehrhart poly for the standard triangle}
\rm{
Let's see by example how we can start with the discrete integer point transform, and end up with an Ehrhart polynomial. 
We will find a formula for the Ehrhart polynomial $L_{\P}(t) := | \Z^2 \cap t\P |$ of the standard triangle, continuing Example \ref{example:standard triangle integer point transform}.    It turns out
that the method we use in this example is universal - it can always be used to find the Ehrhart polynomial of any rational polytope.  We will formalize this method in the ensuing sections.

In this example we are lucky in that we may use brute-force to compute it, since the number of integer points in the $t$-dilate of $\P$ may be computed along the diagonals:   
\[
L_{\P}(t) = 1 + 2 + 3 + \cdots + (t+1) = \frac{(t+1)( t+2)}{2} = 
\frac{1}{2} t^2 + \frac{3}{2}t + 1.
\]
Now we can confirm this lucky answer with our brand new machine, as follows.
Using \eqref{simplified discrete Brion identity},
and the formulation  \eqref{last line} from 
Example \ref{example:standard triangle integer point transform}., 
we have the integer point transform for the dilates of $\P$:
\begin{align} 
& \sum_{n \in t\Delta \cap \Z^d}    e^{ \langle n,  z\rangle}  =      
            \sum_{n \in \K_{tv_1} \cap \Z^d}    e^{ \langle n,  z\rangle}  
       +   \sum_{n \in \K_{tv_1} \cap \Z^d}    e^{ \langle n,  z\rangle} 
       +   \sum_{n \in \K_{tv_3} \cap \Z^d}    e^{\langle n,  z\rangle} \\
         \label{integer point transform for standard triangle}   
 &= \frac{1}{(1- e^{z_1} )(1-  e^{z_2})} + \frac{e^{t z_1}}{(e^{-z_1}-1 ) ( e^{-z_1+z_2} -1)}
+\frac{e^{t z_2}}{(e^{-z_2}-1)(e^{z_1-z_2}-1)}  \\ \label{symmetric rationals}
& := F_1(z) + F_2(z) + F_3(z),      
\end{align}
where we have defined $F_1, F_2, F_3$ by the last equality.   We can let $z\rightarrow 0$ along almost any direction, but it turns out that we can simplify our computations by taking advantage of the symmetry of this polytope, so we will pick $z = \icol{ \ x \\ -x} $, which will simplify our computations
(see Note \ref{Michel Faleiros}).  Here is our plan: 
\begin{enumerate}[(a)]
\item  We pick $z:= \icol{ \ x\\ -x} $.
\item  We  expand all three meromorphic functions $F_1, F_2, F_ 3$ 
in terms of their Laurent series in $x$, giving us Bernoulli numbers.
\item  Finally, we let $x \rightarrow 0$, to retrieve the constant term (which will be
a  polynomial function of $t$) of the resulting Laurent series.
\end{enumerate}

To expand $F_1(z), F_2(z), F_3(z)$ in their Laurent series,  we recall the definition 
\ref{Def. of Bernoulli numbers}
 of the Bernoulli numbers in terms  of their generating function, namely 
$ \frac{t}{e^t-1}   =  \sum_{k =0}^\infty   B_k \frac{t^k}{k!}$:
\begin{align*}     
F_1(x, -x) &=  \frac{-1}{x^2} \sum_{m \geq 0}  B_m \frac{x^m}{m!}  
                                           \    \sum_{n \geq 0} B_n   \frac{(-x)^n}{n!} \\ 
&=   \frac{-1}{x^2}  \left( 1 - \frac{x}{2} + \frac{x^2}{12} + O(x^3)\right)     
                              \left( 1 + \frac{x}{2} + \frac{x^2}{12} + O(x^3)\right)   \\
&=   \frac{-1}{x^2} - \frac{1}{3} + O(x)                                                                     
\end{align*}

Similarly, we have
\begin{align*}     
F_2(x, -x) &=  \frac{1 + t x + \frac{t^2}{2!} x^2 + O(x^3) }{2x^2} 
                       \left(1 +    \frac{x}{2} + \frac{x^2}{12} + O(x^3) \right)   
                       \left(1 + \frac{(2x)}{2} + \frac{(2x)^2}{12} +  O(x^3) \right) \\
&= \frac{1}{2x^2} + \frac{3}{4x} + \frac{2}{3} + \frac{t}{2x} + \frac{3t}{4} + \frac{t^2}{4} +  O(x)                                                                   
\end{align*}

Now, by symmetry we see that $F_3(x,-x) = F_2(-x, x)$, so that by \eqref{symmetric rationals} 
and the latter expansions, we finally have:
\[
\sum_{n \in t\Delta \cap \Z^d}  e^{ \langle n,  \icol{ \ x\\ -x}  \rangle}   = 
    F_1(x, -x) + F_2(x, -x) + F_2(-x,  x) = 1 + \frac{3}{2} t + \frac{1}{2} t^2 + O(x).
\]
Letting $z:= \icol{ \ x \\ -x}  \rightarrow 0$ in the latter 
computation, we retrieve the (Ehrhart) polynomial: 
\[
\sum_{n \in t\Delta \cap \Z^d}  1 =  L_\Delta(t) = 1 + \frac{3}{2} t + \frac{1}{2} t^2, 
\]
as desired. 
}
\hfill $\square$
\end{example}


\bigskip
\section{The Ehrhart polynomial of an integer polytope, and the Ehrhart quasi-polynomial of a rational polytope}

\begin{wrapfigure}{R}{0.39\textwidth}
\centering
\includegraphics[width=0.30\textwidth]{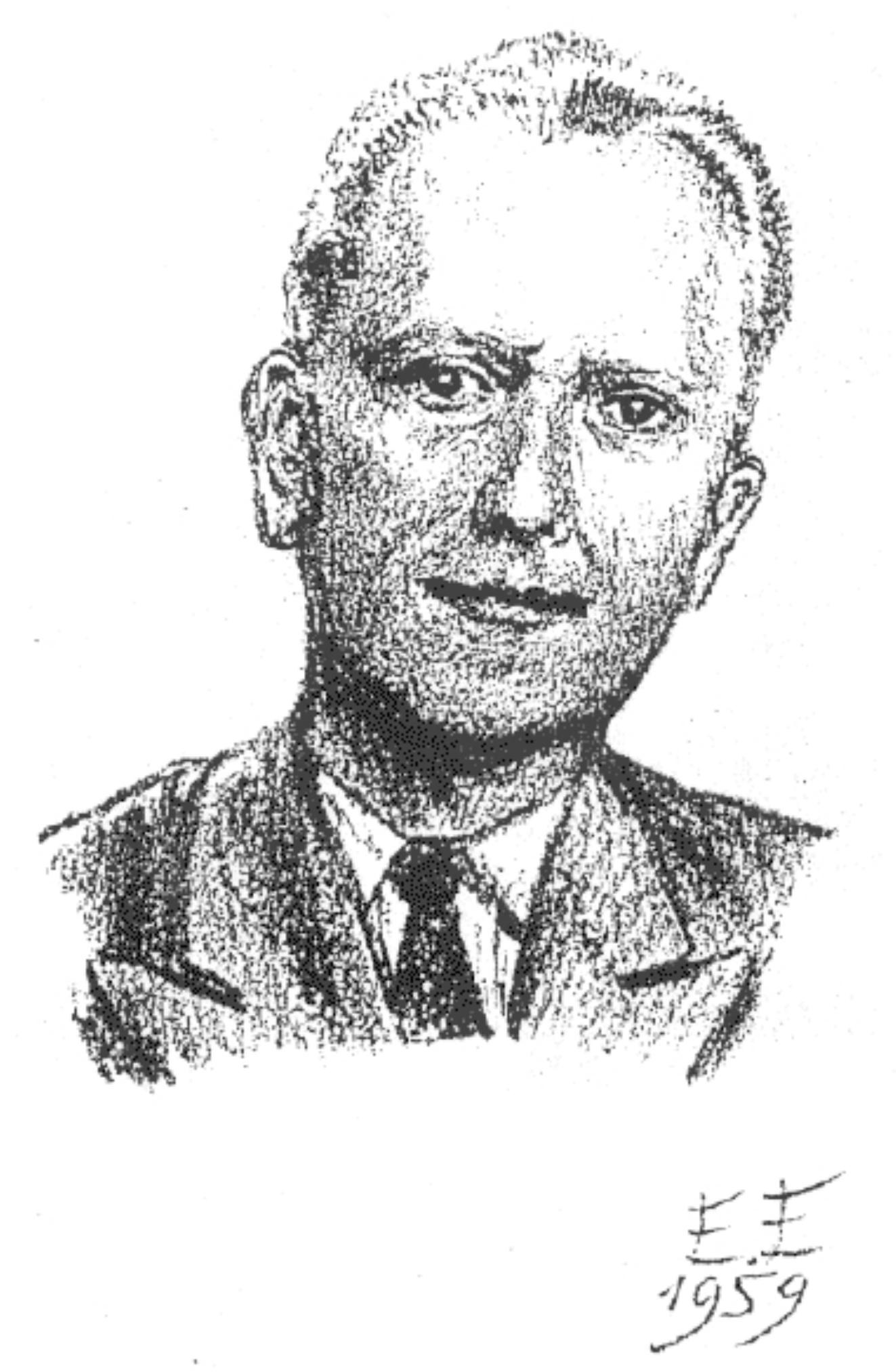}
\caption{Eugene Ehrhart, a self-portrait.}  
\label{Ehrhart}
\end{wrapfigure}

Eugene Ehrhart initiated a systematic study of the integer point enumerator
\[
L_\P(t):= \left| t\P \cap \Z^d  \right|,
\]
for an integer polytope $\P$, which Ehrhart proved was always a polynomial function of the positive integer dilation parameter $t$.
Ehrhart also proved 
that for a  rational polytope $\P \subset \R^d$, the integer point enumerator
$L_\P(t)$ is a {\bf quasi-polynomial} in the positive integer parameter $t$, which means by definition  that
\begin{equation}
 L_\P(t)  =  c_d t^d + c_{d-1}(t) t^{d-1} + \cdots + c_1(t) t + c_0(t),
\end{equation}
where each $c_j(t)$ is a periodic function of $t \in \Z_{>0}$.

The study of Ehrhart polynomials and Ehrhart quasi-polynomials has enjoyed a renaissance in recent years (\cite{BarvinokEhrhartbook}, \cite{BeckRobins}), and has some suprising connections to many branches of science, and even to voting theory, for example.

\begin{thm}[Ehrhart] \label{Ehrhart's main result}
For an integer polytope \\
$\P \subset \R^d$, its discrete volume $L_\P(t)$  
is a polynomial functions of $t$, for all positive integer values of the dilation parameter $t$.    Moreover, we have
\begin{equation}
L_\P(t) = (\vol \P) t^d + c_{d-1} t^{d-1} +  \cdots + c_1 t + 1.
\end{equation}
\hfill $\square$
\end{thm} 

Ehrhart's Theorem \ref{Ehrhart's main result} has an extension to rational polytopes, as follows. 
We will derive the more general Theorem  \ref{Ehrhart's  rational polytope theorem} of Ehrhart, by using 
the discrete Brion Theorem \ref{brion, discrete form}. 

\bigskip
\begin{thm}[Ehrhart] \label{Ehrhart's  rational polytope theorem}
For a rational polytope \\
$\P \subset \R^d$, its discrete volume $L_\P(t)$  
is a quasi-polynomial function of $t$, for all positive integer values of the dilation parameter $t$.    
In particular, we have
\begin{equation}
L_\P(t) = (\vol \P) t^d + c_{d-1}(t) t^{d-1} +  \cdots + c_1(t)  t + c_0(t),
\end{equation}
where each {\bf quasi-coefficient} $c_k(t)$ is a periodic function of $t\in \Z_{>0}$.
\end{thm} 

\noindent
\begin{figure}[htb]
\begin{center}
\includegraphics[totalheight=2.3in]{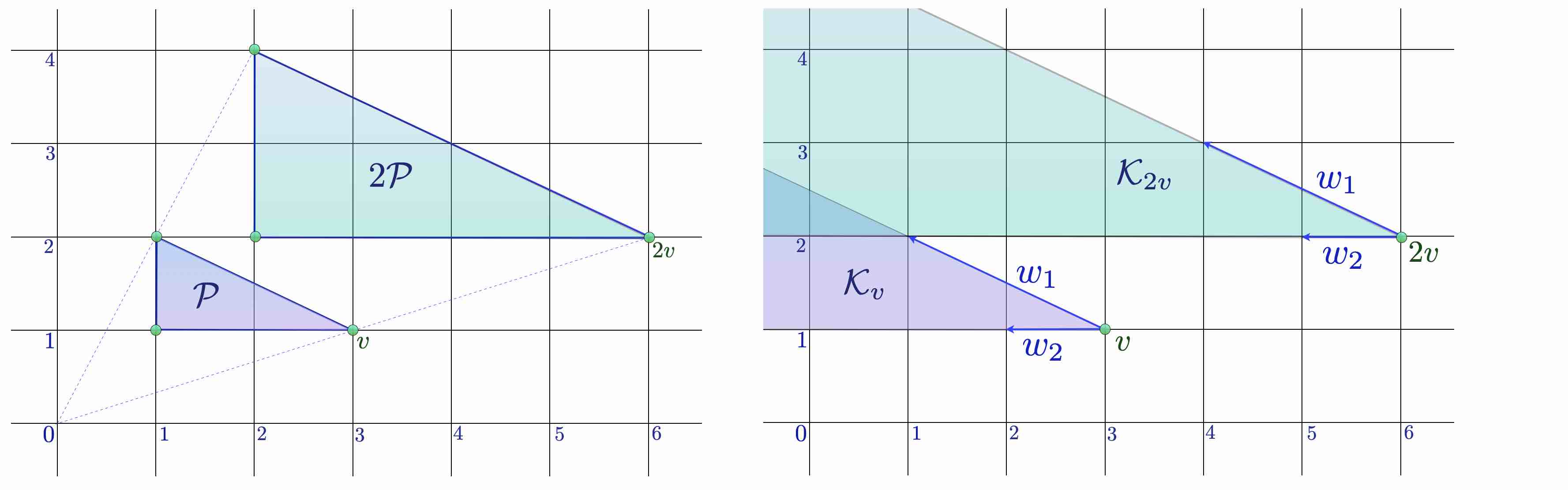}
\end{center}
\caption{Left: A triangle $\P$, and its dilate $2\P$.   \, Right: the vertex tangent cones $\K_v$ and $\K_{2v}$ have 
the same edge vectors  $w_1, w_2$. }    \label{Ehrhart2}
\end{figure}

\begin{proof} 
To begin, suppose that $p$  the least common denominator of the coordinates of all the rational vertices of $\P$. 
We need to show that, for each fixed $0\leq r < p$, the integer point enumerator $L_\P(r+ pk)$ is a polynomial in the parameter $k\in \Z_{>0}$. 
 By definition of a quasi-polynomial, this will prove that $L_\P(t)$ is a quasi-polynomial in $t \in \Z_{>0}$.  In other words, we restrict attention to each fixed arithmetic progression of dilations in $t:= r + pk$.  Now, from the discrete Brion Theorem \ref{brion, discrete form}, we know that 
\begin{equation}\label{Brion, for this proof} 
\sigma_{  \P}(z) =  \sigma_{   \K_{v_1}}(z) + \cdots  + \sigma_{   \K_{v_N}}(z),
\end{equation}
and we also know the elementary relation 
\[
\sigma_{  t\P}(0) := \sum_{n \in \Z^d \cap  t\P} 1 = L_\P(t).
\]
So we'd like to let $z \rightarrow 0$ on both sides of Brion's discrete identity \eqref{Brion, for this proof}:
\begin{equation}\label{first step of Ehrhart in terms of integer point transforms}
 L_\P(t) = \lim_{z\rightarrow 0} 
 \Big(
  \sigma_{   t\K_{v_1}}(z) + \cdots  + \sigma_{  t \K_{v_N}}(z)
\Big).
\end{equation}

The bad news is that
the right-hand-side of    \eqref{Brion, for this proof}
 introduces local singularities in the denominators of each 
rational-exponential function
 \[
\sigma_{   \K_{v_j}}(z) = \sum_{n \in \K_v \cap \Z^d}    e^{ \langle n,  z\rangle}.
\]
 But there is good news too! These singularities must be removable singularities. The reason is easy - 
$ \sigma_{\P}(z)$ is a finite sum of exponentials (by compactness of $\P$), and is therefore an analytic function of $z$, so any singularities on the right-hand side of \eqref{Brion, for this proof} must be removable singularities.
To proceed further, we'll begin by writing each vertex tangent cone $\K_v$ in terms of its vertex $v$, and edge vectors $w_j$:
\[
\K_v :=  \left\{ v + \sum_{j=1}^{M_v} \lambda_j w_j  \mid  \text{ all }  \lambda_j \geq 0  \right\}.
\]
Now we consider the dilates of $\K_v$ a bit more carefully, and we will use the fact that the edge 
vectors $w_1, \dots, w_{M_v}$ of any dilate of a vertex tangent cone $\K_{tv}$ remain invariant, as in Figure \ref{Ehrhart2}:
 \[
\K_{tv} :=  \left\{t v + \sum_{j=1}^{M_v} \lambda_j w_j  \mid  \text{ all }  \lambda_j \geq 0  \right\}.
\]

{\bf Case $1$.}  Suppose $r \not=0$.  Then:
\begin{align*}
t\K_v &:=  (r+pk)\K_v =  \left\{ (r+pk)v + \sum_{j=1}^{M_v} \lambda_j w_j  \mid  \text{ all }  \lambda_j \geq 0 \right\} \\
&= k(pv) +  \left\{ rv + \sum_{j=1}^{M_v} \lambda_j w_j  \mid  \text{ all }  \lambda_j \geq 0  \right\} \\
&= k(pv) + r \K_v.
\end{align*}
The salient feature of this computation is that $pv$ is an integer vector, by definition of $p$.  This implies that 
\begin{align*}
\sigma_{  t \K_{v}}(z)  &:=\sigma_{  (r+pk) \K_{v}}(z) :=  
            \sum_{n \in  \Big( k(pv) + r \K_v \Big) \cap \Z^d}    e^{ \langle n,  z\rangle} \\ 
            &= \sum_{m \in r \K_v  \cap \Z^d}    e^{ \langle  k(pv) + m,  z\rangle} \\
            &= e^{ \langle  k(pv),  z\rangle}    \sum_{m \in r \K_v  \cap \Z^d}    e^{ \langle  m,  z\rangle} \\
             &:= e^{ \langle  k(pv),  z\rangle}  \sigma_{rK_v}(z)  \\
\end{align*}
Summarizing, \eqref{first step of Ehrhart in terms of integer point transforms} gives us:
\[
 L_\P(t)= \lim_{z\rightarrow 0} 
\sum_{v \in V}  e^{ \langle  k(pv),  z\rangle}  \sigma_{rK_v}(z),   
\]
and giving a common denominator to all of the rational functions (of $e^{z_j}$)  $\sigma_{rK_v}(z)$, 
we may apply L'Hospital's rule a finite number of times.  Because the integer variable $k$ only appears in the exponents
$e^{ \langle  k(pv),  z\rangle}$, we see that each time we apply L'Hospital, an extra factor of $k$ comes down, giving us a polynomial function of $k$.

{\bf Case $2$.}  Suppose $r =0$.  Here the situation is slightly easier:  $t = pk$, so
\begin{align*}
t\K_v &:=  pk \K_v = \left\{ pk v + \sum_{j=1}^{M_v} \lambda_j w_j  \mid  \text{ all }  \lambda_j \geq 0 \right\} \\
&= k(pv) +  \left\{  \sum_{j=1}^{M_v} \lambda_j w_j  \mid  \text{ all }  \lambda_j \geq 0   \right\} \\
&= k(pv) + (\K_v - v),
\end{align*}
which is an integer cone because $pv$ is an integer vector, and $\K_v - v$ is an integer cone with apex at the origin.  
Similarly to the computation above, we have
\begin{align*}
\sigma_{  t \K_{v}}(z)  
          &= \sum_{n \in  \Big( kpv +    \left(\K_v-v \right) \Big) \cap \Z^d}    e^{ \langle n,  z\rangle}  
            = \sum_{m \in    \left(\K_v-v \right)   \cap \Z^d}    e^{ \langle  kpv + m,  z\rangle} \\
            &= e^{ \langle  kpv,  z\rangle}    \sum_{m \in  \left(\K_v-v \right)\cap \Z^d}    e^{ \langle  m,  z\rangle} \\
             &:= e^{ \langle  kpv,  z\rangle}  \sigma_{ \left(\K_v-v \right)}(z),  
\end{align*}
and the remaining steps are identital to Case $1$. 
\end{proof}

We note that for an integer polytope $\P$, the same proof gives us Theorem \ref{Ehrhart's main result}, namely
that  $L_\P(t)$ is a polynomial for positive integer dilations $t$;
here we just need Case $2$, with $t:=k$, so that $r=0$ and $p=1$. 

We emphasize again that one of the important steps in the latter computation was the fact that in both cases of the proof above, $k(pv)$ was an integer vector, allowing us to rewrite the integer point transform of the cone in a simpler way. 
As a first application of Theorem \ref{Ehrhart's main result}, we show that the discrete volume of a (half-open) parallelepiped has a particularly elegant and useful form.
\begin{lem}
Let $D$ be any half-open integer parallelepiped in $\R^d$, defined by
\[
D:= \left\{  
\lambda_1 w_1 + \cdots + \lambda_d w_d  \mid  0 \leq \lambda_1, \dots, \lambda_d < 1
\right\},
\]
where $w_1, \cdots w_d \in \Z^d$ are linearly independent.  Then:
\begin{equation}\label{second proof of parallelepiped fact}
\#\{ \Z^d \cap D\} = \vol D,
\end{equation}
and for each positive integer $t$, we also have 
\[
\#\{ \Z^d \cap tD\} = \left(\vol D\right)  t^d.
\] 
\end{lem}
\begin{proof}  We can tile $tD$ by using $t^d$ translates of $D$, because $D$ is half-open.  Therefore
\[
\#\{ \Z^d \cap tD\} = \#\{ \Z^d \cap D\} t^d,
\]
and by definition $\#\{ \Z^d \cap tD\} = L_D(t)$.   On the other hand, we also know by Ehrhart's Theorem \ref{Ehrhart's main result} that $L_D(t)$ is a polynomial for integer values of $t$, whose leading coefficient is $\vol D$.   Since $L_D(t)=\#\{ \Z^d \cap D\} t^d$ for all positive integer values of $t$, we conclude that
\[
\#\{ \Z^d \cap D\}= \vol D.
\]
\end{proof}
Although we've already proved \eqref{second proof of parallelepiped fact}
in Chapter \ref{chapter.lattices} on lattices (see Theorem \ref{sublattice index}, part
 \ref{third part of sublattice thm}), here we see a completely different approach to it, via Ehrhart theory.


\bigskip
\section{The Ehrhart series}

Given an integer $d$-dimensional polytope $\P \subset \R^d$, we now recall one of Ehrhart's tricks, 
which entails building a $(d+1)$-dimensional integer cone from all of the integer dilates of $\P$.

We first place a copy of $\P$ in the $x_{d+1}= 1$ hyperplane, so that this copy of $\P$ has the form
$\left\{  (u, 1) \mid u \in \P \right\}$.    If $v_1, \dots, v_N$ are the vertices of $\P$, we define a $(d+1)$-dimensional cone called Cone($\P$),  by using the edge vectors 
$w_1:= (v_k, 1), \dots, w_N:= (v_N, 1) \in \R^{d+1}$.   By construction, Cone($\P$) has apex at the origin, and is a 
pointed cone in $\R^{d+1}$.  Moreover, Cone($\P$) is naturally built up from all of the simultaneous dilations of $\P$, appearing as slices of Cone($\P$) parallel to the $x_1=0$ hyperplane, as in Figure \ref{Cone Over P}.   Cone($\P$) is sometimes called  {\bf the cone over $\P$} (and this process is sometimes called homogenizing a polytope). 

\begin{figure}[htb]
 \begin{center}
\includegraphics[totalheight=3.8in]{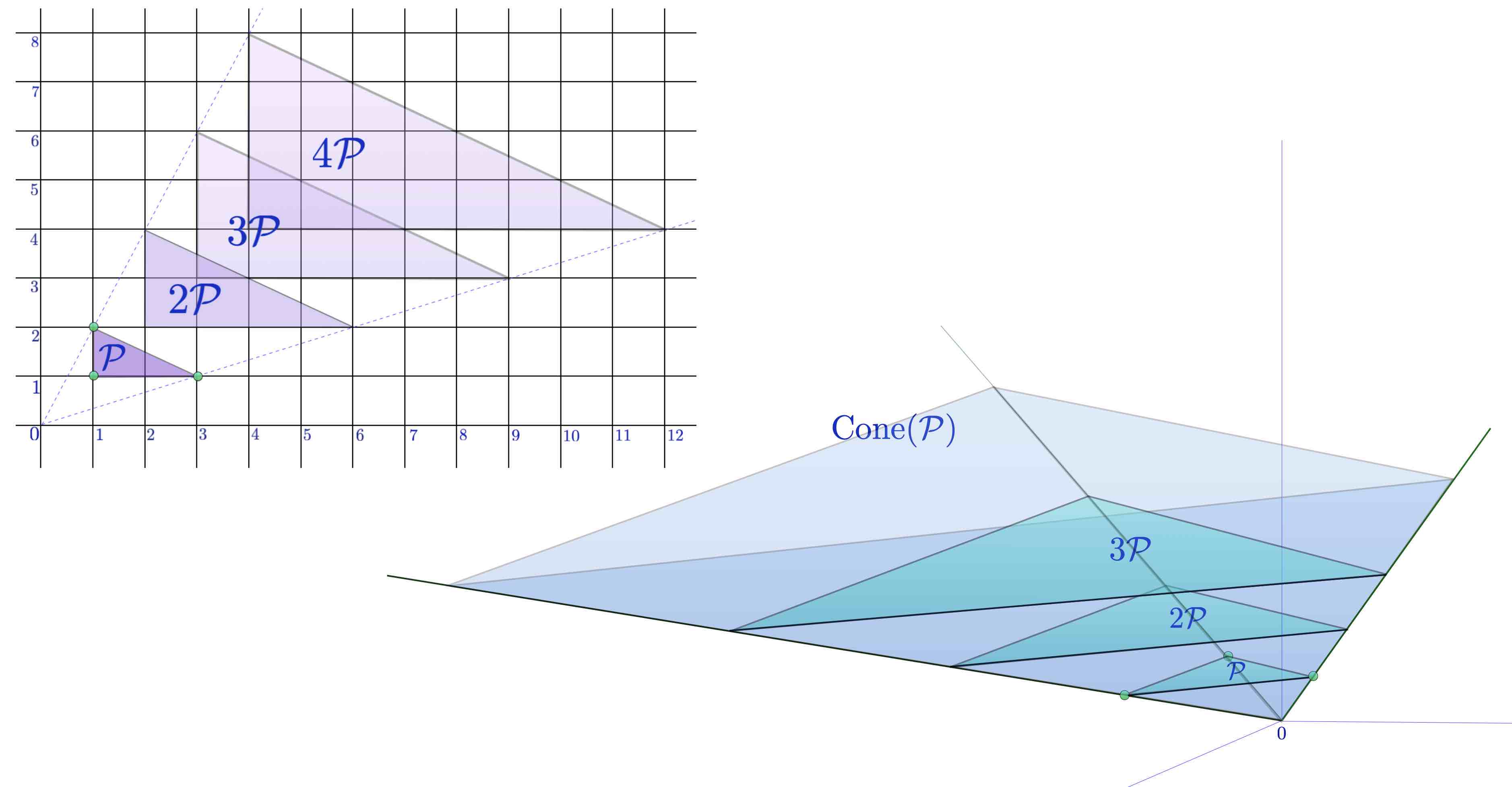}
\end{center}
\caption{Upper left:   a triangle $\P\subset \R^2$, with its dilates in $\R^2$.   Lower right:   a copy of $\P$ contained in the hyperplane $x_3 = 1$, with its dilates as slices of the $3$-dimensional cone called Cone($\P$).}  
\label{Cone Over P}
\end{figure}

Our next goal is to write down the integer point transform of Cone($\P$), so we first parametrize all of its integer points:
\[
\rm{Cone}(\P) \cap \Z^{d+1} = \left\{  \left(u, k \right) \mid u \in k\P\cap \Z^d, \text{ and } k= 0, 1, 2, 3, \dots               \right\},
\]
with the convention that when $k=0$, we define $0 \P$ (the zeroth dilate) to be the origin.  We recall that the integer point transform
of any cone $\K$ is defined by $\sigma_{\K}(z) :=    \sum_{n \in \K \cap \Z^d}    e^{2\pi i  \langle n,  z\rangle}$, but we will now specialize to 
$z:= (0, 0, \dots, 0, z_{d+1})$, and we also define $q:= e^{2\pi i z_{d+1}}$ to get the following
 special case of the integer point transform:
\begin{align*}
\sigma_{\rm{Cone}(\P)}(0, 0, \dots, 0, z_{d+1}) &:=    \sum_{n \in \K \cap \Z^d}    e^{2\pi i  \langle n,  (0, 0, \dots, 0, z_{d+1}) \rangle} \\
&=   \sum_{(u, k) \in  \Z^{d+1}  \atop  u \in k\P \cap \Z^d}    e^{2\pi i  \langle (u, k),  (0, 0, \dots, 0, z_{d+1}) \rangle} \\
&=   \sum_{(u, k) \in  \Z^{d+1}  \atop  u \in k\P \cap \Z^d}    e^{2\pi i  k z_{d+1} } \\
&=  1+  \sum_{k=1}^\infty  \left| k\P \cap \Z^d \right|  q^k = 1 + \sum_{k=1}^\infty L_\P(k)  q^k.
\end{align*}
The latter series is a generating function for the Ehrhart polynomial of $\P$, evaluated at all nonnegative integers, and it is by definition the {\bf Ehrhart series of } $\P$:
\[
{\rm Ehr}_\P(q):=  1 + \sum_{k=1}^\infty L_\P(k)  q^k.
\] 
One of the first and most famous Ehrhart-type theorems, discovered and proved by Richard Stanley \cite{StanleyDecompositions}, is the following characterization for the Ehrhart series of integer polytopes. 
\begin{thm}[Stanley]
\label{thm:Stanley nonnegativity}
Suppose that $\P  \subset \R^d$ is a $d$-dimensional integer polytope.  Then its Ehrhart series is given by
\[
{\rm Ehr}_\P(q) = \frac{\delta_d q^d + \delta_{d-1} q^{d-1} + \cdots + \delta_0}{(1-q)^{d+1}},
\]
and the coefficients $\delta_0, \dots, \delta_d$ are all nonnegative integers. 
\end{thm}
(For a proof see, for example,  \cite{BeckRobins} Theorem $3.12$)


\bigskip
\section{Families of Ehrhart polynomials}

\begin{example}[\bf{The unit cube}]
\rm{The simplest family consists of the $d$-dimensional unit cube in $\R^d$, defined by $P:=[0, 1]^d$.   This family extends Example
\ref{Ehrhart:the square}. 
 Here, the $t$-dilate $t\P$ has $(t+1)^d$ integer points, giving us the polynomial 
 \[
 L_\P(t)= (t+1)^d = t^d + d t^{d-1} + {d \choose 2} t^{d-2} + {d \choose 3} t^{d-3}  + \cdots + d t + 1,
 \]
for each positive integer value of $t$.
}
\hfill $\square$
\end{example}

\bigskip
\begin{example}[\bf{The standard simplex}]
\rm{
For the standard simplex
\index{standard simplex}
  $\Delta$, we consider its $t$-dilate, given by
\[
t\Delta := \{  (x_1, \dots, x_d) \in \R^d \mid \sum_{k=1}^d x_i \leq t, \text{ and  all } x_k \geq 0\}.
\]
We can quickly compute its Ehrhart polynomial by using combinatorics.  By definition, we need to find the number of nonnegative
 integer solutions to 
\[
x_1 + \cdots + x_d \leq t, 
\]
which is equal to $L_\Delta(t)$, for each fixed positive integer $t$.  We can introduce a `slack variable' 
\index{slack variable}
$z$, to transform the
latter inequality to an equality:  $x_1 + \cdots + x_d + z =  t$, where $ 0 \leq z \leq t$.   By a very classical and pretty argument, (involving placing $t$ balls into urns that are separated by $d$ walls) this number is equal to 
${t+d \choose d}$   (Exercise \ref{Ehrhart poly for closure of standard simplex}).
So we find that 
\begin{equation}\label{Ehrhart poly for closed standard simplex}
L_\Delta = {t+d \choose d} =  \frac{ (t+d) (t+d-1) \cdots (t+1)}{d!},
\end{equation}
a degree $d$ polynomial, valid for all positive integers $t$.  

What about the interior of $\Delta$?  Here we need to find the number of {\bf positive} integer solutions to 
$x_1 + \cdots + x_d <  t$, for each positive integer $t$.  
It turns out that by a very similar argument as above (Exercise \ref{Ehrhart poly for interior of standard simplex}), the number of positive integer solutions is ${t-1 \choose d} = L_{\interior \Delta}(t)$.  
So is it really true that 
\[
(-1)^d {d-t \choose d} = {t-1 \choose d}  \ ?
\]
Let's compute, substituting $-t$ for $t$ in    \eqref{Ehrhart poly for closed standard simplex}  to get: 
\begin{align*}
L_\Delta(-t) =  {-t+d \choose d} &=  \frac{ (-t+d) (-t+d-1) \cdots (-t+1)}{d!} \\
 &= (-1)^d\frac{ (t-d) (t-d+1) \cdots (t-1)}{d!} \\
 &=   (-1)^d  {t-1 \choose d} = (-1)^d L_{\interior \Delta}(t),
\end{align*}
confirming  that Ehrhart reciprocity works here as well. 
}
\hfill $\square$
\end{example}

\bigskip
\begin{example}[\bf{A Pyramid}]
\rm{
}
\hfill $\square$
\end{example}

\bigskip
\begin{example}[\bf{Zonotopes}]
\rm{
}
\hfill $\square$
\end{example}

\bigskip
\begin{example}[\bf{The crosspolytope}]
\rm{
}
\hfill $\square$
\end{example}

\bigskip
\begin{example}[\bf{The permutohedron}]
\rm{
We consider the vector $v_1:=(1, 2, 3, \dots, d)^T$, together with all of the $d!$ vectors whose coordinates are permutations of the coordinates of $v_1$, and we take their convex hull.   In other words, if we let $S_d$ be the symmetric group, 
then the {\bf Permutohedron} $\wp_d$ is defined by the convex hull
\[
\wp_d:= \conv  \{  \left(  \sigma(1), \sigma(2), \sigma(3), \dots, \sigma(d) \right)  \mid  \sigma \in S_d \},
\]
where $\sigma$ runs over all permutations in $S_d$.   It may be somewhat surprising to realize that the coefficients of the Ehrhart polynomial of $\wp_d$ count forests in graph theory.  To recall the definitions, a {\bf tree} is a connected graph (undirected) that does not contain any cycles.   A {\bf forest} is a disjoint union of trees.   
\begin{thm} If $\wp_d$ is the permutohedron, then
\[
L_{\wp_d}(t) = \sum_{k=0}^{d-1} f_k(d) t^k,
\]
\end{thm}
where $f_k(d)$ is the number of forests with $k$ edges, on the graph whose nodes are labelled $\{1, 2, \dots, d\}$.
}
\hfill $\square$

For example, when $d=3$, the Permutohedron is the convex hull of the six integer points:
\[
\wp_3:= 
\conv \left\{
\left(\begin{smallmatrix}
1 \\
2 \\
3
\end{smallmatrix}\right),
\left(\begin{smallmatrix}
1 \\
3 \\
2
\end{smallmatrix}\right),
\left(\begin{smallmatrix}
2 \\
1 \\
3
\end{smallmatrix}\right),
\left(\begin{smallmatrix}
2 \\
3 \\
1
\end{smallmatrix}\right),
\left(\begin{smallmatrix}
3 \\
1 \\
2
\end{smallmatrix}\right),
\left(\begin{smallmatrix}
3 \\
2 \\
1
\end{smallmatrix}\right)
\right\},
\]
which is a hexagon sitting in $\R^3$.

Add:   forests with $k$ edges on a graph with $3$ nodes are pictured  in Figure ...............
($k = 0, 1, $ or $2$ - ADD PICTURE HERE).
\end{example}

\begin{question}\label{classify all delta-vectors of polytopes}
Is it possible to somehow classify completely the vectors $(\delta_0, \dots, \delta_d)$ that occur as a
$\delta$-vector for some integer polytope $\P\subset \R^d$?
\end{question}
 
Question \ref{classify all delta-vectors of polytopes} is terribly difficult in general.  However, 
Alan Stapledon was able to classify these $\delta$-vectors of polytopes in low dimensions
 \cite{stapledonadditive}.


\bigskip
\section{Unimodular polytopes}

We recall that the standard basis vectors of $\R^d$ are defined by the set
$\{  e_1, e_2, \dots, e_d \}$, where $e_k$ is the unit vector pointing in the $x_k$ direction. 
A $d$-dimensional integer simplex $\Delta$ is called a {\bf unimodular simplex} if $\Delta$ is the unimodular image of the 
standard simplex $\Delta_{\rm standard}$,   \index{standard simplex}
which we recall is the convex hull  of the points  $\{ 0,    e_1, \dots,  e_d  \} \subset \R^d$.

\begin{figure}[htb]
 \begin{center}
\includegraphics[totalheight=2.7in]{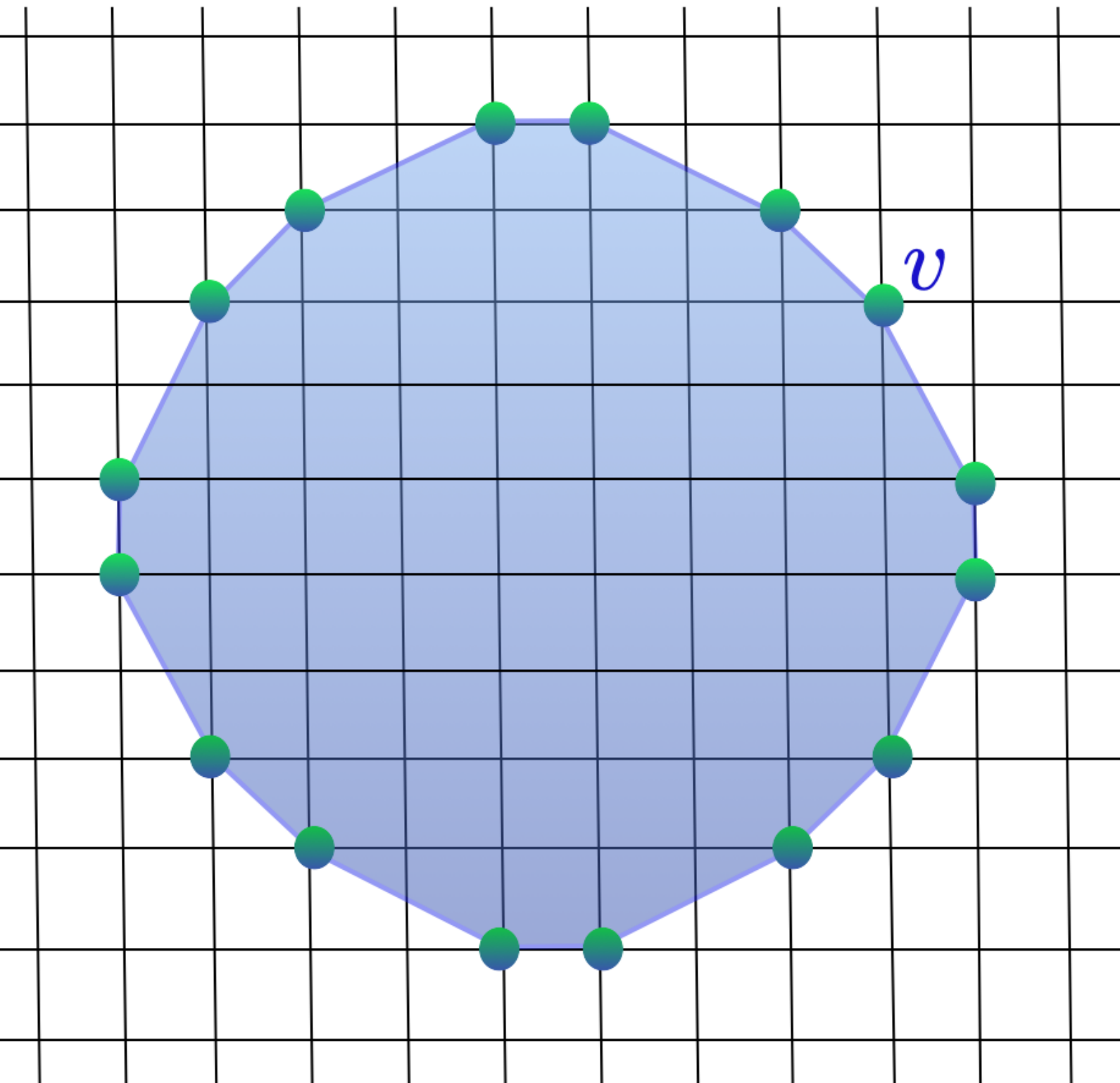}
\end{center}
\caption{A unimodular polygon - each vertex tangent cone is a unimodular cone.  It is clear from the construction in the Figure that we can form arbitrarily large unimodular polygons.}  
\label{unimodular polygon}
\end{figure}

\begin{figure}[htb]
 \begin{center}
\includegraphics[totalheight=2.3in]{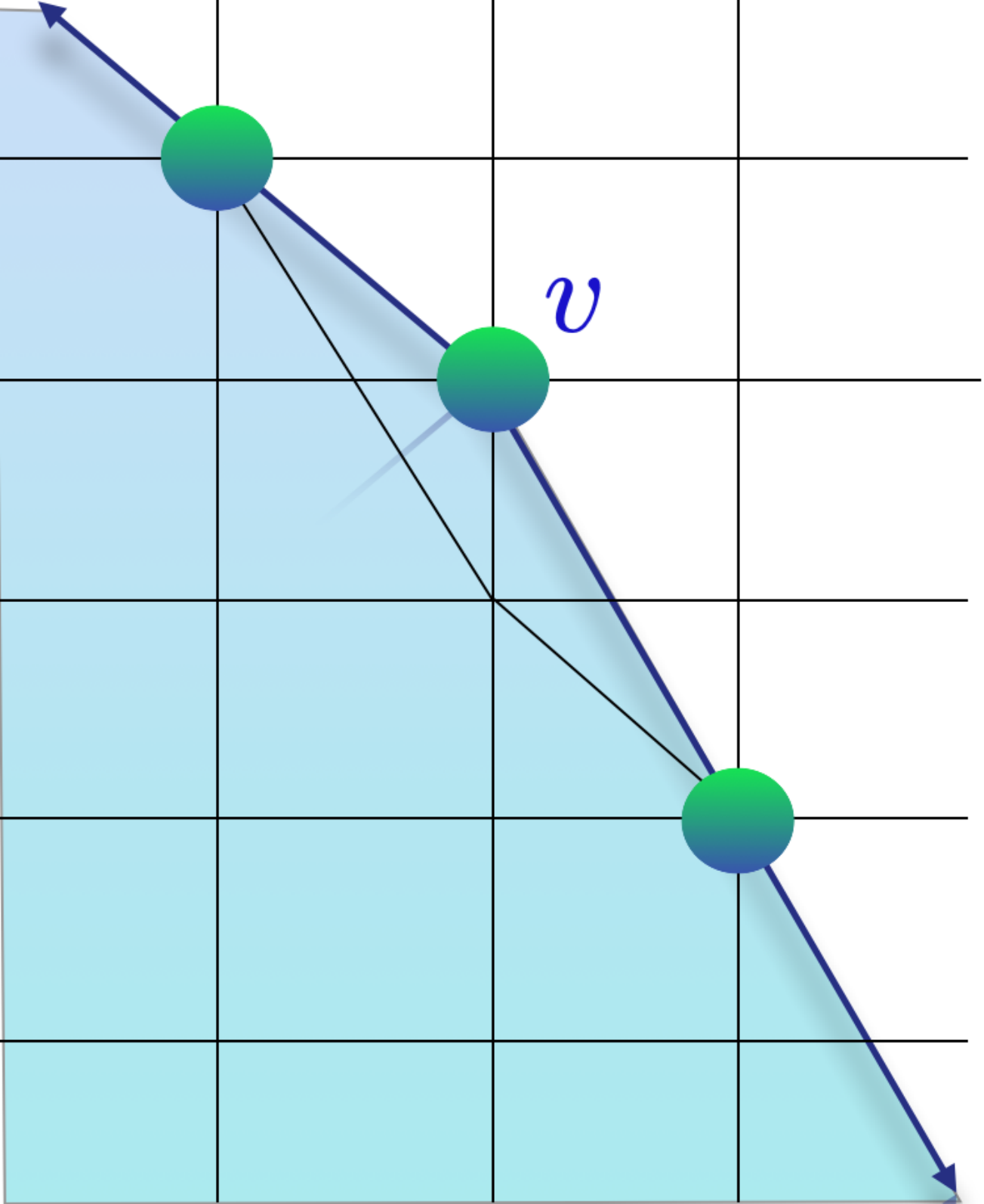}
\end{center}
\caption{A unimodular cone at $v$, appearing as one of the vertex tangent cones in the unimodular polygon of
Figure \ref{unimodular polygon}.  We notice that its half-open fundamental parallelepiped, with vertex
at $v$, does not contain any integer points other than $v$.}  \label{unimodular cone}
\end{figure}

\begin{example}
\rm{
Let $\Delta :=  {\rm conv}\left(   \icol{0\\0\\0},   \icol{1\\0\\0},  \icol{1\\1\\0},  \icol{1\\1 \\ 1}
\right)$, their convex hull.   Then $\Delta$ is a unimodular simplex, because the unimodular matrix 
$\left( 
\begin{smallmatrix}
1 & 1 & 1\\
0 & 1 & 1 \\
0 & 0 & 1 
\end{smallmatrix}
\right)$ maps the standard simplex $\Delta_{\rm standard}$ to $\Delta$. 
}
\hfill $\square$
\end{example}

It is not difficult to show that each of the $d+1$ tangent cones of a unimodular simplex possesses edge vectors that
form a lattice basis for $\Z^d$.   Thus, it is natural to define a {\bf unimodular cone} $\K \subset \R^d$ as a simplicial cone, possessing the additional property that its $d$ edge vectors form a lattice basis for $\Z^d$.

\begin{example}
\rm{
We consider the polygon $\P$ in Figure \ref{unimodular polygon}.  An easy verification shows that
each of its vertex tangent cones is unimodular.  For example, focusing on the vertex $v$, we 
see from Figure  \ref{unimodular cone}, that its vertex tangent cone is
  $\K_v:= v + \{  \lambda_1  \icol{ \ 1\\-2}   +  \lambda_2  \icol{-1\\ \ 1}    \mid 
\lambda_1, \lambda_2 \geq 0 \}$. 
$\K_v$ is a unimodular cone, because the matrix formed by the its two edges 
$\icol{ \  1\\ -2}$ and   $\icol{-1\\   \ 1}$ is  a unimodular matrix. 
}
\hfill $\square$
\end{example}

 More generally, a simple, integer polytope is called a 
{\bf unimodular polytope} if each of its vertex tangent cones is a unimodular cone. 
Unimodular polytopes are the first testing ground for many conjectures in discrete geometry and number theory.  Indeed, we will see later that the number of integer points in a unimodular polytope, namely
 $|\Z^d \cap \P|$, admits a simple and computable formula, if we are given the local tangent cone information at each vertex.  By contrast, it is in general thought to be quite difficult to compute the number of integer points $|\Z^d \cap \P|$, even for (general) simple polytopes, a problem that belongs to the NP-hard class of problems (if the dimension $d$ is not fixed).   The following fact, however, is elementary.

\begin{lem}
Suppose we have two integer polytopes $\P, \mathcal Q \subset \R^d$, which are unimodular images of each other:
\[
\P = U \mathcal Q,
\]
for some unimodular matrix $U$.
Then $L_\P(t) = L_{\mathcal Q}(t)$, for all $t \in \Z_{\geq 0}$.
\end{lem}

Next, we can generalize Exercise \ref{equivalent statements for unimodular triangles}, regarding primitive triangles, 
 to all integer simplices that intersect that integer lattice only in their vertices, as follows.
\bigskip
\begin{thm}\label{NEW: unimodular simplex equivalence}
Suppose that $\Delta \subset \R^d$ is a $d$-dimensional integer simplex.   Then  the following 
 properties are equivalent:
\begin{enumerate}[(a)]
\item  Aside from its integer vertices,  $(d-1)\Delta$ has no other integer points in its interior, or on its boundary.
\item  $\vol \Delta = \frac{1}{d!}$.
\item  $\Delta$ is a unimodular simplex.
\end{enumerate} 
\end{thm}

\section{More examples of rational polytopes and quasi-polynomials}

The following properties for the floor function, the ceiling function, and the fractional part function 
are often useful.
It's convenient to include the following indicator function, for the full set of integers, as well:
\[
 1_{\Z}(x) :=
 \begin{cases}
      1 & \text{if }       x \in \Z \\
     0  & \text{if }        x  \notin \Z \\
     \end{cases},
 \]
 the indicator function for $\Z$.  For all $x\in \R$, we have:
 \begin{enumerate}[(a)]
 \item    $\left\lceil x  \right\rceil  = - \floor{-x}$ \label{fractional part property a}
 \item   $1_{\Z}(x)=  \floor{x} - \left\lceil x  \right\rceil  +1$
\item $ \{ x \} + \{-x\} = 1- 1_{\Z}(x)$
\item $\floor{  x + y } \geq \floor{ x } + \floor{y}$, for all $x, y \in \R$.
\item Let $m \in \Z_{>0}, n \in \Z$.  Then $\floor{ \frac{n-1}{m} } + 1 = \left\lceil   \frac{n}{m} \right\rceil$.
\end{enumerate}
(Exercise  \ref{properties of floor, ceiling, fractional part})

\bigskip
\begin{example}
\rm{
Let's find the integer point enumerator $L_\P(t) := | \Z \cap t\P |$ of the rational line segment 
$ \P := [\frac{1}{3},  \ 1  ]$.   Proceeding by brute-force, for $t \in \Z_{>0}$ we have
\begin{align}
L_\P(t) &= \left|  \left[\frac{t}{3},  \ t  \right]  \cap \Z \  \right|    \label{answer comparison}   
  =\floor{t} - \left\lceil \frac{t}{3} \right\rceil     + 1  \\ 
  &=  t + \floor{   -\frac{t}{3}  }      +1 \\
  &= t +   -\frac{t}{3}    -  \left\{-\frac{t}{3}  \right\}  +1\\
    &= \frac{2}{3} t   -   \left\{   -\frac{t}{3}  \right\}  +1,
\end{align} 
a periodic function on $\Z$ with period $3$.  Here we used property \ref{fractional part property a} in the third equality.
In fact, here we may let $t$ be any positive real number,  and we still obtain the same answer, in this $1$-dimensional case.   

Now we will compare this to a new computation, but this time from the perspective of the vertex tangent cones.
For the cone $\K_{tv_1} := [\frac{t}{3}, + \infty)$, we can parametrize the integer points in this cone by
$\K_{tv_1}  \cap \Z =
    \{   \left\lceil \frac{t}{3} \right\rceil , \left\lceil \frac{t}{3} \right\rceil +1, \dots    \}$, so that
\begin{align*} 
\sigma_{\K_{tv_1}}(z) = e^{ \left\lceil \frac{t}{3} \right\rceil z} \sum_{ n \geq 0 } e^{n z} = 
e^{ \left\lceil \frac{t}{3} \right\rceil z} \frac{1}{1-e^{z}}.
\end{align*}

For the cone $\K_{tv_2} := (-\infty , t]$,  we can parametrize the integer points in this cone by 
$\K_{tv_2}  \cap \Z =  \{ t, t-1, \dots \}$, so that
\begin{align*}
\sigma_{\K_{tv_2}}(z) = e^{ t \cdot z} \sum_{ n \leq 0 } e^{n z} = 
e^{tz} \frac{1}{1-e^{-z}}.
\end{align*}

So by the discrete Brion Theorem (which is here essentially a finite geometric sum), we get:
\begin{align*}
\sum_{n \in [\frac{t}{3}, t] } e^{nz}  &= e^{ \left\lceil \frac{t}{3} \right\rceil z}
                                               \frac{1}{1-e^{z}}     +           e^{tz} \frac{1}{1-e^{-z}}  \\
& =-\left(1 +   \left\lceil \frac{t}{3} \right\rceil   z +   \left\lceil \frac{t}{3} \right\rceil^2 \frac{z^2}{2!} + \cdots    \right)
\left(\frac{1}{z}  -\frac{1}{2} + \frac{1}{12} z + \cdots  \right)  \\
    &+\left(     1 +  (t+1)z +  (t+1)^2 \frac{z^2}{2!} + \cdots    \right) 
    \left(\frac{1}{z}  -\frac{1}{2} + \frac{1}{12} z + \cdots \right)     \\ 
    &=  \frac{1}{2} - \left\lceil \frac{t}{3} \right\rceil  + (t+1) -\frac{1}{2}  +O(z)  \longrightarrow  \ 
        t - \left\lceil \frac{t}{3} \right\rceil +1,
\end{align*}
as $z\rightarrow 0$, recovering the same answer \ref{answer comparison} above.
}
\hfill $\square$
\end{example}

\bigskip
\begin{example}
\rm{
Let's find the integer point enumerator $L_\P(t) := | \Z^2 \cap t\P |$ of the rational triangle
\[
\P:= \conv\left(   \icol{0\\0},   \icol{\frac{1}{2} \\0},  \icol{0\\  \frac{1}{2}  } \right).
\]
First we will proceed by brute-force (which does not always work well), and then we will use the machinery of \eqref{simplified discrete Brion identity}.  

 For the brute-force method, we need to consider separately the even integer dilates and the odd integer dilates.  Letting $t=2n$ be a positive even integer, it's clear geometrically 
that 
\begin{align*}
L_\P(t) &:= | \Z^2 \cap 2n\P | =    1+ 2 + \cdots + n  \\
&=     \frac{n(n+1)}{2} = \frac{   \frac{t}{2}  ( \frac{t}{2} +1)}{2} \\
&= \frac{1}{8} t^2 + \frac{1}{4} t.
\end{align*}
 On the other hand, if $t= 2n-1$, then we 
notice that we never have an integer point on the diagonal face of $\P$, so that in this case
we get:
\[
L_\P(t) := | \Z^2 \cap (2n-1)\P | =  1+ 2 + \cdots + n = \frac{   \frac{t+1}{2}  ( \frac{t+1}{2} +1)}{2}
= \frac{1}{8}t^2 + \frac{1}{2} t + \frac{3}{8}.
\]

Alternatively, we may also rederive the same answer by using the Brion identity 
\eqref{simplified discrete Brion identity}.  We can proceed as in
Example \ref{Ehrhart poly for the standard triangle}.   The only difference now is that the vertex tangent cones have rational apices.   So although we may still use the same edge vectors to parametrize the integer points in 
 $\K_{tv_3} \cap \Z^d$,  we now have a new problem:
the rational vertex $v_3 = \icol{0\\  \frac{1}{2}}$. 
But in any case, we get:
$\K_{tv_3} \cap \Z^d= \left\{  
 \icol{0 \\  \frac{t}{2} }  +    n_1 \icol{ \ 0 \\ -1}   +     n_2\icol{ \ 1\\-1} \mid n_1, n_2 \in \Z_{\geq 0} \right\}$.
}
We invite the reader to complete this alternate derivation of the Ehrhart quasi-polynomial  $L_\P(t)$ in this case. 
\hfill $\square$
\end{example}

\begin{figure}[htb]
 \begin{center}
\includegraphics[totalheight=2.7in]{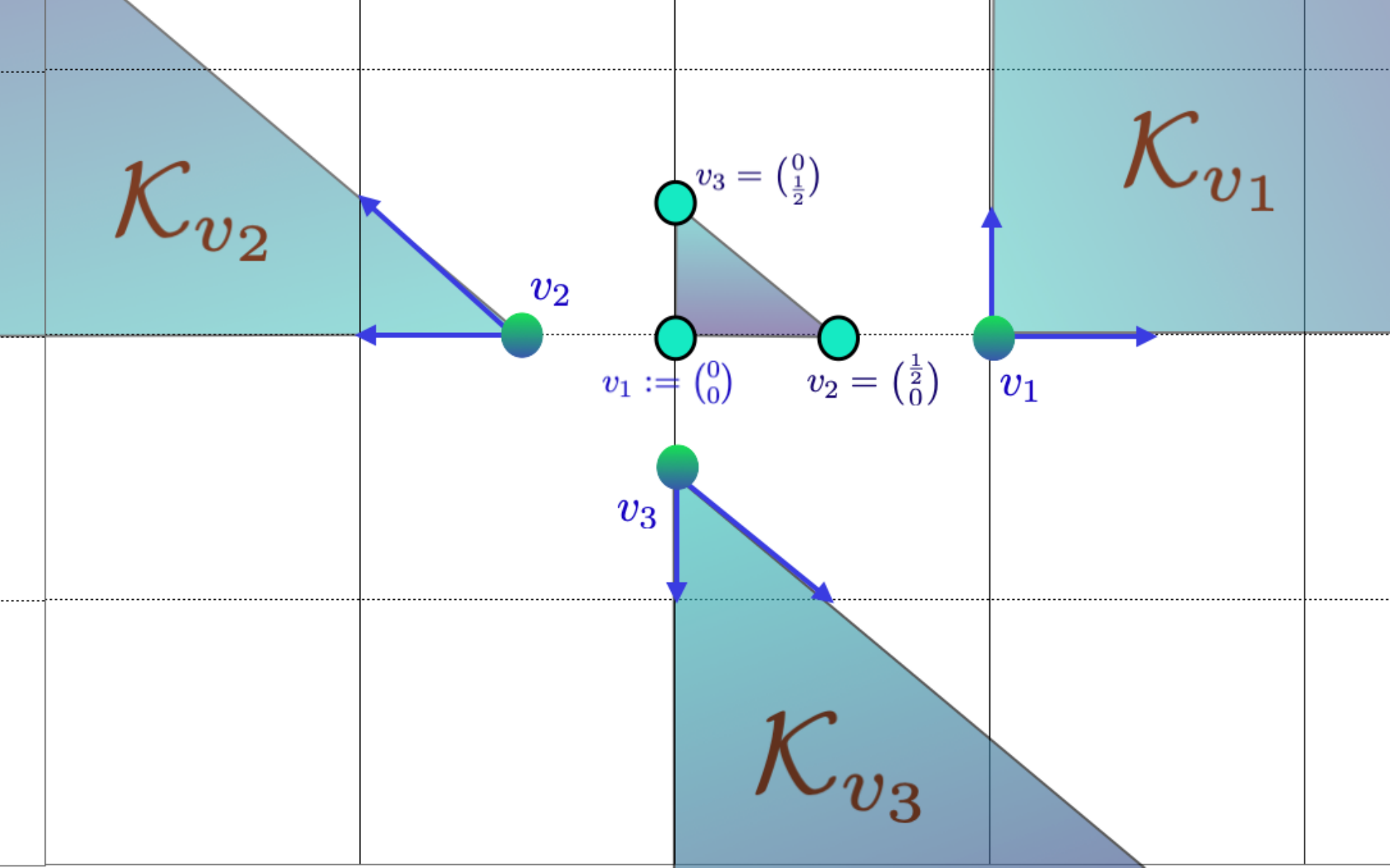}
\end{center}
\caption{A rational triangle, which happens to be a {\bf rational dilate} of the standard simplex.}  \label{rational triangle 1}
\end{figure}

\section{Ehrhart reciprocity} \index{Ehrhart reciprocity}

There is a wonderful, and somewhat mysterious, relation between the Ehrhart polynomial of the (closed) polytope $\P$, and the Ehrhart polynomial 
of its interior, called $\interior \P$.   We recall our convention that all polytopes are, by definition, closed polytopes.  
We first compute $L_\P(t)$, for positive integers $t$, and once we have this polynomial in $t$, we formally replace $t$ by $-t$.  So by definition, we form $L_\P(-t)$ algebraically, and then embark on a search for its new combinatorial meaning.

\bigskip \bigskip

\begin{thm}[Ehrhart reciprocity]
\label{thm:EhrhartReciprocity}
Given a $d$-dimensional rational polytope $\P\subset \R^d$, let 
$
L_{\interior \P}(t) := | \Z^d \cap \interior \P |, 
$
the integer point enumerator of its interior.    Then
\begin{equation}
 L_{ \P}(-t)=  (-1)^d L_{\interior \P}(t),
\end{equation}
for all $t\in \Z$.
\end{thm}

Offhand, this reciprocity law seems like `a kind of magic', and indeed Ehrhart reciprocity is one of the most elegant geometric inclusion-exclusion principles we have.  Some examples are in order.

\begin{example}
\rm{
For the unit cube $\square:= [0, 1]^d$, we can easily compute from first principles
 $L_\square(t) = (t+1)^d= \sum_{k=0}^d {d\choose k} t^k$.
For the open cube $\interior(\square)$ (the interior of $\square$), we can also easily compute 
\begin{align*}
L_{\interior(\square)}(t) &= (t-1)^d =   \sum_{k=0}^d {d\choose k} t^k(-1)^{d-k} \\
&=  (-1)^d \sum_{k=0}^d {d\choose k} (-t)^k \\
&=   (-1)^d L_\square(-t),
\end{align*}
using our known polynomial $L_\square(t) = (t+1)^d$.  
}
\hfill $\square$
\end{example}

A very common question in integer linear programming, as well as convex geometry, is:
\begin{question}  \label{integer point inside}
Given an integer polytope $\P\subset \R^d$, does it contain an integer point in its interior?
\end{question}
As a fun and  rapid consequence of Ehrhart reciprocity, we can give approach  
Question \ref{integer point inside}
in terms of the coefficients of the Ehrhart polynomial.
\begin{lem}
Let $\P\subset \R^d$ be an integer polytope, and let its Ehrhart polynomial be \\
$L_\P(t):= c_d t^d + c_{d-1}t^{d-1} + \cdots + c_1 t + 1$.  Then
\[
c_d - c_{d-1} + c_{d-2 } -  \cdots  +   (-1)^{d-1}c_1   + (-1)^d c_0  \geq 0,
\]
with equality $\iff  \P$ does not contain an integer point in its interior.
\end{lem}
 \begin{proof} By Ehrhart reciprocity, namely Theorem \ref{thm:EhrhartReciprocity}, we have
 \begin{align*}
 L_{\interior \P}(1)&= (-1)^d L_\P(-1)  \\
 &=  (-1)^d  \left(   c_d(-1)^d +  c_{d-1}(-1)^{d-1} + c_{d-2}(-1)^{d-2} +  \cdots  +  c_1(-1) + c_0  \right) \\
 &=  c_d -  c_{d-1} +  c_{d-2} -  \cdots  +  c_1(-1)^{d-1} + c_0(-1)^d.
 \end{align*}
But we also have   $L_{\interior \P}(1) \geq 0$, with equality if and only if 
there are no integer points in the interior of $\P$.
 \end{proof}
 As of this writing, the true complexity of finding even one integer point inside a given integer polytope $\P$ is not known, but this problem is known to be NP-hard.  Theorem \ref{thm:EhrhartReciprocity} also suggests that
  computing coefficients of Ehrhart polynomials appears to be hard in general. 


 \bigskip
\section{The M\"obius inversion formula for the face poset} 
\bigskip

Given a polytope  $\P\subset \R^d$,
the collection of all faces $F$ of $\P$ - including the empty set and $\P$ itself -  is ordered by inclusion.
This ordering forms a partially ordered set, and is called the {\bf face poset}. \index{face poset}
There is a particularly useful inversion formula on this face poset.
\begin{thm}[M\"obius inversion formula for the face poset]
\label{Mobius inversion}
\index{M\"obius inversion formula}
Given any function $g: \P\rightarrow \C$, we may define a sum over the face poset of $\P$:
\begin{equation}
h(\P) := \sum_{F\subseteq \P} g(F).
\end{equation}
We then have the following inversion formula:
\begin{equation} \label{Mobius inversion formula}
g(\P) := \sum_{F\subseteq \P} (-1)^{\dim F}  h(F).
\end{equation}
\hfill $\square$
\end{thm}
To prove (again) that for positive integer  values of $t$, the angle polynomial $A_\P(t)$ is indeed a polynomial in $t$, we 
may use the following useful  little relation between solid angle sums and integer point sums.  
We recall that for any polytope $\F$, the integer point enumerator for the relative interior of $\F$ was defined by   $L_{\interior  \F}(t):= | \Z^d \cap \interior \F |$.

For each face $\F \subseteq \P$, we define the $d$-dimensional {\bf solid angle of the face} $\F$ by picking any point $x$ inside the relative interior of $\F$ and denoting
\[
\om_\P(\F) := \om_\P(x). 
\]

\begin{thm}\label{lemma: relation between angle sum and integer sum}
Let $\P$ be a $d$-dimensional polytope in $\R^d$.  Then we have
\begin{equation}
A_\P(t) =  \sum_{\F \subseteq \P} \om_\P(\F) L_{\interior  \F}(t).
\end{equation}
\end{thm}
\begin{proof}
The polytope $\P$ is the disjoint union of its relatively open faces $\F \subseteq \P$, and similarly
the dilated polytope $t\P$ is the disjoint union of its relatively open faces $t\F \subseteq t\P$.   We therefore have:
\[
A_\P(t)=\sum_{n \in \Z^d} \om_{t\P}(n) =  \sum_{\F \subseteq\P} \sum_{n \in \Z^d}
\om_{t\P}(n)  1_{\interior(t\F)}(n).
\]
But by definition each $\om_{t\P}(n)$ is constant on the relatively open face $\interior(t\F)$ of $t\P$, and we denoted this constant by $\om_\P(\F)$.  Altogether, we have:
\[
A_\P(t)=  \sum_{\F \subseteq\P} \om_\P(\F) \sum_{n \in \Z^d}  1_{\interior(t\F)}(n)
 := \sum_{\F \subseteq \P} \om_\P(\F) L_{\interior  \F}(t).
\]
\end{proof}

\begin{thm}
Given an integer polytope $\P \subset \R^d$, the discrete volume $A_\P(t)$ is a polynomial in $t$, for integer values of the dilation parameter $t$.
\end{thm}
\begin{proof}
\ \ By Ehrhart's Theorem  \ref{Ehrhart's main result}, we know that for each face 
$\F~\subseteq~\P$, 
$L_{\interior  \F}(t)$ is a polynomial function of $t$, for positive integers $t$.  By 
Theorem \ref{lemma: relation between angle sum and integer sum}, we see that $A_\P(t)$ is a finite linear combination of polynomials, with constant coefficients, and is therefore a polynomial in $t$.
\end{proof}

We may apply Theorem \ref{Mobius inversion} to invert the 
relationship in Theorem 
\ref{lemma: relation between angle sum and integer sum}
 between solid angle sums and local Ehrhart polynomials, to get the following consequence of the M\"obius inversion formula. 

\begin{cor}\label{lemma: Mobius mu-function for angle sum and integer sum}
Let $\P \subset \R^d$ be a $d$-dimensional polytope.  Then we have
\begin{equation}
 L_{\interior  \P}(t) =  \sum_{\F \subseteq \P}       (-1)^{\dim F}  A_F(t).
\end{equation}
\end{cor}
\begin{proof}
We begin with the identity of 
Theorem \ref{lemma: relation between angle sum and integer sum}:
\begin{equation}
A_\P(t) =  \sum_{\F \subseteq \P} \om_\P(\F) \,   L_{\interior  \F}(t),
\end{equation}
and we use the M\"obius inversion formula \eqref{Mobius inversion formula}
 to get:
\begin{equation}
 \om_\P(\P)  L_{\interior  \P}(t) =  \sum_{\F \subseteq \P}   (-1)^{\dim F}    A_F(t).
\end{equation}
But $ \om_\P(\P)=1$, by definition, and so we are done.
\end{proof}

\begin{figure}[htb]
 \begin{center}
\includegraphics[totalheight=2.7in]{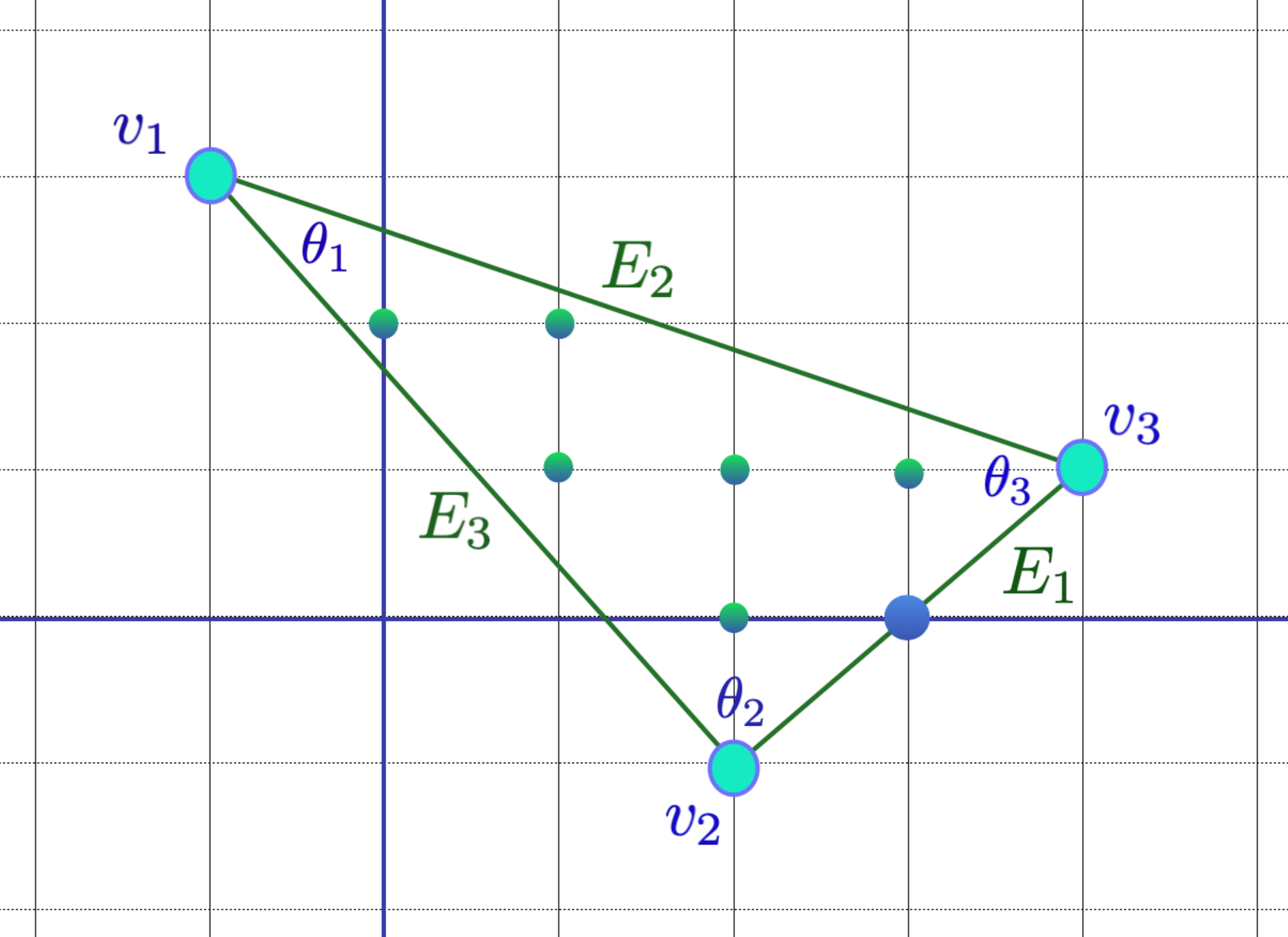}
\end{center}
\caption{An integer triangle, for which we compute $A_F(1)$ for each face $F\subset \P$ in Example \ref{Inverting a solid angle sum using Mobius inversion}, and use M\"obius inversion to find 
$L_{int \P}(1)$} 
\label{Triangle2}
\end{figure}

\begin{example}  \label{Inverting a solid angle sum using Mobius inversion}
\rm{
Let's work out a special case of Corollary \ref{lemma: Mobius mu-function for angle sum and integer sum},
in $\R^2$, for the triangle $\P$ appearing in Figure \ref{Triangle2}, with $t=1$.
$\P$ has vertices $v_1:= \icol{-1\\ \ 3}, v_2:=\icol{ \ 2\\ -1}, v_3:= \icol{4\\1}$, and edges $E_1, E_2, E_3$. 

We have to compute $A_F(1)$ for each face $F \subset \P$.  
At the vertices, we have $A_{v_1}(1) = \theta_1$,    $A_{v_2}(1) = \theta_2$, and $A_{v_3}(1) = \theta_3$.
For the edges of $\P$, we have:
\[
A_{E_1}(1) = \theta_{v_2} + \tfrac{1}{2} + \theta_{v_3},
\]
\[
A_{E_2}(1) = \theta_{v_3} +  \theta_{v_1},
\]
\[
A_{E_3}(1) = \theta_{v_1} +  \theta_{v_2}.
\]
Finally, for $\P$ itself, we have 
\[
A_\P(1) = 6 + \tfrac{1}{2} + \theta_{v_1} +  \theta_{v_2} + \theta_{v_3}
=   7.
\]
Putting everything together, we have:
\begin{align*}
   \sum_{\F \subseteq \P}       (-1)^{\dim F}  A_F(1) &=  
   \Big(    A_{v_1}(1) +  A_{v_2}(1) + A_{v_3}(1)  \Big) 
   -   \Big(    A_{E_1}(1) +  A_{E_2}(1) + A_{E_3}(1)  \Big) 
   +      A_{\P}(1)   \\
   &=   \Big(    \theta_{v_1} +  \theta_{v_2} + \theta_{v_3}    \Big)  
   -   \Big(   
   \theta_{v_2} + \tfrac{1}{2} + \theta_{v_3} +  \theta_{v_3} +  \theta_{v_1} +  \theta_{v_1} +  \theta_{v_2}             
       \Big)  +7 \\
   &=    \frac{1}{2} - \frac{3}{2} + 7 = 6 =  L_{\interior  \P}(1),
\end{align*}
the number of interior integer points in $\P$.
}
\hfill $\square$
\end{example}

 Finally, we mention a fascinating open problem by Ehrhart.
\index{Ehrhart conjecture}
 \begin{question}[Ehrhart, 1964]   \label{Ehrhart conjecture}
 Let $B \subset \R^d$  be a d-dimensional convex body with the origin as its
barycenter. If the origin is the only interior integer point in $B$, then
\[
\vol B \leq  \frac{    (d+1)d}{d!},
\]
and futhermore the equality holds if and only if $B$ is unimodularly equivalent to $(d + 1)\Delta$,
where $\Delta$ is the $d$-dimensional standard simplex.   
\index{standard simplex}
\end{question}
  
Ehrhart proved the upper bound for all $d$-dimensional simplices, and also for all convex bodies in dimension $2$.  But Question \ref{Ehrhart conjecture} remains open in general 
(see \cite{Nill.and.Paffenholz} for more details).

\bigskip \bigskip

\section*{Notes} \label{Notes.chapter.Ehrhart}
\begin{enumerate}[(a)]

\item  Ehrhart theory has a fascinating history, commencing with the fundamental work of Ehrhart \cite{Ehrhart1}, \cite{Ehrhart2}, \cite{Ehrhart3}, 
\cite{Ehrhartbook}, in the 1960's.    Danilov \cite{Danilov} made a strong contribution to the field, but after that
the field of Ehrhart theory lay more-or-less dormant, until it was rekindled by Jamie Pommersheim in 1993 \cite{Pommersheim}, giving it strong connections to Toric varieties.   Using the Todd operators to discretize certain volume deformations of polytopes, Khovanskii and Pukhlikov discovered a wonderful result that helped develop the theory further (see Theorem 12.6 of \cite{BeckRobins}).
In 1993, Alexander Barvinok \cite{Barvinok.algorithm} \index{Barvinok, Alexander}
gave the first polynomial-time algorithm for counting integer points in polytopes in fixed dimension.  

In recent years, Ehrhart theory has enjoyed an enthusiastic renaissance (for example, the books \cite{BarvinokEhrhartbook}, \cite{BeckRobins}, \cite{FultonBook}).  Early connections between Ehrhart theory and Fourier analysis appeared in \cite{diaz} and \cite{brionvergne}.  
 For more relations with combinatorics, the reader may enjoy reading 
Chapter $4$ of the classic book ``Enumerative Combinatorics'',  \cite{StanleyBook} by Richard Stanley. 

\item  Theorem \ref{thm:Stanley nonnegativity} represents one of the many Ehrhart-type results 
 due to Richard Stanley \cite{StanleyDecompositions}.   Stanley's original methods used graded rings in commutative algebra, forming exciting connections between this kind of combinatorial geometry and algebraic geometry \cite{StanleyCommutativeAlgebra}.

\item Regarding the computational complexity of counting integer points in polytopes, Alexander Barvinok settled the problem 
in \cite{Barvinok.algorithm} by showing that for a fixed dimension $d$, there is a polynomial-time algorithm, as a function of the
 `bit capacity' of any given
rational polytope $\P \subset \R^d$, for counting the number of integer points in $\P$.

\item It is also true that for integer polytopes which are not necessarily convex (for example simplicial complexes), 
  the integer point enumerator makes sense as well.  In this more general context, 
the  constant term of the corresponding integer point enumerator equals the (reduced) Euler characteristic of the simplicial complex.

\item   For more information about the rapidly expanding field of Euler-MacLaurin summation over polytopes,  a brief (and by no means complete) list of paper in this direction consists of the work by Berligne and Vergne \cite{BerlineVergne}, Baldoni, Berline, and Vergne  \cite{BaldoniBerlineVergne}, Garoufalidis and Pommersheim  \cite{GaroufalidisPommersheim},   Brandolini, Colzani, Travaglini, and Robins          \cite{BrandoliniColzaniTravagliniRobins2}, 
 Karshon, Sternberg, and Weitsman (\cite{KarshonSternbergWeitsman1}, \cite{KarshonSternbergWeitsman2}), and very recently Fischer and Pommersheim \cite{FischerPommersheim}.
 
\item  There are some fascinating relations between an  integer polytope $\P$ and its polar polytope $\P^o$.  In particular, let $\P \subset \R^2$ be an  integer polygon (convex) whose only interior integer point is the origin. 
Such polygons are called reflexive polygons, and up to unimodular transformations there are only a finite number of them in each dimension.    If we let $B(\P)$ be the number of integer points on the boundary of $\P$, then Bjorn Poonen and Fernando Villegas proved \cite{PoonenVillegas} that
\[
B(\P) + B(\P^o) = 12.
\]
One way to see why we get the number ``12'' is to consider Bernoulli numbers and Dedekind sums, but in  \cite{PoonenVillegas} the authors give 4 different proofs, including Toric varieties and modular forms. 

\item  The book \cite{BeckSanyal} by Matthias Beck and Raman Sanyal covers many classical  instances of combinatorial reciprocity that appear in combinatorial geometry, including the reciprocity for the order polynomial of a poset.

\item  \label{Michel Faleiros}  The trick used in Example \ref{Ehrhart poly for the standard triangle} 
of picking the particular vector $z := (x, -x)$, which turns out to simplify the computations a lot, is due to Michel Faleiros.  

\item There is a fascinating theory that offers an abstract extension of many of the ideas in this chapter, called valuation theory.  The author may wish to consult the excellent introduction to this field, by Katharina Jochemko \cite{KatharinaJochemko}. Allowing the empty set to be included in the collection of all convex bodies, we may define a {\bf valuation} as any mapping
\[
\Phi: \{ \text{convex bodies } \P \subset \R^d\} \rightarrow G,
\]
where $G$ is an abelian group (think of  $G:=\C$), with the following properties:
\[
\Phi(\varphi) = 0, \text{ and } \Phi(\P \cup Q) = \Phi(\P) + \Phi(Q) - \Phi( \P\cap Q),
\]
for all $\P, Q$ such that  $\P \cup Q$ and $\P \cap Q$ are also convex bodies.  The volume of a convex body, as well as the Ehrhart polynomial of an integer polytope, are just two examples of valuations.  Two of the founders of this theory are Peter McMullen \cite{McMullen1} and Jim Lawrence \cite{LawrenceVolume2}.

\item  \label{EM summation note}
In a future version of this book, we will also delve into Dedekind sums, which arise very naturally when considering the Fourier series of certain rational-exponential functions.  To define a general version of these sums, let
$\L$ be a $d$-dimensional lattice in $\R^d$,  let $w_1, \dots, w_d$ be linearly independent vectors from 
$\L^*$, and let $W$  be a matrix with the $w_j$'s  as columns. For any $d$-tuple $e = (e_1, \dots, e_d)$ of positive integers $e_j$, 
define $|e| := \sum_{j = 1}^k e_j$.   Then, for all $x \in \R^d$, a lattice Dedekind sum is defined by
\begin{equation}
L_\L(W, e; x) :=
\lim_{\varepsilon \rightarrow 0} 
\frac{1}{(2\pi i)^{|e|}} 
\sum_{\substack{\xi \in \L   \\ \langle w_j, \xi \rangle \neq 0, \forall j}} 
            \frac{e^{-2\pi i   \langle x, \xi \rangle   }}{\prod_{j = 1}^k \langle w_j, \xi \rangle^{e_j}} e^{-\pi \varepsilon \|\xi\|^2}.
\end{equation}
Gunnells and Sczech \cite{GunnellsSczech} have an interesting reduction theorem for these sums, 
  giving a polynomial-time complexity algorithm for them, for fixed dimension $d$. 
\end{enumerate}


\bigskip
\section*{Exercises}
\addcontentsline{toc}{section}{Exercises}
\markright{Exercises}

\begin{quote}    
If there is a problem you can't solve, then there is an easier problem you can't solve: find it.

--  George Polya  
 \end{quote}
\medskip

\medskip
\begin{prob} \label{exercise:warm up dual lattice}
In $\R$, consider the  $1$-dimensional polytope $\P:= [a,b]$, for any
$a,b \in \Z$.
\begin{enumerate}[(a)]
\item Show that the Ehrhart polynomial of $\P$ is  $L_\P(t) =  (b-a)t + 1$.
\item  Find the Ehrhart quasi-polynomial $L_\P(t)$
for the rational segment $\mathcal Q:= [\frac{1}{3}, \frac{1}{2}]$.
\end{enumerate}
\end{prob}

\medskip
\begin{prob}   
Fix positive integers $a, b$.
Working in $\R^2$, show that the closed line segment $\P \subset \R^2$, whose vertices are the origin and $(a, b)$,
contains exactly $\gcd(a, b) + 1$ integer points of $\Z^2$.   Conclude that we have the lower-dimensional
 Ehrhart polynomial  $L_\P(t) = \gcd(a, b) t + 1$.
\end{prob}

\medskip
\begin{prob} \label{2-d cross polytope Ehrhart}
We recall that the $d$-dimensional cross-polytope was defined by 
\[
\Diamond:=\left\{ \left( x_1, x_2, \dots, x_d \right) \in \R^d \mid
 \, \left| x_1 \right| + \left| x_2 \right| + \cdots +  \left| x_d \right| \leq 1 \right\}.
\]
For $d=2$, find the Ehrhart polynomial $L_\Diamond(t)$.
\end{prob}

\medskip
\begin{prob}
Extending Exercise \ref{2-d cross polytope Ehrhart}, show that the Ehrhart polynomial of $\Diamond$ in $\R^d$ is 
\[
L_{\Diamond}(t) = \sum_{k=0}^d \binom{d}{k} \binom{t-k+d}{d},
\]
for all $t \in \Z_{>0}$.
\end{prob}

\medskip
\begin{prob}
 Let $d=2$, and consider the cross-polytope $\Diamond \subset \R^2$.   Find the Ehrhart quasi-polynomial $L_\P(t)$ for the rational polygon $\P := \frac{1}{2} \Diamond$.
\end{prob}

\medskip
\begin{prob}
Suppose $\Delta$ is the standard simplex in $\R^d$.  Show that the first $d$ dilations of $\Delta$ do not contain
any integer points in their interior: 
\[
t(\interior \Delta) \cap \Z^d = \phi,
\]
for $t = 1, 2, \dots, d$.  
In other words, show that $L_{\interior \P}(1) =  L_{\interior \P}(2) = \cdots = 
L_{\interior \P}(d) = 0$.  Conclude that the same statement is true for any unimodular simplex.
\end{prob}

 \medskip
\begin{prob} \label{Bernoulli polynomial as an Ehrhart polynomial}
Here we show that the Bernoulli polynomial $B_d(t)$, is essentially equal to the Ehrhart polynomial
$L_\P(t)$ for the  ``Pyramid over a cube" (as defined in Exercise
\ref{Pyramid over a square}).  We recall the definition: let $C:=[0,1]^{d-1}$ be the $d-1$-dimensional cube,
considered as a subset of $\R^d$, and let ${\bf e_d}$ be the unit vector pointing in the $x_d$-direction.
Now we define $\P:=  \conv\{   C,   {\bf e_d} \}$, a pyramid over the unit cube.
Show that its Ehrhart polynomial is
\[
L_{\P}(t) = \frac{1}{d}(B_d(t+2) - B_d),
\]
for $t \in \Z_{>0}$.
\end{prob}

 \medskip
 \begin{prob}    \label{Pick's formula, generalization to d dimensions}
 For any integer $d$-dimensional (convex) polytope $\P \subset \R^d$, show that 
 \begin{equation} \label{Volume in terms of forward differences}
 \vol \P = \frac{(-1)^d}{d!} \left(  1 +  \sum_{k=1}^d {d\choose k} (-1)^k L_\P(k) \right),
 \end{equation}
 which can be thought of as yet another {\bf generalization of Pick's formula to $\R^d$}.  
 \index{Pick's formula, generalization}
 
 Note.  \  Using iterations of the forward difference operator 
 \[
 \Delta f(n):=  f(n+1) - f(n), 
 \]
 the latter identity
 may be thought of a {\bf discrete analogue} of the $d$'th derivative of the Ehrhart polynomial.  This idea in fact gives another method of proving \eqref{Volume in terms of forward differences}.  
  \end{prob}

 \medskip
 \begin{prob} \label{Pick's formula from general Ehrhart exercise}
 Show that Pick's formula is the special case of
 Exercise \ref{Pick's formula, generalization to d dimensions} when the dimension $d=2$.
 That is, given an integer polygon $\P \subset \R^2$, we have
 \[
 \rm{Area } \P =  I + \frac{1}{2} B -1,
\]
where $I$ is the number of interior integer points in $\P$, and B is the number of boundary integer points of $\P$. 
 \end{prob}

\medskip
 \begin{prob} \label{convolution of the indicator function with a Gaussian}
 Fix $\epsilon >0$.  Show that the convolution of the indicator function $1_\P$
 with the heat kernel $G_{\varepsilon}$, as in  equation \eqref{Gaussian smoothing},
 is a Schwartz function (of $x \in \R^d$). 
 \end{prob}

\medskip 
 \begin{prob} \label{unimodular triangle}
 Show that any unimodular triangle has area equal to $\frac{1}{2}$. 
 \end{prob}

\medskip
 \begin{prob} \label{unimodular triangle, Ehrhart poly}
 Show that  the Ehrhart polynomial of the standard simplex
 \index{standard simplex}
  $\Delta \subset \R^d$ is 
\[
L_{\Delta}(t) = \binom{t+d}{d}.
\] 
 \end{prob}

\medskip
 \begin{prob} 
 Consulting Figure  \ref{unimodular polygon}:
  \begin{enumerate}[(a)]
\item Find the integer point transform of the unimodular polygon in the Figure. 
\item  Find the Ehrhart polynomial $L_{\P}(t)$ of the integer polygon $\P$ from part (a). 
\end{enumerate}
 \end{prob}

 \medskip
 \begin{prob} $\clubsuit$  \label{two definitions for a solid angle}
Show that  \eqref{def. of solid angle} is equivalent to the following definition, using balls instead of spheres.  
Recall that the unit ball in $\R^d$ is define by  $B_d:= \{  x\in \R^d \mid \| x\| \leq 1 \}$, and similarly
the ball of radius $\varepsilon$, centered at $x\in \R^d$, is denoted by $B_d(x, \varepsilon)$.  Show that
for all sufficiently small $\varepsilon$, we have
\begin{equation*}
    \frac{\vol(S^{d-1}(x,\eps) \cap \P)}{\vol(S^{d-1}(x,\eps))} =
        \frac{\vol(B^{d}(x,\eps) \cap \P)}{\vol(B^{d}(x,\eps))}.
\end{equation*}

 \end{prob}

 \medskip
 \begin{prob} $\clubsuit$ \label{properties of floor, ceiling, fractional part}
Here we gain some practice with `floors', `ceilings', and `fractional parts'.     
First, we recall that by definition, the fractional part of any real number $x$ is  $\{x\} := x - \floor{x}$.
Next, we recall the indicator function of  $\Z$, defined  by:
 $
 1_{\Z}(x) :=
 \begin{cases}
      1 & \text{if }       x \in \Z \\
     0  & \text{if }        x  \notin \Z \\
     \end{cases}.
 $
 
 Show that:
\begin{enumerate}[(a)]
 \item   
 \label{problem:floor and celings, a}
  $\left\lceil x  \right\rceil  = - \floor{-x}$  
 \item          
 \label{problem:floor and celings, b}
 $1_{\Z}(x)=  \floor{x} - \left\lceil x  \right\rceil  +1$
\item           
\label{problem:floor and celings, c}
$ \{ x \} + \{-x\} = 1- 1_{\Z}(x)$
\item  
\label{problem:floor and celings, d}
$\floor{  x + y } \geq \floor{ x } + \floor{y}$, for all $x, y \in \R$.
\item  
\label{problem:floor and celings, e}
Let $m \in \Z_{>0}, n \in \Z$.  Then $\floor{ \frac{n-1}{m} } + 1 = \left\lceil   \frac{n}{m} \right\rceil$.
\end{enumerate}
\end{prob}

\medskip
\begin{prob}   \label{ (-1)^[x] in terms of periodic Beroulli polys}
Show that for $x \in \R \setminus \Z$, we have:
\begin{equation}
(-1)^{\floor{x}} = 2P_1(x) - 4 P_1\left( \frac{x}{2} \right),
\end{equation}
where we recall the definition of the first periodic Bernoulli polynomial $P_1(x):= x - \floor{x} - \tfrac{1}{2}$.
\end{prob} 
 \index{periodic Bernoulli polynomial}

\medskip
\begin{prob}   $\clubsuit$ \label{Ehrhart poly for closure of standard simplex}
\index{standard simplex}
Show that  the number of nonnegative integer solutions $x_1, \dots, x_d, z \in \Z_{\geq 0}$
to 
\[
x_1 + \cdots + x_d + z =  t, 
\]
with $ 0 \leq z \leq t$, equals  ${t+d \choose d}$.
\end{prob}

\medskip
\begin{prob}   $\clubsuit$  \label{Ehrhart poly for interior of standard simplex}
Show that for each positive integer $t$, the number of {\bf positive} integer solutions to 
$x_1 + \cdots + x_d <  t$ is equal to ${t-1 \choose d}$. 
\end{prob}

 \medskip
 \begin{prob}  
 We define the rational triangle whose vertices are $(0, 0), (1, \frac{N-1}{N}), (N, 0)$, where $N \geq 2$ is a fixed integer.
 Prove that the Ehrhart quasi-polynomial is in this case 
 \[
 L_\P(t) = \frac{p-1}{2} t^2 + \frac{p+1}{2} t + 1, 
 \]
 for all $t \in \Z_>0$.
 
 Notes.   So we see here a phenomenon known as `period collapse', where we expect a quasi-polynomial behavior, with some nontrivial period, but in fact we observe a strict polynomial.  
 \end{prob}

 \medskip
 \begin{prob}  
Here we show that the Ehrhart polymomial $L_\P(t)$ remains invariant under the full unimodular group $GL_d(\Z)$.  
In particular, recalling definition \ref{Definition of the unimodular group}, of a unimodular matrix,  show that:
 \begin{enumerate}[(a)]
\item   Every element of $GL_d(\Z)$ acts on the integer lattice $\Z^d$ bijectively. 
\item   \label{invariance of Ehrhart under the unimodular group}
Let $\P$ be an integral polytope, and let $Q := A(\P)$, where $A \in GL_d(\Z)$.   In other words,
$\P$ and $Q$
 are unimodular images of each other, by definition.   Prove that
 \[
L_{\P}(t) = L_Q(t),
 \]
for all $t \in \Z_{>0}$.
\item  Is the converse of part \ref{invariance of Ehrhart under the unimodular group} true?  In other words, given integer polytopes $\P, Q\subset \R^d$, suppose that $L_{\P}(t) = L_Q(t)$, for all positive integers $t$.  Does it necessarily follow that 
$Q~=~A(\P)$, for some unimodular matrix $A \in GL_d(\Z)$?
\end{enumerate} 
 \end{prob}

\medskip
\begin{prob} 
Suppose we are given a general lattice $\L := M(\Z^d)$ for some invertible matrix $M$. Here we 
extend the notion of the Ehrhart polynomial, so that we are counting elements of $\L$:
\[
L_{\P}(t, \L):= \left | t\P \cap \L \right|.
\]
So by definition $L_{\P}(t, \Z^d):= L_\P(t)$.
show that:
\begin{enumerate}
\item $\vol M(\P) = \vol(\P) \det(\L)$.
\item    $L_{\P}(t, \L) = L_{M^{-1} \P}(t)$. 
\end{enumerate}
\end{prob}



\chapter{
\blue{The Fourier transform of a polytope via its hyperplane description:  \\
Stokes' Theorem}
}
\label{Stokes' formula and transforms} 
\index{Stoke's formula}    \index{face poset}

\begin{quote}    
    ``Like a zen koan, Stokes' Theorem tells us that in the end, what happens on the outside
    is purely a function of the change within.''
    
--Keenan Crane  
 \end{quote}

\begin{figure}[htb]
 \begin{center}
\includegraphics[totalheight=2.2in]{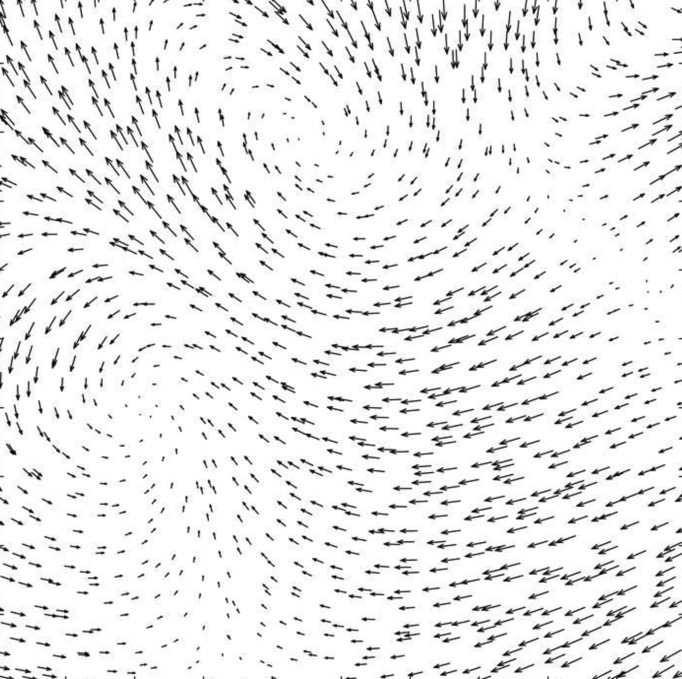}
\end{center}
\caption{A real vector field in $\R^2$}  \label{vector field}
\end{figure}

\red{(Under construction)}

\section{Intuition}

The divergence theorem, which is a special case of Stokes' more general theorem, 
is a multi-dimensional version of ``integration by parts'', a powerful tool from
the $1$-dimensional calculus.   We'll apply the divergence theorem to a polytope, to obtain
a combinatorial version of the divergence theorem.  This discrete version of the divergence theorem will allow us to transfer some of the complexity of computing the Fourier transform of a polytope to the complexity of computing corresponding Fourier transforms of its facets.   This kind of game can be iterated, yielding interesting geometric identities and results for polytopes, as well as for discrete volumes of polytopes.  

In the process, we also obtain another useful way to compute the Fourier transform of a polytope in its own right.

\section{The divergence theorem,   and a combinatorial \\ divergence theorem  for polytopes}
\index{combinatorial divergence theorem}

To warm up, we recall the divergence theorem, with some initial examples. 
A {\bf vector field} on Euclidean space is a 
function $F:\R^d \rightarrow \C^d$ that assigns to each point in $\R^d$ another vector in $\C^d$, which we will denote by
\[
F(x) := (F_1(x), F_2(x), \dots, F_d(x)) \in \C^d.
\]    
If $F$ is a continuous (respectively, smooth)  function, we say that $F$ is a 
{\bf continuous vector field}  (respectively, {\bf smooth vector field}).
If all of the coordinate functions $F_j$ are real-valued functions, we say that we have a {\bf real vector field}. 

We define the  {\bf divergence} of $F$ at  each $x := (x_1, \dots, x_d) \in \R^d$  by 
\[
{\rm div} F(x) :=   \frac{\partial F_1}{{\partial} x_1} + \cdots + \frac{\partial F_d}{\partial x_d},
\]
assuming that $F$ is a smooth (or at least once-differentiable) vector field.  This divergence of $F$ is a measure of the local change (sink versus source) of the vector field at each point $x\in \R^d$.  
Given a surface $S \subset \R^d$, and an outward pointing unit normal vector $\n$, defined at each point $x\in S$,  we also define the {\bf flux} of the vector field $F$ across the surface $S$ by 
\[
\int_S F\cdot \n  \ dS,
\]
where $dS$ denotes the Lebesgue measure of the surface $S$, and where the dot product 
$F\cdot \n$ is the usual inner product $\langle F, \n \rangle := \sum_{k=1}^d  F_k n_k$.  We will apply the divergence theorem (which is technically a special case of Stokes' Theorem) 
to a polytope $\P\subset \R^d$, and its $(d-1)$-dimensional bounding surface $\partial \P$.  Intuitively,
the divergence theorem tells us that the total divergence of a vector field $F$ inside a manifold is equal to the total 
flux of $F$ across its boundary.

\begin{thm}[The Divergence Theorem] 
 Let  $M \subset \R^d$ be a piecewise smooth manifold, and let $F$ be a smooth vector field. 
  Then
\begin{equation}\label{Divergence Theorem} \index{divergence Theorem}
\int_M {\rm div}F(x) dx =  \int_S F\cdot \n  \ dS.
\end{equation}
\hfill $\square$
\end{thm}

\begin{example} \label{Pyramid formula via the divergence theorem}
\rm{
Let $\P \subset \R^d$ be a $d$-dimensional polytope, containing the origin, with defining facets $G_1, \dots, G_N$. 
Define the real vector field 
\[
F(x):= x, 
\]
for all $x\in \R^d$.   First, we can easily compute 
here the divergence of $F$, which turns out to be constant:  
\begin{align*}
{\rm div } F(x) &=   \frac{\partial F_1}{{\partial} x_1} + \cdots + \frac{\partial F_d}{\partial x_d} = 
 \frac{\partial x_1}{{\partial} x_1} + \cdots + \frac{\partial x_d}{\partial x_d} = d.
\end{align*}
If we fix any facet $G$ of $\P$ then, due to the piecewise linear structure of the polytope,
every point $x \in G$ has the same constant outward pointing 
normal vector  to $F$, which we call  $\n_G$.
Computing first the left-hand-side of the divergence theorem, we see that 
\begin{equation}
\int_P {\rm div } F(x) dx = d  \int_P  dx = (\vol \P)d.
\end{equation}
Computing now the right-hand-side of the divergence theorem, we get
\begin{align*}
 \int_S F\cdot \n  \ dS = \int_{\partial \P} \langle x, \n  \rangle \ dS  = \sum_{k=1}^N \int_{G_k} \langle x, \n_G \rangle \ dS. 
\end{align*}
Now it's easy to see that the inner product $ \langle x, n_G \rangle$ is constant on each facet $G \subset \P$, namely
it is the distance from the origin to $G$ (Exercise \ref{distance to a facet}),  denoted by $\rm{dist}(G)$.
So we now have
\begin{align*}
 \int_{\partial \P} F \cdot n  \ dS    &= \sum_{k=1}^N     \int_{G_k}    \langle x, \n_{G_k} \rangle dS \\
&=  \sum_{k=1}^N  {\rm dist}(G_k)   \int_{G_k}  dS =  \sum_{k=1}^N  {\rm dist}(G_k) \vol G_k,
\end{align*}
so that altogether we the following conclusion from the divergence theorem:
\begin{equation}\label{pyramid formula}  \index{pyramid  formula}
\vol \P = \frac{1}{d}  \sum_{k=1}^N  {\rm dist}(G_k) \vol G_k.
\end{equation}
known as ``the pyramid formula'' for a polytope, a classical result in Geometry, which also has a very easy
geometrical proof (Exercise \ref{Pyramid formula}). 
}
\hfill $\square$
\end{example}

\bigskip
\begin{example}
\rm{
Let $\P \subset \R^d$ be a $d$-dimensional polytope with defining facets $G_1, \dots, G_N$, and outward pointing unit vectors $n_{G_1}, \dots, n_{G_N}$.   We fix any constant vector $\lambda \in \C^d$, and  we consider the {\bf constant vector field}  
\[
F(x):= \lambda, 
\]
defined for all $x\in \R^d$.   Here the divergence of 
$F$ is $\rm{div} F(x) = 0$, because $F$ is  constant, and so the left-hand-side of Theorem
\ref{Divergence Theorem}  gives us 
\begin{align*}
\int_P \rm{div } F(x) dx = 0.
\end{align*}
Altogether, the divergence theorem gives us:
\begin{align*}
 0 = \int_{\partial \P}   F \cdot \n  \ dS    &=       \sum_{k=1}^N     \int_{G_k}    \langle \lambda, \n_{G_k} \rangle dS \\
 &= \sum_{k=1}^N     \langle \lambda, \n_{G_k} \rangle        \int_{G_k}   dS \\
&=   \langle    \lambda,        \sum_{k=1}^N    \vol G_k   \n_{G_k}  \rangle,
\end{align*}
and because this holds for any constant vector $\lambda$, we can conclude that 
\begin{equation}\label{Minkowski relation}
 \sum_{k=1}^N    \vol G_k   \n_{G_k}= 0.
\end{equation}
Identity \eqref{Minkowski relation} is widely known as the {\bf Minkowski relation} for polytopes.  There is a marvelous converse to the latter relation, given by Minkowski as well, 
for any convex polytope
(See  Theorem \ref{Minkowski's problem for polytopes}).
}
\hfill $\square$
\end{example}

\bigskip
Now we fix $\xi \in \R^d$, and we want to see how to apply the divergence theorem to the vector-field  
\begin{equation}\label{First vector field}
F(x) :=  e^{- 2\pi i \langle x, \xi \rangle} \xi.
\end{equation}
Taking the divergence of the vector field $F(x)$, we have:
\begin{align*}
{\rm div } F(x)  &=   \frac{\partial  \left(  e^{- 2\pi i \langle x, \xi \rangle} \xi_1  \right)   }{{\partial} x_1}  + \cdots +
                     \frac{\partial  ( e^{- 2\pi i \langle x, \xi \rangle} \xi_d  )}{\partial x_d} \\
&=    (-2\pi i \xi_1^2)   e^{- 2\pi i \langle x, \xi \rangle}    + \cdots +   (-2\pi i \xi_d^2)  e^{- 2\pi i \langle x, \xi \rangle} \\     
&=     -2\pi i \| \xi \|^2    e^{- 2\pi i \langle x, \xi \rangle}.
\end{align*}             
So by the divergence theorem we have
\begin{align}  \label{initial divergence}
\int_{x\in P} - 2\pi i ||\xi||^2  e^{- 2\pi i \langle x, \xi \rangle} dx 
=    \int_{x\in P}  \text{div} F(x) dx =    \int_{\partial P} e^{- 2\pi i \langle x, \xi \rangle} 
\langle \xi,  \n \rangle \ dS,
\end{align}
where $\n$ is the outward-pointing unit normal vector at each point $x\in \partial \P$.   When $\P$ is a polytope, these arguments quickly give the following conclusion.

\begin{thm}   \label{FT of a polytope, first iteration of divergence}
Given any $d$-dimensional polytope $\P\subset \R^d$, with outward pointing normal vector $n_G$ to each facet $G$ of $\P$, its Fourier transform has the form
\begin{equation} \label{the first step of divergence}
\hat 1_\P(\xi) =\frac{1}{-2\pi i } \sum_{G\subset \partial P} \frac{ \langle \xi, \n_G \rangle}{  ||\xi ||^2} \hat 1_G(\xi),
\end{equation}
for all nonzero $\xi \in \C^d$.  Here  the integral that defines each $\hat 1_G$ is taken with respect to Lebesgue measure that matches the dimension of the facet $G \subset \partial P$.   
\end{thm}
\begin{proof}
 \begin{align*}
\hat 1_\P(\xi) &:= \int_{x\in P}  e^{- 2\pi i \langle x, \xi \rangle} dx  \\
&=  \frac{ 1 }{-2\pi i \|\xi\|^2} \int_{\partial P} \langle \xi,  \n \rangle e^{- 2\pi i \langle x, \xi \rangle}  dS 
         \quad (\text{using} \, \eqref{initial divergence}) \\
&=  \frac{ 1}{-2\pi i \|\xi\|^2} 
               \int_{G_1}  \langle \xi,   \n_{G_1} \rangle   e^{- 2\pi i \langle x, \xi \rangle}  dS 
                 + \cdots +  
      \frac{ 1}{-2\pi i \|\xi\|^2}    
               \int_{G_N}  \langle \xi,   \n_{G_N} \rangle   e^{- 2\pi i \langle x, \xi \rangle}  dS  \\
&=  \frac{  \langle \xi,   \n_{G_1} \rangle    }{-2\pi i \|\xi\|^2}    \hat 1_{G_1}(\xi) 
                 + \cdots +  
      \frac{  \langle \xi,  \n_{G_N} \rangle    }{-2\pi i \|\xi\|^2}     \hat 1_{G_N}(\xi),
\end{align*}
where in the third equality we used the fact that the boundary $\partial \P$ of a polytope is a finite union of $(d-1)$-dimensional polytopes (its facets), and hence $\int_{\partial P} = \int_{G_1} + \cdots + \int_{G_N}$, a sum of integrals over the $N$ facets of $\P$. 
\end{proof}
This result allows us to reduce the Fourier transform of $\P$ to a finite sum of Fourier transforms of the facets of $\P$.  This process can clearly be iterated, until we arrive at the vertices of $\P$.  But we will need a few book-keeping devices first.

To simplify the notation that will follow, we can also the {\bf Iverson bracket} notation, defined as follows.  
Suppose we have any 
boolean property $P(n)$, where $n\in \Z^d$;  that is,  $P(n)$ is either true or false. 
Then the Iverson bracket $[ P ]$ is defined by:
\begin{equation}\label{Iverson bracket}
[P] = 
\begin{cases}
1     & \mbox{if P is true }   \\ 
0     & \mbox{if P is false }
\end{cases}
\end{equation}
Now we 	may rewrite the identity of Theorem  \ref{FT of a polytope, first iteration of divergence}
as follows:
\begin{equation} 
\hat 1_\P(\xi) =\vol \P \ [\xi = 0] + 
\frac{1}{-2\pi i } \sum_{G\subset \partial P} \frac{ \langle \xi, \n_G \rangle}{  ||\xi ||^2} \hat 1_G(\xi)\ 
 [\xi \not= 0].
\end{equation}

Later, after Theorem \ref{Combinatorial divergence theorem} below, we will return to the Iverson bracket, and be able to use it efficiently.
To proceed further, we need to define the {\bf affine span} \index{affine span}  of a face $F$ of $\P$:
\begin{equation} \label{affine span of F}
{\rm aff}(F) := \left\{  \sum_{j=1}^k   \lambda_j v_j  \mid     k>0,  v_j \in F,   
\lambda_j \in \R, \text{  and  }  \sum_{j=1}^k   \lambda_j  = 1     \right\}.
\end{equation}

\begin{figure}[htb]
 \begin{center}
\includegraphics[totalheight=3.3in]{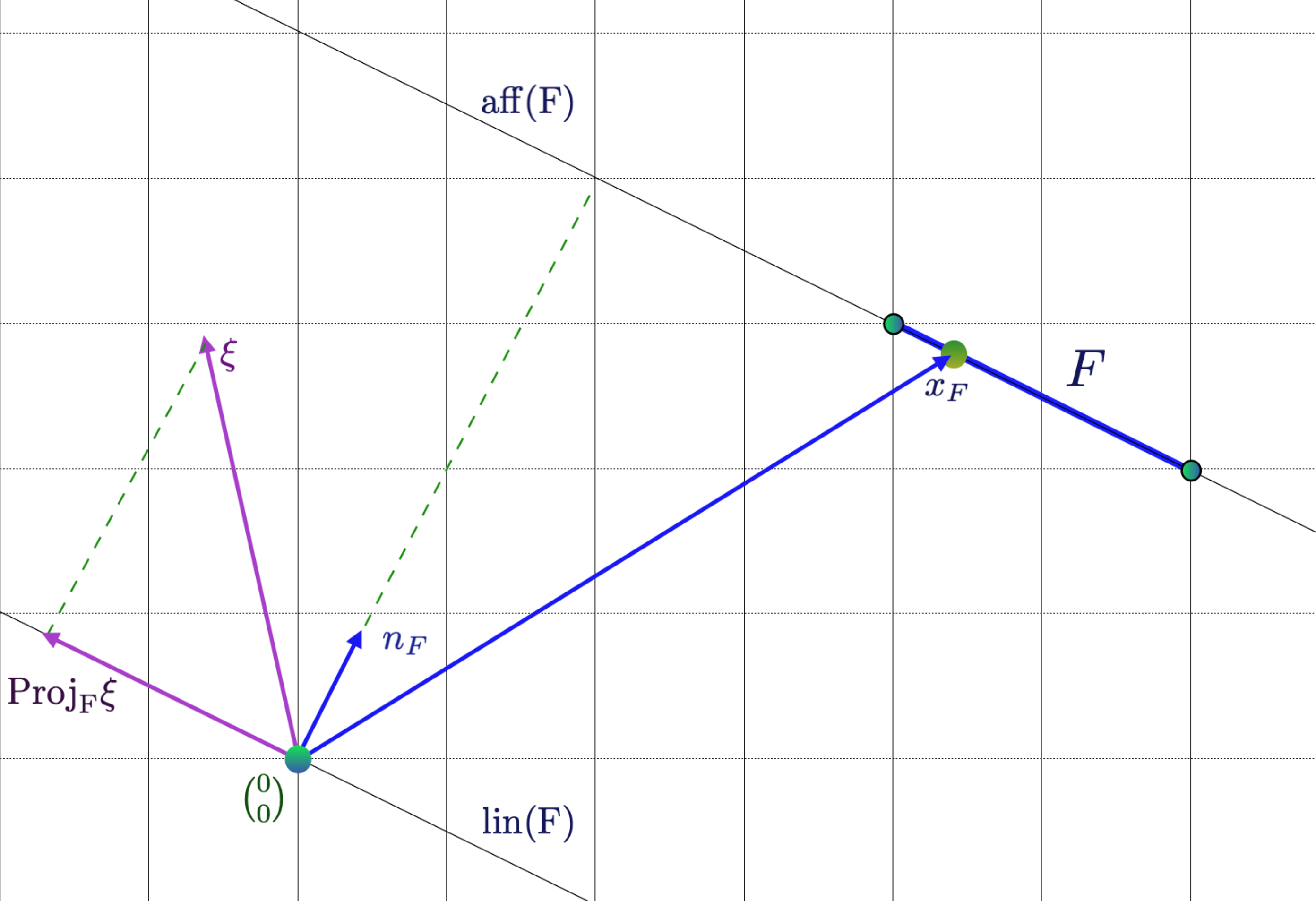}
\end{center}
\caption{The affine span of a face $F$, its linear span                , and the projection of $\xi$ onto $F$.  Here
we note that the distance from the origin to $F$ is $\sqrt{20}$.}  \label{affine span}
\end{figure}

In other words, we may think of the affine span of a face $F$ of $\P$ as follows.  We first translate $F$ by any element $x_0 \in F$.  So this translate, call if $F_0:= F - x_0$, contains the origin.  Then we take all real linear combinations of points of $F_0$, obtaining a vector subspace of $\R^d$, which we call the {\bf linear span } of $F$.  Another way to describe the linear span of a face $F$ of $\P$ is:
\[  
 {\rm lin}(F):= \left\{   x - y    \mid  x, y \in F  \right\}.
 \]
 \index{linear span}
Finally, we may translate this subspace ${\rm lin}(F)$ back using the same translation vector $x_0$, to obtain   ${\rm aff}(F):= {\rm lin}(F) + x_0$ (see Figure \ref{affine span}).

\begin{example}
\rm{
The affine span of two distinct points in $\R^d$ is the unique line in $\R^d$ passing through them.   The affine span of three points in $\R^d$ is the unique $2$-dimensional plane passing through them.  
The affine span of a $k$-dimensional polytope $F \subset \R^d$ is a translate of a $k$-dimensional vector subspace of $\R^d$. 
  Finally, the affine span of a whole $d$-dimensional polytope $\P \subset \R^d$ is all of $\R^d$. 
  }    
\hfill $\square$
\end{example}

In formalizing \eqref{the first step of divergence} further, we will require the notion of the projection of any point 
$\xi \in \R^d$ onto the linear span of any face $F\subseteq \P$, which we abbreviate by $\rm{Proj}_F \xi$:
\begin{equation}
\rm{Proj}_F \xi := \rm{Proj}_{\rm{lin}(F)}(\xi).
\end{equation}
(see Figure \ref{affine span})
We will also need the following elementary fact.  Let $F$ be any $k$-dimensional  polytope in
$\R^d$, and fix the outward-pointing unit normal to $F$, calling it $\n_F$.
 It is straightforward to show that if we take any point  $x_F  \in F$,  then $\langle x_F, \n_F \rangle$ is the distance from the origin to $F$.  Therefore, if $\rm{Proj}_F \xi = 0$, then a straightforward computation shows that
$\langle \xi, x_F \rangle  = \| \xi \|  \rm{dist}(F)$ (Exercise \ref{distance to a facet}).


\medskip
\section{A combinatorial divergence theorem}

We're now ready to extend   \eqref{the first step of divergence} to polytopes whose dimension is lower than the ambient dimension, as follows.
\begin{thm}[Combinatorial Divergence Theorem]  \label{Combinatorial divergence theorem}
Let $F$ be a polytope in $\R^d$, where $1 \leq \dim F \leq d$.  For each facet $G \subseteq F$, we let
$\n(G, F)$ be the unit normal vector to $G$, with respect to  $\rm{lin}(F)$.     
Then for each $\xi \in \R^d$, we have:
\begin{enumerate}[(a)]
\item  If $\rm{Proj}_F \xi = 0$, then 
\begin{equation}
\hat 1_F(\xi) = (\vol F) e^{-2\pi i     \| \xi \|  \rm{dist}(F)      }.
\end{equation}
\item  If $\rm{Proj}_F \xi \not= 0$, then
\begin{equation}
\hat 1_F(\xi) =  \frac{1}{-2\pi i } \sum_{G\subset \partial F} 
\frac{ \langle \rm{Proj}_F \xi, \n(G, F) \rangle}{  ||\rm{Proj}_F \xi ||^2}    \hat 1_G(\xi).
\end{equation}
\end{enumerate}
\end{thm}
\hfill $\square$

\bigskip
We notice that, as before, we are getting rational-exponential functions for the Fourier transform of a polytope.  But Theorem \ref{Combinatorial divergence theorem} gives us the extra freedom to begin with a lower-dimensional polytope $F$, and then find its Fourier transform in terms of its facets.

We are now set up to iterate this process, defined by Theorem \ref{Combinatorial divergence theorem},
reapplying it to each facet $G \subset \partial \P$.  Let's use the Iverson bracket, 
defined in \eqref{Iverson bracket},
 and apply the combinatorial divergence Theorem \ref{Combinatorial divergence theorem} to $\P$ twice:

\begin{align*} 
\hat 1_\P(\xi)   &=\vol \P \ [\xi = 0] + 
\frac{1}{-2\pi i } \sum_{F_1 \subset \partial P} \frac{ \langle \xi, \n_{F_1} \rangle}{  ||\xi ||^2}  \  [\xi \not= 0]  \  \hat 1_{F_1}(\xi) \\
                       &=\vol \P \ [\xi = 0] + 
\frac{1}{-2\pi i } \sum_{{F_1}\subset \partial P} \frac{ \langle \xi, \n_{F_1} \rangle}{  ||\xi ||^2}   [\xi \not= 0]  \\
&   \cdot  \Big(     (\vol {F_1}) e^{-2\pi i \langle \xi, x \rangle} \ [\rm{Proj}_{F_1} \xi = 0 ] + 
\frac{1}{-2\pi i } \sum_{F_2 \subset \partial {F_1}} 
\frac{ \langle \rm{Proj}_{F_2} \xi, \n(F_2, F_1) \rangle}{  ||\rm{Proj}_{F_2} \xi ||^2}    \hat 1_{F_2}(\xi) [\rm{Proj}_{F_1} \xi \not= 0]  \Big)
\\
&=   \vol \P \ [\xi = 0]  +
\frac{1}{-2\pi i } \sum_{F_1 \subset \partial P} 
\frac{ \langle \xi, \n_{F_1} \rangle  (\vol F_1) e^{-2\pi i \langle \xi, x \rangle} }{  ||\xi ||^2}  
 \ [\xi \not= 0][\rm{Proj}_{F_1} \xi = 0 ] \\
& +  \frac{1}{(-2\pi i )^2} \sum_{F_1 \subset \partial P}   \sum_{F_2 \subset \partial F_1} 
\frac{ \langle \xi, \n_{F_1} \rangle}{  ||\xi ||^2}  
\frac{ \langle \rm{Proj}_{F_2} \xi,  \n(F_2, F_1) \rangle}{  ||\rm{Proj}_{F_2} \xi ||^2}    
\hat 1_{F_2}(\xi)  \ [\xi \not= 0]  [\rm{Proj}_{F_1} \xi \not= 0]
\end{align*}

It is an easy fact that the product of two Iverson brackets is the Iverson bracket of their intersection: 
$[ P ] [ Q ] = [ P \text{ and } Q ]$ (Exercise  \ref{Exercise Iverson bracket}).    Hence, if we define
\[
F^\perp := \{  x \in \R^d \mid \langle x, y \rangle = 0 \text{ for all } y \in \rm{lin}F \}, 
\]
Then we see that $\P^\perp = \{ 0 \}$, and we can rewrite the latter identity as
\begin{align*} 
\hat 1_\P(\xi)   &=   \vol \P \ [\xi \in \P^\perp ]  +
\frac{1}{-2\pi i } \sum_{F_1 \subset \partial P} 
\frac{ \langle \xi, \n_{F_1} \rangle  (\vol F_1) e^{-2\pi i \langle \xi, x \rangle} }{  ||\xi ||^2}  
 \ [ \xi \in F_1^\perp - \P^\perp] \\
& +  \frac{1}{(-2\pi i )^2} \sum_{F_1 \subset \partial P}   \sum_{F_2 \subset \partial F_1} 
\frac{ \langle \xi, \n_{F_1} \rangle}{  ||\xi ||^2}  
\frac{ \langle \rm{Proj}_{F_2} \xi,  \n(F_2, F_1) \rangle}{  ||\rm{Proj}_{F_2} \xi ||^2}    
\hat 1_{F_2}(\xi)  \ [ \xi \not\in F_1^\perp].
\end{align*}

In order to keep track of the iteration process, we will introduce another book-keeping device. 
The {\bf face poset} of a polytope $\P$ \index{face poset}  is defined to be the partially ordered set (poset) of all faces of $\P$, ordered by inclusion, including $\P$ and the empty set.

\bigskip
\begin{example}
\rm{
Consider a $2$-dimensional polytope $\P$ that is a triangle.  We have the following picture for the face poset ${\frak F}_P$ of $\P$, as in Figure \ref{FacePoset1}.   It turns out that if we consider a $d$-simplex $\P$, then its  face poset ${\frak F}_P$ has the structure of a ``Boolean poset'' (which is isomorphic to the edge graph of a $(d+1)$-dimensional cube). 
\begin{figure} 
\centering
     \includegraphics[width=.5\textwidth]{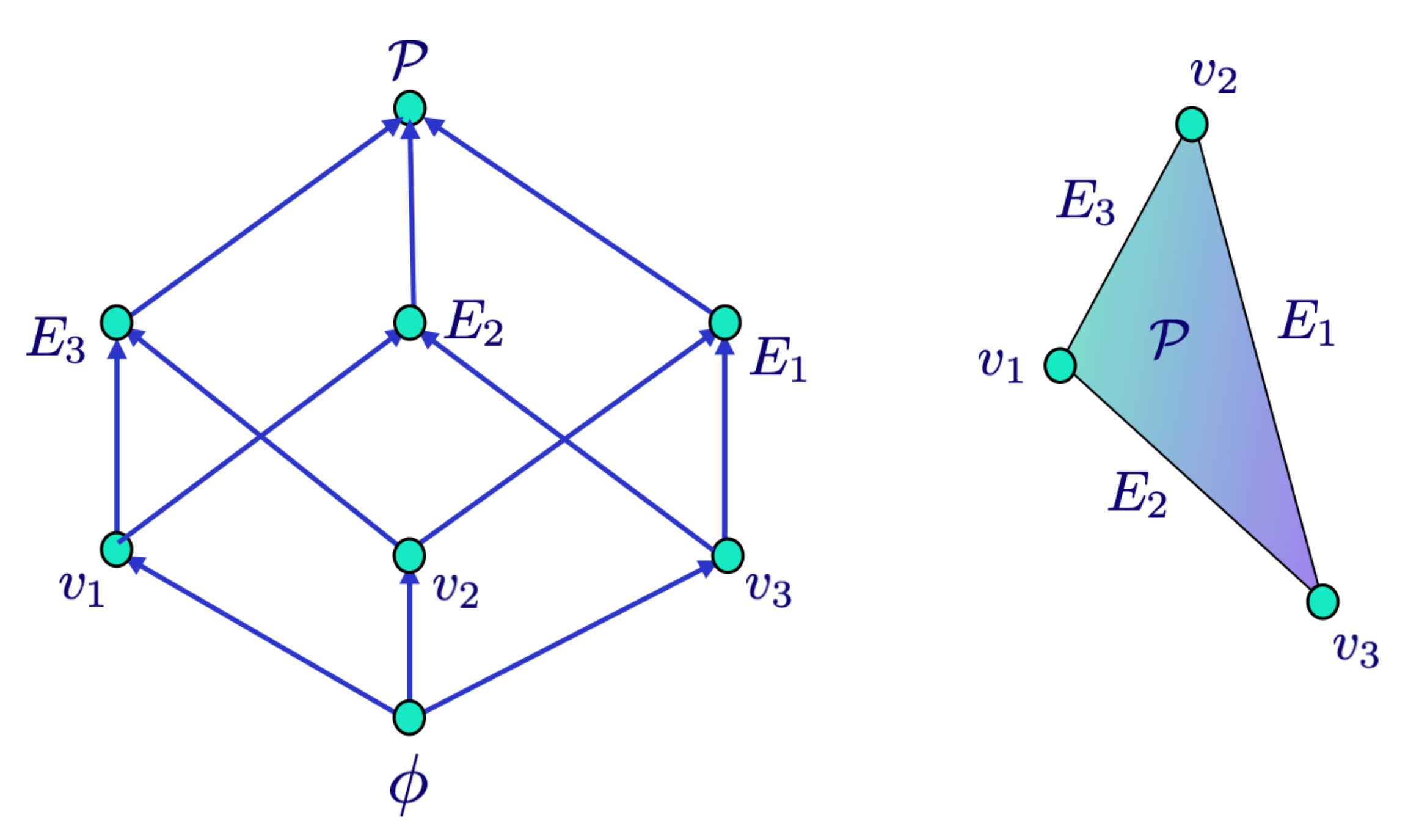}
     \caption{The face poset of a triangle}  \label{FacePoset1}
\end{figure} 
}
\end{example}


We only have to consider rooted chains in the  face poset ${\frak F}_P$, which means chains whose root is $P$.  The only appearance of non-rooted chains are in the following definition. 
 If $G$ is a facet of $F$, we attach the following weight to any  (local) chain $(F,G)$, of length $1$,  in the face poset of $P$:

\begin{equation}\label{weight}
 W_{(F,G)}(\xi):=\frac{-1}{2 \pi i} \frac{\langle \proj_{F} (\xi), \n(G, F) \rangle }{\| \proj_{F} (\xi) \|^2}. 
 \end{equation}

Note that these weights are functions of $\xi$ rather than constants. Moreover, they are all homogeneous of degree $-1$.     Let $\mathbf{T}$ be any rooted chain in ${\frak F}_P$, given by 
\[
T:= (P \to F_1 \to F_2, \dots, \to F_{k-1} \to F_k),
\]
so that by definition $\dim(F_j) = d-j$.
We define the {\bf admissible set} $S(\mathbf{T})$ of the rooted chain
 $\mathbf{T}$ to be the set of all vectors  $\xi\in \R^d$ that are orthogonal to the linear span of $F_k$ but not  orthogonal to the linear span of $F_{k-1}$.  In other words, 
 \begin{align*}
 S(\mathbf{T}) &:= \{  \xi \in \R^d \mid \xi  \perp \rm{lin}(F_k),    \text{  but   }
 \xi  \not\perp \rm{lin}(F_{k-1})    \}   \\
& = \{        \xi \in \R^d \mid \xi   \in F_k^\perp -     F_{k-1}^\perp          \}.
 \end{align*}

Finally, we define the following weights associated to any such rooted chain $\mathbf{T}$:

\begin{figure}
\centering
     \includegraphics[width=.5\textwidth]{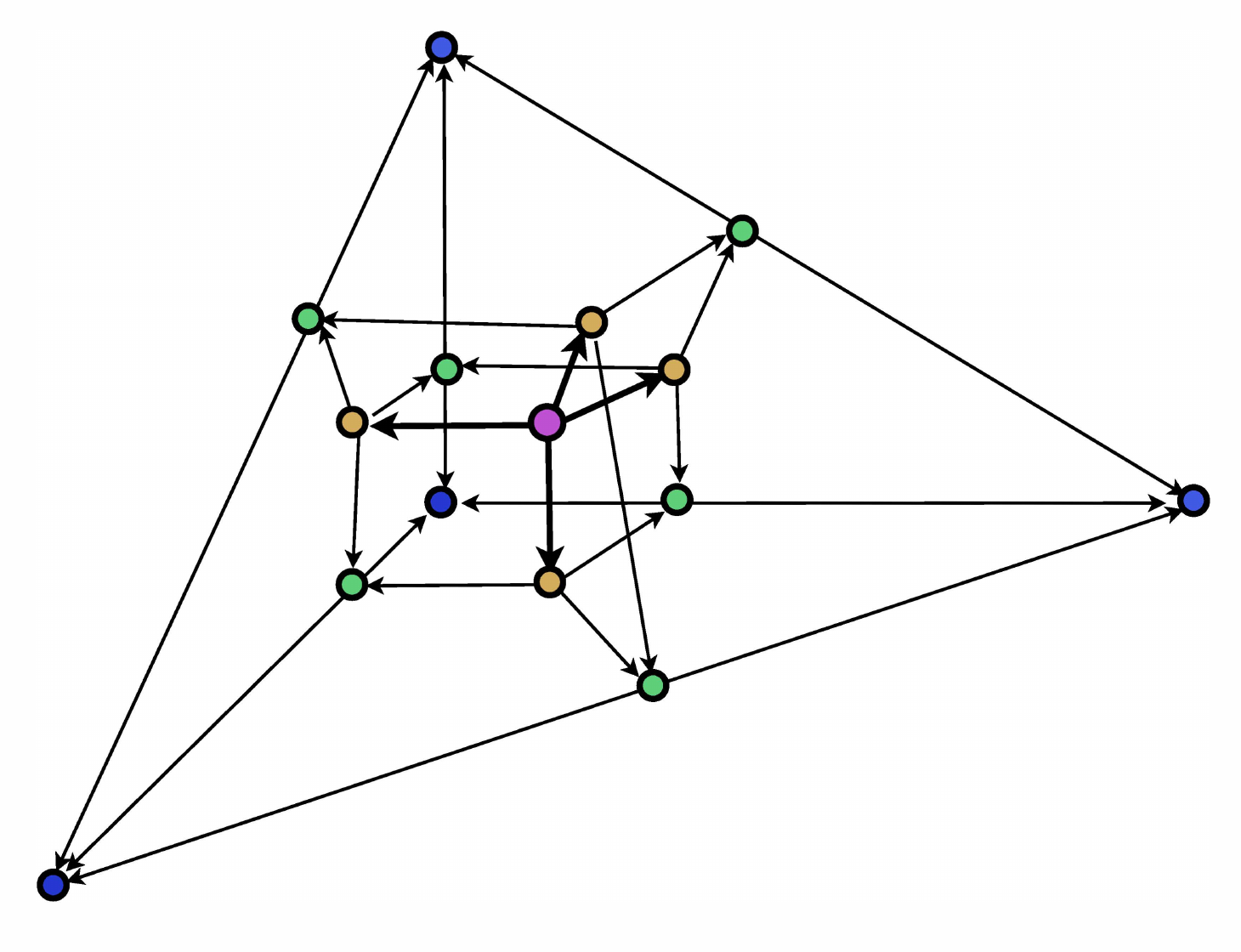}
     \caption{A symbolic depiction of the face poset ${\frak F}_P$, where $P$ is a $3$-dimensional tetrahedron.  Here the points and arrows are drawn suggestively, as a directed graph.  We can see all  the rooted chains, beginning from a symbolic vertex in the center, marked with the color purple.  The rooted chains that terminate with the yellow vertices have length $1$, those that terminate with the green vertices have length $2$, and those that terminate with the blue vertices have length $3$. }
\end{figure}

\begin{enumerate}[(a)]
		\item   The rational weight $\mathcal{R}_{\mathbf{T}}(\xi) = \mathcal{R}_{(P \to ... \to F_{k-1} \to F_k)}(\xi)$ is defined to be the product of weights associated to all the rooted chains  $\mathbf{T}$
 of length $1$, times the Hausdorff volume of $F_k$ (the last node of the chain $\mathbf{T}$).  It is clear from this definition that $\mathcal{R}_{\mathbf{T}}(\xi)$ is a homogenous rational function of $\xi$.
		
\bigskip
\item     The exponential weight 
$\mathcal{E}_{\mathbf{T}}(\xi) = \mathcal{E}_{(P \to ... \to F_{k-1} \to F_k)}(\xi)$ 
is defined to be the evaluation of $e^{-2\pi i\langle\xi,x\rangle}$ at any point $x$ on the  face $F_k$:
\begin{equation} \label{exponential.weight}
\mathcal{E}_{\mathbf{T}}(\xi) := e^{-2\pi i\langle\xi,x_0\rangle},
\end{equation}
for any $x_0 \in F_k$.  We note that the inner product $\langle\xi,x_0 \rangle$ does not depend on the position of $x_0 \in F_k$.

\bigskip		
\item      The total weight of a rooted chain $T$ is defined to be the rational-exponential function
\begin{equation}
W_{\mathbf{T}}(\xi) = W_{(P \to ... \to F_{k-1}  \to F_k)}(\xi):= \mathcal{R}_{\mathbf{T}}(\xi) \mathcal{E}_{\mathbf{T}}(\xi) \mathbf{1}_{S(\mathbf{T})}(\xi), 
\end{equation}

\noindent
where $\mathbf{1}_{S(\mathbf{T})}(\xi)$ is the indicator function of the admissible set $S(\mathbf{T})$ of $\mathbf{T}$.

\end{enumerate} 
	\bigskip
	
\noindent
By repeated applications of the combinatorial divergence 
Theorem \ref{Combinatorial divergence theorem},  
we arrive at  a description of the Fourier transform of $P$ as the sum of weights of all the rooted chains of 
the face poset ${\frak F}_P$, as follows.

\begin{thm}  \label{ingredient1}
\begin{align} \label{explicit Fourier transform of a polytope}
		\hat 1_{P}(\xi) = \sum_{\gT} W_{\gT}(\xi) = \sum_{\gT} \cR_{\gT}(\xi) \cE_{\gT}(\xi) \One_{S(\gT)}(\xi),
\end{align}
valid for any fixed $\xi \in \R^d$.	
\end{thm}
For a detailed proof of Theorem \ref{ingredient1}, see \cite{RicardoNhatSinai}.
Using this explicit description of the Fourier transform of a polytope, we will see an application of it in the following
section,  for the coefficients of Macdonald's angle quasi-polynomial.	In the process, equation 
\eqref{explicit Fourier transform of a polytope}, which gives an explicit description of the Fourier transform of a polytope, using the facets of $\P$ as well as lower-dimensional faces of $\P$,  will become even more explicit with some examples.

\section{Generic frequencies versus special frequencies} 
\label{sec: generic vs. special frequencies}

Given a polytope $\P \subset \R^d$, we call a vector $\xi \in \R^d$ a {\bf generic frequency} (relative to $\P$) \index{generic frequencies}
 if $\xi$ is not orthogonal to any face of $\P$.  All other $\xi \in \R^d$ are orthogonal to some face $F$ of $\P$, and are called {\bf special frequencies}.    

We see from Theorem \ref{ingredient1} that for a generic frequency $\xi$, we have 
\begin{align}
		\hat 1_{P}(\xi) = \sum_{\gT: P \to ... \to F_{1} \to F_0} \cR_{\gT}(\xi)   e^{-2\pi i \langle \xi, F_0 \rangle},
\end{align}
where the $F_0$ faces are the vertices of $\P$.   In other words, for generic frequencies, all of our  rooted chains in the face poset of $\P$ go all the way to the vertices.   The special frequencies, however, are more complex.   But we can collect the special frequencies in `packets', giving
us the following result.


\bigskip
\begin{thm} [Coefficients for Macdonald's angle quasi-polynomial] \cite{RicardoNhatSinai}
  \label{thm:main} 
  
Let $P$ be a $d$-dimensional rational polytope in $\R^d$, and let $t$ be a positive real number.
Then we have the quasi-polynomial
			\[ A_P(t) =\sum_{i = 0}^d a_i(t)t^i, \]
		where, for $0 \leq i \leq d$,
\begin{equation}\label{complicatedcoeff}
 a_i(t) := \lim_{\eps\to 0^+} \sum_{\xi\in\Z^d \cap S(\gT)} 
 \sum_{l(\gT) = d-i} \cR_{\gT}(\xi) \cE_{\gT}(t\xi) \  e^{-\pi\eps\|\xi\|^2},
 \end{equation}
\end{thm}

\noindent
where $l(\gT)$ is the length of the rooted chain $\gT$ in the face poset of $P$, 
$\cR_{\gT}(\xi)$ is the rational function of $\xi$ defined above,  $\cE_{\gT}(t\xi) $ is the complex exponential defined in \eqref{exponential.weight} above, and 
$\Z^d \cap S(\gT)$ is the set of all integer points that are orthogonal to the last node in the chain $T$, 
but not to any of its previous nodes. 
\hfill $\square$

See \cite{RicardoNhatSinai} for the detailed proof of Theorem \ref{thm:main}.

\medskip
\section{The codimension-$1$ Ehrhart coefficient under continuous dilations}
We call the coefficients $a_i(t)$ the {\bf quasi-coefficients} of the solid angle sum $A_P(t)$.   \index{quasi-coefficients}
As a consequence of Theorem \ref{thm:main}, it turns out that there is a closed form for the codimension-$1$ quasi-coefficient, which extends previous special cases of this coefficient.  

We recall our first periodic Bernoulli polynomial, from \eqref{definition of periodic Bernoulli polys}:
\begin{equation}
P_1 (x):=  
 \begin{cases} 
 x - \lfloor x  \rfloor- \frac{1}{2}            & \mbox{if } x \notin \Z  \\ 
0 & \mbox{if } x \in \Z, 
\end{cases}
\end{equation}
where $\lfloor x  \rfloor$ is the integer part of $x$.

\begin{thm}\cite{RicardoNhatSinai}  
 \label{codim1coeff}
Let $P$ be any real polytope.  Then the codimension-1 quasi-coefficient of the solid angle sum $A_P(t)$ has the following closed form:
\begin{equation}\label{codimension-1 formula for angle poly}
a_{d-1} (t) = 
-\sum_{\substack{F \textup{ a facet of } P \\ with \ n_F \neq 0}}  \frac{\vol F}{\|n_F\|} 
P_1 (\langle n_F, x_F \rangle t),
\end{equation}
where $n_F$ is the unique primitive integer vector which is an outward-pointing normal vector  to $F$,  $x_F$ is
 any point lying in the affine span of $F$, and $t$ is any positive real number. 
\hfill $\square$
\end{thm}
We note that  
the latter formula shows in particular that for any rational polytope $\P$,  the quasi-coefficient  $a_{d-1}(t)$ is always  a periodic function of $t > 0$, with a period of $1$. 
For rational polytopes and all of their real dilates, the quasi-coefficients of their quasi-polynomials
are periodic functions of  real dilations $t>0$, as we show below.
\begin{example}
{\rm
To see what  formula \eqref{codimension-1 formula for angle poly} of 
Theorem \ref{codimension-1 formula for angle poly}
tells us for  $2$-dimensional integer polygons, 
let's fix an integer polygon $\P$, whose vertices are $v_1, \dots v_N \in \Z^2$.  (Finish this example)
}
\hfill $\square$
\end{example}


\bigskip
\section{An extension of Pick's theorem to $\R^d$, using solid angles}


\begin{thm}  \label{cs.facets}  \index{Barvinok, Alexander}
Suppose $P$ is a $d$-dimensional integer polytope in $\R^d$ 
all of whose facets are symmetric.  Then
\[
A_\P(t)  = (\vol \P) t^d,
\]
for all positive integers $t$.
\end{thm}
\begin{proof}
We recall the formula for the solid angle polynomial: 
\begin{equation}\label{eq:APtsum}
A_\P(t) = \lim_{\varepsilon \to 0}   \sum_{\xi \in \Z^d} \hat{1}_{t\P}(\xi)e^{-\pi\varepsilon\|\xi\|^2}.
\end{equation}

The Fourier transform of the indicator function of a polytope may be written as follows, after one application of the combinatorial divergence formula:
\begin{align}  
\hat 1_{t\P}(\xi) = t^d \vol \P \, [\xi =0]
+ \left( \frac{-1}{2 \pi i}  \right)      t^{d-1} 
\sum_{{\substack{F \subseteq \P  \\ \dim F = d-1}}} 
 \frac{\langle \xi, \n_F \rangle }{\| \xi \|^2}      \hat 1_F       (t \xi) [\xi \not= 0],
\end{align}
where we sum over all facets $F$ of $\P$.   Plugging this into~\eqref{eq:APtsum} we get
\begin{align}  \label{imaginary}
A_\P(t)  - t^d \vol \P 
= \left( \frac{-1}{2 \pi i} \right) t^{d-1}
\lim_{\varepsilon \to 0} \sum_{\xi \in \Z^d \setminus\{0\}}
     \frac{e^{-\pi\varepsilon\|\xi\|^2}}{\| \xi \|^2} \sum_{{\substack{F \subseteq \P  \\ \dim F = d-1}}}
\langle \xi, \n_F  \rangle   \hat 1_F(t \xi),
\end{align}
so that it is sufficient to show that the latter sum over the facets vanishes.
The assumption that all facets of $\P$ are symmetric implies that $\P$ itself is also centrally symmetric, by Theorem~\ref{cs1}. 
We may therefore combine the facets of $\P$ in pairs of opposite facets $F$ and $F'$. We know that $F' = F + c$, where $c$ is an integer vector, using the fact that the facets are centrally symmetric (see Exercise \ref{half-integer translation for symmetric facets} for this little fact about $c$).


Therefore, since $\n_F' = - \n_F $, we have
\begin{align*}
\langle \xi, \n_F  \rangle  &\hat 1_{F}(t \xi) 
+ \langle \xi, -\n_F  \rangle   \hat 1_{F+ c}(t \xi)\\
&= \langle \xi, \n_F  \rangle   \hat 1_{F}(t \xi) 
- \langle \xi, \n_F  \rangle   \hat 1_{F}(t \xi) 
e^{-2\pi i\langle   t\xi, c  \rangle} \\
&=  \langle \xi, \n_F  \rangle  \hat 1_{F}(t \xi) 
\left( 1  -  e^{-2\pi i\langle   t\xi, c  \rangle} \right) = 0,
\end{align*}
because $\langle t\xi, c \rangle \in \Z$ when both  $\xi \in \Z^d$ and $t \in \Z$.  We conclude that the entire right-hand side of \eqref{imaginary} vanishes, and we are done.
\end{proof}
Theorem \ref{cs.facets}  appeared in  \cite{BarvinokPommersheim}, and here we gave a different proof, using the methods of this chapter.   The result of Alexandrov and Shephard (Theorem \ref{Alexandrov-Shepard thm}) from chapter \ref{Geometry of numbers} came in handy in our proof: if all the facets of $\P$ are symmetric, then $\P$ must be symmetric as well.   

One might wonder if the assumption of Theorem \ref{cs.facets} necessarily implies that 
$\P$ is a zonotope.  That this is not true is shown by the $4$-dimensional polytope called the $24$-cell, depicted in Figure \ref{24-cell}.
Fourier analysis can also be used to give yet more general classes of polytopes that satisfy the formula 
$A_P(t) = (\vol \P) t^d$, for positive integer values of $t$ (See also \cite{FabricioSinai2}, \cite{DeligneTabachnikovRobins}).

There is a wonderful result of Minkowski that gives a converse to the relation 
\eqref{Minkowski relation}, as follows.

\begin{thm}[The Minkowski problem for polytopes]   \label{Minkowski's problem for polytopes}   
\index{Minkowski problem for polytopes}
Suppose that $u_1, \dots, u_k\in \R^d$ are unit vectors that do not lie in a hyperplane.
Suppose further that we are given positive numbers $\alpha_1, \alpha_2, \dots, \alpha_k >0$ that 
satisfy the relation
\[
\alpha_1 u_1 + \cdots + \alpha_k u_k =0.
\]
Then there exists a polytope $\P\subset R^d$, with facet normals $u_1, \dots, u_k\in \R^d$, and 
facet areas $\alpha_1, \alpha_2, \dots, \alpha_k$.  
Moreover, this polytope $\P$ is unique, up to translations.
\hfill $\square$
\end{thm}

There is a large body of work, since the time of Minkowski, that is devoted to extensions of
 Minkowski's Theorem \ref{Minkowski's problem for polytopes}, to other convex bodies, as well as to other manifolds.

\bigskip \bigskip

\section*{Notes} \label{Notes.chapter.Divergence}
\begin{enumerate}[(a)]

\item  We could also define another useful vector field, for our combinatorial divergence theorem, besides 
our vector field in equation \eqref{First vector field}.   Namely, if we define 
$F(x):= e^{2\pi i \langle x, \xi \rangle} \lambda$,
for a fixed $\lambda\in\C^d$, then we would get the analogous combinatorial divergence formula as shown below in 
(Exercise   \ref{alternate combinatorial divergence Theorem}), 
and such vector fields
have been used, for example, by Alexander Barvinok \cite{Barvinok1} in an effective way.  
\index{Barvinok, Alexander}
To the  best of our knowledge,
the first researcher to use iterations of Stokes' formula to obtain lattice point asymptotics was Burton Randol \cite{Randol3}, \cite{Randol4}. 

\item  The Minkowski problem for polytopes can also be related directly to generalized isoperimetric inequalities for mixed volumes, as well as the Brunn-Minkowski inequality for polytopes, 
as done by Daniel Klain in \cite{Klain1}. 

\item We haven't delved into the differential forms perspective of Stokes' theorem, which may be even more appropriate for this line of research, in order to keep the background necessary to a minimum.   However, the differential forms approach is coordinate-free, and therefore has its advantages as well.  
\end{enumerate}

\bigskip 
Here's a true anecdote, which transpired in the $1979$ international symposium on differential geometry, in Berkeley, honoring S. S. Chern:

\begin{quote}   
(Person from the audience) ``What is the most important theorem in Differential geometry, in your opinion, professor Chern?"
 
(Chern's answer)  
\blue{ ``There is only one theorem in Differential Geometry, and that is Stokes' theorem.''}

(Person from the audience)     ``What is the most important theorem in Analysis, professor Chern?"

(Chern's answer) 
\blue{  ``There is only one theorem in Analysis, and that is Stokes' theorem.''}
    
(Another person from the audience)     ``And what is the most important theorem in Complex Analysis, professor Chern?"

 (Chern's answer) 
\blue{“There is only one theorem in complex variables, and that that is Cauchy's theorem. 
But if one assumes the derivative of the function is continuous, then this is just Stokes’ theorem”.}
\label{Chern's answer}

--Shiing-Shen Chern  
 \end{quote}


\bigskip
\section*{Exercises}
\addcontentsline{toc}{section}{Exercises}
\markright{Exercises}

\medskip
\begin{prob}   \label{Stokes implies Cauchy}
{\rm
If you know a bit of complex analysis,  then prove that Chern is correct in the anecdote above. 
 In other words, let $f(z) := u(x, y) + i v(x, y)$, where $z:=x+iy$, suppose that
$\partial f/\partial x$, $\partial f/\partial x$ are continuous
 on the unit ball $B:= \{ z\in \C \mid \|z\| \leq 1\}$.  Prove that Stokes' theorem implies  Cauchy's theorem:
\[
\int_{S^1} f(z) dz = 0.
\]
 }
\end{prob}

\medskip
\begin{prob} $\clubsuit$   \label{Pyramid formula}
{\rm 
We define the distance from the origin to $F$, denoted by $\rm{dist}(F)$,  as the length of the shortest vector of translation between 
${\rm aff}(F)$ and $\rm{lin}(F)$ (the affine span of $F$ and the linear span of $F$ were
defined in 
\eqref{affine span of F}).   
Figure  \ref{affine span} shows what can happen in such a scenario.
\begin{enumerate}[(a)]
\item        Suppose that we consider a facet $F$ of a given polytope $\P \subset \R^d$, and we
 let $\n_F$ be the unit normal vector to $F$.     Show that  the function 
\[
x_F \rightarrow  \langle x_F,   \n_F \rangle
\]
 is constant for $x_F \in F$,  and is in fact equal to the distance from the origin to $F$.  In other words, show that
 \[
\langle x, \n_F \rangle =  \rm{dist}(F).
\]
\item   Show that if $\rm{Proj}_F \xi = 0$, then $\langle \xi, x_F \rangle = \| \xi \| \rm{dist} F$.
\end{enumerate}
}
\end{prob}

\medskip
\begin{prob}   
{\rm
Here we prove the elementary geometric formula for a pyramid over a polytope.  Namely, suppose we are 
given a 
$(d-1)$-dimensional polytope $\P$, lying in the vector space defined by the first $d-1$ coordinates. 
We define a pyramid over $\P$, of height $h > 0$, as the $d$-dimensional polytope defined by 
\[
\rm{Pyr}(\P) := \conv\{ \P,  \ h \cdot e_{d}  \},
\]
where $e_d := (0, 0, \dots, 0, 1) \in \R^d$.  Show that 
\[
\vol \rm{Pyr}(\P) = \frac{h}{d} \vol \P.
\]
}
\end{prob}

\medskip
\begin{prob}    $\clubsuit$     \label{distance to a facet}
{\rm
Prove the Pyramid formula, \eqref{pyramid formula} in 
Example \ref{Pyramid formula via the divergence theorem},
  for a $d$-dimensional polytope $\P$ which contains the origin,  but now using just elementary geometry:
\begin{equation}
\vol \P = \frac{1}{d}  \sum_{k=1}^N  \rm{dist}(G_k) \vol G_k,
\end{equation}
where the $G_k$'s are the facets of $\P$, and $\rm{dist}(G_k)$ is the distance from the origin to $G_k$.
}
\end{prob}

\blue{We note that the next $3$ exercises are meant to be done together}

\medskip
\begin{prob}  $\clubsuit$   \label{alternate combinatorial divergence Theorem}
{\rm
Show that if we replace the vector field in equation \eqref{First vector field} by the 
alternative vector field  $F(x):= e^{-2\pi i \langle x, \xi \rangle} \lambda$, with a constant nonzero vector $\lambda \in \C^d$, then we get:
\begin{equation} \label{our alternate formula for the transform}
\hat 1_\P(\xi) =\frac{1}{-2\pi i } \sum_{G\subset \partial P} 
\frac{   \langle \lambda, \n_G \rangle}{ \langle  \lambda,  \xi \rangle} 
\hat 1_G(\xi),
\end{equation}
valid for all nonzero $\xi \in \R^d$.  Note that one advantage of this formulation
 of the Fourier transform
of $\P$ is that each summand in the right-hand-side of 
\eqref{our alternate formula for the transform}
is free of singularities, assuming the vector
$\lambda$ has a nonzero imaginary part. 
}
\end{prob}

\medskip
\begin{prob}  \label{equivalent identity to the alternate vector field}
Show that the identity \eqref{our alternate formula for the transform}  of Exercise \ref{alternate combinatorial divergence Theorem} is equivalent to the vector identity:
\[
\xi   \hat 1_\P(\xi)    =   \frac{1}{-2\pi i }        \sum_{G\subset \partial P} \n_G    \hat 1_G(\xi),
\]
valid for all  $\xi \in \R^d$. 
\end{prob}

\medskip
\begin{prob} \label{strange proof of Minkowski relation}
Show that the result of Exercise \ref{equivalent identity to the alternate vector field} 
quickly gives us
the Minkowski relation  \eqref{Minkowski relation}:  
\[
\sum_{ \text{facets } G \text{ of } P}    (\vol G)  \n_{G}= 0.
\]
\end{prob}

\medskip
\begin{prob}  
Continuing Exercise \ref{alternate combinatorial divergence Theorem}, show that by iterating this particular version of the Fourier transform of a polytope $\P$,  $k$ times, we get:
\begin{equation}
\hat 1_\P(\xi) =\frac{1}{(-2\pi i )^k} \sum_{G_k \subset G_{k-1} \subset \cdots G_1 \subset \partial P} 
 \prod_{j=1}^k 
 \frac{  
 \langle \lambda, \n_{G_{j}, G_{j-1}}  \rangle}{  \langle  \lambda,  \rm{Proj}_{G_{j-1}}  \xi \rangle } 
\hat 1_{G_k}(\xi),
\end{equation}
valid for all nonzero $\xi \in \R^d$, and where we sum over all chains $G_k \subset G_{k-1} \subset \cdots G_1$ of length $k$ in the face poset of $\P$, with \rm{codim}$(G_j) = j$. 
\end{prob}

\medskip
\begin{prob} \label{Geometric interpretation of Minkowski's relation for d=2}
Show that in the case of polygons in $\R^2$,  the Minkowski relation \eqref{Minkowski relation} 
has the meaning that 
the sum of the pink vectors in Figure \ref{Divergence Exercise} sum to zero.  In other words, the geometric interpretation
of the Minkowski relation in dimension $2$ is that the sum of the boundary (pink) vectors 
wind around the boundary and close up perfectly.
 \end{prob}

\begin{figure}[htb]
 \begin{center}
\includegraphics[totalheight=1.2in]{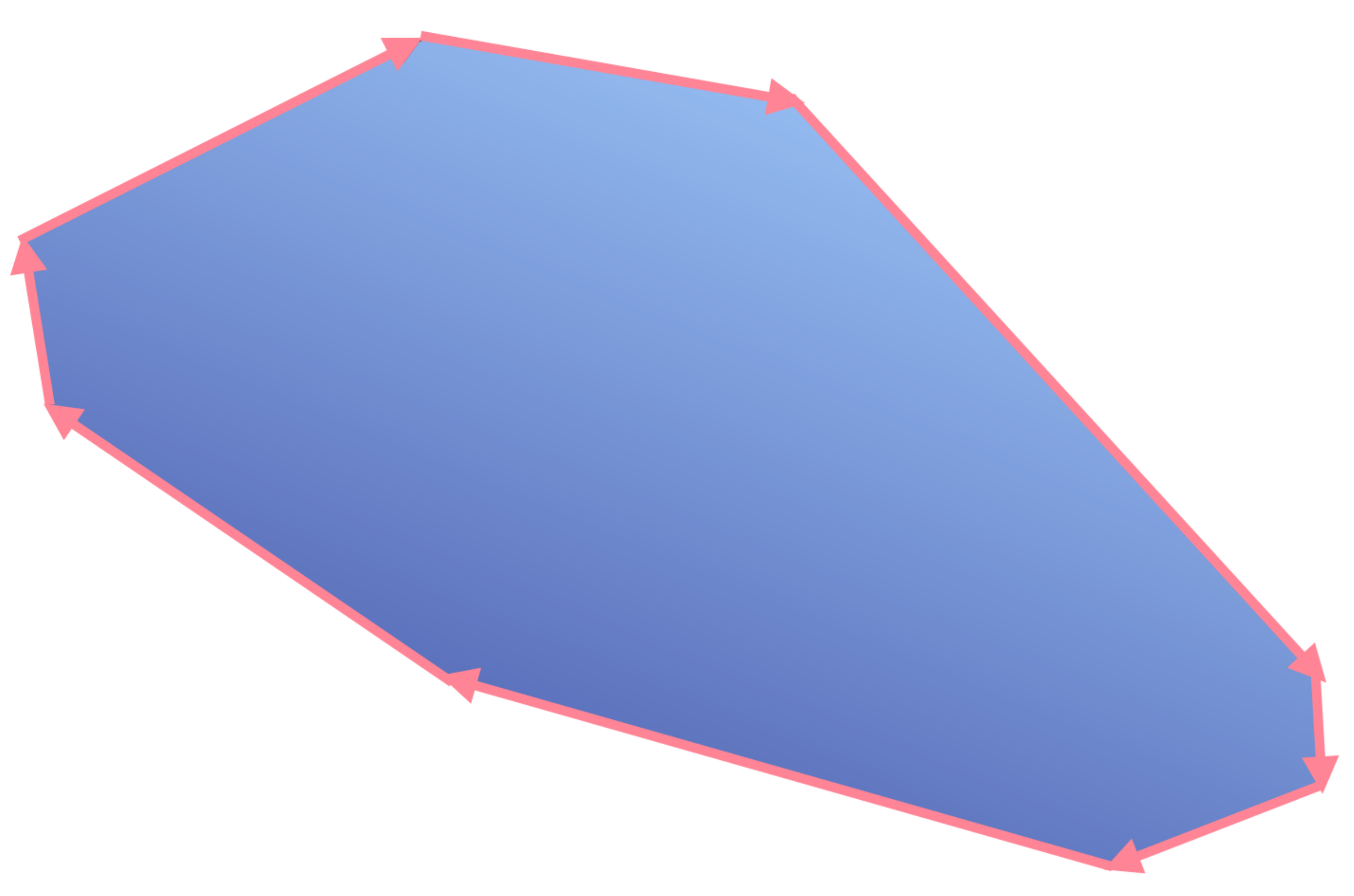}
\end{center}
\caption{The meaning of Minkowski's relation in dimension $2$ - see Exercise
\ref{Geometric interpretation of Minkowski's relation for d=2} }  
\label{Divergence Exercise}
\end{figure}

\medskip
\begin{prob} $\clubsuit$ \label{no simplex is symmetric}
Let's consider a simplex $\Delta \subset \R^d$ whose dimension satisfies $2 \leq \dim \Delta \leq d$. 
Show that $\Delta$ is not a symmetric body.
\end{prob}

\medskip
\begin{prob}  $\clubsuit$ \label{Exercise Iverson bracket}
To get more practice with the Iverson bracket, defined in equation \eqref{Iverson bracket}, 
show that for all logical statements $P$, we have:
\begin{enumerate}[(a)]   
\item $ [ P \rm{ \ and \ } Q ] = [ P ] [ Q ]$.
\item  $[ P \rm{ \ or \ } Q ] =  [P] + [Q] - [P][Q]$.
\item $[\neg P] = 1 - [P]$, where $\neg P$ means the logical negation of $P$. 
\end{enumerate}
\end{prob}

\medskip
\begin{prob} $\clubsuit$ \label{half-integer translation for symmetric facets}
Let $F \subset \R^d$ be a facet of a centrally symmetric, integer $d$-dimensional 
polytope $\P \subset \R^d$. 
Show that the distance from the origin to $F$ is always a half-integer or an integer.   In other words, show that 
\[
\rm{dist}(F) \in  \frac{1}{2} \Z.
\]
(See Exercise \ref{Pyramid formula} above for the definition of distance of $F$ to the origin)
\end{prob}


\chapter{
\blue{Classical geometry of numbers \\
 Part III: \, The covering radius, \\
 the packing radius, and successive minima}     
 } 
\label{Chapter.geometry of numbers III}

\begin{quote}    
My dear Watson, once you eliminate the impossible, then whatever remains -  no matter how improbable - must be the truth.

-- Arthur Conan Doyle (in his book Sherlock Holmes)
 \end{quote}


\bigskip
\section{The successive minima of a lattice}

A very important characteristic of a lattice $\L$ is the {\bf length of its  shortest nonzero vector}:
\index{shortest nonzero vector in a lattice}
\[
\lambda_1(\L):=\min \left\{    \| v \|  \biggm |    v \in \L-\{0\}    \right\}.
\]
Every lattice has at least two shortest nonzero vectors, because if $v \in \L$, then $-v \in \L$.  Therefore, when we use the words `its shortest vector', we always mean that we are free to make a choice between any of its vectors that have the same shortest, nonzero length.  

\begin{example} \label{first example of lattice basis reduction}
\rm{
Consider the following lattice in $\R^2$:
\[
\L := \left\{ m  \icol{102  \\ 11 }    + n  \icol{200\\16}  \bigm |     m,n\in\Z  \right\}.
\]
What is the shortest nonzero vector in this lattice $\L$?   Without using any fancy 
Theorems, we might still try simple subtraction, sort of
mimicking the Euclidean algorithm.  So for example, we might try   $\icol{200\\16} -   2 \icol{102  \\ 11 } =   \icol{-4 \\ -6 }$, which is pretty short.   So we seem to have gotten lucky - we found a relatively short vector.  But here
comes the impending question: 
how do we know whether or not this is really the shortest nonzero vector in our 
lattice $\L$? 
Can we find an even shorter vector in  $\L$?  
}
\hfill $\square$
\end{example}

\begin{figure}[htb]
 \begin{center}
\includegraphics[totalheight=3in]{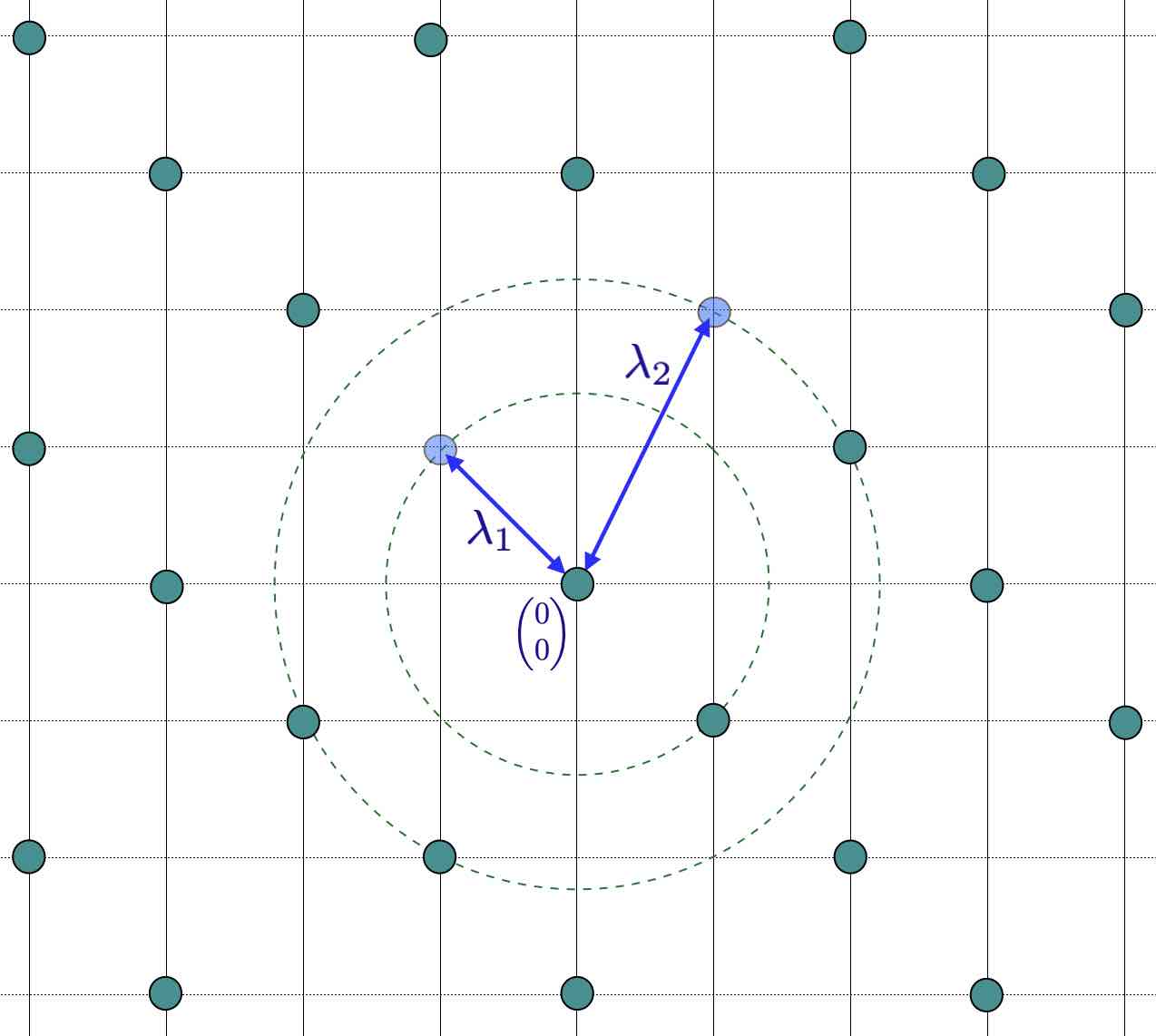}
\end{center}
\caption{The two successive minima for this lattice $\L$ are $\lambda_1(\L) = \sqrt 2$, and $\lambda_2(\L) = \sqrt 5$.}  \label{Successive Minima}
\end{figure}

The question raised in Example \ref{first example of lattice basis reduction} is not easy to answer in general, and we need to learn a bit more theory even to approach it in $\R^2$.   
In dimensions $d\geq 3$, the corresponding problem of finding a shortest nonzero vector in any given lattice is terribly difficult.  It is considered to be one of the most difficult - and one of the most important -  problems in computational number theory.
\begin{question}\rm{[The Shortest Vector Problem (SVP)]}
\label{shortest vector:question}
Given a basis for a lattice $\L \subset \R^d$, find a shortest nonzero vector in $\L$. 
\end{question}
Given that Question \ref{shortest vector:question} is notoriously hard in general (depending on the data we are given for the lattice $\L$), we can ask an easier question.
\begin{question}\label{length of shortest vector:question}
Can we find the {\bf length} of a shortest nonzero vector in $\L$?
\end{question}

Even this seemingly easier question turns out to be quite difficult, and important in many applications.  Minkowski gave an interesting approximation to 
Question \ref{length of shortest vector:question}, as we'll see shortly in Theorem \ref{First successive minima bound}.

To capture the notion of the second-smallest vector in a lattice, and third-smallest vector, etc, we begin by 
imagining balls of increasing radii, centered at the origin, and we can (at least theoretically) keep track of how 
they intersect $\L$.

Given a convex body $K\subset \R^d$,  
we let $r$ be the smallest positive real number such that  the dilated body $rK$ contains at least $j$ linearly independent lattice points of $\L$, for each $1\leq j \leq d$. 
This value of $r$ is called $\lambda_j(K, \L)$, 
 the {\bf  $j$'th successive minima} of the lattice, relative to $K$.

Here's another way of saying the same thing:
\begin{equation}\label{def:general successive minima}
\lambda_j(K, \L) := \min   \left\{ r >0  \bigm |   
  \dim \big(     \text{span }(\L \cap rK)    \big)   \geq j     \right\}.
\end{equation}
In the special case that $K= B$, the unit ball, we'll simply write 
\[
\lambda_j(B, \L):= \lambda_j(\L), 
\]
following the standard conventions in the geometry of numbers.  It follows from the definitions above that 
$0 < \lambda_1(K, \L) \leq  \lambda_2(K, \L) \leq \cdots \leq  \lambda_d(K, \L) < \infty$.  
Figure \ref{Successive Minima} shows an example of the two successive minima 
$\lambda_1(\L), \lambda_2(\L)$  for a $2$-dimensional lattice.

\medskip
\begin{example} 
\rm{
For $\L:= \Z^d$, the shortest nonzero vector  has length $\lambda_1(\Z^d) = 1$, and all of the other successive minima for $\Z^d$ have the same value:
$\lambda_2(\Z^d) = \cdots = \lambda_d(\Z^d) =1$.
One choice for their corresponding vectors is 
$v_1:= {\bf e_1}, \dots,  v_d:= {\bf e_d}$, the standard basis vectors.

\hfill $\square$
}
\end{example}

\medskip
\begin{example}\label{Eisenstein lattice}
 \rm{
In $\R^2$, there is a very special
lattice, sometimes called the {\bf hexagonal lattice} (also known as the  {\bf Eisenstein lattice}): 
\[
\L := \left\{ m  { \frac{1}{2}     \choose \frac{\sqrt{3}}{2}  }    + n {1 \choose 0}  
\mid m,n\in\Z  \right\}.
\]
This lattice has $\det \L = \frac{\sqrt{3}}{2}$  and is generated by the $6$'th roots of unity, as in Figure \ref{Eisenstein Lattice} (Exercise \ref{Eisenstein lattice}).  Given the basis above, we see that here we have $\lambda_1(\L) = \lambda_2(\L) =1$.
  It also turns out to be an 
extremal lattice in the sense that it (more precisely a dilate of it) is the lattice that achieves Hermite's constant $\gamma_2$, below, over all lattices in $\R^2$.
(Exercise \ref{minimal lattice in R^2}).  \hfill $\square$
}
\end{example}

\begin{figure}[htb]
 \begin{center}
\includegraphics[totalheight=2.8in]{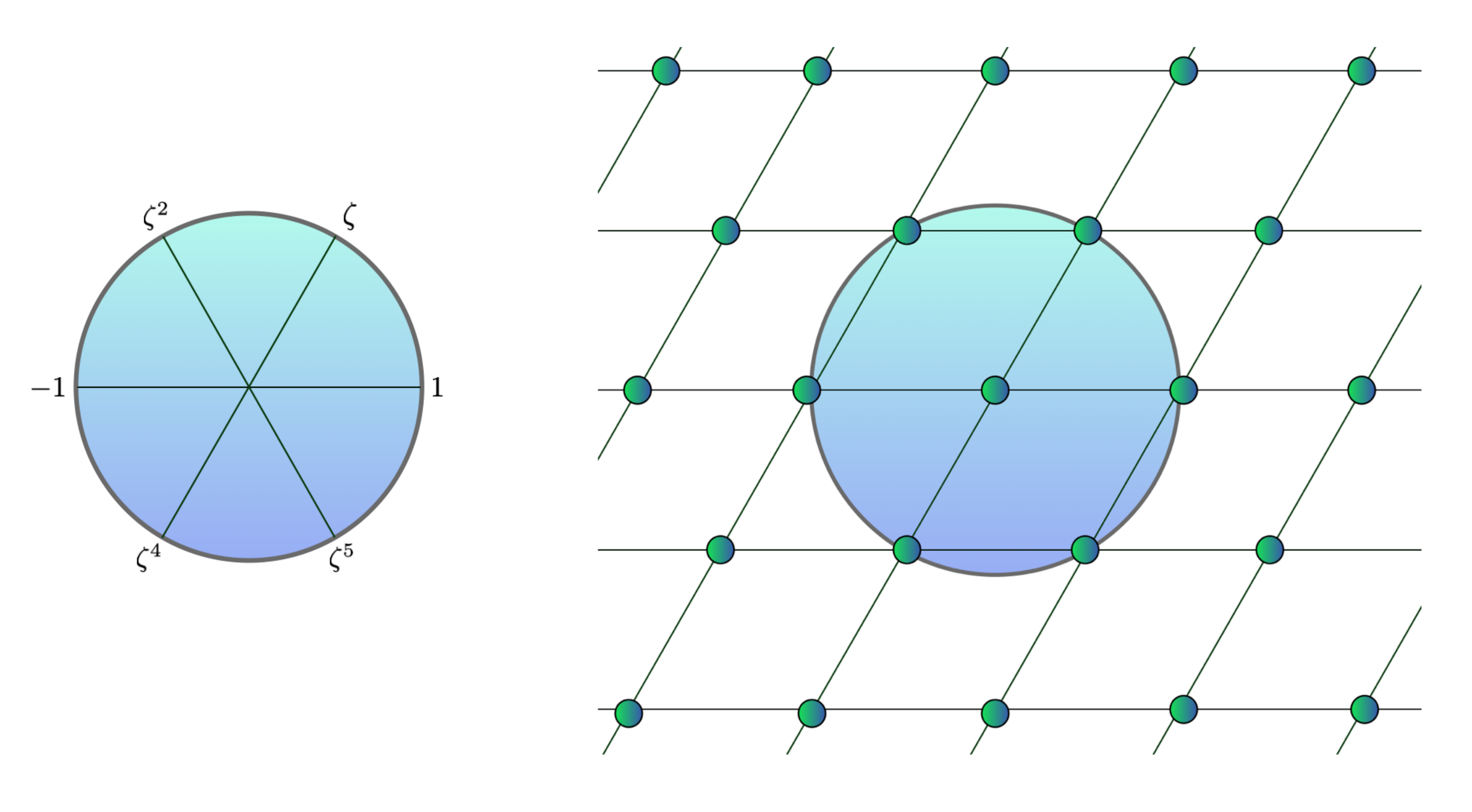}
\end{center}
\caption{Left:  the $6$'th roots of unity.  Right: the hexagonal lattice, also known as the Eisenstein lattice}  
\label{Eisenstein Lattice}
\end{figure}

\medskip
\begin{example} \label{a curve in the space of lattice}
 \rm{
Let's define the following family of $2$-dimensional lattices. For each $t > 0$, we let
\[
 M:=  \begin{pmatrix}
e^t & 0 \\
 &   e^{-t} \end{pmatrix}, \text{ and we let } \L_t:= M(\Z^d),
 \]
 so that we get a parametrized family of lattices.  While all of the lattices in this family have $\det \L=1$, their shortest nonzero vectors approach $0$ as $t\rightarrow \infty$, since $\lambda_1(\L_t) = e^{-t}$.   
 So we see that it does not necessarily make sense to talk about the shortest nonzero vector among a collection of lattices, but it will make sense to consider a ``max-min problem'' of this type 
 (Hermite's constant \eqref{Hermite's constant} below).
 }
 \hfill $\square$
\end{example}

For each dimension $d$, we define {\bf Hermite's constant} as follows:
\begin{equation}\label{Hermite's constant}
\gamma_d :=  \max    \left\{  \lambda_1(\L)^2  \bigm |   \L   \text{ is a full-rank lattice in $\R^d$, with  } \det \L =1       \right\}.
\end{equation}
In words, Hermite's constant is retrieved by varying over all normalized lattices in $\R^d$, which have determinant $1$, picking out the smallest squared norm of any nonzero vector in each lattice, and then
taking the maximum of the latter quantity over all such lattices.  In a later chapter, on sphere packings, we will see an interesting interpretation of Hermite's constant in terms of the densest lattice packing of spheres.

We next give a simple bound, in Theorem \ref{First successive minima bound} below, for the shortest nonzero
vector in a lattice and hence for Hermite's constant.   But first we need to give  
 a simple lower bound for the volume of the unit ball, in  Lemma \ref{volume bound for the ball}.   
 Curiously, Hermite's constant $\gamma_d$ is only known
precisely for $1 \leq d \leq 8$, and $d=24$, as of this writing.

\bigskip
\begin{lem}\label{volume bound for the ball}
\[
\vol B_d(r)  \geq   \left(  \frac{2r}{\sqrt d}  \right)^d.
\]
\end{lem} 
\begin{proof}
The cube $C:= \left\{  x\in \R^d \bigm |    
 \text{ all } |x_k| \leq \frac{r}{\sqrt d}   \right\}$ is contained in the ball $B_d(r)$: 
if $x \in C$ then $\sum_{k=1}^d x_k^2 \leq   d \left(    \frac{r}{\sqrt d}  \right)^2 = r^2$.  So the volume of the ball $B_(r)$ is greater than the volume of the cube, which is equal to $\left(  \frac{2r}{\sqrt d}  \right)^d$.
\end{proof}

\bigskip
The following result of Minkowski gives a bound for the shortest nonzero vector in a lattice.
\begin{thm}[Minkowski]  \label{First successive minima bound}
Suppose that $\L \subset \R^d$ is a full-rank lattice.   Then the shortest nonzero vector $v \in \L$ satisfies
\begin{equation}\label{shortest vector in a lattice bound}
\|v\| \leq    \sqrt{d}      (\det \L)^{\frac{1}{d}}.
\end{equation}
Equivalently, we may write 
\[
\lambda_1(\L) \leq \sqrt{d}      (\det \L)^{\frac{1}{d}}.
\]
\end{thm}
\begin{proof}
The idea is to apply Minkowski's convex body Theorem \ref{Minkowski convex body Theorem for L}
to a ball of sufficiently large radius. 
Let $r:= \lambda_1(\L)$ be the length of the shortest nonzero vector in $\L$, and consider the 
 ball $B_d(r)$  of radius $r$. By definition, $B_d(r)$ does not contain any lattice points of $\L$ in its interior.   
 So by Minkowski's convex body Theorem, and Lemma \ref{volume bound for the ball}, 
 \[
  \left(  \frac{2 \lambda_1(\L)}{\sqrt d}  \right)^d
 \leq
  \vol B_d(r) 
 \leq 
 2^d \det \L.
 \]
 It follows that $\lambda_1(\L) \leq \sqrt{d}  \left( \det \L  \right)^{\frac{1}{d}}$, proving the claim. 
\end{proof}

Despite the  bound \eqref{shortest vector in a lattice bound} on the shortest nonzero vector in a lattice, there are currently no known efficient algorithms to find such a vector for an arbitrary lattice, and it is thought to be one of the most difficult problems we face today.   In practice, researchers often use the LLL algorithm  
  to find a `relatively short' vector in a given lattice, and the same algorithm even finds a relatively short basis for $\L$.   

We already have enough knowledge to relate the length of a shortest nonzero vector of a lattice $\L$ to the length of a
shortest nonzero vector of its dual lattice $\L^*$, as follows. 
\begin{cor}
Let $\L\subset \R^d$ be a full-rank lattice, and let $\L^*$ be its dual lattice.     Then
\begin{equation} \label{cor:transference between shortest vectors in a lattice and its dual}
\lambda_1(\L) \lambda_1(\L^*)  \leq d.
\end{equation}
\end{cor}
\begin{proof}
By Minkowski's bound, namely Theorem \ref{First successive minima bound}, applied to both $\L$ and $\L^*$, we have:
\[
\lambda_1(\L) \lambda_1(\L^*)  \leq           \sqrt{d} (\det \L)^{\frac{1}{d}}  \sqrt{d}  (\det \L^*)^{\frac{1}{d}} = d,
\]
using the relation $(\det \L)( \det \L^*) = 1$.
\end{proof}
Such relations are called {\bf transference theorems}, as they can transfer the complexity of computing a lattice parameter in $\L$ to the complexity of computing a (usually different) parameter in the dual lattice $\L^*$.  In the case of equation 
\eqref{cor:transference between shortest vectors in a lattice and its dual},
we have a quantitative measure of the fact that the shortest vectors in a lattice and its dual lattice cannot both be ``too long".

\bigskip
\section{Minkowski's second theorem}

\begin{thm}[Minkowski's second theorem, for a convex body $K$]  
\label{successive minima for K}
The successive minima of a full-rank lattice $\L$, relative to a convex body $K$, enjoy the property:
\begin{equation}\label{2'nd thm of Minkowski, the inequality}
 \lambda_1(\L, K) \cdots \lambda_d(\L, K)   \vol K  \leq  2^d \det \L.
\end{equation}
\hfill $\square$
\end{thm}
This result, called Minkowski's second theorem, has many proofs, and of course the first was given by Minkowski in 1896 (\cite{Minkowski}, p. 199).  One of the easiest 
(and most clever) ways to see why Theorem \ref{successive minima for K} is true, was given by Henk \cite{Henk4}.  

Let's compare Minkowski's second theorem to Minkowski's first theorem, which may be (easily) rewritten as follows.
\begin{thm}[Minkowski's first theorem for a convex body $K$, equivalent formulation]  
The shortest nonzero lattice point of $\L$, relative to a convex body $K$, enjoys the property:
\begin{equation} \label{1'st thm of Minkowski, equivalent inequality}
 \lambda_1(\L, K)^d  \vol K  \leq  2^d \det \L.
\end{equation}
\hfill $\square$
\end{thm}
We will not deprive the reader of the pleasure of proving this equivalence (Exercise \ref{Rewriting Minkowski's first theorem}).  
Recalling that  $\lambda_1(\L, K)    \leq    \lambda_2(\L, K)  \leq \cdots \lambda_d(\L, K)$,
 it's now apparent that 
\eqref{2'nd thm of Minkowski, the inequality} is in general a huge improvement upon
\eqref{1'st thm of Minkowski, equivalent inequality}.

In this short section, we'll prove a simpler result, for the case of the successive minima of the 
unit ball in $\R^d$, 
namely for $\lambda_j(\L, B):=\lambda_j(\L)$.   
While we may not know explicitly all of the short vectors in a given lattice, it is often still useful to construct 
 an ellipsoid that is based on the successive minima of a lattice. 
 In the spirit of reviewing basic concepts from Linear Algebra, an
 {\bf ellipsoid boundary} \index{ellipsoid} 
 centered at the origin is defined by the $(d-1)$-dimensional body
\begin{equation} \label{ellipsoid}
\left\{    x\in \R^d  \bigm |      \sum_{j=1}^d \frac{{\langle x, b_j\rangle}^2}{c_j^2} =1   \right\},
\end{equation}
for some fixed orthonormal basis
 $\{ b_1, \dots, b_d \}$ of $\R^d$.   Here the vectors $b_j$ are called the {\bf principal axes} of the ellipsoid, and 
 the $c_j$'s are the lengths along the principal axes of the ellipsoid. 
A more geometric way of 
defining an ellipsoid (which turns out to be equivalent to our definition above) is attained 
by applying a linear transformation $M$ to the unit sphere 
$S^{d-1} \subset  \R^d$ (Exercise \ref{Ellipsoid problem}).   For the remainder of this section, we follow the approach taken by Oded Regev \cite{RegevNotes}.

\begin{figure}[htb]
 \begin{center}
\includegraphics[totalheight=2in]{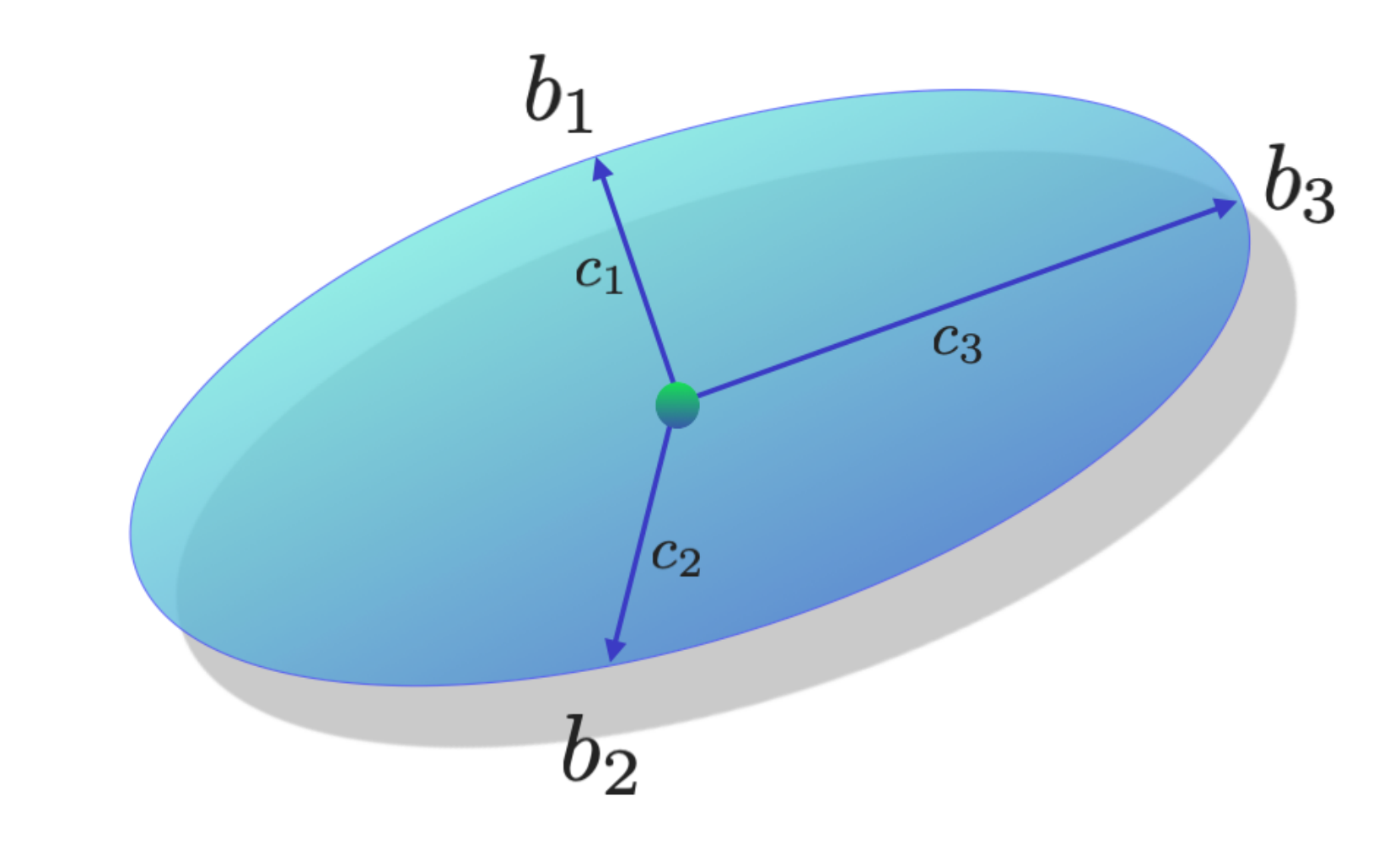}
\end{center}
\caption{An ellipsoid in $\R^3$.}  \label{Ellipsoid}
\end{figure}

Corresponding to the successive minima of a full-rank lattice $\L$, by definition 
we have
 $d$ linearly independent vectors
$v_1, \dots, v_d$, so $\| v_k \| := \lambda_k(\L)$.   We apply the 
Gram-Schmidt algorithm to this set of vectors $\{v_1, \dots, v_d\}$, obtaining a corresponding orthonormal basis 
$\{b_1, \dots, b_d\}$ for $\R^d$.   

Now we define the following open {\bf ellipsoid} by:
\begin{equation}\label{open ellipsoid}
E:=\left\{ x\in \R^d \bigm |     \sum_{k=1}^d     \frac{{\langle x, b_k \rangle}^2}{{\lambda_k}^2}   
< 1
\right\},
\end{equation}
whose axes are the $b_k$'s, and whose radii are the $\lambda_k:= \lambda_k(\L)$. With this notation in hand, we have the following.
\bigskip
\begin{lem}\label{Empty Ellipsoid}
\rm{
 The ellipsoid $E$ does not contain any lattice points of $\L$.
}
\end{lem}
\begin{proof}
We fix any vector $v \in \L$.  
Let $1\leq k \leq d$ be the maximal index such that $\lambda_k(\L) \leq \|v\|$.
We may write $v = \sum_{j=1}^d \langle v, b_j \rangle b_j$, so that 
$\|v\|^2 =  \sum_{j=1}^d    {\langle v, b_j \rangle}^2$.

Now $v$ must lie in $\text{span}\{v_1, \dots v_k\} = \text{span}\{b_1, \dots b_k\}$, for some $1\leq k\leq d$.
Hence we may write $v = \sum_{j=1}^d \langle v, b_j \rangle b_j =  \sum_{j=1}^k \langle v, b_j \rangle b_j $, 
so that 
$\|v\|^2 =  \sum_{j=1}^k | \langle v, b_j \rangle |^2$.
We now check if $v$ is contained in $E$:
\[
\sum_{j=1}^d     \frac{{\langle v, b_j \rangle}^2}{{\lambda_j}^2} =
\sum_{j=1}^k    \frac{{\langle v, b_j \rangle}^2}{{\lambda_j}^2}   \geq
 \frac{1}{ {\lambda_k}^2  }      \sum_{j=1}^k   {\langle v, b_j \rangle}^2
= \frac{  \|v\|^2 }{   {\lambda_k}^2  }  \geq 1,
\]
so that $v\not\in E$.
\end{proof}

More generally, it's easy to prove the following refinement of 
Theorem  \ref{First successive minima bound}, which  gives us a bound for the first $d$ shortest (nonzero) vectors in a lattice.   
\begin{thm} \label{successive minima bound}
The successive minima of a full-rank lattice $\L$ enjoy the property:
\[
\lambda_1(\L)    \cdots \lambda_d(\L) 
          \leq d^{\frac{d}{2}} \det \L.
\]
\end{thm}
\begin{proof}  Using Lemma \ref{Empty Ellipsoid}, the ellipsoid $E$  contains no lattice points belonging to 
$\L$, so that by Minkowski's convex body Theorem, we have $\vol E \leq 2^d \det \L$.  
We also know that 
\[
\vol E = \left( \prod_{j=1}^d \lambda_j \right) \vol B_1      \geq 
\left( \prod_{j=1}^d \lambda_j \right)   \left(\frac{2}{\sqrt{d}}   \right)^d.
\]
Altogether, we have
\[
  2^d \det \L \geq \vol E \geq \left( \prod_{j=1}^d \lambda_j \right)   \left(\frac{2}{\sqrt{d}}   \right)^d,
\]
arriving at the desired inequality.
\end{proof}

\section{The distance function of a body, and its support function}
 \index{support function} \index{distance function}

There is a natural correspondence between norms on $\R^d$ and convex, centrally-symmetric bodies in $\R^d$.   A {\bf norm} \index{norm} 
on $\R^d$ is a function $g: \R^d \rightarrow \R_{\geq 0}$ that enjoys the following properties:

\begin{enumerate}
\item Triangle inequality (Subadditivity): 
$g(x+y)  \leq  g(x)+g(y) $,  for all $x, y \in \R^d$. 
\item
Homogeneity: $g(rx) = |r| g(x)$, for all $r\in \R$ and $x \in \R^d$.
\item
Positive definiteness:  for any $x\in \R^d$,  $g(x) = 0 \iff x=0$. 
\end{enumerate}
Given a centrally symmetric convex body $\K\subset \R^d$,  there is a natural norm we can associate with $\K$.  Namely, we define 
\begin{equation}\label{def:distance function}
|x|_\K:=  \min \left\{ r \geq 0   \mid x \in r\K \right\},
\end{equation}
called the {\bf distance function} \index{distance function}
of $\K$.
We recall that $r\K:= \{ rx \mid x \in \K\}$ is the $r$'th dilation of $\K$. It's straightforward to prove that $|\cdot|_\K$ is indeed a norm (Exercise \ref{proof that the gauge function is a norm}), when $\K$
 is convex and centrally-symmetric.  Moreover, the unit ball of the norm $|x|_\K$ is $\K$ itself.
\begin{example}
When $\K:= B$, the unit ball in $\R^d$, we have $|x|_\K = \| x \| := \sqrt{x_1^2 + \cdots + x_d^2}$, the usual norm on $\R^d$.
\hfill $\square$
\end{example}
In some of the literature, the distance function $|\cdot|_\K$ of a body $\K$ is also called the gauge function of $\K$ (for example, in Siegel's book \cite{SiegelBook}).

\begin{lem}
Suppose $\K\subset \R^d$ is a $d$-dimensional convex body that contains the origin. Then its volume is equal to
\begin{equation}
\vol \K = \frac{1}{d} \int_{S^{d-1}}  \left( \frac{1}{ |x|_\K } \right)^d dx.
\end{equation}
\end{lem}
\begin{proof}
Since $\K$ contains the origin, a ray emanating from the origin and passing through any given point 
$x \in S^{d-1}$ must intersect $\K$ in a unique point $y \in \partial \K$, the boundary of $\K$.  
By definition, we have $y =\frac{1}{ |x|_\K} x $.  Defining $R(x):= \| y\|$, we have 
$R(x) = \frac{1}{ |x|_\K} \|x\| =  \frac{1}{ |x|_\K}$, which is a continuous function on the unit sphere because $R(x)$ does not vanish there.  Using polar coordinates in $\R^d$, namely $y= rx$, where $x \in S^{d-1}$
and $ r \in \R_{\geq 0}$, we now have
\begin{align*}
\vol \K := \int_\K dy &=  \int_{S^{d-1}} \int_0^{R(x)} r^{d-1} dr dx
=\frac{1}{d}  \int_{S^{d-1}} R^d(x) dx \\
&= \frac{1}{d} \int_{S^{d-1}}  \left( \frac{1}{ |x|_\K } \right)^d dx.
\end{align*}
\end{proof}

In general, we also have the following basic fact:  there is a one-to-one correspondence between norms and centrally symmetric convex bodies, given by the following mapping.
\begin{thm}\label{norms and convex bodies}
We let $K_o^d$ be the collection of all convex, centrally-symmetric, $d$-dimensional bodies 
$\K\subset\R^d$.   We define the mapping
\begin{equation}
\Phi: K_0^d \rightarrow \left \{  \text{norms on } \R^d \right\}
\end{equation}
by $\Phi( \K ) := | \cdot |_\K$, the distance function of the body $\K$.  
Then $\Phi$ is a  $1-1$, onto correspondence.  In particular, given any norm $g$ on $\R^d$, 
 $\Phi^{-1}(g)$ is the unit ball for the norm $g$.
\end{thm}
\hfill $\square$

The proof of Theorem \ref{norms and convex bodies} is fun, so we leave it as Exercise 
\ref{exercise:norms and convex bodies}.  It follows from Theorem \ref{norms and convex bodies} that each convex body is uniquely determined by 
its distance function.
Next, there is a related gadget called the {\bf support function} for each convex body $K\subset \R^d$, defined by
\begin{equation}
h_K(x):= \max \left \{ \langle x, k\rangle \mid k \in K \right \}.
\end{equation}

\begin{example}
For the unit ball $B\subset \R^d$, its support function is $h_B(x) = \|x\|$, the length of $x$. 
For the unit cube $\square:= \left[ -\tfrac{1}{2}, \tfrac{1}{2} \right]^d$,
\hfill $\square$
\end{example}
Some of the nice properties of the support function of any convex body, which follow quickly from the definition, include:
\begin{enumerate}[(a)]
\item (homogeneity) $h_K(\lambda x) = \lambda h_K(x)$, for all $\lambda>0$.
\item (subadditivity) $h_K(x+y) \leq h_K(x) + h_K(y)$, for all $x, y \in \R^d$.
\item $h_K$ is a continuous function of $x\in \R^d$. 
\end{enumerate}
In fact, much more is true:
\begin{lem}\label{connection between support function and distance function}
 If $K\subset \R^d$ is a centrally-symmetric convex body containing the origin, then 
 its polar 
$K^o$ is also a centrally-symmetric convex body containing the origin.  Moreover:
\begin{equation}
h_K(x)=|x|_{K^o}.
\end{equation}
\hfill $\square$
\end{lem}
It follows from Lemma \ref{connection between support function and distance function}, together with Theorem \ref{norms and convex bodies}
that every centrally-symmetric convex body is uniquely determined by its support function (although it's also easy to prove it from scratch).


\section{The theta function of a lattice}
\label{section:theta function of a lattice}
\index{theta function of a lattice}

There is a  beautiful analytic function that we can associate to a given full-rank lattice $\L \subset \R^d$, by using Poisson summation. 
For each fixed $t>0$, and each $y\in \R^d$, we define
\begin{equation} \label{Theta function of a lattice}
\theta_\L(t, y):= \sum_{n\in \L} e^{-\pi t \|y-n\|^2},
\end{equation}
called the theta function of the lattice $\L$.   These theta functions have a
 very rich and long history.  In 1859, Bernhard Riemann used $\theta_\Z(t, 0)$, together with its functional equation, to prove the functional equation for the Riemann zeta function $\zeta(s):= \sum_{n\geq 1} \frac{1}{n^s}$ (via the Mellin transform).  In Section \ref{The covering radius of a lattice, and its packing radius}, we will use the theta function \eqref{Theta function of a lattice} to relate the packing radius of a lattice to the covering radius of its dual lattice.  Such relations use the following basic and classical functional equation.

 \begin{thm}\label{Functional equation for theta}
 For any full-rank lattice $\L \subset \R^d$, and any fixed $y \in \R^d$, 
 we have the functional equation
 \begin{equation}\label{thm: functional equation for theta}
t^{\frac{d}{2}}  \theta_\L(t, y) =   \frac{1}{\det \L} \sum_{n\in \L^*}   
e^{-\frac{1}{t} \pi  \|n\|^2 + 2\pi i \langle n, y \rangle}.
 \end{equation}
 \end{thm}
\begin{proof}
The proof of the functional equation for the theta function \eqref{thm: functional equation for theta} is a simple application of Poisson summation for Schwartz functions (Theorem \ref{Poisson.Summation}), applied to the Gaussian  
$G_t(n, y):= e^{ t\pi  \|y-n\|^2}$:
\begin{align}
\theta_\L(t, y) &:= \sum_{n\in \L} e^{-\pi t \|y-n\|^2} := \sum_{n\in \L} G_t(n, y)    \\ 
\label{theta application of Poisson summation}
&=  \frac{1}{\det \L} \sum_{\xi\in \L^*} \hat G_t(\xi, y) \\
\label{theta function eq, FT of Gaussian}
&= \frac{1}{\det \L} \sum_{n\in \L^*}   
t^{-\frac{d}{2}}   e^{-\frac{1}{t} \pi  \|n\|^2 + 2\pi i \langle n, y \rangle},
\end{align}
the desired identity.   In \eqref{theta application of Poisson summation} 
we used Poisson summation, and in 
\eqref{theta function eq, FT of Gaussian} we used the Fourier transform of the Gaussian.
\end{proof}
In the following section, namely Section \ref{The covering radius of a lattice, and its packing radius}, we will use these theta functions to derive bounds on some fundamental
 lattice parameters.

\section{The covering radius of a lattice, 
\\
and its packing radius, via Poisson summation}
\label{The covering radius of a lattice, and its packing radius}

Throughout this section, we are given a convex body $K \subset \R^d$, containing the origin, and a full-rank lattice $\L\subset \R^d$.
The {\bf covering radius of the lattice $\L$} \index{covering radius} relative to $K$, is defined by the
smallest $r>0$ such that every point $x\in \R^d$ is covered by some translate of $rK$ by a vector from the lattice $\L$.
More compactly, we may also give the following description for the covering radius:
\begin{align}
\mu(\L, K)
&:= \min \{ r \geq 0   \mid   r K + \mathcal L = \mathbb{R}^d \}.
\end{align}
The most common scenario occurs when $K:= B$, the unit ball, and in this case it's traditional to use the following shorter notation for the covering radius of a lattice:
\[
\mu(\L, B):= \mu(\L).
\]
In words,  $\mu(\L)$ is the smallest $\mu >0$ such that the collection of open balls of radius $\mu$, 
centered at all lattice points of $\L$, completely covers $\R^d$.  It's useful to define for any set $S\subset \R^d$ and any point $x\in \R^d$, 
the distance \index{distance}
\begin{equation}
\text{ dist}(x, S) := \inf_{s \in S } \| x - s \|.
\end{equation}
It follows directly from the definitions above that 
\begin{equation}
\mu(\L)= \max_{x \in \R^d}\rm{dist}\left( x, \L \right).
\end{equation}

\begin{example}
\rm{
For the integer lattice $\Z^d$, the covering radius $\mu(\Z^d)$ (relative to the unit ball) is
maximum distance between any point $x\in \R^d$ and the nearest integer point.  This distance is clearly achieved by the point 
$\left(\tfrac{1}{2}, \dots, \tfrac{1}{2} \right)^T$, which is the centroid of the unit cube $[0, 1[^d$, and gives us
$\mu(\Z^d) = \frac{\sqrt d}{2}$.
}
\hfill $\square$
\end{example}

In a somewhat dual fashion,  the
{\bf packing radius of the lattice $\L$}  \index{packing radius}
relative to $K$, written as $\rho(\L, K)$,
is the largest $r>0$ such that 
$(r K + l_1) \cap (r K + l_2) \not= \phi$, for all $l_1, l_2 \in \L$.    

   In words, the packing radius $\rho(\L)$ is the largest $r > 0$ such that the collection of 
open balls of radius $r$, centered at all lattice points of $\L$,  do not intersect.  
When $K:= \interior(B)$, the open unit ball, it is traditional to omit $K$ in the notation, and 
 we simply write the packing radius in this case as 
 \[
 \rho(\L, B):= \rho(\L).
 \]

It follows from the definitions that the packing radius of a lattice $\L$ equals precisely half the  distance to the shortest vector of $\L$:
\begin{equation}\label{simple relation between shortest vector and packing radius}
\rho(\L) = \tfrac{1}{2} \lambda_1(\L).
\end{equation}
\begin{example}
\rm{
For the integer lattice $\Z^d$, the packing radius $\mu(\Z^d)$ (relative to the unit ball) is
just $\rho(\Z^d) = \frac{1}{2}$.
}
\hfill $\square$
\end{example}

  There are certain useful dualities between $\mu(\L)$ and $\rho(\L^*)$, known as `transference theorems'.   These results `transfer' the problem of computing certain lattice parameters of $\L$ to the problem of computing certain other parameters of its dual lattice $\L^*$.
 Here we prove such a `transference theorem', discovered by Banaszczyk \cite{Banaszczyk}, which is another application of Poisson summation and theta functions.  
 In this section we'll follow the approach taken in the lecture notes of Oded Regev  \cite{RegevNotes}.

 \begin{thm}[Banaszczyk, 1993] \label{strong transference theorem of Banaszczyk}
 For any full-rank lattice $\L \subset \R^d$, we have
 \begin{equation}
 \mu(\L^*) \lambda_1(L) \leq \frac{d}{2}.
 \end{equation}
 \hfill $\square$
 \end{thm}

 Following Regev's notes \cite{RegevNotes} (with only tiny modifications) we prove here a result that is almost as good,  but with a weaker constant, as follows.
 \begin{thm} \label{transference theorem of Banaszczyk}
 For any full-rank lattice $\L \subset \R^d$, we have
 \begin{equation}
 \mu(\L^*) \lambda_1(L) \leq d.
 \end{equation}
 \hfill $\square$
 \end{thm}
A good (and elementary) exercise is to relate the Voronoi cell $V(\L)$ \index{Voronoi cell}
of a lattice to the packing radius and the covering radius of $\L$.  Namely, the packing radius $\rho(L)$ equals the inradius (radius of largest inscribed sphere) of $V(\L)$, 
and the covering radius $\mu(\L)$ equals the circumradius (radius of smallest circumscribed sphere) of $V(L)$ (Exercise \ref{Voronoi cell and rho and mu}).
From the observation \eqref{simple relation between shortest vector and packing radius}, Theorem \ref{transference theorem of Banaszczyk} can also be trivially restated as 
\begin{equation}
 \mu(\L) \lambda_1(\L^*) \leq d.
\end{equation}

Based on Poisson summation, we already derived a functional equation for the theta function of a lattice, namely Theorem \ref{Functional equation for theta}.  
Following \cite{RegevNotes}, we set up the proof of Theorem \ref{transference theorem of Banaszczyk} by proving a few self-contained and useful lemmas. 
 
\begin{lem}\label{first lemma: section covering radius}
Fix any $x\in \R^d$, and $t >0$. Then we have:
\begin{equation}
 \theta_\L(t, x):=\sum_{n \in \L} e^{-\pi t \| x+n\|^2}   \leq   \sum_{n\in \L} e^{-\pi t \| n\|^2}:= \theta_\L(t, 0),
\end{equation}
with equality if and only if $x \in \L$.
\end{lem}
\begin{proof}
Using Theorem \ref{Functional equation for theta},  the functional equation 
\eqref{thm: functional equation for theta} of the theta function $ \theta_\L(t, x)$ gives us: 
\begin{align*}
t^{\frac{d}{2}}  \theta_\L(t, x) &=   
\frac{1}{\det \L} \sum_{n\in \L^*}      e^{-\frac{1}{t} \pi  \|n\|^2 + 2\pi i \langle n, x \rangle}\\
&=   
\frac{1}{\det \L} \sum_{n\in \L^*}      e^{-\frac{1}{t} \pi  \|n\|^2}
\cos\left(  2 \pi \langle n, x \rangle \right)      
+i \frac{1}{\det \L} \sum_{n\in \L^*}      e^{-\frac{1}{t} \pi  \|n\|^2}
\sin \left(  2 \pi \langle n, x \rangle \right)                            \\
&=   
\frac{1}{\det \L} \sum_{n\in \L^*}      e^{-\frac{1}{t} \pi  \|n\|^2}
\cos\left(  2 \pi \langle n, x \rangle \right)                          \\
&
\leq \frac{1}{\det \L} \sum_{n\in \L^*}      e^{-\frac{1}{t} \pi  \|n\|^2}  \\
&=  t^{\frac{d}{2}}   \sum_{n\in \L} e^{-\pi t \| n\|^2}:= t^{\frac{d}{2}}   \theta_\L(t, 0).                 \\
\end{align*}
We used the functional equation in the first equality and in the penultimate equality above.  The equality condition occurs precisely when $\cos\left(  2 \pi \langle n, x \rangle \right)  =1$, which in turn occurs if and only if 
$\langle n, x \rangle\in \Z$ for all $n \in \L^*$.  Finally, the latter condition holds precisely when 
$x \in \left(\L^*\right)^*$, but we already know that 
$\L^{**} = \L$.
\end{proof}

\begin{lem}\label{second lemma: section covering radius}
Fix any $0\leq t \leq 1$. Then we have:
\begin{equation}
\sum_{n \in \L} e^{-\pi t \| n \|^2}   \leq  t^{-\frac{d}{2}}  
\sum_{n\in \L} e^{-\pi  \| n\|^2}.
\end{equation}
\end{lem}
\begin{proof}
Using the functional equation \eqref{thm: functional equation for theta}
for the theta function again, in both equalities below, we have:
\begin{align*}
 t^{\frac{d}{2}}    \sum_{n \in \L} e^{-\pi t \| n \|^2} &= \frac{1}{\det \L} \sum_{n\in \L^*}  e^{-\frac{1}{t} \pi  \|n\|^2}  \\
 &\leq \frac{1}{\det \L} \sum_{n\in \L^*}  e^{\pi  \|n\|^2} = \sum_{n \in \L} e^{-\pi  \| n \|^2}.
\end{align*} 
We used $e^{-\frac{1}{t} \pi  \|n\|^2} \leq e^{-\pi  \|n\|^2}$ in the inequality above, which is valid for $t \in [0, 1]$.
\end{proof}
To summarize, we now know that in the range $t \in [0, 1]$, the previous two lemmas together give us:
\begin{equation}\label{inequality to be used in third lemma below}
\sum_{n \in \L} e^{-\pi t \| x+n\|^2}   \leq   t^{-\frac{d}{2}} \sum_{n\in \L} e^{-\pi  \| n\|^2}.
\end{equation}

We recall that $B_d(r)\subset \R^d$ is the $d$-dimensional ball of radius $r$, centered at the origin.
The next lemma tells us that for any full-rank lattice $\L$, most of its contribution to its own 
theta function is already contained in those lattice points that belong to $B_d(\sqrt d)$.
\begin{lem}   \label{third lemma: section covering radius}
\begin{equation}
 \sum_{n \in \left(   \L +x  \right)  \text{ and } \|n\| \geq \sqrt d} 
        e^{-\pi \|n\|^2} 
 < 
\frac{1}{2^{\tfrac{3d}{2}} } \sum_{n \in \L}  e^{-\pi \| n \|^2}.
\end{equation}
\end{lem}
\begin{proof}
We'll bound $\theta_{\L}(\frac{1}{2}, x)$ from above and below.   For the upper bound, 
\eqref{inequality to be used in third lemma below} gives us
\begin{equation*}
\theta_{\L}\left(\tfrac{1}{2}, x\right):= 
\sum_{n \in \L} e^{-\pi \tfrac{1}{2} \| x+n\|^2}   \leq   2^{\frac{d}{2}} \sum_{n\in \L} e^{-\pi  \| n\|^2}.
\end{equation*}
For the lower bound, we have
\begin{align*}
\theta_{\L}\left(\tfrac{1}{2}, x\right) 
&\geq
 \sum_{n \in \left(   \L +x  \right)   \setminus B_d\left(\sqrt d \right)} 
           e^{-\pi \tfrac{1}{2} \|n\|^2} 
       = \sum_{n \in \left(   \L +x  \right)  \text{ and } \|n\| \geq \sqrt d} 
           e^{-\pi \tfrac{1}{2} \|n\|^2} \\
&= \sum_{n \in \left(   \L +x  \right)  \text{ and } \|n\| \geq \sqrt d} 
         e^{\frac{\pi}{2}  \|n\|^2}  e^{-\pi \|n\|^2} \\\
&\geq    e^{\frac{\pi}{2} d}  \sum_{n \in \left(   \L +x  \right)  \text{ and } \|n\| \geq \sqrt d} 
        e^{-\pi \|n\|^2} \\
&>    4^d \sum_{n \in \left(   \L +x  \right)  \text{ and } \|n\| \geq \sqrt d} 
        e^{-\pi \|n\|^2}.
 \end{align*}
 Together, these upper and lower bounds give
\begin{equation*}
 \sum_{n \in \left(   \L +x  \right)  \text{ and } \|n\| \geq \sqrt d} 
        e^{-\pi \|n\|^2} 
 < 
2^{\frac{d}{2} - 2d} \sum_{n \in \L}  e^{-\pi \| n \|^2} 
       = \frac{1}{2^{\tfrac{3d}{2}} }
 \sum_{n \in \L}  e^{-\pi \| n \|^2}.
\end{equation*}
\end{proof}

A direct consequence of Lemma \ref{third lemma: section covering radius} is the following simple bound on theta functions.
\begin{lem}   \label{fourth lemma: section covering radius}
If the lattice $\L$ enjoys $\lambda_1 > \sqrt d$, then:
\begin{equation}
 \sum_{n \in \L\setminus \{0\}}   e^{-\pi \| n \|^2} 
 \leq  
 \frac{1}{2^{\frac{3d}{2}} }.
\end{equation}
\end{lem}
\begin{proof}
The assumption that $\lambda_1 > \sqrt d$ means that 
$ \sum_{n \in \left(   \L  \right)  \text{ and } \|n\| \geq \sqrt d} 
        e^{-\pi \|n\|^2} 
     =
\sum_{n \in \L\setminus \{0\}}
        e^{-\pi \|n\|^2}$.
Now Lemma \ref{third lemma: section covering radius} with $x=0$ immediately implies that
\begin{equation}
\sum_{n \in \L\setminus \{0\}}
        e^{-\pi \|n\|^2} 
=
 \sum_{n \in \left(   \L  \right)  \text{ and } \|n\| \geq \sqrt d} 
        e^{-\pi \|n\|^2} 
\leq
        \frac{1}{2^{\tfrac{3d}{2}} } \sum_{n \in \L}  e^{-\pi \| n \|^2}
=
\frac{1}{2^{\tfrac{3d}{2}} } 
\left(
\sum_{n \in \L\setminus \{0\}}  e^{-\pi \| n \|^2} +1
\right)
\end{equation} 
Solving for $z:=\sum_{n \in \L\setminus \{0\}}  e^{-\pi \| n \|^2}$, we have 
$z\leq  \frac{1}{2^{\tfrac{3d}{2}} } ( z+1)$, and hence the desired inequality.
\end{proof}

\begin{lem}   \label{fifth lemma: section covering radius}
If the lattice $\L$ enjoys $\lambda_1 > \sqrt d$, then:
\begin{equation}\label{nearly constant theta function}
  \det \L - \det \L\left(    \frac{1}{2^{\frac{3d}{2}} }  \right)
  \leq
   \sum_{n \in \L^*}   e^{-\pi \| n + x\|^2} 
 \leq  
  \det \L + \det \L\left(    \frac{1}{2^{\frac{3d}{2}} }  \right),
\end{equation}
for all $x\in \R^d$.
\end{lem}
\begin{proof}
Starting with Poisson summation, but this time exchanging the roles of $\L$ and $\L^*$, we have:
\begin{align*}
 \sum_{n \in \L^*}   e^{-\pi \| n + x\|^2} 
 &= \det \L \sum_{n \in \L}  e^{-\pi \| n \|^2} e^{2\pi i \langle n, x \rangle}
 =\det \L \left( 1+   \sum_{n \in \L\setminus \{0\}}e^{-\pi \| n \|^2} \cos{2\pi \langle n, x \rangle}\right) \\
 &\leq 
    \det \L \left( 1+   \sum_{n \in \L\setminus \{0\}}e^{-\pi \| n \|^2} \right) \\
& \leq  
    \det \L + \det \L\left(    \frac{1}{2^{\frac{3d}{2}} }  \right),
\end{align*}
using Lemma \ref{fourth lemma: section covering radius} in the last inequality.

\end{proof}
Intuitively, Lemma \ref{fifth lemma: section covering radius} tells us that under the hypothesis that a shortest nonzero vector
of $\L$ is `not too short', the theta function $ \sum_{n \in \L^*}   e^{-\pi \| n + x\|^2} $ 
is nearly equal to the constant $\det \L$, especially as the dimension $d$ grows.

\bigskip \bigskip
(Proof of Theorem \ref{transference theorem of Banaszczyk}) 

We'll proceed to give a proof by contradiction, so we'll assume there exists a lattice $\L$ such that
$\mu(\L^*) \lambda_1(\L)  > d$.  Because we may rescale, without loss of generality we assume that $\lambda_1(\L)> \sqrt d$ and $\mu(\L^*)> \sqrt d$.

The assumption that $\mu(\L^*)> \sqrt d$ tells us that there exists $x \in \R^d$ such that 
$\rm{dist}(x, \L^*) > \sqrt d$, which is equivalent to $\rm{dist}(0, \L^*-x) > \sqrt d$. 
Now by Lemma \ref{third lemma: section covering radius} (applied to $\L^*$ instead of $\L$), 
we have
\begin{align*}
   \sum_{n \in \left(  \L^* -x  \right) }
        e^{-\pi \|n\|^2} 
&= \sum_{n \in \left(  \L^* -x  \right)  \text{ and } \|n\| \geq \sqrt d} 
        e^{-\pi \|n\|^2}   \\
&<  
\frac{1}{2^{\tfrac{3d}{2}} } \sum_{n \in \L^*}  e^{-\pi \| n \|^2}
=\frac{1}{2^{\tfrac{3d}{2}} }
\left( 
1+  \sum_{n \in \L^*\setminus \{0\}}  e^{-\pi \| n \|^2}
\right) \\
&\leq 
 \frac{1}{2^{\tfrac{3d}{2}} }  \left( 1 + \frac{1}{2^{\tfrac{3d}{2}} } \right),
\end{align*}
where we've used Lemma \ref{fourth lemma: section covering radius} in the last inequality above.
But now as $d \rightarrow \infty$, the theta function 
 $\sum_{n \in \left(  \L^* -x  \right) }
        e^{-\pi \|n\|^2}$ approaches the zero function, contradicting the left-hand side of 
\eqref{nearly constant theta function} in Lemma \ref{fifth lemma: section covering radius}.
\hfill $\square$

Theorem \ref{strong transference theorem of Banaszczyk}
 gives the best-possible inequality of this type, up to a multiplicative constant.   
 Conway and Thompson (\cite{ConwayThompson}, p. 46) 
 have shown that there exist self-dual lattices with the property that 
$\lambda_1(\L)^2 \geq \frac{d}{2\pi e}\left( 1 + o(1) \right)$, as $d \rightarrow \infty$.


\bigskip
\section{Mordell's measure of non-convexity} \label{Mordell: non-convexity}
\index{asymmetry} \index{Mordell}

A lot of the main results in the geometry of numbers assume that a body $\P$ is convex, which is equivalent to $\P + \P \subset 2\P$.  But it is of great interest to study non-convex bodies as well.  Given a body $\P\subset \R^d$, not necessarily convex, Mordell introduced \cite{Mordell.convexity.measure} a natural measure for the lack of convexity of $\P$, as follows.  If there is a constant $r \geq 2$ such that 
\begin{equation}
\P + \P \subset r\P,
\end{equation}
we'll call $\P$ {\bf quasi-convex, at level $r$}.  So by definition a quasi-convex body at level $r=2$ is also convex.   With these definitions, we can state and prove Mordell's result.  Mordell's proof contains such beautiful ideas that we simply could not resist.
\begin{thm}[Mordell, 1935]
Let $\P\subset \R^d$ is a centrally-symmetric body that is not necessarily convex, but that is quasi-convex at level $r$.
If 
\begin{equation}\label{assymption;Mordell}
\vol \P \geq r^d, 
\end{equation}
then $\P$ contains a nonzero integer point in its interior.
\end{thm}
\begin{proof}
We'll prove the case that $\vol \P > r^d$  (for the equality case we refer the reader 
to Mordell's paper \cite{Mordell.convexity.measure}).  
For any positive integer $m$, 
$\vol\left( \frac{m}{c} \P\right) = \frac{m^d}{r^d} \vol \P$.  This implies that the number $N$ of integer points in the interior of 
$\frac{m}{r}\P$ satisfies
\begin{equation}\label{key inequality;Mordell}
N:=\left| \frac{m}{r} \P^{\interior} \cap \Z^d \right|  \sim  \frac{m^d}{r^d} \vol \P > m^d,
\end{equation}
where the last inequality follows by assumption.  We fix a sufficiently large $m$ for which 
$N > m^d$ holds.

Now we reduce each of the  $N$ integer points
(in the interior of $\frac{m}{r} \P$) modulo $m$, in each coordinate, giving us more than $m^d$ integer points in $[0, m-1]^d \cap \Z^d$ (because $N > m^d$). 
By the pigeon-hole principle, there exist at least two of these integer points, say $p, q \in  \frac{m}{r} \P^{\interior}$, with $p \not= q$, that are congruent to each other mod $m$.   So we know that $\frac{1}{m}(p-q)$ is an integer point.  

Next, the quasi-convexity of $\P$ at level $r$ gives us 
$\frac{m}{r}\P^{\interior} + \frac{m}{r}\P^{\interior} 
     \subset \frac{m}{r}r\, \P^{\interior}= m\P^{\interior}$, so that:
\begin{align*}
\frac{1}{m}(p-q) &\in   \frac{1}{m} \left( \frac{m}{r} \P^{\interior} - \frac{m}{r} \P^{\interior}\right) \\
&=   \frac{1}{m} \left( \frac{m}{r} \P^{\interior} + \frac{m}{r} \P^{\interior}\right) \\
&\subset \frac{1}{m} m\P^{\interior} = \P^{\interior}.
\end{align*}
So we've found a nonzero integer point in $\P^{\interior}$, proving the result.
\end{proof}


\bigskip
\section{The Minkowski conjecture} \label{the Minkowski conjecture}

One of the most imporant open problems in number theory is a long-standing conjecture due to Minkowski, known as \emph{the Minkowski conjecture}, regarding products of linear forms.   
 \begin{conjecture}[The Minkowski conjecture]
 \label{Minkowski conjecture}
 For each $1\leq k \leq d$, let 
 \[
 L_k(x):= a_{k, 1}x_1 + \cdots + a_{k, d}x_d
 \]
  be a
 linear form with real coefficients $a_{i, j}$.  Suppose that the matrix formed by these coefficients is invertible:  
 $\Delta:= \det(a_{i,j}) \not=0$.  Then for any given real numbers 
 $c_1, \dots, c_d$, there exists an integer vector $n\in \Z^d$ such that 
 \begin{equation}\label{Minkowski conjecture, inequality in coordinates}
 \left| (L_1(n) + c_1)\cdots(L_d(n)+c_d) \right|
  \leq \frac{\Delta}{2^d}.
 \end{equation}
  \hfill $\square$
 \end{conjecture}
 Minkowski proved the case $d=2$,  Robert Remak proved the case $d=3$ \cite{Remak2},  
 Freeman Dyson \cite{Dyson} proved the case $d=4$,   Skubenko \cite{Skubenko} and Bambah and Woods \cite{BambahWoods} proved the case $d=5$,   
 Curtis McMullen \cite{CurtMcMullen} proved the case $d=6$, and by now Minkowski's conjecture has been proved up to dimension $d=10$ \cite{Kathuria}. 
C. McMullen's approach is more modern in the sense that he used  ideas from homogeneous dynamics. 
But Minkowski's Conjecture \ref{Minkowski conjecture} 
still appears to be beyond the reach of current methods.

Let's give an equivalent formulation of Conjecture \ref{Minkowski conjecture} 
in terms of the geometry of coverings.
 Following the philosophy of Section \ref{The covering radius of a lattice, and its packing radius}, suppose we are given a body $K\subset \R^d$.  We call a lattice $\L \subset \R^d$ a 
 {\bf covering lattice for $K$} if 
 \[
 K + \L = \R^d.
 \]
 In other words, $\L$ is a covering lattice for $K$ if all of the lattice translates of $K$ 
 cover the full space $\R^d$.  Now suppose we're given any function 
 $f:\R^d \rightarrow \R$, and we consider the set 
 \begin{equation}\label{particular choice of K for covering lattice}
 K:= \{ x\in \R^d \mid  |f(x)| \leq 1 \}.
 \end{equation}
 We have the following easy equivalence, which we leave as 
 Exercise \ref{proof of lemma for covering lattice of K}.
 
 \begin{lem}\label{lem: equivalence for covering lattice of K}
For the set $K$ defined above,  the following conditions are equivalent:
\begin{enumerate}[(a)]
 \item $\L$ is a covering lattice for $K$. 
 \item   \label{second part of lattice covering of K}
 For any $x\in \R^d$, there exists $n\in \L$ such that 
 $\left | f(n+x) \right | \leq 1$.
 \end{enumerate}
 \hfill $\square$
 \end{lem}
 Restating part \ref{second part of lattice covering of K} using coordinates, we see that an equivalent condition for $\L$ to be a covering lattice for the set $K$ in 
 \eqref{particular choice of K for covering lattice} 
 is that for any $x \in \R^d$, and any basis $v_1, \dots, v_d$ of $\L$, there are integers 
 $n_1, n_2, \dots, n_d$ such that $\left | f(n_1  v_1 + \cdots + n_d v_d + x) \right | \leq 1$.
It follows that we may recast Minkowski's conjecture in the following form.

 \begin{conjecture}[The Minkowski conjecture, an equivalent geometric formulation]
Any full-rank lattice $\L \subset \R^d$ is a covering lattice for the set 
 \[
 K:= \left \{ x \in \R^d \mid \left | x_1 x_2 \cdots x_d \right | \leq \frac{ \det \L}{2^d} \right \}.
 \]
  \hfill $\square$
\end{conjecture}
 
It's also easy to see that the equality conditions in Minkowski's conjecture are achieved by the diagonal linear forms $L_k(x_1, \dots, x_d):= 2 c_k x_k$ (Exercise \ref{equality condition for Minkowski's conjecture}).   Minkowski further conjectured that such diagonal linear forms should be 
 the only case of equality in \eqref{Minkowski conjecture, inequality in coordinates}.


\bigskip
\section{Quadratic forms and lattices} \label{quadratic forms}

The study of lattices is in a strong sense equivalent to the study of positive definite quadratic forms, 
over integer vector inputs, for the following 
simple reason.  Any positive definite quadratic form $f:\R^d \rightarrow \R$
 is defined by $f(x):= x^T A x$, where 
$A$ is a positive definite matrix, so the image of $\Z^d$ under $f$  is 
\[
\{ x^T A x \mid x \in \Z^d\}.
\]
On the other hand, any full-rank lattice in $\R^d$ is  by definition
$\L := M(\Z^d)$, for some real non-singular matrix $M$.  By definition, this implies that the square of the norm of any vector
in $\L$ has the following shape:  $\| v\|^2 = v^T v = x^T M^T M x$, for some $x \in \Z^d$.  We notice that
$M^T M$ in the last identity is positive definite.  

We may summarize this discussion as follows.  Given any lattice $\L:= M(\Z^d)$, we have
\begin{equation}\label{norms of lattice vectors, and pds forms}
\left\{ \| v\|^2   \bigm |    v \in \L \} = \{ x^T A x   \bigm |    x \in \Z^d   \right\},
\end{equation}
where $A:= M^T M$ is positive definite.

So the distribution of the (squared) norms of all vectors in a given lattice is equivalent to the image of $\Z^d$
under a positive 
definite quadratic form.   

Interestingly, despite this equivalence, for an arbitrary given lattice $\L$ it is not known in general whether the knowledge of
 the norms of all vectors in $\L$ uniquely determines the lattice $\L$.  In very small dimensions it is true,
 but for dimensions $\geq 4$ there are some counterexamples due to Alexander Schiemann (\cite{Schiemann1},  \cite{Schiemann2}).

The above equivalence between lattices in $\R^d$ and quadratic forms is straightforward but often useful, because it allows both algebraic and analytic methods to come to bear on important problems involving lattices.

 Gauss initiated the systematic study of finding  the minimum value of positive definite, binary quadratic forms $f(x, y) := a x^2 + 2b xy + cy^2$, over 
 all integer inputs $(x, y) \in \Z^2$.  Gauss' theory is also known as a reduction theory for positive definite binary quadratic forms, and is now a popular topic that can be found in many standard Number Theory books.  
 
 By the discussion of this short section, in particular 
 \eqref{norms of lattice vectors, and pds forms},
 it's clear that minimizing positive definite quadratic forms is essentially equivalent
  to finding a vector of smallest nonzero length in a lattice. 

 
 We close with a result of Mahler, regarding sequences of lattices.
So far we worked with one lattice at a time, but it turns out to be fruitful to work with infinite collections of lattices simultaneously. 
But what does it mean for a sequence of lattices to converge?  Luckily this notion is not difficult to define.  Suppose that we have a sequence of lattices $\L_n\subset \R^d$, and a fixed lattice $\L \subset \R^d$.  We say that
\[
\lim_{n\rightarrow \infty} \L_n = \L  
\]
if there exists sequences of bases $\beta_n$ of the  lattices $\L_n$ that converge to a basis $\beta$ of $\L$, in the sense that the $j$'th basis vector of $\beta_n$ converges to the $j$'th basis vector of $\beta$. 
In this direction, the following result, often called Mahler's compactness theorem,   is due to Kurt Mahler, \index{Mahler, Kurt} 
who was one of the main contributors to the development of the Geometry of Numbers.    
\begin{thm}[Mahler] \label{Mahler}
Fix $\rho >0, C>0$.  Then any infinite sequence of lattices 
$\L \subset \R^d$ such that
\[
\min \left\{  \|x\|  \bigm |    x  \in \L-\{0\}   \right\}    \geq     \rho,  \text{  and  }  \det \L \leq C,
\]
has an infinite convergent subsequence of lattices.
\hfill $\square$
\end{thm}
In other words, Mahler realized that among all lattices that have determinant equal to $1$, 
if a sequence of lattices diverges, then it must be true that the lengths of the shortest nonzero vectors
of these lattices tend to zero.


\bigskip
\section*{Notes}

 \begin{enumerate}[(a)]
\item  There is a well-known meme in Mathematics:  ``Can one hear the shape of a drum?", which is the title of Mark Kac's famous paper regarding the desire to discern the shape of a drum from its `frequencies'. 
An analogous question for lattices, studied by John Conway, is
``which properties of quadratic forms are determined by their representation numbers?''. 
 For further reading,  there is the lovely little book by Conway called ``The sensual quadratic form'', which draws connections between quadratic forms and 
  many different fields of Mathematics \cite{Conway.Book.SensualForm}.    

Of course, no library is complete without the important and biblical  ``Sphere Packings, Lattices and Groups", by John H. Conway and Neil Sloane \cite{ConwaySloan.book}.

\item  The idea of periodicity, as embodied by any lattice in $\R^d$,  also occurs on other manifolds, besides Euclidean space.    If we consider a closed geodesic on a manifold, then it's intuitively clear that 
as we flow along that geodesic, we have a periodic orbit  along that geodesic.
One important family of manifolds where this type of periodicity occurs naturally is the family of Hyperbolic manifolds.    Following the philosophy that `if we have periodicity, then we have Fourier-like series', 
it turns out that there is also an hyperbolic analogue of the Poisson summation formula, known as the Selberg trace formula \index{Selberg trace formula}, and this type of number theory has proved extremely fruitful. 

\item  It's natural to try to extend Minkowski's geometry of numbers to discretized volumes. 
In other words, we replace $\vol K$ in Minkowski's inequalities (and their extensions) by the lattice point enumerator 
 $\left | K \cap \L\right |$, which we will call its discretized volume. 
   In this direction, Betke, Henk, and Wills \cite{BetkeHenkWills}
 extended Minkowski's first theorem, obtaining good upper bounds on the lattice point enumerator, as follows.  For this note, $K\subset \R^d$ is always a centrally-symmetric, convex body.
 \begin{thm} \cite{BetkeHenkWills}
 With the notation above, we have:
 \[
 \left | K \cap \L\right | \leq  \left( \frac{2}{\lambda_1(K, \L) } + 1 \right)^d.
 \]
 \end{thm}
  \hfill $\square$
  
 It's easy to see that Minkowski's first theorem may be rewritten as
 \begin{equation}
 \frac{\vol K}{\det \L} \leq  \left(
 \frac{2}{\lambda_1(K, \L)}
 \right)^d,
 \end{equation}
as you showed in Exercise \ref{Rewriting Minkowski's first theorem}.  
Minkowski's second theorem may also be trivially rewritten as
 \begin{equation}
 \frac{\vol K}{\det \L} \leq 
 \prod_{k=1}^d
 \left(
 \frac{2}{\lambda_k(K, \L)}
 \right).
 \end{equation}
Therefore for a natural analogue of Minkowski's second theorem 
(Theorem \ref{successive minima for K}), the authors of \cite{BetkeHenkWills} gave the following conjecture for the discretized volume of $K$. 
 \begin{conjecture} \label{Conjecture:BetkeHenkWills}
 With the notation above, we have:
 \[
 \left | K \cap \L\right |  \leq  
  \prod_{k=1}^d 
\left \lfloor 
\frac{2}{\lambda_k(K, \L) } + 1 
 \right \rfloor.
 \]
  \hfill $\square$
 \end{conjecture}
In \cite{Malikiosis1}, Romanos-Diogenes Malikiosis proved 
Conjecture \ref{Conjecture:BetkeHenkWills} in dimension $3$.  As of this writing,
Conjecture \ref{Conjecture:BetkeHenkWills} is open in dimensions $d \geq 4$.

\item A strong bound for Hermite's constant in dimension $d$ was given by Blichfeldt \cite{Blichfeldt1}:  \index{Blichfeldt}
\[
\gamma_d \leq  \left( \frac{2}{\pi} \right)  \Gamma \left( 2 + \frac{ d}{2} \right)^{\frac{2}{d}}.
\]

\item  Related to the Hermite Normal Form is another extremely important reduction, called the Smith Normal Form (see the classic reference to integer matrices \cite{MorrisNewman}, by Morris Newman, as well as Exercise \ref{Ex:smith normal form}).

\item
Our proof of Theorem \ref{Prepping for main result on lattices} follows the approach taken in Barvinok's lecture notes \cite{BarvinokNotes}.   For Theorem \ref{sublattice index}, see \cite{FukshanskyBook} for a more standard approach, using the Hermite-normal form.

\item Some authors use the word \emph{gauge function} of $K$ instead of its distance function $|x|_K$. \index{gauge function} 

\item The family of diagonal matrices in Example \ref{a curve in the space of lattice} is very important in the study of homogeneous dynamics, because it  acts by multiplication on the left,
on the space of all lattices that have $\det \L = 1$.   This fascinating action is sometimes
called the ``modular flow'', and was studied intensively by Etienne Ghys.  
  A beautiful result in this direction is that the periodic orbits of the modular flow are in bijection with the conjugacy classes of hyperbolic elements in the modular group $GL_2(\Z)$, and furthermore that
these periodic orbits 
produce incredible knots in the complement of the trefoil knot.  


\item Almost all of the founders of the geometry of numbers, as well as many current mathematicians, have tried to prove Minkowski's conjecture, and the proofs in each dimension ($1 \leq d \leq 10$) have brought new ideas to the table.
 
 \item 
 It is clear that because lattices offer a very natural way to discretize $\R^d$, they continue to be of paramount importance to modern research.  In particular,  the theory of modular forms, with their Hecke operators that are defined using lattices and their fixed finite index sublattices, is crucial for modern number theory.   Euclidean lattices are also the bread-and-butter of crystallographers. 
 
\end{enumerate}

\newpage
\section*{Exercises}
\addcontentsline{toc}{section}{Exercises}
\markright{Exercises}

 \begin{quote}   
``No one ever wrote five pages of mathematics without a mistake."               

--  G. H. Hardy 
\end{quote}

\medskip
\begin{prob}  $\clubsuit$ 
\label{translating the Voronoi cell around}
Given a full rank lattice $\L\subset \R^d$, and any $m \in \L$, show that  
\[
{\rm Vor}_0(\L) + m =  {\rm Vor}_m(\L).
\]
\end{prob}

\medskip
\begin{prob}  $\clubsuit$ 
\label{facts about Vor cell}
Show that $\rm{Vor}_0(\L)$  is symmetric about the origin, convex, and compact.
\end{prob}

\blue{
\noindent Problems \ref{Voronoi cell and rho and mu} - 
\ref{Exercise:E_8} develop practice with the packing radius of a lattice, and with its covering radius.
}
\begin{prob}  $\clubsuit$ 
\label{Voronoi cell and rho and mu}
Given a full-rank lattice $\L \subset \R^d$, we have the following relations with its Voronoi cell $V(\L)$.
\begin{enumerate}[(a)]
\item Prove the packing radius $\rho(\L)$ equals the inradius of $V(\L)$.
\item The covering radius $\mu(\L)$ equals the circumradius of $V(\L)$.
\end{enumerate}
\end{prob}

\medskip
\begin{prob}
Suppose $\L_0 \subset \L\subset \R^d$ are two lattices of rank $d$.  Prove the following two inequalities:
\begin{equation}
\rho(\L) \leq \rho(\L_0)  \leq 
\rho(\L)  \left| \L / \L_0                \right|.
\end{equation}
\end{prob}

\medskip
\begin{prob}
For the lattice $D_n$, show that its covering radius is 
$\mu(D_n) =  \tfrac{1}{2}\sqrt d$ for $d \geq 4$, but $\mu(D_3) =1$.
\end{prob}

\medskip
\begin{prob} \label{Exercise:E_8}
Suppose $\L\subset \R^d$ is a full-rank lattice.    
Prove that for the $E_8$ lattice in $\R^8$, we have:
\[
\mu(E_8) = 1.
\]
\end{prob}


\medskip
\begin{prob} 
 $\clubsuit$
\label{proof that the gauge function is a norm}
\rm{
Given a centrally symmetric convex body $\K\subset \R^d$,  prove that the distance function
$|x|_\K$ that we defined in \eqref{def:distance function} is a norm on $\R^d$. 
}
\end{prob}

\medskip
\begin{prob} \label{elementary theta function exercise}
Suppose $\L\subset \R^d$ is a full-rank lattice.   Prove or disprove: 
\begin{equation}
\sum_{n \in \L} e^{-\pi  \| x-n\|^2} \geq
 e^{-\pi \| x\|^2} 
 \sum_{n\in \L} e^{-\pi  \| n\|^2},
\end{equation}
for each $x \in \R^d$.
\end{prob}

\medskip
\begin{prob} \label{successive minima of K are equal to those of the difference body of K}
Let $K\subset \R^d$ be a convex body, and fix a full-rank lattice $\L\subset \R^d$.   Here we show that the successive minima of $K$ are equal to those of its difference body.  That is, show that:
\[
\lambda_k(K, \L) = \lambda_k\left(\tfrac{1}{2} K-\tfrac{1}{2} K, \L\right),
\]
for all $1\leq k \leq d$.
\end{prob}

\medskip
\begin{prob}\label{Minkowski's little fact about the symmetrized body}
We say that a lattice $\L\subset \R^d$ is a {\bf packing lattice for $K$} if
\[
(K+n_1) \cap (K+n_2) = \phi, 
\]
for all $n_1, n_2 \in \L$, with $n_1 \not= n_2$.  In other words all of the distinct
translates of $K$, using vectors from $\L$, are disjoint.  Given a convex 
body $K\subset \R^d$, prove the following are equivalent.
\begin{enumerate}[(a)]
\item  $\L$ is a packing lattice for $K$ 
\item $\L$ is a packing lattice for the symmetrized body $\tfrac{1}{2} K-\tfrac{1}{2}K$.
\end{enumerate}
\end{prob}

\medskip
\begin{prob}\label{exercise:norms and convex bodies}
Prove Theorem \ref{norms and convex bodies}.
\end{prob}

\medskip
\begin{prob} \label{Extension of Erdos to dimension d}
\rm{
(hard)  Erd\H os' question, given in Exercise \ref{Erdos lattice partition problem}, possesses a natural extension to dimension $d$, as follows.
\begin{question}\label{tiling the lattice with translated sublattices}
Suppose that the integer lattice  $\Z^d$  is partitioned into a disjoint union of a  finite number of translates of  integer sublattices, say:
\[
\Z^d = \{ \L_1 +v_1 \} \cup \{  \L_2 + v_2  \}    \cup \dots \cup \{    \L_N +v_N\}.  
\]
 Is it true that there are at least two integer sublattices, say  $\L_j, \L_k$, 
that enjoy the property that $\L_k = \L_j + w$, for some integer vector $w$?  
\end{question}

Here we prove that in $\R^3$, Question \ref{tiling the lattice with translated sublattices} 
has a negative answer.  In particular, find a partition  of $\Z^3$ into $4$ integer sublattices, such that no two of them are integer translates of one another.  Using an easy extension to $d >3$, also show that the answer to the question above 
 is `no', if $d \geq 3$.

Notes.  Question \ref{tiling the lattice with translated sublattices}  remains unsolved in dimension $d=2$ \cite{FeldmanProppRobins}.
}
 \end{prob}

\medskip
\begin{prob} \label{proof of lemma for covering lattice of K}
$\clubsuit$
\rm{
Prove Lemma  \ref{lem: equivalence for covering lattice of K}.
}
\end{prob}

\medskip
\begin{prob}\label{Rewriting Minkowski's first theorem}
\rm{ 
Show that Minkowski's first theorem (namely Theorem 
\ref{Minkowski convex body Theorem for L}) 
has the following equivalent formulation.

Let $K\subset \R^d$ be a $d$-dimensional convex body, symmetric about the origin,
and let $\L$ be a (full rank) lattice in $\R^d$.
If $K$ is a convex body, symmetric about the origin, then
\begin{equation}
\lambda_1(K, \L)^d  \vol K  \leq  
 2^d \det \L.
\end{equation}
}
\end{prob}
Notes.   \ Although this reformulation is really just an observation, it does lead us to think more  carefully about the interior of $K$ (recalling that by definition a body is compact) versus $K$ itself.

\medskip
\begin{prob} \index{support function}
\rm{
Given any measurable sets $A, B \subset\R^d$, and all $x \in \R^d$, 
prove the following properties for their support functions:
\begin{enumerate}[(a)]
\item $h_{\lambda A}(x) = \alpha h_A(x)$, for any  $\lambda \geq 0$.
\item For any translation vector $v \in \R^d$, we have $h_{A+v}(x) = h_A(x) + \langle x, v \rangle$.
\item $h_{A+B}(x) = h_A(x) + h_B(x)$, where $A+B$ is the Minkowski sum of $A$ and $B$.
\index{Minkowski sum}
\end{enumerate}
}
\end{prob}

\medskip
\begin{prob} \label{equality condition for Minkowski's conjecture}
$\clubsuit$
\rm{
Show that the equality condition in Minkowski's conjecture \ref{Minkowski conjecture} 
is achieved by the diagonal linear forms $L_k(x_1, \dots, x_d):= 2 c_k x_k$.
}
\end{prob}


 \chapter{
 \blue{
 Sphere packings}
 } \label{Sphere packings}
 \index{sphere packings}

\begin{quote}
The problem of packing, as densely as possible, an unlimited number of equal nonoverlapping circles in a plane 
was solved millions of years ago by the bees, who found that the best arrangement consists of circles inscribed 
in the hexagons of the regular tessellation.  \ -- \  H. S. M.  Coxeter  \index{Coxeter}
\end{quote}

 \begin{quote}    
 There is geometry in the humming of the strings. There is music in the spacing of the spheres.
 \ -- \ Pythagoras   \index{Pythagoras}
 \end{quote}
                     

\begin{figure}[htb]
 \begin{center}
\includegraphics[totalheight=1.6in]{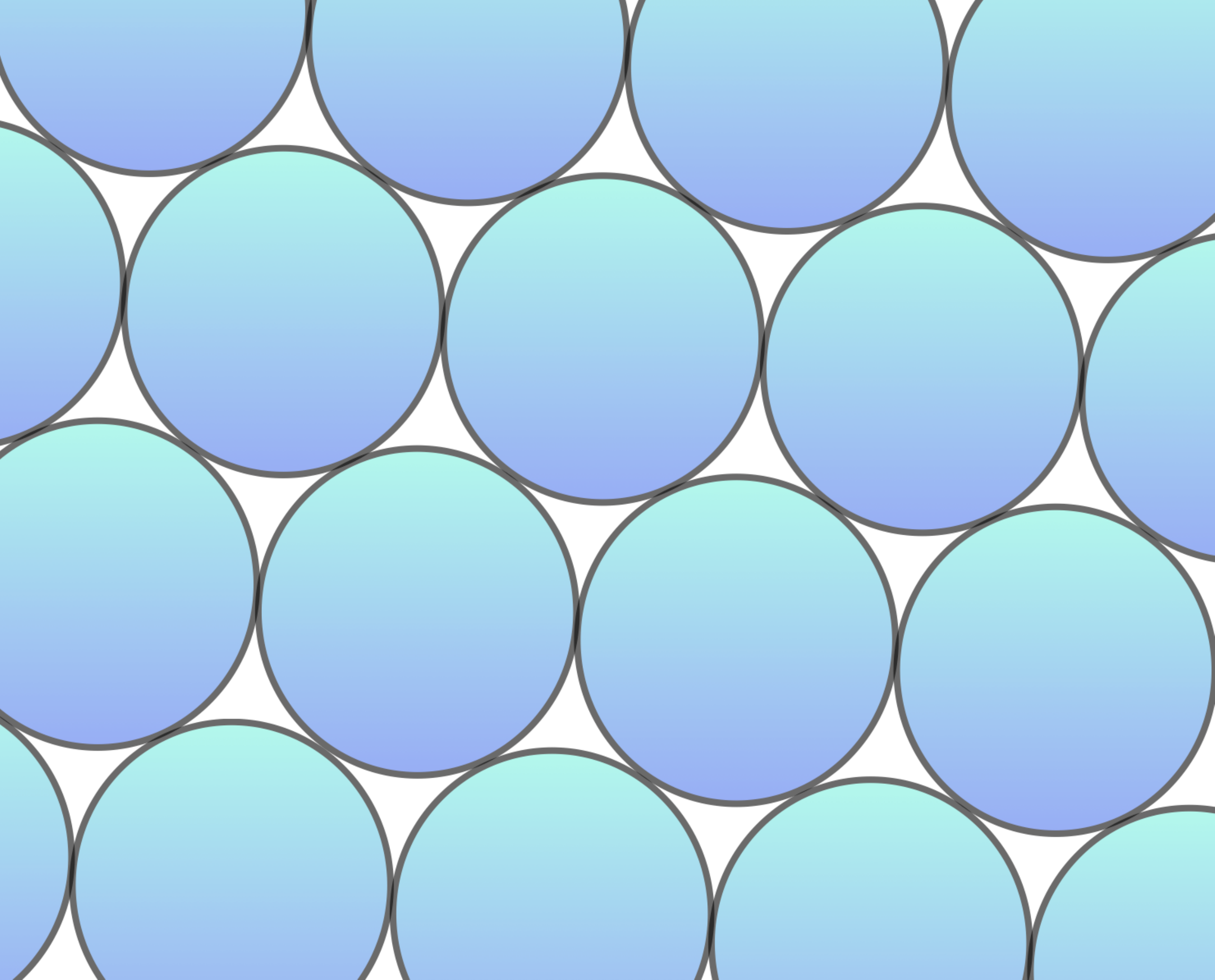}
\end{center}
\caption{A lattice sphere packing, using the hexagonal lattice, which gives the densest sphere packing in 2 dimensions.}  \label{periodic packing, Eisenstein}
\end{figure}

\section{Intuition}

\begin{wrapfigure}{R}{0.4\textwidth}
\centering
\includegraphics[width=0.24\textwidth]{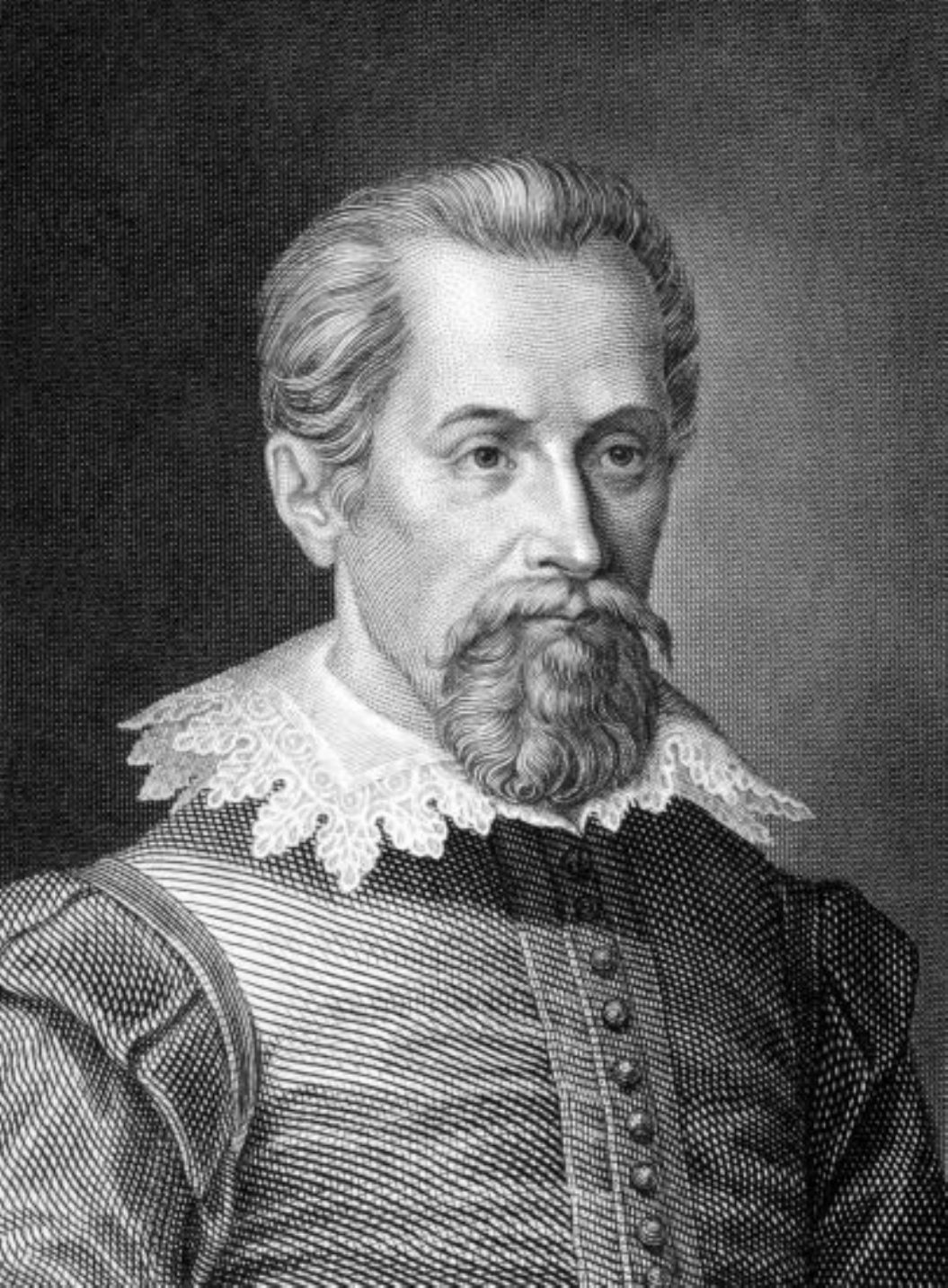}
\caption{Johannes Kepler}  \label{Kepler.pic}
\end{wrapfigure}

 The sphere packing problem traces its roots back to Kepler, and it asks for a packing 
 of solid spheres in Euclidean space that achieves the maximum possible density.  In all of the known cases,
 such optimal configurations - for the centers of the spheres - form a lattice.   It's natural, therefore, that Fourier analysis comes into the picture.   We prove here a result of Cohn and Elkies, from $2003$, which is a beautiful application of Poisson summation, and  gives certain upper bounds for the maximum densities of sphere packings
 in $\R^d$.

At this point it may be wise to define carefully all of the terms - what is a packing?  what is density?  
Who was Kepler?


\section{Lattice sphere packings}
A {\bf lattice sphere packing} in $\R^d$  
is a packing of balls, all having the same radius, with the property that the centers of the balls
 are located precisely at the  
points of some fixed lattice $\L \subset \R^d$, as in Figure \ref{Sphere Packing 1}.


A {\bf densest lattice sphere packing} is a lattice sphere packing 
with the additional property that no spheres with a larger radius will form a packing. 

\begin{example}
{\rm
In Figure \ref{Sphere Packing 1}, we've used the same lattice $\L\subset \R^2$ for two distinct sphere packings: on the left, we have a sphere packing that is not optimal, and on the right we have a sphere packing that is optimal.  But there are other lattices in $\R^2$ that give us a 
 sphere packing with a higher density. A natural question comes to mind: can we find a lattice that achieves the largest possible packing density?
\begin{figure}[htb]
 \begin{center}
\includegraphics[totalheight=2.3in]{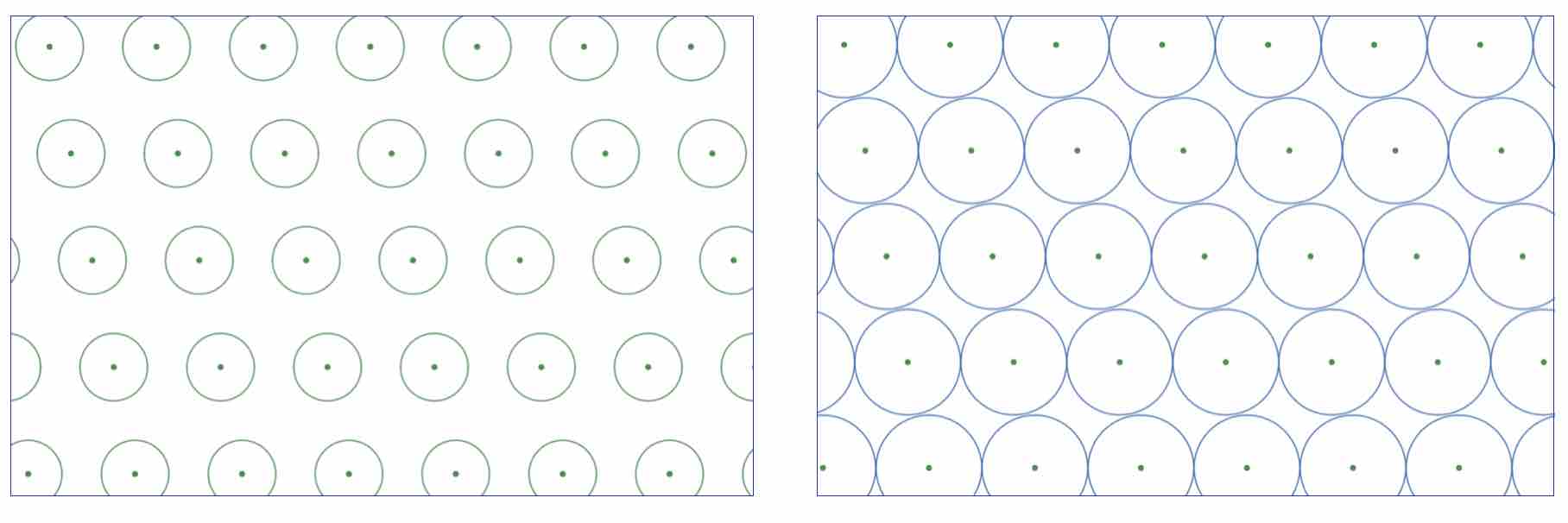}
\end{center}
\caption{Left: a lattice sphere packing, with a lattice $\L$, but with a small packing density. 
Right: the densest lattice sphere packing, for the same lattice $\L$.}  
\label{Sphere Packing 1}
\end{figure}
}
\hfill $\square$
\end{example}

Suppose we fix a lattice $\L$ and we have an optimal (densest) lattice packing with 
respect to $\L$. 
To quantify the packing density for a lattice sphere packing, we may think of ``how many balls
 of radius $r$ do we have per unit volume of the lattice $\L$''? 

We recall that $B_d$ is the ball of radius $1$ in $\R^d$, and that $B_d(r) \subset \R^d$ is the ball of  radius $r>0$. If we have a densest lattice sphere packing for a fixed lattice $\L$,   then 
we define its  {\bf lattice packing density}, relative to $\L$, by the expression:
\begin{equation}\label{def:lattice packing density}
\delta^*(B_d, \L):= \frac{\vol B_d(r)}{\det \L}.
\end{equation}
It's easy to see that to achieve the densest possible lattice sphere packing with this particular $\L$, we should use
 $r = \frac{1}{2}\lambda_1(\L)$, which is half-way from the origin to a nearest lattice point of $\L$.  
 In other words, we may rewrite the expression \eqref{def:lattice packing density} for the 
densest sphere packing \emph{of a fixed lattice} $\L$ as follows:
\begin{equation}\label{def:densest lattice packing density}
\delta^*(B_d, \L):=  \frac{1}{2^d}\lambda_1^d(\L)\frac{\vol B_d}{\det \L},
\end{equation}
where we used $\frac{\vol B_d(r)}{\det \L} = r^d \frac{\vol B_d}{\det \L}$.  

\begin{figure}[htb]
\begin{center}
\includegraphics[totalheight=1.73in]{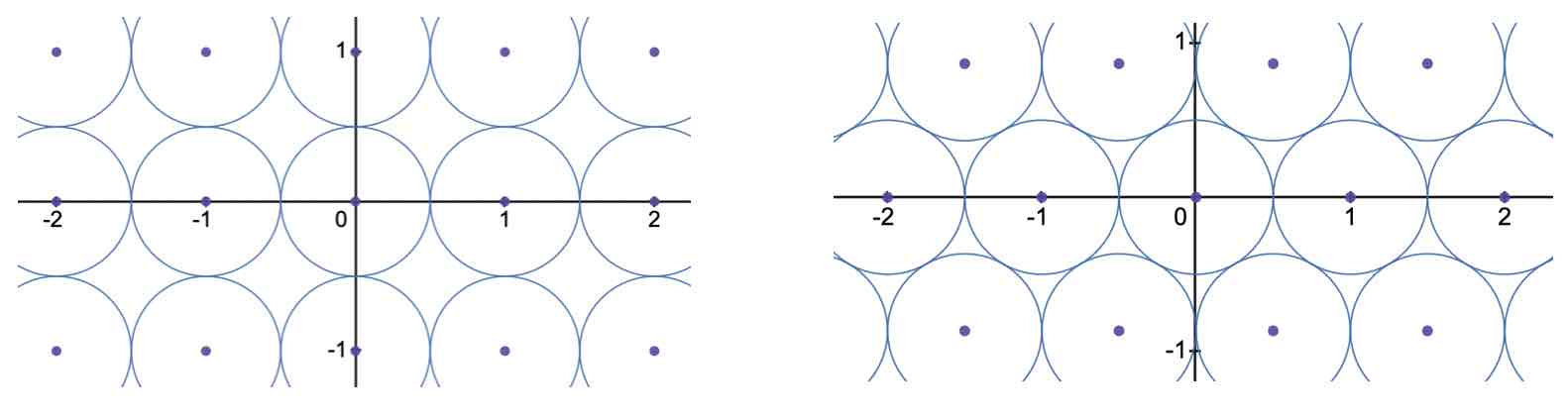}
\end{center}
      \caption{Left: the densest sphere packing for the lattice $\L_1:=\Z^2$, with a packing density of 
$\delta^*(B_2, \L_1)=\frac{\pi}{4} \approx .7854$, so that approximately $78.5 \%$ of the plane is covered by this configuration of balls.  Right: the densest sphere packing for the hexagonal lattice $\L_2$, with a packing density of 
$\delta^*(B_2, \L_2)=\frac{\pi\sqrt 3}{6} \approx .9068$.} 
         \label{densest integer lattice packing}
\end{figure}

 \begin{example}
 \rm{
 Consider the integer lattice $\L:= \Z^2$. It is clear that we can place non-overlapping 
 spheres of radius 
 $r = \frac{1}{2}$ at each integer point, as in Figure \ref{densest integer lattice packing} on the left.  It is also clear that any larger radius for our spheres will incur overlapping spheres. 
 So this particular packing gives us a sphere packing density of 
 \[
\frac{  \vol B_2(r)  }{\det \L}:=  \frac{  \frac{\pi}{4}  }{   \det \Z^2 } =  \frac{\pi}{4}
 \approx .7854.
 \]
 }
 \hfill $\square$
 \end{example}
 
\begin{example}
\label{densest lattice sphere packing in R^2 - Lagrange}
 \rm{
Now consider the hexagonal lattice 
$
\L:=  
     \begin{pmatrix}  1 &  \frac{1}{2}   \\  
                        0  & \frac{\sqrt 3}{2}   
        \end{pmatrix}\Z^2, 
$
as in Figure \ref{densest integer lattice packing} on the right.
It's clear that we can still place non-overlapping spheres of radius 
 $r = \frac{1}{2}$ at each of these lattice points.
This particular arrangement of spheres gives us a packing density of 
 \[
\frac{  \vol B_2(r)  }{\det \L}
=  \frac{  \left( \frac{\pi}{4} \right) }{ \left( \frac{\sqrt 3}{2} \right) } 
=\frac{ \pi \sqrt 3}{6}   \approx  .9068.
 \]
 As it turns out, this is the densest lattice sphere packing, and Lagrange \cite{Lagrange} was the first to prove it, in 1773.  
 }
 \hfill $\square$
 \end{example}

More generally, suppose we are given a convex, compact, centrally-symmetric set $K\subset \R^d$. We may similarly  define a densest packing density for $K$, relative to any fixed lattice $\L\subset \R^d$:
\begin{equation}\label{def:densest lattice packing density for K}
\delta^*(K, \L):=  \frac{1}{2^d}\lambda_1^d(K, \L)\frac{\vol K}{\det \L},
\end{equation}
where we used definition \eqref{def:general successive minima} 
for the first successive minimima of $K$ relative to $\L$, namely
$\lambda_1(K, \L)$. 

But which lattice achieves the densest lattice sphere packing? This is clearly an extremely difficult problem in general, because it is a discrete optimization problem over 
an infinite-dimensional space of lattices.  Nevertheless, we may define 
the densest lattice packing of $K\subset \R^d$, as we vary over all full-rank lattices 
$\L\subset \R^d$:
\begin{equation}\label{def:densest lattice packing density for K, over all lattices}
\delta^*(K):=  \sup_{\L \subset \R^d} \delta^*(K, \L),
\end{equation}
where the supremum is taken over all full-rank lattices in $\R^d$.

\begin{question}[The lattice sphere-packing problem]
\label{The sphere packing problem}
For each dimension $d\geq 2$, find the value of $\delta^*(B_d)$, and find a lattice
$\L\subset \R^d$ that achieves it.
\end{question}

\begin{example}\label{densest lattice sphere packing in R^3 - Gauss}\index{Gauss}
 \rm{
In $\R^3$, we consider the ``face-centered-cubic'' lattice, defined by: \index{face-centered-cubic lattice}
\[
\L:=  
     \begin{pmatrix} \ \ 1 & \ \ 1 &\ \  0    \\  
                           \ \    1 &    -1 &\ \  1    \\  
                           \ \    0 & \ \ 0 &    -1    \\  
        \end{pmatrix}\Z^3.
\]
This lattice sphere packing gives us a packing density of 
 \begin{equation}\label{Gauss lattice sphere density in R^3}
\frac{  \vol B_3(r)  }{\det \L}
=\frac{\pi}{ \sqrt{18}} \approx  .7405.
 \end{equation}
 As it turns out, this is the densest lattice sphere packing, a statement that was first proved in 1831, by Gauss.
In the notation of Question \ref{The sphere packing problem}, 
Gauss proved that  $\delta^*(B_3)= \frac{\pi}{ \sqrt{18}}$.
 }
 \hfill $\square$
 \end{example}
 
\begin{example}\label{densest lattice sphere packing in R^4}
 \rm{
In 1872,  Korkin and Zolotarev  \cite{KorkinZolotarev} discovered the answer to the lattice sphere packing problem in $\R^4$. Namely, the ``checkerboard'' lattice $D_4$  (defined in Example \ref{D_n lattices}) gives 
\[
\delta^*(B_4) = \frac{\pi^2}{16}.
\]
In 1877, Korkin and Zolotarev also discovered that $D_5$ gives the 
densest lattice packing in $\R^5$ \cite{KorkinZolotarev5}. 
 \hfill $\square$
 }
 \end{example}

Putting some of the definitions above together, we arrive at the following elementary but useful equivalence.
\begin{lem}\label{lem:equivalence of lattice packing density}
Suppose we are given a convex, compact, centrally-symmetric set $K\subset \R^d$, 
with $\vol K \leq  C$.  
The following are equivalent:
\begin{enumerate}[(a)]
\item \label{part a of the equivalence of lattice packing density}
There exists a 
lattice $\L\subset \R^d$ with $\det \L = 1$,  such that $K$ contains no nonzero 
points of $\L$.
\item \label{part b of the equivalence of lattice packing density}
There exists a lattice $\L\subset \R^d$ such that 
$ \delta^*(K, \L) \geq \frac{1}{2^d}C$. \hfill $\square$
\end{enumerate}
\end{lem}

We'll let the reader enjoy writing out the straightforward proof of 
Lemma \ref{lem:equivalence of lattice packing density} 
(Exercise \ref{proof of Lemma, equivalence of packing density}).  When $K$ is any centrally-symmetric convex body,  
and $C := 2 \zeta(d)$, Hlawka proved \cite{Hlawka} that 
part \ref{part a of the equivalence of lattice packing density}
Lemma \ref{lem:equivalence of lattice packing density} is in fact true (see \eqref{The Minkowski-Hlawka bound} below).

\bigskip
\section{More general sphere packings}

In order to allow more general packings of spheres, we will relax the restriction of using just one lattice, as follows.  
A {\bf densest periodic sphere packing} is a packing of spheres of radius $r$, with a lattice $\L$, but also with a finite collection of its translates, 
say $\L + v_1, \dots, \L + v_N$, such that the differences $v_i-v_j \notin L$, and such that no larger radius will form a packing.
   Figure \ref{periodic packing} gives such an example, by using $3$ lattices, rather than just one.

\begin{figure}[htb]
 \begin{center}
\includegraphics[totalheight=2in]{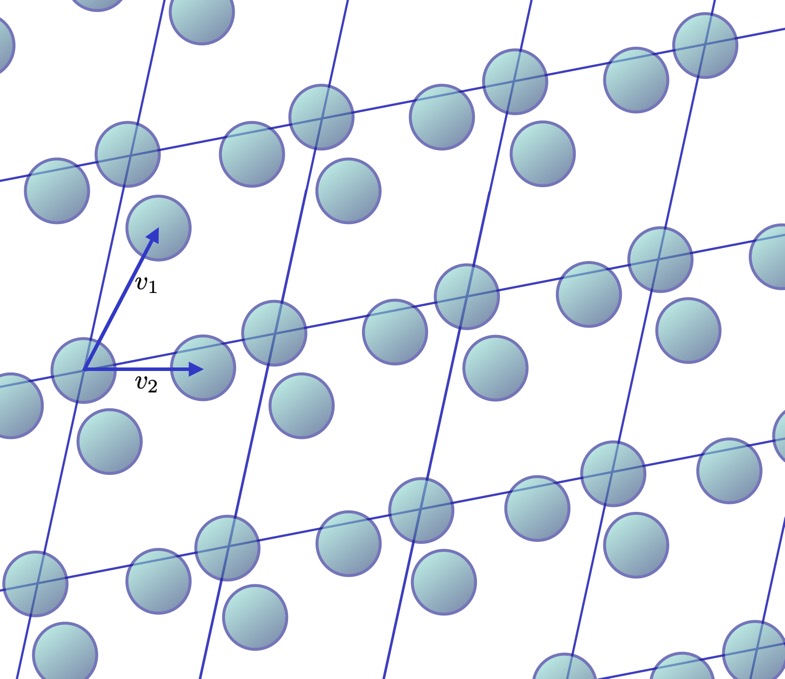}
\end{center}
\caption{A periodic sphere packing with two translates of the same lattice.  This packing is not a lattice packing.}  \label{periodic packing}
\end{figure}

Given a period sphere packing with a lattice $\L$ and a set of translates $v_1, \dots, v_N$,
 we define its {\bf periodic sphere packing density} by 
  \begin{equation} \label{periodic packing density}
 \delta(B_d, \L) := \frac{ N \vol B_d(r)   }{\det \L},
 \end{equation}
 corresponding to placing a sphere of radius $r$ at each point of $\L$, and also at each point
 of its translates $\L + v_1, \dots, \L + v_N$.  It's not hard to prove that the latter definition 
 \ref{periodic packing density} matches our intuition that  any fixed fundamental parallelepiped of $\L$
 intersects this configuration of spheres in a set whose measure is exactly 
 $N  \vol B_d(r) $ (Exercise  \ref{volume of periodic packing}).

 The densest periodic sphere packing is defined by varying over 
 all lattices, as follows. 
\begin{equation}\label{def:densest periodic sphere packing density, over all lattices}
\delta(B_d):=  \sup_{\L \subset \R^d} \delta(B_d, \L),
\end{equation}
where the supremum is taken over all full-rank lattices in $\R^d$.

\begin{question}[The general sphere-packing problem]
\label{The sphere packing problem}
For each dimension $d\geq 2$, find the value of $\delta(B_d)$, and find a lattice
$\L\subset \R^d$ that achieves it.
\end{question}
Many other questions naturally arise. 
\begin{question} \label{lattice or a few lattices?}
Fixing the dimension $d$, is $\delta(B_d)= \delta^*(B_d)$?
\end{question}
 In other words, may we always use just one lattice to achieve a densest sphere packing?
  Or are there dimensions $d$ for which we need to use some translates of a lattice as well?

\begin{question}  \label{unique lattice?}
 If the answer to Question \ref{lattice or a few lattices?} is affirmative, then is such an optimal lattice unique in each dimension?
 \end{question}

It turns out that periodic sphere packings come arbitrarily close to arbitrary sphere packings, as shown in Appendix A of \cite{Cohn-Elkies} (see also Greg Kuperberg's paper \cite{Kuperberg}). 
So there is no loss of generality in merely considering periodic sphere packings for the general problem of sphere packings, as done by Cohn and Elkies \cite{Cohn-Elkies}.

With these news ideas in mind, we can revisit our densest lattice sphere packings in dimensions $2$ and $3$, from Examples 
\ref{densest lattice sphere packing in R^2 - Lagrange} 
and \ref{densest lattice sphere packing in R^3 - Gauss}.  In 1910, Axel Thue \cite{Thue} was the first to prove that in $\R^2$ we do indeed have
 \[
 \delta(B_2)= \delta^*(B_2), 
 \]
 finally settling completely the question of the general sphere packing problem for $d=2$ 
 (for any configuration of spheres of the same radius), a question that eluded even Gauss. 
 In 1950, Laszlo F. Toth \cite{Toth1} extended the work of Thue, by allowing incongruent circles of different radii, among other extensions.

For $\R^3$, the famous astronomer Johannes Kepler conjectured in $1611$ that the densest possible sphere packing was given by a lattice, namely the face-centered cubic lattice of Example \ref{densest lattice sphere packing in R^3 - Gauss}. 
In 2005, Thomas C. Hales published a proof  \cite{Hales} that indeed
\[
\delta(B_3) = \frac{\pi}{ \sqrt{18}} \approx  .7405,
\]
confirming that indeed Kepler was correct, and that $\delta(B_3) = \delta^*(B_3)$.  Hales' proof, part of which was done with Sam Ferguson,
 was a huge breakthrough,  even though it was a computer-aided proof.  
 
In $\R^4$,  is it true that the densest periodic sphere packing given by the $D_4$  lattice? It seems reasonable to conjecture that it is, but this is still open as well.
The only dimensions $d$ for which we know the answers to Question \ref{lattice or a few lattices?} and 
Question \ref{unique lattice?} are $d=1, 2, 3, 8, 24$, and in these known cases the answer is affirmative.  Nevertheless, it seems reasonable to think that in higher dimensions the answer will become negative - there is more 'freedom to move around'.  That is, it is widely believed that in higher dimensions, we might discover that often $\delta(B_d) \not= \delta^*(B_d)$.
Even in dimension $10$, the densest known sphere packing is not a lattice packing.
The sphere packing problem continues to intrigue, and it is a very important 
problem in Geometry, Number theory, Coding theory, and information theory.

\bigskip
\section{Upper bounds for sphere packings, via Poisson summation}

Here we give an exposition of the ground-breaking result of Henry Cohn and Noam Elkies on the sphere packing problem \cite{Cohn-Elkies}.
This result sets up the machinery for finding certain {\bf magical functions} $f$, as defined
 in Theorem \ref{Cohn-Elkies}  below, 
 that allow us to give precise upper bounds on 
$\Delta_{periodic}{\L}$.   The main tool is Poisson summation again, for arbitrary lattices.  We recall that we defined a function $f$ to be \emph{nice} if $f$ satisfies the Poisson summation formula 
\[
\sum_{n \in \L} f(n+v) = \frac{1}{\det \L}  \sum_{\xi \in \L^*}  
                           \hat f(\xi)  e^{2\pi i \langle v, \xi \rangle}, 
\]
pointwise for all  $v\in \R^d$.

\begin{thm}[Cohn-Elkies]  \label{Cohn-Elkies}
Let $f:\R^d \rightarrow \R$ be a nice function, not identically zero, which enjoys the following 
three conditions:
\begin{enumerate}
\item $f(x) \leq 0$, for all $\|x\| \geq r$.      \label{condition 1}
\item  $\hat f(\xi) \geq 0$, for all $\xi \in \R^d$.  \label{condition 2}
\item  $f(0) >0$, and $\hat f(0) > 0$.     \label{condition 3}
\end{enumerate}
Then the periodic sphere packing density  has the upper bound:
\[
\delta(B_d) \leq \frac{       f(0)     }{      \hat f(0)    } \vol  B_d(r) .
\]
\end{thm}
\begin{proof}
Suppose we have a periodic packing with spheres of radius $r$, a lattice $\L$, and translation vectors
$v_1, \dots, v_N$, so that by definition the packing density is 
$\delta(B_d, \L) := \frac{ N \vol  B_d(r)  }{\det \L}$.   \\
\noindent
By Poisson summation,  \index{Poisson summation formula}
we have
\begin{equation}
\sum_{n \in \L} f(n+v) = \frac{1}{\det \L}  \sum_{\xi \in \L^*}  
                           \hat f(\xi)  e^{2\pi i \langle v, \xi \rangle},
\end{equation}
converging absolutely for all $v \in \R^d$.  Now we form the following finite sum and rearrange 
the right-hand-side of Poisson summation:
\begin{align}\label{fancy Poisson}
 \sum_{1\leq i \leq j \leq N}  
 \sum_{n \in \L} f(n+v_i - v_j) &= \frac{1}{\det \L}  \sum_{\xi \in \L^*}  
                           \hat f(\xi)     \sum_{1\leq i \leq j \leq N} e^{2\pi i \langle v_i - v_j, \xi \rangle}  \\  
&=   \frac{1}{\det \L}  \sum_{\xi \in \L^*}  
                           \hat f(\xi)     \Big|  \sum_{1\leq k \leq N} e^{2\pi i \langle v_k,  \xi \rangle} \Big|^2.    
                           \label{RHS of Poisson}      
\end{align}
Every summand on the right-hand-side of \eqref{RHS of Poisson}   is nonnegative, because by the second
assumption of the Theorem, we have
$\hat f (\xi) \geq 0$, so that the whole series can be bounded from below by its constant term, which for $\xi = 0$ 
gives us the bound  $  \frac{ \hat f(0) N^2}{\det \L}$.

On the other hand, let's ask what the positive contributions are, from the left-hand-side of \eqref{fancy Poisson}.
Considering the vectors  $n+v_i - v_j$ on the left-hand-side of \eqref{fancy Poisson},
suppose we have $\| n+v_i - v_j \| \geq r$. Then the first hypothesis of the Theorem guarantees 
that $f(n+v_i - v_j ) \leq 0 $.  So we may restrict attention to those vectors that satisfy $\| n+v_i - v_j \| < r$.  Here
the vector $n+v_i - v_j$ is contained in the sphere of radius $r$, centered at the origin, but this means (by the packing assumption) that it must be the zero vector:  $n  + v_i - v_j = 0$.   By assumption, the difference between any two translations $v_i-v_j$ is never a nonzero element of $\L$, so we have $i=j$, and  now $v_i = v_j \implies n=0$.  We conclude that the only positive contribution from the left-hand-side of \eqref{fancy Poisson}  is the $n=0$ term, and so the
 left-hand-side of \eqref{fancy Poisson} has an upper bound of $N f(0) >0$.  
 
 Altogether, Poisson summation gave us the bounds:
\[
 N f(0) 
 \geq
 |  \sum_{1\leq i \leq j \leq N}    \sum_{n \in \L} f(n+v_i - v_j) |
=
 \frac{1}{\det \L}  \sum_{\xi \in \L^*}  
                \hat f(\xi)     \Big|  \sum_{1\leq k \leq N} e^{2\pi i \langle v_k,  \xi \rangle} \Big|^2
\geq  
\frac{ \hat f(0) N^2}{\det \L}.
\]
Simplifying, we have
\[
\frac{ f(0) }{ \hat f(0) } \geq \frac{N}{\det \L} := 
\frac{ \delta(B_d, \L)}{  \vol   B_d(r)    }.
\]
Since the upper bound  $\frac{ f(0) }{ \hat f(0) } \vol   B_d(r)$ does not depend on the lattice $\L$, we get the desired result.
\end{proof}

\bigskip
\begin{example} [The trivial bound]
\rm{
Let $\L$ be a full-rank lattice in $\R^d$, whose shortest nonzero vector has length $r>0$. 
  We  define the function 
\[
f(x):= 1_{K}(x) * 1_{K}(x), 
\]
where $K$ is the ball
of radius $r$, centered at the origin.   
 We claim that $f$ satisfies all of the conditions of 
Theorem \ref{Cohn-Elkies}. Indeed, by the convolution Theorem,
\[
\hat f(\xi)  = \widehat{\left(1_{K} * 1_{K}\right)}(\xi)
= \Big(  \hat 1_{K}(\xi) \Big)^2 \geq 0, 
\]
for all $\xi \in \R^d$, verifying condition \ref{condition 2}.    Condition \ref{condition 1} is also easy to verify, because the support of $f$ is equal to the Minkowski sum (by Exercise  \ref{support of convolution}) 
\index{Minkowski sum}
 $K + K = 2K$, a sphere of radius $2r$. 
It follows that $f$ is identically zero outside a sphere of radius $2r$.   For condition \ref{condition 3}, by the definition of convolution we have $f(0) = 
\int_{\R^d} 1_{K}(0-x)  1_{K}(x)dx =  \int_{\R^d}  1_{K}(x)dx =
\vol K >0$.  Finally, $\hat f(0)  = \Big(  \hat 1_{K}(0) \Big)^2 =
\vol^2( K) > 0$.

By the Cohn-Elkies Theorem \ref{Cohn-Elkies}, we know that the packing density of such a lattice is therefore bounded above by
\[
 \frac{ f(0) }{  \hat f(0)} \vol  B_d(r) 
= \frac{\vol K }{  \vol^2( K)} \vol K = 1,
\]
the trivial bound.   So we don't get anything interesting, but all this tells us is that our particular choice of function $f$ above was a poor choice, as far as density bounds are concerned.   We need to be more clever in picking our magical $f$.
}
\hfill $\square$
\end{example}

Although it is far from trivial to find magical functions $f$ that satisfy the hypothesis of the 
Cohn-Elkies Theorem, and  simultaneously give a strong upper bound, 
there has been huge success recently in finding exactly such functions -  in dimensions $8$ and $24$.   These recent magical functions gave the densest sphere packings in these dimensions, knocking off the whole sphere packing problem in dimensions $8$ and $24$. 

Another observation that is useful is that if we have a magical function $f$ that enjoys all three 
hypotheses of the Cohn-Elkies
Theorem \ref{Cohn-Elkies}, then  $f\circ \sigma$ also satisfies the same hypotheses, 
for any $\sigma \in SO_d(\R)$ (Exercise \ref{invariance of magical functions under orthogonal transformations}).
We may therefore take certain radial functions as candidates for magical functions.

As of this writing, the provably densest packings are known only in dimensions
$1, 2, 3, 8$,  and $24$.  Each dimension seems to require slightly different methods, and sometimes wildly different methods, 
such as $\R^3$.     
Somewhat surprisingly, the sphere packing problem is  still open in all other dimensions.

\bigskip
\section{Lower bounds for sphere packings}

Here we discuss lower bounds for the optimal sphere packing problem.  As always, things began with Minkowski.  


\begin{thm}[The Minkowski-Hlawka theorem, 1943]
\label{thm:Minkowski-Hlawka}
 Let $K\subset \R^d$ be a $d$-dimensional centrally symmetric convex body.
  Then we have:
\begin{equation}\label{The Minkowski-Hlawka bound}
\delta^*(K) \geq \frac{\zeta(d)}{2^{d-1}},
\end{equation}
where $\zeta(s)$ is the Riemann zeta function.
\hfill $\square$
\end{thm}
 Minkowski proved \eqref{The Minkowski-Hlawka bound} in the special case that $K:= B$, the unit ball in $\R^d$, and he conjectured that the same lower holds for all centrally symmetric convex bodies.  This particular conjecture of Minkowski was finally proved by Hlawka \cite{Hlawka}. 
We recall that a unimodular lattice $\L$ simply enjoys $\det \L = 1$.  
To prove Theorem \ref{thm:Minkowski-Hlawka}, Hlawka first proved 
the following interesting fact (in a straightforward manner) concerning the existence of a certain unimodular lattice \cite{Hlawka}.
\begin{lem}[Hlawka, 1943]
\label{lemma:Hlawka}
Let $f(x)\in L^1(\R^d)$ be a bounded, real-valued, compactly supported function.  For each $\varepsilon >0$, there exists a unimodular lattice $\L\subset \R^d$ such that 
\begin{equation}
\sum_{n\in \L\setminus\{0\}} f(n) < \int_{\R^d} f(x) dx + \varepsilon.
\end{equation}
\hfill $\square$
\end{lem}

In 1945, Carl Ludwig Siegel extended Hlawka's Lemma \ref{lemma:Hlawka} with the following fundamental fact \cite{Siegel1}.  To state this result, we will need to assume a bit more background of the reader, just for the remainder of this section. 
By way of introductions, the group  $G:= \SL_d(\R) / \SL_d(\Z)$ 
(which is a very important Lie group in Number theory and physics),
 may be thought of as the space of all unimodular lattices in $\R^d$. 
\begin{thm}[C. L. Siegel, 1945]
\label{thm:Siegel mean values}
Let $f(x)\in L^1(\R^d)$ be a bounded, compactly supported function.  Then:
\begin{equation} \label{Siegel's equation for mean values}
\int_{G} \left( \sum_{n \in \Z^d \setminus \{0\} } f(An) \right) d\mu
=
\int_{\R^d} f(x) dx,
\end{equation} 
where $\mu$ is the unique normalized Haar measure on $G:= \SL_d(\R) / \SL_d(\Z)$, 
and the matrix $A$ varies over $G$.
\hfill $\square$
\end{thm}
So Hlawka's Lemma \ref{lemma:Hlawka} follows immediately from 
Siegel's Theorem \ref{thm:Siegel mean values}, even with $\varepsilon = 0$. 
As another example, let's consider the special
case of $f(x):= 1_S(x)$, where $S$ is any bounded, measurable subset of $\R^d$.   Here \eqref{Siegel's equation for mean values} gives us the intuitively compelling conclusion that 
as $\L$ varies over all unimodular lattices in $\R^d$, the number of nonzero lattice points in 
$\L \cap S$ is ``on average'' equal to the volume of $S$ (where this ``average'' is really the integral over the space of all unimodular lattices).
The reader may also consult Zong \cite{Zong.SpherePackingsBook}, 
for a proof of Theorem \ref{thm:Siegel mean values}.

In 1992, Keith Ball \cite{KeithBall.1} improved upon the Minkowski-Hlawka theorem, 
and gave the following lower bounds, which included a new linear term.
\begin{thm}[Ball, 1992]
\[
\delta^*(B)  \geq  \frac{(d-1)}{2^{d-1}} \zeta(d), 
\]
where $B$ is the $d$-dimensional unit ball.  
\hfill $\square$
\end{thm}
For large dimensions, $\zeta(d)$ is very close to $1$, so some authors omit the factor of $\zeta(d)$. 

In 2013, Akshay Venkatesh \cite{Venkatesh}  has given an improvement over the known lower bounds, by using a variation of Siegel's Theorem \ref{thm:Siegel mean values} above.
\begin{thm}[Venkatesh, 2013]
There exist infinitely many dimensions $d$ for which
\[
\delta^*(B) > \frac{ \log \log d}{2^{d+1}}d.
\]
In addition, for all sufficiently large dimensions, we have $\delta^*(B) > \frac{65, 963 }{2^d}d$.
\hfill $\square$
\end{thm} 

\bigskip \bigskip

\section*{Notes} \label{Notes.chapter.SpherePackings}
\begin{enumerate}[(a)]
\item  Each dimension $d$ appears to have a separate theory for sphere packings.  This intuition is sometimes tricky to conceptualize, but there are facts that help us do so.  For example, it is a fact that the Gram matrix (see  \ref{Gram matrix positive semidefinite}) of a lattice $\L \subset \R^d$ consists entirely of integers, with even diagonal elements  $\iff   \ d$ is divisible by $8$.  
  For this reason, it turns out that the theta series of a  lattice possesses certain functional equations (making it a modular form)  if and only if $8 \mid d$, which in turn allows us to build some very nice related `magical' functions $f$ that are sought-after in Theorem \ref{Cohn-Elkies}, at least for $d=8$ and $d=24$ so far.   
  
  In dimension $2$, it is an open problem to find such magical functions, even though there is
 an independent proof that the hexagonal lattice is the optimal sphere packing lattice. 

\item Johannes Kepler  (1571 --1630)  \index{Kepler, Johannes}
 was a German astronomer and mathematician.  
Kepler's laws of planetary motion motivated Sir Isaac Newton to develop further the theory of gravitational attraction and planetary motion.  Kepler conjectured that the densest packing of sphere is given by the 
``face-centered cubic''  packing.   It was Gauss (1831) \index{Gauss}
who first proved that, if we assume the packing to be a lattice packing, then Kepler's conjecture is true. 
In $1998$ Thomas Hales (using an approach initiated by L. Fejes T\'oth (1953)), gave an unconditional proof of the Kepler conjecture.

\item  It is also possible, of course, to pack other convex bodies.   One  such variation is to pack regular tetrahedra in $\R^3$.   
The interesting article by Jeffrey Lagarias and Chuanming Zong \cite{LagariasZong} gives a nice account of this story.   For a classical introduction to sphere packings, the reader may also consult the book by Chuanming Zong \cite{Zong.SpherePackingsBook}.   There is also a nice survey paper on many other aspects of packings, coverings, and tilings by Zong \cite{Zong.SpherePackingsSurveyPaper}.

\item We mention some of the recent spectacular applications of the 
Cohn-Elkies Theorem.   In 2016, 
Maryna Viazovska was able to find these magical functions for $\R^8$, thereby proving that the $E_8$ lattice gives the densest sphere packing in dimension $8$.  Shortly afterwards, professor Viazovska, joined with the team effort of 
Henry Cohn, Abhinav Kumar, Stephen D. Miller,  and Danylo Radchenko, managed to also find magical functions in
$\R^{24}$ \cite{Cohn.etal}.
Here is a synopsis of some of  their results.

\begin{thm}[\cite{Cohn.etal}]
The lattice $E_8$ is the densest periodic packing in $\R^8$.   The Leech lattice is the densest
periodic packing in $\R^{24}$.   In addition, these lattices are unique, 
in the sense that there do not exist any
other periodic packings that achieve the same density. 
\end{thm}

\end{enumerate}


\newpage
\section*{Exercises}
\addcontentsline{toc}{section}{Exercises}
\markright{Exercises}

 \begin{quote}    
 ``It is better to do the right problem the wrong way, than the wrong problem the right way.''
 
 --  Richard Hamming    \index{Hamming, Richard}
 \end{quote}

\medskip
\begin{prob} \label{volume of periodic packing}
\rm{
Given a periodic lattice packing, by $N$ translates of a lattice $\L \subset \R^d$, show that any fixed fundamental 
parallelepiped of $\L$
 intersects the union of all the spheres in a set of measure $N  \vol B_d(r) $, where $r:= \frac{1}{2}\lambda_1(\L)$.
 Thus, we may compute the density of a periodic sphere packing by just considering the
 portions of the spheres that lie in one fundamental parallelepiped.
 }
\end{prob}

 \medskip
\begin{prob}  
\rm{
Here we show that the integer lattice $\Z^d$ is a very poor choice for sphere packing.
\begin{enumerate}[(a)]
\item  Compute the packing density of the integer lattice $\Z^d$.
\item Compute the packing density of the lattices $D_3$ and $D_4$.
\item Compute the packing density of the lattices $D_n$, for $n\geq 5$.
\end{enumerate}
}
\end{prob}

\medskip
\begin{prob}  \label{radial function transform}
\rm{
If $f\in L^1(\R^d)$ is a radial function, then prove that its Fourier transform $\hat f$
is also a radial function.
}
\end{prob}

 \medskip
\begin{prob} 
\rm{
Suppose we pack equilateral triangles in the plane, by using only translations of a fixed equilateral triangle $\Delta$.
To make the problem easier, we'll restrict attention to lattice packings here.
What is the maximum lattice packing density of $\Delta$?  Do you think it may be the worst possible density among 
 lattice packings of any convex body in $\R^2$?
}
\end{prob}

\medskip
\begin{prob}  \label{invariance of magical functions under orthogonal transformations}
\rm{
 Show that if we have a magical function $f$ that enjoys all $3$ hypotheses of Theorem~\ref{Cohn-Elkies},
 then   $f\circ \sigma$ also satisfies the same hypotheses, for any orthogonal transformation $\sigma \in SO_d(\R)$.
 }
\end{prob}
 
 \medskip
\begin{prob} 
\rm{
We define a rigid motion of a compact set $K$ to be any orthogonal transformation of $K$, composed with any translation of $K$.  
\begin{enumerate}[(a)]
\item When $d=1$, find a continuous function $f:\R \rightarrow \C$ (other than the zero function) such that:
\[
\int_c^{c+R}  f(x) dx = 0,
\]
for all constants $c, R>0$. 
\item   \label{ex:part b, for Pompeiu}
More generally, in any dimension $d$, find a (nontrivial) continuous function $f:\R^d \rightarrow \C$ that allows the following integrals (taken over any ball of radius $r$) to vanish:
\[
\int_{B_d(r)+c} f(x) dx = 0,
\]
for all constants $c, R>0$. 
\end{enumerate}
Notes.  For part \ref{ex:part b, for Pompeiu}, 
it's advisable to think about the Fourier transform of the ball. 
 It is conjectured that for any bounded set $K$ with nonempty interior, the balls in this example are the only examples of objects that allow such nonzero continuous functions $f$ to exist.  This is known as the {\bf Pompeiu problem} - see also Question \ref{Pompeiu conjecture}.
 }
\end{prob}
 
\medskip
\begin{prob} \label{the integral of the example in terms of factorials}
\rm{
Show that when $p$ is a positive integer, the identity 
\eqref{identity of J-Bessel example} of Example \ref{ex:integral using Bessel functions} simplifies to a ratio of factorials. 
}
\end{prob}

\medskip
\begin{prob} \label{the integral of the example in terms of factorials}
\rm{
Show that the Hermite constant is very simply related to 
the densest lattice sphere packing density by:
\begin{equation}
\gamma_d = \left( \frac{ \delta^*(B_d) }{ \vol B_d }\right)^{\frac{2}{d}}.
\end{equation}
}
\end{prob}

\medskip
 \begin{prob}  
 \rm{  
 \label{positive FT over R^d}   
Using the idea of Exercise \ref{positive FT over R} in Chapter \ref{Fourier analysis basics}, 
and using the sum of two 
indicator functions of balls (with incommensurable radii) in $\R^d$, 
show that there exists a compactly supported function $f:\R^d \rightarrow \C$ such that 
 \[
 \hat f(\xi) >0,
 \]
for all $\xi \in \R^d$.
}
\end{prob}

\medskip
\begin{prob} $\clubsuit$
\label{proof of Lemma, equivalence of packing density}
Prove Lemma \ref{lem:equivalence of lattice packing density}.
\end{prob}


\bigskip
\chapter{\blue{
Shannon sampling, in one and several dimensions}
}
 \label{Chapter:Shannon sampling}

 \begin{quote}    
 ``It is easy to argue that real signals must be band-limited.  \\
 \quad     It is also easy to argue that they cannot be so.'' 
 
-- David Slepian   \index{Slepian, David}
 \end{quote}

\red{(Under construction)}

\section{Introduction}

Sampling theory consists in the  reconstruction of a continuous function with only a discrete or finite amount of data and has many applications in signal processing and other  engineering applications. At a first glance this task sounds impossible, however it can be done well in practice. One of the reasons for this success comes from the Fourier analysis, which deals with the representation of a function in terms of its ``frequencies", and functions without high frequencies (bandlimited)  represents very well the real-world signals. 

In one dimension, the classical example is a sound signal, and since typical humans can only hear sounds with frequencies smaller than $20$  kHz, the bandlimited assumption is appropriate.  
Examples in higher dimensions include images or MRI exams where higher frequencies are associated with random noises and measurement errors, more connected to the physical apparatus than the object being measured \cite{EpsteinBook}.  In this sense it is even desirable to remove the high 
 frequency information. 

More recently, the interest in bandlimited functions increased in the machine learning community, because
 it was observed that neural networks learns low frequencies faster and this might explain why they often generalize quickly from the training sets. 

On the other hand, by a basic uncertainty principle of Theorem \ref{basic uncertainty principle}, we know that a function with compact support can never be bandlimited, so representing an arbitrary function using
this class of functions is in general not exact.  It is therefore desirable to also give 
some theoretical results concerning the error of such approximations.


Here we introduce the classical sampling theorem by Shannon and Whittaker 
for one dimensional sampling, and then we study some of its generalizations to higher dimensions, 
where much less is known.  
An excellent introduction to Sampling Theory, from an expository as well as a rigorous perspective, 
is the book of J. R. Higgins \cite{Higgins1996}.


\bigskip
\section{The Shannon-Whittaker sampling Theorem}
\label{sec:one-dimensional-stuff}

Claude Shannon~\cite{Shannon1} showed how to reconstruct a complete signal $f$ by sampling it only discretely,
in a classical paper  that gave rise to the field of information theory.  
To accomplish this, Shannon used an interesting assumption, namely that the Fourier transform of $f$ 
vanishes outside of some interval.  

One of the main characters of this story is our old friend, the $\sinc$ function:
\begin{equation}\label{eq:sinc-def}
\sinc(x):= \int_{-\frac{1}{2}}^{\frac{1}{2}} e^{2\pi i \xi x} d\xi = \begin{cases}  
\frac{\sin(\pi x)}{\pi x},     &\mbox{if } x \not= 0 \\ 
1  & \mbox{if } x= 0,
\end{cases}
\end{equation}
which plays a central role in the sampling theory for functions in $\R$, because it turns out to be a building block for a basis of the Paley-Wiener space $PW_c$, as the Shannon-Whittaker sampling theorem shows. 

Reviewing some of the Fourier facts that we learned in Chapter  \ref{Fourier analysis basics}, 
we recall that if $f \in L^1(\R)$, then  $\hat f$ is uniformly continuous and $\hat f(\xi) \to 0$ as $|\xi| \to \infty$. 
So not every function can be the Fourier transform of some other function in $L^1(\R)$. 


In practice, we are often interested in functions that are not absolutely integrable, and yet possess a (conditionally convergent) Fourier transform, such as the important $\sinc$ function.  To resolve this issue, the theory progresses by first defining the  transform in the space $L^1(\R)$, and then extending the definition of $\hat f$ to all of $L^2(\R)$, by taking the limit $\lim_{n \to \infty} \int_{|x| < n}f(x)e^{-2\pi i x \xi}dx$.   This unique extension of the Fourier transform, 
from the $L^1(\R)$ space to the $L^2(\R)$ space, is sometimes called the Plancherel-Fourier transform.  From now on, we'll follow the usual Fourier convention and simply call both transforms ``the Fourier transform''.

For a given number $c>0$, a function $f\in L^2(\R)$ is called {\bf c-bandlimited} if 
\begin{equation}\label{first def of c-bandlimited}
\hat{f}(x) = 0 \text{ for all } x \not\in [-c, c].
\end{equation}
We will sometimes just say `bandlimited' if the $c$ is not contextually important.  A bandlimited function $f$ has a Fourier transform that decays at the `best possible rate', in the sense that its Fourier transform is identically zero outside the interval $[-c, c]$.  
It's easy to notice that any $c$-bandlimited function $f$ must be equal (almost everywhere) 
to an infinitely smooth 
 function, because by Fourier inversion, we have:
\begin{equation}\label{eq:bandlimit-inversion}
f(x) = \int_{\R} \hat f(\xi) e^{2\pi i \xi x} d\xi
=\int_{-c}^c \hat f(\xi) e^{2\pi i \xi x} d\xi.
\end{equation}
This identity implies that we can differentiate the last expression with respect to $x$ as many times as we like under the integral sign, because the integrand is a smooth function of $x$, and we are integrating over a compact domain.  Therefore $f$ is infinitely smooth. 
For simplicity, when considering a bandlimited function $f$, 
we will always assume that $f$ is also continuous, which is
consistent with the equality in~\eqref{eq:bandlimit-inversion}.

Given any $c>0$, we define the space of all $c$-bandlimited functions in 
$L^2(\R)$ by 
\[
PW_c:= \{ f \in L^2(\R) \mid \,  \hat f(\xi) = 0 \text{ for } \xi \notin  [-c, c],    \text{  and }  f
\text{ is continuous} \},
\]
called {\bf the Paley-Wiener space}  \cite{Higgins1996}.

\begin{figure}[htb]
 \begin{center}
\includegraphics[totalheight=4.6in]{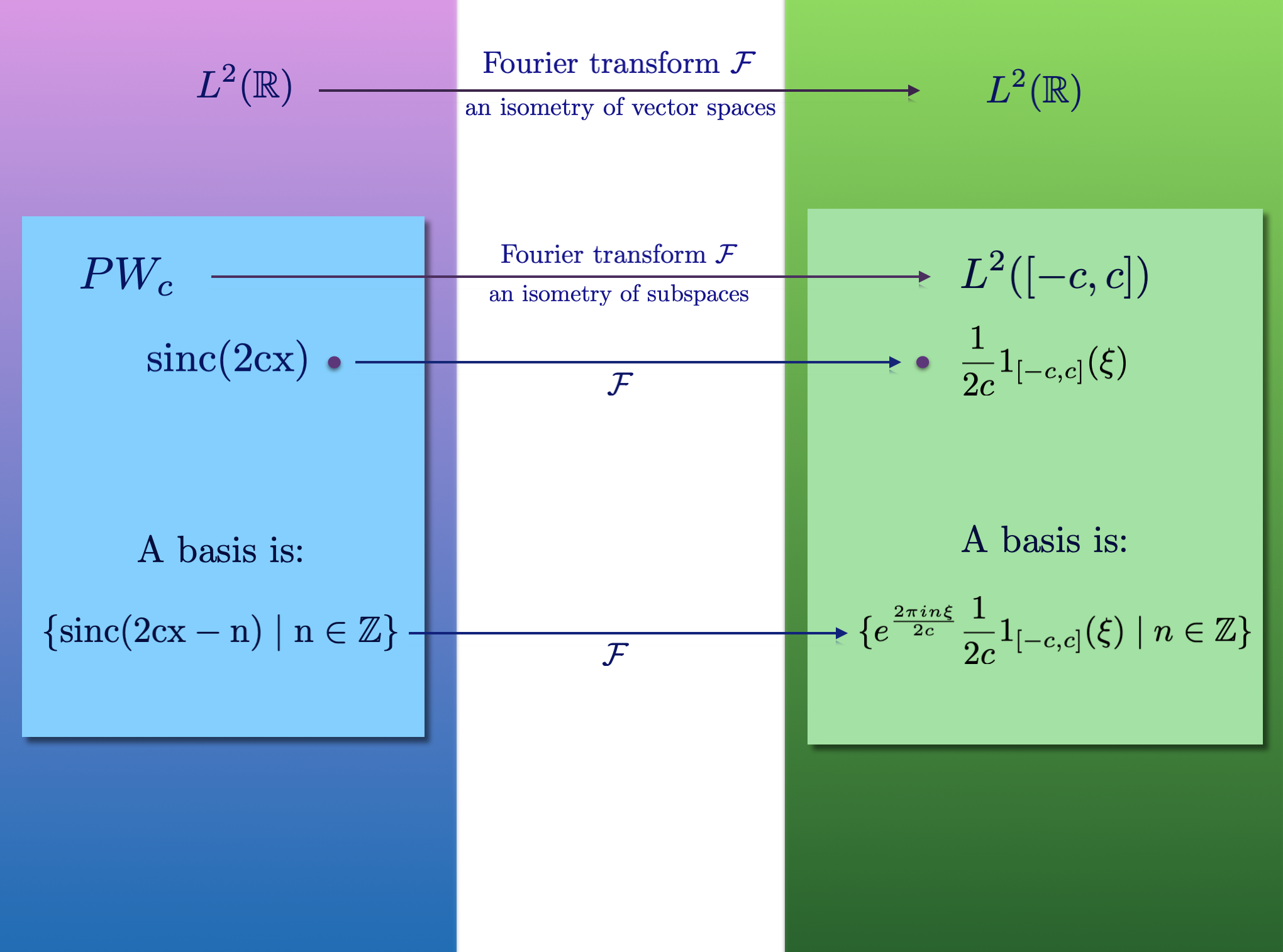}
\end{center}
\caption{We can visualize the isometries of the various spaces of functions, given by the Fourier transform. 
First, the Fourier transform $\mathcal{F}$ gives us an isometry from $L^2(\R)$ onto itself.  
Second, restricting attention to the subspace $PW_c \subset L^2(\R)$ of bandlimited functions, $\mathcal{F}$ also gives us an isometry from $PW_c$ to $L^2([-c, c])$, carrying the basis of translated $\sinc$ functions to the basis of exponentials. }
 \label{Shannon1}
\end{figure}


\begin{thm}[Shannon-Whittaker] \label{Shannon}
Suppose that $f \in PW_c$.  Then we have
\begin{equation}
\label{shannon_sampling}
f(x) = \sum_{n \in \Z} f\left(  \frac{n}{2c} \right)   {\rm sinc } \left( 2c  x - n  \right),
\end{equation}
and the series converges absolutely and uniformly over $\R$.
\hfill $\square$
\end{thm}
In other words, if we sample a $c$-bandlimited function $f$ at only the discrete set of points $\{ \frac{n}{2c} \mid n \in \Z\}$, we may reconstruct the whole function $f(x)$ for all $x\in \R$!  In the next two sections, we give two different proofs of Theorem \ref{Shannon}. 

The quantity $\frac{1}{2c}$ is called the {\bf sample spacing} and its reciprocal $2c$ is the  {\bf sampling rate}. 
So what the theorem says is that to reconstruct a function with bandlimit $c$, one has to sample at a rate~$2c$. Offhand, it seems rather incredible that some (non-periodic)  functions  $f:\R \rightarrow \C$ may be completely recovered by knowing only their values $f(n)$ at a discrete set of points. 
This phenomenon shows, in a sense, how the Paley-Wiener space is a very special subspace of $L^2(\R)$.


\bigskip
\section{The approach of G. H. Hardy}\label{sec:Hardy}

G. H. Hardy's proof \cite{Hardy41} of Theorem \ref{Shannon} is particularly interesting because it also answers the following informal question: 
\begin{question}\rm{[Rhetorical]}
How large is the space of bandlimited functions?
\end{question}
Hardy's approach also clarifies some of the underlying structure of  bandlimited functions. It relies on the following isometry.

\begin{lem}\label{isometry}
The Fourier transform $\F$ gives a bijection between the following two Hilbert spaces:
\begin{equation*}
    \F: PW_{c} \rightarrow 
    L^2(\left[-c, c\right]).
\end{equation*}
Moreover, this bijection is an isometry. 
\end{lem}
\begin{proof}
First, given any $f \in PW_{c} \subset L^2(\R)$, we need to show that $\hat f \in  L^2(\left[-c, c\right])$.
By definition $f$ is $c$-bandlimited, hence its Fourier transform can be naturally viewed as a function with domain $[-c,c]$. 
We need to show that $\hat f$ has a finite $L^2(\left[-c, c\right])$ norm. 
The following computation uses Parseval's identity, namely that 
$\| f \|^2_{L^2(\R)} = \|\hat f \|^2_{L^2(\R)}$:
\begin{align}
{\| \hat f \|}^2_{L^2(\left[-c, c\right])}
&:=\int_{-c}^{c} 
 |\hat f(\xi)|^2 d\xi 
 =\int_{\R} 
 |\hat f(\xi)|^2 d\xi 
=\int_{\R} 
 |f(x)|^2 dx 
 =: \| f \|^2_{L^2(\R)} 
 < \infty,
\end{align}
proving that
$\hat f \in  L^2(\left[-c, c\right])$.  

Conversely, given any $g \in  L^2(\left[-c, c\right])$, we need to show that $\F^{-1}(g) \in PW_c$.
We may extend $g$ to be equal to $0$ outside the interval $\left[-c, c\right]$, so that now $g \in L^2(\R)$.  By construction of $g$, we also 
have $g \in L^1(\R)$, so that now Lemma \ref{uniform continuity} guarantees that  $\F^{-1}(g)$ is uniformly continuous on $\R$.
Because the Fourier transform is an isometry of $L^2(\R)$, and $g \in L^2(\R)$, we also have $\F^{-1}(g) \in L^2(\R)$.   So now we have
$\F^{-1}(g) \in PW_c$.

Finally, the Fourier transform is invertible, and since we have 
$\|f-g\|_{L^2(\R)} = \|\hat f-\hat g\|_{L^2(\R)}$ by Parseval again, we have an isometry between the two Hilbert spaces $PW_c$ and $L^2([-c, c])$.
\end{proof}

\medskip
Hardy's insight is to consider an orthonormal basis for $L^2([-c,c])$ and then pull it back to an orthonormal basis for $PW_c$.  We recall the classical fact (Theorem \ref{Fourier series for periodic functions})  that the set of exponentials 
\[ \Big\{ \frac{1}{\sqrt{2c}}e^{\frac{2\pi i n x}{2c}} \mid n \in \Z\Big\}\]
form a complete orthonormal basis for the Hilbert space $L^2(\left[-c, c\right])$.  
Moreover, any $g \in L^2(\left[-c, c\right])$ has a unique representation in this basis (which we called its Fourier series), that converges in the $L^2$-norm on $[-c, c]$:
\begin{equation}\label{good series for now}
g(\xi) 
\underset{L^2([-c,c])}{=}
\sum_{m\in \Z} \hat g_m e^{\frac{2\pi i m \xi}{2c}},
\end{equation}
with coefficients equal to
\begin{equation} \label{good series for now: coefficients}
     \hat g_m = \frac{1}{2c}\big\langle g(\xi), e^{\frac{2\pi i m \xi}{2c}} \big\rangle 
     := \frac{1}{2c}\int_{-c}^c  g(\xi) e^{-\frac{2\pi i m \xi}{2c}} d\xi.
\end{equation}

Using the Fourier series~\eqref{good series for now}, we may expand  $\|g\|_{L^2([-c,c])} = \langle g, g\rangle^{1/2}$, obtaining
\begin{equation*}
\| g \|_{L^2([-c,c])} =  \Big( 2c \sum_{m \in \Z}  |\hat g_m|^2  \Big)^{1/2}.
\end{equation*}

Although we've only scratched the surface, we've already scratched it enough in order to prove 
Theorem \ref{Shannon}.

\bigskip
\begin{proof}[Proof of Theorem~\ref{Shannon}]
 
For any $f \in PW_{c}$, we know by Lemma~\ref{isometry} that $\hat f \in L^2(\left[-c, c\right])$, so that $\hat f$ has a Fourier series that converges in the $L^2$-norm on $[-c, c]$:
\begin{equation*}
\hat f(\xi) 
\underset{L^2([-c,c])}{=}
\sum_{m\in \Z} c_m e^{\frac{2\pi i m \xi}{2c}},
\end{equation*}
with coefficients equal to
$c_m = \frac{1}{2c}\int_{-c}^c \hat f(\xi) e^{-\frac{2\pi i m \xi}{2c}} d\xi = \frac{1}{2c}f(-\frac{m}{2c})$ by the Fourier inversion formula \eqref{eq:bandlimit-inversion}. It follows that
\begin{equation}
\label{eq:fhat-L2series}
\hat f(\xi) \underset{L^2([-c,c])}{=} \frac{1}{2c} \sum_{m\in \Z} f\left(\frac{m}{2c}\right)  e^{-\frac{2\pi i m \xi}{2c}}.
\end{equation}
From the orthonormality of the exponentials and Parseval's identity (Lemma~\ref{isometry}) we have
\begin{equation}\label{eq:series-fnorm}
\| f \|_{L^2(\R)} = \| \hat f \|_{L^2([-c,c])} =  \Big( \frac{1}{2c} \sum_{m \in \Z}  \Big|f\left(\frac{m}{2c}\right)\Big|^2  \Big)^{1/2}.
\end{equation}

We recall that the equality in the norm in~\eqref{eq:fhat-L2series} means that
\[
\lim_{N \to \infty} \bigg\| \hat f(\xi) - \frac{1}{2c}\sum_{|m| < N} f\left(\frac{m}{2c}\right)  e^{-\frac{2\pi i m \xi}{2c}} \bigg \|_{L^2([-c,c])} = 0.
\]
Using the Fourier inversion and the isometry stated in Lemma~\ref{isometry},
\begin{align*}
0 &= \lim_{N \to \infty} \bigg\| \int_{-c}^{c} \Big( \hat f(\xi) - \frac{1}{2c}\sum_{|m| < N} f\left(\frac{m}{2c}\right)  e^{-\frac{2\pi i m \xi}{2c}}\Big) e^{2\pi i \xi x}d\xi \bigg \|_{L^2(\R)}\\
&= \lim_{N \to \infty} \bigg\| \int_{-c}^{c} \hat f(\xi)e^{2\pi i \xi x}d\xi - \frac{1}{2c}\sum_{|m| < N} f\left(\frac{m}{2c}\right) \int_{-c}^{c} e^{-\frac{2\pi i m \xi}{2c}} e^{2\pi i \xi x}d\xi \bigg \|_{L^2(\R)}\\
&= \lim_{N \to \infty} \bigg\| f(x) - \sum_{|m| < N} f\left(\frac{m}{2c}\right) \frac{\sin(\pi(2cx-m))}{\pi(2cx-m)} \bigg \|_{L^2(\R)}\\
&= \lim_{N \to \infty} \bigg\| f(x) - \sum_{|m| < N} f\left(\frac{m}{2c}\right) \sinc(2cx-m) \bigg \|_{L^2(\R)},
\end{align*}
and therefore
\begin{equation}\label{eq:formulal2norm}
f(x) \underset{L^2(\R)}{=} \sum_{m \in \Z} f\left(\frac{m}{2c}\right) \sinc(2cx-m).
\end{equation}

To pass from the convergence in the norm to pointwise convergence, we need to show that the latter 
series converges uniformly, so that we can conclude that it represents a continuous function and hence by Lemma \ref{norm convergence plus absolute convergence implies equality} it is equal to $f$ everywhere. 

To prove the uniform convergence, we  make use of  the Cauchy-Schwartz inequality for infinite series, namely
\begin{equation}\label{C-S for infinite series}
\sum_{m=N}^\infty \Big|f\left(\frac{m}{2c}\right) \sinc(2cx-m)\Big| \leq 
\Big(\sum_{m=N}^\infty \Big|f\left(\frac{m}{2c}\right)\Big|^2\Big)^{1/2}
\Big(\sum_{m=N}^\infty \sinc^2(2cx-m)\Big)^{1/2}.
\end{equation}
The rest of the proof consists in showing that the right-hand side of \eqref{C-S for infinite series}
goes to zero as $N \rightarrow \infty$, uniformly for $x \in \R$. 
The same proof will also work for the series defined from $-N$ to $-\infty$. Together these results show that the expression in~\eqref{eq:formulal2norm} converges absolutely and uniformly over~$\R$, giving the result stated in the theorem.
From~\eqref{eq:series-fnorm}, we see that 
\begin{equation}\label{eq:fgoestozero}
\Big(\sum_{m = N}^\infty \Big|f\left(\frac{m}{2c}\right)\Big|^2\Big)^{1/2} \to 0,\quad \text{as } N \to \infty.
\end{equation}
Clearly
\[
\sum_{m = N}^\infty \sinc^2(2cx - m) \leq \sum_{m \in \Z} \sinc^2(2cx - m),
\]
and since the latter series are periodic function of $x$, with period $\frac{1}{2c}$, 
we may assume that $0 \leq x < \frac{1}{2c}$. 
For $m = 0$ and $1$, we note that $\sinc^2(2 c x - m) \leq 1$. For $m \geq 2$, we use the estimate
\[
\sinc^2(2cx -m) = \frac{\sin^2(2\pi cx - \pi m)}{(2\pi cx - \pi m)^2} \leq \frac{1}{\pi^2(m-1)^2},
\]
so that
\[
\sum_{m = 2}^\infty \sinc^2(2cx - m) \leq \sum_{m=2}^\infty \frac{1}{\pi^2(m-1)^2} 
= \frac{1}{6}.
\]
Similarly, for $m \leq -1$,
\[
\sum_{m = -\infty}^{-1} \sinc^2(2cx - m) \leq \sum_{m=1}^\infty \frac{1}{\pi^2 m^2} =\frac{1}{6}.
\]
We conclude that  for all $x \in \R$, 
$\sum_{m = N}^\infty \sinc^2(2cx-m) \leq \sum_{m \in \Z} \sinc^2(2cx-m) \leq 2 + \frac{1}{6}+ \frac{1}{6} = \frac{7}{3}$,
and therefore the series in~\eqref{C-S for infinite series} converges uniformly to $0$ as $N\rightarrow \infty$.
\end{proof}


\bigskip
It follows from this approach of G.H. Hardy, that despite being a subspace of $L^2(\R)$, the Paley-Wiener space $PW_{c}$ has 
a concrete, countable basis, which we record as follows. 
\begin{cor}\label{cor:sinc-orthonormal}
The set of translated $\sinc$ functions 
\begin{equation}
\{(\sqrt{2c}) \,{\rm sinc}(2cx-n) \mid n \in \Z\}
\end{equation}
is a complete orthonormal basis for the Hilbert space 
$PW_{c}$ of $c$-bandlimited functions.
\hfill $\square$
\end{cor}
\noindent
It is also worthwhile recording here the orthonormality of the $\sinc$ functions explicitly. For each 
$n, m \in \Z$, we have:
\begin{equation}\label{eq:sinc-orthonormal}
2c\int_{\R} \sinc(2cx-n) \sinc(2cx-m) dx = 
\begin{cases}
1 &\mbox{if } n = m,\\
0 &\mbox{otherwise.}
\end{cases}
\end{equation}

\bigskip
\section{An alternative proof, using Poisson summation}\label{sec:proof-Poisson}
\index{Poisson summation formula}

Here we give Shannon's proof of the classical Shannon-Whittaker sampling theorem (Theorem~\ref{Shannon}), with some added details.  
This proof uses the Poisson summation formula.  As we've seen several times before,  Poisson summation often simplifies proofs in surprising ways. 

To state the formula more precisely,  we use $\underset{L^1(\R)}{=}$ and $\underset{L^1([-c,c])}{=}$ to denote convergence in the $L^1$-norm, so that equality between functions holds almost everywhere but cannot be assumed at a specific point, unless we have an additional assumption like continuity.  

Assuming only that $f \in L^1(\R)$, the Poisson summation formula (See \cite{SteinWeiss})
states that the periodized function defined by the series
$
\sum_{n \in \Z} f(x+2cn)
$
converges in the norm of $L^1([-c,c])$ to a function whose Fourier expansion is
\begin{equation}\label{first Poisson summation-one dimension}
\sum_{n \in \Z} f(x+2cn)
\underset{L^1([-c,c])}{=}
\frac{1}{2c}  \sum_{m \in \Z}  \hat f\left(\frac{m}{2c}\right)  e^{\frac{2\pi i  m x}{2c} }.
\end{equation}

\begin{proof}[Proof of Theorem~\ref{Shannon}]
We begin with the Fourier series \eqref{first Poisson summation-one dimension}, which converges in the $L^1([-c, c])$ norm.

Step $1$. \   Our first goal will be to  
exchange the roles of $f$ and $\hat f$.  To justify this,
we begin by noting that our assumption  that $f \in PW_c$ implies 
$\hat f \in L^2([-c,c])$ by Lemma~\ref{isometry}, and $ L^2([-c,c]) \subset L^1([-c,c])$ by  
Lemma \ref{proper containment of L^2 in L^1 for torus}.    So we have $\hat f  \in L^1([-c,c] \subset L^1(\R)$,
allowing us to apply the same Poisson summation formula as above, together with Fourier inversion:
\begin{equation}\label{Poisson strikes again-one dimension}
\sum_{n \in \Z} \hat f(\xi +   2cn)
\underset{L^1([-c,c])}{=}  
\frac{1}{2c}   \sum_{m \in \Z} f\left( \frac{m}{2c} \right) e^{-\frac{2\pi i m \xi  }{2c}}.
\end{equation}

Step $2$. \  Since we are assuming that $f \in PW_c$ and thus $\hat f$ is supported on $[-c, c]$, we may use the indicator function $1_{[-c, c]}$, defined as $1_{[-c, c]}(\xi) = 1$ if $\xi \in [-c,c]$ and $1_{[-c, c]}(\xi) = 0$ otherwise, and write the trivial identity
\begin{equation} \label{step 2}
\hat f( \xi ) = 1_{[-c, c]}(\xi) \sum_{n\in \Z} \hat f(\xi +    2c n    ),
\end{equation}
for all real $\xi \not= c,  \xi \not= -c$.  The reason is that the series on the right-hand-side contains only one term, namely
the $n=0$ term $\hat f(\xi)$.

Step $3$.  \  Using the Poisson summation formula \eqref{Poisson strikes again-one dimension} above, together with \eqref{step 2}, we see that
\begin{equation}\label{step 3}
\hat f( \xi ) = 1_{[-c, c]}(\xi) \sum_{n\in \Z} \hat f(\xi + 2c n) 
\underset{L^1(\R)}{=} \sum_{n \in \Z} f\left(  \frac{n}{2c} \right)   \Big(
 \frac{1}{2c}  1_{[-c, c]}(\xi) e^{ - \frac{2\pi i n \xi }{2c}} \Big).
\end{equation}
We recall that the inverse Fourier transform of the interval $[-c, c]$ is
\[
 \int_\R 1_{[-c, c]}(\xi) e^{-2\pi i \xi x} d\xi  =2c \,  \sinc(2c x),
\]
so that after composing the $\sinc$ function with a translation, we know that the Fourier transform of  
$\sinc \left( 2c \left(x - \frac{n}{2c} \right)\right)$
is $\frac{1}{2c}  1_{[-c, c]}(\xi)  
e^{ - \frac{2\pi i n \xi }{2c}}$.
Multiplying both sides of  \eqref{step 3} by $e^{-2\pi i \xi x}$ and integrating term-by-term  over $x \in [-c, c]$, we get:
\[
f(x) \underset{L^1(\R)}{=}  \sum_{n \in \Z} f\left(  \frac{n}{2c} \right)  \sinc \left(  2c x -  n \right),
\]
applying Fourier inversion again on the left-hand side. 

Finally, we recall that we are assuming  $f$ is continuous, since $f\in PW_c$.  So to pass from the convergence in the norm to the pointwise convergence we may apply the same procedure from the first proof to conclude that the series on the right converges uniformly in $\R$ and hence also represents a continuous function.
\end{proof}

\bigskip
There are many different possible extensions of the Shannon-Whittaker sampling theorem to higher dimensions, and below we glimpse some of them below.



\section{Special properties of bandlimited and sinc functions} 
\label{sec:special}

Here we focus on some special properties of bandlimited functions.  We've already seen  in Section~\ref{sec:Hardy} that the space $PW_c$ is isometric to $L^2([-c,c])$, 
so many of its special properties comes from $L^2([-c,c])$,  and they are then pulled back via inverse Fourier transform.

The special case $c=\frac{1}{2}$ of Theorem \ref{Shannon} is worth pointing out:
\begin{equation}\label{Shannon simplified}
f(x) = \sum_{n \in \Z} f(n) \  {\rm sinc } \left( x - n  \right),
\end{equation}
a classical version of the Shannon-Whittaker formula.
The choice of $c = \frac{1}{2}$ means that we begin with the interval $[-\frac{1}{2}, \frac{1}{2}]$ in the frequency space;
this interval is a Voronoi cell of the integer lattice $\Z$. 

\begin{example}  \label{ex-1}
\rm{
What happens if we apply the Shannon-Whittaker formula 
\eqref{Shannon simplified} to the $\sinc$ function itself?  Let's try it!   

With  $f(x):= \sinc(y - x)$, and any fixed $y\in \R$, we have:
\begin{equation}
\sinc(y - x) =  \sum_{n \in \Z}  \sinc(y - n) \, \sinc( x - n ).
\end{equation}
As a special case, if we let $x=y$, we get:
\begin{equation}
1= \sinc(0) =  \sum_{n \in \Z}   \sinc^2(x - n).
\end{equation}
\hfill $\square$
}
\end{example}

In Corollary \ref{cor:sinc-orthonormal} we showed that the functions $\sinc(x-n)$ with $n \in \Z$ 
form an orthonormal basis for the space of bandlimited functions.  This has some nice consequences. 

\begin{thm}   \label{thm:discrete-inner-product}
If $f$ and $g$ are $\tfrac{1}{2}$-bandlimited,
then
\[
\int_\R f(x) \overline{g(x)} dx
= \sum_{n \in \Z}  f(n) \overline{g(n)}.
\]
\end{thm}
[{\bf Intuitive proof}] \ If we work formally, then we can use the orthonormality of the sinc functions \eqref{eq:sinc-orthonormal},
together with \eqref{Shannon simplified} to quickly see that:
\begin{align*}
\langle f(x), g(x) \rangle &= 
\Big\langle \sum_{m \in \Z} f(m) \sinc(x-m), \sum_{n \in \Z} g(n)\sinc(x-n) \Big\rangle \\
&= \sum_{m, n \in \Z} 
f(m)\overline{g(n)} \Big\langle   \sinc(x-m),  \sinc(x-n) \Big\rangle \\
&=\sum_{n \in \Z}f(n) \overline{g(n)},
\end{align*}
and we're done.  Although this intuitive proof may seem `fast and loose', these steps can
 be made rigorous if we would prove just a bit more about Hilbert spaces, 
 because the Paley-Wiener space $PW_c$ is a Hilbert space, and the translated 
 functions  $\sinc(x-n)$ are a basis
 for this Hilbert space. 

One important case of the previous theorem is when $g(x) := \sinc(x - y)$, which combined again with the Shannon-Whittaker formula
\eqref{Shannon simplified} results in the next theorem.

\begin{thm}     \label{thm:ReproducingSinc}
The space $PW_{\frac{1}{2}}$ is a space with a reproducing kernel $\sinc(x-y)$, which means by definition 
that any $f \in PW_{\frac{1}{2}}$ can be written as
\[
f(x) = \int_\R f(y) \sinc(x-y) d y.
\]
\end{thm}
\begin{proof}
For $x \in \R$, take $g(y) := \sinc(x - y)$ in Theorem~\ref{thm:discrete-inner-product}:
\[
\int_R f(y) \sinc(x-y) d y = \sum_{n \in \Z} f(n) \sinc(x-n) = f(x),
\]
where the second equality follows from the Shannon-Whittaker formula.
\end{proof}

\bigskip
We also have the following somewhat surprising properties of bandlimited functions on $\R$.  First we recall
that we called  $f:\R^d \rightarrow \C$ a \emph{nice} function  if  $f, \hat f \in L^1(\R^d)$, and $f$ satisfies the Poisson summation formula:
\begin{equation}  \label{Poisson summation again} 
\sum_{n \in \Z^d}  f(n+ x)  =  \sum_{\xi \in \Z^d} \hat f(\xi) e^{2\pi i \langle \xi, x \rangle},
\end{equation}
valid pointwise for each $x \in \R^d$. 

\begin{thm} 
\label{bandlimited, Poisson summation} 
Let $f:\R\rightarrow \C$ be a nice function, such that $f$ is $c$-bandlimited.   Then we have:
\begin{enumerate}[(a)]
   \item   \label{firstpart} 
        We have, for each $k>c$,
        \begin{equation}\label{exact Riemann approximation}
          \frac{1}{k}\sum_{n\in\Z} 
          f\left(\frac{n}{k}\right)
           = \int_{\R} f(x) dx,
        \end{equation}
        \blue{We note that the identity \eqref{exact Riemann approximation}
        can be interpreted to mean that the Riemann approximation to the integral is always exact for such an $f$, provided that the step size is $\Delta x:=\frac{1}{k}<\frac{1}{c}$.}
    \item  \label{secondpart} 
        For all $a, k$ with $a>c$ and $k > a+c$, we have
    \begin{equation} 
     \sum_{n\in\Z} f\left(\frac{n}{k}\right)
       e^{\frac{2\pi i n a}{k}}=0.
    \end{equation}
\end{enumerate}
\end{thm}
\begin{proof}
To prove \ref{firstpart}, we use Poisson summation \eqref{Poisson summation again}, with $\L:= \frac{1}{k} \Z$:
\begin{equation} \label{equality of sum and integral}
        \sum_{n\in\Z} 
          f\left(\frac{n}{k}\right)=
          \sum_{\xi \in \L} f(\xi) = 
          \frac{1}{\det \L}  
         \sum_{m \in \L^*} \hat f(m) =
         k \sum_{m \in \Z} \hat f(mk) = k\int_\R f(x) dx,
\end{equation}
which is the desired identity.  In the last equality we used the assumption that the indices of summation satisfy $ |mk| >  c$, 
for $m\not=0$, 
so that $\hat f(mk) = 0$ because $f$ is $c$-bandlimited by assumption. 
We also used the fact that 
$\hat f(0) = \int_\R f(x) dx$. 


\noindent
To prove \ref{secondpart}, we apply the following small variation of Poisson summation:
 \begin{equation}\label{magic2}
       \frac{1}{k} 
     \sum_{n \in \Z} f\left(  \frac{n}{k} \right) e^{\frac{2\pi i n a }{k}}
       =\sum_{n\in \Z} \hat f(-a+kn),
\end{equation}
which follows quickly from the Poisson summation formula given above in \eqref{Poisson summation again} (Exercise \ref{Exercise:PoissonSummation1}).
But by the assumption that $f$ is $c$-bandlimited, we also have
\begin{equation}
       \sum_{n\in \Z} \hat f(-a+kn)=0,
\end{equation}
provided that 
\begin{equation}\label{constraint}
|-a+kn|>c, \ \text{ for all } n\in \Z. 
\end{equation}
For $n=0$, we see that a necessary condition for \eqref{constraint} is $|a|>c$.  Geometrically, \eqref{constraint} tells us that the arithmetic progression $\{ kn -a\}_{n \in \Z}$ does not intersect the interval $[-c, c]$.  It is easily checked that the additional constraint  $k>c+a$ gives us a sufficient condition for \eqref{constraint} to hold.
\end{proof}

As is easily observed, sums and products of bandlimited functions are again bandlimited.  In particular, more precise statements such as the following are possible.  

\begin{lem}\label{lemma:product}
Suppose that $f$ is $c$-bandlimited, and $g$ is $d$-bandlimited. \\
Then  $fg$ is $(c+d)$-bandlimited.
\end{lem}
\begin{proof}  By assumption, $\hat f(\xi) = 0$ outside of $[-c, c]$, and 
$\hat g(\xi) = 0$ outside of $[-d, d]$.   We need to show that $\widehat{fg}$ vanishes outside the interval $[-c-d, c+d]$.
We also have, by assumption, that $f, g \in L^2[-c, c]$, and since $L^2[-c, c] \subset L^1[-c, c]$ (Lemma \ref{proper containment of L^2 in L^1 for torus}), it follows that $f, g \in L^1[-c, c]$.  Therefore the convolution 
  Lemma \ref{convolution theorem} \ref{convolution under FT}
   applies: 
   \[
   \widehat{(fg)}(\xi)=(\hat f*\hat g)(\xi).
   \]
   We know that the support of the convolution $\hat f(\xi)*\hat g(\xi)$ is contained in the closure of the Minkowski sum \index{Minkowski sum}
of the individual supports of $f$ and $g$ (by Exercise \ref{support of convolution}), which is equal to $[-c, c]+[-d, d] = [-c-d, c+d]$. 
\end{proof}

\bigskip
\begin{example}
\rm{
Here are some fun consequences of Theorem \ref{bandlimited, Poisson summation}.  Let's fix any $\epsilon >0$.
By Theorem  \ref{bandlimited, Poisson summation}, part \rm{\ref{secondpart}} , we can pick $k=1, c=\frac{1}{2}-\epsilon$, and $a=\frac{1}{2}$, all of 
which satisfy the hypothesis, so that $f$ is $(\frac{1}{2}-\epsilon)$-bandlimited by definition.  We then have
\begin{equation} 
 0=  \sum_{n\in\Z} f(n)
       e^{\pi in}
  = \sum_{n\in\Z} f(n)
       (-1)^n.
\end{equation}
Seperating the lattice sum into $n$ even and $n$ odd, we have 
\begin{equation}
\sum_{m\in \Z} f(2m) 
=\sum_{m \in \Z} f(2m+1).
\end{equation}
Generalizing the latter identity, we fix any positive integer $N$, and we  let  $k=\frac{2}{N}, c=\frac{1}{N}-\epsilon$, and $a=\frac{1}{N}$,  so that $f$ is $(\frac{1}{N}-\epsilon)$-bandlimited by definition.
By Theorem 
\ref{bandlimited, Poisson summation}, part \rm{\ref{secondpart}}:
\begin{equation} 
0=     \sum_{n\in\Z} f\left(\frac{n}{k}\right)
       e^{\frac{2\pi i n a}{k}}
=\sum_{n\in\Z} f\left(\frac{Nn}{2}\right)
      (-1)^n,
\end{equation}
so that we get the identity
\begin{equation} 
\sum_{m\equiv 0\text{ mod N }} f(m)
 =\sum_{m\equiv 0\text{ mod N }} 
 f\Big(m+ \frac{N}{2}\Big).
\end{equation}
}
\hfill $\square$
\end{example}

\bigskip
\begin{example}
With  $k=1, c=\frac{1}{3}-\epsilon$, and $a=\frac{1}{3}$, part \rm{\ref{secondpart}} of Theorem 
\ref{bandlimited, Poisson summation}
gives:
\begin{align*} 
0&=\sum_{n\in\Z} f(n) e^{\frac{2\pi i n}{3}}
=\sum_{n\equiv 0\text{ mod 3 }} f(n) 
+ \omega\sum_{n\equiv 1\text{ mod 3 }} f(n)
+ \omega^2\sum_{n\equiv 2\text{ mod 3 }} f(n),
\end{align*}
where $\omega:= e^{2\pi i / 3}$.
\hfill $\square$
\end{example}

\bigskip
\begin{example}
\rm{
Consider $\sinc^N(x / \pi)$, which has the bandlimit $c:= \frac{N}{2\pi}$. 
By Theorem~\ref{bandlimited, Poisson summation}, the strange relation
\begin{equation} \label{example:strange 1}
\sum_{n \in \Z} \sinc^N\left(\frac{n}{\pi}\right) = \int_{\R} \sinc^N\left(\frac{x}{\pi}\right) dx,
\end{equation}
holds for $N = 2, \dots, 6$, because in this range we have
$c=\frac{N}{2\pi} \leq \frac{6}{2\pi} < 1=: k$.  It's also true for $N=1$, with some care:
\[
\lim_{M\rightarrow \infty}  \sum_{|n| < M \atop n \in \Z}
 \sinc \left(\frac{n}{\pi}\right) = \lim_{M\rightarrow \infty}   \int_{-M}^M  \sinc \left(\frac{x}{\pi}\right) dx.
\]

It turns out that this identity fails, however, for $N \geq 7$. Indeed, by 
Poisson summation~\eqref{Poisson summation again}, for a nice function $f$ we have:
\[
\sum_{n \in \Z}f(n) = \sum_{m \in \Z}\hat f(m) = \int_\R f(x) dx + \sum_{m \in \Z \setminus \{0\}} \hat f(m).
\]
Taking $f(x) := \sinc^N\left( \frac{x}{\pi}\right)$, we see that 
that the last sum is zero when $N \leq 6$ and positive when $N \geq 7$, since $\hat f$ has support  $[-\frac{N}{2\pi}, \frac{N}{2\pi}]$ and is positive inside this interval.

In a similar manner to eq. \eqref{example:strange 1}, we have:
\begin{equation} \label{example:strange 2}
\sum_{n \in \Z}\, \prod_{k = 0}^N \sinc\left(\frac{n}{(2k+1)\pi}\right) = \int_{\R}\, \prod_{k = 0}^N \sinc\left(\frac{x}{(2k+1)\pi}\right) dx,
\end{equation}
holds for $N = 0, \dots, 40248$, since for these $N$ we have $1 + \frac{1}{3} + \dots + \frac{1}{2N+1} < 2\pi$. 
It can be also checked that the equality above fails for $N = 40249$. These facts are easy corollaries of Theorem~\ref{bandlimited, Poisson summation}, but may seem surprising when taken out of this context.  The identities \eqref{example:strange 1}
and \eqref{example:strange 2} appeared in \cite{Baillie}.
}
\hfill $\square$
\end{example}


\bigskip
\section{Shannon sampling in higher dimensions}

The first research into higher-dimensional Shannon-type sampling theorems, as far as we know,
 was 
the work of Petersen and Middleton~\cite{PetersenMiddleton62}.   We'll also follow a bit of 
Chapter $14$ from Higgins~\cite{Higgins1996}.

For a convex body $\P$, we say that a function $f$ is {\bf $\P$-bandlimited} if $\hat f$ vanishes outside of $\P$.  
We note that this does not preclude the possibility that $\hat f$ may  only be nonzero on some proper subset of $\P$. 
Assuming that $f$ is real-valued, we know that the image of $\hat f$ is  symmetric about the origin 
(Lemma \ref{symmetric iff FT is real}); 
so the assumption that $\P$ is symmetric is natural.

By analogy with the $1$-dimensional Paley-Wiener space $PW_c$, we define the 
{\bf Paley-Wiener space} of $\P$-bandlimited functions  in $L^2(\R^d)$:
\begin{equation}
    PW_\P:= \{ f \in L^2(\R^d) \mid 
    f \text{ is continuous and  $\P$-bandlimited} \}. 
\end{equation}
     
A new twist in higher dimensions is the strong distinction between packing and tiling, 
so the following question motivates some of these research directions.

\begin{question}
Given a convex $d$-dimensional body $\P \subset \R^d$, suppose we want to have a sampling theorem for functions that are 
$\P$-bandlimited.   Does $\P$ have to tile $\R^d$ by translations, or is it sufficient to consider a packing of $\P$ by some lattice $\L$? 
\end{question}

Interestingly, we don't observe this distinction in dimension $1$, because optimal packing and tiling are equivalent.  But they are quite different in dimensions $d \geq 2$. 
Luckily, our elementary $1$-dimensional Lemma~\ref{isometry} does extend directly to our new $d$-dimensional setting. 
\begin{lem}
\label{isometry lemma for R^d}
Let $\P$ be a bounded convex body in $\R^d$. The Fourier transform $\F$ is an isometry between the two Hilbert spaces:
\begin{equation*}
    \F: PW_{\P} \rightarrow 
    L^2(\P).
\end{equation*}
\end{lem}
\begin{proof}
Given $f \in PW_{\P} \subset L^2(\R^d)$, by definition $\mathrm{supp}(\hat f) \subseteq \P$, 
so using Parseval's identity we have: 
\begin{align}
{\| \hat f \|}^2_{L^2(\P)} 
&:=\int_\P 
 |\hat f(\xi)|^2 d\xi =
 \int_{\R^d} 
 |\hat f(\xi)|^2 d\xi 
=\int_{\R^d} 
 |f(\xi)|^2 d\xi =: \| f \|^2_{L^2(\R^d)} 
 < \infty,
\end{align}
which shows that $\hat f \in  L^2(\P)$.
Using the fact that the Fourier transform is an isometry of $L^2(\R^d)$, and is invertible, we are done.
\end{proof}

Now this chapter comes full circle with the goals of the previous chapters:  to better understand the Fourier transform of a convex body.  The following result helps.
\bigskip
\begin{thm}[Higher-dimensional sampling formula]
\label{Higher dimensional sampling}
Suppose we have a lattice packing for a symmetric convex body $\P$, with a lattice~$\L^*\subset \R^d$.
If $f \in PW_\P$, then
then we can reconstruct the function $f$ completely by sampling it only at the  lattice points of $\L$:
\[
f(x) = \det \L   \sum_{n \in \L} f(n)   \hat 1_P(x-n),
\]
and the series converges absolutely and uniformly over $\R^d$.
\end{thm}
\begin{proof}
[Proof of Theorem~\ref{Higher dimensional sampling}]
The assumption that $f \in PW_\P$, together with $\P$ being compact, implies that 
 $\hat f \in L^2(\P) \subseteq L^1(\P) \subseteq L^1(\R^d)$.   Now we may use the Poisson summation formula
 \eqref{Poisson summation again}, but with $f$ replaced by $\hat f$, and with $\L$ replaced by $\L^*$:
 
\begin{equation}\label{Poisson strikes again}
\sum_{m \in \L^*} \hat f(\xi + m)
\underset{L^1(\R^d/\L^*)}{=}  
\det \L \sum_{n \in \L} f\left( n \right) e^{-2\pi i \langle \xi,  n \rangle},
\end{equation}
where we also have used that $\hat{\hat{f}}(m) = f(-m)$.
Because $\hat f$ is supported on $\P$, we have by definition $\sum_{m \in \L^*} \hat f(\xi + m) = \hat f(\xi)$, 
so that we may write
\begin{equation} \label{step 2.1}
\hat f( \xi ) = 1_{\P}(\xi) \sum_{m \in \L^*} \hat f(\xi + m),
\end{equation}
for all  $\xi\in \R^d$ that do not lie on the boundary of $\P$.  

Because of our packing assumption, all of the translated supports of $ \hat f(\xi + m)$ are disjoint, 
as $m$ varies over the lattice   $\L^*$. In other words, these supports are 
\begin{equation*}
\{     \supp( \hat  f ) + m   \mid \, m  \in \L^*     \}   \subseteq     \{\P+ m  \mid \, m  \in \L^*\}, 
\end{equation*}
a disjoint collection of translates of $\P$.
This means that the latter identity \eqref{step 2.1} holds
 because the series on the right-hand-side contains only one term, namely
the $m=0$ term $\hat f(\xi)$.
Next, we combine~\eqref{Poisson strikes again} with~\eqref{step 2.1} to get
\begin{equation}\label{step 3.1}
\hat f( \xi ) \underset{L^1(\R^d)}{=}  \det \L       \sum_{n \in \L} f\left( n \right)  1_{\P}(\xi)  e^{-2\pi i \langle \xi,  n \rangle}. 
\end{equation}
Now we'd like to take  the inverse Fourier transform of both sides of \eqref{step 3.1}.  
We'll use the following elementary identity, for a fixed $n$:
\[
\F^{-1} \big(  1_{\P}(\xi)  e^{-2\pi i \langle \xi,  n \rangle}  \big)(x) =   \F^{-1}(1_{\P})(x - n) 
= \int_\P e^{2\pi i \langle  \xi,  x-n  \rangle}d\xi = \hat 1_\P(x-n).
\]
We finally arrive at
\begin{align*}
f(x) \underset{L^1(\R^d)}{=}  \det \L \sum_{n \in \L} f\left( n \right) \hat 1_{\P}(x - n).
\end{align*}
To finish the proof, we just to pass from the $L^1$-convergence of the latter identity, 
to its pointwise and uniform convergence.
The series
$\sum_{n \in \L} f\left( n \right) \hat 1_\P(x-n)$ converges
uniformly on $\R^d$ by  a standard argument, as in \cite{Higgins1996} for example.
\end{proof}

Let's conduct a sanity check and verify that  Theorem \ref{Higher dimensional sampling}  indeed generalizes
Theorem~\ref{Shannon}, 
 the classical $1$-dimensional Shannon-Whittaker sampling theorem.
   In the one dimensional case, the lattice $\L$ is just the sampling domain $\{\frac{n}{2c} \mid n \in \Z\}$ and 
 hence $\det \L = \frac{1}{2c}$, while $\P:= [-c,c]$ is simply an interval.  
Therefore:
\[
\phi(x) = \int_{-c}^c e^{2\pi i \xi x}d\xi 
= \frac{1}{2\pi i x}e^{2\pi i xc} - \frac{1}{2\pi i x}e^{-2\pi i xc} 
= \frac{1}{\pi x} \sin(2 \pi c x)  
= 2c \, \sinc(2cx).
\]
The formula from Theorem~\ref{Higher dimensional sampling}  reduces to the formula from Theorem~\ref{Shannon}:
\begin{align*}
f(x) = \det \L \sum_{n \in \L}f(n)\phi(x-n) 
&= \frac{1}{2c} \sum_{n \in \Z} f\Big(\frac{n}{2c}\Big) 2c \, \sinc\Big(2c\big(x-\frac{n}{2c}\big)\Big)  \\
&=  \sum_{n \in \Z} f\Big(\frac{n}{2c}\Big)  \sinc(2cx-n),
\end{align*}
which is the Shannon-Whittaker sampling formula.

Intuitively, if we pick a  larger set $\P$, then the vectors from $\L^*$ will  have to be more widely spaced  in order to satisfy the packing condition
for $\P$.
Therefore in Theorem \ref{Higher dimensional sampling}, we will need to sample from a denser lattice $\L$, due to the reciprocal relation $(\det \L )(\det \L^*) = 1$. 

On the other hand, for a given sampling lattice $\L$, in this multidimensional case we can consider infinitely many different
 bodies $\P$ that form a packing of $\R^d$ with the lattice $\L$. 
  One of the most natural choices for such a convex set $\P$ is the Voronoi cell of $\L^*$.

In closing, we note that it is impossible for both $f$ and $\hat f$ to be simultaneously bandlimited, by the basic uncertainty principle,
Theorem \ref{basic uncertainty principle}.    However, in practice, if we are given a function $f\in L^2(\R^d)$ that is not bandlimited, 
we can form a sequence of bandlimited functions that approach $f$ as $n\rightarrow \infty$, as follows.   To make $\hat f$ compactly supported, we'll multiply $\hat f$ by $1_{[-n, n]^d}$, the indicator function of the cube.  
 Pulling things back to the space domain, we have:
\begin{align}
\F^{-1} \left(    
1_{[-n, n]^d} \, \hat f   
\right)
&=  \F^{-1} \left(      1_{[-n, n]^d}     \right) * \F^{-1} (\hat f)  \\   \label{inverting sinc^d}
&= \sinc^d* f.
\end{align}
So if we define  $f_n:=  \sinc^d* f$, then $\F(f_n) = 1_{[-n, n]^d} \, \hat f$, a compactly supported function that  is bandlimited to the cube $ [-n, n]^d$.   

The careful reader might notice that in \eqref{inverting sinc^d}, we are really 
applying the Fourier inversion formula on $L^2(\R^d)$ (as opposed to the Fourier inversion formula on $L^1(\R^d$)).
We do this because although  $\sinc^d(x) \notin L^1(\R^d)$,  we do have $\sinc^d(x) \in L^2(\R^d)$.


\bigskip \bigskip

\section*{Notes} \label{Notes.chapter.Shannon}
\begin{enumerate}[(a)]
\item  John Higgins' book \cite{Higgins1996}, Chapter $6$, has an excellent account of the Paley-Wiener space, and its connections to the Paley-Wiener theorem, which also answers the question: ``how may we analytically continue bandlimited functions of a real variable, to $\C$?''   Moreover, Higgin's book has more mathematical rigor than many other books that treat sampling.

\item  For further reading, the following articles are of interest:  \cite{Entezari09}, \cite{Unser00}, \cite{Ye12}. 

\item Interesting relations between rates of convergence of bandlimited functions, 
Nikol'skij type functions spaces,  and Plancherel-Polya type inequalities is given in \cite{SchmeisserSickel2000}.
 
\end{enumerate}


\section*{Exercises}
\addcontentsline{toc}{section}{Exercises}
\markright{Exercises}

 \begin{quote}    
 ``The only true wisdom is in knowing you know nothing.''
 
 --  Socrates   \index{Socrates}
 \end{quote}

\begin{prob}   
{\rm
By using an example, show that a bandlimited function $f\in L^2(\R)$ may not be in $L^1(\R)$. 
}
\end{prob}

\medskip
\begin{prob}
{\rm
Let $f \in PW_c$, and fix any $x_0 \in \R$.  Prove  that $f$ is completely determined
by the samples 
\[
\left\{     f \left(x_0 +\frac{\pi n}{c}   \right)  \mid   n \in \Z \right\}.
\]
}
\end{prob}

\medskip
\begin{prob}  $\clubsuit$    \label{Exercise:PoissonSummation1}
{\rm
Here we give another small variation on Poisson summation, namely that for any $a, k \in \R$, we have
 \begin{equation}
       \frac{1}{k} 
     \sum_{n \in \Z} f\left(  \frac{n}{k} \right) e^{\frac{2\pi i n a }{k}}
       =\sum_{n\in \Z} \hat f(-a+kn),
\end{equation}
where $f:\R\rightarrow \C$ is a nice function (in the sense of \eqref{Poisson summation again}). 
}
\end{prob}

\medskip
\begin{prob} 
{\rm
Consider the function $f(x):= \rm{sinc}^2(x) := \left( \frac{\sin(\pi \xi)}{\pi \xi} \right)^2$, when $x\not\in \Z$.
\begin{enumerate}[(a)]
\item Show that $f$ is $1$-bandlimited.
\item Show that for each $k > 1$, we have:
\[
          \frac{1}{k} \sum_{n\in\Z} 
         \rm{sinc}^2\left(\frac{n}{k}\right)
           = \int_{\R} \rm{sinc}^2(x) dx, 
\]
using any results from this book.
\item Show that $\int_{\R} \rm{sinc}^2(x) dx=1$.
\end{enumerate}
}
\end{prob}

\appendix

\bigskip
\chapter{The dominated convergence theorem, and other goodies}
 \label{Appendix A}

 A frequent question that comes up in proofs is ``when may we take the limit inside the integral''? 
 A general tool that allows us to do so is the Dominated convergence theorem.   Here we remind the reader of 
 some of the basic results from real analysis, but we skip the proofs and give references for them.   For our purposes, we only need these results in Euclidean spaces, although all of these theorems have extensions to arbitrary measure spaces.   All functions here are assumed to be measurable.

 \begin{thm}[Fatou's lemma]  \label{Fatou}
 \index{Fatou's lemma}
Fixing any subset $E\subset \R^d$, let  $f_n:E \rightarrow [0, \infty)$ be a sequence of  nonnegative functions.  Then we have:
 \begin{equation}
 \int_{E} \lim \inf \, f_n(x) dx \leq     \lim \inf  \int_{E}  \, f_n(x) dx.
 \end{equation}
\hfill $\square$
 \end{thm}  
The inherent flexibility in {\bf Fatou's lemma} allows it to be useful in many different contexts, because the lim inf $f_n$ always exists, and are even allowed to be equal to $\pm$ infinity.   In fact, Fatou's lemma is the main tool in proving Lebesgue's dominated convergence theorem, below.

 Another essential fact for us is {\bf Fubini's theorem}, which allows us to interchange integrals with integrals, and series with integrals, for product spaces.
 If we write $\R^d = \R^m \times \R^n$, and we denote a point $z\in \R^d$ by $z:= (x,y)$, then we may also write
 $f(z):= f(x, y)$.

 \begin{thm}[Fubini]  \label{Fubini}
 \index{Fubini's theorem}
 Let $f \in L^1(\R^d)$.  Then:
 \begin{equation}
 \int_{\R^d} f(z) dz = \int_{\R^n}  \left(  \int_{\R^m} f(x, y) dx \right) dy,
 \end{equation}
 and 
  \begin{equation}
 \int_{\R^d} f(z) dz = \int_{\R^m}  \left(  \int_{\R^n} f(x, y) dy \right) dx.
 \end{equation}
\hfill $\square$
 \end{thm}  
There is also a version of Fubini's theorem that uses the counting measure in one of the
factors of  $\R^m \times \R^n$, giving us:
\begin{equation}\label{Fubini for sums and integrals}
\sum_{\xi \in \Z^n}  \left(  \int_{\R^m} f(x, \xi) dx \right) =  \int_{\R^m}    \left(   \sum_{\xi \in \Z^n}    f(x, \xi) \right)   dx.
\end{equation}
(See \cite{PaulSally1}, p. 220,  for a proof of Theorem \ref{Fubini})


\medskip
\section{The Dominated Convergence Theorem}
\begin{thm}[Dominated convergence theorem]  \label{Dominated convergence theorem}
\index{Lebesgue dominated convergence theorem}
 \ Suppose that we have a sequence of  functions $ f_n(x):\R^d \rightarrow \C$, for $n = 1, 2, 3, \dots $, and suppose 
 there exists a limit function $f(x) =\lim_{n\rightarrow \infty} f_n(x)$, valid for all $x\in \R^d$. 
 
 If there exists a function  $g \in L^1(\R^d)$ such that for all $x \in \R^d$, we have:
 \[
\left| f_n(x) \right| \leq g(x),  \quad n = 1, 2, 3, \dots
 \]
  then:
 \begin{enumerate}[(a)]
  \item  $f \in L^1(\R^d)$.
  \item   $ \lim_{n\rightarrow \infty} \int_{\R^d} | f_n(x) - f(x) | dx = 0$.
  \item    And finally, we may interchange limits and integrals:
  \[
  \lim_{n \rightarrow \infty}  \int_{\R^d} f_n(x) dx =  \int_{\R^d} f(x) dx.
  \]
 \end{enumerate}
\hfill $\square$
 \end{thm}  
 
 Theorem \ref{Dominated convergence theorem} is sometimes called the \emph{Lebesgue dominated convergence theorem}, honoring the work of Lebesgue.  
 There is a useful application of Lebesgue's dominated convergence theorem, which allows us to interchange summations with integrals as follows. 
 
 \medskip
\begin{thm} \label{Application of dominated convergence}
 \ Suppose that we have a sequence of  functions $ f_n(x):\R^d \rightarrow \C$, such that 
 \[
 \sum_{n=1}^\infty \int_{\R^d}  | f_n(x) | dx < \infty.
 \]    
 Then the series   $ \sum_{n=1}^\infty f_n(x)$
 converges for all $x\in \R^d$, and we have:
 \[
\sum_{n=1}^\infty      \int_{\R^d}    f_n(x)   dx =    \int_{\R^d}    \sum_{n=1}^\infty   f_n(x)   dx.
 \]
\hfill $\square$
 \end{thm}

(See \cite{RudinGreenBook}, p. 26 for a proof of Theorem \ref{Dominated convergence theorem}, and 
\cite{RudinGreenBook}, p. 29 for a proof of Theorem \ref{Application of dominated convergence})


\section{Big-O and Little-o}

Very often we'd like to compare, in a quick-and-dirty way that avoids uncountably many details, 
how fast two functions grow.   We review here two of the most common ways to do this.   

Suppose we are given two functions $f, g:\R^d \rightarrow \C$.   
We say that $f(x) = O( g(x) )$ (pronounced ``Big o''), as $x \rightarrow x_0$,
if {\bf there exists a positive constant} $C$ such that
\begin{equation}
|f(x)| \leq C |g(x)|,
\end{equation}
for all $x$ that are sufficiently close to $x_0$.   Here we allow $x_0$ to be any real vector, 
and we also allow the very common case $x_0 = \pm \infty$. 
Equivalently, we may  say  
\[
\left|   \frac{f(x)}{g(x)}  \right| \text{ is eventually bounded above}. 
\]

\begin{example}
\rm{
We write
$e^x = 1 + x + \frac{1}{2} x^2 + O(x^3)$,
as $x \rightarrow 0$.     We could, of course, also write $e^x -(1 + x + \frac{1}{2} x^2 ) = O(x^3)$, though the
former way of writing it is much more common. 
In this case, we can give a `better' Big-O estimate by adding more terms of the Taylor series:
$e^x = 1 + x + \frac{1}{2} x^2 + \frac{1}{6} x^3 + O(x^4)$,
as $x \rightarrow 0$.
\hfill $\square$
}  
\end{example}


\begin{example}
\rm{
Given $f(x) := x \sin \left(\frac{1}{x}  
\right)$, and $g(x) := x^2 - 12$,  we have
\[
f(x) = O( g(x) ), \text{ as } x \rightarrow \infty.
\]
In other words, for all sufficiently large $x$,  $|f(x)| \leq C g(x)$, despite the fact that
this statement is false for these particular functions, for some small positive values of $x$. 
\hfill $\square$
}  
\end{example}
Claim. \ Big-O enjoys transitivity: 
\[
\text{ If } f = O(g), \text{ and }  g = O(h), \text{ then } f= O(h). 
\]
\begin{proof}
For all $x$ sufficiently close to $x_0$, there exists positive constants $C_1, C_2$ such that 
$|f(x)| \leq C_1 |g(x)|$ and $|g(x)| \leq C_2 |h(x)|$, implying that
\[
|f(x)| \leq C_1 |g(x)| \leq C_1 C_2  |h(x)|. 
\]
\end{proof}


There is another very useful comparison technique, for any two given functions $f, g:\R^d \rightarrow \C$. 
We say that $f(x) = o( g(x) )$ (pronounced ``little o''), as $x \rightarrow x_0$,
if {\bf for all positive constants} $C$, we have:
\begin{equation}
|f(x)| \leq C | g(x) |,
\end{equation}
for all $x$ that are sufficiently close to $x_0$.   Again 
 we allow $x_0$ to be any real vector, and we also allow the very common case $x_0 = \pm \infty$. 
Equivalently, we may also write
\[
\lim_{x \rightarrow x_0}   \left|   \frac{    f(x) }{   g(x)  }  \right| = 0,
\]
which intuitively means that $g$ approaches $x_0$ faster than $f$ does. 

\begin{example}
\rm{
Given $f(x) :=  \sqrt x$,   and $g(x) :=x$, where we restrict the domain of 
both functions to be $(0, +\infty)$.
We claim
$f(x) = o( g(x) ), \text{ as } x \rightarrow 0$.
\begin{proof}
\[
\lim_{x \rightarrow 0}  \left|  \frac{    f(x) }{   g(x)  }  \right| = 
\lim_{x \rightarrow 0}   \left|  \frac{ \sqrt x }{ x }  \right| =
\lim_{x \rightarrow 0}   \left|  \frac{ 1}{ \sqrt x }  \right| 
= 0.
\]
So $g$ approaches $0$ much faster than $f$.
\end{proof}
}  
\end{example}   

\medskip
Claim.  Little-o also enjoys transitivity:
\[
\text{ If } f = o(g), \text{ and }  g = o(h), \text{ then } f= o(h). 
\]
\begin{proof}
The two given limits
$\lim_{x \rightarrow x_0}   \left|   \frac{    f(x) }{   g(x)  }  \right| = 0$ and 
$\lim_{x \rightarrow x_0}   \left|   \frac{    g(x) }{   h(x)  }  \right| = 0$ together imply   that
\[
\lim_{x \rightarrow x_0}   \left|   \frac{    f(x) }{   h(x)  }  \right| 
=\lim_{x \rightarrow x_0}   \left|   \frac{    f(x) }{   g(x)  }  \right|   \left|   \frac{    g(x) }{   h(x)  }  \right|  
= 0.
\]
\end{proof}


\bigskip
\chapter{Various forms of convergence}
 \label{Appendix B}
\section{Weierstrass M-test}

How can we quickly conclude that certain series converge uniformly? The following criterion, discovered by Karl Weierstrass, comes to the rescue. 
For the proofs of these  basic real analysis results, see for example the classic \cite{RudinGreenBook}.

\begin{thm} \rm{[Weierstrass M-test]}
Suppose that $f_n(x)$ is a sequence of complex-valued functions defined on a set $E\subset \R$, 
such that there exists a sequence of numbers $M_n\geq 0$ satisfying the following conditions:
\begin{enumerate}
\item   $ |f_{n}(x)|   \leq M_{n},    \  \forall n \in \Z \text{ and  all } x \in E$.
\item    $\sum_{n \in \Z} M_n < \infty$.
\end{enumerate}
Then the series
$\sum _{n\in \Z}  f_n(x) $
converges absolutely and uniformly on $E$.
\hfill $\square$
\end{thm}

In practice, the Weierstrass $M$-test gets used together with the following test, which allows us to partially answer the question:
\begin{question}
When does a series $\sum_{n\in \Z} f_n(x)$ converge to a continuous function of $x$? 
\end{question}

\begin{thm}  \label{uniform limit test}      \rm{[Uniform limit]}
Suppose that $s_n(x):E  \rightarrow \C$ is a sequence of continuous functions defined on a set $E\subset \R$, 
and that $s_n$ converges uniformly to $s(x)$, on $E$.   Then $s(x)$ is continuous on $E$. 
\hfill $\square$
\end{thm}


\bigskip
\section{Some things you wanted to know about convergence but were afraid to ask}

It's often useful to pass from $L^2$ convergence to pointwise convergence, under some additional hypothesis on $f$.    Throughout, we fix a real number $1\leq p < \infty$.   Given a measurable subset $E\subset \R^d$, and a sequence of functions $f_n:E \rightarrow \C$, we say that
$f_n(x) \rightarrow f(x)$ in the $p$-norm  if
 \begin{equation} \label{convergence in L^p norm, take 2}
\lim_{n \rightarrow \infty} 
\int_{E} \left|  f_n(x) - f(x) \right |^p dx =0,
\end{equation}
for which we will also use here the notation $\lim_{n\rightarrow \infty} \| f_n - f \|_{L^p(E)} = 0$.
Sometimes, if the constant $p$ is not specified, it is common to simply call \eqref{convergence in L^p norm, take 2}
{\bf convergence in norm}.
The two most common subsets are $E:= \R^d$, and $E:= [0, 1]^d$.
A natural question arises:
\begin{question}
When can we pass from convergence in norm to pointwise convergence?
\end{question}

Given a series $\sum_{n \in \Z^d} f_n(x)$, we consider the sequence of partial sums
$S_N(x):= \sum\limits_{|n| < N} f_n(x)$.    
By definition, we say the series converges
\begin{enumerate}[(a)]
\item {\bf pointwise on $E$}  if the sequence $\{S_N(x)\}_{N=1}^\infty$ converges, for each $x \in E$.
\item  {\bf absolutely on $E$}  if the series $\sum_{n \in \Z^d} |f_n(x)|$ converges pointwise, 
for each $x \in E$.
\item {\bf uniformly on $E$}  if the sequence of partial sums $S_N(x)$ converge uniformly on $E$.
\item {\bf in the $p$-norm} on $E$ if $\lim_{n\rightarrow \infty} \| f_n - f \|_{L^p(E)} = 0$.\end{enumerate}

\medskip
\begin{lem}  \label{technical equality a.e.}
Consider the partial sums 
$
S_N(x):= \sum\limits_{|n| < N \atop n \in \Z^d} f_n(x), 
$
for all $x$ in a given subset $E\subset \R^d$.  Suppose we have the following two properties:
\begin{enumerate}
\item There exists a function $f:\R^d \rightarrow \C$ such that $S_N(x) \rightarrow f(x)$ in the $p$-norm, on  $E$.    
\item $S_N(x)$ converges uniformly to the series $S(x):= \sum_{n\in \Z^d} f_n(x)$ on $E$.
\end{enumerate}
Then $ S(x) = f(x)  \text{ a.e. on } E$.
 \hfill $\square$
\end{lem}

\begin{lem}   \label{absolute convergence of Fourier series implies continuity}
Let $f\in L^1([-c,c])$, and suppose we already know that its Fourier series  converges pointwise:
\begin{equation} \label{Fourier series converging pointwise}
f(x) =          \sum_{n\in \Z} \hat f_n   e^{\frac{2\pi i n x}{c}},  
\end{equation}
 If the series \eqref{Fourier series converging pointwise}
converges absolutely, then $f$ is a continuous function on $[-c,c]$, and $f(-c) = f(c)$.
\end{lem}
\begin{proof}
The idea is to use the uniform limit Theorem \ref{uniform limit test}, together with the fact that the summands $\hat f_n   e^{\frac{2\pi i n x}{c}}$ are continuous functions of $x$.   So it remains to show that the convergence of the series is uniform.  
\[
|S_N(x)| := 
\left|  \sum_{|n|<N}  \hat f_n   e^{\frac{2\pi i n x}{c}}  \right|
\leq  \sum_{|n|<N} \left| \hat f_n   e^{\frac{2\pi i n x}{c}}  \right|
= \sum_{|n|<N} | \hat f_n |  < \infty,
\]
where the penultimate equality holds because $| e^{\frac{2\pi i n x}{c}}  |=1$, and  the last inequality holds by assumption. 
So by the $M$-test, with $M_n := | \hat f_n |$, we have uniform convergence of the series.   
 Finally, the claim $f(-c) = f(c)$ is trivial, because $f(\pm c) :=  \sum_{n\in \Z} \hat f_n   e^{\pm 2\pi i n} = \sum_{n\in \Z} \hat f_n$.
\end{proof}

In the previous lemma, we could have also used the alternate notation of the circle $ \R/c\Z$, and rewrite everything 
in terms of it, which automatically incorporates periodicity. 
The following {\bf passage from convergence in the $ L^2([-c,c])$ norm, to pointwise convergence}, is often useful. 
\begin{lem}   \label{norm convergence plus absolute convergence implies equality}
Let $f\in L^2([-c,c])$ be a continuous function, and write its Fourier series as 
\begin{equation} \label{L^2 convergent series}
f(x)  \underset{L^2([-c,c])}{=}
\sum_{n\in \Z} \hat f_n   e^{\frac{2\pi i n x}{c}},  
\end{equation}
which by definition means that this series converges in the $L^2([-c,c])$-norm.   

If the series \eqref{L^2 convergent series}
converges absolutely, then it also converges pointwise and uniformly to $f(x)$, for all $x \in [-c,c]$.
\end{lem}
\begin{proof}
Repeating the computation of the previous proof, we have:
\[
|S_N(x)| := 
\left|  \sum_{|n|<N}  \hat f_n   e^{\frac{2\pi i n x}{c}}  \right|
\leq  \sum_{|n|<N} \left| \hat f_n   e^{\frac{2\pi i n x}{c}}  \right|
\leq \sum_{|n|<N} | \hat f_n |  < \infty,
\]
 Therefore by the $M$-test again, $S_N(x)$ converges
 uniformly to the series $S(x)~:=~\sum\limits_{n\in \Z} \hat f_n   e^{\frac{2\pi i n x}{c}}$,  for all  $x \in [-c, c]$.
 We also know, by Lemma \ref{absolute convergence of Fourier series implies continuity}, 
  that $S(x)$ is continuous on $[-c, c]$. 
We still need to prove that the series converges to $f$, but now we at least know that both hypotheses of Lemma \ref{technical equality a.e.} are satisfied (with $p=2$ and $E:= [-c, c]$), 
and therefore  $S(x) = f(x)  \text{ a.e. on } E$.

To prove that $S(x) = f(x)$ for all $x \in [-c, c]$, we observe that the summands 
$\hat f_n   e^{\frac{2\pi i n x}{c}}$ are continuous functions of $x$, and hence by the uniform limit theorem (Theorem \ref{uniform limit test}), the series $S(x)$ is itself a continuous function of $x$.   Since $f$ is also continuous on $[-c, c]$, and $S(x) = f(x)$ almost everywhere, they must agree everywhere.  
\end{proof}

\section{Bump functions}\label{Appendix:bump functions and inner products}

Perhaps the easiest bump function to define is the function (\cite{SteinShakarchi}, page 209):
\begin{equation*}
    \varphi(x):=
    \begin{cases}
        c\, e^{-\frac{1}{1-\|x\|^2}} & \text{ if } \|x\|<1, \\
        0 & \text{ if } \|x\|   \geq1. \\
    \end{cases}
\end{equation*}
where the constant $c$ is chosen so that
$\int_{\R^d} \varphi(x) dx = 1$.  By definition, $\varphi$ is compactly supported, on the unit ball. 
   It turns out that $\varphi$ is infinitely smooth.  
As usual, using $\varphi$ we can build a family of integrable functions:
\begin{equation*}
    \varphi_{\varepsilon}(x):=\varepsilon^{-d}\varphi(x\varepsilon^{-1}), \text{ for all } 0<\varepsilon\leq1.
\end{equation*}
Thus, the family $\{\varphi_\varepsilon\}$ is an approximate identity.

More generally, a {\bf bump function} is defined to be any infinitely smooth function $\varphi:\R^d\rightarrow \C$ that is compactly supported.  
By Lemma \ref{useful Schwartz fact}, we know that any such bump function
$\varphi$ lies in the Schwartz class $S(\R^d)$.   Clearly finite linear combinations of bump functions are again bump functions, making the space of bump functions  a vector subspace of the space
of Schwartz functions.


Often, a slightly more general sort of space than a Hilbert space is required.   
 Suppose there exists 
a function called $\langle \cdot, \cdot \rangle$, defined from $V \times V\rightarrow \R$, that enjoys the following
properties:
\begin{enumerate}[(a)]
\item  (Strict positivity)      $\langle x, x \rangle >0$, for all nonzero $x \in V$.   
\item  (Symmetry)        $\langle x, y \rangle =   \langle y, x \rangle$, for all 
$x, y \in V$.
\item   (Linearity)    for any fixed $a \in V$,   the function $x \rightarrow \langle x, a \rangle$ is linear, which means that
\[
 \langle x+y, a \rangle = \langle x, a \rangle + \langle y, a \rangle,
 \]
  for all $x, y \in V$.
\end{enumerate}
Then $\langle \cdot, \cdot \rangle$ is called an {\bf inner product on $V$}, and 
$V$ is an {\bf inner product space}  (over $\R$).  
 Inner products also 
interact well with continuity, in the following precise sense (\cite{EinsiedlerWardBook}, p. 74):
\begin{lem}
If $x_n \rightarrow x$ in $V$, and $y_n \rightarrow y$ in $V$, then $\langle x_n, y_n \rangle \rightarrow \langle x, y \rangle$.
\end{lem}
We say that an inner product space $V$ is {\bf complete}, relative to the distance function 
$d(x_n, x_m):= 
\langle x_n-x_m, x_n-x_m \rangle^{\tfrac{1}{2}}$,
if every Cauchy sequence $\{x_n\}_{n = 1}^\infty$ in $V$ converges to a point of $V$.
Finally, we mention a basic fact about linear functions acting on complete  inner product space.
\begin{thm}[Riesz Representation Theorem]
Let $V$ be a complete inner product space (finite or infinite dimensional), and suppose that $f: V \rightarrow \R$ is a continuous linear functional on $V$.
Then there exists a unique $w\in V$ such that 
\[
f(x) = \langle x, w \rangle,
\]
for all $x \in V$.
\end{thm}


\chapter{Solutions and hints}

\begin{quote}    
 ``There are no problems - just pauses between ideas.''
 
 -- David Morrell, Brotherhood of the Rose  
 \end{quote}

{\bf \Large Chapter \ref{Chapter.Tiling.A.Rectangle}}
\bigskip
\index{tiling}

\bigskip
 Exercise \ref{TrivialExponential} 
 \quad
  By Euler, we have $1 = e^{i \theta} = \cos\theta  + i\sin\theta$, which holds
  if and only if
 $\cos\theta = 1$, and $\sin\theta = 0$.  The latter two conditions hold simultaneously if and only if 
 $\theta \in 2 \pi k$, with 
 $k \in \Z$.

 \medskip
Exercise \ref{bound of the exponential function}
\quad   Let $z := a + bi$, so that $|e^z| = |e^{a+bi}|  =|e^a| | e^{bi} | = e^a \cdot 1
\leq e^{\sqrt{ a^2 + b^2}} = e^{|z|}$, using the fact that $| e^{bi} |= 1$ for all real numbers $b$.

 \bigskip
Exercise \ref{orthogonality for exponentials} 
\quad
  In case $a \not= b$, we have
\begin{equation*}
\int_0^1 e_a(x)  \overline{e_b(x)} dx =  \int_0^1 e^{2\pi i (a-b) x} dx 
=  \frac{e^{2\pi i (a-b)}}{2\pi i(a-b)} - 1=0,
\end{equation*}
because we know that $a-b \in \Z$. In case $a = b$, we have
\begin{equation*}
\int_0^1 e_a(x)  \overline{e_a(x)} dx =  \int_0^1 dx =  1.
\end{equation*}

  \bigskip
Exercise \ref{definition of complex integral} \ 
\quad  By definition, 
\begin{align*}
\int_{[0,1]}   e^{-2\pi i \xi x}  dx  &:=   \int_{[0,1]}  \cos(2\pi \xi x) dx + i \int_{[0,1]}  \sin(2\pi \xi x) dx\\
&= \frac{\sin(2\pi \xi)}{2\pi\xi} +i \frac{-\cos(2\pi \xi)+1}{2\pi\xi}\\
&= \frac{i\sin(2\pi \xi)}{2\pi i \xi} + \frac{\cos(2\pi \xi)-1}{2\pi i \xi}\\
&= \frac{e^{2\pi i \xi}-1}{2\pi i \xi}.
\end{align*}

\bigskip
 Exercise \ref{SumOfRootsOfUnity} \ 
\quad   Let $S:= \sum_{k = 0}^{N-1}    e^{\frac{2\pi i k}{N}}$, and note that we may write
 \[
 S = \sum_{k\text{ mod } N}     e^{\frac{2\pi i k}{N}}.
 \]
   Now, pick any $m$ such that $e^{\frac{2\pi i m}{N}} \not=1$.   Consider
 \begin{align*}
 e^{\frac{2\pi i m}{N}} S &= \sum_{k\text{ mod } N}    e^{\frac{2\pi i (k + m)}{N}} \\
 &= \sum_{n\text{ mod } N}    e^{\frac{2\pi i n}{N}} =  S,
 \end{align*}
 so that $0=( e^{\frac{2\pi i m}{N}} -1)S$, and since
 by assumption $e^{\frac{2\pi i m}{N}} \not=1$, we have $S = 0$.
 
 \bigskip
Exercise \ref{DivisibilityUsingExponentials} 
\quad
    We use the finite geometric series:  
$1+x+ x^2 + \cdots + x^{N-1} = \frac{x^{N} - 1}{x-1}$.  Now, if $N \not{\mid}  M$, then $x:=  e^{\frac{2\pi i M}{N}} \not=1$, so we may substitute this value of $x$ into the finite geometric series to get:
\begin{align*}
\frac{1}{N}  \sum_{k = 0}^{N-1}    e^{\frac{2\pi i kM}{N}} &=  
\frac{ e^{\frac{2\pi i MN}{N}}- 1}{ e^{\frac{2\pi i M}{N}}-1} \\
&= \frac{0}{e^{\frac{2\pi i M}{N}}-1}=0.
\end{align*}
On the other hand, if $N \mid M$, then $ \frac{1}{N}  \sum_{k = 0}^{N-1}    e^{\frac{2\pi i kM}{N}} =
\frac{1}{N}  \sum_{k = 0}^{N-1}  1 = 1$.

 \bigskip
Exercise \ref{Orthogonality.for.roots.of.unity} \ 
\begin{align*} 
\frac{1}{N} \sum_{k = 0}^{N-1}    e^{\frac{2\pi i ka}{N}} e^{-\frac{2\pi i kb}{N}} 
&=   \frac{1}{N} \sum_{k = 0}^{N-1}    e^{\frac{2\pi i k(a-b)}{N}}.  \\
\end{align*}
Therefore, using Exercise \ref{DivisibilityUsingExponentials}, we see that the latter sum equals
$1$ exactly when $N \mid a-b$, and vanishes otherwise.

\bigskip
Exercise \ref{trick-write an integer as a product with roots of unity} \ 
\quad   We begin with the factorization of the polynomial
 $x^n-1= \prod_{k=1}^n (x - \zeta^k)$,  with $\zeta:= e^{2\pi i /n}$.  Dividing both sides by 
 $x-1$, we obtain $1+ x + x^2 + \cdots + x^{n-1} =  \prod_{k=1}^{n-1} (x - \zeta^k)$.  Now
 substituting $x=1$, we have $n =  \prod_{k=1}^{n-1} (1 - \zeta^k)$.

\bigskip
Exercise \ref{PrimitiveRootsOfUnity}  \ 
\quad   Suppose to the contrary, that a primitive  $N$'th root of unity is of the form $e^{2\pi i m/N}$, where 
$\gcd(m,N) > 1$.   Let $m_1 := \frac{m}{\gcd(m, N)}$, and $k:=\frac{N}{\gcd(m, N)}$, so that by assumption both
$m_1$ and $k$ are integers.     Thus   $e^{2\pi i m/N} = e^{2\pi i m_1/k}$, a $k$'th root of unity, with $k < N$, a contradiction.

\bigskip
Exercise \ref{zeros of the sin function} \
\quad   We recall Euler's identity:  
\[
e^{iw} = \cos w + i \sin w,
\]
 which is valid for all $w \in \mathbb C$.   Using Euler's identity first with $w:= \pi z$,  and then with $w := -\pi z$, we have the two identities  $e^{\pi i z} =  \cos \pi z + i \sin \pi  z$, and $e^{-\pi i z} =  \cos \pi z - i \sin \pi z$.    Subtracting the second identity from the first, we have 
\[
\sin(\pi z) = \frac{1}{2i}\left(  e^{\pi i z} - e^{-\pi i z}  \right).
\]
Now it's clear that $\sin(\pi z) = 0  \iff e^{\pi i z} = e^{-\pi i z}  \iff e^{2\pi i z}  = 1 \iff  z \in \Z$, by Exercise~\ref{TrivialExponential}.

\bigskip \medskip
Exercise \ref{Erdos lattice partition problem} \
\quad   We will assume, to the contrary, 
 that we only have one arithmetic progression with a common difference of $a_N$, the largest of the common differences. 
 We hope to obtain a contradiction.  To each arithmetic progression $\{ a_k n + b_k \mid n \in \Z\}$, we associate the generating function
\[
f_k(q):= \sum_{a_k n + b_k \geq 0, \ n\in \Z} q^{a_k n + b_k },
\]
where $|q| < 1$, in order to make the series converge.  The hypothesis that we have a tiling
\index{tiling}
 of the integers by these $N$ 
arithmetic progressions translates directly into an identity among these generating functions:
\[
 \sum_{a_1 n + b_1 \geq 0, \ n\in \Z} q^{a_1 n + b_1 } + \cdots +
  \sum_{a_N n + b_N \geq 0, \ n\in \Z} q^{a_N n + b_N } =  \sum_{n=0}^\infty q^n.
\]
Next, we use the fact that we may rewrite each generating function in a `closed form' of the following kind, because they are geometric series:
$f_k(q):= \sum_{a_k n + b_k \geq 0, \ n\in \Z} q^{a_k n + b_k } = \frac{q^{b_k}}{1-q^{a_k}}$. Thus, we have:
\[
 \frac{q^{b_1}}{1-q^{a_1}} +     \cdots +       \frac{q^{b_N}}{1-q^{a_N}} =  \frac{1}{1-q}.
\]
Now we make a `pole-analysis' by observing that each rational function $f_k(q)$ has  poles at precisely all of the $k$'th roots of unity.
The final idea is that the `deepest' pole, namely $e^{   \frac{2\pi i}{N}  }$, cannot cancel with any of the other poles.  To make this idea precise, 
we  isolate the only rational function that has this pole (by assumption): 
\[
 \frac{q^{b_N}}{1-q^{a_N}} = \frac{1}{1-q} - \left(  \frac{q^{b_1}}{1-q^{a_1}}  + \cdots +  \frac{q^{b_{N-1}}}{1-q^{a_{N-1}}}  \right).
\]
Finally, we let  $q\rightarrow e^{   \frac{2\pi i}{N}  }$, to get a finite number on the right-hand-side, 
and infinity on the left-hand-side of the latter identity, a contradiction.


 \bigskip \bigskip 
{\bf \Large Chapter \ref{Chapter.Examples}}

 \medskip
Exercise \ref{transform.of.interval.a.to.b} \ 
 \quad   If $\xi = 0$, we have $\hat 1_{[a,b]}(0) := \int_a^b e^0 dx = b-a$.  
If $\xi \not=0$,  we can compute the integral:
\begin{align*}
\hat 1_{[a,b]}(\xi) &:= \int_a^b e^{-2\pi i \xi x} dx \\
&=\frac{e^{-2\pi i \xi b}    -  e^{-2\pi i \xi a}    }{-2\pi i \xi}.
\end{align*}

\medskip
Exercise \ref{transform.of.unit.cube}  \ 
\quad  Beginning with the definition of the Fourier transform of the unit cube $[0,1]^d$, we have:
\begin{align*}
\hat 1_{\square}(\xi) &=  \int_{\square}  e^{2\pi i \langle x, \xi \rangle}dx \\
&= \int_0^1 e^{2\pi i \xi_1 x_1} dx_1   \int_0^1  e^{2\pi i \xi_2 x_2} dx_2  \cdots   
         \int_0^1  e^{2\pi i \xi_d x_d} dx_d \\
&= \frac{1}{(-2\pi i)^d}  \prod_{k=1}^d \frac{    e^{-2\pi i \xi_k} -1    }{   \xi_k         },
\end{align*}
valid for all $\xi \in \R^d$, except for the finite union of hyperplanes defined by \\
$H := \{  x \in \R^d \mid \xi_1 = 0 \text{ or } \xi_2 = 0  \dots  \text{ or }  \xi_d = 0 \}$.

\medskip
Exercise \ref{brute force Bernoulli polys}  \ 
\quad   To see that the generating-function definition of the Bernoulli polynomials in fact gives polynomials, we first write the Taylor series of the following two analytic functions:
\[
\frac{t}{e^t - 1} = \sum_{k=0}^\infty  \frac{B_k}{k!} t^k
\]
\[
e^{xt} = \sum_{j=0}^\infty  \frac{ x^j t^j}{j!}.
\]
Multiplying these series together by brute-force gives us:
\begin{align}
\frac{t}{e^t - 1} e^{xt}    &= \left( \sum_{k=0}^\infty  \frac{B_k}{k!}    t^k  \right)  
\left( \sum_{j=0}^\infty  \frac{ x^j}{j!}  t^j  \right) \\
&= \sum_{n=0}^\infty \left(     \sum_{j+k = n}  \frac{B_k}{k!}     \frac{ x^j}{j!}    \right) t^n \\
&=  \sum_{n=0}^\infty \left(     \sum_{k = 0}^n  \frac{B_k}{k!}     \frac{ x^{n-k}}{(n-k)!}    \right) t^n.
\end{align}
The coefficient of $t^n$ on the LHS is by definition $\frac{1}{n!} B_n(x)$, and by uniqueness of Taylor series, this must also be the coefficient on the RHS, which is seen here to be a polynomial in $x$.   In fact, we see more, namely that 
\[
\frac{1}{n!} B_n(x) =  \sum_{k = 0}^n  \frac{B_k}{k!}     \frac{ x^{n-k}}{(n-k)!}, 
\]
which can be written more cleanly as
$ B_n(x) = \sum_{k = 0}^n  {n\choose k}    B_k x^{n-k}$.

\medskip
Exercise  \ref{Reflection property for B_n(x)} \ 
\quad  Commencing with the generating-function definition of the Bernoulli polynomials, 
equation
\ref{generating function for Bernoulli polynomials}, 
we replace $x$ with $1-x$ in order to observe the coefficients $B_k(1-x)$:
\begin{align*}
\sum_{k=0}^\infty  \frac{B_k(1-x)}{k!} t^k  &=   \frac{te^{t(1-x)}}{e^t - 1} \\
&= \frac{te^t e^{-tx}}{e^t - 1} \\
&= \frac{t e^{-tx}}{1-  e^{-t}  } \\
&=  \frac{-t e^{-tx}}{e^{-t}-1  } \\
&=  \sum_{k=0}^\infty  \frac{B_k(x)}{k!} (-t)^k,
\end{align*}
where the last equality follows from the definition of the same generating function, namely equation
\ref{generating function for Bernoulli polynomials}, but with the variable $t$ replaced by $-t$. 
Comparing the coefficient of $t^k$ on both sides, we have $B_k(1-x) = (-1)^k B_k(x)$.

\medskip
 Exercise \ref{difference of Bernoulli polys}  \ 
 To show that
$B_n(x+1) - B_n(x) = n x^{n-1}$, we play with:
\begin{align*}
\sum_{k=0}^\infty  \left(\frac{B_k(x+1)}{k!} t^k - \frac{B_k(x)}{k!} t^k  \right)  &=   
                        \frac{te^{t(x+1)}}{e^t - 1}   -    \frac{te^{t(x)}}{e^t - 1}  \\
                   &=      e^t \frac{te^{tx}}{e^t - 1}   -    \frac{te^{t(x)}}{e^t - 1}  \\
                   &=   (e^t - 1) \frac{te^{tx}}{e^t - 1}    \\
                   &=   te^{tx}  \\
               &=    \sum_{k=0}^\infty   \frac{x^k}{k!}  t^{k+1}   \\
               &=    \sum_{k=1}^\infty   \frac{x^{k-1}}{(k-1)!}  t^{k}   \\
               &=    \sum_{k=1}^\infty   \frac{k x^{k-1}}{k!}                   t^{k}.
\end{align*}
Therefore, again comparing the coefficients of $t^k$ on both sides, we arrive at the required identity.

\medskip
Exercise \ref{derivative of Bernoulli polys}   \ 
\quad We need to show that  $\frac{d}{dx}  B_n(x) = n B_{n-1}(x)$.  Well,
\begin{align*}
\sum_{k=0}^\infty  \frac{d}{dx} \frac{B_k(x)}{k!} t^k 
                   &=     \frac{d}{dx} \frac{te^{tx}}{e^t - 1}   \\
              &=   t \sum_{k=0}^\infty   \frac{B_k(x)}{k!} t^k      \\
          &=   \sum_{k=0}^\infty   \frac{B_k(x)}{k!} t^{k+1}      \\
          &=   \sum_{k=1}^\infty   \frac{B_{k-1}(x)}{(k-1)!} t^{k}      \\
                    &=   \sum_{k=1}^\infty  k \frac{B_{k-1}(x)}{k!} t^{k},
\end{align*}
so that comparing the coefficient of $t^n$ on both sides, the proof is complete.

\medskip
Exercise \ref{Elementary bounds for sin(x), sinc(x)}
We'll prove part \ref{Elementary trig bounds, part b}.   To begin, we have:
\begin{align*}
\big| e^{i\theta} - 1 \big|^2 &= \big| \cos \theta -1 + i \sin \theta  \big|^2 = (\cos \theta -1)^2  + \sin^2\theta \\
&= 2 - 2 \cos \theta  = 4 \sin^2\left(\frac{\theta}{2}\right).
\end{align*}
So it suffices to show that
 $4 \sin^2\left(\frac{\theta}{2}\right) \leq  \theta^2$, for all $0 \leq \theta \leq 2\pi$.
 In other words, the problem is reduced to the Calculus I problem of showing that  
 $ \sin\left(\frac{\theta}{2}\right)      \leq      \frac{\theta}{2}$, for $\theta \in [0, 2\pi]$.
To prove this, we let $y(x)=x- \sin x $, so that it suffices to prove that $y \geq 0$ on $[0, \pi]$.
Computing its derivative, 
 $y'(x) = 1-\cos x  \geq 0$ on $[0, \pi]$, and since $y(0) = 0$, we conclude that $y$ is an increasing function.  This proves $y \geq 0$ on $[0, \pi]$.



\medskip
Exercise \ref{polytope with 5 vertices}
For part (a), suppose to the contrary that none of the vertices of $\P$ have degree $4$.  Because each of its vertices must have degree at least $3$, then all of vertices have degree $3$.  By the `handshanking lemma' of elementary graph theory, we have:
\[
2 |E| =  \sum_{\text{vertices }v\in\P} \text{deg}(v) = 3\cdot 5,
\]
a contradiction.
To prove part (b), consider the unit cube in $\R^3$, and take $4$ vertices that belong to one facet, with one vertex from an oppostive facet.  It's clear that all of its vertices have degree $4$.

\medskip
Exercise \ref{Dirichlet's convergence test}
\quad  Considering the partial sum $S_n:= \sum_{k=1}^n a_k b_k$, we know by Abel summation that
\[
S_n = a_n B_n + \sum_{k=1}^{n-1} B_k(a_k - a_{k+1}),
\]
for each $n \geq 2$, where $B_n :=  \sum_{k=1}^n b_k$. 
 By assumption, $|B_n|:= | \sum_{k=1}^n b_k |  \leq M$, and the $a_k$'s are going to $0$, 
 so we see that  the first part of the right-hand-side approaches 
 zero, namely: 
 $|a_n B_n| :=   |a_n| | \sum_{k+1}^n b_k|   \rightarrow 0$, as $n \rightarrow \infty$. 
 
Next, we have
 \[
 |  \sum_{k=1}^{n-1} B_k(a_k - a_{k+1}) |    \leq \sum_{k=1}^{n-1} | B_k| |  a_k - a_{k+1} |
 \leq    M  \sum_{k=1}^{n-1} |  a_k - a_{k+1} | =  M  \sum_{k=1}^{n-1}   (a_k - a_{k+1}), 
 \]
 where the last equality holds because by assumption the $a_k$'s are decreasing.  But the last finite
sum equals  $ -M a_{n} + M  a_1$, and we have $ \lim_{n\rightarrow \infty} (-M a_{n} +M  a_1)=  M a_1$, a finite limit.
Therefore    $\sum_{k=1}^{n-1} B_k(a_k - a_{k+1})$ converges absolutely, and so $S_n$ converges, as desired.

\medskip
Exercise \ref{exponential sum bound}
\quad  We fix  $x \in \R-\Z$, and let $z:= e^{2\pi i x}$, which lies on the unit circle, and by assumption $z \not= 1$. 
 Then 
\begin{equation}
 \left |  \sum_{k= 1}^n     e^{2\pi i k x} \right |   
= \left | \sum_{k= 1}^n   z^k  \right |    = \left |  \frac{z^{n+1} -1}{z-1}  \right |  \leq \frac{2}{z-1},
\end{equation}
because $|z^{n+1} -1| \leq |z^{n+1}|  +1 = 2$.  We also have
\begin{align*}
 |z-1|^2 &= |e^{2\pi ix}-1||e^{-2\pi ix}-1| = |2-2\cos(2\pi x)|  =  4\sin^2(\pi x),
\end{align*}
so that we have the equality
$\left | \frac{2}{z-1} \right | =\left |    \frac{1}{\sin(\pi x)}  \right |$.
Altogether, we see that 
\begin{equation}
 \left |  \sum_{k= 1}^n     e^{2\pi i k x} \right |  \leq     \frac{1}{  |  \sin(\pi x)  |  }.
\end{equation}

\medskip
Exercise \ref{rigorous convergence of P_1(x)}
\quad We fix $a \in \R - \Z$ and need to prove that $\sum_{m = 1}^\infty  \frac{e^{2\pi i m a}}{m}$ converges.
Abel's  summation formula \eqref{actual Abel summation}   gives us
\[
\sum_{k = 1}^n  \frac{e^{2\pi i k a}}{k} = \frac{1}{n}\sum_{r=1}^n e^{2\pi i r a} 
+  \sum_{k=1}^{n-1}      \Big(    \sum_{r=1}^k   e^{2\pi i r a}   \Big)     \frac{1}{k(k+1)},
\]
so that 
\[
\sum_{k = 1}^\infty    \frac{e^{2\pi i k a}}{k} = 
 \sum_{k=1}^{\infty}      \Big(    \sum_{r=1}^k   e^{2\pi i r a}   \Big)     \frac{1}{k(k+1)}.
\]
and the latter series in fact converges absolutely.

 \bigskip \bigskip 
{\bf \Large Chapter \ref{Fourier analysis basics}}

\medskip
Exercise   \ref{elementary norm relations}
\quad  For all four inequalities, we will use an arbitrary vector $a \in \R^d$. 
 For the first inequality, $a_1^2 + \cdots + a_d^2 \geq \max\{ |a_1|, \dots, |a_d| \}^2 :=  \|a\|_\infty^2$. 

The second inequality $\|a\|_2 \leq \|a\|_1$   means that $\sqrt{a_1^2 + \cdots + a_d^2} \leq |a_1| + \cdots +  |a_d|$, which is clear by squaring both sides.

To prove the third and most interesting inequality here, we use the Cauchy-Schwarz inequality, with the two vectors
$x:= (a_1, \dots, a_d)$ and $(1, 1, \dots, 1)$:
\[
\| a \|_1:=   |a_1| \cdot 1 +   \cdots + |a_d|  \cdot 1     \leq     \sqrt{ a_1^2 + \cdots + a_d^2}   \sqrt{  1 + \cdots + 1 }
    =  \sqrt{d} \  \| a \|_2,
 \]
which also shows that we obtain equality if and only if $(a_1, \dots, a_d)$ is a
scalar multiple of $(1, 1, \dots, 1)$.

For the fourth inequality, we have:
\[
\sqrt{a_1^2 + \cdots + a_d^2} \leq 
 \sqrt{  d  \max\{ |a_1|, \dots, |a_d| \}^2   }:=
  \sqrt{d}    \|a\|_\infty.
\]

\medskip
Exercise \ref{exercise:hyperbolic cosine and sine} 
\quad  
To prove part (a), we compute:  
\begin{align*}
 \left(    \frac{e^{t} + e^{-t}}{2}    \right)^2      -    \left(    \frac{e^{t} - e^{-t}}{2}    \right)^2 
 &= \frac{e^{2t} + 2 + e^{-2t} - \left(    e^{2t}  -   2 +   e^{-2t}     \right)  }{4} = 1.
\end{align*}
To prove part (b), we begin with the definition of the hyperbolic cotangent:
\begin{align*}
 t \coth t &=  t\frac{ e^t + e^{-t}}{e^t - e^{-t}} = t \frac{ e^t }{e^t - e^{-t}} + t \frac{ e^{-t}}{e^t - e^{-t}} \\
  &=  \frac{ t }{1 - e^{-2t}} +  \frac{ t }{  e^{2t}-1}. 
\end{align*}
 Recalling the definition of the Bernoulli numbers, namely
$
        \frac{t}{e^t-1}   =  \sum_{k =0}^\infty   B_k \frac{t^k}{k!},
$
we see that 
\begin{align*}
 t \coth t &=    \frac{1}{2} \left( \frac{ -2t }{ e^{-2t}-1} \right) + \frac{1}{2} \left( \frac{ 2t }{  e^{2t}-1}  \right) \\
 &=    \frac{1}{2}  \sum_{k =0}^\infty   B_k \frac{(-2t)^k}{k!} + \frac{1}{2}  \sum_{k =0}^\infty   B_k \frac{(2t)^k}{k!}  \\
 &=  \sum_{k =0}^\infty           \tfrac{1}{2}  \left(   (-1)^k + 1  \right)    B_k \frac{(2t)^k}{k!},
\end{align*}
so the only surviving terms in the latter series are the terms  whose index $k$ is
an even integer.  This yields
$
t \coth t = \sum_{n=0}^\infty  \frac{2^{2n}}{(2n)!} B_{2n} t^{2n}.
$

\medskip  
Exercise \ref{compute FT for exponential of abs value}
\quad 
We know, by equation \eqref{FT of the abs value exponential},  that the Fourier transform of
$f(x):=e^{-2\pi t |x|}$ is equal to $\hat f(\xi) =  \frac{  t   }{\pi (t^2 + \xi^2)}$. 
So using Poisson summation, we have:
 \index{Poisson summation formula}
\begin{align*}
\sum_{n \in \Z} e^{-2\pi t |n|} = \sum_{n \in \Z} f(n) 
= \sum_{\xi \in \Z}\hat  f(\xi) 
=  \frac{t}{\pi} \sum_{\xi \in \Z} \frac{1}{\xi^2 + t^2}. 
\end{align*}

\medskip
Exercise  \ref{continuity of convolution for two L^2 functions}
We are given that $f, g \in L^2(\R^d)$, and we wish to prove that
$f*g$ is always continuous on $\R^d$.   By definition, we need to show that
$\lim\limits_{h\rightarrow 0}\int_{\R^d} f(x-h - y)g(y) dy = (f*g)(x)$.

To prove part \ref{first part of continuity, L^2 convolution}, we fix any sequence of functions $f_n \in L^2(\R^d)$ with the property that
$f_n \rightarrow f$ in $L^2(\R^d)$.  We'll prove that 
$\lim\limits_{n \rightarrow \infty} \big(  (  f_n - f   )* g  \big)(x)= 0$,  for each $x \in \R^d$.  Well, we have:
\begin{align}
\big|  \left((f_n - f)*g\right)(x) \big| 
& \leq  \int_{\R^d} \left| f_n(x-y) - f(x-y)\right| |g(y)| dy    \label{triangle ineq for part a} \\  
\label{C-S for part b}
&\leq  \int_{\R^d} \left| f_n(x-u) - f(x-u)\right|^2 du 
 \int_{\R^d} |g(v)|^2 dv    \\
 &= \int_{\R^d} \left| f_n(u) - f(u)\right|^2 du 
 \int_{\R^d} |g(v)|^2 dv 
\end{align}
using the triangle inequality for integrals in \eqref{triangle ineq for part a}, and the Cauchy-Schwartz inequality in \eqref{C-S for part b}.  Since $f_n \rightarrow f$ in $L^2(\R^d)$, we are done.
For part \ref{second part of continuity, L^2 convolution},  we must show that 
\[
\int_{\R^d} \big| {T_h f}(x) - f(x) \big|^2 dx := \int_{\R^d} \big| f(x-h) - f(x) \big|^2 dx
\rightarrow 0, 
\]
as $h \rightarrow 0$.  First we'll show that the latter integral converges for each fixed nonzero vector $h$:
\begin{align*}
\int_{\R^d} \big| f(x-h) - f(x) \big|^2 dx &\leq \int_{\R^d} \left | f(x-h) \right |^2 dx
+ 2 \int_{\R^d} \left | f(x-h) f(x) \right | dx     +    \int_{\R^d} \left | f(x) \right |^2 dx \\
&\leq 2 \| f \|_{L^2(\R^d)}^2 
+ 2 \left( \int_{\R^d} \left | f(x-h) \right |^2 dx \right)^{\frac{1}{2}}
\left( \int_{\R^d} \left | f(x) \right |^2 dx \right)^{\frac{1}{2}}\\
&\leq 4 \| f \|_{L^2(\R^d)}^2.
\end{align*}
using $\int_{\R^d} \left | f(x-h) \right |^2 dx = \int_{\R^d} \left | f(x) \right |^2 dx:=  \| f \|_{L^2(\R^d)} < \infty$, and the
Cauchy-Schwartz inequality for $ \int_{\R^d} \left | f(x-h) f(x) \right | dx$.  We therefore have convergence of the integral for each nonzero $h$.  Next, we'll separate the integral into two pieces, one of which is a `neighborhood of infinity':
\[
\int_{\R^d} \big| f(x-h) - f(x) \big|^2 dx =
\int_{\|x\| > R} \big| f(x-h) - f(x) \big|^2 dx +  \int_{\|x\| \leq R} \big| f(x-h) - f(x) \big|^2 dx.
\]
By the convergence of the integral, we know that given any $\varepsilon >0$, there exists $R>0$ such that 
\[
\int_{\|x\| > R} \big| f(x-h) - f(x) \big|^2 dx < \frac{\varepsilon}{2}.
\]
It remains to handle the remaining integral, 
where we'll label the remaining compact set $E:=\{ x \in \R^d \mid \|x\| \leq R \}$:
\[
\int_{E} \big| f(x-h) - f(x) \big|^2 dx < \vol E \cdot
\sup_{x\in E}\{ \big| f(x-h) - f(x) \big|^2\}.
\]
Although $f$ may not necessarily be continuous, we may still conclude that as $h\rightarrow 0$, the latter expression tends to $0$ (otherwise the integral would diverge), finishing part \ref{second part of continuity, L^2 convolution}.

To prove part \ref{third part of continuity, L^2 convolution},
we must show that $\lim\limits_{h\rightarrow 0} (T_h f * g)(x) = (f*g)(x)$, for each $x \in \R^d$.  We pick the sequence of functions $f_n:= T_{h_n} f$,  with some sequence of vectors $h_n \rightarrow 0$. 
By part \ref{second part of continuity, L^2 convolution} we know that $T_{h_n} f \rightarrow f$ in $L^2(\R^d)$.  So $f_n \rightarrow f$ in $L^2(\R^d)$.  Now we may invoke part \ref{first part of continuity, L^2 convolution} to  conclude that 
$\lim\limits_{n \rightarrow \infty} \big(  (  f_n - f   )* g  \big)(x)= 0$.  In other words, we've shown that
$\lim\limits_{h \rightarrow 0} \big(  (T_h f* g)  \big)(x)= (f*g)(x)$, for each $x \in \R^d$, meaning that
 $f*g$ is continuous on $\R^d$.

\medskip
Exercise  \ref{strictly less than the FT at zero}
To prove part (a), suppose we are given $f\in L^1(\R^d)$ with $f(x) >0$ for all $x\in \R^d$.  
By the triangle inequality, we know that 
$|\hat f(\xi)| \leq \int_{\R^d} |f(x) e^{2\pi  i \langle x,\xi \rangle}  | dx 
= \int_{\R^d} f(x) dx:= \hat f(0)$, where $|f(x)| = f(x)$ follows from our assumption that 
$f(x) > 0$ for all $x\in \R^d$.  
To prove the strict inequality
$|\hat f(\xi)| < \hat f(0)$, for all nonzero $\xi \in \R^d$, suppose to the contrary that there exists 
a nonzero $\xi \in \R^d$ such that $|\hat f(\xi)| = \hat f(0)$.  Then
\[
\left | \int_{\R^d} f(x) e^{2\pi  i \langle x,\xi \rangle}  dx \right | =  \int_{\R^d} f(x) dx 
=  \int_{\R^d} |  f(x) e^{2\pi  i \langle x,\xi \rangle}   | dx,
\]
and we can now invoke Corollary \ref{cor:equality case of the triangle inequality}, which allows us to conclude that 
\[
\alpha \left(  f(x) e^{2\pi  i \langle x,\xi \rangle}  \right)
 =  \left | f(x) e^{2\pi  i \langle x,\xi \rangle} \right | 
 = |f(x)| = f(x),
\]
for some complex constant $\alpha:= e^{2\pi i \theta}$ on the unit circle, and for almost all $x\in \R^d$ ($\alpha=1$ is also allowed and poses no problems).  
In other words, we have $f(x) \left( 
e^{2\pi i \theta}    e^{2\pi  i \langle x,\xi \rangle}  -1
\right)=0$ almost everywhere.  Now our assumption that $f(x) >0$ for all $x\in \R^d$ implies that
$e^{2\pi i \theta}    e^{2\pi  i \langle x,\xi \rangle}  =1$ for almost all $x$. 
But this is a contradiction because $e^{2\pi i (\theta+ \langle x,\xi \rangle )} = 1$ precisely when 
$\theta+ \langle x,\xi \rangle \in \Z$.  That is, the latter condition occurs exactly when 
$x$ belongs to the discrete union of hyperplanes 
\[
\{ x\in \R^d \mid \langle x,\xi \rangle = -\theta + \Z\}, 
\]
a set of measure $0$ (for the $d$-dimensional meassure in $\R^d$).


Part (b) is almost identical.  Again arguing by contradiction, we suppose that 
there exists 
a nonzero $\xi \in \R^d$ such that $|\hat 1_\P(\xi)| = \hat 1_\P(0) := \vol \P$. 
We proceed in exactly the same manner, where the only difference is that we
 replace all the integrals over $\R^d$ by integrals over $\P$.  
 We  arrive at the following conclusion:  $e^{2\pi i (\theta+ \langle x,\xi \rangle )} = 1$ for almost all $x \in \P$.  This is again a contradiction, because the solution set to the latter equality is precisely the finite union of hyperplanes
$ \{ x\in \P \mid \langle x,\xi \rangle = -\theta + \Z\}$, which has measure $0$ (as a $d$-dimensional subset of $\R^d$).

\medskip
Exercise \ref{positive FT over R}  
We need to show that there exist two real numbers $r, s$ such that  
\[
f:= 1_{[-r, r]}*1_{[-r, r]}  + 1_{[-s, s]}*1_{[-s, s]}
\]
enjoys the property:
 \[
 \hat f(\xi) >0,
 \]
for all $\xi \in \R$.  Let's pick any two real numbers $r, s$ that are incommensurable, meaning that 
$\frac{r}{s} \notin \Q$.  Using \eqref{Stretch lemma for the sinc function}, we compute $\hat f$:
\[
\hat f(\xi):= \Big( \hat 1_{[-r, r]}(\xi) \Big)^2 +  \Big( \hat 1_{[-s, s]}(\xi) \Big)^2    = 
\left(  \frac{ \sin(2r\pi \xi)  }{   \pi \xi     }  \right)^2  +  \left(  \frac{ \sin(2s\pi \xi)  }{   \pi \xi     }  \right)^2 \geq 0.
\]
To prove positivity, suppose to the contrary that there exists a nonzero $\xi\in \R$ such that $\hat f(\xi)=0$.    Then
$\left(  \sin(2r\pi \xi)  \right)^2  +  \left(   \sin(2s\pi \xi)  \right)^2 = 0$,  but the vanishing of a sum of two squares (of real numbers) implies that
they must both equal $0$:
\[
 \sin(2r\pi \xi) =0, \text{ and }     \sin(2s\pi \xi)   = 0.
\]
Therefore $2r \pi \xi = m \pi$ and $2s \pi \xi = n \pi$, for some integers $m, n$.   We conclude that $\xi = \frac{m}{2r} = \frac{n}{2s}$,
so $\frac{r}{s}=\frac{m}{n} \in \Q $, a contradiction that proves $\hat f(\xi) >0$ for all nonzero real $\xi$.

\bigskip
Exercise \ref{tricky application of Poisson summation}
By assumption, $g:\R^d\rightarrow \C$ is infinitely smooth, and compactly supported. 
By Corollary \ref{cor: f smoother implies FT of F decays faster},  $\hat g$ is a rapidly decreasing function.  Because $g$ has compact support, we also know that $\hat g$ is infinitely smooth.  So $\hat g$ is a Schwartz function (and
$g$ is also a Schwartz function - in fact $g$ is a  `bump function', by definition).  Therefore we may apply the Poisson summation formula for Schwartz functions (Theorem \ref{Poisson.Summation})  to $\hat g$:
\[
\sum_{\xi \in \Z^d}  \hat g(\xi)  = \sum_{n \in \Z^d}  g(n),
\]
which is a finite sum due to the compact support of $g$.



 \bigskip \bigskip 
{\bf \Large Chapter \ref{Geometry of numbers}}

\medskip
Exercise \ref{easy first problem}
\quad
We're given a symmetric convex body $K\subset \R^2$ of area $4$, which contains only the origin.  
By Theorem \ref{thm:extremal bodies}, 
 $\tfrac{1}{2}K$ must tile $\R^2$ by translations with vectors from $\Z^2$, because 
 $2^2\det \Z^2 = 4 = \vol K$ (and $\tfrac{1}{2}K$ is therefore an extremal body).  But since
  $\tfrac{1}{2}K$ tiles $\R^2$ by translations, so does $K$ itself.

\medskip
Exercise \ref{convexity of K-K}
\quad 
We're given $d$-dimensional compact convex  sets $K, L \subset \R^d$. 
To prove that $K+L$ is convex, pick any $x, y \in K+L$, and we must show that 
$ \lambda_1 x + \lambda_2 y \in  K+L$ for all nonnegative  $\lambda_1, \lambda_2$ with 
$\lambda_1+ \lambda_2 =1$.   By assumption $x = k_1 + l_1$ and $y = k_2 + l_2$, with $k_1, k_2 \in K$, $l_1, l_2 \in L$.  We have:
\[
\lambda_1 x + \lambda_2 y = \lambda_1 \left( k_1 + l_1 \right) + \lambda_2 \left( k_2 + l_2 \right) 
= \left( \lambda_1  k_1 + \lambda_2 k_2 \right) +  \left( \lambda_1  l_1 + \lambda_2 l_2 \right) \in K + L,
\]
where we used the convexity of $K$ and of $L$ in the very last step above.  The same conclusion holds for $K-L$, because the convexity of $L$ implies the convexity of $-L$.

\medskip
Exercise \ref{practice with Minkowski sums}
\quad
We are given  $d$-dimensional compact, convex  sets $A, B \subset \R^d$. 
To prove that $ A \cap B \subseteq \tfrac{1}{2} A + \tfrac{1}{2} B$, we pick any
 $x \in A \cap B$.  Noticing that  $x = \tfrac{1}{2} x + \tfrac{1}{2}x$, where $\tfrac{1}{2} x \in \tfrac{1}{2}A$ and $\tfrac{1}{2} x \in \tfrac{1}{2}B$, we're done.   

To prove the second containment $ \tfrac{1}{2} A + \tfrac{1}{2} B \subseteq \conv\left( A \cup B \right)$, we pick $y \in \tfrac{1}{2} A + \tfrac{1}{2} B$.   So we may write $y = \tfrac{1}{2} a + \tfrac{1}{2} b$, where $a \in A, b \in B$, which is a convex linear combination of elements from $A$ and $B$, hence belongs to  
$\conv\left( A \cup B \right)$.  We'll leave the equality cases for the reader.

\medskip
Exercise  \ref{Convexity: A+A = 2A}
\quad
We are given a $d$-dimensional convex  set $A \subset \R^d$. To prove that
$A + A = 2A$, we pick any $x, y \in A$.  By the convexity of $A$, we know that 
$\frac{1}{2}x + \frac{1}{2}y \in A$, so that $x + y \in 2A$, proving that $A  + A \subseteq 2A$.
For the reverse inclusion $A + A \supseteq 2A$, we just notice that for any $a \in A$, 
$2a = a + a \in A+A$.

\medskip
Exercise \ref{c.s. C equals its symmetrized body}  
\quad  
For part (a), we suppose that 
\begin{equation} \label{symmetrized toy}
   \frac{1}{2}C - \frac{1}{2}C = C.
\end{equation}
For any $x \in C$, we need to
show  that $-x \in C$.   Since $x  \in \frac{1}{2}C - \frac{1}{2}C$,
we know that there must exist $y, z\in C$ such that $x=\frac{1}{2}y - \frac{1}{2}z$.  
This implies that $-x =   \frac{1}{2}z  -\frac{1}{2}y \in  \frac{1}{2}C - \frac{1}{2}C\subseteq C$. Therefore $C$ is centrally symmetric. 

To show part (b), first let's suppose that $C$ is convex and centrally symmetric (cs). Then
$\frac{1}{2}C - \frac{1}{2}C =\frac{1}{2}C + \frac{1}{2}C$.  Now using convexity, we claim that
$\frac{1}{2}C + \frac{1}{2}C = C$.  The convexity assumption implies that
$\frac{1}{2}C + \frac{1}{2}C \subseteq C$, because
for any $x, y \in C$, we have $ \frac{1}{2}x + \frac{1}{2}y \in C$. On the other hand, we always have 
$\frac{1}{2}C + \frac{1}{2}C \supseteq C$, because we can write each $x\in C$ as 
$x=\frac{1}{2}x + \frac{1}{2}x \in \frac{1}{2}C + \frac{1}{2}C$. 

So altogether we have $\frac{1}{2}C - \frac{1}{2}C =\frac{1}{2}C + \frac{1}{2}C = C$, proving the first direction.  

For the other direction of part (b), we  assume that 
\begin{equation}\label{hypothesis for exercise convexity}
\frac{1}{2}C - \frac{1}{2}C = C,
\end{equation}
 and we need to prove that $C$ is convex and cs.  By part (a), we already know that 
 $\frac{1}{2}C - \frac{1}{2}C$ is cs, hence the hypothesis 
 \eqref{hypothesis for exercise convexity} shows that $C$ must also be cs.

To prove convexity, let $x, y \in C$.   Using the hypothesis \eqref{hypothesis for exercise convexity}, together with the central symmetry of $C$, we have $C = \frac{1}{2}C - \frac{1}{2}C =
\frac{1}{2}C + \frac{1}{2}C$, so in particular $\frac{1}{2}x + \frac{1}{2}y \in C$.

For part (c), a compact counter-example is given by 
 $C:= [-2, -1] \cup [1, 2]$, a nonconvex set in $\R$.  Here $C$ is centrally symmetric, yet
$C - C = [-3, 3] \not=       [-4, -2] \cup [2, 4]   =  2C$.

Another (non-compact) counter-example is $\Z$, which is not convex, yet clearly centrally-symmetric.

\bigskip
Exercise   \ref{support of convolution}  
\quad   To prove part (a),   we are given two convex bodies $A, B \subset \R^d$, so by definition we have 
\[
\supp( 1_A * 1_B)  :=  \closure \left\{   y \in \R^d \bigm |   \int_{\R^d}  1_A(x) 1_B(y-x) dx \not= 0 \right\},
\]
and we must prove that $\supp( 1_A * 1_B)  = A + B$, their Minkowski sum.  \index{Minkowski sum}
In general, we have:
\begin{align}  \label{equivalences for Minkowski sums}
1_A(x) 1_B(y-x)>0 
& \iff   1_A(x) =1 \text{ and }   1_B(y-x) = 1  \\
& \iff  x\in A \text{ and }   y-x \in B   \\
& \iff   y \in A+B.
\end{align}
If we fix any $y \notin \supp( 1_A * 1_B)$,
then  $\int_{\R^d}  1_A(x) 1_B(y-x) dx = 0$, which implies that  $1_A(x) 1_B(y-x)=0$ for all $x\in \R^d$. 
But by the equivalences \eqref{equivalences for Minkowski sums}
above, we see that
$1_A(x) 1_B(y-x)=0 \iff   y\notin A+B$, proving that  $  A + B \subset \supp( 1_A * 1_B)$.  

Conversely, suppose that $y \in  \supp( 1_A * 1_B)$, meaning that there exists a sequence $y_n\in \R^d$ with
$ \int_{\R^d}  1_A(x) 1_B(y_n-x) dx \not= 0$.   This implies that for each such $y_n$, there exists at least
one $x \in \R^d$ with  $1_A(x) 1_B(y_n-x) >0$.  This last inequality, using our 
equivalences \eqref{equivalences for Minkowski sums}, implies that the sequence $y_n \in A+B$.
 Because $A+B$ is a closed set, we finally have $y := \lim_{n\rightarrow \infty} y_n \in A+B$.

\bigskip

To prove part (b), we must show that 
$ \supp(f*g)\subseteq  C$, where 
\[
C:= \closure\left(   \supp(f) +  \supp(g)  \right).
\]
We'll prove the contrapositive: if $x \notin C$, then $x \notin \supp(f*g)$ .   So we suppose $x \notin C$, and we have to prove that $(f*g)(x)=0$.    
By our assumption on $x$, for each $y \in \supp(g)$,  we have that $x-y \notin \supp(f)$.   The last assertion means that $f(x-y) = 0$, so we now know that $f(x-y) g(y)=0$ for all $y\in \R^d$.  
Finally,  we have $(f*g)(x) := \int_{\R^d} f(x-y) g(y) dy = 0$.

\bigskip
Exercise  \ref{exercise:unimodular simplex}
Show that in $\R^d$, an integer simplex $\Delta$ is unimodular $\iff \vol \Delta =  \frac{1}{d!}$.

\bigskip
Exercise \ref{non-unimodular but empty simplex}
\quad 
Define $\Delta:= \conv\{ (0, 0, 0), (1, 1, 0), (1, 0, 1), (0, 1, 1)\}$, an integer $3$-simplex.  It's clear that 
$\Delta$ is subset of the unit cube $[0, 1]^3$, and therefore $\Delta$ has no integer points in its interior. 
To see that $\Delta$  is not a unimodular simplex, its sufficient to consider its tangent $K_0$ cone at the origin, and show that this tangent cone is not unimodular.   $K_0$ has primitive integer edge 
vectors  $(1, 1, 0), (1, 0, 1), (0, 1, 1)$, so that the determinant of $K_0$ is 
equal to $\left| \det
\begin{pmatrix}
1 & 1 & 0  \\
1 & 0 & 1  \\
0 & 1 & 1
\end{pmatrix}
\right| = 2 > 1.
$

\medskip
Exercise   \ref{FT of a polytope is not Schwartz}
\quad
Suppose to the contrary, that for some polytope $\P$ we have $\hat 1_{\P}(\xi)= g(\xi)$, a Schwartz function. Taking the Fourier transform of both sides of the latter equality, and using the fact that the Fourier transform takes Schwartz functions to Schwartz functions, we would have $1_{\P}(-x) = \hat g(-x)$ is a Schwartz function.   But this is a contradiction, because the indicator function of a polytope is not even continuous.

\medskip
Exercise \ref{an application of Cauchy-Schwartz 1}
\quad We use the Cauchy-Schwartz inequality:  
\[
{\Big \langle      \icol{a \\ b} ,   \icol{ \sin x \\ \cos x}     \Big \rangle}^2
   := ( a \sin x+ b \cos x )^2   
    \leq       \big( a^2 + b^2  \big)   \big( \sin^2 x+ \cos^2 x \big) = a^2 + b^2. 
\]
By the equality condition of Cauchy-Schwartz, we see that the maximum is obtained when 
the two vectors are linearly dependent, which gives
$\tan x = \frac{a}{b}$.

\medskip
Exercise   \ref{elementary condition on the covariogram}
\quad
\rm{
Well, we have $x \in \tfrac{1}{2}K \cap \left(     \tfrac{1}{2}K +   n    \right) \iff x = \tfrac{1}{2}y$ and $x = \tfrac{1}{2}z + n$, where $y, z \in K$.  The latter conditions hold 
$\iff n = \tfrac{1}{2}y-\tfrac{1}{2}z \in \tfrac{1}{2}K-\tfrac{1}{2}K$.  Because of its convexity and central symmetry, we know that the latter condition is equivalent to $n \in K$, 
by Exercise \ref{c.s. C equals its symmetrized body}.
}


 \bigskip \bigskip 
{\bf \Large Chapter \ref{chapter.lattices}}

\medskip
Exercise  \ref{distance between hyperplanes}
\quad We are given the hyperplanes $H_1:= \{ x\in \R^d \mid   c_1 x_1 + \cdots + c_d x_d =  k_1\}$,
 and $H_2:= \{ x\in \R^d \mid   c_1 x_1 + \cdots + c_d x_d =  k_2\}$.  First we'll pick a point $x \in H_1$, 
 and then we'll  `walk along its normal vector', until we get to $H_2$.  With this `walk' in mind, we may assume WLOG
 that $k_2 > k_1$, and that the normal vector is pointing from $H_1$ towards $H_2$.  
 
 For simplicity, we'll let  $L:=  \sqrt{ c_1^2 + \cdots + c_d^2}$, and with this definition 
the unit normal vector to $H_1$ is 
 $n:= \frac{1}{ L   }(c_1, \dots, c_d)^T$.  We want to find $\delta>0$ such that
 $x + \delta n \in H_2$.  
 Unraveling the definition of the latter statement, we must have
 \begin{align*}
  & c_1 ( x_1 + \delta  \tfrac{1}{L} c_1) + \cdots + c_d (x_d + \delta  \tfrac{1}{L} c_d)  =  k_2 \\
\iff & (c_1 x_1 +   \cdots  + c_d x_d)   +  \frac{\delta}{L}( c_1^2  + \cdots + c_d^2)   =  k_2 \\
\iff & k_1 +  \delta\sqrt{c_1^2  + \cdots + c_d^2}  =  k_2 \\ 
\iff &   \delta   =  \frac{k_2 - k_1}{\sqrt{c_1^2  + \cdots + c_d^2}}.   
\end{align*}


\medskip
Exercise \ref{Hadamard's inequality, exercise}
\quad We consider each $k$'th row of $M$ as a vector, call it $v_k$.  By assumption, 
the norm of $v_k$ is bounded by $\|v\| \leq \sqrt{  B^2 + \cdots B^2    } = B\sqrt d$. 
Using Hadamard's inequality \ref{Hadamard inequality}, we have:
\begin{align*}
|\det M| \leq   \|v_1\| \cdots \|v_d\| \leq  \left(B\sqrt d \right)^d. 
\end{align*}

\medskip
Exercise  \ref{Ellipsoid problem}
\quad It's easy to see that the inverse matrix for $M$ is 
\[
M^{-1} := \begin{pmatrix} |  &  |  &  ... & |  \\  
                        \frac{1}{c_1}   b_1  &    \frac{1}{c_2} b_2     &  ... &     \frac{1}{c_d} b_d   \\  
                          |  &  |  &  ... & |  \\ 
        \end{pmatrix}^T.
\] 
The image of the unit sphere under the matrix $M$ is, by definition:
 \begin{align*}
 M(S^{d-1}) &:=  \{ u \in \R^d \mid u = Mx,  x \in S^{d-1}  \} \\
 &=  \{ u \in \R^d \mid   M^{-1}u  \in S^{d-1}  \} \\
  &=  \{ u \in \R^d \mid     \frac{1}{c_1^2}  \langle  b_1, u \rangle^2 + \cdots +    
                                        \frac{1}{c_d^2}  \langle  b_d, u \rangle^2 = 1 \},
\end{align*}
using our description of $M^{-1}$ above. 

For part (b), we begin with the definition of volume, and we want to compute the 
volume of the region $M(B):=  \{ u \in \R^d \mid u = My,  \text{ with } \| y \| \leq 1  \}$,
where $B$ is the unit ball in $\R^d$.
 \begin{align*}
 \vol(Ellipsoid_M)  &:= \int_{M(B)} du \\
 &= | \det M |  \int_{B} dy \\
 &= | \det M | \vol(B).
\end{align*}
using the change of variable $u = My$, with $y \in B$.  We also used the Jacobian,
 which gives $du = | \det M | dy$.
 
 Finally, we note that the matrix $M^T M$ is a diagonal matrix, with diagonal entries $c_k^2$, due
 to the fact that the $b_k$'s form an orthonormal basis.  Thus we use:  $| \det M |^2 = | \det M^T M | 
 =  \prod_{k=1}^d c_k^2$, so taking the positive square root, we arrive at $| \det M | = \prod_{k=1}^d c_k$, 
 because all of the $c_k$'s are positive by assumption.

\medskip
Exercise \ref{exercise:2by2 positive definite matrix}
\quad  Let $A:=
 \left(\begin{smallmatrix} a & b \\  b & d \end{smallmatrix}\right) $ be an invertible, symmetric matrix.
 Because $A$ is symmetric, we know both of its eigenvalues $\lambda_1, \lambda_2$ are real. 
The characteristic polynomial of $A$, namely
 $(a-\lambda)(d-\lambda) - b^2 $, may also be factored and  rewritten as 
 \[
 \lambda^2 - (a+ d) \lambda + (ad-b^2) = (\lambda - \lambda_1)(\lambda - \lambda_2) =  \lambda^2 - (\lambda_1 + \lambda_2)\lambda + \lambda_1 \lambda_2.
 \]
 Equating coefficients of the latter identity between polynomials, we therefore have $\lambda_1 + \lambda_2 =  {\rm Trace } A$, and  $\lambda_1 \lambda_2= \det A$.
 From these last two relations, we see that if both eigenvalues are positive, then ${\rm Trace } A>0$ and $\det A>0$.   
 
 Conversely, suppose that  ${\rm Trace } A>0$ and $\det A>0$. 
 Then $\lambda_1 \lambda_2 >0$, so either both eigenvalues are positive, or both eigenvalues are negative.  
But the eigenvalues cannot both be negative, for this would contradict our assumption that $\lambda_1 + \lambda_2>0$.

\medskip
Exercise \ref{permuting the elements of a group}
We're given any group $G$ (not necessarily finite), and any element $g \in G$.  We note that 
$gG\subseteq G$ 
by definition of closure in $G$: for any $h \in G$, we have $gh \in G$.  
To show $gG\supseteq G$,  we fix any $a \in G$ and we must find some $x\in G$ such that 
$gx = a$.  Since inverses exist in $G$, we find that $x = g^{-1} a$, and we're done.

 \bigskip \bigskip 
{\bf \Large Chapter \ref{Chapter.geometry of numbers II}}

\medskip
Exercise  \ref{one iteration of Minkowski symmetrization}
\quad
We're given $Q:= \tfrac{1}{2}K - \tfrac{1}{2}K$, where $K\subset \R^d$ is compact and convex. 
We already know that $Q$ is centrally symmetric.  Moreover, the convexity of $K$ implies that 
$\tfrac{1}{2}K +\tfrac{1}{2} K = K$.  So we have:
\begin{equation*} 
\tfrac{1}{2}Q - \tfrac{1}{2}Q =  
\frac{1}{2} \left( 
\tfrac{1}{2}K - \tfrac{1}{2}K 
\right) 
-\frac{1}{2} \left( 
\tfrac{1}{2}K - \tfrac{1}{2}K 
\right) =\tfrac{1}{4}K - \tfrac{1}{4}K - \tfrac{1}{4}K + \tfrac{1}{4}K = \tfrac{1}{2}K - \tfrac{1}{2}K = Q.
\end{equation*}

\medskip
Exercise \ref{elementary inequality, concavity 1}
\quad 
We are given  $r>1$ a fixed constant.  Dividing the inequality
$
(x+y)^r \geq x^r + y^r
$,
by $x^r$, it suffices to prove that $(1+t)^r \geq 1+ t^r$, for all positive $t$.  But this follows, for example, from the consideration of the function $f(t):=(1+t)^r  -1- t^r$ and the fact that its derivative $f'(t) =r(1+t)^{r-1} - r t^{r-1}$ is positive on $(0, \infty)$.

\medskip
Exercise  \ref{convergence of basic lattice sum}
We prove the claim by induction on the dimension $d$.   
For $d=1$, the claim is simply the usual test for convergence of the `p-series' $\sum\limits_{n\geq 1}\frac{1}{n^p}$, and therefore holds.  
Now we fix any $d\geq 2$, and we assume that $r > d$.  We must prove that  
$
\sum\limits_{n\in \Z^d} \frac{1}{ \| n \|^r} 
$
converges.  
The first step below makes use of the inequality $(x+y)^r \geq x^r + y^r$ for $r>1$ and $x, y > 0$ (see Exercise \ref{elementary inequality, concavity 1}).
We have:
\begin{align}
\| n \|^r    \label{elementary c>1}
&:= \left( n_1^2 + \cdots + n_d^2 \right)^\frac{r}{2} 
\geq     
     \left( n_1^2 \right)^\frac{r}{2} + \cdots +  \left( n_d^2 \right)^\frac{r}{2}    \\  
&=  |n_1|^r + |n_2|^r \cdots  + |n_d|^r  \\
& \geq
\left( \left( |n_1|^r \cdot |n_2|^r \cdots |n_d|^r \right) \right)^\frac{1}{d} d,  \label{AG inequality step}
\end{align}
using the Arithmetic-Geometric mean inequality in \eqref{AG inequality step}.   When considering the series
 $\sum\limits_{n\in \Z^d} \frac{1}{ \| n \|^r}$, we notice that by induction on the dimension it is sufficient to only prove convergence of the sub-series 
 $\sum\limits_{n\in \Z^d \atop n_1 n_2 \cdots n_d \not=0} \frac{1}{ \| n \|^r}$ with the property that \emph{none} of the coordinates of $n\in \Z^d$ vanish.  
From  \eqref{AG inequality step}, we have
\begin{align}
\sum_{n\in \Z^d \atop n_1 n_2 \cdots n_d \not=0} \frac{1}{ \| n \|^r}  
&<   \frac{1}{d}     
 \sum_{n_1\in \Z \atop n_1\not=0}  \frac{1}{ |n_1|^\frac{r}{d}  } 
 \cdots 
  \sum_{n_d\in \Z \atop n_d\not=0}  \frac{1}{ |n_d|^\frac{r}{d}  }  
   = \frac{2^d}{d}  \zeta^d \left(\frac{r}{d}\right),
\end{align}
which converges because $\frac{r}{d}>1$.

To prove the converse, we fix any $r\leq d$, and we must show that 
$\sum\limits_{n \in \Z^d-\{0\}}    \frac{1}{\|n\|^r} $  
diverges.
We recall the norm  $\|x\|_\infty:= \max\{ |x_1|, \dots, |x_d| \}$, for any $x \in \R^d$.
By Exercise \ref{elementary norm relations}, we had
$
  \| x \|    \leq   \sqrt{d} \,  \|x\|_\infty,  
$   
for all $x \in \R^d$. We therefore have 
\begin{equation}\label{comparing norms to get divergence}
\sum_{n \in \Z^d-\{0\}}    \frac{1}{\|n\|^r} 
\geq \frac{1}{d^{\frac{r}{2} }}  \sum_{n \in \Z^d -\{0\}} 
    \frac{1}{  \|n\|_\infty^r},  
\end{equation}
and the point is that now it is easy to count the number of integer points that have a fixed $\|n\|_\infty$ norm.  In fact, to count the number of integer points $n\in \Z^d$ such that $\|n\|_\infty =k$, we realize that this equals the number of integer points that lie on the boundary of the cube $[-k, k]^d$.   Thus, we may compute the number of these boundary integer points easily:
\[
(2k+1)^d - (2k-1)^d = \sum\limits_{j=0}^d \binom{d}{j} (2k)^{d-j}
- \sum\limits_{j=0}^d \binom{d}{j} (-1)^j (2k)^{d-j} = 2 \sum\limits_{m=0}^d \binom{d}{2m+1} (2k)^{d-2m-1},
\]
where the upper summation limit is never achieved, but takes care of both parity cases of $d$.  We notice that the latter finite sum is a sum of strictly positive terms, and to prove divergence we'll only keep the
 leading term $2 d (2k)^{d-1} = d 2^d k^{d-1}$.  Continuing from \eqref{comparing norms to get divergence}, we have:
 \begin{equation}
\sum_{n \in \Z^d-\{0\}}    \frac{1}{\|n\|^r} 
\geq  \frac{1}{d^{\frac{r}{2} }}  \sum_{n \in \Z^d -\{0\}} 
    \frac{1}{  \|n\|_\infty^r} 
>  \frac{d 2^d}{d^{\frac{r}{2} }}  \sum_{k=1}^\infty  \frac{ k^{d-1}}{ k^r},
\end{equation}
 which diverges precisely when $r - d + 1 \leq 1$.  We conclude that we have divergence when $r \leq d$.

\medskip
Exercise \ref{Arithmetic-geometric inequality application}
\quad
By the Arithmetic-Geometric mean inequality, we know that $\frac{1+ a_j}{2} \geq \sqrt a_j$, for each  $1\leq j \leq d$, and now we multiply all of these together:
\[
\frac{1}{2^d}(1+a_1)(1+a_2) \cdots (1+a_d) \geq \sqrt {a_1 a_2 \cdots a_d} = 1.
\]

 \bigskip \bigskip 
{\bf \Large Chapter \ref{chapter.Brion}}

\medskip
Exercise  \ref{independent of edge vectors}
We are given $\alpha >0$, and a simplicial cone $\K_v$,   with edge vectors $w_1, \dots, w_d \in \R^d$.
By definition, $\det \K_v$ is the determinant of the matrix whose columns are the $w_k$'s. 
Replacing each $w_k(v)$ by $\alpha_k w_k(v)$, we see that  the determinant $|\det \K_v|$
gets multiplied  by $\alpha^d$, and so
\[
\frac{ \alpha^d |\det \K_v|  }{\prod_{k=1}^d  \langle  \alpha w_k(v), z  \rangle}
= \frac{ |\det \K_v|  }{\prod_{k=1}^d  \langle  w_k(v), z  \rangle}.
 \]

\medskip
Exercise \ref{duality of dual cone}
We have to show that if we have the inclusion of cones
 $\K_1 \subset \K_2$, then $\K_2^* \subset \K_1^*$.   So we let 
 $x \in \K_2^*:=  \{  x \in \R^d \mid \langle x, u \rangle < 0 \text{ for all } u\in \K_2  \}$, implying that in particular  
 $\langle x, u \rangle < 0 \text{ for all } u\in \K_1$, because  $\K_1 \subset \K_2$.  But by definition this means that $x \in \K_1^*$ as well.

\medskip
Exercise \ref{polytope from pentagons}
\quad 
Euler's formula gives us 
\[
V-E+F =2,
\]
and the hypotheses also imply that:
\begin{align}
5F&=2E \\  
5F &\geq 3V.   
\end{align}
Altogether, we get 
\begin{equation*}
2=V-E+F \leq \frac{5}{3} F  - \frac{5}{2} F + F = \frac{1}{6} F,
\end{equation*}
so that $F \geq 12$.

 \bigskip \bigskip 
{\bf \Large Chapter \ref{chapter:Discrete Brion}}

Exercise  \ref{unimodular cone, integer point transform}
\quad  
The main point here is that at each vertex $v\in V$, the edge vectors form a basis for $\Z^d$, and therefore 
the only integer point in the (half-open) fundamental parallelepiped $\Pi_v$ is $v$ itself.   So we see that 
its integer point transform of $\Pi$ is $\sigma_{\Pi_v}(x) = e^{\langle v, z  \rangle}$. 
Now we use Theorem \ref{brion, discrete form}, followed by Theorem \ref{closed form for integer point transform of a cone}: 
\begin{equation*}
\sigma_\P(z) = \sum_{v \in V}  \sigma_{\K_v}(z) = 
\frac{e^{\langle v, z \rangle} }{\prod_{k=1}^d \left( 1 - e^{\langle w_k, z\rangle}  \right) }.
\end{equation*}

Exercise \ref{bound for integer point transform}
\quad 
Because $|e^{2\pi i \langle x, n \rangle}|=1$ for all $x\in \R^d$, we have:
\[
\left |  \sigma_\P(2\pi i x)  \right |        
    \leq         \sum_{n\in \P \cap \Z^d } \left|   e^{2\pi i \langle x, n \rangle}  \right|
= \sum_{n\in \P \cap \Z^d } 1 =   \left |\Z^d \cap \P \right |.
\]

 \bigskip \bigskip 
{\bf \Large Chapter \ref{Ehrhart theory}}

 \medskip
Exercise  \ref{Bernoulli polynomial as an Ehrhart polynomial}
\quad  
Here $\P:=  \conv\{   C,   {\bf e_d} \}$, where $C$ is the $(d-1)$-dimensional unit cube $[0, 1]^{d-1}$.  To compute the Ehrhart
polynomial $\L_{\P}(t)$ here, we use the fact that a `horizontal' slice of $\P$, meaning a slice parallel to $C$, and orthogonal
to $e_d$, is a dilation of $C$.   Thus, each of these slices counts the number of points in a $k$-dilate of $C$, as $k$ varies from $0$ to $t+1$.   Summing over these integer dilations of $C$, we have
\[
\L_{\P}(t) = \sum_{k=0}^{t+1}  (t+1 - k)^{d-1} =  \sum_{k=0}^{t+1}  k^{d-1} =
 \frac{1}{d}(B_d(t+2) - B_d),
\]
where the last step holds thanks to Exercise \ref{historical origin of Bernoulli poly}.

\medskip 
Exercise \ref{unimodular triangle}
Using Pick's formula, the unimodular triangle $\P$ has area:
\[
\rm{Area } \P  = I + \frac{1}{2} B -1 = 0 + \frac{1}{2} 3 -1 = \frac{1}{2}.
\]

 \medskip
Exercise  \ref{properties of floor, ceiling, fractional part}

 Throughout, we first write $x = n + \alpha$, with $\lfloor x \rfloor := n \in \Z$ and $0 \leq \alpha < 1$. 
We prove part \ref{problem:floor and celings, a}, namely that $- \floor{-x} = \left\lceil  x \right\rceil$.
  
Case $1$:   $x \in \Z$.    Here $\alpha = 0$ and $x=n$, so that $- \floor{-x} = - (-n)= n =  \left\lceil  x \right\rceil$. 

Case $2$:   $x \notin \Z$.  In this case  
$\left\lceil  x \right\rceil = n+1$.   We have   $-x = -n - \alpha = - n-1 + (1-\alpha)$,
 from which we see that  $- \floor{-x} = - (-n-1) = n + 1 =  \left\lceil  x \right\rceil$.

To prove part \ref{problem:floor and celings, b}, we need to show that $\floor{x} - \left\lceil x  \right\rceil  +1=1_{\Z}(x)$.

Case $1$:   $x \in \Z \implies  \floor{x} - \left\lceil x  \right\rceil  +1 = n - n + 1 = 1 = 1_{\Z}(x)$.

Case $2$:   $x \notin \Z  \implies \floor{x} - \left\lceil x  \right\rceil  +1 
= n - (n + 1)  + 1  = 0 = 1_{\Z}(x)$.

To prove part \ref{problem:floor and celings, c}, we need to show that 
$ \{ x \} + \{-x\} = 1- 1_{\Z}(x)$.
This follows from part \ref{problem:floor and celings, b} if we use the definitions 
 $\floor{x} := x - \{x\}, \lceil x \rceil:= x + \{x\}$.   Using the identity of part \ref{problem:floor and celings, b}, we have
 \[
 1- 1_{\Z}(x) =  \left\lceil x  \right\rceil - \floor{x}  =  x + \{x\} - ( x - \{x\}) =   \{x\} + \{x\}.
 \]

To prove part \ref{problem:floor and celings, d}, we must show that
$\floor{  x + y } \geq \floor{ x } + \floor{y}$, for all $x, y \in \R$.  So we let $x := n + \alpha$ and $ y:= m + \beta$, where $n:= \floor x, m:= \floor y$, and by definition $0\leq \alpha< 1$, $0\leq \beta< 1$. 
Now $\floor{  x + y } = \floor{  n + m + \alpha + \beta }  \geq n+m := \floor{x}+\floor{y}$.

Finally, for part  \ref{problem:floor and celings, e}, we have to prove that 
if  $m \in \Z_{>0}, n \in \Z$, then $\floor{ \frac{n-1}{m} } + 1 = \left\lceil   \frac{n}{m} \right\rceil$.
We begin by using the division algorithm, which gives us $n=qm+r$, with integers $q$ and   $0\leq r < m$.

Case $1$:   $r=0$.  Here $n = qm$, and we have $\floor{ \frac{n-1}{m} } + 1 = \floor{ q- \frac{1}{m} } + 1=q 
= \frac{n}{m} =  \left\lceil   \frac{n}{m} \right\rceil$.

Case $2$:   $0 < r   < m$.  Here 
$\floor{ \frac{n-1}{m} } + 1 = \floor{ \frac{qm + r - 1}{m} } + 1
= \floor{ q + \frac{ r - 1}{m} } + 1 = \floor{ \frac{ r - 1}{m} } + 1 = 1$.
On the other hand, $  \left\lceil   \frac{n}{m} \right\rceil =  \left\lceil   \frac{qm + r}{m} \right\rceil   
=   \left\lceil   q + \frac{ r}{m} \right\rceil     =  \left\lceil    \frac{ r}{m} \right\rceil =1   $.

\medskip
Exercise \ref{ (-1)^[x] in terms of periodic Beroulli polys}
To show that for $x \in \R \setminus \Z$, we have
$(-1)^{\floor{x}} = 2\, P_1(x) - 4\,  P_1\left( \frac{x}{2} \right)$,  we expand the right-hand side:
\begin{align*}
2\, P_1(x) - 4 \, P_1\left( \frac{x}{2} \right) &:= 
2\left(  x - \floor{x} - \tfrac{1}{2}   \right) 
- 4 \left(  \frac{x}{2}  - \floor{\frac{x}{2} } - \frac{1}{2}\right) \\
&= 1 - 2\floor{x} + 4\floor{ \frac{x}{2} } = 1 - 2\floor{x} + 4\floor{ \frac{\floor{x}}{2} } \\
&=\begin{cases}  
1  -2\floor{x} + 4\left(  \frac{\floor{x}}{2} \right)  &      \mbox{if }  \floor{x}   \mbox{ is even}, \\ 
1  -2\floor{x} + 4\left(  \frac{\floor{x}-1}{2} \right)  &        \mbox{if } \floor{x}   \mbox{ is odd}.
\end{cases}
=\begin{cases}  
1    &      \mbox{if } \floor{x}   \mbox{ is even}, \\ 
-1  &        \mbox{if } \floor{x}   \mbox{ is odd}
\end{cases} \\
&= (-1)^{\floor{x}}.
\end{align*}


 \bigskip \bigskip 
{\bf \Large Chapter \ref{Stokes' formula and transforms} }

\medskip
Exercise  \ref{Stokes implies Cauchy}
Here $f(z) := u(x, y) + i v(x, y)$, with $z:=x+iy$, and we suppose that
$f'(z)$ is continuous on the unit ball $B:= \{ z\in \C \mid \|z\| \leq 1\}$.  
We must show that Stokes' theorem implies Cauchy's theorem:
$
\int_{S^1} f(z) dz = 0.
$
 We know by Goursat's Lemma that the continuity of the partial derivatives implies that 
 $f$ has a complex derivative in $B$.  Now the Cauchy-Riemann equations follow:
 $ \frac{ \partial u }{ \partial y} =  \frac{\partial v}{ \partial x }$.  So we conclude that:
 \begin{equation}
 \int_{S^1} f(z) dz =  \int_{S^1} \big( u(z) + i v(z) \big)  dz 
 = \int_B  \left( 
 \frac{ \partial u }{ \partial y} -  \frac{\partial v}{ \partial x }
  \right) dx dy
 = 0,
 \end{equation}
 where we've used Stoke's theorem in the penultimate equality.

\medskip
Exercise \ref{alternate combinatorial divergence Theorem}
We have to show that if   $F(x):= e^{-2\pi i \langle x, \xi \rangle} \lambda$, with a constant nonzero vector $\lambda \in \C^d$, then:
\begin{equation} \label{Ex:desired identity for alternate divergence thm}
\hat 1_\P(\xi) =\frac{1}{-2\pi i }
 \sum_{G\subset \partial P} 
\frac{   \langle \lambda, \n_G \rangle}{ \langle  \lambda,  \xi \rangle} 
\hat 1_G(\xi),
\end{equation}
valid for all nonzero $\xi \in \R^d$.  
Taking the divergence of the vector field $F(x)$, we have:
\begin{align}
{\rm div } F(x)  
&=   \frac{\partial  \left(  e^{- 2\pi i \langle x, \xi \rangle} \lambda_1  \right)   }{{\partial} x_1}  
   + \cdots +
        \frac{\partial  ( e^{- 2\pi i \langle x, \xi \rangle} \lambda_d  )}{\partial x_d} \\
&=    -2\pi i \langle \xi, \lambda \rangle  e^{- 2\pi i \langle x, \xi \rangle}.
\end{align}
By the divergence theorem, we now have
\begin{align}  \label{initial divergence}
\int_{x\in P} - 2\pi i \langle \xi, \lambda \rangle   e^{- 2\pi i \langle x, \xi \rangle} dx 
&=    \int_{x\in P}  \text{div} F(x) dx =    \int_{\partial P} e^{- 2\pi i \langle x, \xi \rangle} 
\langle \lambda,  \n \rangle \ dS \\
&= 
 \sum_{G\subset \partial P} 
 \langle \lambda, \n_G \rangle
\hat 1_G(\xi), \label{Exercise:FT of a polytope, using a constant vector field}
\end{align}
where $\n_G$ is the outward-pointing unit normal vector at each point of the facet $G\in \partial \P$,
and where we've used  $1_{\partial \P}=  \sum_{G\in \partial P} 1_G$ (since $\P$ is a polytope).
In other words, we have:
\begin{equation} \label{Ex: final desired identity for constant vector field}
 \langle  \lambda,  \xi \rangle  \hat 1_\P(\xi) =\frac{1}{-2\pi i }
 \sum_{G\subset \partial P} 
  \langle \lambda, \n_G \rangle 
\hat 1_G(\xi),
\end{equation}
the desired identity.

\medskip
Exercise \ref{equivalent identity to the alternate vector field}
Revisiting \eqref{Ex: final desired identity for constant vector field}
in  Exercise \ref{alternate combinatorial divergence Theorem}, we have
\begin{equation} 
 \langle  \lambda,  \xi \rangle  \hat 1_\P(\xi) =\frac{1}{-2\pi i }
 \sum_{G\subset \partial P} 
  \langle \lambda, \n_G \rangle 
\hat 1_G(\xi),
\end{equation}
for each constant vector $\lambda\in \C^d$ with nonzero imaginary part, and for each
$\xi \in \R^d$, including $\xi=0$. In other words,
\begin{equation} 
 \langle  \lambda,  \xi \hat 1_\P(\xi) \rangle   =\frac{1}{-2\pi i }
  \langle \lambda, 
  \sum_{G\subset \partial P} 
 \n_G 
\hat 1_G(\xi)  \rangle
\end{equation}
gives us the vector identity
\begin{equation}
\xi   \hat 1_\P(\xi)    =   \frac{1}{-2\pi i }        \sum_{G\subset \partial P} \n_G    \hat 1_G(\xi),
\end{equation}
valid for all  $\xi \in \R^d$.

\medskip
Exercise \ref{strange proof of Minkowski relation}
To show that  Exercise \ref{equivalent identity to the alternate vector field} 
easily gives us
the Minkowski relation  \eqref{Minkowski relation}, we simply evaluate both sides
of $
\xi   \hat 1_\P(\xi)    =   \frac{1}{-2\pi i }        \sum_{G\subset \partial P} \n_G    \hat 1_G(\xi),
$
 at $\xi=0$:
\[
0 =  \frac{1}{-2\pi i }        \sum_{G\subset \partial P} \n_G    \hat 1_G(0)
=\sum_{ \text{facets } G \text{ of } P}    (\vol G)  \n_{G}.
\]

\medskip
Exercise \ref{Exercise Iverson bracket}
Each of these identities is easily proved using a truth table.  For example, to prove that
 $[ P \rm{ \ or \ } Q ] =  [P] + [Q] - [P][Q]$, we compare two truth tables:

Truth table for $ [ P \rm{ \ or \ } Q ]$:

\begin{tabularx}{0.8\textwidth} { 
  | >{\raggedright\arraybackslash}X 
  | >{\centering\arraybackslash}X 
  | >{\raggedleft\arraybackslash}X | }
 \hline
  & Q is true & Q is false \\
 \hline
P is true & $1$  & $1$  \\
\hline
P is false  & $1$ & $0$  \\
\hline
\end{tabularx}

Truth table for  $[P] + [Q] - [P][Q]$:

\begin{tabularx}{0.8\textwidth} { 
  | >{\raggedright\arraybackslash}X 
  | >{\centering\arraybackslash}X 
  | >{\raggedleft\arraybackslash}X | }
 \hline
  & Q is true & Q is false \\
 \hline
P is true & $1+1-1 = 1$  & $1+0-0=1$  \\
\hline
P is false  & $0+1 - 0 = 1$ & $0+0-0=0$  \\
\hline
\end{tabularx}

 \bigskip \bigskip 
{\bf \Large Chapter \ref{Chapter.geometry of numbers III}}

\medskip
Exercise \ref{translating the Voronoi cell around} 
\quad 
Given a full rank lattice $\L\subset \R^d$, and any $m \in \L$, we have:
\begin{align}
{\rm Vor}_0(\L) + m 
&:=   
 \left\{ x +m \in \R^d \bigm |     \|x\|  \leq   \|x - v\|, \ \text{ for all  } v \in \L    \right\}  \\  \label{last expression for Voronoi}
&=  \left\{ y \in \R^d \bigm |     \|y-m\|  \leq   \|y-m - v\|, \ \text{ for all  } v \in \L    \right\}.
\end{align}
But as $v$ varies over $\L$, so does $m + v$, because $m \in \L$.  Hence the expression 
\eqref{last expression for Voronoi} above is equal to  ${\rm Vor}_m(\L)$.

\medskip
Exercise \ref{Rewriting Minkowski's first theorem} 
Let $K\subset \R^d$ be a $d$-dimensional convex body, symmetric about the origin,
and let $\L$ be a (full rank) lattice in $\R^d$.
 We have to show that the following two statements are equivalent. 
\begin{enumerate}[(a)]
\item \label{first version of Minkowski's first theorem}
\begin{equation} 
\text{ If }   \vol K > 2^d (\det \L),  \text{  then }   K
\text{  must contain a nonzero point of } \L \text{ in its interior}. 
\end{equation}
\item  \label{second version of Minkowski's first theorem}
\begin{equation}
\lambda_1(K, \L)^d  \vol K  \leq  
 2^d \det \L.
\end{equation}
\end{enumerate}
First we assume that part \ref{second version of Minkowski's first theorem} is true.  So if 
 $\vol K > 2^d (\det \L)$, then
 \[
  2^d (\det \L) < \vol K \leq \frac{1}{\lambda_1(K, \L)^d} 2^d \det \L,
 \]
giving us $\lambda_1(K, \L)^d \leq 1$, or $\lambda_1(K, \L) \leq 1$.  By definition of the first successive minima, this means that $K$ contains a nonzero lattice point of $\L$.  So we've proven part \ref{first version of Minkowski's first theorem}.

Now we assume part \ref{first version of Minkowski's first theorem} is true. Suppose to the contrary that part \ref{second version of Minkowski's first theorem} is false.  Here the
 main idea is that $\vol(\lambda K) = \lambda^d \vol K$ for any positive number $\lambda$.
So we have 
\[
\vol\left( \lambda_1 K\right)  > 2^d \det \L.
\]
Now applying part \ref{first version of Minkowski's first theorem} to the new body $\lambda_1 K$, we may conclude that  $\lambda_1 K$ must contain a nonzero point of $\L$ 
in its interior.   But this contradicts the definition of $\lambda_1$, and we're done.

\medskip
Exercise \ref{equality condition for Minkowski's conjecture}
We need to
 prove that the equality conditions in Minkowski's conjecture \ref{Minkowski conjecture} 
are achieved by the diagonal linear forms $L_k(x_1, \dots, x_d):= 2 c_k x_k$. 
Letting $M$ denote the matrix formed by the coefficients of these diagonal linear forms,
we see that $\det M = 2^d c_1 \cdots c_d$.  To satisfy Minkowski's conjecture in this case, we need to find an integer point $n \in \Z^d$ such that 
$\left| (L_1(n) + c_1)\cdots(L_d(n)+c_d) \right|
  \leq \frac{|\det M|}{2^d}$, which in our case translates to the requirement that 
  \[
  \left| (2c_1n_1 + c_1)\cdots(2c_d n_d+c_d) \right|
  \leq \frac{2^d c_1 \cdots c_d }{2^d} = c_1 \cdots c_d.
  \]
But the latter inequality becomes an equality when $n=0$, so we're done.

\printindex
\backmatter

\bigskip
\begin{thebibliography}{100}


\bibitem{AkopyanKarasev}
A. V. Akopyan and R. N. Karasev, 
\emph{Bounding minimal solid angles of polytopes}, (2015) 
(\url{https://arxiv.org/abs/1505.05263}).

\bibitem{Alexandrov}
A.  D.  Alexandrov,  \emph{A theorem on convex polyhedra},
 Trudy Mat. Int. Steklov, Sect. Math, 4:87, (1933).

\bibitem{ApostolBook}
Tom M. Apostol, \emph{Introduction to Analytic Number Theory},
Springer Unergraduate texts in Mathematics, (1976),
1--350.

\bibitem{DavidAustin}
David Austin, \emph{Fedorov's five parallelohedra}, Notices of the American Math. Society, Feature column, 2013.   

\url{http://www.ams.org/publicoutreach/feature-column/fc-2013-11}


\bibitem{Averkov}
Gennadiy Averkov,  \emph{Equality Case in Van der Corput's Inequality and Collisions in Multiple Lattice Tilings},
Discrete \& Computational Geometry \textbf{65},  (2021), 212--226.

\bibitem{Babai}
L\'aszl\'o Babai,   \emph{On Lovász'  lattice reduction and the nearest lattice point problem},   \newline
Combinatorica \textbf{6},  (1986), 1--13.





\bibitem{Baillie}
Robert Baillie, David Borwein, and Jonathan M. Borwein,
\emph{Surprising sinc sums and integrals}, The American Mathematical Monthly, \textbf{115}(10), (2008), 888--901.




\bibitem{BaldoniBerlineVergne}
Velleda Baldoni, Nicole Berline, and Mich{\`e}le Vergne, \emph{Local
{E}uler--{M}aclaurin expansion of {B}arvinok valuations and {E}hrhart
coefficients of a rational polytope}, Integer points in polyhedra---geometry,
number theory, representation theory, algebra, optimization, statistics,
Contemp. Math., vol. 452, Amer. Math. Soc., Providence, RI, 2008, pp.~15--33.

\bibitem{KeithBall.1}
Keith Ball, \emph{A lower bound for the optimal density of lattice packings},
International Mathematics Research Notices, Vol 1992, Issue 10, 1992), 217--221.




\bibitem{Banaszczyk1}
W.  Banaszczyk, \emph{New bounds in some transference theorems in the geometry of numbers}, 
Math. Ann. \textbf{296} (1993),  625--635.

\bibitem{Barany}
Imre B\'ar\'any,  \emph{Random points and lattice points in convex bodies}, Bull. Amer. Math. Soc. (N.S.) \textbf{45} (2008), no. 3, 339--365.

\bibitem{BaranyAkopyanRobins}
Imre B\'ar\'any,  Arseniy Akopyan, and  Sinai Robins, Algebraic vertices of non-convex polyhedra, Advances in Math, \textbf{308}, (2017), 627-644.

\bibitem{Barvinok1}
Alexander Barvinok, \emph{Exponential integrals and sums over convex
polyhedra}, Funktsional. Anal. i Prilozhen. \textbf{26} (1992), no.~2,
64--66.

\bibitem{Barvinok.algorithm}
Alexander Barvinok,  \emph{A polynomial time algorithm for counting integral
points in polyhedra when the dimension is fixed}, Math. Oper. Res.
\textbf{19} (1994), no.~4, 769--779.

\bibitem{BarvinokNotes}
Alexander Barvinok,  \emph{Combinatorics, Geometry, and Complexity of integer points}, Online lecture notes:

\url{http://www.math.lsa.umich.edu/~barvinok/latticenotes669.pdf}

\bibitem{Barvinok.A.course.in.convexity}
Alexander Barvinok, \emph{A course in convexity},  Graduate Studies in Mathematics, 54. American Mathematical Society, Providence, RI, 2002. 

\bibitem{BarvinokEhrhartbook}
Alexander Barvinok, \emph{Integer points in polyhedra}, Zurich Lectures in Advanced
Mathematics, European Mathematical Society (EMS), Zurich, 2008.

\bibitem{BarvinokPommersheim}
Alexander Barvinok and James~E. Pommersheim, \emph{An algorithmic theory of
lattice points in polyhedra}, New {P}erspectives in {A}lgebraic
{C}ombinatorics (Berkeley, CA, 1996--97), Math. Sci. Res. Inst. Publ.,
vol.~38, Cambridge Univ. Press, Cambridge, 1999, pp.~91--147.

\bibitem{Batyrev}
Victor~V. Batyrev, \emph{Dual polyhedra and mirror symmetry for Calabi--Yau
hypersurfaces in toric varieties}, J. Algebraic Geom. \textbf{3} (1994),
no.~3, 493--535, {\tt arXiv:alg-geom/9310003}.

\bibitem{BatirevHofscheier}
Victor Batyrev and Johannes Hofscheier,
\emph{A generalization of a theorem of G. K. White},
Moscow Journal of Combinatorics and Number Theory  vol. 10(4) (2021), 281--296.

\bibitem{BeckSanyal}
Matthias Beck and Raman Sanyal, 
\emph{Combinatorial reciprocity theorems,  an invitation to enumerative geometric combinatorics},
Grad. Stud. Math., vol. 195, Providence, RI: American Mathematical Society (AMS), 2018.

\bibitem{BerlineVergne}
Nicole Berline and Mich\'ele Vergne,  \emph{Local Euler-Maclaurin formula for polytopes}, 
 Mosc. Math. J., \textbf{7} (3) (2007), 355--386.

\bibitem{JozsefBeck}
J\'ozsef Beck, \emph{Probabilistic Diophantine approximation, Randomness in lattice point counting}, Springer Monographs in Mathematics, Springer, Cham, (2014), 1--487.





\bibitem{BeckRobins}
Matthias Beck and Sinai Robins,
\emph{Computing the continuous discretely: integer-point enumeration in polyhedra}, 
$2$'nd edition, Springer, New York, (2015), 1--285.

\bibitem{BeckRobins.polygons}
Matthias Beck and Sinai Robins, \emph{Explicit and efficient formulas for the
lattice point count in rational polygons using {D}edekind--{R}ademacher
sums}, Discrete Comput. Geom. \textbf{27} (2002), no.~4, 443--459, {\tt
arXiv:math.CO/0111329}.




\bibitem{BeckDiaz}
Matthias Beck, Ricardo Diaz, and Sinai Robins, \emph{The Frobenius problem,
rational polytopes, and Fourier--Dedekind sums},    
J. Number Theory \textbf{96} (2002), no.~1, 1--21.
{\tt arXiv:math.NT/0204035}.



\bibitem{BeckRobinsSam}
Matthias Beck, Sinai Robins, and Steven~V Sam, \emph{Positivity theorems for
solid-angle polynomials}, Beitr\"age Algebra Geom. \textbf{51} (2010), no.~2,
493--507, {\tt arXiv:0906.4031}.


 \bibitem{BetkeHenkWills}
 Ulrich Betke, Martin Henk, and J{\"o}rg~M. Wills, 
 \emph{Successive-minima-type inequalities},
Discrete Comput. Geom. \textbf{9} (1993), 165--175. 

\bibitem{HenkRootsOfEhrhart}
Christian Bey, Martin Henk, and J{\"o}rg~M. Wills, \emph{Notes on the roots of
{E}hrhart polynomials}, Discrete Comput. Geom. \textbf{38} (2007), no.~1,
81--98, {\tt arXiv:math.MG/0606089}.


\bibitem{BambahWoods}
R. P. Bambah and A. C. Woods, 
\emph{Minkowski's conjecture for $n=5$;  a Theorem of Skubenko}, 
Journal of Number Theorey \textbf{12}, (1980) 27--48.

\bibitem{Banaszczyk}
W. Banaszczyk, \emph{New bounds in some transference theorems in the geometry of numbers}, 
Math. Annalen, \textbf{296} (1993), 625--635.

\bibitem{Banaszczyk2}
W. Banaszczyk,  \emph{Inequalities for convex bodies and polar reciprocal lattices in $\R^n$}, Discrete Comput. Geom, (1995), 217--231.

\bibitem{Bianchi1}
Gabriele Bianchi,
\emph{The covariogram and Fourier–Laplace transform in $C^n$},
 Proc. London Math. Soc. (3) 113 (2016) 1--23.






\bibitem{BisgardLinearAlgebraBook}
James Bisgard, 
\emph{Analysis and linear algebra: the singular value decomposition and applications},
Student Mathematical Library, 94. American Mathematical Society, Providence, RI, (2021), 1--217.

\bibitem{Blichfeldt2}
Hans F. Blichfeldt,  \emph{A new principle in the geometry of numbers, with some applications},
Trans. Amer. Math. Soc. \textbf{15} (1914), no. 3, 227--235.


\bibitem{Blichfeldt1}
Hans F. Blichfeldt,  \emph{The minimum value of quadratic forms, and the closest packing of spheres},
Math. Annalen \textbf{101} (1929), 605--608.




\bibitem{boas72}
R. P. Boas,   \emph{Summation formulas and band-limited signals},  Tohoku Math. J.,
\textbf{24} (1972), no.~2,  121--125.




\bibitem{bockerliptak}
Sebastian B{\"o}cker and Zsuzsanna Lipt{\'a}k, \emph{The money changing problem
revisited: computing the {F}robenius number in time {$O(ka_1)$}}, Computing
and combinatorics, Lecture Notes in Comput. Sci., vol. 3595, Springer,
Berlin, 2005, ~965--974.

\bibitem{BochnerBook}
Salomon Bochner, \emph{Lectures on Fourier integrals}, Princeton University Press, translated from the original by Morris Tenenbaum and Harry Pollard, (1959), 1--338.

\bibitem{BokowskiHadwigerWills}
J\"urgen  Bokowski,  Hugo Hadwiger, and J\"org M. Wills, 
\emph{Eine Ungleichung zwischen Volumen,
Oberfl\"ache und Gitterpunktanzahl konvexer K\"orper im n-dimensionalen euklidischen
Raum}, Math. Z. \textbf{127} (1972), 363--364.

\bibitem{Bombieri}
Enrico Bombieri, 
\emph{Sulla dimostrazione di C. L. Siegel del teorema fondamentale di Minkowski nella geometria dei numeri},
Bollettino dell’Unione Matematica Italiana, Serie 3, Vol. 17 (1962), n.3,  283--288.

\bibitem{BorodzikNguyenRobins}
Maciej Borodzik, Danny Nguyen, and Sinai Robins, \emph{Tiling the integer lattice with translated sublattices}, 
Moscow Journal of Combinatorics and Number Theory, Vol 6, issue 4, (2016), 3--26.

\bibitem{Butzer.etal}
P. L. Butzer, P. J. S. G. Ferreira, G. Schmeisser and R. L. Stens, 
\emph{The Summation Formulae of Euler-Maclaurin,  Abel-Plana,  Poisson, and their Interconnections with the Approximate Sampling Formula of Signal Analysis},
Results. Math. \textbf{59} (2011), 359--400.




\bibitem{BrandoliniColzaniTravagliniRobins1}
Luca Brandolini, Leonardo Colzani, Sinai Robins, and Giancarlo Travaglini,
\emph{Pick's Theorem and Convergence of multiple Fourier series}, The American Mathematical Monthly, \textbf{128}, Issue 1, (2021), 41--49.

\bibitem{BrandoliniColzaniTravagliniRobins2}
Luca Brandolini, Leonardo Colzani, Sinai Robins, and Giancarlo Travaglini,  
\emph{An Euler-MacLaurin formula for polygonal sums}, to appear in Transactions of the AMS, 2021.

\bibitem{Brion}
Michel Brion, \emph{Points entiers dans les poly\`edres convexes}, Ann. Sci. Ecole Norm. Sup.
\textbf{4} 21 (1988), no.~4, 653--663.

\bibitem{brionvergne}
Michel Brion and Mich{\`e}le Vergne, \emph{Residue formulae, vector partition
functions and lattice points in rational polytopes}, J. Amer. Math. Soc.
\textbf{10} (1997), no.~4, 797--833.

\bibitem{BrownSchreiberTaylor}
L. Brown,  B. Schreiber, and B.A. Taylor, \emph{Spectral synthesis and the Pompeiu problem},
 Ann. Inst. Fourier,  \textbf{23} (3) (1973),  125--154. 

\bibitem{brualdigibson1}
Richard~A. Brualdi and Peter~M. Gibson, \emph{Convex polyhedra of doubly
stochastic matrices. {I}. {A}pplications of the permanent function}, J.
Combinatorial Theory Ser. A \textbf{22} (1977), no.~2, 194--230.



\bibitem{Buhmann}
M. D. Buhmann,  \emph{Radial functions on compact support},  
Proceedings of the Edinburgh Mathematical Society. Series II,
\textbf{41}, (1998), no.~1, 33--46. 
{https://doi.org/10.1017/S0013091500019416}




\bibitem{camengasolidangles}
Kristin~A. Camenga, \emph{Vector spaces spanned by the angle sums of
polytopes}, Beitr\"age Algebra Geom. \textbf{47} (2006), no.~2, 447--462,
{\tt arXiv:math.MG/0508629}.


\bibitem{Car}
Leonard Carlitz, \emph{The reciprocity theorem for Dedekind-Rademacher sums},
Acta Math. XXIX, (1976), 309--313.

\bibitem{CasselsBook}
J. W. S.  Cassels,  \emph{An introduction to the geometry of numbers}, 
Corrected reprint of the 1971 edition, Classics in Mathematics, Springer-Verlag, Berlin, (1997), 
1--344.

 \bibitem{BillChen}
 William Y. C. Chen  and Peter L. Guo, \emph{Equivalence classes of full-dimensional $0/1$-polytopes with many vertices}, 
 Discrete and Computational Geometry \textbf{52} (2) (2014), 630--662.   

\bibitem{Cohn.etal}
Henry Cohn, Abhinav Kumar, Stephen D. Miller, Danylo Radchenko and Maryna Viazovska, 
\emph{The sphere packing problem in dimension $24$}, Annals of Mathematics 
SECOND SERIES, vol. 185, No. 3 (2017), 1017--1033.

\bibitem{Cohn-Elkies}
Henry Cohn and Noam Elkies, \emph{New upper bounds on sphere packings I},
Annals of Mathematics, \textbf{157} (2003), 689--714.

\bibitem{ConwaySloan.book}
John H. Conway and Neil  J. A. Sloane, \emph{Sphere Packings, Lattices and Groups}, 
Third Edition, with additional contributions by E. Bannai, R. E. Borcherds, J. Leech, S. P. Norton,  A. M. Odlyzko, R. A. Parker, L. Queen and B.  B. Venkov Grundlehren der Mathematischen Wissenschaften [Fundamental Principles of Mathematical Sciences], 290, Springer-Verlag, New York, (1999), 1--703.

\bibitem{Conway.Book.SensualForm}
John H. Conway, \emph{The sensual quadratic form}, The Carus Mathematical Monographs, published by the MAA, (1997), 1--152.

\bibitem{ConwayThompson}
John Milnor and Dale Husemoller, 
\emph{Symmetric bilinear forms},
in the series Ergebnisse der Mathematik und ihrer Grenzgebiete,
volume 73, 1973, 1--155.

\bibitem{Cordoba1}
Antonio Cordoba, 
\emph{La formule sommatoire de Poisson},
 C.R. Acad Sci. Paris,  306, Serie I,  (1988), 373--376.

\bibitem{CLS}
Dan Cristofaro-Gardiner, Teresa Xueshan Li,  and Richard Stanley,
\emph{New examples of period collapse}, (2015).
arXiv:1509.01887v1



\bibitem{Danilov}
Vladimir~I. Danilov, \emph{The geometry of toric varieties}, Uspekhi Mat. Nauk
\textbf{33} (1978), 85--134, 247.

\bibitem{DeBruijn.Book}
N. G. De Bruijn,  \emph{Filling boxes with bricks},  Amer. Math. Monthly \textbf{76} (1969), 37--40.

\bibitem{DeligneTabachnikovRobins}
Pierre Deligne, Sergei Tabachnikov,  and Sinai Robins, \emph{The Ice Cube Proof}, 
The Mathematical Intelligencer, Vol 36, no. 4,  2014, 1--3. 

\bibitem{deloerahemmeckekoeppe}
Jes{\'u}s~A. De~Loera, Raymond Hemmecke, and Matthias K{\"o}ppe,
\emph{Algebraic and {G}eometric {I}deas in the {T}heory of {D}iscrete
{O}ptimization}, MOS-SIAM Series on Optimization, vol.~14, Society for
Industrial and Applied Mathematics (SIAM), Philadelphia, PA; Mathematical
Optimization Society, Philadelphia, PA, 2013.

\bibitem{DRS}
Jes\'us de Loera, J\"org Rambau, and Francisco Santos,
\emph{ Triangulations.  Structures for algorithms and applications},
Algorithms and Computation in Mathematics, 25. 
Springer-Verlag, Berlin, (2010), 1--535.
 
\bibitem{DesarioRobins}
David Desario and Sinai Robins , \emph{Generalized solid-angle theory for real
polytopes}, The Quarterly Journal of Mathematics, \textbf{62} (2011), no.~4, 1003--1015, {\tt
arXiv:0708.0042}.

\bibitem{diaconisgangoli}
Persi Diaconis and Anil Gangolli, \emph{Rectangular arrays with fixed margins},
Discrete Probability and Algorithms (Minneapolis, MN, 1993), Springer, New
York, 1995, pp.~15--41.

\bibitem{diaz}
Ricardo Diaz and Sinai Robins , \emph{The {E}hrhart polynomial of a lattice polytope}, Annals of Math.
(2) \textbf{145} (1997), no.~3, 503--518.

\bibitem{RicardoNhatSinai}
Ricardo Diaz, Quang-Nhat Le and Sinai Robins, \emph{Fourier transforms of polytopes, solid angle sums, and discrete volumes}, preprint.
\url{https://drive.google.com/file/d/0B223XJaVpyE_MU16UER2VnFQRHc/view}

\bibitem{Dutour et al.parallelohedra}
 Dutour Sikiri\'c M,    Garber A,    Sch\"urmann A,    Waldmann C,
 \emph{The complete classification of five-dimensional Dirichlet-Voronoi polyhedra of translational lattices},  Acta Crystallogr A,  Found Adv.   Nov 1,  (2016) 72 (Pt 6),  673--683.
 
 
 
 \bibitem{DymMcKean}
H.  Dym and  H. P. McKean, \emph{Fourier Series and Integrals},  Academic Press, (1972), 1--295.
 
\bibitem{Dyson} 
Freeman Dyson, \emph{On the product of four non-homogeneous linear forms}, 
 Annals of Math. (2), \textbf{49} (1948), 82--109.
 
\bibitem{Ehrhart1}
Eug{\`e}ne Ehrhart, \emph{Sur les poly\`edres rationnels homoth\'etiques \`a
{$n$}\ dimensions}, C. R. Acad. Sci. Paris \textbf{254} (1962), 616--618.

\bibitem{Ehrhart2}
Eug{\`e}ne Ehrhart,  \emph{Sur un probl\`{e}me de g\'{e}om\'{e}trie diophantienne lin\'{e}aire I},  J. reine. angew. Math. 226, (1967), 1--29.

\bibitem{Ehrhart3}
Eug{\`e}ne Ehrhart,  \emph{Sur un probl\`{e}me de g\'{e}om\'{e}trie diophantienne lin\'{e}aire II},  J. reine. angew. Math. 227, (1967), 25--49.

\bibitem{Ehrhartbook}
Eug{\`e}ne Ehrhart, \emph{Polyn\^omes arithm\'etiques et m\'ethode des} \emph{poly\`edres en
combinatoire}, Birkh\"auser Verlag, Basel, 1977, International Series of
Numerical Mathematics, Vol. 35.


\bibitem{EinsiedlerWardBook}
Manfred Einsiedler and Thomas Ward, \emph{Functional Analysis, spectral theory, and applications},  Springer GTM series, (2017), 1--614. 


\bibitem{Entezari09}
Alireza Entezari, Ramsay Dyer, and Torsten M\"{o}ller, 
\emph{From sphere packing to the theory of optimal lattice sampling}, in
{Mathematical foundations of scientific visualization, computer
   graphics, and massive data exploration},  Series Math. Vis., Springer, Berlin, (2009),  227--255. 


\bibitem{EpsteinBook}
Charles L. Epstein, \emph{Introduction  to  the  mathematics  of  medical  imaging}, Society for Industrial and Applied Mathematics (SIAM), 
Philadelphia, PA, second edition, 2008.

\bibitem{KlainFeldman}
David  V.  Feldman  and Daniel A. Klain,
\emph{Angles as probabilities}, 
American Mathematical Monthly,  \textbf{116} (2009), no. 8, 732--735.

\bibitem{FeldmanProppRobins}
David Feldman, Jim Propp, and Sinai Robins, \emph{Tiling lattices with sublattices I}, 
Discrete \& Computational Geometry, Vol. 46, No. 1, (2011), 184--186.  

\bibitem{FischerPommersheim}
Benjamin Fischer and Jamie Pommersheim,  \emph{An algebraic construction of sum-integral interpolators}, preprint, 2021. 
\url{https://arxiv.org/abs/2101.04845v1}

\bibitem{Folland}
Gerald Folland, \emph{Fourier analysis and its applications}, Wadsworth \& Brooks/Cole Advanced Books \& Software, (1992), 1--433.

\bibitem{Fuglede74}
Bent Fuglede,
\emph{Commuting self-adjoint partial differential operators and a group theoretic problem},
J. Functional Analysis,
\textbf{16}, (1974), 101--121. 




\bibitem{FukshanskyBook}
Lenny Fukshansky and Stephan Ramon Garcia,  
 \emph{Geometric Number Theory}, Cambridge University Press, to appear in 2023. 

\bibitem{fukshanskysolidangle}
Lenny Fukshansky and Sinai Robins, \emph{Bounds for solid angles of lattices of
rank three}, J. Combin. Theory Ser. A \textbf{118} (2011), no.~2, 690--701,
{\tt arXiv:1006.0743}.


\bibitem{FultonBook}
William Fulton, \emph{Introduction to {T}oric {V}arieties}, Annals of Mathematics Studies, vol. 131, Princeton University Press, Princeton, NJ, 1993.

\bibitem{GaroufalidisPommersheim}
Stavros Garoufalidis and James Pommersheim, 
\emph{Sum-Integral interpolators and the Euler-MacLaurin formula for polytopes},
Transactions of the AMS, \textbf{364}, Number 6, June 2012, 2933--2958

\bibitem{GrepstadLev}
Sigrid Grepstad and Nir Lev,  \emph{Multi-tiling and Riesz bases}, Advances in Mathematics,
 \textbf{252}(2014),  1--6.

\bibitem{polymake}
Ewgenij Gawrilow and Michael Joswig, \emph{polymake: a framework for analyzing
convex polytopes}, Polytopes---combinatorics and computation (Oberwolfach,
1997), DMV Sem., vol.~29, Birkh\"auser, Basel, 2000, pp.~43--73, Software
{\tt polymake} available at  \url{https://www.polymake.org/doku.php}.




\bibitem{GravinKolountzakisRobinsShiryaev}
Nick Gravin, Mihail Kolountzakis, Sinai Robins, and Dmitry Shiryaev, 
\emph{Structure results for multiple tilings in 3D}, Discrete and Computational Geometry, (2013), Vol. 50, 1033--1050.  

\bibitem{Golub}
Gene H. Golub and Charles F. Van Loan, 
\emph{Matrix computations}, 4th ed. JHU Press, (2013).

\bibitem{GourionSeeger}
D.  Gourion, D. and A.  Seeger,   
\emph{Deterministic and stochastic methods for computing volumetric moduli of convex cones},  Comput. Appl. Math. \textbf{29},  (2010), 215--246.

\bibitem{GruberBook}
Peter M. Gruber, 
\emph{Convex and Discrete Geometry}, SpringerVerlag, Berlin, (2007), 1--590.

\bibitem{GruberLekkerkerker}
Peter M. Gruber and Cornelis G. Lekkerkerker, 
\emph{Geometry of numbers},
North-Holland mathematical library,  Elsevier science publishers,
(1987), 1--732.

\bibitem{GravinShiryaevRobins}
Nick Gravin, Sinai Robins, and Dmitry Shiryaev, 
\emph{Translational tilings by a polytope, with multiplicity}, 
Combinatorica \textbf{32} (2012), no.~6,
629--649, {\tt arXiv:1103.3163}.








\bibitem{Grunbaum}
 Branko Gr\"unbaum, \emph{Convex Polytopes}, Graduate Texts in Mathematics, vol. 221, Springer-Verlag, New York, 2003. Second edition, prepared by V. Kaibel, V. Klee, and G. M. Ziegler (original
edition: Interscience, London, 1967).

\bibitem{Grunbaum2}
 Branko Gr\"unbaum, \emph{Are your polyhedra the same as my polyhedra?}, 
 in Aronov, Boris; Basu, Saugata; Pach, J\'anos; Sharir, Micha (eds.), Discrete and Computational Geometry: The Goodman-Pollack Festschrift,  Algorithms and Combinatorics, (2003).

\bibitem{GunnellsSczech}
Paul Gunnells and Robert Sczech, \emph{Evaluation of Dedekind sums, Eisenstein cocycles, and special values of L-functions}, 
Duke Math. J. \textbf{118} (2003), no.~2, 229--260.

\bibitem{Gutierrez.Jimenez.Villa}
David Alonso-Guti\'errez, C. Hugo Jim\'enez, Rafael Villa, 
\emph{Brunn--Minkowski and Zhang inequalities for convolution bodies},
 Advances in Mathematics, \textbf{238} (2013),  50--69.
 
\bibitem{Zong2022}
Mei Han, Kirati Sriamon, Qi Yang, Chuanming Zong, 
\emph{Characterization of the three-dimensional multiple translative tiles},
Advances in Mathematics,  Volume 410, Part B, (2022).

\bibitem{Hales}
Thomas C. Hales, \emph{A proof of the Kepler conjecture}, 
 Ann. of Math. (2) 162 (2005), no. 3, 1065--1185. 

 \bibitem{Hardy41}
Godfrey~H. Hardy, \emph{Notes on special systems of orthogonal functions (IV): the orthogonal functions of Whittaker's cardinal series},  
Mathematical Proceedings of the Cambridge Philosophical Society, \textbf{37} (1941), 331--348.

https://doi:10.1017/S0305004100017977

\bibitem{Hardy.uncertainty}
Godfrey~H. Hardy, \emph{A Theorem Concerning Fourier Transforms},  
Journal of the London Mathematical Society, s1-8(3), (1933), 227--231.

\bibitem{HardyLittlewood}   
Godfrey~H. Hardy and John~E. Littlewood,  \emph{Some problems of Diophantine approximation: The lattice-points of a right-angled triangle (Second memoir)},   Abh. Math. Sem. Univ. Hamburg, 
no. 1 (1922), no. 1, 211--248.


\bibitem{Hardy.Wright.book}
Godfrey~H. Hardy  and Edward M.  Wright, \emph{Introduction to the theory of numbers}, 
Sixth edition, revised by D. R. Heath-Brown and J. H. Silverman, Oxford University Press (2008), 1--621.

\bibitem{Henk2}
Martin Henk, \emph{Inequalities between successive minima and intrinsic volumes of a convex body}, 
Monatsh. Math. \textbf{110} (1990) 279--282.

\bibitem{Henk4}
Martin Henk, \emph{Successive minima and lattice points},
Rend. Circ. Mat. Palermo (2) Suppl. no. 70, part I (2002), 377--384. 

\bibitem{HenkSchurmannWills}
Martin Henk, Achill Sch{\"u}rmann, and J{\"o}rg~M. Wills, \emph{Ehrhart
polynomials and successive minima}, Mathematika \textbf{52} (2005), no.~1--2,
1--16.  {\tt arXiv:math.MG/0507528} 
%

\bibitem{HenkWills}
Martin Henk and J{\"o}rg~M. Wills,
\emph{A Blichfeldt-type inequality for the surface area},
Monatsh Math \textbf{154},  (2008), 135--144.



\bibitem{HenkHenzWills}
Martin Henk, Matthias Henze, and J{\"o}rg~M. Wills, 
\emph{Blichfeldt-type inequalities and central symmetry},
Advances in geometry, \textbf{11} (2011), 731--744.
 
\bibitem{Henk3}
Martin Henk, \emph{An introduction to geometry of numbers}, lecture notes, preprint. 

\bibitem{Hensley}
Douglas Hensley,  \emph{Lattice vertex polytopes with interior lattice points},
Pacific Journal of Mathematics,  \textbf{105} (1983), no. 1,  183--191.


\bibitem{Herstein}
Israel Nathan Herstein, \emph{Topics in Algebra},
 (2nd ed.) 1975, Wiley \& sons,  1--388.

\bibitem{Higgins1996}
 John R. Higgins,  \emph{Sampling Theory in Fourier and Signal Analysis}, Clarendon Press, Oxford,
 (1996) 1-- 222.  
     
\bibitem{HlawkaSchoissengeierTaschner}
Edmund Hlawka , Rudolf Taschner , and Johannes Schoißengeier,
\emph{Geometric and Analytic Number Theory}, 
Springer Universitext (UTX), (1991), 1--238.

\bibitem{Hlawka}
Edmund Hlawka, \emph{Zur Geometrie der Zahlen}, Math. Z. \textbf{49} (1943), 285--312.
     
\bibitem{IosevichKatzTao}
Alex Iosevich, Nets Katz, and Terence Tao, \emph{The Fuglede spectral conjecture holds for convex planar domains},  Mathematical Research Letters \text{10},  (2003) 559--569.

\bibitem{KatharinaJochemko}
Katharina Jochemko,
\emph{A Brief Introduction to Valuations on Lattice Polytopes}, 
Algebraic and Geometric Combinatorics on Lattice Polytopes. Proceedings of the Summer Workshop on Lattice Polytopes. Hibi, T., Tsuchiya A. (eds), 
38--55, (2019), World Sci. Publ.

\bibitem{kannan}
Ravi Kannan, \emph{Lattice translates of a polytope and the {F}robenius problem}, 
Combinatorica \textbf{12} (1992), no.~2, 161--177.







\bibitem{kantorkhovanskii}
Jean-Michel Kantor and Askold~G. Khovanski{\u\i}, \emph{Une application du
th\'eor\`eme de {R}iemann--{R}och combinatoire au polyn\^ome d'{E}hrhart des
polytopes entiers de $\mathbf {R}\sp d$}, C. R. Acad. Sci. Paris S\'er. I
Math. \textbf{317} (1993), no.~5, 501--507.

\bibitem{Karasev.et.al}
Roman Karasev, Jan Kyncl, Pavel Paták, Zuzana Patáková, and Martin Tancer, 
\emph{Bounds for Pach’s Selection Theorem and for the
Minimum Solid Angle in a Simplex},
Discrete Comput Geom (2015) \textbf{54}, 610--636.

\bibitem{KarshonSternbergWeitsman1}
Yael Karshon, Shlomo Sternberg, and Jonathan Weitsman, \emph{The
{E}uler--{M}aclaurin formula for simple integral polytopes}, Proc. Natl.
Acad. Sci. USA \textbf{100} (2003), no.~2, 426--433.

\bibitem{KarshonSternbergWeitsman2}
Yael Karshon, Shlomo Sternberg,  and Jonathan Weitsman,  
\emph{Exact Euler-Maclaurin formulas for simple lattice polytopes}, Advances in Applied Mathematics, (2007), Vol 39 (1), 1--50.

\bibitem{Kathuria}
Leetika Kathuria and Madhu Raka, 
\emph{On conjectures of Minkowski and Woods for n=10},
Proc. Indian Acad. Sci. (Math. Sci.) 132:45 (2022).

\bibitem{Katznelson}
Yitzhak Katznelson,  \emph{An introduction to harmonic analysis}, Third edition, Cambridge Mathematical Library, Cambridge University Press, Cambridge, (2004), 1--314.

\bibitem{Katz.Stapledon}
Eric Katz and Alan Stapledon, \emph{Local h-polynomials, invariants of subdivisions, and mixed Ehrhart theory}, Adv. Math., \textbf{286} (2016), 181--239.

\bibitem{Knuth}
Donald Knuth,  \emph{Notes on generalized Dedekind sums},  Acta Arith. (1977), 297--325.




\bibitem{Klain1}
Daniel A.  Klain,   \emph{The Minkowski problem for polytopes},  Advances in Mathematics,  \textbf{185}  (2004), no. 2, 270--288.         

\bibitem{CarolineKlivans}
Caroline J. Klivans, 
\emph{The Mathematics of Chip-firing},
Discrete Mathematics and its Applications, Taylor \& Francis Group, LLC, (2019), 1--295.

\bibitem{Kobayashi1}
 T. Kobayashi, 
 \emph{The null variety of the Fourier transform of the characteristic function of a bounded
domain}, Semin. Rep. Unitary Represent. 6 (1986), 1--18.

\bibitem{Kobayashi2}
 T. Kobayashi, 
\emph{Asymptotic behaviour of the null variety for a convex domain in a non-positively curved
space form}, 
J. Fac. Sci. Univ. Tokyo Sect. IA Math. 36 (1989), 389--478.


\bibitem{KoldobskyBook}
Alexander Koldobsky, \emph{Fourier analysis in convex geometry}, Mathematical Surveys and Monographs, American Mathematical Society, Providence, RI, 1--170.

\bibitem{Kolountzakis1}
Mihalis  N. Kolountzakis, \emph{On the structure of multiple translational tilings by polygonal regions}, 
 Discrete and Computational Geometry,  \textbf{23} (4),  (2000),  537--553.

\bibitem{Kolountzakis2}
Mihalis  N. Kolountzakis, \emph{The study of translational tilings with Fourier analysis},  in
 Fourier analysis and convexity, 
Appl. Numer. Harmon. Anal., Birkhäuser Boston, Boston, MA, (2004), 131--187.

\bibitem{KorkinZolotarev}
A. N. Korkin and E. I. Zolotarev, 
\emph{Sur les formes quadratiques positives quaternaires}, 
Math. Ann. \textbf{5} (1872), 581--583.

\bibitem{KorkinZolotarev5}
A. N. Korkin and E. I. Zolotarev, 
\emph{Sur les formes quadratique positives},
Math. Ann. \textbf{11} (1877), 242--292.

\bibitem{Kuperberg}
Greg Kuperberg, \emph{Notions of denseness}, Geom. Topol. 4 (2000) 277--292.

\bibitem{Kraft}
Roger L. Kraft, 
\emph{What's the difference between Cantor sets?},
American Mathematical Monthly, 
\textbf{101} (7), (1994), 640--650.

\bibitem{LagariasZiegler}
Jeffrey C. Lagarias and G\"unter M. Ziegler, 
\emph{Bounds for lattice polytopes containing a fixed number of interior points in a sublattice},
Canadian J. Math. \textbf{43}, (1991), no. 5, 1022--1035.




\bibitem{LagariasZong}
Jeffrey C.  Lagarias and  Chuanming  Zong, 
\emph{Mysteries in packing regular tetrahedra}, 
Notices Amer. Math. Soc. \textbf{59} (2012), no. 11, 1540--1549.

\bibitem{Lagrange}
Joseph-Louis Lagrange,  \emph{Recherches d’arithm\'etique}, 
Nouveaux M\'emoires de L'Acad\'emie royal des Sciences et Belles-Lettres de Berlin (1773), 265--312.


\bibitem{LasserreZeron}
Jean B. Lasserre and Eduardo S. Zeron, \emph{On counting integral points in a convex rational polytope}, 
Math. Oper. Res., 28(4)  (2003), 853--870.

\bibitem{LawrenceVolume2}
Jim Lawrence,
\emph{Rational-function-valued Valuations on Polyhedra},
 Discrete and Computational Geometry (1990), 199--208.

\bibitem{LawrenceVolume}
Jim Lawrence, \emph{Polytope volume computation}, Math. Comp. \textbf{57} (1991),
no.~195, 259--271.

\bibitem{LevLiu}
Nir Lev and  Bochen Liu,   Multi-tiling and equidecomposability of polytopes by lattice translates, Bulletin of the London Math. Society, Vol 51, issue 6, (2019), 1079--1098.

\bibitem{LevMatolcsi}
Nir Lev and M\'at\'e  Matolcsi,  \emph{The Fuglede conjecture for convex domains is true in all dimensions}, preprint, 2021.

\bibitem{Lighthill}
Lighthill, M. J., \emph{Introduction to Fourier analysis and generalised functions}, 
Cambridge University Press, New York (1960), 1--79. 

\bibitem{linke}
Eva Linke, \emph{Rational {E}hrhart quasi-polynomials}, J. Combin. Theory Ser.
A \textbf{118} (2011), no.~7, 1966--1978, {\tt arXiv:1006.5612}.

\bibitem{Lions}
J. L. Lions, \emph{Supports de produits de composition I}, Comptes Rendus 232, (1951) 1530--1532; 11, Comptes Rendus 232,  (1951), 1622--1624.

\bibitem{Liu}   
Bochen Liu, \emph{Periodic structure of translational multi-tilings in the plane}, 
American Journal of Mathematics, Volume 143, Number 6, December
(2021), 1841--1862.

\bibitem{macdonald1}
Ian~G. Macdonald,
\emph{The volume of a lattice polyhedron},  Proc. Cambridge Philos. Soc., \textbf{59} (1963), 719--726.

\bibitem{macdonald2}
Ian~G. Macdonald, \emph{Polynomials associated with finite cell-complexes}, J. London Math. Soc. (2) \textbf{4} (1971), 181--192.



\bibitem{FabricioSinai1}
Fabricio Caluza Machado and Sinai Robins,  \emph{The null set of a polytope and the Pompeiu property for polytopes}, to appear in Journal d'Analyse Mathematique, 2022.

\bibitem{FabricioSinai2}
Fabricio Caluza Machado and Sinai Robins, \emph{Coefficients of the solid angle and Ehrhart quasi-polynomials}, preprint, 2019.

\bibitem{Malikiosis1}
Romanos-Diogenes Malikiosis, 
\emph{A discrete analogue for Minkowski’s second theorem on successive minima},
Advances in Geometry \textbf{12} (2012), 365--380.

\bibitem{Martinet}
J. Martinet, \emph{Perfect lattices in Euclidean spaces}, 
Grundlehren der Mathematischen Wissenschaften, 
Fundamental Principles of Mathematical Sciences, vol. 327, Springer-Verlag,
Berlin, (2003).

\bibitem{MartinsRobins}
Michel Faleiros Martins and Sinai Robins,
\emph{The covariogram and extensions of the Bombieri-Siegel formula}, (2023), 
https://arxiv.org/abs/2204.08606

\bibitem{McAllisterWoods}
Tyrrell B. McAllister and Kevin M. Woods, \emph{The minimum period of the Ehrhart quasi-polynomial of a rational polytope}, Journal of Combinatorial Theory, Series A \textbf{109} (2005) 345--352.


\bibitem{CurtMcMullen}
Curtis McMullen, \emph{Minkowski’s conjecture, well-rounded lattices and topological dimension}, Journal of the American Mathematical Society \textbf{18}, (2005), revised 2007, 
711--734.



\bibitem{McMullen1}
Peter McMullen,
\emph{Lattice invariant valuations on rational polytopes},  Arch. Math., 31, (1978), 509--516.


\bibitem{McMullen2}
Peter McMullen, 
\emph{Non-linear angle-sum relations for polyhedral cones and polytopes},   Math. Proc. Cambridge Phil. Soc., 78, (1975), 247--261.


\bibitem{McMullen3} 
Peter McMullen, \emph{Angle-sum relations for polyhedral sets},  Mathematika \textbf{33} (1986), 
no. 2, 173--188.


\bibitem{McMullen5}
Peter McMullen. \emph{Valuations and Euler-Type Relations on Certain Classes of Convex Polytopes},
Proceedings of the London Mathematical Society 3.1 (1977), 113--135.


\bibitem{McMullen4} 
Peter McMullen, \emph{Polytopes with centrally symmetric faces}, Israel J. Math., 8 (1970), 194--196.

\bibitem{Mercer}
James Mercer, 
 \emph{Functions of positive and negative type and their connection with the theory of integral equations}, 
Philosophical Transactions of the Royal Society A, 209 (441–458): (1909), 415--446.

\bibitem{Minkowski} 
Hermann Minkowski,  \emph{Geometrie der Zahlen}, Teubner, Leipzig, 1896.

\bibitem{Minkowski1897}
Hermann Minkowski, \emph{Allgemeine Lehrsatze iiber konvexen Polyeder}, 
Nachr. K. Akad. Wiss. Gottingen, Math.-Phys. Kl. ii (1897), 198--219.

\bibitem{DraismaMcAllisterNill}
Jan Draisma, Tyrrell B. McAllister, and Benjamin Nill, 
\emph{Lattice-Width Directions and Minkowski's $3^d$-Theorem}, 
SIAM Journal on Discrete Mathematics, Vol. 26, No. 3 (2012),  1104--1107.

\bibitem{morelli}
Robert Morelli, \emph{Pick's theorem and the {T}odd class of a toric variety},
Adv. Math. \textbf{100} (1993), no.~2, 183--231.

\bibitem{Mordell.convexity.measure}
Louis  J. Mordell, \emph{On some arithmetical results in the geometry of numbers},
Compositio Mathematica, \textbf{1} (1935), 248--253.


\bibitem{MorrisNewman}
Morris Newman, \emph{Integral Matrices}, Academic press, (1972), 1--223.


\bibitem{Randol1}
Marina Nechayeva and Burton Randol,  \emph{Asymptotics of weighted lattice point counts inside dilating polygons}, 
Additive number theory,  Springer, New York, (2010), 287--301.

\bibitem{Nill.and.Paffenholz}
Benjamin Nill and Andreas Paffenholz,
\emph{On the equality case in Ehrhart's volume conjecture},  Adv. Geom. 14 (2014), no. 4, 579--586.

\bibitem{Nosarzewska}
Maria Nosarzewska, 
\emph{\'Evaluation de la diff\'erence entre l'aire d'une r\'egion plane convexe et le nombre des points aux coordonn\'ees enti\`eres couvertes par elle}, Colloq. Math. \textbf{1} (1948), 305--311.

\bibitem{OldsBook}
C. D. Olds,  Anneli Lax, Giuliana P. Davidoff,  \emph{The geometry of numbers}, 
Mathematical Association of America, (2000), 1--193.

\bibitem{OsgoodBook}
Brad G. Osgood,  \emph{Lectures on the Fourier transform and its applications},
Pure and Applied Undergraduate Texts, \textbf{33}, 
American Mathematical Society, Providence, RI, 2019. 1--693.

\bibitem{payneehrharttriang}
Sam Payne, \emph{Ehrhart series and lattice triangulations}, Discrete Comput.
Geom. \textbf{40} (2008), no.~3, 365--376, {\tt arXiv:math/0702052}.


\bibitem{PetersenMiddleton62}
Daniel Petersen, and David Middleton,
\emph{Sampling and reconstruction of wave-number-limited functions in {$N$}-dimensional Euclidean spaces},
Information and Control,  \textbf{5}    (1962),  279--323.

\bibitem{perlesshephardangles}
Micha~A. Perles and Geoffrey~C. Shephard, \emph{Angle sums of convex
polytopes}, Math. Scand. \textbf{21} (1967), 199--218.

\bibitem{MarkPinsky}
Mark A. Pinsky, 
\emph{Introduction to Fourier Analysis and Wavelets}, 
Brooks/Cole, Pacific Grove, California, (2002), 1--376.


\bibitem{Pikhurko}
O. Pikhurko, 
\emph{Lattice points in lattice polytopes}, Mathematika 48 (2001), 
no. 1-2, 15--24.



\bibitem{Podkorytov}
A. N. Podkorytov and Mai Van Minh, 
\emph{The Fourier formula for discontinuous functions of several variables},
Journal of Mathematical Sciences, Vol. 124, No. 3, (2004), 5018--5025.
 
\bibitem{Pommersheim}
James~E. Pommersheim, 
\emph{Toric varieties, lattice points and {D}edekind
sums}, Math. Ann. \textbf{295} (1993), no.~1, 1--24.

\bibitem{PoonenVillegas}
Bjorn Poonen and Fernando Rodriguez-Villegas, 
\emph{Lattice polygons and the number 12}, Amer. Math. Monthly 107 (2000), no. 3, 238--250.






\bibitem{Postnikov.Permutahedra}
Alexander Postnikov, \emph{Permutohedra, associahedra, and beyond}, Int. Math. Res. Not. (2009), no.~6, 1026--1106, {\tt arXiv:math/0507163}.



\bibitem{Nhat}
Quang-Nhat Le,  \emph{A discrete Stokes formula and
the solid-angle sum of polytopes}, undergraduate dissertation, (2009).

\bibitem{Nhat.Sinai}
Quang-Nhat Le and Sinai Robins,  \emph{Macdonald's solid-angle sum for real dilations of rational polygons}, preprint.

\bibitem{Ramanujan1}  Srinivasa Ramanujan,  \emph{Some definite integrals}, Messenger of Mathematics \textbf{44} (1915), 10--18.

\bibitem{JorgeRami­rez}
Jorge Luis Ramirez Alfonsin,    \emph{Complexity of the Frobenius problem}, Combinatorica, 
\textbf{16} (1), (1996), 143--147/








\bibitem{Randol3}
Burton Randol, \emph{On the Fourier transform of the indicator function of a planar set},
Trans. Amer. Math. Soc., \textbf{139} (1969), 271--276.
 
\bibitem{Randol4}
Burton Randol, \emph{On the asymptotic behavior of the Fourier transform of a convex set},
Trans. Amer. Math. Soc., \textbf{139} (1969), 279--285.


\bibitem{Randol2}
Burton Randol, \emph{On the number of integral lattice-points in dilations of algebraic polyhedra}, 
Internat. Math. Res. Notices (1997) no. 6, 259--270.

\bibitem{RegevNotes}
Oded Regev,  \emph{Lattices in Computer Science}, 
Online lecture notes:  

\url{https://cims.nyu.edu/~regev/teaching/lattices_fall_2009/index.html}

\bibitem{Remak1}
Robert Remak, 
\emph{Vereinfachung eines Blichfeldtschen Beweises aus der Geomtrie der Zahlen},
Math. Zeitschr. \textbf{26}, (1927), 694--699.

\bibitem{Remak2}
Robert Remak, \emph{Verellgemeinerung eines Minkowskischen Satzes, I, II},  
Mathematische Zeitschrift, \textbf{17} (1923), 1--34;  \textbf{18} (1923), 173--200.

\bibitem{Reznick}
Bruce Reznick, \emph{Lattice point simplices}, Discrete Math. \textbf{60}  (1986), 219--242.



\bibitem{ribando}
Jason~M. Ribando, \emph{Measuring solid angles beyond dimension three},
Discrete Comput. Geom. \textbf{36} (2006), no. 3, 479--487.

\bibitem{RogersBook}
C.  A. Rogers, Packing and covering,  Cambridge Tracts in Mathematics and Mathematical Physics, No. 54, Cambridge University Press, New York, 1964.





\bibitem{RogersShephard}
C. A. Rogers and  G. C. Shephard, 
\emph{Convex bodies associated with a given convex body}, 
Journal of the London Math. Soc., 1 (1958), no. 3, 270--281.

\bibitem{RudinGreenBook}
Walter Rudin, \emph{Real and complex analysis},  Third edition,  McGraw-Hill Book Co., New York, (1987), 
1--416.

\bibitem{RudinGroups}
Walter Rudin, \emph{Fourier analysis on groups}, Wiley Classics Library, (1990), 1--285.
 
\bibitem{Tiago1}
Tiago Royer, \emph{Reconstruction of rational polytopes from the real-parameter Ehrhart function of its translates}, preprint, 2017.
\url{https://arxiv.org/abs/1712.01973} 

\bibitem{Tiago2}
Tiago Royer, \emph{Reconstruction of symmetric convex bodies from Ehrhart-like data}, preprint, 2017.
\url{https://arxiv.org/abs/1712.03937}

\bibitem{GervasioSantos}
Gerv\'asio Prot\'asio dos Santo Neto, \emph{The theory and computation of solid angles}, Master's thesis, IME, Universidade de S\~ao Paulo, (2021), 1--85. 

\bibitem{PaulSally1}
Paul Sally, \emph{Fundamentals of Mathematical Analysis}, AMS (The Sally series), Pure and Applied Undergraduate texts (20), (2013), 1--384.


\bibitem{Schiemann1}
Alexander Schiemann,   \emph{Ein Beispiel positiv definiter quadratischer Formen der Dimension 4 mit gleichen Darstellungszahlen},
  Arch. Math. 54 (1990), 372--375.

\bibitem{Schiemann2}
 Alexander Schiemann,  \emph{Temare positiv defInite quadratische Fonnen mit gleichen Darstellungszahlen}, Dissertation, Bonn, 1993. 
 
\bibitem{schlafli}
Ludwig Schl\"afli, \emph{Theorie der vielfachen {K}ontinuit\"at}, Ludwig
{S}chl\"afli, 1814--1895, {G}esammelte {M}athematische {A}bhandlungen,
{V}ol.~I, Birkh\"auser, Basel, 1950, pp.~167--387.

\bibitem{SchleimerSegerman}
Saul Schleimer  and Henry Segerman,  \emph{Puzzling the 120-cell},  Notices Amer. Math. Soc. \textbf{62} (2015), no. 11, 1309--1316.

\bibitem{SchmeisserSickel2000}
Hans--Jürgen Schmeisser, and Winfried Sickel,
\emph{Sampling theory and function spaces}, 
Applied Mathematics Reviews,  (2000), 205--284.

\bibitem{Schneider.book}
Rolf Schneider, \emph{Convex Bodies: The Brunn--Minkowski Theory}, 2nd edition, Encyclopedia of Mathematics and its Applications, Cambridge University Press, 2013. 

\bibitem{Schneider.book2}
Rolf Schneider and Wolfgang Weil, \emph{Stochastic and integral geometry}, 
Springer Science \& Business Media,  2008.



\bibitem{schrijver}
Alexander Schrijver, \emph{Combinatorial {O}ptimization. {P}olyhedra and
{E}fficiency. {V}ol. {A}--{C}}, Algorithms and Combinatorics, vol.~24,
Springer-Verlag, Berlin, 2003.


\bibitem{SenechalGaliulin}
Marjorie Senechal and  R.V. Galiulin,  \emph{An Introduction to the Theory of Figures: the Geometry of E.S. Fedorov}, Structural Topology, \textbf{10},  (1984),  5--22.

\bibitem{shallitfrobgeneralizations}
Jeffrey Shallit, \emph{The {F}robenius problem and its generalizations},
Developments in language theory, Lecture Notes in Comput. Sci., vol. 5257,
Springer, Berlin, 2008, pp.~72--83.


\bibitem{Shannon1}
Claude E. Shannon,   \emph{Communication in the Presence of Noise}, 
   Proceedings of the IRE,  \textbf{37}, number 1, {1949}, 10--21.

\bibitem{ShephardSymmetricPolytopes}
Geoffrey C. Shephard, \emph{Polytopes with centrally symmetric faces}, 
Canadian J. Math., \textbf{19} (1967), 1206--1213.

\bibitem{ShephardGramRelations}
Geoffrey~C. Shephard, \emph{An elementary proof of {G}ram's theorem for convex
polytopes}, Canad. J. Math. \textbf{19} (1967), 1214--1217.


\bibitem{SiegelBook}
Carl~Ludwig Siegel, \emph{Lectures on the {G}eometry of {N}umbers},
Springer-Verlag, Berlin, 1989, Notes by B. Friedman, rewritten by Komaravolu
Chandrasekharan with the assistance of Rudolf Suter, with a preface by
Chandrasekharan.

\bibitem{Siegel1}
Carl~Ludwig Siegel, \emph{A Mean Value Theorem in Geometry of Numbers}, 
Annals of Mathematics,  Second Series, Vol. 46, No. 2 (1945), 340--347.

\bibitem{Skriganov1}
Maxim M.  Skriganov,  
\emph{Ergodic theory on homogeneous spaces and the enumeration of lattice points in
polyhedra (Russian)}, Dokl. Akad. Nauk \textbf{355} (1997), no. 5, 609--611.

\bibitem{Skriganov2}
Maxim M.  Skriganov, 
\emph{Ergodic theory on SL(n), Diophantine approximations and anomalies in the lattice point problem}, Invent. Math. \textbf{132} (1998), no. 1, 1--72.


\bibitem{Skubenko}
B. F. Skubenko, 
\emph{A proof of Minkowski’s conjecture on the product of $n$ linear inhomogeneous forms in $n$ variables 
for $n \leq 5$},  J. Soviet Math. 6 (1976), 627--650; Proc. Steklov Inst. Math. 33 (1973),4--36.

\bibitem{sommerville}
Duncan M.~Y. Sommerville, \emph{The relation connecting the angle-sums and
volume of a polytope in space of $n$ dimensions}, Proc. Roy. Soc. London,
Ser. A \textbf{115} (1927), 103--119.


\bibitem{Stanleyreciprocity}
Richard~P. Stanley,   \emph{Combinatorial reciprocity theorems}, Advances in Math.
\textbf{14} (1974), 194--253.

\bibitem{StanleyBook}
Richard~P. Stanley, \emph{Enumerative Combinatorics, Volume 1}, Second edition, 
Cambridge Studies in Advanced Mathematics, vol.~49, 
Cambridge University Press, Cambridge, 2012.

\bibitem{StanleyDecompositions}
Richard~P. Stanley, \emph{Decompositions of rational convex polytopes}, Ann. Discrete Math. \textbf{6} (1980), 333--342.

\bibitem{StanleyCommutativeAlgebra}
Richard~P. Stanley,    \emph{Combinatorics and Commutative Algebra}, Second edition, 
Progress in Mathematics, vol. 41, Birkha\"user Boston Inc., Boston, MA, 1996.

\bibitem{stapledonadditive}
Alan Stapledon, \emph{Additive number theory and inequalities in Ehrhart
theory}, 
International Mathematics Research Notices, No. 5, (2016), 1497--1540.

\bibitem{SteinWeiss}
Elias Stein and Guido Weiss, \emph{Introduction to Fourier analysis on Euclidean spaces}, 
Princeton University Press, Princeton Mathematical Series, No. 32, Princeton, N.J., 1971.

\bibitem{SteinShakarchi}
Elias 
Stein and Rami Shakarchi,     
\emph{Fourier analysis, an introduction}, 
Princeton Lectures in Analysis, 1. Princeton University Press, Princeton, NJ, 2003, 1--311.

\bibitem{Sturmfels}
Berndt Sturmfels,  \emph{On vector partition functions}, 
Journal of Combinatorial Theory Series A. \textbf{72} (2)  (1995), 302--309.



\bibitem{Thue}
Axel Thue,  \emph{\"Uber die dichteste Zuzammenstellung von kongruenten Kreisen in der Ebene}, Norske Vid. Selsk. Skr. 1 (1910), 1--9.

\bibitem{Titchmarsh2}
E. C. Titchmarsh,  \emph{The zeros of certain integral functions}, 
Proc. Lond. Math. Soc. \text{25},  (1926), 283--302.



\bibitem{Titchmarsh}
E. C. Titchmarsh, 
\emph{Introduction to the Theory of Fourier Integrals},
Oxford University Press, Oxford, (1937).

\bibitem{Toth1}
L\'aszl\'o Fejes T\'oth,   \emph{Some packing and covering theorems},
 Acta Sci. Math. 12A (1950), 62--67.

\bibitem{Travaglini}
Giancarlo  Travaglini,  
\emph{Number theory, Fourier analysis and geometric discrepancy}, 
London Mathematical Society Student Texts, 81. Cambridge University Press, Cambridge (2014), 1--240.

\bibitem{Terras.Harmonic1}
Audrey Terras, \emph{Harmonic Analysis on Symmetric Spaces and Applications I},
Springer New York, NY,  (1985), 1--341.

\bibitem{Terras}
Audrey Terras, \emph{Fourier Analysis on Finite Groups and Applications}, London Mathematical Society,
Student Texts, vol. 43, Cambridge University Press, Cambridge, (1999).

\bibitem{Unser00}
Michael Unser,  
\emph{Sampling - 50 years after Shannon}, Proceedings of the IEEE, \textbf{88}, issue 4, (2000), 569--587.

\bibitem{Van.der.Corput.1}
J. G. van der Corput, 
\emph{Verallgemeinerung einer Mordellschen Beweis-methode in der Geometrie der Zahlen}, 
Acta Arithmetica \textbf{2}  (1936 (a)), 145--146.

\bibitem{Venkatesh}
 Akshay Venkatesh, \emph{A note on sphere packings in high dimension},  Int. Math. Res. Not.  IMRN (2013), 1628--1642.

\bibitem{verdoolaege}
Sven Verdoolaege, \emph{Software package {\tt barvinok}}, (2004),
electronically available at

 \url{http://freshmeat.net/projects/barvinok/}.

\bibitem{StanWagon}
Stan Wagon, \emph{Fourteen Proofs of a result about tiling a rectangle},
The American Mathematical Monthly \textbf{94} (1987), 601--617.

\bibitem{White.G.K}
G. K. White,  
\emph{A Refinement of Van Der Corput’s Theorem on Convex Bodies},
American Journal of Mathematics vol. 85, no. 2 (1963), 320--26. 

\bibitem{Wills1}
J. M. Wills, \emph{\"Uber konvexe Gitterpolygone}, Comment. Math. Helv. 
\textbf{48} (1973) 188--194.





\bibitem{Woods}
Kevin Woods,  \emph{The unreasonable ubiquitousness of quasi-polynomials}, 
 Electronic Journal of Combinatorics \textbf{21} (1), 
Paper 1.44,   (2014), 1--23.


\bibitem{YangZong}
Qi Yang and Chuanming Zong,  \emph{Multiple lattice tilings in Euclidean spaces}, Canad. Math. Bull. 
\textbf{62}  (2019), no. 4, 923--929. 

\bibitem{YauZhang}
Stephen T. Yau and Letian Zhang, \emph{An upper estimate of integral points in 
real simplices with an application to singularity theory}, Math. Res. Lett. \textbf{13} (2006), no. 6, 911--921.





\bibitem{Ye12}
Wenxing Ye,  and Alireza Entezari,  \emph{A geometric construction of multivariate sinc functions}, 
IEEE Transactions on Image Processing,
\textbf{21} (2012), no. 6,  2969--2979.

\bibitem{Ziegler}
G\"unter M. Ziegler,  \emph{Lectures on polytopes}, Graduate Texts in Mathematics, Volume 152, 
Springer-Verlag, New York, 1995.

\bibitem{Zong.lattices}
Chuanming Zong, 
\emph{Classification of the sublattices of a lattice}, 
Aust. Math. Soc. (2020), 1--12.

\bibitem{Zong.book}
Chuanming Zong, \emph{The cube - a window to convex and discrete geometry}, Cambridge University Press, (2006), 
1--174.

\bibitem{Zong.SpherePackingsBook}
Chuanming Zong, \emph{Sphere Packings}, Springer Universitext, (1999), 
1--245.

\bibitem{Zong.SpherePackingsSurveyPaper}
Chuanming Zong, 
\emph{Packing, covering and tiling in two-dimensional spaces},
Expo. Math. 32 (2014) 297--364.

\bibitem{Zygmund}
Antoni Zygmund,  \emph{Trigonometric Series},  Third edition, Volumes I $\&$ II combined, with a forward by Robert Fefferman, Cambridge University Press, (2002).




\end{thebibliography}
\end{document}